\pdfoutput=1
\documentclass[a4paper]{report}

\usepackage[utf8]{inputenc}
\usepackage{amsmath,amssymb}
\usepackage{graphicx}
\usepackage{xcolor}
\usepackage{geometry}
\usepackage{microtype}
\usepackage{tikz}
\usetikzlibrary{calc,angles,quotes,arrows.meta}
\usepackage{makeidx}
\usepackage[square,sort,comma,numbers]{natbib}
\usepackage[hidelinks,bookmarksnumbered]{hyperref}
\counterwithout{section}{chapter}

\newcommand{\eof}{{\unskip\nobreak\hfil\penalty50
          \hskip2em\hbox{}\nobreak\hfil\mbox{\rule{1ex}{1ex} \qquad}
   \parfillskip=0pt
   \finalhyphendemerits=0\par\medskip}}

\def\a{\alpha}
\def\abs#1{\vert  #1 \vert}        
\def\act#1#2{\langle#1,#2\rangle}  
\def\ae{ {\textrm{~a.e.~}} }

\def\Ab{\tilde{A}}
\def\AT{{\mathcal{A}_T}}

\def\b{\beta}
\def\ba{\left\{\begin{array}{llll}   }
\def\bac{\begin{array}{c}}
\def\bal{\begin{array}{l}}
\def\bari{\begin{array}{r}}
\def\ban{\begin{array}{llll}}
\def\be{\begin{equation}}
\def\benum{\begin{enumerate}}

\def\card{\mathrm{card}}
\def\cedge{c_\edge}
\def\centercv{x_\cv}                         
\def\centercvv{x_\cvv}                       
\def\centercvw{x_\cvw}                       
\def\centeredge{y_\edge}                     
\def\centeredget{y_\edget}                     
\def\comp{{\bar{\o}}}
\def\cv{K}                                    
\def\cvbar{{\overline \cv}}                   
\def\cvedgem{{\cv_{\edge,-}}}
\def\cvedgep{{\cv_{\edge,+}}}
\def\cvv{{L}}                                 
\def\cvw{{M}}                                 
\def\cvvbar{{\overline \cvv}}                 
\def\cvedge{{\cv_\edge}}

\def\d{\delta}
\def\dcvedge{ {d_{\cv,\edge}} }    
\def\dcvvedge{ {d_{\cvv,\edge}} }  
\def\dedge{d_\edge}                
\def\dcvcvv{d_{\edgecvcvv}}        
\def\dedgecvcvv{d_{\edgecvcvv}}    
\def\dfrontiere{{d \gamma}}
\def\diam{\mathrm{diam}}
\def\diamO{\mathrm{diam}(\O)}        
\def\dint{\int_\R \int_{\R_+}}    
\def\dim{d}
\def\div{\mathrm{div}}
\def\dO{\partial \Omega}
\def\dr{\partial}
\def\dsp{\displaystyle} 
\def\dsurdn#1{\frac {\dr #1} {\dr n}} 
\def\dxi#1{\frac {\dr #1} {\dr x_i}}
\def\dxj#1{\frac {\dr #1} {\dr x_j}}

\def\D{\mathcal{D}}
\def\Dcvcvv{\mathcal{D}_{\cv, \cvv}}   
\def\Dcvedge{\mathcal{D}_{\cv, \edge}} 

\def\ea{\end{array}\right.}
\def\ean{\end{array}}
\def\ecv{{e_{\cv}}}
\def\edge{\sigma}
\def\edget{{\tilde \sigma}}
\def\edgecvcvv{{\cv\vert\cvv}}
\def\edges{\mathcal{E}}              
\def\edgesatyp{\mathcal{E_a}}        
\def\edgescv{{\edges}_\cv}        
\def\edgescvv{{\edges}_\cvv}        
\def\edgescvedge{{\edges}_\cvedge} 
\def\edgesint{{\edges}_\mathrm{int}} 
\def\edgesext{{\edges}_\mathrm{ext}} 
\def\ee{\end{equation}}
\def\eedgep{e_{\edge,+}}
\def\eenum{\end{enumerate}}
\def\ein{{e_i^{n}}}
\def\eimn{{e_{i-1}^{n}}}
\def\einp { {e_i^{n+1} } }
\def\eps{\varepsilon}
\def\epsol{entropy process solution } 
\def\error{e}
\def\errorcv{\error_\cv}               
 
\def\E{\mathcal{E}}

\def\fcvn{f_{\cv}^n} 
\def\fcvnm{{f_{\cv}^{n-1}}}
\def\fim{F_{i-\frac 1 2}^n}
\def\fip{F_{i+\frac 1 2}^n}
\def\fipt{F_{i+1/2}^n}
\def\fipb{\tilde{F}_{i+\frac 1 2}^n}
\def\fipbt{\tilde{F}_{i+1/2}^n}

\def\F{\mathbf{F}}
\def\Fast{F^\star}
\def\Fastcvedge{F^\star_{\cv,\edge}}
\def\Fbar{\overline F}
\def\Fbarcvedge{\Fbar_{\cv,\edge}}
\def\Fcvcvvn{F_{\cv,\cvv}^n}
\def\Fcvedge{F_{\cv,\edge}}
\def\Fcvvedge{F_{\cvv,\edge}}

\def\Fcvvcvn{F_{\cvv,\cv}^n}

\def\g{\gamma} 

\def\gammat{\overline{\gamma}} 
\def\gcv{g_{\cv}^n}

\def\Gcvedge{G_{\cv,\edge}}

\def\half{{\frac 1 2}}
\def\halft{{1/2}}
\def\him{{h_i^-}}
\def\hip{{h_i^+}}
\def\hipum{{h_{i+1}^-}}
\def\himud{ {h_\imud} }

\def\Hast{H^\star}
\def\Hbar{\overline H}
\def\Hastde{H^{\star,+}}
\def\Hastun{H^{\star,-}}

\def\Hunz{H^1_0}

\def\imud{{i-\frac 1 2}}         
\def\iptd{{i+\frac 3 2}}         
\def\ipud{{i+\frac 1 2}}         
\def\imudt{{i-1/2}}              
\def\imu{{i- 1 }}                
\def\ipu{{i+1} }                 

\def\jpud{{j+\frac 1 2}}         

\def\l{\lambda}

\def\lapl{{\Delta}}
\def\lbar{{\overline \lambda}}

\def\lund{{\underline \lambda}}

\def\La{{\Lambda}}

\def\m{\mu}
\def\mapp{\mu_{\mesh,k}}
\def\mapz{\mu_{\mesh}}
\def\mcv{\meas(\cv)}                      
\def\mcvbord{{\meas(\dr\cv \cap \dO )}}   
\def\mcvv{\meas(\cvv)}                    
\def\meas{\mathrm{m}}
\def\medge{\meas(\edge)}                  
\def\medgecvcvv{\meas(\edgecvcvv)}
\def\mesh{\mathcal{T}}
\def\meshatyp{\mathcal{T_a}}
\def\meshprim{\mathcal{M}}
\def\mesO{{\meas(\O)}}
\def\mnn{{m \in \N}}

\def\M{\mathcal{M}}
\def\Mesh{\mathcal{T}}

\def\n {\mathbf{n}}
\def\nnn{{n \in \N}}
\def\ncvcvv{\n_{\cv,\cvv}}
\def\ncv{\n_{\cv}}  
\def\ncvedge {\n_{\cv,\edge}}  
\def\nedge {\n_{\edge}}  
\def\norm#1#2{\Vert{#1}\Vert_{#2}}    
\def\normef{\|f\|_{L^\infty(\O\times(0,T))}}  
\def\normefd{\|f\|_{L^\infty(\O\times(0,2T))}}
\def\normuz{\Vert u_0\Vert_{L^\infty(\R^\dim)}}
\def\normuzr{{\Vert u_0 \Vert_{L^\infty(\O)}}}

\def\NN{\{0,\ldots,\NTk\}}
\def\NNpar{\{0,\ldots,\Nk\}}
\def\Nk{N_{k}}
\def\NTk{N_{T,k}}
\def\NNparmun{\{0,\ldots,\Nk-1\}}
\def\NNparpun{\{0,\ldots,\Nk+1\}}
\def\NNparun{\{1,\ldots,\Nk+1\}}
\def\Nv{\mathcal{N}(\cv)}

\def\O{\Omega}
\def\o{\omega} 
\def\obar{\overline \omega}

\def\Obar{\overline \Omega}

\newcommand{\partd}{\eth}

\def\pcv{{p_\cv}}
\def\pcvv{{p_\cvv}}
\def\pa{\partial}
\def\points{\mathcal{P}}
\def\proof{\textsc{Proof}\ }

\def\P{\mathcal{P}}

\def\r{\rho} 
\def\rbarbar{\overline{\overline r}} 
\def\rcvedge{r_{\cv,\edge}}          
\def\rcvvedge{r_{\cvv,\edge}}        
\def\redge{r_\edge}

\def\refep#1{(\ref{#1}) page \pageref{#1}}
\def\refp#1{\ref{#1}    page \pageref{#1}}

\def\Redge{R_\edge}
\def\Rpast{\R_+^\star}
\def\Rcvedge{R_{\cv,\edge}}          
\def\Rcvvedge{R_{\cvv,\edge}}        
\def\RdR{\R^\dim \times (0,T)}
\def\Rin   { {R_i^{n } } }

\def\RNRp{{\R^\dim \times \R_+}}
\def\Rn{\R^\dim}
\def\Rp{\R_+}

\def\s{\sigma} 
\def\separ{\vspace*{0.5cm}}
\def\sgn{\mathrm{sign}} 
\def\size{\mathrm{size}}
\def\sizemesh{\size(\mesh)}
\def\sk{{\kappa}}
\def\snorm#1#2{\vert #1 \vert_{#2}}
\def\sumedgecv{\dsp{\sum_{\edge \in \edgescv}}}
\def\sumedges{\dsp{\sum_{\edge \in \edges}}}
\def\sumedgesint{\dsp{\sum_{\edge \in \edgesint}}}
\def\sumedgesints{\sum_{\edge \in \edgesint}}
\def\sumedgesintkl{\dsp{\sum_{\edge=\edgecvcvv \in \edgesint}}}
\def\sumedgesext{\dsp{\sum_{\edge \in \edgesext}}}
\def\summesh{\dsp{\sum_{\cv \in \mesh}}}
\def\summeshs{\sum_{\cv \in \mesh}}
\def\summeshR{\dsp{\sum_{\cv \in \mesh_R}}}
\def\sumE{\dsp{\sum_{\edgecvcvv \in \edgesint}}}
\def\sumEll{\dsp{\sum_{\edgecvcvv \in \edgesint}}}
\def\sumER{\dsp{\sum_{(\cv,\cvv) \in \edges_R^n}}}
\def\sumERi{\dsp{\sum_{(\cv,\cvv) \in \edges_R^{n+1}}}}
\def\sumIN{\dsp{\sum_{n\in\N}}}
\def\sumN{\sum_{n=0}^N}
\def\sumparTsk{\dsp{\sum_{n=0}^{\Nk}}}
\def\sumTsk{\dsp{\sum_{n=0}^{\NTk}}}
\def\sumund{\dsp{\sum_{i=1}^d}}
\def\sumv{\dsp\sum_{\cvv\in\Nv}}

\def\S{\mathcal{S}}

\def\transedge{\tau_{\edge}}            
\def\transcvcvv{\tau_{{\edgecvcvv}}}    
\def\tends{\rightarrow} 

\def\tn{{t_n}}
\def\tnp{t_{n+1}}

\def\T{\mathcal{T}}
\def\Tbar{\overline T}

\def\uapp{{u_{\Mesh,k}}}
\def\uappm{{u_{\Mesh_m,k_m}}}
\def\ubar{\overline{u}}
\def\ucv{{u_{\cv}}}
\def\ucvn{{u_{\cv}^{n}}}
\def\ucvnm{{u_{\cv}^{n-1}}}
\def\ucvnp{{u_{\cv}^{n+1}}}
\def\ucvr{{u_{\cv}^{(r)}}}
\def\ucvv{{u_{\cvv}}}
\def\ucvvn{{u_{\cvv}^n}}
\def\ucvvnm{{u_{\cvv}^{n-1}}}
\def\ucvvnp{{u_{\cvv}^{n+1}}}   
\def\ucvvr{{u_{\cvv}^{(r)}}}
\def\ucvz{{u_{\cv}^0}}

\def\uedgep{u_{\edge,+}}
\def\uedgem{u_{\edge,-}}
\def\uexi{u(x_i)}
\def\uexip{u(x_{i+1})}
\def\uexin   { {u(x_i,nk)} }
\def\ueximn{ {u(x_{i-1},nk)} }
\def\uexinp{ {u(x_i,(n+1)k)} }
\def\ug{{\gammat(u)}}
\def\uimn{ {u_{i-1}^{n } } }
\def\uin{u_i^{n}}
\def\uinp{u_i^{n+1}}
\def\uipn{u_{i+1}^n}

\def\util{{\tilde u}}
\def\u0{u_0}

\def\Ubar{\overline U}
\def\Ubarbar{\overline{\overline U}}

\def\vareps{\varepsilon}
\def\vbf{\mathbf{v}}
\def\vedge{v_{\edge}} 
 
\def\vcvedge{v_{\cv,\edge}} 
\def\vcvcvvn{v_{\cv,\cvv}^{n}}
\def\vcvcvv{v_{\cv,\cvv}}
\def\vcvvcv{v_{\cvv,\cv}}
\def\vcvvcvn{v_{\cvv,\cv}^{n}}
\def\vphi{\varphi}
\def\vphig{\varphi_{\gamma}}

\def\V{\mathcal{V}}
\def\Vastcvedge{V^\star_{\cv,\edge}}
\def\Vbarcvedge{\Vbar_{\cv,\edge}}
\def\Vcvedge{V_{\cv,\edge}}
\def\Vertices{\mathcal{V}}

\def\wcv{{w_{\cv}}}
\def\wcvv{{w_{\cvv}}}

\def\Wcvedge{W_{\cv,\edge}}

\def\Xmesh{X(\mesh)}
\def\zetamesh{\zeta_{\mesh}}

\newcommand{\bfd}{\boldsymbol{d}}
\newcommand{\bfn}{\boldsymbol{n}}
\newcommand{\bfu}{\boldsymbol{u}}
\newcommand{\bfv}{\boldsymbol{v}}

\newcommand{\bfx}{\boldsymbol{x}}
\newcommand{\bfvarphi}{\boldsymbol{\varphi}}

\newcommand{\ybar}{\overline{y}}
\newcommand{\zbar}{\overline{z}}
\newcommand{\Vbar}{\overline{V}}

\newcommand{\dive}{\mathrm{ div}}
\newcommand{\grad}{\nabla}

\def\R{\ifmmode\mathbb{R}\else$\mathbb{R}$\ \fi} 

\def\Z{\ifmmode\mathbb{Z}\else$\mathbb{Z}$\ \fi} 

\def\Q{\ifmmode\mathbb{Q}\else$\mathbb{Q}$\ \fi} 

\def\N{\ifmmode\mathbb{N}\else$\mathbb{N}$\ \fi} 

\def\C{\ifmmode\mathbb{C}\else$\mathbb{C}$\ \fi}

\def\smallskip{\vspace{.1cm}} 
\def\medskip{\vspace{.3cm}}

\newtheorem{remark}{Remark}[section]
\newtheorem{example}{Example}[section]
\newtheorem{prop}{Proposition}[section]
\newtheorem{defi}{Definition}[section]
\newtheorem{theo}{Theorem}[section]
\newtheorem{lemma}{Lemma}[section]
\newtheorem{corol}{Corollary}[section]
\newtheorem{hyp}{Assumption}[section]

\numberwithin{equation}{section}

\def\hpage#1{\markboth{}{}}

\makeindex

\begin{document}

\title{\Large\bfseries Finite Volume Methods}
\author{\large\bfseries Robert Eymard\thanks{Ecole Nationale des Ponts et Chauss\'ees, 
Marne-la-Vall\'ee, et Universit\'e de Paris XIII}, 
Thierry Gallou\"{e}t\thanks{Ecole Normale Sup\'erieure
de Lyon} and Rapha\`ele
Herbin\thanks{Universit\'{e} de Provence, Marseille}}
\date{
January 2019. This manuscript is an update of the preprint \\ 
 n0 97-19 du LATP, UMR 6632, Marseille, September 1997\\ 
which appeared in
Handbook of Numerical Analysis, \\ P.G. Ciarlet, J.L. Lions eds, vol 7, pp 713-1020 \\
}

\maketitle

\tableofcontents

\chapter{Introduction}

\hpage{Introduction}

\hpage{Introduction}

The finite volume\index{Finite volume principles} method is a discretization 
method which is well suited 
for the numerical simulation of various types (elliptic, parabolic or 
hyperbolic, for instance)  of conservation laws; 
it has been extensively used in several engineering 
fields, such as fluid mechanics, heat and mass transfer or 
petroleum engineering. Some of the important features of the finite volume method 
are similar to those of the finite element method, see \citet{oden}: it may  be used on 
arbitrary geometries, using structured or unstructured meshes, and it leads
to robust schemes. An additional feature is the local 
conservativity\index{Conservativity} of the
numerical fluxes, that is the numerical flux is conserved from one 
discretization cell to its neighbour. This  last feature makes the finite volume method
quite attractive when modelling problems for which the flux is of importance, such
as in fluid mechanics, semi-conductor device simulation, heat and mass transfer\ldots
The finite volume method is locally conservative because it
 is based on a `` balance" approach: a local 
balance is written on each discretization cell which is often called  
``control volume"\index{Control volume}; 
by the divergence formula, an integral formulation of the fluxes 
over the boundary of the control volume 
is then obtained. 
The fluxes on the 
boundary are discretized with respect to the discrete unknowns. 

\smallskip

Let us 
introduce the method more precisely on simple examples, and then give a 
description of the discretization of general conservation laws.

\smallskip

\section{Examples}

\hpage{Introduction}

Two basic  examples can be used to introduce the finite volume method. 
They will be developed in details in the following chapters.  

\begin{example}[Transport equation] \label{introte} \normalfont 
Consider first the linear transport\index{Equation!Transport --} equation
\be
\ba
u_t(x,t) + \div(\vbf u) (x,t) & = & 0, x \in \R^2,  t \in \R_+,\\
u(x,0) = u_0(x), x \in \R^2
\ea
\label{introtransp}
\ee
where $u_t$ denotes the time derivative of $u$, 
$\vbf \in C^1(\R^2,\R^2)$,  and $u_0 \in L^\infty(\R^2)$. 
 Let $\mesh$ be a  mesh of $\R^2$ consisting of polygonal bounded convex 
 subsets of $\R^2$ and let $\cv \in \mesh$ be a ``control 
 volume"\index{Control volume}, that is 
 an element of the mesh $\mesh$. Integrating the first equation of 
 \eqref{introtransp} over $\cv$ yields the following ``balance equation" 
 over $\cv$:
 \be
 \int_\cv u_t (x,t) dx + \int_{\partial \cv} \vbf(x,t) \cdot \ncv(x) u(x,t)
 \dfrontiere(x) 
 = 0, \forall t \in \R_+,
 \label{introbaltr}
 \ee
where $\ncv$ denotes the normal vector to $\partial \cv$, outward to 
$\cv$.
Let $k\in \R_+^\ast$ be a constant time discretization step  and let
 $t_n = nk$, for $n \in \N$. Writing equation \eqref{introbaltr} at time 
 $t_n$, $\nnn$ and discretizing the time partial derivative by the Euler 
 explicit scheme \index{Scheme!Explicit Euler --} suggests to find an approximation $u^{(n)}(x)$  of the 
 solution of \eqref{introtransp} at time $t_n$  which   satisfies the 
 following semi-discretized equation:
 \be
 \frac 1 k \int_\cv (u^{(n+1)} (x) - u^{(n)}(x) ) dx + \int_{\partial \cv}
 \vbf(x,t_n) \cdot \ncv(x)  u^{(n)}(x) \dfrontiere(x) = 0, \forall \n \in \N, \forall 
 \cv \in \mesh,
 \label{introsemidisc}
 \ee 
 where  $\dfrontiere$ denotes the one-dimensional Lebesgue measure on 
 $\partial \cv$ and  $
 u^{(0)}(x) = u(x,0) = u_0(x).
 $
\index{Finite volume scheme!-- for hyperbolic problems!-- in two or three 
space dimensions|(}
 We need to  define the discrete unknowns for the (finite volume) space 
discretization. We shall be concerned here principally with the so-called 
``cell-centered" finite volume method in which each discrete unkwown is 
associated with a control volume. Let $(u^{(n)}_\cv)_{\cv \in \mesh, n \in \N}$
denote the discrete unknowns. For $\cv \in \mesh$, 
let $\edgescv$ be the set of edges which are included in $\partial \cv$, 
\index{Edge}
and for $\edge \subset \partial \cv$, let  $\ncvedge$ denote the unit normal 
to $\edge$ outward to 
$\cv$.   
The second integral in \eqref{introsemidisc} may then be split as:
\be
 \int_{\partial \cv}
 \vbf(x,t_n) \cdot \ncv(x)  u^{(n)}(x)  \dfrontiere(x) = \sumedgecv \int_\edge
 \vbf(x,t_n) \cdot  \ncvedge u^{(n)}(x) \dfrontiere(x);
\label{introflux1}
\ee
for $\edge \subset \partial \cv$, let 
\[
\vcvedge^{(n)} = \int_\edge \vbf(x,t_n) \ncvedge(x) \dfrontiere(x).
\]
Each term of the sum in the right-hand-side of \eqref{introflux1}
is then discretized as
\be
F_{\cv,\edge}^{(n)} = \ba \vcvedge^{(n)} u^{(n)}_{\cv} & \text{ if 
}\vcvedge^{(n)} \ge 0, \\ \vcvedge^{(n)} u^{(n)}_{\cvv} & \text{ if } \vcvedge^{(n)} 
< 0, \ea
\label{introflux2} 
\ee
where $\cvv$ denotes the neighbouring control volume to $\cv$ with 
common edge $\edge$. 
This ``upstream"\index{Upstream} or  ``upwind" choice is classical for transport equations;
it may be seen, from the mechanical point of view, as the choice of the
``upstream information"   with respect to   the location of $\edge$. This 
choice is  crucial in the mathematical analysis;
it ensures the stability properties of the
finite volume scheme (see chapters \ref{hyper1d} and \ref{hypmd}). 
We have therefore derived the following finite volume 
scheme for the 
discretization of \eqref{introtransp}:
\be
\ba
\dsp{\frac \mcv k }(u_\cv^{(n+1)} - u_\cv^{(n)}) + 
\sumedgecv F_{\cv,\edge}^{(n)} = 
0, \forall \cv \in \mesh, \forall n \in \N, \\
u^{(0)}_\cv = \dsp{\int_\cv u_0(x) dx},
\ea 
\ee
where $\mcv$ denotes the measure of the control volume \index{Control volume}
 $\cv$
 and $ F_{\cv,\edge}^{(n)}$ is defined in \eqref{introflux2}.
This scheme is locally conservative  in the sense that if $\edge$ is a 
common edge to the control volumes $\cv$  and $\cvv$, then
$F_{\cv,\edge} = - F_{\cvv,\edge}$. This property is important in 
several application fields; it will later be shown to be
a key ingredient in the 
mathematical proof of convergence.
Similar schemes for the discretization of linear or nonlinear
hyperbolic equations will be studied in chapters 
\ref{hyper1d} and \ref{hypmd}.
\index{Finite volume scheme!-- for hyperbolic problems!-- in two or three 
space dimensions|)}
\end{example}

\begin{example}[Stationary diffusion equation]\label{introde} \normalfont
Consider the basic diffusion equation\index{Equation!Diffusion --}
\be
\ba
-\lapl u = f \text{ on } \Omega = ]0,1[ \times ]0,1[, \\
u = 0 \text{ on } \partial \Omega.
\ea 
\label{introdiffeq}
\ee
 Let $\mesh$ be a rectangular mesh. Let us integrate 
the first equation of  \eqref{introdiffeq} over a control volume
\index{Control volume} $\cv$
of the mesh; with  the same notations as 
in the previous example, this yields:
\be
 \sumedgecv \int_\edge - \grad u(x) \cdot \ncvedge \dfrontiere(x)  = 
\int_\cv f(x) dx.
\ee
For each control
volume $\cv \in \mesh$, let  $\centercv$ be the center of $\cv$.
Let $\edge$ be the common edge between the control volumes $\cv$ 
and $\cvv$.  
One way to approximate the flux $- \int_{\edge} \grad u(x) \cdot \ncvedge \dfrontiere(x)$ (although clearly not the only one), 
 is to use a centered  finite difference approximation:
\be
F_{\cv,\edge} = - \frac \medge \dedge (u_\cvv - u_\cv),
\label{introfluxint}
\ee
where $(u_\cv)_{\cv \in \mesh}$ are the discrete unknowns and $\dedge$  is the distance between $\centercv$ and $\centercvv$. This finite difference approximation of the first order derivative $\grad u \cdot \n$ on the edges of the mesh (where 
$\n$ denotes the unit normal vector) is consistent\index{Consistency!-- of the fluxes}: the truncation error on the flux is of order $h$, where $h$ is the maximum length of the edges of the mesh. We may note that the 
consistency of the flux holds because for any $\sigma = K \vert L$ common to the control volumes $K$ and $L$, the line segment $[\centercv \centercvv]$ is perpendicular to  $\sigma = K \vert L$.
Indeed, this is the case here since the control volumes are rectangular. 
This property is satisfied by other meshes which will bestudied hereafter. It is crucial for the discretization of diffusion operators. 

In the case where the edge $\edge$ is part of the boundary, then $\dedge$ denotes the distance 
between the center $\centercv$ of the control volume $\cv$ to which  $\edge$ belongs and  the boundary. 
The flux $- \int_{\edge} \grad u(x) \cdot \ncvedge \dfrontiere(x)$, is then approximated by
\be
F_{\cv,\edge} =  \frac \medge \dedge  u_\cv,
\label{introfluxext}
\ee
 Hence the finite volume scheme for the discretization of 
 \eqref{introdiffeq} is:
\be
 \sumedgecv \Fcvedge = \mcv f_\cv, \forall \cv \in \mesh,
\ee
where $\Fcvedge$ is defined by \eqref{introfluxint} and 
\eqref{introfluxext}, and $f_\cv$ denotes (an approximation of) the mean 
value of $f$ on $\cv$. 
We shall see later (see chapters \ref{ellund},  \ref{ellmd} 
and \ref{parabolic})
that the finite volume scheme is easy to generalize  to a triangular
mesh, whereas the finite 
difference\index{Finite difference method} method is not. As in the previous 
example, the finite volume scheme is locally conservative, 
since for any 
edge $\edge$ separating $\cv $ from $\cvv$, one has $\Fcvedge = 
-\Fcvvedge$.

\end{example}

\section{The finite volume principles for general conservation 
laws}\index{Finite volume principles}

\hpage{Introduction}

The finite volume method is used for the discretization of conservation 
laws. We gave in the above section   two  examples of such conservation laws. 
Let us now present the discretization of general conservation laws by
finite volume schemes. 
As suggested by its name, a conservation law expresses the conservation 
of a  quantity $q(x,t)$. For instance, the conserved
quantities may be the energy, the mass, or the number of moles
of some chemical species. 
Let us first assume that the local form of 
the conservation equation\index{Equation!Conservation --} may be written as 

\be
q_t(x,t) + \div \F (x,t)  = f(x,t),
\label{conserv.loc} 
\ee

at each point $x$ and each time $t$ where the conservation  
of   $q$ is to be 
written. In equation  (\ref{conserv.loc}), $(\cdot)_t$ denotes the time 
partial derivative of the entity  within the parentheses, $\div$
represents the space divergence operator: $\div \F = \dr F_1 /{\dr
x_1}+\dots + \dr F_d / {\dr
x_d}$, where $\F=(F_1,\ldots,F_d)^t$ denotes a 
vector function depending  on the space variable $x$ and on the time $t$, 
$x_i$ is the $i$-th space 
coordinate, for $i = 1,\ldots,d$,
and $\dim$ is the space dimension, i.e. $\dim = 1,2$ or $3$; 
the quantity $\F$ is a flux which expresses a transport mechanism of $q$; 
the ``source term" $f$ expresses a possible 
volumetric exchange, due for instance  to chemical reactions between
the conserved quantities.

\medskip

Thanks to the physicist's work, the problem can be closed
by introducing   constitutive laws
which relate $q$, $\F$, $f$ with some scalar or vector unknown 
$u(x,t)$, function of the space variable $x$ and of the time $t$.
For example, the components of $u$ can be pressures, concentrations, 
molar fractions of the various chemical species by unit volume\ldots
The quantity $q$ is often given by means of a 
known function $\bar{q}$ of $u(x,t)$, of the space
variable $x$ and of the time $t$, that is $q(x,t)=\bar{q}(x,t,u(x,t))$. 
The quantity $\F$ may also be given by means of a function of the space 
variable $x$, the time variable $t$ and of the unknown $u(x,t)$ and (or) 
by means of the gradient of $u$ at point $(x,t)$\ldots.  
The transport equation of Example \ref{introte} is a particular case of  
\eqref{conserv.loc} with $q(x,t) = u(x,t)$, $\F(x,t) = \vbf  u(x,t)$
and $f(x,t) = f(x)$; so is the stationary
diffusion equation of Example \ref{introde} with $q(x,t) = u(x)$, 
 $\F(x,t)= - \grad  u(x)$, and $f(x,t) = f(x)$. 
 The source term $f$ may  also be  given by means of a function 
of $x$, $t$ and $u(x,t)$. 

\begin{example}[The one-dimensional Euler equations] \normalfont 
Let us consider as an example of a system of conservation laws 
the 1D Euler 
equations\index{Euler equations} for 
equilibrium real gases; these equations  may be written under 
the form \eqref{conserv.loc}, 
with \[  q  = \Bigl(\ban \rho\\ \rho u\\ E \ean \Bigr)  \text{ and }
\F = \Bigl(\ban \rho u \\ \rho u^2  +p \\  u (E + p )\ean \Bigr), \]
where $\rho,  u, E$ and $p$ are functions of the space variable $x$ 
and the time $t$, and
refer respectively to the density, 
the velocity, the total energy and  the pressure of the particular gas under
consideration. The system of equations is closed by introducing the constitutive laws
which relate $p$ and $E$ to the specific volume $\tau$, with $\tau = \frac 1 \rho$
and the entropy $s$, through the constitutive laws:
\[ p = \frac {\partial \eps } {\partial \tau } (\tau,s) \text{ and } E = \rho(\eps(\tau,s)  
+ \frac {u^2} 2), \]
where $\eps$ is the internal energy per unit mass, which is a given function of $\tau$ and $s$.
\end{example}

Equation (\ref{conserv.loc}) may be seen as the expression of the conservation of $q$ 
in an infinitesimal domain; it is formally equivalent to the equation

\be
\bari 
\dsp{\int_\cv q(x,t_2) dx - \int_\cv q(x,t_1) dx} + 
 \dsp{\int_{t_1}^{t_2} \int_{\partial \cv} 
\F(x,t)\cdot \bfn_{\cv}(x) \dfrontiere(x) dt} \\
 = \dsp{\int_{t_1}^{t_2} \int_\cv f(x,t) dx dt},
\ean
\label{conserv.vol}
\ee

for any subdomain $\cv$ and for all times $t_1$ and $t_2$, where  $\bfn_{\cv}(x)$ 
is the unit normal vector to the boundary $\dr \cv$, at point $x$, outward to 
$\cv$. Equation \eqref{conserv.vol} expresses the conservation law in 
subdomain $\cv$ between times $t_1$ and $t_2$. Here and in the sequel, 
unless otherwise mentionned, $dx$ is the integration symbol for the 
$\dim$-dimensional Lebesgue measure in $\R^\dim$ and $\dfrontiere$ is the 
integration symbol for the $(\dim-1)$-dimensional Hausdorff measure on the considered boundary.

\subsection{Time discretization}

\hpage{Introduction}

The time  discretization of Equation  (\ref{conserv.loc}) is performed by 
introducing  an increasing sequence $(t_n)_\nnn$ with $t_0 = 0$. For the 
sake of simplicity, only constant time steps will be considered here, 
keeping in mind that the generalization to variable time steps is 
straightforward. Let  $k \in \R_+^\star$ denote the time step, and let  
$t_n=nk$, for $n \in \N$. It can be noted that
Equation~(\ref{conserv.loc}) could be written with the use of a  
space-time divergence.  Hence, 
Equation~(\ref{conserv.loc}) could be either discretized using a  
space-time finite volume
discretization or a space finite volume discretization with a time finite difference scheme 
(the explicit Euler scheme, for instance).  In the first case, the 
conservation law is integrated over a time interval and a space ``control 
volume" as in the formulation \eqref{conserv.loc}. 
In the latter case, it is only integrated space wise, and the 
time derivative is approximated by a finite difference scheme; 
with the explicit Euler scheme\index{Scheme!Explicit Euler --}, the term $(q)_t $ is therefore approximated
by the differential quotient 
$(q^{(n+1)} - q^{(n)}) / k$, and $q^{(n)}$ is computed with an approximate value of
$u$ at time $t_n$, denoted by $u^{(n)}$. Implicit and higher order schemes 
may also be used.

\subsection{Space discretization}

\hpage{Introduction}

In order to perform a space finite volume discretization of equation~(\ref{conserv.loc}), a 
mesh $\mesh$ of the domain $\O$ of $\R^\dim$, over which the conservation law is to be 
studied, is introduced. The mesh is such that $\Obar = \cup_{\cv \in 
\mesh} \cvbar$, where an element of $\mesh$, denoted by $\cv$, 
is an open subset of $\O$ and 
is called  a 
control volume. Assumptions on the meshes will be needed for the 
definition of the schemes; they also depend on the type of equation to 
be discretized. 

\smallskip

For the finite volume schemes considered here, 
the discrete unknowns at time $t_n$ 
are denoted by $u^{(n)}_\cv$, $\cv \in \mesh$. The value $u^{(n)}_\cv$ is expected to 
be some approximation  of $u$ on the cell
$\cv$ at time $t_n$. The basic principle of the classical finite
volume method is to integrate equation~(\ref{conserv.loc}) over each
cell $\cv$ of the mesh $\mesh$.  
One obtains a conservation law under a nonlocal form 
(related to equation~(\ref{conserv.vol})) written
for the volume $\cv$. Using  the Euler time discretization, 
this yields 

\be
\int_\cv \frac {q^{(n+1)}(x) - q^{(n)}(x)} k dx +
 \int_{\partial  \cv} \F(x,t_n) \cdot \bfn_{\cv}(x) \dfrontiere(x) = 
 \int_{\cv} f(x,t_n)dx,
\label{eVF3} 
\ee

where $ \bfn_{\cv}(x) $ is the unit normal vector to $\dr \cv$
at point $x$, outward to $\cv$. 

\smallskip

The remaining step in order to define the finite volume scheme
is therefore the approximation of the ``flux", $\F(x,t_n) \cdot 
\bfn_{\cv}(x)$, across the
boundary $\partial  \cv$ of each control volume, in terms of $\{u^{(n)}_{\cvv},
\cvv \in \mesh\}$ (this flux approximation has to be done in terms of 
$\{u^{n+1}_{\cvv}, \cvv \in \mesh\}$  
if one chooses the implicit Euler scheme  instead of the explicit Euler scheme for 
the time discretization). More 
precisely, omitting the terms on the boundary of $\O$, let
$\edgecvcvv=\cvbar \cap \cvvbar$, with $\cv$, $\cvv \in 
\mesh$, the exchange term (from $\cv$ to $\cvv$), $\int_\edgecvcvv 
\F(x,t_n) \cdot \bfn_{\cv}(x)
\dfrontiere(x)$,  between the control volumes $\cv$ and $\cvv$ during the time
interval $[t_n,t_{n+1})$ is approximated by some quantity, 
$F^{(n)}_{\cv,\cvv}$, which is a function of $\{u^{(n)}_{M}, M \in \mesh\}$ (or a 
function of $\{u^{n+1}_{\cvw}, \cvw \in \mesh\}$  for the implicit Euler 
\index{Scheme!Implicit Euler --} 
scheme, or  more generally  a function of $\{u^{(n)}_{M}, M \in \mesh\}$ and 
$\{u^{n+1}_{\cvw}, \cvw \in \mesh\}$ if the time discretization is a one-step
 method). Note that $F^{(n)}_{\cv,\cvv} = 0$ if the Hausdorff dimension of $\cvbar
\cap \cvvbar$ is less than $d-1$ (e.g.  $\cvbar
\cap \cvvbar$ is a point in the case $d=2$ or a line segment in the case 
$d=3$). 

\smallskip

Let us point out that two important features of the classical finite volume method are
\begin{enumerate}
\item the conservativity\index{Conservativity}, that is $F^{(n)}_{\cv,\cvv} = - F^{(n)}_{\cvv,\cv}$,
 for all $ \cv$ and $ \cvv$ $\in \mesh$ and for all $\nnn$.
\item the ``consistency" of the approximation of 
$\F(x,t_n) \cdot \bfn_{\cv}(x)$, which has to
be defined for each relation type between $\F$ and the unknowns.  
\end{enumerate}

These properties, together with adequate stability 
properties which are obtained by
estimates on the approximate solution, will give some convergence
properties of the finite volume scheme.

\section{Comparison with other discretization techniques}

\hpage{Introduction}

The finite volume method is quite different from (but sometimes related 
to)  the
finite difference method or the finite element method. On these classical 
methods see e.g. \citet{DB}, \citet{finitedifference}, \citet{ciarlet},
\citet{finitelement},
\citet{mixedfinitelement}.

\smallskip

Roughly speaking, the principle of the finite difference 
method\index{Finite difference method} is, 
given a  number
of discretization points which may be defined by a mesh, to assign
one discrete unknown per discretization point, and to write one
equation per discretization point. At each discretization point, the
derivatives of the unknown are replaced by finite differences through
the use of Taylor expansions. 
The finite difference method becomes difficult to use when the coefficients involved in the equation are discontinuous\index{Discontinuous
coefficients} 
(e.g. in the case of heterogeneous media). With the finite volume method, 
discontinuities of
the coefficients will not be any problem if the mesh is chosen such
that the discontinuities of the coefficients occur on the boundaries
of the control volumes (see sections \ref{ell1D} and \ref{ell2Dgen}, for 
elliptic problems).  Note that
the finite volume scheme is often called ``finite difference scheme" or  ``cell
centered difference scheme". Indeed, in the finite volume method, the finite 
difference approach can be used for the approximation of the fluxes on the
boundary of the control volumes. Thus, the finite volume scheme   differs  from
the finite difference scheme in that the finite difference  approximation is used
for the flux rather than for the operator itself. 

\smallskip

The finite element method\index{Finite element method|(}  (see e.g. \citet{ciarlet}) is based 
on a variational formulation, which is  written for both the continuous  
and the discrete
problems, at least in the case of conformal finite element methods  which are 
considered here.  The variational formulation 
is obtained by multiplying the original equation by a ``test function". The
continuous unknown is then approximated by a linear combination of 
``shape" functions; these shape functions are the test functions 
for the discrete variational formulation
(this is the so called ``Galerkin expansion\index{Galerkin expansion}"); 
the resulting equation is integrated
over the domain. The finite volume method is sometimes called a 
``discontinuous finite
element method" since the original equation is multiplied by 
the characteristic function
of each grid cell  which is defined by 
$1_\cv(x) = 1$, if $ x \in \cv$,  $1_\cv(x) = 0$, if $ x \notin \cv$, 
and the discrete unknown may be considered as a linear combination of 
shape functions. However, the techniques used to prove the convergence of 
finite element methods do not  generally apply for this choice of test 
functions. In the following chapters, the finite
volume method will be compared in more detail with the classical and 
the mixed finite element
methods\index{Finite element method|)}.

\medskip

From the industrial point of view, the finite volume method is known as  
a robust and cheap method for the discretization of 
conservation laws (by robust, we mean a scheme which behaves well even 
for particularly
difficult equations, such as nonlinear systems of hyperbolic equations    
and which  can easily 
be extended to more realistic and physical contexts than the classical academic
problems). The finite volume method is cheap thanks to 
short and  reliable computational coding for complex problems.
It may be more adequate than the finite difference method\index{Finite 
difference method} (which in 
particular requires a simple 
geometry). 
However, in some cases, it is difficult to design schemes which give 
enough  precision. 
Indeed, the finite element method can be much more  precise than the finite
volume method when using higher order polynomials, but it  requires an adequate
functional framework which is not always available in industrial problems. 
Other more precise methods are, for instance, particle methods  
or spectral methods
but these methods can be more expensive and less robust than the finite 
volume method.

\section{General guideline}

\hpage{Introduction}

\smallskip 

The mathematical theory of finite volume schemes has recently been 
undertaken.  Even though we choose here to refer to 
the class of scheme which is the object of our study  as
the "finite volume'' method, we must point out that there are several 
methods with 
different names (box method, control volume finite element\index{Control 
volume finite element method} methods, 
balance method to cite only a few) which may be viewed as finite volume 
methods. The name "finite difference" has also often been used 
referring to  the finite volume method. We shall mainly  quote here the
works regarding the mathematical analysis of the finite volume method, 
keeping in mind that there exist numerous works on applications of the finite volume methods in 
the applied sciences, some references to which may be found in the books
which are cited below.

\smallskip

Finite volume methods for convection-diffusion 
equations  seem to have been first introduced in the early sixties  by 
 \citet{samarskitich}, \citet{samarski1}  and
\citet{samarski2}.  
 
The convergence theory of such schemes in several space dimensions has 
only recently been undertaken. 
In the case of vertex-centered finite volume schemes, studies were carried out by
 \citet{samarskilaz} in the case of Cartesian meshes, \citet{Heinrich}, 
  \citet{BR}, \citet{Cai}, \citet{CMMC} and 
\citet{Vanselow} in the case of unstructured meshes; see also 
\citet{mortonsuli}, \citet{suli}, \citet{mackenzie}, \citet{mortonss} 
 and \citet{shashkov}  in the case of quadrilateral meshes.
Cell-centered finite volume schemes are addressed in \citet{WM}, \citet{FS}, \citet{weiser} and \citet{lazarovsiam}  in the case of Cartesian meshes  and in \citet{vass}, \citet{VF4}, \citet{rouenrh}, \citet{rouenlazarov}, \citet{mishev} in the case of  triangular or Vorono\"{\i} meshes; let us also mention \citet{rouencoudiere} and  \citet{coudiere} where more general meshes 
are treated, with, however, a somewhat technical geometrical condition. 
In the pure diffusion case,the cell centered finite volume method has also been analyzed with finite element tools:
\citet{ABMO}, \citet{angermann}, \citet{BMO}, \citet{arbogast}, \citet{angermann}.
Semilinear convection-diffusion are studied in \citet{FFL2} with a combined finite element-finite volume method, \citet{cvnl} with a pure finite volume scheme. 
\smallskip

Concerning nonlinear hyperbolic conservation laws, the one-dimensional  case is now classical; let us mention the following books on   numerical methods for hyperbolic problems:  \citet{godlewski-ellipses}, \citet{leveque}, \citet{godlewski-springer}, \citet{kronerbook},  and references therein. 
 In the multidimensional case, let us mention the convergence results which where obtained in
 \citet{cgh}, \citet{kroner1}, \citet{Coquelconv}  and the error estimates of \citet{Coquel} and \citet{Vi} in the case of an explicit scheme and \citet{eggh} in the case of explicit and implicit schemes. 
The proof of the error estimate of \citet{eggh}, which is concerned with a flux of the form $\bfv(\bfx,t) f(\bfu)$ can easily be adapted for general fluxes of the form $F(\bfx,t,\bfu)$ \citet{Chainais2}.

\medskip

The purpose of the following chapters  is to lay out a mathematical framework for the convergence and error analysis of the  finite volume method for the discretization of elliptic, parabolic or hyperbolic  partial differential equations under conservative form, following the philosophy of the works of \citet{cgh}, \citet{VF4}, \citet{eggh} and \citet{cvnl}. 
In order to do so, we shall describe the implementation of the finite volume method on some simple (linear or non-linear) 
academic problems, and develop the tools which are needed for the mathematical analysis. 
This approach helps determine  the properties of finite volume schemes which lead to ``good" schemes for complex applications.

 Chapter \ref{ellund} introduces the finite volume discretization of an elliptic 
operator in one space dimension. The resulting numerical scheme is
compared to finite difference, finite element and mixed finite element 
methods in this particular case. An error estimate is given; this estimate 
is in fact contained in results shown later in the multidimensional case; 
however, with the one-dimensional case, one can already understand the 
basic principles of the convergence proof, and understand the difference 
with the proof of \citet{WM}  or \citet{FS}, which does not seem to generalize  to the unstructured meshes. 
In particular, it is made clear that, although the finite volume scheme is not consistent \index{Consistency!-- 
in the finite difference sense} in the finite difference sense since  the truncation error does not tend to 0, the conservativity  of the scheme, together with a consistent approximation of the fluxes and some ``stability" allow the proof of convergence.
The scheme and the error estimate are then generalized to the case of a more general elliptic operator allowing discontinuities in the diffusion coefficients. 
Finally, a semilinear problem is studied, for which a convergence result is proved. 
The principle of the proof of this result may be used for nonlinear problems in several space dimensions.   
It is used  in Chapter \ref{ellmd} in order to prove convergence results for linear problems when no  
regularity on the exact solution is known.

\smallskip 

In Chapter \ref{ellmd}, the discretization of elliptic problems in several space 
dimensions by the finite volume method is presented. Structured meshes 
are shown to be an easy generalization of the one-dimensional case; unstructured meshes are then considered, for Dirichlet and Neumann conditions on the boundary of the domain. In both cases, admissible meshes are defined, and, 
following \citet{cvnl}, convergence results (with no  regularity on the data) and error estimates  assuming a $C^2$ or $H^2$ regular solution to the continuous problems are proved.  
As in the one-dimensional case, the conservativity of the scheme, together with a consistent approximation of the fluxes   and some ``stability" are used for the proof of convergence.
In addition to the properties already used in the one-dimensional case, the multidimensional estimates require the  use of a ``discrete Poincar\'e" inequality  which is proved in both Dirichlet and Neumann cases, along with some  compactness properties  which are also used and are given in the last section.  
It is then shown how to deal with matrix diffusion coefficients and more general boundary conditions. 
Singular sources and mesh refinement are also studied. 

\smallskip

Chapter \ref{parabolic} deals with the discretization of parabolic problems.
Using the same concepts as in Chapter \ref{ellmd}, an error estimate is 
given in the linear case. A nonlinear degenerate parabolic problem is 
then studied, for which a convergence result is proved, thanks to a 
uniqueness result which is proved at the end of the chapter.

\smallskip

Chapter \ref{hyper1d} introduces the finite volume discretization of a hyperbolic  operator in one space dimension. 
Some basics on entropy weak solutions to nonlinear hyperbolic equations are recalled. 
Then the concept of stability  of a scheme is  explained on a simple linear advection problem, for which both finite difference and finite volume schemes are considered.
Some well known schemes are presented with a  finite volume formulation in the nonlinear case. 
A proof of convergence using a ``weak $BV$ inequality" which was found to be crucial in the multidimensional case (Chapter  \ref{hypmd})
is given in the one-dimensional case for the sake of clarity. 
For the sake of completeness, the proof of convergence based on ``strong $BV$ estimates"  and the Lax-Wendroff theorem 
is also recalled, although it is not used  for general meshes in the multidimensional case.

\smallskip

In Chapter \ref{hypmd},  finite volume schemes  for the discretization of multidimensional nonlinear hyperbolic conservation equations are studied. 
Under suitable assumptions, which are satisfied by several well known schemes, it is shown that the considered schemes are $L^\infty$ stable (this is classical) but also satisfy some ``weak $BV$ inequality". 
This ``weak $BV$" inequality is the key estimate to the proof of convergence of the schemes.
Following \citet{eggh}, both time implicit and explicit discretizations are considered. 
In the case of the implicit scheme, the existence of the solution must first be proved.
The approximate solutions are shown to satisfy some discrete entropy inequalities. 
Using the weak $BV$ estimate, the approximate solution is also shown to satisfy some continuous entropy inequalities. 
Introducing the concept of ``entropy process solution" to the nonlinear hyperbolic equations (which is similar to the notion of measure valued solutions of \citet{DP}), the approximate solutions are proved to converge towards an entropy process solution as the mesh size tends to 0. 
The \epsol is shown to be unique, and is therefore equal to the entropy weak solution, which concludes the convergence of the approximate solution
towards the entropy weak solution. 
Finally error estimates are proved for both the explicit and implicit schemes. 

\smallskip 

The last chapter is concerned with systems of equations. In the case of 
hyperbolic systems 
which are considered in the first part, 
little is known concerning the continuous problem, 
so that the schemes which are introduced are only shown to be efficient 
by numerical experimentation. These ``rough" schemes seem to be efficient 
for complex cases such 
as the Euler equations \index{Euler equations} for real gases. 
The incompressible Navier-Stokes\index{Navier-Stokes equations}   equations are then considered; after 
recalling the classical staggered grid finite volume 
formulation\index{Staggered grid} (see
e.g. \citet{Pat}), a finite volume scheme defined on a triangular mesh 
for the Stokes equation is studied. In the case of equilateral triangles, 
 the tools of Chapter \ref{ellmd} allow to show that the approximate velocities   
converge to the exact velocities.
Systems arising from  modelling  multiphase flow in porous media are
then considered. The convergence of the approximate finite volume 
solution for a simplified case is then proved with the tools introduced 
in Chapter \ref{hypmd}.

\medskip

More precise 
references to recent works on the convergence of finite volume methods 
will be made in the following chapters. However, we shall not quote here 
the numerous works on applications of the finite volume methods in 
the applied sciences.

 \chapter{A one-dimensional elliptic problem}\label{ellund}

The purpose of this chapter is to give some developments of the example 
\ref{introde} of the introduction in the one-dimensional case. The 
formalism needed to define admissible finite volume meshes is first given 
and applied to the Dirichlet problem. After some comparisons with other 
relevant schemes, convergence theorems and error estimates are provided.  
Then, the case of general linear elliptic equations is handled and 
finally, a first approach of a nonlinear problem is studied and introduces
some compactness theorems in a quite simple framework; 
these compactenss theorems will be useful in further 
chapters.

\section{A finite volume method for the Dirichlet problem} \label{Diric1D}

\hpage{Elliptic 1d}

\subsection{Formulation of a finite volume scheme}\label{fvse1d}

The principle of the finite volume method will be shown here
on the academic Dirichlet  problem, namely a second order differential operator 
without time dependent terms and with homogeneous Dirichlet boundary 
conditions\index{Dirichlet boundary condition|(}. 
Let $f$ be a given function from $(0,1)$ to $\R$, consider the following differential 
equation:
 
\be
\begin{array}{lcl}  
- u_{xx}(x) & = & f (x), \ \  x \in (0,1),   \\
u(0) & = & 0,   \\
u(1) & = & 0. 
\end{array}
\label{elliptic1D}
\ee

If $f \in C([0,1],\R)$, there exists a unique solution $u \in  C^2([0,1],\R)$ to Problem
(\ref{elliptic1D}). In the sequel, this exact solution will
be denoted by $u$. Note that the equation $-u_{xx}=f$ can be written in the conservative
form $\div (\F)= f$ with $\F=-u_x$. 

In order to compute a numerical approximation to the solution of this equation, 
let us define a mesh,  denoted by $\mesh$, of the interval $(0,1)$ consisting of $N$ 
cells (or control volumes\index{Control volume}), denoted by $K_i$, $i = 1,\ldots,N$, and $N$ 
points of $(0,1)$, denoted by $x_i$, $i = 1,\ldots,N$, satisfying the 
following
assumptions:

\begin{defi}[Admissible one-dimensional mesh]  \normalfont \label{meshund} 
\index{Mesh!Admissible --!--in the one-dimensional elliptic case}
An admissible mesh of $(0,1)$, denoted by $\mesh$, is given by a family 
$(K_i)_{i = 1,\cdots,N}$, $N \in \N^\star$, such that
$K_i = (x_{i-\frac 1 2}, x_{i+\frac 1 2})$, and a family $(x_i)_{i = 0,\cdots,N+1}$ such that

\[
x_0=x_{\frac 1 2}=0 <x_1<x_{\frac 3 2} < \cdots < x_{i-\frac 1 2}<x_i<x_{i+\frac 1 2}
<\cdots < x_N<x_{N+\frac 1 2}=x_{N+1}=1.
\]
One sets

\[
\bal
h_i = \meas(\cv_i)=x_{i+\frac 1 2}-x_{i-\frac 1 2}, \,  i=1,\ldots,N, 
\text{~and therefore } \dsp{\sum_{i=1}^N h_i = 1},
\\
h_i^- = x_i-x_{i-\frac 1 2}, h_i^+ = x_{i+\frac 1 2}-x_i,  \, i=1,\ldots,N,
\\
h_{i+\frac 1 2}=x_{i+1}-x_i, \,  i=0,\ldots,N, 
\\
\sizemesh=h=\max \{h_i, \, i = 1,\ldots, N\}.
\ean
\]

\end{defi}


The discrete unknowns are denoted by $u_i$, $i = 1,\ldots, N$, and are expected to be some approximation of $u$ in the cell $K_i$ (the discrete unknown $u_i$ can be viewed as an approximation of the mean value of $u$ over $K_i$, or of the value of $u(x_i)$, or of other values of $u$ in the control volume $K_i$\ldots). 
The first equation of~(\ref{elliptic1D}) is integrated over each cell $K_i$, as in \eqref{eVF3} and yields

\[
-u_x(x_{i+\frac 1 2}) +  u_x(x_{i-\frac 1 2}) = \int_{K_i} f (x) dx,
\qquad i = 1,\ldots, N.
\]

A reasonable choice for the approximation of $- u_x(x_{i+\frac 1 2})$ (at least, for
$i=1,\ldots,N-1$) seems to be the differential quotient

\[
F_\ipud = - \dsp{\frac {u_{i+1} - u_i} {h_{i+ \frac 1 2}}}.
\]
 This approximation is consistent \index{Consistency!-- of the fluxes} in the sense that, if $u \in
C^2([0,1],\R)$, then there exists $C \in \R_+$ only depending on $u$ such that
\be
\vert R_{\ipud} \vert=  \vert \Fast_{\ipud} +u_x(x_{i+\frac 1 2}) \vert  \le Ch, \mbox{ where }  \Fast_{\ipud} =  - \frac {u(x_{i+1}) - u(x_i)}{h_{i+\frac 1 2}}.
\label{consistflux}
\ee
The quantity $R_{\ipud}$ is called the consistency error \index{Consistency! -- error}.
\begin{remark}[Using the mean value] \normalfont
Assume that $x_i$ is the center of $K_i$. Let $\tilde{u}_i$ denote the mean value over $K_i$  of the exact solution $u$ to Problem~(\ref{elliptic1D}). 
One may then remark that  $\vert \tilde{u}_i - u(x_i) \vert \le C h_i^2$, with some $C$ only depending on $u$;  it follows easily that  $(\tilde{u}_{i+1} - \tilde{u}_i) / h_{i+\frac 1 2}
= u_x(x_{i+\frac 1 2}) + 0(h)$ also holds, for $i=1,\ldots,N-1$ (recall 
that $h = \max\{h_i,i=1,\ldots,N\}$). Hence the approximation 
of the flux is also consistent  if the discrete unknowns $u_i$, $i = 
1,\cdots,N$, are viewed as approximations of the mean value of $u$ in 
the control volumes.
\end{remark}

The Dirichlet boundary conditions are taken into account by using the
values imposed at the boundaries to compute the fluxes on these boundaries.
Taking these boundary conditions into consideration and setting $f_i = \frac 1 {h_i} \int_{K_i} f(x) dx$ for $i=1,\ldots,N$ (in an actual computation, an approximation of $f_i$ by numerical integration can be used), the finite volume scheme  \index{Finite volume scheme!-- for elliptic problems!-- in one space dimension|(} for problem~(\ref{elliptic1D}) reads
\be
F_{i+ \frac 1 2} - F_{i- \frac 1 2}  =  h_i f_i, \, i =
1,\ldots,N
\label{VF1}
\ee
\be
F_{i+ \frac 1 2}  =   - \frac {u_{i+1} - u_i} {h_{i+ \frac 1 2}}, \, i =
1, \ldots, N-1,
\label{VF2}
\ee
\be
F_{\frac 1 2}  =  - \frac {u_{1}} {h_{\frac 1 2}},
\label{VF3}
\ee
\be
F_{N+ \frac 1 2} =  \frac {u_{N}}{h_{N+\frac 1 2}}.
\label{VF4}
\ee

Note that (\ref{VF2}), (\ref{VF3}), (\ref{VF4}) may also be written

\be
F_{i+ \frac 1 2} = -\frac {u_{i+1} - u_i} {h_{i+ \frac 1 2}}, \, i =
0,\ldots,N,
\label{VF2m}
\ee

setting

\be
u_0 =u_{N+1}= 0.
\label{VF3m}
\ee

\vspace{.1cm}

The numerical scheme (\ref{VF1})-(\ref{VF4}) may be
written under the following  matrix form:

\be
A U=b,
\label{AUb}
\ee
where $U = (u_1,\ldots,u_N)^t$, $b = (b_1,\ldots,b_N)^t$, with
\eqref{VF3m} and with $A$ and $b$ defined by

\be
(AU)_i = \frac 1 {h_i} \Bigl(- \frac {u_{i+1} - u_i}
{h_{i+ \frac 1 2}} + \frac {u_{i} - u_{i-1}} {h_{i- \frac 1 2}} \Bigr), \,
i = 1,\ldots,N,
\label{matrix}
\ee
\be
b_i =  \frac 1 h_i \int_{K_i} f(x) dx, \,
i = 1,\ldots,N,
\label{rhsdf}
\ee

\index{Dirichlet boundary condition|)}
\begin{remark} \normalfont\label{angot1}
There are other finite volume schemes for problem \eqref{elliptic1D}.
\benum
\item For instance, it is possible, in Definition \ref{meshund}, to take $x_1
\ge 0$, $x_N \le 1$ and, for the definition of the scheme (that is 
\eqref{VF1}-\eqref{VF4}), to write \eqref{VF1} only for $i=2,\ldots,N-1$ and 
to replace \eqref{VF3} and \eqref{VF4} by $u_1=u_N=0$ (note that \eqref{VF2} 
does not change). For this so-called ``modified finite volume" scheme, it 
is also possible to obtain an error estimate as for the scheme 
\eqref{VF1}-\eqref{VF4} (see Remark \ref{h2conv}). Note that, with this 
scheme, the union of all control volumes for which the ``conservation 
law" is written is slightly different from $[0,1]$ (namely $[x_{3/2},x_{N-1/2}] \ne [0,1]$) \index{Finite volume scheme!-- for elliptic problems!-- in one space dimension|)}.

\item Another possibility is to take (primary) unknowns associated to the boundaries of the control volumes \citet{keller}, \citet{croisille}. 
We do not consider this case here.
\eenum
\end{remark}

\subsection{Comparison with a finite difference scheme}\index{Finite difference method} \label{fdvf}

With the same notations as in Section \ref{fvse1d}, consider that $u_i$ is now an approximation of $u(x_i)$. 
It is interesting to notice that the expression 
\[
\partd_i^2 u = \frac 1 h_i (F_{\ipud}-F_{\imud})= \frac 1 h_i \Bigl(- \frac {u_{i+1} - u_i} {h_{i+ \frac 1 2}} + \frac {u_{i} - u_{i-1}} {h_{i- \frac 1 2}}\Bigr)
\]
is not a consistent approximation  of $- u_{xx}(x_i)$ in the finite difference sense, that is the error made by replacing the derivative by a difference quotient (the truncation error \citet{DB}) does not tend to 0 as $h$ tends to 0. 
Indeed, let $\Ubar = \bigl(u(x_1),\ldots, u(x_N)\bigr)^t$; with the notations of (\ref{AUb})-(\ref{rhsdf}), the truncation error may be defined as
\[
r = A \Ubar - b,
\]

with $r=(r_1,\ldots,r_N)^t$.
Note that for $f$ regular enough, which is assumed in the sequel, $b_i = f(x_i) + 0(h)$. 
An estimate of $r$ is obtained by using Taylor's expansion:
\[
u(x_{i+1}) = u(x_i) + h_{i+ \frac 1 2} u_x(x_i) + \frac 1 2 h_{i+ \frac 1 2}^2 u_{xx}(x_i) + \frac 1 6 h_{i+ \frac 1 2}^3u_{xxx}(\xi_i),
\]
for some $\xi_i \in (x_i,x_{i+1})$, which yields   

\[
r_i = - \frac 1 h_i \frac
{h_{i+ \frac 1 2} + h_{i- \frac 1 2}} 2 u_{xx}(x_i) + u_{xx}(x_i) + 0(h),  \ \ 
i=1,\ldots,N, 
\]

which does not, in general tend to 0 as $h$ tends to 0  (except in
particular cases) as  may be seen on the simple following example: 

\begin{example}\normalfont
Let $f \equiv 1$ and consider a mesh of $(0,1)$, in the sense of Definition \ref{meshund}, satisfying $h_i = h$ for even $i$, $h_i = h / 2$ for odd $i$ and $x_i = (x_{i+1/2} + x_\imudt)/ 2$, for $i=1,\ldots,N$. 
An easy computation shows that the truncation error $r$ is such that

\[ 
r_i  = \begin{cases}
 	 -\frac 1 4,  \text{ for even  } i
		\\
 	+\frac 1 2,  \text{ for odd  } i. 
	\end{cases}
\]
Hence $\sup \{ \vert r_i \vert, i = 1,\ldots,N\}   \not \tends 0 $ as  $h\tends 0$. 
\end{example}

Therefore, the scheme  obtained from ~(\ref{VF1})-(\ref{VF4}) is not consistent \index{Consistency!-- in the finite difference sense} in the finite difference sense, even though it is consistent in the finite volume sense, that is,  the numerical approximation of the fluxes is conservative and the truncation error on the fluxes tends to 0 as $h$ tends to 0. 

If, for instance, $x_i$ is  the center of $K_i$,  for $i=1,\ldots,N$,
it is well known that for problem~(\ref{elliptic1D}), the consistent
finite difference scheme\index{Finite difference method} would be, omitting boundary conditions, 

\be \frac 4 {2 h_i + h_{i-1} + h_{i+1}}
\Bigl[- \frac {u_{i+1} - u_i}
 {h_{i+ \frac 1 2}} + \frac {u_{i} - u_{i-1}} {h_{i- \frac 1 2}} \Bigr] = f(x_i), \,
i=2,\ldots,N-1, \label{FD}
\ee

\begin{remark} \normalfont
Assume that  $x_i$ is, for $i=1,\ldots,N$, the center of $K_i$ and that the discrete 
unknown  $u_i$ of the finite volume scheme is considered as an approximation of 
the
mean value  $\tilde{u}_i$ of $u$ over $K_i$  
(note that 
$\tilde{u}_i =  u(x_i) + (h_i^2/24) u_{xx}(x_i) + 
0(h^3)$, if
$u \in C^3([0,1],\R)$) instead of $u(x_i)$, then again, the finite volume scheme, 
considered once
more as a finite difference scheme, is not consistent
 in the finite difference sense. 
Indeed, 
let
$\tilde{R}=A\tilde{U}-b$, with $\tilde{U}=(\tilde{u}_1,\ldots,\tilde{u}_N)^t$, and
$\tilde{R}=(\tilde{R}_1,\ldots,\tilde{R}_N)^t$, then, in general, 
$\tilde{R}_i$ does not go to $0$ as $h$ goes to $0$. In fact, it will be shown later 
that the finite volume scheme, when seen as a finite difference scheme, is
consistent  in the
finite difference sense if $u_i$ is considered as an approximation of
$u(x_i)-(h_i^2 / 8) u_{xx}(x_i)$. This is the idea
 upon which  the first proof of convergence by Forsyth and
Sammon in 1988 is based, see \citet{FS} and Section \ref{fdconvergence}.
 \end{remark} 

In the case of Problem~(\ref{elliptic1D}), both the finite volume and finite difference schemes  are convergent.  
The finite difference scheme (\ref{FD}) is convergent since it is stable, in the sense that
$\Vert X \Vert_\infty \le  C \Vert AX \Vert_\infty$, for all $X \in
\R^N$, where $C$ is a constant and $\Vert X \Vert_\infty=
\mathrm{sup}(\vert X_1 \vert,\ldots,\vert X_N \vert)$,
$X=(X_1,\ldots,X_N)^t$, and consistent 
 in the usual finite difference
sense. 
Since  $A(\overline{U}-U)=R$, the stability\index{Stability} property implies
that $\Vert \overline{U}-U \Vert_\infty \le  C \Vert R \Vert_\infty$ 
which goes to $0$, as $h$ goes to $0$, by definition of the
consistency 
in the finite difference sense. The convergence of the finite volume
scheme (\ref{VF1})-(\ref{VF4}) needs some more work and is described
in Section \ref{fvconvergence1D}.

\subsection{Comparison with a mixed finite element method}\index{Finite element method!-- mixed|(}  
The finite volume method has often be thought of as a kind of mixed finite element method, since both methods involve the fluxes.
However, we show here that, on the simple Dirichlet problem (\ref{elliptic1D}), the two methods yield two different schemes.
For Problem~(\ref{elliptic1D}), the discrete unknowns of the finite volume method are the values$u_i$, $i = 1,\ldots, N$.
The finite volume method also introduces one discrete unknown at each of the  control volume  extremities,
namely the numerical flux between the corresponding control volumes. 
And so indeed, the finite volume method for elliptic problems may appear closely related to the mixed finite element method.   
Recall that the mixed finite element method consists in introducing in Problem~(\ref{elliptic1D}) the auxiliary variable $q =  -u_x$, which yields the following system: 
\[\bac
q + u_x  =  0,
\\
q_x  =  f; 
\ean
\]
assuming  $f \in L^2((0,1))$,  a variational formulation of this system  is:
\be
q \in H^1((0,1)), \;  u \in L^2((0,1)), \label{efm0}
\ee
\be
\int_0^1  q(x)  p(x) dx  = \int_0^1 u(x)  p_x(x) dx, \, \forall \, p \in
H^1((0,1)), \label{efm1}
\ee
\be
\int_0^1  q_x(x) v(x) dx  =  \int_0^1 f(x) v(x)  dx, \, \forall \, v
\in L^2((0,1)). \label{efm2}
\ee
Considering an admissible mesh of $(0,1)$ (see Definition \ref{meshund}), the usual discretization of this variational formulation consists in taking the classical  piecewise linear  finite element functions for the approximation $H$ of $H^1((0,1))$ and the   piecewise
constant  finite element for the approximation  $L$ of $L^2((0,1))$.
Then, the discrete unknowns are $\{u_i,i=1,\ldots,N\}$ and $\{q_{i+ 1/2},i=0,\ldots,N\}$ ($u_i$ is an approximation of $u$ in 
$K_i$ and $q_{i+ 1/2}$ is an approximation of $-u_x(x_{i+ 1/2}$)). 
The discrete equations are obtained by performing a Galerkin expansion \index{Galerkin expansion} of $u$ and $q$ with respect to the  natural basis functions  $\psi_l$, $l = 1,\ldots,N$ (spanning $L$), and $\varphi_{j+1/2}$, $j= 0,\ldots,N$ (spanning $H$) and by taking $p = \varphi_{i+1/2}$, $i= 0,\ldots, N$ in (\ref{efm1}) and  $v = \psi_k, k= 1,\ldots,N$ in (\ref{efm2}). Let $h_0=h_{N+1}=0$, $u_0=u_{N+1}=0$ and  $q_{-1/2}=q_{N+3/2}=0$.
Then the discrete system obtained by the mixed finite element method has $2N  +1$ unknowns and reads 
\begin{align*}
 q_{i+\frac 1 2}(\frac {h_i+h_{i+1}} 3) + q_{i-\frac 1 2}(\frac {h_i} 6) + q_{i+\frac 3 2}(\frac {h_{i+1}} 6) = u_i - u_{i+1}, \, i=0,\ldots,N,  
\\
q_{i+\frac 1 2}-q_{i-\frac 1 2}=\int_{K_i} f (x) dx, \, i=1,\ldots,N.
\end{align*}
 Note that the unknowns $q_{i+1/2}$ cannot be eliminated from the system. 
The resolution of this system of equations does not give the same values 
$\{u_i,i=1,\ldots,N\}$ than those obtained by using the finite volume scheme~(\ref{VF1})-(\ref{VF4}). In fact it is  easily 
seen that,
in this case, the finite volume scheme can be  obtained from the mixed finite element 
scheme by
using the following numerical integration for the left handside of
(\ref{efm1}):  
\[
\int_{K_i} g(x) dx = \frac {g(x_{i+1}) + g(x_{i})} 2 h_i.
\]

This is also true for some two-dimensional elliptic problems and therefore the finite volume
 error estimates for these problems may be obtained via the mixed finite
element theory, see \citet{ABMO},   \citet{BMO}\index{Finite element method!-- mixed|)}.

\section{Convergence and error analysis for the Dirichlet problem}

\hpage{Elliptic 1d}

\index{Dirichlet boundary condition|(}

\subsection{Error estimate with $C^2$ regularity}\label{fvconvergence1D}

We shall now prove the following error estimate, which will be 
generalized to more general elliptic problems and in higher  space dimensions. 

\begin{theo}\label{eelapl1D}\index{Error estimate!-- for an elliptic equation!-- in the one dimensional case}~ 

Let $f \in C([0,1],\R)$ and let $u\in C^2([0,1],\R)$ be the (unique) solution of Problem
(\ref{elliptic1D}). Let $\mesh = (K_i)_{i=1,\ldots,N}$ be an admissible 
mesh in the sense of Definition \ref{meshund}.
Then, there exists a unique vector $U=(u_1,\ldots,u_N)^t \in \R^N$ solution to
(\ref{VF1})  -(\ref{VF4}) and there exists $C \ge 0$, only depending   on $u$, such
that

\be
\sum_{i=0}^{N} \frac {(e_{i+1}-e_i)^2} {h_{i+\frac 1 2}} \le C^2 h^2, \label{eel1}
\ee

and

\be
\vert e_i \vert \le C h, \; \forall i \in
\{1,\ldots,N\},  \label{eel2} 
\ee

with $e_0=e_{N+1}=0$ and $e_i=u(x_i)-u_i$, for all $i \in \{1,\ldots,N\}$.

\end{theo}

%
%

\proof

First remark that there exists a unique vector $U=(u_1,\ldots,u_N)^t \in \R^N$ solution to 
\eqref{VF1}-\eqref{VF4}. Indeed, multiplying \eqref{VF1}
by $u_i$ and summing for $i=1,\ldots,N$ gives

\[
\frac {u_1^2} {h_{\frac 1 2}} + \sum_{i=1}^{N-1} \frac {(u_{i+1}-u_i)^2} 
{h_{i+\frac 1 2}} + \frac {u_N^2} {h_{N + \frac 1 2}}= \sum_{i=1}^N 
u_ih_if_i.
\]

Therefore, if $f_i=0$ for any $i \in \{1,\ldots,N\}$, then the unique solution to (\ref{VF1}) is obtained by taking $u_i=0$, for any
$i \in \{1,\ldots,N\}$. This gives existence and uniqueness of $U=(u_1,\ldots,u_N)^t \in \R^N$ solution to 
\eqref{VF1} (with \eqref{VF2}-\eqref{VF4}).

\vspace{0.5cm}

One now proves \eqref{eel1}.
Let 

\[
 \Fbar_{i+\frac 1 2} = -  u_x(x_\ipud), \; i = 0,\ldots,N,
\]

Integrating the equation $-u_{xx} = f$ over $K_i$ yields

\[
 \Fbar_{i+\frac 1 2} - \Fbar_{i-\frac 1 2} = h_i f_i,
\; i =1, \ldots, N.
\]

By (\ref{VF1}), the numerical fluxes $F_{i+\frac 1 2}$ satisfy 

\[
 F_{i+\frac 1 2} - F_{i-\frac 1 2} = h_i f_i,
\; i = 1, \ldots, N.
\]

Therefore, with $G_{i+\frac 1 2}=\Fbar_{i+\frac 1 2} -  F_{i+\frac 1 2}$,

\[
 G_{i+\frac 1 2} - G_{i-\frac 1 2} = 0, \, i  =1, \ldots, N.
\]

Using the consistency of the fluxes \eqref{consistflux}, there exists $C > 0$, 
 only depending  on $u$, such that

\be
\Fast_{i+\frac 1 2} =   \Fbar_{i+\frac 1 2} +   R_{i+\frac 1 2} \text{ 
and }   
\abs{R_{i+\frac 1 2}}.
\le C h, \label{errconst} 
\ee

Hence with $e_i = u(x_i) -u_i$, for $i = 1,\ldots,N$, and $e_0=e_{N+1}=0$, one has 

\[
G_{i+\frac 1 2} = - \frac {e_{i+1} -e_i} {h_\ipud} - R_{\ipud}, \, i =
0,\ldots,N,
\]

so that $(e_i)_{i=0,\ldots,N+1}$ satisfies
 
\be 
- \frac {e_{i+1} -e_i} {h_\ipud}  -  R_{\ipud} +   \frac {e_{i} -e_{i-1}} 
{h_\imud} + R_{\imud}=0, \, \forall i \in \{1,\ldots,N\}.
\label{eqerror}
\ee

Multiplying (\ref{eqerror}) by $e_i$ and summing over
$i=1,\ldots, N$ yields 

\[
-\sum_{i=1}^N \frac {(e_{i+1}-e_i)e_i} {h_{i + \frac 1 2}} +\sum_{i=1}^N \frac {(e_i-
e_{i-1})e_i}
{h_{i - \frac 1 2}} = -\sum_{i=1}^N R_{\imud} e_i + \sum_{i=1}^N R_{\ipud} 
e_i.
\]

Noting that $e_0 = 0$, $e_{N+1} = 0$ and reordering by parts, this yields 
(with (\ref{errconst}))

\be
\sum_{i=0}^N \frac {(e_{i+1}-e_i)^2} {h_{i + \frac 1 2}} \le
Ch\sum_{i=0}^N  \abs{e_{i+1}-e_i}. 
\label{caspie1}
\ee

The Cauchy-Schwarz  inequality applied to the right hand side gives 

\be
\sum_{i=0}^N \vert e_{i+1}-e_i \vert  \le
\Bigl(\sum_{i=0}^N \frac {(e_{i+1}-e_i)^2} {h_{i + \frac 1 2}}\Bigr)^{\frac 1
2}  \Bigl(\sum_{i=0}^N  {h_{i + \frac 1 2}}\Bigr)^{\frac 1 2}. 
\label{caspie2}
\ee

Since $\dsp{\sum_{i=0}^N  {h_{i + \frac 1 2}} = 
1}$ in \eqref{caspie2} and from \eqref{caspie1}, one deduces  \eqref{eel1}.

Since, for all $i \in \{1,\ldots,N\}$, $e_i=\dsp{\sum_{j=1}^i(e_j-e_{j-1})}$, one can deduce,
from \eqref{caspie2} and \eqref{eel1} that \eqref{eel2} holds. 
\eof

\begin{remark} \normalfont
The error estimate given in this section
does not use the discrete maximum principle (that is the fact that $f_i \ge 0$, 
for all $i=1,\ldots,N$, implies $u_i \ge 0$, 
for all $i=1,\ldots,N$), which is used in the proof of error estimates by 
the finite difference techniques, but the coerciveness of the elliptic 
operator, as in the proof of error estimates by 
the finite element techniques.
\end{remark}

\begin{remark} \normalfont\label{h2conv}~

\begin{enumerate}

\item The above proof of convergence gives an error estimate of order $h$. 
It is sometimes possible to obtain an error estimate of order $h^2$. Indeed, 
this is the case, at least if $u \in
C^4([0,1],\R)$, if $x_i$ is the center of $K_i$ for all $i=1,\ldots,N$. 
One obtains, in this case, $\vert e_i \vert \le C h^2$, for all $i \in
\{1,\ldots,N\}$, where $C$ only depends on $u$ (see \citet{FS}).

\item It is also possible to obtain an error estimate for the modified 
finite volume scheme described in the first item of Remark \refp{angot1}. 
It is even possible to obtain an error estimate of order $h^2$ in the case $x_1=0$, 
$x_N=1$ and assuming that $x_{i + 1/2}= (1/2)(x_i+x_{i+1})$, for all 
$i=1,\ldots,N-1$. In fact, in this case, one obtains $\vert R_{i+1/2} 
\vert \le C_1 h^2$, for all $i=1,\ldots,N-1$. Then, the proof of Theorem 
\ref{eelapl1D} gives \eqref{eel1} with $h^4$ instead of $h^2$ which yields 
$\vert e_i \vert \le C_2 h^2$, for all $i \in \{1,\ldots,N\}$ (where 
$C_1$ and $C_2$ are only depending  on $u$). 
Note that this modified finite volume scheme is also consistent in the finite difference sense. 
Then, the finite difference techniques yield also an 
error estimate on $\vert e_i \vert$, but only of order $h$. 

\item  It could be tempting to try and find error estimates with respect 
to the mean value of the exact solution on the control volumes 
rather than with respect to its 
value at some point of the control volumes. This is not such a good idea: 
indeed, if   
$x_i$ is not the center of $K_i$ (this will be the general case in 
several space dimensions), then  one does not have 
(in general) $\vert \tilde{e}_i \vert \le C_3 h^2$ (for some $C_3$ only 
depending on $u$) with $\tilde{e}_i=\tilde{u}_i-u_i$ where 
$\tilde{u}_i$ denotes the mean value of $u$ over $K_i$.
\end{enumerate}
\end{remark}

\begin{remark} \normalfont~ \label{rem:integral-remainder}

\begin{enumerate}

\item  If the assumption $f \in C([0,1],\R)$ is replaced by the assumption $f \in L^2((0,1))$ in Theorem \ref{eelapl1D}, then $u \in H^2((0,1))$ instead of $C^2([0,1],\R)$, but the estimates of  Theorem \ref{eelapl1D} still hold. 
In this case, the consistency  of the fluxes  must be obtained with a Taylor expansion with an integral remainder. 
This is feasible for $C^2$ functions, and since the remainder only depends on the $H^2$ norm, a density argument allows to conclude; see also Theorem \refp{testh2} below and \citet{cvnl}. 

\item  If the assumption $f \in C([0,1],\R)$ is replaced by the assumption 
$f \in L^1((0,1))$ in Theorem \ref{eelapl1D}, then $u \in C^2([0,1],\R)$ no longer holds and neither does $u \in H^2((0,1))$, but the convergence still holds; 
indeed there exists $C(u,h)$,  only depending on $u$ and $h$, such that $C(u,h) \tends 0$, as $h \tends 0$, and  $\vert e_i \vert \le C(u,h)$, for all $i= 1,\ldots,N$. 
The proof is similar to the one above, except that the estimate (\ref{errconst}) is replaced by $ \vert R_{i+1/2} \vert
\le C_1(u,h)$, for all $i= 0,\ldots,N$, with some $C_1(u,h)$,  only
depending on $u$ and $h$, such that $C(u,h) \tends 0$, as $h \tends 0$. 
\end{enumerate} 
\end{remark}

\begin{remark} \normalfont
Estimate (\ref{eel1}) can be interpreted as a ``discrete
$H^1_0$" estimate on the error. A theoretical result which underlies the
$L^\infty$  estimate  (\ref{eel2}) is the fact that if $\O$ is an open bounded
subset of $\R$, then 
$H^1_0(\O)$ is imbedded in $L^\infty(\O)$. This is
no longer true in higher dimension. In two space dimensions, for instance,  a
discrete version of the imbedding of $H^1_0$ in $L^p$ allows to obtain (see
e.g. \citet{thesejmf}) $\Vert e \Vert_p \le Ch$,  for all finite $p$, 
which in turn yields  $\Vert e \Vert_\infty \le Ch \ln h$ for convenient 
meshes (see Corollary \refp{norminf}).
\end{remark}

The important features needed for the above proof 
seem to be the consistency  of the approximation of the fluxes and the
conservativity 
 of the scheme; this conservativity is natural
 the fact that the scheme is obtained by integrating the
equation over each cell, and the approximation of the flux on any interface is obtained by
taking into account the flux balance (continuity of the flux in the case of no source term 
on
the interface).

The above proof generalizes to other elliptic problems, such as a 
convection-diffusion\index{Equation!Convection-diffusion --}
 equation of
the form $-u_{xx} + a u_x +b u =f$, and to equations of the form $-(\l 
u_x)_x = f$ where
$\l \in L^\infty$ may be discontinuous, and is such that there exist $\a$ and $\b$
in $\R_+^\star$ such that $\a \le \l \le \b$. These generalizations are studied in the next section. 
Other generalizations include similar problems in 2 (or 3) space dimensions, with meshes 
consisting of rectangles (parallepipeds), triangles (tetrahedra), or general meshes of
Vorono\"{\i} type, and the corresponding evolutive (parabolic) problems.  These 
generalizations
will be addressed in further chapters. 

Let us now give a proof of Estimate \eqref{eel2}, under slightly different conditions, which
uses finite difference techniques. 

\subsection{An error estimate using a finite difference technique}\label{fdconvergence}

Convergence can be obtained via a method
similar to that of the finite difference\index{Finite difference method} proof of convergence (following, for
instance, \citet{FS}, \citet{WM}, \citet{Fa1}). Most of these methods, are, however,
limited to the finite volume method for 
Problem (\ref{elliptic1D}). Using the  notations of Section \ref{fdvf}  
(recall that $\Ubar = (u(x_1),\ldots,u(x_N))^t$, and $r = A \Ubar - b = 0(1)$), 
the idea is to find
$\Ubarbar$ ``close" to $\Ubar$, such that 
\[
A\Ubarbar = b + \rbarbar, \text{ with }\rbarbar = 
0(h).
\] 
This value of $\Ubarbar$ was found in \citet{FS} and is such that 
$\Ubarbar = \Ubar - V$, where 
\[
V =  (v_1,\ldots,v_N)^t \text{ and } v_i =  \frac {h_i^2 u_{xx}(x_i)} 8, 
\, i=1,\ldots,N.
\] 
Then, one may decompose the 
truncation error as
\[
r = A (\Ubar-U) = A V  + \rbarbar \text{ with } \Vert V
\Vert_\infty = 0(h^2) \text{ and } \rbarbar  = 0(h).
\]

The existence of
such a $V$ is given in Lemma \ref{V}.  In order to prove the convergence of the scheme,
a stability property  is established in Lemma \ref{stab}.   

\begin{lemma} \label{V}
Let $\mesh = (K_i)_{i=1,\cdots,N}$ be an admissible mesh of $(0,1)$, in the sense of 
Definition \refp{meshund}, such that $x_i$ is the center of $K_i$ for all 
$i=1,\ldots,N$. Let $\a_\mesh >0$ be such that $h_i > \a_\mesh h$ for all 
$i=1,\ldots,N$ (recall that $h = \max \{h_1,\ldots,h_N\}$). Let $\Ubar$
$=$ $(u(x_1),\ldots,u(x_N))^t$ $\in$ $\R^N$, where $u$ is the solution to
(\ref{elliptic1D}), and assume $u \in C^3([0,1],\R)$. 
 Let $A$ be the matrix defining the numerical scheme, given in
\refep{matrix}. 
Then there exists a unique $U=(u_1,\ldots,u_N)$ solution of
(\ref{VF1})-(\ref{VF4}) and there exists $\rbarbar$ and
$V \in \R^N$ such that
\[
r = A (\Ubar-U) = A V  + \rbarbar, \text{ with } \Vert V
\Vert_\infty \le Ch^2 \text{ and } \Vert \rbarbar \Vert_\infty \le Ch,
\] 

where $C$ only depends on $u$ and $\a_\mesh$.
\end{lemma}

\proof of Lemma \ref{V}

\smallskip

The existence and uniqueness of $U$ is classical (it
is also proved in Theorem \ref{eelapl1D}). 

For $i = 0,\ldots N$, define

\[
R_{i+\frac 1 2} = - \frac {\uexip -\uexi} {h_{i+\frac 1 2}} + 
u_x (x_{i+\frac 1 2}).
\]

Remark that 

\be
r_i = \frac 1 {h_i}(R_{i+\frac 1 2} -R_{i-\frac 1
2}), \text{ for } i=0,\ldots,N,
\label{aieaie}
\ee

where $r_i$ is the $i-$th component of  $r = A (\Ubar-U) $. 

The computation of  $R_{i+\frac 1 2}$ yields

\[
\bal
R_{i+\frac 1 2}    = -  \frac 1 4 (h_{i+1}-h_i)  u_{xx}(x_{i+\frac 1 2}) + 0(h^2) , \,
i=1,\ldots, N-1, \\
R_{\frac 1 2}   =  - \frac 1 4 h_1  u_{xx}(0) + 0(h^2), \, R_{ N + \frac 1 2}=
\frac 1 4 h_N u_{xx}(1) + 0(h^2). 
\ean
\]

Define $V=(v_1,\ldots,v_N)^t$ with $v_i=  \frac {h_i^2 u_{xx}(x_i)} {8}$,
$i=1,\ldots,N$. Then,

\[
- \frac {v_{i+1} - v_i} {h_{i+ \frac 1 2}}  =  R_{i+\frac 1 2} + 0(h^2), \,
i=1,\ldots,N-1,
\]
\[
- \frac {2 v_1} {h_1}  = R_{\frac 1 2} + 0(h^2),
\]
\[
\frac {2 v_N} {h_N}  =  R_{N + \frac 1 2} + 0(h^2).
\]

Since $h_i \ge \a_\mesh h$, for $i = 1,\ldots,N$, replacing $R_{i+\frac 1 2}$ in
\eqref{aieaie}  gives that
$r_i = (AV)_i + 0(h)$, for $i = 1,\ldots,N$, and  $\Vert V 
\Vert_\infty = 0(h^2)$.  Hence the lemma
is proved. 
\eof

\begin{lemma}[Stability]\label{stab}\index{Stability}
Let $\mesh = (K_i)_{i=1,\cdots,N}$  be an admissible mesh of $[0,1]$ in the sense of Definition 
\ref{meshund}. Let  $A$ be the matrix
defining the finite volume scheme given in
\eqref{matrix}. Then $A$ is invertible and  
\be \Vert A \Vert_\infty^{-1} \le \frac 1 4.
\label{eqstab} \ee 
\end{lemma}

\proof of Lemma \ref{stab}

\smallskip

First we prove a discrete maximum principle; indeed  if $b_i  \ge 0$,
for all $i = 1, \ldots,N$, and if $U$ is solution of $A U = b$ then we prove that $u_i \ge 0$ 
for all $i = 1,\ldots,N$.

Let $a= \min\{u_i, i= 0,\ldots, N+1\}$ (recall that $u_0=u_{N+1}=0$) and 
$i_0 = \min\{i \in \{0,\ldots,N+1\}$; $u_i = a\}$.

If $i_0 \ne 0$ and $i_0 \ne N+1$, then

\[
\frac 1 {h_{i_0}} \Bigl( \frac { u_{i_0} - u_{i_0 -1}} {h_{i_0-\frac 1 2}} - 
 \frac  { u_{i_0+1} - u_{i_0} } {h_{i_0+\frac 1 2}} \Bigr) =b_{i_0} \ge 0,
\]

this is impossible since $u_{i_0+1} - u_{i_0} \ge 0$ and $u_{i_0} - 
u_{i_0 -1} < 0$, by definition of $i_0$. Therefore,
$i_0 = 0$ or $N+1$. Then, $a=0$ and $u_i \ge 0$ for all $i =1,\ldots,N$. 

Note that, by linearity, this implies that $A$ is invertible. 

Next, we shall prove that there exists $M >0$ such that $\Vert A^{-1} \Vert_\infty \le M$ (indeed, $M=1/4$ is convenient). 
Let $\phi$ be defined on $[0,1]$ by $\phi(x) = \frac 1 2 x(1-x)$. 
Then $- \phi_{xx}(x) =1$ for all $x \in [0,1].$ 
Let $\Phi= (\phi_1, \ldots,\phi_N)$ with $\phi_i = \phi(x_i)$; if $A$ 
represented the
usual finite difference approximation of the second order derivative, then we would have 
$A \Phi =
1$, since the difference quotient approximation of the second order derivative 
of a second 
order
polynomial is exact ($\phi_{xxx}=0$). Here, with the finite volume scheme
(\ref{VF1})-(\ref{VF4}), we have $A \Phi - \mathbf{1}= A W$ 
(where $\mathbf{1}$  denotes the vector of $\R^N$ the components of which are all equal to 1), with  $W 
=
(w_1,\ldots,w_N) \in \R^N$ such that $W_i = -\frac {h_i^2} 8$ (see proof 
of Lemma \ref{V}). Let $b \in \R^N$ and 
$A U = b$,
since  $A (\Phi -  W) = \mathbf{1}$, we have 
\[A(U - \Vert b \Vert_\infty(\Phi -  W)) \le 0,\] 
this  last
inequality being meant componentwise. Therefore, by the above maximum principle, 
assuming, without loss of generality, that $h \le 1$, one has
\[
u_i \le 
{\Vert b \Vert_\infty} (\phi_i - w_i), \text{ so that } u_i \le \frac {\Vert b
\Vert_\infty} 4.
\]
(note that $\phi(x) \le \frac 1 8)$.  But we also have 
\[A(U + \Vert b \Vert_\infty (\Phi -  W)) \ge  0,\]
and again by the maximum principle, we obtain
\[
u_i \ge  - \frac {\Vert b
\Vert_\infty}  4.
\] 
Hence
$\Vert U \Vert_\infty \le \frac 1 4 \Vert b \Vert_\infty$. 
This shows that $\Vert A^{-1}  \Vert_\infty \le \frac 1 4$. 
\eof

This stability result, together with
the existence of $V$ given by Lemma \ref{V}, yields the convergence of the finite volume 
scheme, formulated in the next theorem.

\begin{theo}\index{Error estimate!-- for an elliptic 
equation!-- in the one dimensional case}
 Let  $\mesh = (K_i)_{i=1,\cdots,N}$  be an admissible mesh of $[0,1]$ in the sense of 
Definition \refp{meshund}. 
Let $\a_\mesh \in \R_+^\star$ be such that $h_i \ge \a_\mesh h$, for all $i = 1,\ldots,N$ (recall that $h = \max \{h_1,\ldots,h_N\}$). 
Let $\Ubar = (u(x_1),\ldots,u(x_N))^t \in \R^N$, and assume $u \in  C^3([0,1],\R)$
(recall that $u$ is the solution to (\ref{elliptic1D})).
Let $U=(u_1,\ldots,u_N)$ be the solution given by the numerical scheme
(\ref{VF1})-(\ref{VF4}). 
Then there exists $C > 0$, only depending on $\a_\mesh$ and $u$, such that $\Vert U - \Ubar \Vert_\infty \le C h$. 
\end{theo}

\begin{remark} \normalfont\label{rubarbare}
In the proof of Lemma \ref{stab}, it was shown that 
$A(\overline{U}-V)=b+0(h)$; 
therefore, if,
once again,  the finite volume scheme is considered as a finite difference
scheme\index{Finite difference method}, it is consistent, in the finite difference sense, 
when $u_i$ is
considered to be an approximation of $u(x_i)- (1/8)h_i^2 
u_{xx}(x_i)$. 
\end{remark}

\begin{remark} \normalfont\label{wc} With the notations of Lemma \ref{V},
let ${r}$ be the function defined by 
\[
{r}(x)=r_i, \ \  \text{ if } x \in K_i, \ \ i=1,\ldots,N,
\]
the function ${r}$ does not necessarily go to $0$ (as $h$ goes to
$0$) in the $L^\infty$ norm (and even in the $L^1$ norm), but, thanks
to the conservativity 
 of the scheme, it goes to $0$ in
$L^\infty((0,1))$ for the weak-$\star$ topology\index{Convergence!-- for the weak star topology}, that is
\[
\int_0^1 r(x)\varphi(x)dx \rightarrow 0, \ \ \mathrm{as} \ \  h \rightarrow 0, \ \ \forall
\varphi \in L^1((0,1)).
\]
This property  will be called ``weak consistency"\index{Consistency!Weak --} in the sequel and may also be used to 
prove the
convergence of the finite volume scheme (see \citet{Fa1}).
\end{remark} 

The proof of convergence described above may be easily generalized to the two-dimensional Laplace
equation $-\Delta u = f$ in two and three space  dimensions if  a rectangular or a
parallepipedic mesh is used, provided that the solution $u$ is of class $C^3$. However, it does
not seem to be easily generalized to other types of meshes.
\index{Dirichlet boundary condition|)}

\section{General 1D elliptic equations} \label{ell1D}

\hpage{Elliptic 1d}

\subsection{Formulation of the finite volume scheme}
 This section is devoted to the formulation and to the proof of convergence  of a finite volume scheme  for a one-dimensional linear 
convection-diffusion equation, with a discontinuous diffusion coefficient.  The scheme can be  generalized
in the two-dimensional and three-dimensional cases (for a space discretization 
which uses, for instance, simplices or
parallelepipedes or   a ``Vorono\"{\i} mesh", see Section 
\refp{vfquatre}) and to other boundary conditions.

~

Let $\l \in L^\infty((0,1))$ such that there exist $\lund$ and $\lbar \in \R^\star_+$ 
with $\lund  \le \l \le \lbar$ a.e. and let $a, b, c, d \in \R$, with $b \ge 0$,  and $f \in
L^2 ((0,1))$. The aim, here, is to find an approximation to the solution, $u$, 
of the following  problem: 
\be
 -(\l u_x)_x (x) + a u_x(x) + bu(x) = f(x), \; x \in [0,1],
\label{eq1}
\ee
\be
 u(0)=c, \; u(1)=d. 
\label{cl}
\ee

The discontinuity of the
coefficient $\l$ may arise for instance for  the permeability of a porous medium,
the ratio between the permeability of sand and the permeability of clay being of an order of
$10^3$;  heat conduction in a heterogeneous medium can also yield such discontinuities,
since the conductivities of the different components of the medium may be quite
different.  Note that  the assumption $b \ge 0$ ensures the existence of the
solution to the problem. \vspace*{.1cm}
 
\begin{remark} \normalfont
Problem (\ref{eq1})-(\ref{cl}) has a unique solution $u$ in the Sobolev
space  $H^1((0,1))$. This solution is continuous (on $[0,1]$) but is not, in
general, of class $C^2$ (even if $\l(x)=1$, for all $x \in [0,1]$). Note 
that one has $-\l u_x(x) = \int_0^x g(t) dt + C$,
where $C$ is some constant and $g=f-au_x-bu \in L^1((0,1))$, so that $\l u_x$ is a continuous function and $u_x  \in
L^\infty((0,1))$. 
\end{remark}

\vspace*{.1cm}

Let $\mesh = (K_i)_{i=1,\cdots,N}$ be an admissible mesh, in the sense of 
Definition \refp{meshund}, such that the
discontinuities of $\l$ coincide with the interfaces of the mesh.  

The notations being the
same as in section \ref{Diric1D}, integrating Equation (\ref{eq1})
over $K_i$ yields

\[
-(\l u_x)(x_{i+\frac 1 2}) + (\l u_x)(x_{i-\frac 1 2})  + a u(x_\ipud) - a u(x_\imud) +
\int_{K_i} b u(x) dx = \int_{K_i} f(x) dx, \ \ i=1,\ldots,N.
\]

Let $(u_i)_{i = 1,\cdots,N}$ be the discrete unknowns. In the case $a \ge 0$, which will be
considered in the sequel, the convective term
$a u(x_{i+1/2})$ is approximated by $a u_i$ (``upstream\index{Upstream}") because of
stability considerations\index{Stability}. Indeed, this choice always yields a stability 
result whereas the approximation of $a u(x_{i+1/2})$ by 
$(a/2)(u_i+u_{i+1})$ (with the approximation of the other terms as it is 
done below) yields a stable scheme if $ah \le 2\l$, for a 
uniform mesh of size $h$ and a constant diffusion coefficient $\l$. The case $a  \le 0$ is
easily handled in the same way by approximating $a u(x_{i+1/2})$ by $a u_{i+1}$. The term
$\int_{K_i} b u(x) dx$ is  approximated by
$ b h_i u_i$. Let us now turn to the approximation $H_{i+1/2}$ of $- \l
u_x(x_{i+1/2})$.  Let $\l_i = \frac 1 {h_i} \int_{K_i} \l(x) dx$;  since
$\l\vert_{K_i} \in C^1(\bar{K_i})$, there exists $c_\l \in \R_+$, only 
depending on $\l$, such that  $\vert \l_i -
\l(x) \vert \le c_\l h$, $\forall \, x \in K_i$. In order that the scheme be
conservative, the discretization of the flux at $x_{i+1/2}$ should have the same
value on $K_i$ and $K_{i+1}$.  To this purpose, we introduce the auxiliary unknown
$u_{i+1/2}$
(approximation of $u$ at $x_{i+1/2}$). Since on $K_{i}$ and  $K_{i+1}$, $\l$
is continuous, the approximation of $- \l u_x$  may be performed on each side of 
 $x_{i+1/2}$ by using the finite difference principle:
 
\[
H_{i+\frac 1 2}   =  - \l_i  \frac {u_{i+ \frac 1 2} - u_i} {\hip} \text{
on } K_i, \, i=1,\ldots,N,
\]
\[
H_{i+\frac 1 2}  = -  \l_{i+1} \frac {u_{i+1} - u_{i+ \frac 1 2}} {\hipum}
\text{ on } K_{i+1}, \, i=0,\ldots,N-1,
\]

with $u_{1/2}=c$, and $u_{N +1/2}=d$,
for the boundary conditions. (Recall that $\hip = x_{i+1/2}-x_i$ and $\him =
x_i-x_{i - 1/2}$). Requiring the two above approximations of  $\l 
u_x(x_{i+1/2})$ to be equal (conservativity  of the flux) yields the value of 
$u_{i+1/2}$ (for $i=1,\ldots,N-1$):

\be
u_{i+ \frac 1 2} = \frac 
{  \dsp{ u_{i+1} \frac{\l_{i+1}} {\hipum} + u_i \frac{\l_{i}} {\hip}} }
{  \dsp{  \frac{\l_{i+1}} {\hipum} +   \frac{\l_{i}} {\hip} } }
\label{uipud}
\ee

which, in turn, allows to give the expression of the approximation
$H_{i+ \frac 1 2}$ of $\l u_{x}(x_{i + \frac 1 2})$: 

\be 
\bal
H_{i+ \frac 1 2}  = - \tau_\ipud 
(u_{i+1} - u_i), \, i=1,\ldots,N-1, \\
H_{ \frac 1 2}  = -  \dsp{\frac {\l_{1}} { h_1^-} ( u_1 -c),}  \\
H_{ N+\frac 1 2}  = -  \dsp{\frac {\l_{_N}} { h_N^+} (d-u_N) } 
\ean
\label{vfsig}
\ee

with

\be \tau_\ipud  = \frac { \l_{i} \l_{i+1}}
 { \hip \l_{i+1} + \hipum  \l_{i}}, \, i=1,\ldots,N-1. \label{tau1D}
\ee 

\begin{example}\normalfont If $h_i = h$, for all $i \in \{1,\ldots,N\}$, and $x_i$ is assumed to be
the center of $K_i$, then $\hip = \him = \frac h 2$, so that 
\[ 
H_{i+ \frac 1 2}   =  - \frac { 2 \l_{i} \l_{i+1}}
 {\l_{i} + \l_{i+1}} \frac{ u_{i+1} - u_i} h,  
\] 
and therefore the mean harmonic value of $\lambda$ is involved.  
\end{example}
The numerical scheme \index{Finite volume scheme!-- for 
elliptic problems!-- in one space dimension} 
for the approximation of Problem (\ref{eq1})-(\ref{cl}) is
therefore,
\be
 F_{i+\frac 1 2} - F_{i-\frac 1 2} + b h_i u_i = h_i f_i, \; \forall i \in
\{1, \ldots, N\}, \label{eqd1}
\ee

with $f_i = \frac 1 {h_i}\int_{x_{i-\frac 1 2}}^{x_{i+\frac 1 2}}f(x) dx$,  for
$i=1,\ldots, N$, and where $(F_{i+\frac 1 2})_{i \in \{0,\ldots,N\}}$ 
is defined by
the following expressions 

\be
 F_{i+\frac 1 2}= -  \tau_\ipud (u_{i+1}-u_i)  + a u_i, \; \forall i \in \{1, \ldots, N-1\},
\label{eqd2}
\ee
\be
F_{\frac 1 2} = - \frac {\l_{1}} { h_1^-} (u_1-c)  + a
c, \;
 F_{N + \frac 1 2} = -  \frac {\l_{_N}} { h_N^+} (d-u_N) + a u_N.
\label{cld}
\ee

\begin{remark} \normalfont In the case $a \ge 0$, the choice of the approximation 
 of $a u(x_{i + 1/2})$ by  $a u_{i+1}$ would yield an
unstable scheme, except for $h$ small enough (when $a \le 0$, the unstable scheme is $a u_i$). 
\end{remark} 

Taking \eqref{tau1D}, (\ref{eqd2}) and (\ref{cld}) into account, the numerical scheme (\ref{eqd1})
yields a system of $N$ equations with $N$ unknowns $u_1,\ldots,u_N$.

\subsection{Error estimate}\label{ce1r}

\begin{theo}\label{ee1dell}\index{Error estimate!-- for an elliptic 
equation!-- in the one dimensional case}~

Let $a, b \ge 0$, $c, d \in \R$, $\l \in L^\infty((0,1))$ such that 
$\lund \le \l \le
\lbar$ a.e. with some $\lund$,  $\lbar \in \Rpast$ and $f \in L^1((0,1))$. Let $u$ be
the (unique) solution of (\ref{eq1})-(\ref{cl}). Let 
$\mesh=(K_i)_{i=1,\cdots,N}$ be an admissible mesh, in the sense of Definition \ref{meshund}, such 
that $\l \in C^1(\cvbar_i)$
and $f \in C(\cvbar_i)$, for all $i = 1,\cdots,N$. Let $\gamma = \max\{\Vert u_{xx}
\Vert_{L^\infty(K_i)}, i = 1,\cdots,N\}$ and $\d = \max\{\Vert \l
\Vert_{L^\infty(K_i)}, i = 1,\cdots,N\}$. Then, 

\begin{enumerate}
\item there exists a unique vector
$U=(u_1,\ldots,u_N)^t \in \R^N$ solution to (\ref{tau1D})-(\ref{cld}), 
\item there
exists $C$, only depending  on $\lund, \lbar, \gamma$ and $\d$, such
that
\be
\sum_{i=0}^{N} \tau_\ipud (e_{i+1}-e_i)^2  \le C h^2, \label{ee1} 
\ee
where $\tau_\ipud$ is defined in (\ref{tau1D}), and

\be
\vert e_i \vert \le C  h, \, \forall i \in
\{1,\ldots,N\},  \label{ee2} 
\ee

with $e_0=e_{N+1}=0$ and $e_i=u(x_i)-u_i$, for all $i \in \{1,\ldots,N\}$.

\end{enumerate}
\end{theo}

\proof of Theorem \ref{ee1dell}
 
\smallskip

\textbf{Step 1. Existence and uniqueness of the solution to 
(\ref{tau1D})-(\ref{cld}).}

Multiplying (\ref{eqd1})
by $u_i$ and summing for $i=1,\ldots,N$ yields that if $c=d=0$
and $f_i=0$ for any $i \in \{1,\ldots,N\}$, then the unique solution to
(\ref{tau1D})-(\ref{cld}) is obtained by taking $u_i=0$, for any
$i \in \{1,\ldots,N\}$. This yields existence and uniqueness of the solution to 
(\ref{tau1D})-(\ref{cld}).

\medskip

\textbf{Step 2. Consistency of the fluxes.\index{Consistency!-- of the fluxes}}

Recall that $h=\max\{h_1,\ldots,h_N\}$. 
Let us first show  the consistency of the fluxes.

Let $\Hbar_{i+1/2} = - 
(\l u_{x})(x_{i+1/2})$ and $\Hast_{i+1/2} = - 
\tau_{i+1/2} (u(x_\ipu) - u(x_i))$, for $i=0,\ldots,N$, with 
$\tau_{1/2}=\l_1/h_1^-$ and $\tau_{N+1/2}=\l_N/h_N^+$. Let us first
show that there exists $C_1 \in \R_+^\star$,  only depending  on $\lund, \lbar, \gamma$ and
$\d$, such that

\be 
\bal
\Hast_\ipud = \Hbar_\ipud + T_\ipud,\\
\vert T_\ipud \vert \le C_1 h, \, i=0,\ldots,N. 
\ean 
\label{lambdaconsist}
\ee

In order to show this, let us introduce 
\be
\Hastun_\ipud = - \l_i \frac {u(x_\ipud) - u(x_i)} \hip \text{ and } 
\Hastde_\ipud =  - \l_\ipu \frac {u(x_\ipu)- u(x_\ipud) } \hipum;
\label{Gast}
\ee

since $\l \in C^1(\bar{K_i})$, one has $u \in C^2(\bar{K_i})$; hence, 
there exists $C \in
\Rpast$,  only depending  on $\g$ and $\d$,  such that 

\be
\Hastun_\ipud  = \Hbar_\ipud  + R_{\ipud}^{-}, \text{ where }   \vert 
R_{\ipud}^{-} \vert
\le C h, \, i=1,\ldots,N, 
\label{eqRun}
\ee 

and 

\be
\Hastde_\ipud  =  \Hbar_\ipud  + R_{\ipud}^{+}, \text{ where }  \vert 
R_{\ipud}^{+} \vert
\le C h, \, i=0,\ldots,N-1. 
\label{eqRde}
\ee  

This yields \eqref{lambdaconsist} for $i=0$ and $i=N$.

The following equality:
\be
 \Hbar_{i+\frac 1 2}=  \Hastun_\ipud -   R_{\ipud}^{-} =  \Hastde_\ipud 
 -   R_{\ipud}^{+}, \, i=1,\ldots,N-1,
\label{Gbarast}
\ee
yields that 

\be
u(x_\ipud) = \frac { \dsp 
{\frac {\l_\ipu} {\hipum} u(x_\ipu) +  \frac {\l_i} {\hip} 
u(x_i)}  } { \dsp 
{ \frac {\l_i} {\hip}  + \frac {\l_\ipu} {\hipum}}  } + S_\ipud, \, 
i=1,\ldots,N-1,
\label{lambdauipud}
\ee
where 
\[
S_\ipud = \frac { R_{\ipud}^{+}  -  R_{\ipud}^{-}} 
{ \frac {\l_i} {\hip}  + \frac {\l_\ipu} {\hipum} }
\]
so that 
\[
\vert S_\ipud \vert \le \frac 1 \lund \frac  { \hip \hipum} {\hip + \hipum}
 \vert R_{\ipud}^{+}  -  R_{\ipud}^{-} \vert. 
\]

Let us replace  the
expression (\ref{lambdauipud}) of $u(x_{i+1/2})$ in $\Hastun_{i+1/2}$ defined by
(\ref{Gast}) (note  that the computation is similar to that performed in
(\ref{uipud})-(\ref{vfsig})); this  yields 
\be 
\Hastun_\ipud = - \tau_\ipud   (u(x_\ipu) - u(x_i))  - \frac  {\l_i} \hip 
S_\ipud, \, i=1,\ldots,N-1.
\label{computGast}
\ee
Using (\ref{Gbarast}), this implies that
$ \Hast_{\ipud} = \Hbar_{\ipud} +  T_\ipud $
where

\[
\vert T_\ipud  \vert \le \vert R_{\ipud}^{-} \vert  + \vert R_{\ipud}^{+}  -  R_{\ipud}^{-} 
\vert  \frac \lbar
{2\lund}.
\]

Using (\ref{eqRun}) and \eqref{eqRde}, this last inequality yields that there exists 
$C_1$,  only depending  on
$\lbar,\lund,\g,\d$, such that

\[
\vert \Hast_{\ipud} - \Hbar_\ipud \vert=  \vert T_\ipud \vert \le 
C_1 h, \, i=1,\ldots,N-1.
\]

Therefore (\ref{lambdaconsist}) is proved. 

\smallskip

Define now the total exact fluxes;

\[
\Fbar_{i+\frac 1 2}= - (\l u_{x})(x_\ipud)  +
a u(x_\ipud), \; \forall i \in \{0, \ldots, N\},
\]

and define

\[
 \Fast_{i+\frac 1 2}= -\tau_\ipud (\uexip - \uexi)  + a u(x_i), \; \forall i \in \{1,
\ldots, N-1\},
\]

\[
\Fast_{\frac 1 2} = - \frac {\l_{1}} { h_1^-} (u(x_1)-c) + a
c, \; 
\Fast_{N + \frac 1 2} = -  \frac {\l_{_N}} { h_N^+} (d-u(x_N)) + a u_N.
\]

Then, from (\ref{lambdaconsist})  and the regularity of
$u$, there exists $C_2$,  only depending  on
$\lund,\lbar,\g$ and $\d$, such that 

\be
\Fast_{i+\frac 1 2} = \Fbar_{i+\frac 1 2} + R_\ipud, \text{ with }
\vert R_\ipud \vert \le  C_2  h, \, i=0,\ldots,N. \label{consistflux1D}
\ee

Hence the numerical approximation of the flux is consistent\index{Consistency!-- of the fluxes}.

\medskip
\textbf{Step 3. Error estimate.}

Integrating Equation (\ref{eq1}) over each control volume yields that 

\be 
\Fbar_{i+\frac 1 2} - \Fbar_{i-\frac 1 2} + b h_i (\uexi + S_i) = h_i f_i, \; \forall i
\in \{1, \ldots, N\}, \label{eqdc1} \ee

where $S_i \in \R$ is such that there exists $C_3$  only depending  on $u$ such that 
$\vert S_i \vert \le C_3 h$, for $i = 1, \ldots,N$.  

Using (\ref{consistflux1D}) yields that

\be 
 \Fast_{i+\frac 1 2} - \Fast_{i-\frac 1 2} + b h_i (\uexi + S_i) = h_i f_i
+ R_\ipud - R_\imud, \; \forall i \in \{1, \ldots, N\}.
\label{eqfast1D} 
\ee

Let $e_i = u(x_i) - u_i$, for $i = 1, \ldots,N$, and $e_0 = e_{N+1} = 0.$
Substracting (\ref{eqd1}) from (\ref{eqfast1D}) yields

\[
- \tau_{i+\frac 1 2} (e_\ipu  - e_i) + \tau_{i-\frac 1 2} (e_i - e_\imu)
+ a (e_i - e_\imu) + b h_i e_i = - b h_i  S_i +
R_\ipud - R_\imud, \; \forall i \in \{1, \ldots, N\}.
\]

Let us multiply  this equation by $e_i$, sum  for $i = 1,\ldots,N$,
reorder  the summations.  Remark that

\[
\sum_{i=1}^{N} e_i(e_i - e_{i-1}) = \frac 1 2
\sum_{i=1}^{N+1}(e_i-e_{i-1})^2
\]

and therefore  

\[
\sum_{i=0}^{N} \tau_{i+\frac 1 2} (e_\ipu  - e_i)^2 
+  \frac a 2 \sum_{i=1}^{N+1} (e_i-e_{i-1})^2  + \sum_{i=1}^{N} b h_i 
e_i^2 =  -\sum_{i=1}^{N} b h_i 
S_i e_i - \sum_{i=0}^{N} R_\ipud (e_\ipu - e_i).
\]

Since  $\vert S_i \vert \le C_3 h$ and thanks to (\ref{consistflux1D}), one has

\[
\sum_{i=0}^{N} \tau_{i+\frac 1 2} (e_\ipu  - e_i)^2 
\le  \sum_{i=1}^{N} b C_3 h_i 
h \vert e_i \vert + \sum_{i=1}^{N} C_2 h  \vert e_\ipu  - e_i \vert.
\]

Remark that $\vert e_i \vert \le \sum_{j = 1}^N  \vert e_j - e_{j-1} 
\vert$. Denote by
$A = \Bigl(\sum_{i=0}^{N} \tau_{i+\frac 1 2} (e_\ipu  - e_i)^2 \Bigr)^{\frac 1 2}$ and 
$B =  \Bigl(\sum_{i=0}^{N} \frac 1 {\tau_\ipud} \Bigr)^{\frac 1 2}$. 
The Cauchy-Schwarz  inequality yields

\[
A^2 
\le  \sum_{i=1}^{N} b C_3 h_i 
h A B  +  C_2 h  A B.
\]

Now, since 

\[
\frac 1 {\tau_\ipud}  \le \frac {\lbar}
{\lund^2} (\hipum + \hip), \, \sum_{i=0}^{N}  (\hipum + \hip) = 1,
\text{~with } h_0^+=h_{N+1}^-=0, \text{~and } \sum_{i=1}^{N} h_i =
1,
\]

one obtains that $ A \le C_4h$, with $C_4$  only depending  on $\lund, \lbar, 
\gamma$ and $\d$,  which yields  Estimate
(\ref{ee1}). 

Applying once again   the Cauchy-Schwarz  inequality yields  Estimate
(\ref{ee2}).
\eof

\subsection{The case of a point source term}
In many physical problems, some discontinuous or point source terms appear.  In the
case where a source term exists at the interface $x_{i + 1/2}$, the fluxes
relative to $K_i$ and  $K_{i+1}$ will differ because of this source term. The
computation of the fluxes is carried out in a similar way, writing that the sum of the
approximations of the fluxes must be equal to the source term at the interface.
Consider again the one-dimensional conservation problem \eqref{eq1}, 
\eqref{cl} (with, for the sake of simplification, $a=b=c=d=0$, we use below 
the notations of the previous section), but assume
now that at $\underline x \in (0,1)$, a point source of intensity $\a$ exists. In this case, the
problem may be written in the following way: 
\be
- (\l u_{x}(x))_x  =  f (x), \ \ x \in (0,\underline x) \cup (\underline x,1), \label{eqdisc}
\ee
\be
u(0) =  0, \label{cldisc1}
\ee
\be
u(1)  =  0, \label{cldisc2}
\ee
\be
(\l u_{x})^+ (\underline x) - (\l u_{x})^ - (\underline x) = - \a, \label{fluxdisc}
\ee

where   \[(\l u_{x})^+  (\underline x) = \lim_{x \tends \underline x, x >\underline x} 
(\l u_{x}) (x) \text{ and  }
(\l u_{x})^- (\underline x) =  \lim_{x \tends \underline x, x < \underline x} (\l u_{x}) (x).\]
 Equation  (\ref{fluxdisc}) states that the flux is discontinuous
at point $\underline x$. Another formulation of the problem is the following:

\be
- (\l u_{x})_x  =  g   \text{ in } \D'((0,1)), \label{eqmatdisc}
\ee
\be
u(0)  =  0, \label{cldisc21}
\ee
\be
u(1)  =  0, \label{cldisc22}
\ee

where $g = f + \a \d_{\underline x}$, where $ \d_{\underline x}$ denotes the Dirac measure, 
which is defined by
$< \d_{\underline x},\varphi>_{\D',\D} = \varphi(\underline x)$, for any $\varphi \in \D((0,1)) =
C^\infty_c((0,1),\R)$, and $\D'((0,1))$ denotes  the set of
distributions on (0,1), i.e. the set of continuous linear forms on
$\D((0,1))$. 

Assuming the mesh \index{Finite volume scheme!-- for 
elliptic problems!-- in one space dimension} 
 to be such that $\underline x = x_{i+ 1/2}$ for some $i \in
{1,\ldots,N-1}$, the equation corresponding to the unknown $u_i$ is
$F_{i+ 1/2}^- -F_{i- 1/2} = \int_{K_i} f(x) dx$, while the
equation   corresponding to the unknown $u_{i+1}$ is $F_{i+ 3/2}  -F_{i + 
1/2}^+ = \int_{K_{i+1}} f(x) dx$. In order to compute
the values of the numerical fluxes $F_{i+ 1/2}^\pm $,  one must take the source term into 
account while writing the
conservativity  of the flux; hence at $x_{i+ 1/2}$, the two numerical fluxes at 
$x = \underline x$, namely $F^+_{i+ \frac
1 2}$ and $F^-_{i+  \frac 1 2}$, must satisfy, following Equation 
(\ref{fluxdisc}),

\be 
F^+_{i+ \frac 1 2} - F_{i+ \frac 1 2}^- = \a.
\label{ddf}
\ee

Next, the fluxes  $F^+_{i+1/2}$ and $F^-_{i+1/2}$ must be expressed
in terms of the discrete variables $u_k$, $k = 1, \ldots,N$; in order to do so, introduce 
the auxiliary  variable
$u_{i+ 1/2}$ (which will be eliminated later), and write

\[
F^+_{i+ \frac 1 2}  =  - \l_{i+1} \frac {u_{i+1} - u_{i+\frac 1 2}} {
h_{i+1}^- }
\]

\[
F^-_{i+ \frac 1 2}  =  - \l_{i} \frac { u_{i+\frac 1
2} - u_{i}} {h_{i}^+ }.
\]

Replacing these expressions in (\ref{ddf}) yields

\[
u_{i+ \frac 1 2}  =  \frac {h_{i}^+ h_{i+1}^-}{(h_{i+1}^- \l_{i} + h_{i}^+ 
\l_{i+1})} 
[\frac {\l_{i+1}}  {h_{i+1}^-} u_{i+1} +  \frac {\l_{i}}  {h_{i}^+} u_{i} + \a].
\]

and therefore 

\[
F^+_{i+ \frac 1 2}  =  \frac { h_i ^+ \l_{i+1}} { h_{i+1}^-  \l_{i} +  h_i^+
\l_{i+1}} \a  - \frac {  \l_{i} \l_{i+1}} { h_{i+1}^- \l_{i} +  h_i^+
\l_{i+1}} (u_{i+1} - u_{i})
\]

\[
F^-_{i+ \frac 1 2}  =  \frac { -h_{i+1}^-\l_{i}} { h_{i+1}^- \l_{i} +  h_i ^+\l_{i+1}} \a  -  \frac {  \l_{i}
\l_{i+1}} { h_{i+1}^- \l_{i} +  h_i ^+\l_{i+1}} (u_{i+1} - u_{i}). 
\]

Note that the source term $\a$ is distributed on either
side of the interface  proportionally to the coefficient $\l$, and
that, when $\a = 0$, the above expressions lead to 

\[
F^+_{i+ \frac 1 2}
= F^-_{i+ \frac 1 2}  = - \frac {\l_{i} \l_{i+1}} { h_{i+1}^- \l_{i} + 
h_i^+ \l_{i+1}} (u_{i+1} - u_{i}).
\] 

Note that the error estimate given in Theorem \ref{ee1dell} still holds in this case
(under adequate assumptions).

\section{A semilinear elliptic problem}

\hpage{Elliptic 1d}

\subsection{Problem and Scheme}

This section is concerned with the proof of  convergence for
some nonlinear problems.  We are interested, as an example, by  the
following problem:

\be
-u_{xx}(x)  =  f(x,u(x)), \, x \in (0,1), 
\label{enl1d}
\ee

\be 
u(0)=u(1)=0,
\label{enl1dcl}
\ee

with a function $f :\, (0,1) \times \R \tends \R$ such that
\be
\ban f(x,s) \text{ is measurable with respect to } x \in (0,1) \text{ for all } s \in \R
\\
\text{ and continuous with respect to } s \in \R  \text{ for a.e. } x \in (0,1),
\ean 
\label{carat}
\ee
\be
f \in L^{\infty}((0,1) \times \R).
\label{bornee}
\ee

It is possible to prove that there exists at least one weak solution to \eqref{enl1d},
\eqref{enl1dcl},  that is a function $u$ such that

\be
u \in H^1_0((0,1)), \, \int_0^1 u_{x}(x)v_x(x)dx = \int_0^1f(x,u(x))v(x)dx, \, \forall v \in 
H^1_0((0,1)).
\label{enl1dw}
\ee

Note that \eqref{enl1dw} is equivalent to ``$u \in H^1_0((0,1))$ and $-u_{xx}=f(\cdot,u)$ in the
distribution sense in $(0,1)$".

The proof of the existence of such a solution is possible by using, for instance, the 
Schauder's fixed point theorem (see e.g. \citet{deimling}) or 
by using the convergence theorem  \ref{enl1dcr}  which is proved in the sequel.

\medskip

Let $\mesh$ be an admissible mesh of $[0,1]$ in the sense of Definition \ref{meshund}.
In order to discretize \eqref{enl1d}, \eqref{enl1dcl}, let us consider the following  (finite volume) scheme\index{Finite volume scheme!-- for elliptic problems!-- in  the one  dimensional case} 
 
 \be
 F_{i+\frac 1 2} - F_{i-\frac 1 2} =h_if_i(u_i), \, i=1,\ldots,N,
 \label{enl1dd}
 \ee 
 
 \be
 F_{i+\frac 1 2}= -\frac {u_{i+1}-u_i} {h_{i+\frac 1 2}}, \,
i=0,\ldots,N,
 \label{enl1df}
 \ee
 
 \be
 u_0=u_{N+1}=0,
 \label{enl1dcld}
 \ee
  with $f_i(u_i)=\frac 1 {h_i} \int_{K_i} f(x,u_i)dx$, $i=1,\ldots,N$.

 The discrete unknowns are therefore $u_1,\ldots,u_N$.
 In order to give a convergence result for this scheme (Theorem \ref{enl1dcr}), one first proves the existence of a solution to 
\eqref{enl1dd}-\eqref{enl1dcld}, a stability result, that is, an  estimate on the solution of \eqref{enl1dd}-\eqref{enl1dcld} (Lemma \ref{enl1dl1}) and a compactness lemma (Lemma \ref{enl1dl2}).

\begin{lemma}[Existence and stability result]\label{enl1dl1}\index{Stability}
\index{Estimate on the approximate solution!-- for an elliptic equation}
 Let $f:\, (0,1) \times \R \tends \R$  satisfying \eqref{carat}, \eqref{bornee} and $\mesh$ be an admissible mesh of 
$(0,1)$ in the sense of Definition \ref{meshund}. 
Then, there exists $(u_1,\ldots,u_N)^t \in \R^N$ solution of \eqref{enl1dd}-\eqref{enl1dcld} and which satisfies: 
\be
 \sum_{i=0}^N \frac {(u_{i+1}-u_i)^2} {h_{i+\frac 1 2}} \le C,
 \label{h1dest}
 \ee
for some $C \ge 0$  only depending  on $f$.

\end{lemma}
 
\proof of Lemma \ref{enl1dl1}
 
\smallskip

Define $M=\Vert f \Vert_{L^{\infty}((0,1) \times \R)}$. 
The proof of estimate (\ref{h1dest}) is given in a first step, and the existence of a solution to \eqref{enl1dd}-\eqref{enl1dcld} in a second step.
 
~

\textit{Step 1} (Estimate) 
 
 Let $V=(v_1,\ldots,v_N)^t \in \R^N$, there exists a unique  $U=(u_1,\ldots,u_N)^t \in \R^N$  solution of \eqref{enl1dd}-\eqref{enl1dcld} with $f_i(v_i)$ instead of $f_i(u_i)$ in the right hand-side (see  Theorem \refp{eelapl1D}). 
 One sets $U=F(V)$, so that $F$ is a continuous mapping from $\R^N$ to $\R^N$, and ($u_1,\ldots,u_N$) is a solution   to \eqref{enl1dd}-\eqref{enl1dcld} if and only if $U=(u_1,\ldots,u_N)^t$  is a fixed point to $F$.
Multiplying \eqref{enl1dd} by $u_i$ and summing over $i$ yields
  
\be
 \sum_{i=0}^N \frac {(u_{i+1}-u_i)^2} {h_{i+\frac 1 2}} \le M \sum_{i=1}^N 
 h_i \vert u_i \vert,
\label{preest}
\ee

and from   the Cauchy-Schwarz  inequality, one has
 
\[
 \vert u_i \vert \le \bigl(\sum_{j=0}^N \frac {(u_{j+1}-u_j)^2} {h_{j+\frac 1
2}}\bigr)^{\frac 1  2},
  \, i=1,\ldots,N,
\]

then  \eqref{preest} yields,  with $C=M^2$,

\be 
\sum_{i=0}^N \frac {(u_{i+1}-u_i)^2} {h_{i+\frac 1 2}} \le C.
\label{h1destbis}
\ee

This gives, in particular, Estimate \eqref{h1dest} if $(u_1,\ldots,u_N)^t \in \R^N$ is a solution of \eqref{enl1dd}-\eqref{enl1dcld} (that is $u_i=v_i$ for all $i$).
 
~

\textit{Step 2} (Existence)

The mapping $F:\, \R^N \tends \R^N$ defined above is continuous and, taking in $\R^N$ the norm 
\[
\Vert V \Vert= \big(\sum_{i=0}^N \frac
{(v_{i+1}-v_i)^2} {h_{i+\frac 1 2}}\big)^{\frac 1 2}, \text{~for } 
V=(v_1,\ldots,v_N)^t,\text{~with } v_0=v_{N+1}=0,
\]
one has $F(B_M) \subset B_M$, where $B_M$ is the closed ball of radius $M$ and center $0$ in $\R^N$. 
Then, $F$ has a fixed point in $B_M$ thanks to the Brouwer fixed point theorem (see e.g. \citet{deimling}).  
This fixed point is a solution to \eqref{enl1dd}-\eqref{enl1dcld}.
\eof

\subsection{Compactness results}

\index{Compactness!Discrete -- results|(}.

\begin{lemma}[Compactness]\label{enl1dl2}~
 
For an admissible mesh $\mesh$ of $(0,1)$ (see definition  
\ref{meshund}), let $(u_1,\ldots,u_N)^t\in \R^N$
satisfy \eqref{h1dest} for some $C \in \R$ (independent of $\mesh$) and 
let $u_{\mesh}:
\, (0,1) \tends \R$ be defined by $u_{\mesh}(x)=u_i$ if $x \in K_i$, 
$i=1,\ldots,N$.

Then, the set  $\{u_{\mesh}$, $\mesh$ admissible mesh of $(0,1)\}$ is
relatively  compact in $L^2((0,1))$. Furthermore, if $u_{\mesh_n}\tends u$ in
$L^2((0,1))$ and $\size(\mesh_n) \tends 0$, as $n \tends \infty$, then, $u \in
H^1_0((0,1))$.
\end{lemma}
 
\proof of Lemma \ref{enl1dl2} 
 
\smallskip

A possible proof is to use ``classical" compactness results, 
replacing $u_\mesh$ by a continuous function, say $\ubar_{\mesh}$,
piecewise affine, such that $\ubar_{\mesh}(x_i)=u_i$ for
$i=1,\ldots,N$, and $\ubar_{\mesh}(0)=\ubar_{\mesh}(1)=0$. The set
$\{\ubar_{\mesh}$, $\mesh$ admissible mesh of $(0,1)\}$ is
then bounded in  $H^1_0((0,1))$, see Remark \refp{Liebig}.

Another proof is given here, the
interest of which is its simple generalization to multidimensional
cases (such as the case of one unknown per triangle in $2$ space
dimensions, see Section \refp{vfquatre} and Section \refp{rellichd}) when the
construction of such a function, $\ubar_{\mesh}$, ``close" to
$u_{\mesh}$ and  bounded in $H^1_0((0,1))$ (independently of $\mesh$), is 
not so easy.

~

In order to have $u_{\mesh}$ defined on $\R$, one sets $u_{\mesh}(x)=0$ for $x \notin
[0,1]$. The proof may be decomposed  into four steps.

~

\textit{Step 1}. First remark that the set $\{u_{\mesh}$, $\mesh$  an admissible mesh of 
$(0,1)\}$ is bounded in $L^2(\R)$. Indeed, this an easy consequence of
\eqref{h1dest}, since one has, for all $x \in [0,1]$ (since $u_0=0$ and by
 the Cauchy-Schwarz  inequality), 

\[
\vert u_\mesh(x) \vert \le \sum_{i=0}^N \vert u_{i+1}-u_i \vert \le 
(\sum_{i=0}^N \frac {(u_{i+1}-u_i)^2} {h_{i+\frac 1 2}})^{\frac 1 2} \le C.
\]

~

\textit{Step 2}. Let $0 < \eta < 1$. One proves, in this step, that

\be
\Vert u_{\mesh}(\cdot+\eta)-u_{\mesh} \Vert_{L^2(\R)}^2 \le C\eta(\eta+2h).
\label{cleest}
\ee

(Recall that $h=\size(\mesh)$.)

Indeed, for $i \in \{0,\ldots,N\}$ define $\chi_{i+1/2}:$ $\R \tends \R$, by
$\chi_{i+1/2}(x)=1$, if $x_{i+1/2} \in [x,x+\eta]$ and 
$\chi_{i+1/2}(x)=0$, if $x_{i+1/2} \notin [x,x+\eta]$.
Then, one has, for all $x \in \R$,

\[
 (u_{\mesh}(x+\eta)-u_{\mesh}(x))^2 \le  \big(\sum_{i=0}^N 
\vert u_{i+1}-u_i \vert
\chi_{i+\frac 1 2}(x) \big)^2
\]
\be
\le \Bigl(\sum_{i=0}^N \frac {(u_{i+1}-u_i)^2} {h_{i+\frac 1
2}}\chi_{i+\frac 1 2}(x)\Bigr)\Bigl(\sum_{i=0}^N \chi_{i+\frac 1 2}(x)
h_{i+\frac 1 2}\Bigr).
\label{clest}
\ee

Since $\sum_{i=0}^N \chi_{i+1/2}(x) h_{i+1/2} \le \eta + 2h$, for all $x 
\in \R$, and $\int_\R\chi_{i+1/2}(x)dx=\eta$, for all $i \in 
\{0,\ldots,N\}$, integrating \eqref{clest} over $\R$ yields \eqref{cleest}.

~

\textit{Step 3}. For $0 < \eta < 1$,  Estimate \eqref{cleest} implies that 

\[
\Vert u_{\mesh}(\cdot+\eta)-u_{\mesh} \Vert_{L^2(\R)}^2 \le 3C\eta.
\]

This gives (with Step 1), by the Kolmogorov compactness theorem (recalled in Section \ref{rellichd}, see Theorem \refp{Kolm}), the relative compactness of the set $\{u_{\mesh}$,
$\mesh$  an admissible mesh of $(0,1)\}$ in $L^2((0,1))$ and also in $L^2(\R)$ (since $u_{\mesh}=0$ on $\R \setminus [0,1]$).

~

\textit{Step 4}. In order to conclude the proof of Lemma \ref{enl1dl2}, one may use Theorem \refp{kolmh10}, which we prove here in the one-dimensional case for the sake of clarity. 
Let  $(\mesh_n)_{n  \in \N}$  be a sequence of admissible meshes of $(0,1)$ such that $\size(\mesh_n) \tends 0$ and $u_{\mesh_n} \tends u$, in $L^2((0,1))$, as $n \tends \infty$. 
Note that $u_{\mesh_n} \tends u$, in $L^2(\R)$, with $u=0$ on $\R \setminus [0,1]$. 
For a given $\eta \in (0,1)$, let $n \to +\infty$ in \eqref{cleest}, with $u_{\mesh_n}$ instead of $u_{\mesh}$ (and
$\size(\mesh_n)$ instead of $h$); one  obtains 

\be
\Vert \frac {u(\cdot+\eta)-u} {\eta} \Vert_{L^2(\R)}^2 \le C.
\label{preh1}
\ee

Since $(u(\cdot+\eta)-u)/\eta$ tends to $Du$ (the distribution derivative of $u$) in the distribution sense, as $\eta \tends 0$,
Estimate \eqref{preh1} yields
that $Du \in L^2(\R)$. Furthermore, since $u=0$ on $\R \setminus
[0,1]$, the restriction of $u$ to $(0,1)$ belongs to $H^1_0((0,1))$. The proof of Lemma
\ref{enl1dl2} is complete. \eof

\index{Compactness!Discrete -- results|)}.

\subsection{Convergence}

\index{Convergence of the finite volume method!-- for a semilinear elliptic 
equation|(}
The following convergence
result follows from  lemmata  \ref{enl1dl1} and \ref{enl1dl2}.

\begin{theo}\label{enl1dcr}
 Let $f :\, (0,1) \times \R \tends \R$  satisfying \eqref{carat}, 
 \eqref{bornee}. For an admissible mesh, $\mesh$, of $(0,1)$ (see 
 Definition \ref{meshund}), let 
 $(u_1,\ldots,u_N)^t \in \R^N$ be a solution to \eqref{enl1dd}-\eqref{enl1dcld} (the existence of which is given by Lemma
\ref{enl1dl1}), and let
$u_{\mesh}:
\, (0,1) \tends \R$ by
$u_{\mesh}(x)=u_i$, if $x \in  K_i$, $i=1,\ldots,N$.

Then, for any sequence $(\mesh_n)_{n \in \N}$  of admissible meshes  such that
$\size(\mesh_n) \tends 0$, as $n \tends \infty$, there exists a subsequence, still denoted
by $(\mesh_n)_{n \in \N}$,  such that $u_{\mesh_n} \tends u$, in $L^2((0,1))$, as $n 
\tends
\infty$, where $u \in H^1_0((0,1))$ is a weak solution to  \eqref{enl1d}, \eqref{enl1dcl} 
(that
is, a solution to \eqref{enl1dw}).
\end{theo} 

\proof of Theorem \ref{enl1dcr} 

\smallskip

Let $(\mesh_n)_{n \in 
\N}$ be
a sequence of admissible meshes of $(0,1)$ such that
$\size(\mesh_n) \tends 0$, as $n \tends \infty$. By lemmata \ref{enl1dl1} 
and \ref{enl1dl2}, there
exists a subsequence, still denoted by $(\mesh_n)_{n \in \N}$,  such that $u_{\mesh_n} 
\tends u$,
in $L^2((0,1))$, as $n \tends \infty$, where  $u \in H^1_0((0,1))$. In order to 
conclude, it only remains to prove that $-u_{xx}=f(\cdot,u)$ in the
distribution sense in $(0,1)$. 

To prove this, let $\varphi \in C_c^\infty((0,1))$. Let
$\mesh$ be an admissible mesh of $(0,1)$,  and $\varphi_i=\varphi(x_i)$, 
$i=1,\ldots,N$, and
$\varphi_0=\varphi_{N+1}=0$. If $(u_1,\ldots,u_N)$ is a solution to
\eqref{enl1dd}-\eqref{enl1dcld}, multiplying \eqref{enl1dd} by $\varphi_i$ and summing 
over
$i=1,\ldots,N$ yields

\be
\int_0^1 u_{\mesh}(x)\psi_{\mesh}(x)dx=\int_0^1f_{\mesh}(x)\varphi_{\mesh}(x)dx,
\label{yenamar}
\ee

where 

\[
\psi_{\mesh}(x)=\frac 1 {h_i} (\frac {\varphi_{i}-\varphi_{i-
1}} {h_{i-\frac 1 2}}-
\frac {\varphi_{i+1}-\varphi_{i}} {h_{i+\frac 1 2}}), \, 
f_{\mesh}(x)=f(x,u_i) \text{~and } \varphi_{\mesh}(x)=\varphi_i, \text{ 
if } x \in K_i.
\]

Note that, thanks to the regularity of the function $\varphi$, 

\[
\frac {\varphi_{i+1}-\varphi_{i}} {h_{i+\frac 1 
2}}=\varphi_x(x_{\ipud})+R_{\ipud}, 
\, \vert R_{\ipud} \vert \le C_1 h,
\]

with some $C_1$  only depending   on $\varphi$, and therefore

\[
\ban
\dsp{\int_0^1 u_{\mesh}(x)\psi_{\mesh}(x)dx }& = \dsp{\sum_{i = 1}^N \int_{K_i} 
\frac {u_i} {h_i} \Bigl(\varphi_x (x_{\imud}) - \varphi_x(x_{\ipud}) 
\Bigr) dx +  \sum_{i = 1}^N  u_i (R_{\imud} - R_{\ipud}) }\\ 
& =  \dsp{ \int_0^1 - u_\mesh (x)  \theta_\mesh (x) dx +  \sum_{i = 0}^{N} R_\ipud (u_{i+1} - u_i),}
\ean
\]
with $u_0=u_{N+1}=0$, where the piecewise constant function 
\[ \theta_\mesh =\sum_{i=1,N} \frac 
{\varphi_x (x_{\ipud}) - \varphi_x(x_{\imud})}  {h_i}  1_{\cv_i} \]
tends to $\varphi_{xx}$ as $h$ tends to 0.

Let us consider  \eqref{yenamar} with ${\mesh_n}$ instead of  ${\mesh}$; 
thanks to the Cauchy-Schwarz  inequality, a passage to the limit as $n \tends \infty$ 
gives, thanks to \eqref{h1dest}, 

\[
-\int_0^1u(x)\varphi_{xx}(x)dx=\int_0^1f(x,u(x))\varphi(x)dx,
\]

and therefore $-u_{xx}=f(\cdot,u)$ in the
distribution sense in $(0,1)$. This concludes the proof of Theorem 
\ref{enl1dcr}. Note that the crucial idea of this proof is to use the 
property of consistency of the fluxes on the regular test function 
$\varphi$.
\eof

\begin{remark} \normalfont\label{sublinear}
It is possible to give some extensions of the results of this section. For instance,  
Theorem \ref{enl1dcr} is true with an assumption of ``sublinearity" on $f$ instead of
\eqref{bornee}. Furthermore, in order to have both existence and uniqueness of the
solution to \eqref{enl1dw} and a rate of convergence (of order $h$) in Theorem
\ref{enl1dcr}, it is sufficient to assume, instead of \eqref{carat} and \eqref{bornee}, 
that $f \in C^1([0,1]\times \R,\R)$ and that there exists $\gamma
<1$, such that $(f(x,s)-f(x,t))(s-t) \le \gamma (s-t)^2$, for all $(x,s) \in
[0,1]\times \R$.
\end{remark}

\index{Convergence of the finite volume method!-- for a semilinear elliptic 
equation|)}
\chapter{Elliptic problems in two or three dimensions} \label{ellmd}

The topic of this chapter is the discretization of elliptic problems in several space 
dimensions  by the 
finite volume method. The one-dimensional case which was studied in 
Chapter \ref{ellund} is easily generalized to nonuniform rectangular or
parallelipedic meshes. However, for general shapes of  control volumes, 
the definition of the scheme (and the proof of convergence) requires some 
assumptions which define an ``admissible mesh".   Dirichlet and Neumann 
boundary conditions are both considered. In both cases,  a discrete 
Poincar\'e inequality is used, and  the stability
of the scheme is proved 
by establishing estimates on the approximate solutions. 
The convergence of the scheme without any assumption on the 
regularity of the exact solution is proved; this result may be 
generalized, under adequate assumptions, to nonlinear equations. 
Then, again in both the Dirichlet and Neumann cases, an error estimate 
between the finite volume approximate solution and the  $C^2$ or $H^2$
regular exact solution to the continuous 
problems are proved. 
The results are generalized to the case of  matrix 
diffusion coefficients and more general boundary conditions. 
Section \ref{scvfe} is devoted to finite volume schemes written with 
unknowns located at the vertices. Some links between the finite element 
method, the ``classical" finite volume method and the ``control volume 
finite element" method  introduced by \citet{F1} are given. 
Section \ref{meshref} is devoted to the treatment of  singular sources 
and to mesh refinement; under suitable assumption, it can be shown that 
error estimates still hold for ``atypical" refined meshes. 
Finally, Section 
\ref{rellichd} is  devoted to the proof of compactness results which are
used in the proofs of convergence of the schemes.

\section{Dirichlet boundary conditions}\label{sdbc}

\hpage{Elliptic multid}
\index{Dirichlet boundary condition|(}

Let us consider here the following elliptic 
equation\index{Equation!Convection-diffusion --}

\be        
- \Delta u(x) + \div(\vbf u) (x) + b u(x) = f(x), \ \   x \in \O,        
\label{elldifstaf} \ee  
with   Dirichlet boundary condition:
\be        
u(x) = g(x), \ \   x \in \dr\O,      
\label{elldifstab}\ee 

where

\begin{hyp}\label{ellH}~
\benum
\item $\O$ is an open bounded polygonal subset of $\R^d$, $d=2$ or $3$, 
\item $b\ge 0$,
\item $f \in L^2(\O)$,
\item $\vbf \in C^1(\Obar, \R^d); \div \vbf \ge 0$,
\item $g \in C (\dr \O, \R)$ is such that there exists $\tilde g \in
H^1(\O)$ such that $\gammat(\tilde g) = g$ a.e. on $\dr \O$.
\eenum
\end{hyp}

Here, and in the sequel, ``polygonal" is used for both $d = 2$
and $d= 3$  (meaning polyhedral in the latter case) and $\gammat$ denotes 
the trace operator from $H^1(\O)$ into $L^2(\dr \O)$. Note also that 
``a.e. on $\dr \O$" is a.e. for the $d-1$-dimensional Lebesgue measure on 
$\dr \O$.

Under Assumption \ref{ellH}, by the
Lax-Milgram theorem,  there exists a unique variational solution $u\in H^1(\O)$
of Problem  \eqref{elldifstaf}-\eqref{elldifstab}. (For the study of elliptic 
problems and their discretization by finite 
element methods, see e.g. \citet{ciarlet} and references therein).
This solution satisfies $u = w + \tilde
g$, where $\tilde g \in H^1(\O)$ is such that $\gammat(\tilde g) = g$, 
a.e. on $\dr \O$, and $w$
is the unique function of $H^1_0(\O)$ satisfying

\be
\bac
\dsp{\int_\O \Bigl( \grad w(x) \cdot \grad \psi(x) + \div (\vbf w) (x)
\psi(x)  +b w(x) \psi(x) \Bigr) dx = }\\

\dsp{ \int_\O \Bigl(- \grad \tilde  g(x) \cdot \grad \psi(x)
- \div (\vbf  \tilde g) (x) \psi(x) -b \tilde g(x) \psi(x)  + f(x)
\psi(x)\Bigr) dx, \;     \forall \psi\in H^1_0(\O).}
\ean
\label{ellsolvar}
\ee

\subsection{Structured meshes}\label{smesh}

If $\O$ is a rectangle ($d=2$) or a parallelepiped  ($d=3$), 
it may then be meshed 
with  
rectangular or parallelepipedic control volumes\index{Control volume}. In this case, the one-dimensional
scheme may easily be generalized\index{Mesh!Structured --|(}.

\subsubsection{Rectangular meshes for the Laplace operator}

Let us for instance consider the case $d=2$, let 
$\O = (0,1) \times (0,1)$, and $f \in C^2(\O,\R)$ (the three dimensional case is
similar). Consider Problem \eqref{elldifstaf}-\eqref{elldifstab} and
assume here that $b = 0$, $\vbf = 0$ and $g = 0$ (the general case is 
considered later, on general unstructured meshes).  The problem reduces 
to the pure diffusion equation\index{Equation!Diffusion --}:

\be
\bal
-\Delta u(x,y)  =  f(x,y), \, (x,y)  \in \O, \\
u (x,y)  =  0, \, (x,y) \in \partial \O. 
\ean
\label{lapl2d}
\ee  

In this section, it is convenient to denote by $(x,y)$ the current point 
of $\R^2$ (elsewhere, the notation $x$ is used for a point or a vector of 
$\R^\dim$).
 
Let $\mesh = (K_{i,j})_{i=1,\cdots,N_1;j=1,\cdots,N_2}$ be an admissible mesh of $(0,1) \times 
(0,1)$, that is, satisfying the following  assumptions
(which generalize  Definition \ref{meshund})

\begin{hyp} \label{hyprect}\index{Mesh!Rectangular --}
Let $N_1 \in \N^\star, N_2 \in \N^\star, \, h_1,\ldots,h_{N_1}> 0,  \,
k_1,\ldots,k_{N_2}> 0$ such that 

\[ 
\dsp{ \sum_{i=1}^{N_1} h_i = 1,
\sum_{i=1}^{N_2} k_i = 1, }
\] 
and let $h_0 = 0, h_{N_1+1} =0, k_0 = 0,
k_{N_2+1} =0$. For $ i=1,\ldots,N_1,$  let $
x_{\frac 1 2}=0, \, x_{i+\frac 1 2}=x_{i-\frac 1 2} + h_i,$ (so
that $x_{N_1+ \frac 1 2}=1),$
and for $j=1,\ldots,N_2,$ 
$y_{\frac 1 2}=0, \, y_{j+\frac 1 2}=y_{j-\frac 1 2} + k_j,$ 
(so that $y_{N_2+\frac 1 2}=1)$
and 
\[
K_{i,j}= [x_{i-\frac 1 2}, x_{i+\frac 1 2}] 
\times [y_{j-\frac 1 2}, y_{j+\frac 1 2}]. 
\]
 Let $(x_i)_{i=0,N_1+1},$ and $(y_j)_{j=0,N_2+1},$ such that
\[
x_{i-\frac 1 2}<x_i<x_{i+\frac 1 2}, \, \text{ for } i=1,\ldots,N_1, \, x_0=0, \, 
x_{N_1+1}=1,
\]
\[
y_{j-\frac 1 2}<y_j<y_{j+\frac 1 2}, \, \text{ for } j=1,\ldots,N_2, \, y_0=0, \, 
y_{N_2+1}=1, 
\]
 and let $x_{i,j} = (x_i, y_j)$, for $i=1,\ldots,N_1,$, $j=1,\ldots,N_2;$ 
set 

\[
{h_i}^- = x_i - x_{i-\frac 1 2}, \,  {h_i}^+ =  
x_{i+\frac 1 2} - x_i ,  \, \text{ for } i=1,\ldots, N_1, \,  
h_{i+\frac 1 2}=x_{i+1}-x_i, \, \text{ for } i=0,\ldots, N_1,
\]

\[
{k_j}^- =  y_j - y_{j-\frac 1 2}, \, {k_j}^+ =  
y_{j+\frac 1 2} - y_j, \, \text{ for } j=1,\ldots, N_2, \,  k_{j+\frac 1 2}=  
y_{j+1}-y_j, \, \text{ for } j=0,\ldots, N_2.
\]

Let $h = \max \{(h_i,{i=1,\cdots,N_1}),(k_j, {j=1,\cdots,N_2})\}.$ 
\end{hyp}


As in the 1D case, the finite volume scheme is found by integrating the
first equation  of (\ref{lapl2d}) over each control volume $K_{i,j}$,  which
yields 

\[ 
\ba
 - \dsp{ \int_{y_{j-\frac 1 2}}^{y_{j+\frac 1 2}} u_x(x_{i +
\half},y) dy  +  \int_{y_{i-\frac 1 2}}^{y_{i+\frac 1 2}}
u_x(x_{i -\half},y) dy  } 
  \\
  +\dsp{  \int_{x_{i-\frac 1 2}}^{x_{i+\frac 1 2}}
u_y(x,y_{j-\half}) dx  - \int_{x_{i-\frac 1 2}}^{x_{i+\frac 1 2}}
u_y(x,y_{j+\half}) dx = \int_{K_{ij}} f(x,y) dx \, dy. }
\ea 
\]

The fluxes are then approximated by differential quotients with respect to the 
discrete unknowns 
$(u_{i,j}, i=1,\cdots,N_1, j=1,\cdots,N_2)$ in a similar manner to the 1D case; 
hence the numerical scheme \index{Finite volume scheme!-- for elliptic problems!-- in two or three space dimensions}  reads 
\be 
F_{i + \half,j} - F_{i -\half,j} + F_{i ,j+ \half} - F_{i ,j- \half} = 
h_{i,j} f_{i,j}, \,  \forall \ (i,j) \in \{1,\ldots, N_1\}\times\{1,\ldots, N_2\},
\label{vfrect} 
\ee
where $h_{i,j}=  h_i \times k_j$, $f_{i,j}$ is the mean value of $f$ over 
$K_{i,j}$, and

\be 
\bal
F_{i + \half,j}  =  - \dsp{ \frac  {k_j} {h_\ipud}
 (u_{i+1,j} - u_{i,j}), } \, \text{ for }  i = 0, \cdots, N_1, j =
1,\cdots,N_2, \\
 F_{i ,j+ \half}   =  - \dsp{ \frac  {h_i} {k_\jpud}
(u_{i,j+1} - u_{i,j}), }  \,
\text{ for }  i = 1,\cdots, N_1, j = 0,\cdots,N_2,
\ean
\label{fluxrect}\ee

\be
u_{0,j}=u_{N_1+1,j}=u_{i,0}=u_{i,N_2+1}=0, \, \text{ for }  i = 1,\ldots, 
N_1, j = 1,\ldots,N_2.
\label{clrect}
\ee

The numerical scheme \eqref{vfrect}-\eqref{clrect} is therefore clearly 
conservative  and the numerical approximations of the fluxes can easily be 
shown to be consistent.

\begin{prop}[Error estimate]\index{Error estimate!-- for an elliptic 
equation!-- in the general case}
Let $\O = (0,1) \times  (0,1) $ and $f \in L^2(\O)$. Let $u$ be the 
unique variational solution to \eqref{lapl2d}. Under Assumptions 
\ref{hyprect}, let $\zeta >0$ be such that $h_i \ge \zeta h$  for $i = 1,\ldots,
N_1$ and  $k_j \ge \zeta h$  for $ j= 1,\ldots,
N_2$. Then, there exists a unique solution $(u_{i,j})_{i = 1,\cdots, N_1,  j = 1,\cdots,N_2}$ to 
\eqref{vfrect}-\eqref{clrect}. Moreover, there exists $C >0$  only depending   on
$u$, $\O$ and $\zeta$ such that 

\be \sum_{i,j} \frac {( e_{i+1,j} -e_{i,j})^2} {h_{i+\frac 1 2}} k_j +
\sum_{i,j} \frac {( e_{i,j+1} -e_{i,j})^2} {k_{j+\frac 1 2}} h_i  \le Ch^2
\label{hunzrect} 
\ee

and 

\be
\sum_{i,j} (e_{i,j})^2 h_i  k_j \le C h^2,
\label{lderect}
\ee

where  $e_{i,j} = u(x_{i,j}) - u_{i,j}$, for $ i = 1,\cdots, N_1, j =
1,\cdots,N_2$.
\end{prop}

In the above proposition, since $f \in L^2(\O)$ and $\O$ is convex, 
it is well known that the variational solution $u$ to \eqref{lapl2d} 
belongs to $H^2(\O)$. We do not 
give here the proof of this proposition since it is in fact included in
Theorem \refp{testh2} (see also \citet{lazarovsiam} where the case $u \in 
H^s$, $s \ge \frac 3 2$ is also studied). 

\smallskip

In the case $u \in C^2(\Obar)$, the estimates (\ref{hunzrect}) and 
(\ref{lderect}) can be shown with the same technique as
in the 1D case (see e.g. \citet{thesejmf}). If  $u \in C^2$ then the above 
estimates are a consequence of 
Theorem \refp{ellesterr}; in this case,  the value  $C$ in \eqref{hunzrect}
and \eqref{lderect} 
independent of $\zeta$, and therefore the assumption  $h_i \ge \zeta h$  for $i = 1,\ldots,
N_1$ and  $k_j \ge \zeta h$  for $ j= 1,\ldots,
N_2$ is no longer needed.

\smallskip

Relation (\ref{hunzrect}) can be seen as an estimate of a ``discrete $H^1_0$ norm" of
the error,  while relation (\ref{lderect}) gives an estimate of the $L^2$ norm of the
error.

\begin{remark} \normalfont\label{h2conv2d}
Some slight modifications of the scheme \eqref{vfrect}-\eqref{clrect} are 
possible, as in the first item of Remark \refp{angot1}. It is also possible to obtain, 
sometimes, an ``$h^2$" estimate on the $L^2$ (or $L^{\infty}$) norm of 
the error (that is ``$h^4$" instead of ``$h^2$" in \eqref{lderect}), 
exactly as in the 1D case, see Remark \refp{h2conv}. In the case 
equivalent to the second case of Remark \ref{h2conv}, the point $x_{i,j}$ 
is not necessarily the center of $K_{i,j}$.
\end{remark}
\index{Mesh!Structured --|)}

When the mesh is no longer rectangular, the scheme (\ref{vfrect})-(\ref{fluxrect})
is not easy to generalize if keeping to a 5 points scheme. In particular,  the
consistency  of the fluxes or the 
conservativity can be lost, see \citet{Fa1}, which
yields a bad numerical behaviour of the scheme. One way to keep both properties is
to introduce  a 9-points scheme.

\subsubsection{Quadrangular meshes: a nine-point scheme}

Let $\Omega$ be an open bounded polygonal subset of $\R^2$, and $f$ be a regular function from
$\Obar$ to $\R$. We still consider Problem \ref{lapl2d}, turning back to 
the usual notation $x$ for the current point of $\R^2$, 

\be 
\bal
- \Delta u (x)= f(x), \, x \in  \Omega, \\ 
u(x) = 0, \,  x \in  \partial \Omega.
\ean 
\label{Lapl-2D} 
\ee

\index{Finite volume scheme!-- for 
elliptic problems!-- in two or three 
space dimensions|(} Let $\mesh$ be a mesh defined over $\O$;  then, integrating the first equation of
(\ref{Lapl-2D}) over any cell $\cv$ of the mesh yields
\[ 
- \int_{\partial \cv} \mathbf{grad} u \cdot \mathbf{n}_\cv  =  \int_{\cv} f,
\]
where  $\mathbf{n}_\cv$  is the normal to the boundary  $\partial \cv$, outward to $\cv$.
Let $u_\cv$ denote the discrete unknown associated to the control volume 
$\cv \in \mesh$.
In order to obtain a numerical scheme, if $\edge$ is a common edge to 
$\cv \in \mesh$ and $\cvv \in \mesh$ (denoted by $\edgecvcvv$) or if 
$\edge$ is an edge of $\cv \in 
\mesh$ belonging to $\dr \O$, the expression $\mathbf{grad} u \cdot \mathbf{n}_\cv$ must be
approximated on $\edge$ by using the discrete
unknowns.
The study of the finite volume scheme in dimension 1 and the above 
straightforward generalization to the rectangular case showed that the fundamental 
properties of the method seem to be
\begin{enumerate}

\item conservativity\index{Conservativity}: in the absence of any source term on  $\edgecvcvv$,
the approximation of $\mathbf{grad} u \cdot \mathbf{n}_\cv$ on  $\edgecvcvv$ which is used
in the  equation
associated with cell $\cv$ is equal to the approximation of $-\mathbf{grad} u \cdot \mathbf{n}_{\cvv} $  which is used in the equation associated with cell  $\cvv$. This
property is naturally obtained when using a finite volume scheme. 

\item consistency of the fluxes\index{Consistency!-- of the fluxes}:  taking for $u_\cv$ the value of $u$ 
 in a fixed point of $\cv$ (for instance, the center of gravity of 
 $\cv$), where $u$ is a regular function, the difference between $\mathbf{grad} u \cdot
\mathbf{n}_\cv $ and the chosen approximation of $\mathbf{grad} u \cdot \mathbf{n}_\cv$ is of
an order less or  equal
to that of the mesh size. This need of consistency will be discussed in 
more detail: see remarks \refp{cafe} and \refp{creme}

\end{enumerate}

Several computer codes use the following 
 ``natural" extension of \eqref{fluxrect} for the approximation  of  $\mathbf{grad} u \cdot
\mathbf{n}_\cv$ on $\cvbar \cap \cvvbar$:

\[
\mathbf{grad} u \cdot \mathbf{n}_\cv = \frac  {u_{\cvv}
- u_\cv} {\dcvcvv},
\]

where $\dcvcvv$ is the distance between the center of the
cells  $\cv$ and 
$\cvv$.
This choice, however simple, is far from optimal, at least in the case of a general (non
rectangular) mesh, because the fluxes thus obtained are not consistent; this yields 
important errors, especially in the case where the mesh cells are all oriented in the
same direction, see \citet{Fa1}, \citet{Fa2}. This problem may be avoided by modifying 
the
approximation of  $\mathbf{grad} u \cdot \mathbf{n}_\cv$ so as to make it 
consistent.  
However, one must be careful, in doing so, to maintain the conservativity 
 of the
scheme. To this purpose, a 9-points scheme was developped, which is denoted by
FV9. 

\begin{figure}[htbp]

\begin{center}
\begin{tikzpicture}[x=0.9mm, y=0.9mm, line cap=round, line join=round]
  \draw[line width=0.8pt] (60,40) -- (100,30);
  \draw[line width=0.8pt] (75,15) -- (87,63);
  \draw (100,30) -- (106.8,70.8);
  \draw (60,70) -- (107,70);
  \draw (60,40) -- (60,70);
  \draw (60,40) -- (50,10);
  \draw (100,30) -- (105,10);
  \draw (50,10) -- (105,10);
  \draw (105,10) -- (130,16.25);
  \draw (100,30) -- (130,35);
  \draw (130,35) -- (130,16);
  \draw (130,35) -- (125,65);
  \draw (125,65) -- (107,69.5);
  \draw (60,40) -- (30,35);
  \draw (30,35) -- (35,65);
  \draw (35,65) -- (60,70);
  \draw (50,10) -- (25,20);
  \draw (25,20) -- (30,35);
  \draw (82,54) -- (115,48.5);
  \draw (42,27.5) -- (79,21.333);
  \draw[-{Stealth[length=4pt]}] (70,33) -- (76,36);
  \node[inner sep=0pt] at (64,34) {${\sigma_{i,j+1/2}}$};
  \draw[-{Stealth[length=4pt]}] (81,63) -- (86,61.333);
  \node[inner sep=0pt] at (80,65) {${\eta_{i,j+1/2}}$};
  \draw[-{Stealth[length=4pt]}] (74,43) -- (84,53);
  \node[inner sep=0pt] at (72,41) {$U_{i,j+1/2}$};
  \draw[-{Stealth[length=4pt]}] (91,14) -- (77,21);
  \node[inner sep=0pt] at (97,12) {$D_{i,j+1/2}$};
  \node[inner sep=0pt] at (40,23) {$C_{i-1,j}$};
  \node[inner sep=0pt] at (84,21) {$C_{i,j}$};
  \node[inner sep=0pt] at (88,29) {$\mathbf{K}_{i,j}$};
  \draw[-{Stealth[length=4pt]}] (63,16) -- (66,22);
  \node[inner sep=0pt] at (65,14) {$\zeta_{i-1/2,j}$};
  \draw[-{Stealth[length=4pt]}] (90,50) -- (94,51.333);
  \node[inner sep=0pt] at (94,48) {$\zeta_{i+1/2,j+1}$};
  \node[inner sep=0pt] at (77,57) {$C_{i,j+1}$};
  \node[inner sep=0pt] at (95,36) {$\mathbf{K}_{i,j+1}$};
  \node[inner sep=0pt] at (114,52) {$C_{i+1,j+1}$};
\end{tikzpicture}
\end{center}
\caption{FV9 scheme}
\label{vf9fig}
\end{figure}
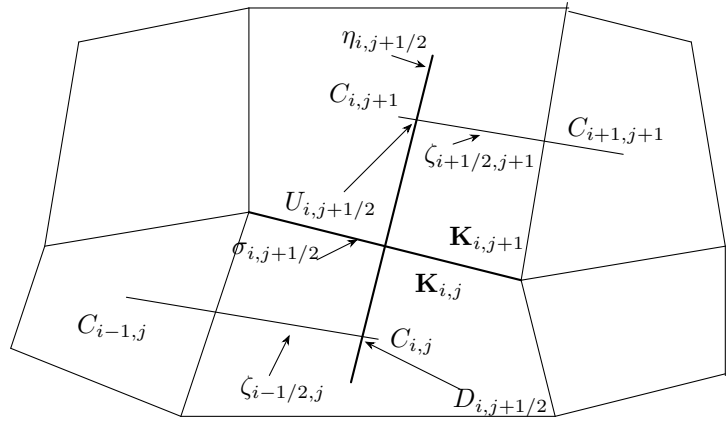

Let us describe now how the flux $\mathbf{grad} u \cdot \mathbf{n}_\cv$ is approximated by 
the
FV9 scheme. Assume here, for the sake of clarity, that the mesh $\mesh$ is structured;
indeed, it  consists in a set of quadrangular cells  $\{K_{i,j}, i= 1,\ldots,N; 
j=1,\ldots,
M\}$. As shown in Figure \ref{vf9fig}, let  $C_{i,j}$  denote the center of gravity of the cell $K_{i,j}$, $\edge_{i,j-1/2}$,
$\edge_{i+1/2,j}$, $\edge_{i,j+1/2}$, $\edge_{i-1/2,j}$ the four edges to $K_{i,j}$ 
and
$\eta_{i,j-1/2}$, $\eta_{i+1/2,j}$, $\eta_{i,j+1/2}$, $\eta_{i-1/2,j}$ their
respective orthogonal bisectors. Let $\zeta_{i,j-1/2}$, (resp. $\zeta_{i+1/2,j}$,
$\zeta_{i,j+1/2}$, $\zeta_{i-1/2,j}$) be the lines joining points $C_{i,j}$ and $C_{i,j-
1}$ (resp. $C_{i,j}$ and  $C_{i+1,j}$, $C_{i,j}$ and $C_{i,j+1}$, $C_{i-1,j}$ and
$C_{i,j}$). 

Consider for instance the edge $\edge_{i,j+1/2}$ which lies between the cells
$K_{i,j}$ and $K_{i,j+1}$ (see Figure \ref{vf9fig}). In order to find an approximation 
of
$\mathbf{grad}  u \cdot \mathbf{n}_\cv$, for $K=K_{i,j}$, at the center of this edge, we shall first derive an 
approximation of $u$ at the two points $U_{i,j+1/2}$ and $D_{i,j+1/2}$ which are 
located
on the  orthogonal bisector $\eta_{i,j+1/2}$ of the edge $\edge_{i,j+1/2}$, on each
side of the edge. Let  $\phi_{i,j+1/2}$  be the approximation of $- \mathbf{grad} u \cdot \mathbf{n}_\cv$ at the center of the edge $\edge_{i,j+1/2}$. A natural choice for
$\phi_{i,j+1/2}$ consists in taking 

\be
 \phi_{i,j+1/2} =  - \frac {u^U_{i,j+1/2} - u^D_{i,j+1/2}} {d(U_{i,j+1/2},
D_{i,j+1/2})}, \label{flux}
\ee

where $ u^U_{i,j+1/2} $ and $ u^D_{i,j+1/2}$ are approximations of $u$  at $
U_{i,j+1/2}$ and $D_{i,j+1/2}$, and $d(U_{i,j+1/2}$, $D_{i,j+1/2})$ is the
distance between points $U_{i,j+1/2}$ and $D_{i,j+1/2}$.

The  points $U_{i,j+1/2}$ and $D_{i,j+1/2}$ are chosen so that they are located on the 
lines
$\zeta$ which join the centers of the neighbouring cells. The
points $U_{i,j+1/2}$ and $D_{i,j+1/2}$  are therefore located at the intersection of the
orthogonal bisector $\eta_{i,j+1/2}$ with the adequate $\zeta$ lines, which are chosen
according to the geometry of the mesh. More precisely,  

\[
\left.\begin{array}{lll} U_{i,j+1/2} &
= \eta_{i,j+1/2} \cap \zeta_{i- 1/2,j+1}
\ \ & \mathrm{if} \ \ \eta_{i,j+1/2} \ \  \text{is} \text{ to the left of } \ \
C_{i,j+1} \\ & = \eta_{i,j+1/2} \cap \zeta_{i+1/2,j+1} \ \ & \mathrm{otherwise} 
\\
D_{i,j+1/2} & = \eta_{i,j+1/2} \cap \zeta_{i-1/2,j} \ \ & \mathrm{if} \ \
\eta_{i,j+1/2} \ \  \text{is} \text{ to the left of} \ \ C_{i,j} \\ & 
=
\eta_{i,j+1/2} \cap \zeta_{i+1/2,j} \ \ & \mathrm{otherwise} \end{array}\right. 
\]

In order to satisfy the property of consistency \index{Consistency!-- of the fluxes} of the fluxes, a second order
approximation of $u$ at points  $U_{i,j+1/2}$ and  $D_{i,j+1/2}$ is required. In the
case  of
the geometry which is described in Figure \ref{vf9fig}, the following linear approximations
of  $u^U_{i,j+1/2}$ and $u^D_{i,j+1/2}$ can be used in (\ref{flux}); 
 
\[
\left.\begin{array}{lll} u^U_{i,j+1/2} & = \alpha u_{i+1,j+1} + (1-\alpha)
u_{i,j+1} & \ \ \text{where} \ \ \alpha = \displaystyle{\frac
{d(C_{i,j+1},U_{i,j+1/2})}{d(C_{i,j+1},C_{i+1,j+1})}} \\  
 u^D_{i,j+1/2} & = \beta
u_{i-1,j} + (1- \beta) u_{i,j} & \ \  \text{where} \ \  \beta =
\displaystyle{ \frac
{d(C_{i,j},D_{i,j+1/2})} {d(C_{i-1,j},C_{i,j})}} \end{array}\right.
\]

The approximation of  $ \mathbf{grad} u \cdot \mathbf{n}_\cv $ at the center of a ``vertical" edge
$\edge_{i+1/2,j}$ is performed in a similar way, by introducing the points
$R_{i+1/2,j}$  intersection of the orthogonal bisector $\eta_{i+1/2,j}$ and,
according to the geometry, of the line $\zeta_{i,j-1/2}$ or $\zeta_{i,j+1/2}$, and 
$L_{i+1/2,j}$ intersection of $\eta_{i+1/2,j}$ and $\zeta_{i+1,j-1/2}$ or
$\zeta_{i+1,j+1/2}$. 

Note that the outmost grid cells  require a particular treatment 
(see \citet{Fa1})\index{Finite volume scheme!-- for 
elliptic problems!-- in two or three 
space dimensions|)}.

The scheme which is described above is stable under a geometrical condition on the family
of meshes which is considered. Since the fluxes are consistent 
 and the scheme is
conservative, it also satisfies a property of 
``weak consistency"\index{Consistency!Weak --}, that is, as in the one
dimensional case (see remark \refp{wc}  of Section \ref{ell1D}), the exact solution of
(\ref{Lapl-2D}) satisfies the numerical scheme with an error which tends to 0 in 
$L^\infty(\Omega)$ for the weak-$\star$ topology. Under adequate restrictive 
assumptions,
the convergence of the scheme can be deduced, see \citet{Fa1}.

Numerical tests were performed for the Laplace operator  and for operators  of the
type $ - \div ($ $\Lambda$ $\mathbf{grad} .)$, where $\Lambda$ is a variable and
discontinuous matrix (see \citet{Fa1});  the discontinuities of $\Lambda$ are
treated in a similar way  as in the 1D case (see Section~\ref{ell1D}). Comparisons
with solutions which were obtained by the bilinear finite element method\index{Finite element method}, 
and with
known analytical solutions, were performed. The results given by the VF9 
scheme and by the finite element scheme were very similar. 

The two drawbacks of this method are the fact that it is a 9-points scheme, and
therefore computationally expensive, and that it yields a nonsymmetric matrix even
if the original continuous operator is symmetric.  Also, its generalization to three
dimensions is somewhat complex.

\begin{remark} \normalfont \label{cafe}
The proof of convergence of this scheme is hindered by the lack of 
consistency for the discrete adjoint operator (see Section \ref{ellcvf4}). 
An error estimate is also difficult to obtain because the numerical 
flux at an interface $\edgecvcvv$ cannot be written under the form 
$\transcvcvv (u_\cv- u_\cvv)$ with $\transcvcvv >0$. Note, however, 
that under some geometrical
assumptions on the mesh, see \citet{Fa1} and \citet{coudiere},
error estimates may be obtained. 
\end{remark}

\subsection{General meshes and schemes}  \label{vfquatre}  

Let us now turn to the discretization of convection-diffusion problems on general structured or non structured grids,  consisting of
any polygonal (recall that we shall call ``polygonal" any polygonal 
domain of $\R^2$ or  polyhedral domain or
$\R^3$) control volumes \index{Control volume} (satisfying adequate geometrical conditions which are stated in the sequel)
and not necessarily  ordered in a Cartesian grid.  The advantage of  finite volume
schemes using    non structured  meshes is clear for convection-diffusion equations.    
On one hand, the stability and convergence properties of the finite volume 
scheme (with an upstream choice for the convective flux)        ensure a robust scheme 
for any admissible mesh as defined in Definitions \refp{meshdirichlet} and
\refp{meshneuman} below, without any need for refinement in the areas of a
large convection flux. On the other hand,  the use of a non structured  mesh
allows
the computation of a solution for any shape of the physical domain.

\medskip

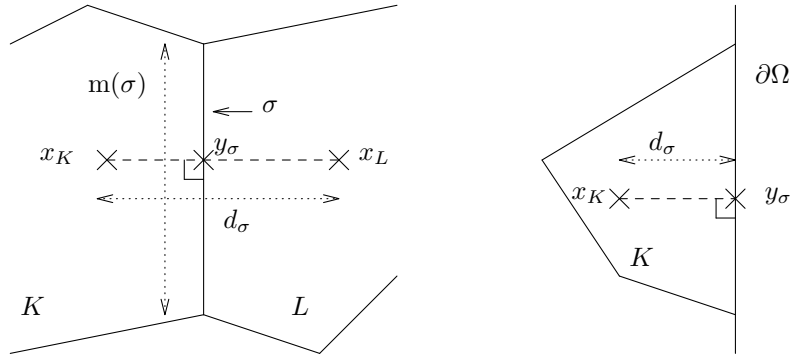
\begin{figure}[ht]
\begin{center}
\begin{tikzpicture}[x=0.00066667in, y=0.00066667in, line cap=round, line join=round]
  \draw (5487,1062) -- (5637,1062);
  \draw (4137,1512) -- (4737,612);
  \draw (1437,1587) -- (1587,1437);
  \draw (1587,1587) -- (1437,1437);
  \node[inner sep=0pt, anchor=base west] at (237,1512) {$x_K$};
  \draw (5487,1212) -- (5487,1062);
  \draw (5637,2712) -- (5637,12);
  \draw (5637,2412) -- (4137,1512);
  \draw (4737,612) -- (5637,312);
  \draw[dashed] (4737,1212) -- (5637,1212);
  \draw[dotted] (4737,1512) -- (5637,1512);
  \draw (4857,1542) -- (4737,1512) -- (4857,1482);
  \draw (5517,1482) -- (5637,1512) -- (5517,1542);
  \node[inner sep=0pt, anchor=base west] at (2712,1512) {$x_L$};
  \node[inner sep=0pt, anchor=base west] at (4362,1212) {$x_K$};
  \node[inner sep=0pt, anchor=base west] at (612,2037) {$\medge$};
  \node[inner sep=0pt, anchor=base west] at (1662,987) {$\dedge$};
  \node[inner sep=0pt, anchor=base west] at (1962,1887) {$\edge$};
  \node[inner sep=0pt, anchor=base west] at (1587,1587) {$y_\edge$};
  \node[inner sep=0pt, anchor=base west] at (5787,2112) {$\partial\O$};
  \node[inner sep=0pt, anchor=base west] at (87,312) {$K$};
  \node[inner sep=0pt, anchor=base west] at (2187,312) {$L$};
  \node[inner sep=0pt, anchor=base west] at (4812,687) {$K$};
  \node[inner sep=0pt, anchor=base west] at (4962,1587) {$\dedge$};
  \node[inner sep=0pt, anchor=base west] at (5862,1212) {$y_\edge$};
  \draw (1362,1362) -- (1512,1362);
  \draw (4662,1287) -- (4812,1137);
  \draw (4662,1137) -- (4812,1287);
  \draw (5562,1287) -- (5712,1137);
  \draw (687,1437) -- (837,1587);
  \draw (2487,1587) -- (2637,1437);
  \draw (2487,1437) -- (2637,1587);
  \draw (687,1587) -- (837,1437);
  \draw (5562,1137) -- (5712,1287);
  \draw[dashed] (762,1512) -- (2562,1512);
  \draw[dotted] (687,1212) -- (2562,1212);
  \draw (807,1242) -- (687,1212) -- (807,1182);
  \draw (2442,1182) -- (2562,1212) -- (2442,1242);
  \draw[dotted] (1212,2412) -- (1212,312);
  \draw (1242,2292) -- (1212,2412) -- (1182,2292);
  \draw (1182,432) -- (1212,312) -- (1242,432);
  \draw (1587,1887) -- (1887,1887);
  \draw (1707,1917) -- (1587,1887) -- (1707,1857);
  \draw (1362,1512) -- (1362,1362);
  \draw (2412,12) -- (3012,612);
  \draw (12,2412) -- (612,2712);
  \draw (612,2712) -- (1512,2412);
  \draw (1512,2412) -- (1512,312);
  \draw (12,12) -- (1512,312);
  \draw (1512,2412) -- (3012,2712);
  \draw (1512,312) -- (2412,12);
\end{tikzpicture}
\caption{Admissible meshes }
\label{CCfig1}
\end{center}
\end{figure}

We saw in the previous section that a consistent discretization of the normal flux
$- \grad u \cdot \bfn$ over the interface of two control volumes $\cv$ and $\cvv$ may be
performed with  a differential  quotient involving values of the unknown located on
the orthogonal  line to the interface between $\cv$ and $\cvv$, on either side of
this interface. This remark suggests the following definition of admissible finite
volume meshes for the discretization of diffusion problems. We shall only consider 
here, for the sake of simplicity,  the case of polygonal domains. The case of
domains with a regular boundary does not introduce any supplementary difficulty
other than complex notations. The definition of admissible meshes 
 and notations introduced in this definition are illustrated in 
 Figure \ref{CCfig1}

\index{Dirichlet boundary condition|)}

\begin{defi}[Admissible meshes]\label{meshdirichlet}  \normalfont
\index{Mesh!Admissible --!--for Dirichlet boundary conditions}
Let $\O$ be an open bounded  polygonal subset of $\R^d$, $d =2,$ or $3$. 
An admissible finite volume mesh of $\O$, denoted by $\mesh$, is given by 
a family of ``control volumes"\index{Control volume}, which are  
 open polygonal convex subsets of 
$\O$ , a family of subsets of $\Obar$ contained in 
hyperplanes of $\R^d$, denoted by  
$\edges$ (these are the edges\index{Edge} (two-dimensional) or  
sides  (three-dimensional) of the control volumes), with strictly positive 
($d-1$)-dimensional measure, and a family of points of $\O$ denoted by 
$\points$ satisfying the  following properties (in fact, we shall denote, somewhat 
incorrectly, by $\mesh$ the family of control volumes):

\begin{enumerate}
 
\item[(i)] The closure of the union of all the control volumes is
$\overline{\O}$;  

\item[(ii)]  For  any $\cv \in \mesh$, there exists  a subset $\edgescv$ of $\edges$
such that $\dr \cv  = \cvbar\setminus \cv = 
\dsp{\cup_{\edge \in \edgescv}}
\overline{\edge} $. Furthermore, $\edges=\dsp{\cup_{\cv \in \mesh}\edgescv}$.

\item[(iii)] For any $(\cv,\cvv)\in \mesh^2$ with $\cv \neq \cvv$, either the 
($d-1$)-dimensional Lebesgue  measure of 
$\cvbar \cap \cvvbar$ is 0  or 
 $\cvbar \cap \cvvbar = \overline\edge$ for some $\edge \in \edges$, which will then be
denoted by $ \edgecvcvv$.

\item[(iv)] The family  $\points = (\centercv)_{\cv \in \mesh}$ is  such that
$\centercv \in \overline \cv$ (for all $\cv \in \mesh$) and, if $\edge =
\edgecvcvv$, it is assumed  that  $\centercv \ne \centercvv$, and
that the straight line $\Dcvcvv$ going through $\centercv$ and $\centercvv$ is
  orthogonal to $\edgecvcvv$.

\item[(v)] For any $\edge \in \edges$ such that $\edge \subset \dr \O$, let 
$\cv$ be the control volume such that $\edge \in \edgescv$. If $\centercv 
\notin \edge$, let $\Dcvedge$ be the straight line going through
$\centercv$  and orthogonal to $\edge$,  then the condition $\Dcvedge \cap \edge 
\ne \emptyset$ is assumed; let $\centeredge =\Dcvedge \cap \edge$.
\end{enumerate}

In the sequel, the following notations are used. 

The mesh size is defined 
by:
$\size(\Mesh)= \sup\{\diam(\cv),$ $\cv\in \mesh\}$.

For any $\cv \in \mesh$ and $\edge \in \edges$, $\mcv$ is the $d$-dimensional 
Lebesgue measure of $\cv$ (it is the area of   $\cv$ in the two-dimensional case and
the volume in the three-dimensional case)
and $\medge$ the ($d-1$)-dimensional measure of
$\edge$. 

The set of interior (resp. boundary) edges is denoted by $\edgesint$
(resp. $\edgesext$), that is $\edgesint = \{\edge \in \edges;$ $\edge \not \subset
\partial \O \}$ (resp.  $\edgesext = \{\edge \in \edges;$ $\edge \subset
\partial \O \}$).

The set of neighbours of $\cv$ is denoted by
$\Nv$, that is $\Nv = \{\cvv \in \mesh;$ $\exists \edge \in \edgescv$, 
$\overline{\edge} =
\cvbar \cap \cvvbar\}$. 

If $\edge = \edgecvcvv$, we denote by
$\dedge$ or $d_{\edgecvcvv}$ the Euclidean distance  between
$\centercv$ and $\centercvv$ (which is positive) and by $\dcvedge$ the 
distance from $\centercv$ to $\edge$. 

If
$\edge \in \edgescv \cap \edgesext$, let
$\dedge$ denote the Euclidean distance between $\centercv$ and 
$\centeredge$ (then, $\dedge=\dcvedge$).

For any $\edge \in \edges$; the ``transmissibility"\index{Transmissibility} 
through $\edge$ is defined by 
$\transedge =  {\medge} /{\dedge}$ if $\dedge \ne 0$.

In some results and proofs given below, there are summations over $\edge \in 
\edges_0$, with $\edges_0=\{\edge \in \edges;$ $\dedge \ne 0\}$. For 
simplicity,  (in these results and proofs) $\edges=\edges_0$ is assumed.

\end{defi}

\begin{remark} \normalfont 

(i) The definition of $\centeredge$ for $\edge \in \edgesext$ 
requires  that $\centeredge \in \edge$. However, In many cases, this condition may
be relaxed. The condition $\centercv \in \cvbar$ may also be relaxed as 
described, for instance, in Example \ref{tm} below.

(ii) The condition  $\centercv \ne \centercvv$  if $\edge =
\edgecvcvv$, is in fact quite easy to satisfy:  two neighbouring control 
volumes $\cv, \cvv$ which do not satisfy it just have to be collapsed into a new 
control volume $\cvw$ with $\centercvw = \centercv = \centercvv$, and the 
edge $\edgecvcvv$ removed from the set of edges.  The new mesh thus 
obtained is admissible. 
\end{remark}

\begin{example}[Triangular meshes] \label{tm} \index{Mesh!Triangular --}
\normalfont Let $\O$ be
an open bounded polygonal subset of $\R^2$. Let $\mesh$ be a family of
open triangular disjoint subsets of $\O$ such that two triangles having a 
common edge have also two common vertices. Assume that all angles of the triangles
are less than $\pi / 2$.  This last condition is sufficient for the orthogonal
bisectors to intersect inside each triangle, thus naturally defining the points 
$\centercv \in \cv$. One  obtains an admissible mesh. In the case of an
elliptic operator, the finite volume scheme defined on such a grid using
differential quotients for the approximation of the normal flux yields a 4-point
scheme \citet{VF4}. This scheme does not lead to a finite 
difference scheme consistent  with the
continuous diffusion operator (using a Taylor expansion).
The consistency is only verified for the approximation of the
fluxes, but this, together with the conservativity  of the scheme yields the
convergence of the scheme, as it is proved below.

Note that the condition that all angles of the triangles
are less than $\pi / 2$ (which yields $\centercv \in \cv$) may be relaxed 
(at least for the triangles the closure of which are in $\O$) to the so called 
``strict Delaunay 
condition"\index{Delaunay condition}
 which is that the closure of the circumscribed circle to each triangle of the 
mesh does not contain any other triangle of the mesh. For such a mesh, 
the point $\centercv$ (which is the intersection of the orthogonal 
bisectors of the edges of $\cv$) is not always in $\cv$, but
the scheme \eqref{easyvf1}-\eqref{easyvf3} is convenient since 
\eqref{easyvf2} yields a consistent  approximation of the diffusion fluxes 
and since the transmissibilities \index{Transmissibility} 
(denoted by $\transcvcvv$) are positive.
\end{example}

\begin{example}[Vorono\"{\i} meshes]\label{evoronoi} \normalfont
\index{Vorono\"{\i}} \index{Mesh!Vorono\"{\i} --} Let
$\O$ be an open bounded polygonal subset of $\R^d$. An admissible 
finite volume mesh can be built by using the so called ``Vorono\"{\i}"
technique. Let $\P$ be a family of points of $\Obar$. For example, this 
family may be chosen as $\P=\{(k_1
h,\ldots,k_d h)$, $k_1,\ldots k_d \in \Z \} \cap \O$, for a given $h>0$. 
The control volumes \index{Control volume} of the Vorono\"{\i} mesh are defined with respect to each point  $x$ of $\P$ by 

\[
\cv_x = \{ y \in \O, \vert x - y \vert < \vert z - y \vert, \; \forall z \in
\P, \,  z \neq x \},  
\]
where $\vert x- y \vert$ denotes the Euclidean distance between 
$x$ and $y$. 
Vorono\"{\i} meshes are admissible in the sense of Definition
\ref{meshdirichlet} if the assumption ``on the boundary", namely part 
$(v)$ of Definition \ref{meshdirichlet}, 
is satisfied. Indeed, this is true, in particular, if the number of 
points $x \in \P$ which are located on $\dr \O$ is ``large enough". 
Otherwise, the
assumption $(v)$ of Definition \ref{meshdirichlet} may be replaced by the 
weaker assumption ``$d(\centeredge,\edge) \le
\size(\mesh)$ for any $\edge \in \edgesext$" which is much easier to satisfy. Note
also that a slight modification of the treatment of the boundary 
conditions in the finite volume
scheme \eqref{ellschema}-\refep{ellschemad} allows us to obtain convergence and error
estimates results (as in theorems \refp{ellcvgce} and \refp{ellesterr}) for 
all Vorono\"{\i} meshes. 
This modification is the obvious generalization of the scheme described in the first item of Remark \refp{angot1} for the 1D case. It 
consists in replacing, for $\cv \in \mesh$ such that $\edgescv \cap \edgesext \ne \emptyset$, the equation \eqref{ellschema}, associated to this control volume, by the equation $u_\cv=g(z_\cv)$, where $z_\cv$ is some point on $\dr \O \cap \dr \cv$. In fact,  Vorono\"{\i} meshes often satisfy  the following property:
\[\edgescv \cap 
\edgesext \ne \emptyset  \Rightarrow \centercv 
\in \dr \O\]
 and the mesh is therefore admissible in the sense of 
Definition \ref{meshdirichlet} (then, the scheme 
\eqref{ellschema}-\refep{ellschemad} yields $u_\cv=g(\centercv)$ if $\cv 
\in \mesh$ is such that $\edgescv \cap 
\edgesext \ne \emptyset$).

An advantage of the Vorono\"{\i} method is that it easily leads to meshes on non polygonal domains $\O$.
\end{example}

Let us now introduce the space of piecewise constant functions associated to an admissible mesh and some ``discrete $H^1_0$" norm  for this space. 
This discrete norm will be used to obtain stability properties  which are given by some estimates on the approximate solution of  a finite volume scheme.

\begin{defi}[Discrete space and norm] \label{Xmesh} \label{defh10d} \index{Discrete!-- $H_0^1$ norm} \normalfont 
Let $\O$ be an open bounded polygonal subset of $\R^d$, $d=2$ or $3$, and $\Mesh$  be an admissible finite volume mesh in the sense of Definition \refp{meshdirichlet}. .
Let $X(\Mesh)$ as the set of functions from $\O$ to $\R$ which are constant over each control volume of the mesh.

For  $u \in X(\Mesh)$, define the discrete  $H_0^1$ norm by 
\begin{equation} 
\Vert u \Vert_{1,\Mesh} = \Bigl(\sumedges \transedge (D_\edge u)^2 \Bigr)^{\frac 1 2},
\label{normh1d}
\end{equation}
where $\transedge = \medge / \dedge$ and $D_\edge u = \abs {u_\cv - u_\cvv}$ if $\edge  \in \edgesint$, $\edge=\edgecvcvv$, $D_\edge u = \abs {u_\cv} $ if $\edge  \in \edgesext \cap \edgescv$,
and where $u_\cv$ denotes the value taken by $u$ on the control volume $\cv$ and the sets $\edges$, $\edgesint$, $\edgesext$ and $\edgescv$ are defined in Definition
\refp{meshdirichlet}.
\end{defi}

The discrete $H_0^1$ norm is used in the following sections to prove the congergence of finite volume schemes and, under some regularity conditions, to give error estimates. 
It is related to the  $H_0^1$ norm, see the convergence of the norms in Theorem \ref{ellcvgce}.
One of the tools used below is the following ``discrete Poincar\'e inequality" which may also be found in \citet{temam}:

\begin{lemma}[Discrete Poincar\'e inequality]
\label{thinpoin}\index{Poincar\'e inequality!-- Dirichlet boundary conditions}
Let $\O$ be an open bounded polygonal subset of $\R^d$, $d=2$ or $3$, $\Mesh$ an admissible
finite volume mesh in the sense of Definition \ref{meshdirichlet} and $u \in
X(\Mesh)$  (see Definition \ref{Xmesh}), then 
\be
\Vert u \Vert_{L^2(\O)} \le \diamO \Vert u  \Vert_{1,\Mesh}, 
\label{ellinpoin}
\ee
where $\Vert \cdot  \Vert_{1,\Mesh}$ is the discrete $\Hunz$ norm defined in
Definition \refp{defh10d}.
\end{lemma}

\begin{remark}[Dirichlet  condition on part  of the boundary]
\label{dirich-part}
This lemma gives a discrete Poincar\'e inequality for Dirichlet boundary conditions on the boundary $\partial \Omega$. 
In the case of a Dirichlet condition on part of the boundary only, it is still
possible to prove a Discrete boundary condition provided that the polygonal bounded open set $\Omega$ is also connex, thanks to Lemma 
\refp{thinpoin} proven in the sequel. 
\end{remark}

\proof  of Lemma \ref{thinpoin}

\smallskip

For $\edge \in \edges$, define $\chi_\edge$ from $\R^\dim \times \R^\dim$ 
to $\{0,1\}$ by $\chi_\edge(x,y) = 1$ if $\edge \cap [x,y] \not = 
\emptyset$ and $\chi_\edge (x,y) = 0$
otherwise.

~

Let $u \in X(\Mesh)$. Let $\bfd$ be a given unit 
vector. For all $x \in \O$, let $\D_x$ be the
semi-line defined by its origin, $x$, and the vector $\bfd$. Let $y(x)$ such 
that $y(x) \in \D_x \cap \dr \O$ and $[x,y(x)] \subset \Obar$, where 
$[x,y(x)]=\{tx+(1-t)y(x),$ $t\in [0,1]\}$ (i.e. $y(x)$ is the first point 
where $\D_x$ meets $\dr \O$). 

~

Let $\cv \in \mesh$. For a.e. $x \in \cv$, one has

\[
\vert u_\cv \vert  \le \sumedges  D_\edge u  \, \chi_\edge (x,y(x)),
\]

where the notations $D_\edge u$ and $u_\cv$ are defined in Definition 
\refp{defh10d}.  We write the above inequality  for a.e $x \in \Omega$ and
not for all $x \in \Omega$ in order to  account for the cases where an edge 
or a vertex of the mesh is included in the semi-line  $[x,y(x)]$; in both 
cases one may not write the above inequality, but there are only a
finite number of edges and vertices, and since $\bfd$ is fixed,
the above inequality may be 
written almost everywhere.

Let
$c_\edge = \vert \bfd \cdot \bfn_\edge \vert$ (recall that $\xi \cdot \eta$ denotes the usual scalar product 
of $\xi$ and $\eta$ in $\R^\dim$). By the Cauchy-Schwarz  inequality, the 
above inequality yields:

\be
\vert u_\cv \vert^2 \le \sumedges \frac {(D_\edge u)^2} {\dedge \cedge}
\chi_\edge(x,y(x)) \sumedges \dedge \cedge \chi_\edge(x,y(x)), \text{ for a.e. }x \in \cv.
\label{bon}
\ee

Let us show that, for a.e. $x \in \O$,

\be
\sumedges \dedge\cedge \chi_\edge(x,y(x)) \le \diamO.
\label{prelim}
\ee

Let $x \in \cv$, $\cv \in \mesh$, such that $\edge 
\cap [x,y(x)]$ contains at most one point, for all $\edge \in \edges$, 
and $[x,y(x)]$ does not contain any vertex of $\mesh$ (proving \eqref{prelim} 
for such points $x$ leads to \eqref{prelim} a.e. on $\O$, since $\bfd$ is fixed).

There exists $\edge_x \in \edgesext$ such that $y(x) \in \edge_x$.
Then, using the fact that the control volumes are convex, one has:

\[
\sumedges \chi_{\edge}(x,y(x)) d_\edge c_\edge=\vert 
(\centercv-x_{\edge_x}) \cdot \bfd  
\vert.
\]

Since $\centercv$ and  $x_{\edge_x} \in \Obar$ (see Definition 
\ref{meshdirichlet}), this gives \eqref{prelim}.

\smallskip

Let us integrate \eqref{bon} over $\O$;   \eqref{prelim} gives

\[
\summesh \int_\cv \vert u_\cv\vert^2 dx \le \diamO \sumedges 
\frac {(D_\edge u)^2} {\dedge \cedge} \int_\O 
\chi_\edge(x,y(x))  dx.
\]

Since $\int_\O \chi_\edge(x,y(x))  dx \le \diamO \medge \cedge$, this last inequality
yields

\[  
\summesh \int_\cv \vert u_\cv \vert^2 dx \le (\diamO)^2 \sumedges 
\vert D_\edge u \vert^2 \frac {\medge} {\dedge}
dx.
\]

Hence the result. 
\eof

Let $\mesh$ be an admissible mesh. Let us now define a
finite volume scheme \index{Finite volume scheme!-- for 
elliptic problems!-- in two or three 
space dimensions|(}
 to discretize \eqref{elldifstaf}, \refep{elldifstab}. 

Let
 
\be
f_\cv = \frac 1 {\mcv} \int_\cv f(x) dx,  \forall \cv\in
\mesh. \label{frhs}
\ee
  
Let $(u_\cv)_{\cv \in \mesh}$  denote the discrete unknowns.
In order to describe the scheme in the most general way, one
introduces some auxiliary unknowns (as in the 1D case, see Section
\ref{ell1D}), namely the fluxes $F_{\cv,\edge}$, for all $\cv \in
\mesh$ and $\edge \in \edgescv$, and some (expected) approximation of
$u$ in $\edge$, denoted by $u_{\edge}$, for all $\edge \in \edges$. 
For $\cv \in \mesh$ and $\edge \in \edgescv$, let $\n_{\cv,\edge}$ denote 
 the normal unit vector to
$\edge$ outward to $\cv$ and $\vcvedge = \int_\edge \vbf(x) \cdot 
\n_{\cv,\edge} \dfrontiere (x)$. Note  that $\dfrontiere$ is the 
integration symbol for the $(\dim-1)$-dimensional Lebesgue measure on the 
considered hyperplane. 
In order to discretize the convection 
term\index{Convection term} $\div(\vbf(x) u(x))$ in a stable way (see Section \refp{ell1D}),  
let us define the upstream\index{Upstream} choice $\uedgep$ of $u$ on an edge $\edge$ with respect 
to $\vbf$ in the following way. If  $\edge =  \edgecvcvv$, then  $\uedgep  =  u_\cv$ if 
$\vcvedge \ge 0$,  and  $\uedgep  =  u_\cvv$ otherwise;  if   $\edge
\subset \cv \cap \dr \O$, then  $\uedgep  =  u_{\cv}$ if  $\vcvedge \ge 0$  
 and   $\uedgep  = g(\centeredge)$ otherwise. 

\smallskip

Let us first assume that the points $x_\cv$ are located in the interior of each
control volume, and are therefore not located on the edges, hence  
$\dcvedge > 0$ for any $\edge \in \edgescv$, where $\dcvedge$ is the 
distance from $\centercv$ to $\edge$. A finite volume scheme can be
defined by the following set of equations:

\be        
\sum_{\edge \in \edgescv} F_{\cv,\edge} 
+ \sumedgecv \vcvedge \uedgep + b  \mcv u_\cv    = 
\mcv f_\cv, \ \forall \cv\in\mesh, 
\label{easyvf1}
\ee 

\be
F_{\cv,\edge}=  - \transcvcvv  (u_\cvv - u_\cv), \, \forall  \edge \in \edgesint, \; \mathrm{if} \;
\edge = \edgecvcvv,
\label{easyvf2}
\ee

\be
F_{\cv,\edge}   = - \transedge (g(\centeredge)-u_{\cv}), \, \forall 
\edge \in \edgesext \text{ such that }\edge \in \edgescv. 
\label{easyvf3}
\ee

In the general case, the center of the cell may be located on an edge. This is the 
case for instance when constructing Vorono\"{\i} meshes with some of the original points 
 located on the  boundary $\partial \O$. In this case, the following formulation of the finite volume scheme
is valid, 
and is equivalent to the above scheme if no cell center is located on an edge:

\be        
\sum_{\edge \in \edgescv} F_{\cv,\edge} 
+ \sumedgecv \vcvedge \uedgep + b  \mcv u_\cv    = 
\mcv f_\cv, \ \forall \cv\in\mesh, 
\label{ellschema}
\ee 

\be
F_{\cv,\edge}= -F_{\cvv,\edge}, \, \forall  \edge \in \edgesint, \; \mathrm{if} \;
\edge = \edgecvcvv,
\label{ellschemab}
\ee

\be
F_{\cv,\edge} \dcvedge = - \medge (u_{\edge}-u_{\cv}), \, \forall 
\edge \in \edgescv, \, \forall \cv \in \mesh,
\label{ellschemac}
\ee

\be        
u_\edge  =  g(\centeredge),
\ \forall \edge \in \edgesext.  
\label{ellschemad}
\ee

Note that \eqref{ellschema}-\eqref{ellschemad} always lead, after an easy elimination of the
auxiliary unknowns, to a linear system of
$N$ equations  \index{Finite volume scheme!-- for 
elliptic problems!-- in two or three 
space dimensions|)} with
$N$ unknowns, namely the $(u_\cv)_{\cv \in \mesh}$, with $N=$
$\card(\mesh)$.

\begin{remark} \normalfont~
\benum
\item Note that one may have, for some $\edge \in \edgescv$,
$\centercv \in \edge$, and therefore, thanks to \eqref{ellschemac},
$u_{\edge}=u_{\cv}$. 
\item The choice $u_\edge  =  g(\centeredge)$ in \eqref{ellschemad} needs 
some discussion. Indeed, this choice is possible since $g$ is assumed to 
belong to $C(\dr \O,\R)$ and then is everywhere defined on $\dr \O$. In 
the case where the solution to \eqref{elldifstaf}, \refep{elldifstab} 
belongs to $H^2(\O)$ (which yields $g \in C(\dr \O,\R)$), it is clearly 
 a good choice since it yields the consistency  of fluxes (even 
though an error estimate also holds with other choices for $u_\edge$, 
the choice given below is, for instance, possible). If $g \in H^{1/2}$ (and not continuous), the 
value $g(\centeredge)$ is not necessarily defined. Then, another choice 
for $u_\edge$ is possible, for instance,

\[
u_{\edge}=\frac 1 {\medge} \int_\edge g(x)\dfrontiere(x).
\]

With this latter choice for $u_\edge$, a convergence result also holds, see Theorem 
\ref{ellcvgcediraf}. 

\eenum
\end{remark}

For the sake of simplicity, it is  assumed  in Definition \ref{meshdirichlet} that
$\centercv \neq \centercvv$, for all $\cv,$ $\cvv \in \mesh$. This condition may be 
relaxed; it simply allows an easy expression of the numerical flux 
$F_{\cv,\edge}= - \transcvcvv  (u_\cvv-u_\cv)$ if $\edge = \edgecvcvv$.

\subsection{Existence and estimates}
 Let us first  prove the existence of the
approximate solution and an estimate on this solution. This estimate 
ensures the stability of the scheme and will be
obtained by using  the discrete Poincar\'{e} inequality\index{Discrete!-- Poincar\'e inequality} \eqref{ellinpoin} and will
yield convergence thanks to a compactness theorem given in Section
\refp{rellichd}.

\begin{lemma}[Existence and estimate] \label{ellexistu} Under Assumptions 
\ref{ellH}, let $\Mesh$ be an
admissible mesh in the sense of Definition \refp{meshdirichlet}; there exists
a unique solution
$(u_\cv)_{\cv \in \mesh}$ to equations \eqref{ellschema}-\eqref{ellschemad}. 

Furthermore, assuming $g = 0$ and   defining 
$u_{\mesh} \in X(\Mesh)$ (see Definition \ref{Xmesh}) by $u_\mesh(x) = u_{\cv}$
for $a.e.$ $x \in \cv$, and for any
$\cv \in \mesh$, the following estimate\index{Estimate on the approximate  
solution!-- for an elliptic equation} holds: 

\be
\Vert u_\mesh \Vert_{1,\Mesh} \le \diamO \Vert f \Vert_{L^2(\O)},
\label{ellestimx}
\ee

where $\Vert \cdot  \Vert_{1,\Mesh}$ is the discrete $\Hunz$ norm defined in
Definition \ref{defh10d}. 
\end{lemma}
    
\proof  of Lemma \ref{ellexistu}

\smallskip

Equations \eqref{ellschema}-\eqref{ellschemad} lead, after an easy elimination of the
auxiliary unknowns, to a linear system of
$N$ equations with
$N$ unknowns, namely the $(u_\cv)_{\cv \in \mesh}$, with $N=$
$\card(\mesh)$.

\medskip

\textit{Step 1 (existence and uniqueness)}

Assume that $(u_\cv)_{\cv \in \mesh}$ satisfies this linear system with
$g(\centeredge)= 0$ for any $\edge \in \edgesext$, and $f_{\cv}=0$ for all $\cv \in
\mesh$. Let us multiply  \eqref{ellschema} by $u_\cv$ and  sum over $\cv$; from 
\eqref{ellschemab} and \eqref{ellschemac} one deduces

\be
b \summesh \mcv u_\cv^2 +
\summesh \sumedgecv \Fcvedge u_\cv + 
\summesh \sumedgecv \vcvedge  \uedgep  u_\cv = 0,
\label{ell1e}
\ee

which gives, reordering the summation over the set of edges
\be
b \summesh \mcv u_\cv^2 + \sumedges \transedge (D_\edge u)^2 + 
\sumedges \vedge \Bigl(\uedgep - \uedgem \Bigr) \uedgep = 0,
\label{ell7e} 
\ee

where

 $\vert D_\edge u \vert=\vert u_{\cv}-u_{\cvv} \vert$, if
$\edge=\edgecvcvv$ and $\vert D_\edge u \vert=\vert u_{\cv}\vert$, if $\edge \in
\edgescv \cap \edgesext$;

 $\vedge = \vert \int_\edge \vbf (x) \cdot \n \dfrontiere(x) \vert$, $\n$ being a unit
normal  vector to $\edge$;

  $\uedgem $ is the downstream value to $\edge$ with respect to $\vbf$,
i.e. if  $\edge =  \edgecvcvv$, then  $\uedgem  =  u_\cv$ if  $\vcvedge \le 0$, 
and  $\uedgem  =  u_\cvv$ otherwise;  if   $\edge \in \edgescv \cap \edgesext$,
then  $\uedgem  =  u_{\cv}$ if  $\vcvedge \le 0$ and $\uedgem  =  u_{\edge}$ if  $\vcvedge > 0$. 
  
Note that $u_{\edge}=0$ if $\edge \in \edgesext$.

\medskip

Now, remark that 
\begin{equation}
\sumedges  \vedge \uedgep (\uedgep - \uedgem)  
= \frac 1 2 \sumedges \vedge \Bigr( (\uedgep - \uedgem)^2 +
(\uedgep^2-\uedgem^2) \Bigl)  \label{vfqconv1}
\end{equation}

and, thanks to the assumption $\div \vbf \ge 0$,  
\begin{equation}
\sumedges  \vedge  (\uedgep^2-\uedgem^2)   
=  \summesh \Bigl(\int_{\dr \cv} \vbf(x) \cdot \n_{\cv} \dfrontiere(x) \Bigr) 
u_\cv^2 = \int_\O (\div \vbf(x)) u_\mesh^2(x)dx \ge 0.
\label{vfqconv11}
\end{equation}
Hence,

\be
b\norm{u_\mesh}{L^2(\O)}^2 + \norm{u_\mesh}{1,\mesh}^2 = b \summesh \mcv u_\cv^2 + \sumedges \transedge (D_\edge u)^2 \le 0,
\label{ell7f} 
\ee

One deduces, from \eqref{ell7f}, that $u_{\cv}=0$  for all $\cv \in
\mesh$.

This proves the existence and the uniqueness of the solution
$(u_\cv)_{\cv \in \mesh}$, of the linear system given by 
\eqref{ellschema}-\eqref{ellschemad}, for any 
 $\{g(\centeredge),$ $\edge \in \edgesext\}$ and $\{f_{\cv},$ $\cv \in \mesh\}$.

\medskip
 
\textit{Step 2 (estimate)}

Assume $g=0$. Multiply 
\eqref{ellschema} by $u_\cv$, sum  over $\cv$; then, thanks to
\eqref{ellschemab},  \eqref{ellschemac}, \eqref{vfqconv1} and \eqref{vfqconv11}  
one has

\[
b \Vert u_\mesh \Vert_{L^2(\O)}^2 + \Vert u_\mesh \Vert_{1,\Mesh}^2 \le \summesh 
\mcv f_\cv
u_\cv.
\]

By  the Cauchy-Schwarz  inequality, this inequality yields

\[
\Vert u_\mesh \Vert_{1,\Mesh}^2 \le (\summesh \mcv u_\cv^2)^{\frac 1 2} 
(\summesh \mcv f_\cv^2)^{\frac 1 2} \le 
\norm{f}{L^2(\O)}\norm{u_\mesh}{L^2(\O)}.
\]

Thanks to the discrete Poincar\'e  inequality \eqref{ellinpoin}, this yields  
$ \Vert u_\mesh \Vert_{1,\Mesh} \le  \Vert f \Vert_{L^2(\O)} \diamO$, which concludes
the proof of  the lemma.
\eof

Let us now state a discrete maximum principle which is satisfied by the
scheme \eqref{ellschema}-\eqref{ellschemad}; this is an interesting stability
property, even though it will not be used in the proofs of the convergence
and error estimate. 

\begin{prop} \label{maxppe}\index{Stability}\index{Discrete!-- maximum 
principle} Under Assumption \refp{ellH}, let $\mesh$ be an admissible mesh
in the sense of Definition \refp{meshdirichlet}, let $(f_\cv)_{\cv \in \mesh}$ be
defined by \eqref{frhs}. If $f_\cv \ge 0$ for all
$\cv \in \mesh$, and $g(\centeredge) \ge 0$, for all $\edge \in \edgesext$, then the
solution $(u_\cv)_{\cv\in \mesh}$ of \eqref{ellschema}-\eqref{ellschemad}  satisfies
$u_\cv \ge 0$ for all $\cv \in \mesh$.
\end{prop}

\proof of Proposition \ref{maxppe}

\smallskip

Assume that
$f_\cv \ge 0$ for all $\cv \in \mesh$ and $g(\centeredge) \ge
0$ for all $\edge \in \edgesext$. Let $a=\min\{u_\cv, \cv
\in \mesh\}$.
Let $\cv_0$ be a control volume such that $u_{\cv_0} =a$. Assume first 
that $\cv_0$ is an ``interior" control volume, in the
sense that $\edgescv \subset \edgesint$, and that
$u_{\cv_0} \le 0$. Then, from \eqref{ellschema}, 

\begin{equation}
  \sum_{\edge \in \edges_{\cv_0}}   F_{\cv_0,\edge}  +  \sum_{\edge \in \edges_{\cv_0}} 
v_{{\cv_0},\edge} \uedgep \ge 0; 
\end{equation}

since for any neighbour $\cvv$ of $\cv_0$ one has 
$u_{\cvv} \ge u_{\cv_0}$, then, noting that $\div \vbf \ge 0$,  one must
have $u_{\cvv} = u_{\cv_0}$ for any neighbour $\cvv$ of $\cv$. Hence, setting 
$B=\{\cv \in \mesh,$ $u_\cv=a\}$, there exists $\cv \in B$ such that 
$\edgescv \not \subset \edgesint$, that is $\cv$ is a control volume 
``neighbouring the boundary".

Assume then that $\cv_0$ is a control volume neighbouring  the
boundary and that $u_{\cv_0}=a < 0$. Then, for an edge $\edge \in \edgesext\cap
\edgescv$,  relations \eqref{ellschemac}
 and \eqref{ellschemad} yield  $g(\centeredge) <0$, which is
in contradiction with the assumption.  Hence Proposition \ref{maxppe} is proved. 
\eof

\begin{remark} \normalfont The maximum principle immediately yields the existence and uniqueness
of the solution of the numerical scheme  \eqref{ellschema}-\eqref{ellschemad}, which
was proved directly in Lemma \ref{ellexistu}.
\end{remark}

\subsection{Convergence} \label{ellcvf4}

Let us now show the  convergence of approximate solutions obtained  by 
the above  finite volume scheme when the size of the mesh tends to 0. 
One uses Lemma \ref{ellexistu} together with the compactness theorem \ref{kolmh10} given at the end of this chapter to
prove the convergence result. In order to use Theorem  \ref{kolmh10}, one needs the following lemma.

\begin{lemma}\label{h1dtot} Let $\O$ be an open bounded set of $\R^d$,
$d=2$ or $3$. Let $\Mesh$ be an admissible mesh in the sense of Definition
\refp{meshdirichlet}  and $u \in
X(\Mesh)$ (see Definition  \ref{Xmesh}). One defines $\tilde{u}$ by $\tilde{u}=u$ 
a.e. on $\O$, and $\tilde{u}=0$  a.e. on $\R^d \setminus \O$. Then there 
exists $C > 0$, only depending  on $\O$, such that

\be
\norm{\tilde{u}(\cdot+\eta)-\tilde{u}}{L^2(\R^d)}^2 \le \norm{u}{1,\Mesh}^2 \vert \eta
\vert (\vert \eta \vert +
C \, \size(\mesh)), \forall \eta \in \R^d.
\label{itransh10}
\ee

\end{lemma}

\proof of Lemma \ref{h1dtot}

\smallskip

For $\edge \in \edges$, define
$\chi_{\edge}$ from $\R^d \times \R^d$ to $\{0,1\}$ by $\chi_{\edge}(x,y)=1$ if
$[x,y] \cap \edge \neq \emptyset$ and $\chi_{\edge}(x,y)=0$ if
$[x,y] \cap \edge = \emptyset$.

Let $\eta \in \R^d$, $\eta \neq 0$. One has

\[
\vert \tilde u(x+\eta)-\tilde u(x) \vert \le \sumedges \chi_{\edge}(x,x+\eta)
\vert D_{\edge}u \vert, \; \text{ for } \ae x \in \O
\]

(see Definition \refp{defh10d} for the definition of $D_{\edge}u$).

This gives, using  the Cauchy-Schwarz  inequality,

\be
\vert \tilde u(x+\eta)-\tilde u(x) \vert^2 \le \sumedges \chi_{\edge}(x,x+\eta)
\frac {\vert D_{\edge}u \vert^2}{d_\edge c_\edge} \sumedges
\chi_{\edge}(x,x+\eta) d_\edge c_\edge, \;
\text{ for }  \ae x \in \R^d,
\label{astucecos}
\ee

where $c_\edge=\vert \bfn_\edge \cdot \frac {\eta}{\vert \eta \vert} \vert$, and
$\bfn_\edge$ denotes a unit normal vector to $\edge$.

~

Let us  now prove that there exists $C > 0$, only depending  on $\O$, such that 

\be
\sumedges \chi_{\edge}(x,x+\eta) d_\edge c_\edge \le \vert \eta \vert + C \,
\size(\mesh),
\label{aie}
\ee
for a.e. $x \in \R^d$.

~

Let $x \in \R^\dim$ such that $\edge 
\cap [x,x+\eta]$ contains at most one point, for all $\edge \in \edges$, 
and $[x,x+\eta]$ does not contain any vertex of $\mesh$ (proving 
\eqref{aie} for such points $x$ gives \eqref{aie} for a.e. 
$x \in \R^d$, since $\eta$ is 
fixed). Since $\O$  is not assumed to be convex, it may happen that the 
line segment  $[x,x+\eta]$  is not included in $\overline\O$. In order to deal with this, 
let $y, z \in [x,x+\eta]$ such that  $y \ne z$ and $[y,z] \subset \Obar$;
there exist $\cv$, $\cvv \in \mesh$ such that $y \in \cvbar$ and $z \in 
\cvvbar$. Hence, 

\[
\sumedges 
\chi_{\edge}(y,z) d_\edge c_\edge=\vert (y_1-z_1) \cdot \frac {\eta}{\vert 
\eta \vert} \vert,
\]


where $y_1=\centercv$ or $\centeredge$ with $\edge \in \edgesext \cap \edgescv$ 
and $z_1=\centercvv$ or 
$\centeredget$ with $\edget \in \edgesext \cap \edgescvv$, depending on 
the position of $y$ and $z$ in $\cvbar$ or $\cvvbar$ respectively.

Since $y_1=y+y_2$, with $\vert y_2 \vert \le \sizemesh$, and $z_1=z + 
z_2$, with $\vert z_2 \vert \le \sizemesh$, one has

\[
\vert (y_1-z_1) \cdot \frac {\eta}{\vert 
\eta \vert} \vert \le \vert y-z \vert + \vert y_2 \vert + \vert z_2 \vert 
\le \vert y-z \vert + 2 \, \sizemesh
\]

and

\be
\sumedges \chi_{\edge}(y,z) d_\edge c_\edge \le \vert y-z \vert + 2 \, \sizemesh.
\label{aier}
\ee

Note that this yields \eqref{aie} with $C=2$ if $[x,x+\eta] \subset \Obar$.

Since $\O$ has a finite number of sides, the line segment  $[x,x+\eta]$ intersects $\dr \O$ a finite number of times; hence  there exist $t_1, \ldots, t_n$ such that $0 \le t_1 < t_2 < \ldots < t_n \le 1$, $n \le N$, where $N$ only depends on $\O$ (indeed, it is possible to take  $N=2$ if $\O$ is convex and $N$ equal to the number of sides of $\O$ for 
a general $\O$) and such that

\[
\sumedges \chi_{\edge}(x,x+\eta) d_\edge c_\edge =\sum_{\substack{i=1,n-1\\\mathrm{odd } i}} \sumedges  \chi_{\edge}(x_i,x_{i+1}) d_\edge c_\edge,
\]
with $x_i=x + t_i \eta$, for $i=1,\ldots,n$,   $x_i \in \dr \O$ if $t_i \notin \{0,1\}$  and $[x_i,x_{i+1}] \subset  \Obar$ if $i$ is odd. 


Then, thanks to \eqref{aier} with $y=x_i$ and $z=x_{i+1}$, for $i=1,\ldots, 
n-1$,  one has \eqref{aie} with $C=2(N-1)$ (in particular, if $\O$ is 
convex, $C=2$ is convenient for \eqref{aie} and therefore for 
\eqref{itransh10} as we shall see below).

~

In order to conclude the proof of Lemma \ref{h1dtot}, remark that, for all 
$\edge \in \edges$,

\[
\int_{\R^d}\chi_{\edge}(x,x+\eta)dx \le \medge c_\edge \vert \eta \vert.
\]

Therefore, integrating \eqref{astucecos} over $\R^d$ yields, with 
\eqref{aie},

\[
\norm{\tilde u(\cdot+\eta)-\tilde u}{L^2(\R^d)}^2 \le \big(\sumedges \frac
{\medge}{d_\edge} \vert D_{\edge}u \vert^2 \big) \vert \eta \vert (\vert \eta
\vert   + C \, \size(\mesh)).
\]
\eof
We are now able to state the convergence theorem.  
We shall first prove the convergence result in the case of homogeneous Dirichlet boundary conditions, i.e. $g = 0$; thenonhomogenous case is then considered (see Theorem \refp{ellcvgcediraf}), following  \citet{cvnl}.

\begin{theo}[Convergence, homogeneous Dirichlet boundary conditions]
\label{ellcvgce}
\index{Convergence of the finite volume method!-- for an elliptic equation}
Under Assumption \refp{ellH} with $g=0$, let $\Mesh$ be an admissible mesh (in the sense of Definition \refp{meshdirichlet}). 
Let$(u_\cv)_{\cv \in \mesh}$ be the solution of the system given by equations \eqref{ellschema}-\eqref{ellschemad} (existence and uniqueness of $(u_\cv)_{\cv \in \mesh}$ are given in Lemma \ref{ellexistu}). 
Define $u_{\Mesh} \in X(\Mesh)$ by $u_\Mesh(x) = u_{\cv}$ for a.e. $x \in \cv$, and for any $\cv \in \mesh$. 
Then $u_{\Mesh}$ converges in $L^2(\O)$ to the unique variational solution $u\in H^1_0(\O)$ of Problem \eqref{elldifstaf}, \eqref{elldifstab} as $\size(\Mesh)\tends 0$. 
Furthermore $\norm{u_{\Mesh}}{1,\mesh}$ converges to $\norm{u}{H^1_0(\O)}$ as $\size(\Mesh)\tends 0$.
\end{theo}
\begin{remark} \normalfont\label{convplus}~\\
\begin{enumerate}
\item In Theorem \ref{ellcvgce}, the hypothesis $f \in L^2(\O)$ is not necessary.
It is used essentially to obtain a bound on $\norm{u_\mesh}{1,\mesh}$. In order to pass to the limit, the hypothesis ``$f \in L^1(\O)$" is sufficient. 
Then, in Theorem \ref{ellcvgce}, the hypothesis $f \in L^2(\O)$ can be replaced by $f \in L^p(\O)$ for some $p > 1$, if $d=2$, and for $p \ge \frac 6 5$, if $d=3$, provided that the meshes satisfy, for some fixed $\zeta >0$, $\dcvedge \ge \zeta \dedge$, for all $\edge \in \edgescv$ and for all control volumes $\cv$.
Indeed, one obtains, in this case, a bound on $\norm{u_\mesh}{1,\mesh}$ by using a ``discrete Sobolev inequality" (proved in Lemma \refp{lq}). 

It is also possible to obtain convergence results, towards a ``very weak solution" of Problem \eqref{elldifstaf}, \eqref{elldifstab}, with only
$f \in L^1(\O)$, by working with some discrete equivalent of the $W^{1,q}_0$-norm, with $q < \frac d {d-1}$. 
This is not detailed here.

\item 
In Theorem \ref{ellcvgce}, it is also possible to prove convergence results
when $f(x)$ (resp. $\vbf(x) $) 
is replaced by some nonlinear function $f(x,u(x))$, 
(resp. $\vbf(x,u(x))$ under adequate assumptions,   see \cite{cvnl}.
\end{enumerate}
\end{remark}

\proof of Theorem \ref{ellcvgce}

\smallskip

Let $Y$ be the set of approximate solutions, that is the set of $u_\mesh$ 
where $\mesh$ is an admissible mesh in the sense of Definition 
\refp{meshdirichlet}. First, we want to prove that $u_\mesh$ tends to the 
unique solution (in $H^1_0(\O)$) to \eqref{ellsolvar} as $\sizemesh \tends 0$.

\smallskip

Thanks to  Lemma \ref{ellexistu}  and to the discrete  Poincar\'e inequality
\eqref{ellinpoin}, there exists $C_1 \in \R$, only depending on $\O$ and 
$f$, such that $\norm{u_{\Mesh}}{1,\mesh} \le C_1$ and $\Vert u_{\Mesh} 
\Vert_{L^2(\O)} \le C_1$ for all $u_\mesh \in Y$. Then, thanks to Lemma 
\ref{h1dtot} and to the
compactness result given in Theorem
\refp{kolmh10}, the set $Y$ is relatively compact in $L^2(\O)$ 
 \index{Compactness!-- of the approximate solutions in $L^2$} and any 
possible limit (in $L^2(\O)$) of a sequence $(u_{\mesh_n})_{\nnn} \subset 
Y$ (such that $\size(\mesh_n) \tends 0$) belongs to $H^1_0(\O)$.
Therefore, thanks to the uniqueness of the solution (in $H^1_0(\O)$) of 
(\ref{ellsolvar}), it is sufficient to prove that if $(u_{\mesh_n})_{\nnn} \subset 
Y$ converges towards some $u \in H^1_0(\O)$, in $L^2(\O)$, and 
$\size(\mesh_n) \tends 0$  (as $n \tends \infty$), 
then $u$ is the solution to (\ref{ellsolvar}). We prove this result below, 
omiting the index $n$, that is assuming $u_\mesh \tends u$ in $L^2(\O)$ 
as $\sizemesh \tends 0$.

\separ

Let $\psi\in C^\infty_c(\O)$ and let $\size(\mesh)$ be small enough so 
that $\psi(x)=0$ if $x \in \cv$ and $\cv \in \mesh$ is such that $\dr \cv 
\cap \dr \O \ne \emptyset$. Multiplying  (\ref{ellschema}) by 
$\psi(\centercv)$, and summing the result over
$\cv\in\mesh$  yields

\be
T_{1} + T_{2} + T_{3} = T_{4}, 
\label{ell1g}
\ee

with

\[
T_{1} = b \summesh  \mcv u_{\cv} \psi(\centercv), 
\]
\[
T_{2} = - \summesh \sumv\transcvcvv (u_{\cvv} - u_{\cv})
\psi(\centercv),  \]

\[
T_{3} = \summesh \sumedgecv \vcvedge \uedgep \psi(\centercv),
\]

\[
T_{4} = \summesh \mcv \psi(\centercv) f_{\cv}.
\]

First remark that,  since  $u_{\mesh}$ tends to $u$ in $L^2(\O)$,   
\[
T_{1} \tends  b\dsp{\int_\O} u(x) \psi(x) dx \text{ as 
}\size(\Mesh)\tends 0.
\] 
Similarly, 
\[
T_{4} \tends \dsp{\int_\O} f(x) \psi(x) dx  \text{ as 
}\size(\Mesh)\tends 0.
\] 

\separ

Let us now turn  to the study of $T_{2}$;

\[
T_{2} =  -   \sumEll   \transcvcvv  
(u_\cvv - u_\cv) (\psi(\centercv) - \psi(\centercvv)).
\]
Consider the following auxiliary expression:

\[
\ban
T'_{2}  & = \dsp{\int_\O} u_{\mesh}(x) \Delta \psi(x) dx \\
& = \summesh u_{\cv} \dsp{\int_\cv} \Delta \psi(x) dx \\
& =   \sumEll (u_{\cv} - u_{\cvv})
\dsp{\int_{\edgecvcvv}} \grad\psi(x) \cdot \bfn_{\cv, \cvv} \dfrontiere(x).
\ean 
\]

Since $u_{\mesh}$ converges to $u$ in $L^2(\O)$, it is clear that $T'_{2}$ tends
to $\dsp{\int_\O} u(x) \Delta \psi(x)\ dx $ as $\size(\Mesh)$ tends to $0$. 

Define

\[
\left.\begin{array}{ll}
R_{\cv ,\cvv} = \dsp{\frac 1 {\medgecvcvv} \int_{\edgecvcvv} }
\grad\psi(x) \cdot \bfn_{\cv,\cvv} \dfrontiere(x)
- \frac {\psi(\centercvv) - \psi(\centercv)} {\dedgecvcvv},
\end{array}\right.
\]

where $\bfn_{\cv,\cvv}$ denotes the unit normal vector to $\edgecvcvv$, outward 
to $\cv$, then 

\[
\left.\begin{array}{ll}
\vert T_{2}+T'_{2} \vert & = 
\vert  \sumEll  \medgecvcvv  (u_\cv - u_\cvv) R_{\cv ,\cvv} \vert \\
& \le 
 \Bigl[ \sumEll  \medgecvcvv \frac {(u_\cv - u_\cvv)^2} {\dedgecvcvv}
\sumEll \medgecvcvv \dedgecvcvv (R_{\cv ,\cvv})^2 \Bigr]^{1/2}, 
\end{array}\right.
\]

Regularity properties of the function $\psi$ give
the existence of $ C_2 \in \R$,  only depending  on $\psi$, 
such that $\vert R_{\cv ,\cvv} \vert \le C_2  \size(\Mesh)$.
Therefore, since 
\[\sumEll \medgecvcvv \dedgecvcvv  \le d \meas(\O),\]
from Estimate (\ref{ellestimx}), we conclude that $ T_{2}+T'_{2} \tends 
0$ as $\size(\Mesh) \tends 0$.

Let us now show that $T_3$ tends to $-\int_\O \vbf(x) u(x) \grad \psi (x) dx$ as
$\size(\Mesh)\tends 0$. Let us decompose  $T_3 =  T'_3 + T''_3$  where

\[
T'_{3} = \summesh \sumedgecv \vcvedge (\uedgep - u_\cv) \psi(\centercv)
\]
and
\[
T''_{3} = \summesh \sumedgecv \vcvedge  u_\cv \psi(\centercv) = \int_\O \div \vbf(x)
u_\mesh(x)   \psi_\mesh(x) dx,
\]
where $\psi_\mesh$ is defined by $\psi_\mesh(x)=\psi(\centercv)$ if $x 
\in \cv$, $\cv \in \mesh$. 
Since $u_\mesh \tends u$ and $\psi_\mesh \tends \psi$ in $L^2(\O)$ as 
$\sizemesh \tends 0$ (indeed, $\psi_\mesh \tends \psi$ uniformly on $\O$ 
as $\sizemesh \tends 0$) and since $\div \vbf \in L^\infty(\O)$, one has 
\[
T''_{3} \tends \int_\O \div \vbf(x) u(x)   \psi(x) dx \text{ 
as }\size(\Mesh) \tends 0.
\]

\smallskip

Let us now rewrite $T'_{3}$ as $T'_{3}=T'''_{3} +r_3$ with
\[
T'''_{3} = \summesh \sumedgecv (\uedgep - u_\cv) \int_{\edge} \vbf(x) \cdot 
\ncvedge \psi(x) \dfrontiere(x)
\]
and
\[
r_3= \summesh \sumedgecv (\uedgep - u_\cv) \int_{\edge} \vbf(x) \cdot 
\ncvedge (\psi(\centercv)-\psi(x)) \dfrontiere(x).
\]

Thanks to the regularity of $\vbf$ and $\psi$, there exists $C_3$  only depending  
 on $\vbf$ and $\psi$ such that

\[
\vert r_3 \vert \le C_3 \sizemesh \sumEll \vert u_\cv - u_\cvv \vert \medgecvcvv,
\]

which yields, with    the Cauchy-Schwarz  inequality,

\[
\vert r_3 \vert \le C_3 \sizemesh (\sumEll \transcvcvv \vert u_\cv - u_\cvv 
\vert^2)^{\frac 1 2} (\sumEll \medgecvcvv \dcvcvv)^{\frac 1 2},
\]

from which one deduces, with Estimate \eqref{ellestimx}, that $r_3 \tends 
0$ as $\sizemesh \tends 0$.

Next, remark that

\[
T'''_{3} = - \summesh u_\cv \sumedgecv \int_{\edge} \vbf(x) \cdot 
\ncvedge \psi(x) \dfrontiere(x) = - \summesh u_\cv \int_\cv \div 
(\vbf(x)\psi(x))dx.
\]

This implies (since $u_\mesh \tends u$ in $L^2(\O)$) that $T'''_{3} 
\tends -\int_\O \div (\vbf(x)\psi(x)) u(x) dx$, so that $T'_3$ has the same 
limit and $T_3 \tends - \int_\O \vbf(x) \cdot \grad  \psi (x) u(x) dx$.

\medskip 

Hence, letting $\size(\Mesh)\tends 0$ in \eqref{ell1g}
yields that the function $u \in H^1_0(\O)$ satisfies

\[
\int_\O \Bigl(b u(x) \psi(x)  - u(x) \Delta \psi(x)   - \vbf(x) u(x) \grad \psi
(x) - f(x) \psi(x)\Bigr) dx = 0, \; \forall \psi\in C^\infty_c(\O),
\]

which, in turn, yields  \eqref{ellsolvar} thanks to the fact  that  $u \in
H^1_0(\O)$, and to the density of  $C^\infty_c(\O)$ in $H^1_0(\O)$.

This concludes the proof of $u_\mesh \tends u$ in $L^2(\O)$ as $\sizemesh 
\tends 0$, where $u$ is the unique solution (in $H^1_0(\O)$) to \eqref{ellsolvar}.

\smallskip

S Let us now prove that $\norm{u_\mesh}{1,\mesh}$ tends to 
$\norm{u}{H^1_0(\O)}$
in the pure diffusion case, i.e. assuming $b = 0$ and $\vbf= 0$. 
Since 

\[
\norm{u_\mesh}{1,\mesh}^2=\int_\O f_\mesh(x) u_\mesh(x) dx 
\tends \int_\O f(x) u(x)dx \text{ as } \sizemesh 
\tends 0,
\]
where $f_\mesh$ is defined from $\O$ to $\R$ by $f_\mesh(x)=f_\cv$ 
a.e. on $\cv$ for all $\cv \in \mesh$, it is easily seen that
\[
\norm{u_\mesh}{1,\mesh}^2 \tends \int_\O f(x) u(x)dx= 
\norm{u}{H^1_0(\O)}^2 \text{ as } \sizemesh 
\tends 0.
\] 

This concludes the proof of Theorem \ref{ellcvgce}.
\eof

\begin{remark}[Consistency for the adjoint operator]\label{creme}  \normalfont
The proof of Theorem \ref{ellcvgce} uses the property of consistency  of 
the (diffusion) fluxes on 
the test functions. This property consists in writing the consistency of 
the fluxes for the adjoint operator to the discretized Dirichlet 
operator. This consistency is achieved thanks to that of fluxes for the 
discretized Dirichlet operator and to the fact that this operator is self 
adjoint. In fact, any discretization of the Dirichlet operator giving 
``$L^2$-stability" and 
consistency of fluxes on its adjoint, yields a convergence result (see also 
Remark \refp{cafe}). On the contrary, the error estimates proved in 
sections \ref{errordirichlet} and  \ref{sh2ee} directly use the 
consistency for the discretized 
Dirichlet operator itself.
\end{remark}

\begin{remark}[Finite volume schemes and $H^1$ approximate solutions] 
\label{Liebig} \normalfont~

In the above proof, we showed that a sequence of approximate solutions 
(which are piecewise constant functions)
converges in  $L^2(\O)$  to a limit which is in $H^1_0(\O)$. An alternative 
to the use of Theorem \ref{kolmh10} is the construction of 
a bounded sequence in $H^1(\R^d)$  from the sequence of 
approximate solutions. This can be performed by convoluting the 
approximate solution  with a
mollifier ``of size $\size(\mesh)$". Using  Rellich's compactness theorem 
and the weak sequential  compactness of the bounded sets of $H^1$, one 
obtains that the limit of the sequence of approximate solutions is in $H^1_0$. 
\end{remark}

Let us now deal with the case of non homogeneous Dirichlet boundary 
conditions, in which case  $g\in H^{1/2}(\dr\O)$  is no longer assumed to 
be $0$. The proof uses the following preliminary result:

\begin{lemma}\label{thuh1}
Let $\O$ be an open bounded polygonal subset of $\R^2$, $\tilde g\in 
H^1(\O)$ and $g = \gammat 
(\tilde g)$ (recall that $\gammat$ is the ``trace" operator from 
$H^1(\O)$ to $H^{1/2}(\dr \O)$). Let $\Mesh$ be  an admissible mesh (in the sense of
Definition \refp{meshdirichlet}) such that, for some $\zeta>0$, the inequality   $\dcvedge \ge
\zeta\diam(\cv)$  holds for all control volumes $\cv\in\mesh$ and for all $\edge \in
\edgescv$, and let  $M\in\N$ be such that $\card (\edgescv) \le M$ for all $\cv\in\mesh$. Let us define $\tilde 
g_\cv$ for all $\cv \in \mesh$ by

\[
\tilde g_\cv  = \frac 1
{\mcv} \int_\cv \tilde g(x) dx 
\]

and $\tilde 
g_\edge$ for all $\edge \in \edgesext$ by

\[
\tilde g_{\edge}  = \frac 1 {\medge} \int_\edge g(x) \dfrontiere(x). 
\]

Let us define
\be
\mathcal{N}(\tilde g, \mesh) = \Bigl(\sumedgesintkl \transcvcvv (\tilde g_\cv -  \tilde g_\cvv)^2  
+  \sumedgesext  \transedge (\tilde g_{\cv(\edge)} -  \tilde 
g_\edge)^2\Bigr)^{\half}, \label{Ngtilde}
\ee

where $\cv(\edge)=\cv$ if $\edge \in 
\edgesext \cap \edgescv$.
Then there exists $C \in \Rp$, only depending  on $\zeta$ and $M$, such that

\be
\mathcal{N}(\tilde g, \mesh) \le C 
\Vert \tilde g \Vert_{H^1(\O)}.
\label{ellthuh1}
\ee
\end{lemma}

\proof of Lemma \ref{thuh1}

\smallskip

Lemma \ref{thuh1} is given in the two dimensional case, an analogous 
result is possible in the three dimensional case. 
Let $\O$, $\tilde g$, $\Mesh$, $\zeta$, $M$ satisfying the hypotheses
of Lemma \ref{thuh1}. By a
classical argument of density, one may assume that
$\tilde g\in C^1(\Obar,\R)$.

A first step consists in proving that there exists $C_1\in\Rp$, only depending  
on $\zeta$, such that

\be
(\tilde g_\cv  - \tilde g_\edge)^2 \le C_1 \frac {\diam(\cv)} {\medge} \int_\cv
\vert \grad \tilde g(x) \vert^2 dx, \forall \cv \in \mesh, \forall 
\edge \in \edgescv,
\label{thuh11}
\ee
where  $\tilde g_\cv$  (resp. $\tilde g_\edge$) is the mean value of 
$\tilde g$ on $\cv$ (resp. $\edge$), for $\cv \in \mesh$ (resp. $\edge 
\in \edges$).
Indeed, without loss of generality, one assumes that
$\edge =\{0\}\times J_0$, with $J_0$ is a closed interval of $\R$ and $\cv
\subset \R_+\times\R$.

\smallskip

Let $\a =\max\{ x_1, x=(x_1,x_2)^t \in \cvbar\}$ and
$a=(\a,\b)^t \in \cvbar$. In the following, $a$ is fixed. For all
$x_1\in(0,\a)$, let
$J(x_1) = \{ x_2\in\R,$ such that $(x_1,x_2)^t\in\cvbar \}$, so that $J_0 = J(0)$.

For a.e. $x=(x_1,x_2)^t \in \cv$ and a.e., for the 1-Lebesgue measure,
$y=(0,\ybar)^t \in \edge$ (with $\ybar \in J_0$), one sets
$z(x,y)=ta+(1-t)y$ with $t=\frac {x_1} {\a}$. Note that, since $\cvbar$ is 
convex, $z(x,y) \in \cvbar$ and $z(x,y)=(x_1,z_2(x_1,\ybar))^t$, with 
$z_2(x_1, \ybar) =  \frac {x_1} {\a} \beta + (1 - \frac {x_1} {\a})\ybar$.

One has, using  the Cauchy-Schwarz  inequality,

\be
(\tilde g_\cv  - \tilde g_\edge)^2 \le \frac 2 {\mcv\medge} (A + B),
\label{thuh111}\ee

where

\[
A= \int_\cv \int_\edge \big( \tilde g(x) -  \tilde g(z(x,y)) \big)^2
\dfrontiere(y)dx,
\]

and

\[
B= \int_\cv \int_\edge \big( \tilde g(z(x,y)) - \tilde g(y) \big)^2
\dfrontiere(y)dx.
\]

Let us now obtain a bound of $A$. Let $D_i \tilde g$, $i=1$ or $2$, denote
the partial derivative of $\tilde g$ w.r.t. the components of 
$x=(x_1,x_2)^t \in \R^2$. Then,

\[
A= \int_{0}^{\a} \int_{J(x_1)} \int_{J(0)} \big( \int_{z_2(x_1,\ybar)}^{x_2}
D_2 \tilde g(x_1,s) ds \big)^2
d\ybar dx_2 dx_1.
\]

The Cauchy-Schwarz  inequality yields

\[
A \le \diam(\cv) \int_{0}^{\a} \int_{J(x_1)} \int_{J(0)} \int_{J(x_1)} \big(
D_2 \tilde g(x_1,s)\big)^2
ds d\ybar dx_2 dx_1
\]

and therefore

\be
A \le \diam(\cv)^3 \int_{\cv} \big( D_2 \tilde g(x)\big)^2  dx.
\label{thuh12}\ee

One now turns to the study of $B$, which can be rewritten as

\[
B= \int_{0}^{\a} \int_{J(x_1)} \int_{J(0)}
\big( \int_{0}^{x_1} [ D_1 \tilde g(s,z_2(s,\ybar)) + \frac {\b - \ybar} {\a}
D_2 \tilde g(s,z_2(s,\ybar))] ds \big)^2
d\ybar dx_2 dx_1.
\]

The Cauchy-Schwarz  inequality and the fact that $\a\ge\zeta \diam(\cv)$ 
give that

\be
B \le 2 \diam(\cv) ( B_1 + \frac 1 {\zeta^2} B_2),
\label{thuh112}\ee

with

\[
B_i = \int_{0}^{\a} \int_{J(x_1)} \int_{J(0)}
\int_{0}^{x_1} \big( D_i \tilde g(s,z_2(s,\ybar)) \big)^2 ds d\ybar dx_2 
dx_1, \, i=1, \, 2.
\]

First, using Fubini's theorem, one has

\[
B_i =   \int_{J(0)} \int_{0}^{\a}
\big( D_i \tilde g(s,z_2(s,\ybar)) \big)^2 \int_{s}^{\a} \int_{J(x_1)} dx_2
dx_1 ds d\ybar.
\]

Therefore

\[
B_i \le \diam(\cv) \int_{0}^{\a} \int_{J(0)}
\big( D_i \tilde g(s,z_2(s,\ybar)) \big)^2 (\a - s) d\ybar ds.
\]

Then, with the change of variables $z_2 = z_2(s,\ybar)$, one gets

\[
B_i \le \diam(\cv) \int_{0}^{\a} \int_{J(s)}
\big( D_i \tilde g(s,z_2) \big)^2 \frac {\a - s} {1 - \frac {s} {\a} } dz_2 ds.
\]

Hence
\be
B_i \le \diam(\cv)^2 \int_{\cv} \big( D_i \tilde g(x) \big)^2 dx.
\label{thuh13}\ee

Using the fact that $\mcv \ge \pi \zeta^2 \big( \diam(\cv)\big)^2$,
\eqref{thuh111}, \eqref{thuh12},
\eqref{thuh112} and \eqref{thuh13}, one concludes \eqref{thuh11}.

\medskip

In order to conclude the proof of \eqref{ellthuh1}, one remarks that

\[
\Bigl( \mathcal{N}(\tilde g, \Mesh) \Bigr)^2 \le
2 \summesh \sumedgecv \transedge (\tilde g_\cv  - \tilde g_\edge)^2.
\]

Because, for all $\cv\in\mesh$ and $\edge\in\edgescv$, $\dedge\ge
\zeta\diam(\cv)$, one gets 
thanks to \eqref{thuh11}, that

\[
\Bigl( \mathcal{N}(\tilde g, \Mesh) \Bigr)^2 \le
2 \summesh \sumedgecv \frac {C_ 1} \zeta \int_\cv \vert \grad \tilde 
g(x) \vert^2 dx.
\]

The above inequality shows that

\[
\Bigl( \mathcal{N}(\tilde g, \Mesh) \Bigr)^2 \le
2 M  \frac {C_ 1} \zeta \int_\O \vert \grad \tilde g(x) \vert^2 dx,
\]

which implies \eqref{ellthuh1}.\eof

\begin{theo}[Convergence, non homogeneous  Dirichlet boundary
condition] \label{ellcvgcediraf} ~
\index{Convergence of the finite volume method!-- for an elliptic equation}

Assume items 1, 2, 3 and 4 of Assumption \refp{ellH} and $g\in H^{1/2}(\dr\O)$.
Let $\zeta\in\Rp$ and $M\in\N$ be given values.
Let $\Mesh$ be an admissible mesh
(in the sense of Definition \refp{meshdirichlet}) such that $\dcvedge \ge
\zeta\diam(\cv)$ for all control volumes $\cv\in\mesh$ and for all $\edge \in
\edgescv$,
and $\card(\edgescv) \le M$ for all $\cv\in\mesh$.
Let $(u_\cv)_{\cv \in \mesh}$ be the solution of the system given by
equations \eqref{ellschema}-\eqref{ellschemac} and

\be
u_\edge  =  \frac 1 \medge \int_\edge g(x) \dfrontiere(x),
\ \forall \edge \in \edgesext.
\label{ellschemadi}
\ee

(note that the proofs of existence and uniqueness of
$(u_\cv)_{\cv \in \mesh}$ which were given in Lemma \refp{ellexistu} 
remain valid).
Define
$u_{\Mesh} \in X(\Mesh)$ by $u_\Mesh(x) = u_{\cv}$ for
a.e. $x \in \cv$ and for any
$\cv \in \mesh$. Then, $u_{\Mesh}$
converges, in $L^2(\O)$, to the unique variational solution $u\in H^1(\O)$ of
Problem \eqref{elldifstaf}, \eqref{elldifstab} as $\size(\Mesh)\tends
0$.

\end{theo}

\proof of Theorem \ref{ellcvgcediraf}

The proof is only detailed for the case $b=0$ and $\vbf=0$ (the extension of
the proof
to the general case is straightforward using the proof of Theorem
\refp{ellcvgce}).
Let $\tilde g\in H^1(\O)$ be such that the trace of $\tilde g$ on $\dr\O$
is equal to $g$. One defines $\util_\Mesh\in X(\Mesh)$ by
$\util_\Mesh = u_\mesh - {\tilde g}_\mesh$ where  ${\tilde g}_\Mesh\in X(\Mesh)$
is defined by $\tilde g(x) = \frac 1 {\mcv} \int_\cv \tilde g(y)
dy$ for all
$x\in\cv$ and all $\cv\in\mesh$. Then $(\util_\cv)_{\cv\in\mesh}$ satisfies

\be
\sum_{\edge \in \edgescv} \tilde F_{\cv,\edge}     =
\mcv f_\cv - \sum_{\edge \in \edgescv} G_{\cv,\edge}, \ \forall
\cv\in\mesh,
\label{ellschemada}
\ee

\be
\tilde F_{\cv,\edge}=  - \transcvcvv  (\util_\cvv - \util_\cv), \, \forall
\edge \in \edgesint, \; \mathrm{if} \;
\edge = \edgecvcvv,
\label{ellschembda}
\ee

\be
\tilde F_{\cv,\edge}   = \transedge (\util_{\cv}), \, \forall
\edge \in \edgesext \text{ such that }\edge \in \edgescv.
\label{ellschemcda}
\ee

\be
G_{\cv,\edge}=  - \transcvcvv
({\tilde g}_\cvv   - {\tilde g}_\cvv), \,
\forall  \edge \in \edgesint, \; \mathrm{if} \;
\edge = \edgecvcvv,
\label{ellschemdda}
\ee

\be
G_{\cv,\edge}   = -
\transedge ({\tilde g}_\edge- {\tilde g}_\cvv), \, \forall
\edge \in \edgesext \text{ such that }\edge \in \edgescv,
\label{ellschemeda}
\ee
where ${\tilde g}_\edge =\frac 1 \medge \int_\edge g(x) \dfrontiere(x)$
Multiplying \eqref{ellschemada} by $\util_\cv$, summing over $\cv\in\mesh$,
gathering by edges in the right hand side and using   the Cauchy-Schwarz  
inequality yields

\[
\Vert \util_\Mesh \Vert_{1,\Mesh}^2 \le \summesh \mcv f_\cv \util_\cv +
 \mathcal{N}(\tilde g, \Mesh)  \Vert \util_\Mesh \Vert_{1,\Mesh},
\]

from the definition \refep{Ngtilde} of $ \mathcal{N}(\tilde g, \Mesh)$ and 
Definition \refp{defh10d} of $\Vert \cdot \Vert_{1,\Mesh}$.
Therefore, thanks to   Lemma
\refp{thuh1}  and the discrete Poincar\'e inequality
\eqref{ellinpoin}, there exists $C_1 \in \R$, only depending  on $\O$, $\Vert
\tilde g \Vert_{H^1(\O)}$,
$\zeta$, $M$ and $f$, such
that $\Vert \util_\Mesh \Vert_{1,\mesh} \le C_1$ and $\Vert \util_\Mesh 
\Vert_{L^2(\O)} \le C_1$. 
Let us now prove that $\util_\Mesh$ converges in $L^2(\O)$, as $\sizemesh 
\tends 0$, towards the unique solution in $H^1_0(\O)$ to 
\eqref{ellsolvar}. We proceed as in Theorem \refp{ellcvgce}. Using Lemma 
\ref{h1dtot}, the
compactness result given in Theorem
\refp{kolmh10} and the uniqueness of the solution (in $H^1_0(\O)$) of 
(\ref{ellsolvar}), it is sufficient to prove that if $\util_{\mesh}$ 
converges towards some $\util \in H^1_0(\O)$, in $L^2(\O)$ 
as $\size(\mesh) \tends 0$, 
then $\util$ is the solution to (\ref{ellsolvar}). 
In order to prove this result, let us introduce the function $\tilde g_\Mesh$ defined by

\[
\tilde g_\Mesh(x) = \frac 1 {\mcv} \int_\cv \tilde g(y) dy, \, \forall x 
\in \cv, \, \forall \cv \in \mesh,
\]

which converges to $\tilde g$ in $L^2(\O)$, as
$\size(\mesh) \tends 0$. Then the function $u_\Mesh$ converges in $L^2(\O)$, as
$\size(\mesh) \tends 0$ to $u = \util + \tilde g \in H^1(\O)$ and the proof 
that $\util$
is the unique solution of (\ref{ellsolvar}) is identical to the corresponding 
part in the proof of Theorem \refp{ellcvgce}. This completes the proof of 
Theorem \ref{ellcvgcediraf}. \eof

\begin{remark}[Lipschitz continuous boundary data]
\normalfont A  simpler proof of convergence for the finite volume scheme with non homogeneous Dirichlet boundary condition is possible if $g$ is the trace of a Lipschitz-continuous function $\tilde g$. 
In thiscase, $\zeta$ and $M$ do not have to be introduced and Lemma \ref{thuh1} is not used.
The scheme is defined with $u_\edge = g(y_\edge)$ instead of the average value of $g$ on $\edge$, and the
proof uses $\tilde g(\centercv)$ instead of the average value of $\tilde g$ on $\cv$.
\end{remark}

\subsection{$C^2$ error estimate} \label{errordirichlet}

Under adequate regularity assumptions on the solution of Problem
\eqref{elldifstaf}-\eqref{elldifstab}, one may prove that the error between the
exact solution and the approximate solution given by the finite volume scheme 
\eqref{ellschema}-\eqref{ellschemad} is of order $\sizemesh = \sup_{\cv \in
\mesh}\diam(\cv)$, in a certain sense which we give in the following theorem:

\begin{theo}\label{ellesterr}\index{Error estimate!-- for an elliptic  equation!-- in the general case}
Under Assumption \refp{ellH}, let $\Mesh$ be an admissible mesh as defined in
Definition \refp{meshdirichlet}
and $u_\mesh \in X(\Mesh)$ (see 
Definition \refp{Xmesh})  be defined a.e.in $\O$ by $u_\mesh (x)  = u_\cv$ 
for a.e.
$x \in \cv$, for all $\cv \in \mesh$, where $(u_\cv)_{\cv \in \mesh}$ is the
solution to \eqref{ellschema}-\eqref{ellschemad}. 
 Assume that the unique variational solution $u$ of
Problem \eqref{elldifstaf}-\eqref{elldifstab} satisfies $u\in
C^2({\overline \O})$. Let, for each $\cv \in \mesh$,  $e_\cv = u(\centercv) -
u_\cv$, and $e_\Mesh \in X(\Mesh)$ defined by $e_\Mesh (x) = e_\cv$ for a.e.
$x \in \cv$, for all $\cv \in \mesh$. 

Then, there exists  $C >0$  only depending  on $u$, $\vbf$ and $\O$ 
such that

\be
\Vert e_\Mesh \Vert_{1,\Mesh} \le C  \size(\Mesh),
\label{eqesterrh1}
\ee

where  $\Vert \cdot \Vert_{1,\Mesh} $ is the discrete $\Hunz$ norm defined in
Definition \ref{defh10d},

\be
\Vert  e_\Mesh  \Vert_{L^2(\O)} \le C \size(\Mesh)
\label{eqesterr}\ee

and

\be
\bac
\dsp{\sum_{\substack{\edge \in \edgesint \\ \edge = \edgecvcvv} }   \medge    \dedge 
\Bigl( \frac {u_\cvv  - u_\cv} \dedge  - \frac 1 \medge 
  \int_{\edge} \grad u(x) \cdot \ncvedge \dfrontiere(x)\Big)^2  + }
\\
\dsp{\sum_{\substack{\edge \in \edgesext \\ \edge \in \cvbar \cap \dr \O} }   \medge   
\dedge  \Bigl( \frac {g(\centeredge) - u_\cv} \dedge  - \frac 1 \medge 
  \int_{\edge} \grad u(x) \cdot \ncvedge \dfrontiere(x)\Big)^2\le C \size(\Mesh)^2.}
\ean
\label{eqesterr1}
\ee

\end{theo}

\begin{remark} \normalfont ~

\begin{enumerate}

\item Inequality \eqref{eqesterrh1} (resp. \eqref{eqesterr}) yields  an 
estimate of order 1 for the discrete $H^1_0$ norm (resp. $L^2$ norm) of the error
on the solution. Note also that, since $u \in C^1(\Obar)$, one deduces, from
\eqref{eqesterr}, the existence of $C$ only depending  on $u$ and $\O$ such that
$\norm{u-u_\mesh}{L^2(\O)} \le C\sizemesh$. Inequality 
\eqref{eqesterr1}  may be seen as an estimate of order 1 for the $L^2$ norm of
the flux. 

\item In  \citet{BMO}, finite element \index{Finite element method} tools are used to obtain error
estimates of order $\sizemesh^2$ in the case  $d=2$, 
$\vbf=b=g=0$ and if 
the elements of $\mesh$ are triangles of a finite element mesh satisfying 
the Delaunay condition \index{Delaunay condition} (see section \ref{scvfe} page \pageref{Delaunay}). 
Note that this result is quite 
different of those of the remarks \refp{h2conv} and \refp{h2conv2d}, which are obtained by 
using a higher order approximation of the flux.

\item The proof of Theorem \ref{ellesterr} given below is close to that of error
estimates for finite element schemes in the sense that it uses the 
coerciveness of
the operator (the discrete Poincar\'e inequality) instead of the discrete
maximum principle of Proposition \refp{maxppe} 
(which is used for error estimates with finite difference
schemes). 
\end{enumerate}
\end{remark}

\proof  of Theorem \ref{ellesterr}

\smallskip

Let $u_\mesh \in X(\Mesh)$ be defined a.e. in $\O$ by $u_\mesh (x)  = u_\cv$ for
a.e. $x \in \cv$, for all $\cv \in \mesh$, where $(u_\cv)_{\cv \in \mesh}$ is the
solution to \eqref{ellschema}-\eqref{ellschemad}. Let us
write the flux balance for any $\cv\in\mesh$;

\be
\sumedgecv  \Bigl(\Fbarcvedge + \Vbarcvedge \Bigr) + b
\int_\cv u(x) dx= \int_\cv f(x)
dx, \label{eller1}
\ee
where
$\Fbarcvedge = - \int_\edge \grad u(x) \cdot \ncvedge \dfrontiere(x)$,  and
$\Vbarcvedge =   \int_\edge u(x)  \vbf(x) \cdot \ncvedge  \dfrontiere(x)$  are
respectively the diffusion and convection fluxes through $\edge$ outward to $\cv$. 

Let
$\Fastcvedge$ and $\Vastcvedge$ be defined by
 
\[
\Fastcvedge = - \transcvcvv
(u(\centercvv)-u(\centercv)), \, \forall  \edge = \edgecvcvv \in \edgescv \cap
\edgesint, \, \forall \cv \in \mesh,
\]

\[
\Fastcvedge d(\centercv, \edge) =
-  \medge  (u(\centeredge)-u(\centercv)) , \, \forall  \edge \in \edgescv \cap
\edgesext,
\,
\forall \cv \in \mesh,
\]

\[
\Vastcvedge   = \vcvedge u(x_{\edge,+}), \forall  \edge \in \edgescv, \,
\forall \cv \in \mesh,
\]
where $x_{\edge,+} = \centercv$ (resp. $\centercvv$) if $\edge \in \edgesint$, 
$\edge = \edgecvcvv$ and $\vcvedge \ge 0$ (resp. $\vcvedge \le 0$) and   
$x_{\edge,+} = \centercv$ (resp. $\centeredge$) if $\edge = \edgescv \cap 
\edgesext$ and $\vcvedge \ge 0$ (resp. $\vcvedge \le 0$). Then, the 
consistency \index{Consistency!-- error} error
on the diffusion and convection fluxes may be defined as 

\be
\Rcvedge  = \dsp{\frac 1 {\medge}} (\Fbarcvedge - \Fastcvedge), \label{eller2}
\ee
\be
\rcvedge  = \dsp{\frac 1 {\medge}} (\Vbarcvedge - \Vastcvedge),\label{eller2b}
\ee

Thanks to the regularity of $u$ and $\vbf$, there exists $C_1 \in \R$,  only 
depending 
on $u$ and $\vbf$, such that $\vert \Rcvedge \vert  + \vert \rcvedge \vert  \le C_1
\size(\Mesh)$  for any $\cv \in \mesh$ and $\edge \in \edgescv$. 
For $\cv \in \mesh$, let 
\[\rho_\cv=u(\centercv) - (1/\mcv) \int_\cv 
u(x)dx,
\]
so that   $\abs {\rho_\cv} \le C_2\size(\Mesh)$ with some $C_2 \in 
\R_+$  only depending  on $u$.

\medskip

Substract 
\eqref{ellschema} to  \eqref{eller1}; thanks to  \eqref{eller2} and
  \eqref{eller2b}, one has

\be
 \sumedgecv \Bigl(\Gcvedge + \Wcvedge \Bigr) + b \mcv e_\cv =  b \mcv
\rho_\cv  - \sumedgecv \medge   (\Rcvedge + \rcvedge), \label{eller3}
\ee

where

\smallskip

   $\Gcvedge = \Fastcvedge - \Fcvedge$ is such that 
\[
\Gcvedge = -\transcvcvv (e_\cvv -e_\cv), \, \forall \cv \in \mesh, \, 
\forall  \edge \in \edgescv \cap
\edgesint, \edge=\edgecvcvv,
\]  

\[
\Gcvedge d(\centercv,\edge) = \medge e_\cv, \, \forall \cv \in \mesh, \, \forall  \edge \in \edgescv \cap \edgesext,
\]
with $e_\cv = u(\centercv) - u_\cv$, and $\Wcvedge = \Vastcvedge - \Vcvedge = \vcvedge (u(x_{\edge,+}) - u_{\edge,+})$

\smallskip

Multiply  \eqref{eller3} by $e_\cv$, sum  for $\cv\in\mesh$, and note
that
\[  
 \summesh \sumedgecv  \Gcvedge e_\cv = \sumedges \vert D_\edge e \vert^2 \frac 
\medge \dedge = \Vert e \Vert_{1, \Mesh}^2.
\]
Hence
\be
 \Vert e_\Mesh \Vert_{1, \Mesh}^2  + \summesh  \sumedgecv  \vcvedge  e_{\edge,+}e_\cv +  b  \Vert e_\Mesh \Vert_{L^2(\O)}^2 \le b \summesh \mcv \rho_\cv e_\cv -  \summesh \sumedgecv \medge (\Rcvedge + \rcvedge) e_\cv,
\label{ellraslebol1}
\ee
where 
\begin{align*}
  &e_\Mesh \in X(\Mesh), e_\Mesh (x) = e_\cv \mbox{ for a.e. }x \in \cv\mbox{ and for all } \cv \in \mesh,  \\
 & \vert D_\edge e \vert=\begin{cases}\vert e_{\cv}-e_{\cvv} \vert \mbox{ if }\edge \in \edgesint, \edge=\edgecvcvv, \\
                         \vert e_{\cv}\vert\vert \mbox{ if }\edge\in \edgescv \cap \edgesext,
                        \end{cases}
\\
&  e_{\edge,+}  = u(x_{\edge,+}) - u_{\edge,+}.
\end{align*}

By Young's inequality, the first term of the left hand side satisfies:
\be
 \vert \summesh \mcv \rho_\cv e_\cv \vert \le \half \Vert e_\Mesh 
 \Vert^2_{L^2(\O)} + \half C_2^2 (\size(\Mesh))^2 \meas(\O).
\label{ellyoung}
\ee

\smallskip

Thanks to the assumption $\div \vbf \ge 0$, one obtains, through a computation similar to \eqref{vfqconv1}-\refep{vfqconv11} that
\[
\summesh  \sumedgecv  \vcvedge e_{\edge,+}  e_\cv \ge 0.
\]  
Hence,  \eqref{ellraslebol1} and \eqref{ellyoung} yield that there exists $C_3$  only depending   on $u, b$ and $\O$ such that
\be
\Vert e_\Mesh \Vert_{1, \Mesh}^2 + \half b \Vert e_\Mesh \Vert_{L^2(\O)}^2
\le  C_3 (\size(\Mesh))^2
 -  \summesh \sumedgecv \medge (\Rcvedge +
\rcvedge) e_\cv, 
\label{ellraslebol2}
\ee

Thanks to  the property of conservativity, one has $\Rcvedge = 
-\Rcvvedge$ and $\rcvedge = -\rcvvedge$ for $\edge \in \edgesint$ 
such that $\edge = \edgecvcvv$. Let
$\Redge =
\abs{\Rcvedge}$ and $ \redge  =
\abs{\rcvedge} $  if $\edge \in \edgescv$. Reordering the 
summation over the edges and 
from  the Cauchy-Schwarz  inequality,  one then obtains

\be
\bac
 \vert \summesh \sumedgecv \medge (\Rcvedge +
\rcvedge) e_\cv \vert \le  \sumedges \medge (
D_\edge e)( \Redge +  \redge ) \le  \\
\Bigl(\sumedges \frac \medge \dedge
(  D_\edge e)^2\Bigr)^\half   \Bigl(\sumedges  \medge \dedge
( \Redge + \redge)^2\Bigr)^\half.
\ean
\label{ellraslebol3}
\ee

Now, since $\vert \Redge + \redge  \vert \le C_1 \size(\Mesh)$ and
since  $\sumedges  \medge \dedge = d \, \meas(\O)$,  
\eqref{ellraslebol2} and \eqref{ellraslebol3} yield the existence of $C_4 \in
\R_+$  only depending  on $u,\vbf$ and $\O$ such that
 
\[
\Vert e_\Mesh \Vert_{1, \Mesh}^2  + \half b \Vert e_\Mesh \Vert_{L^2(\O)}^2
\le  C_3 (\size(\Mesh))^2
+  C_4 \size(\Mesh)\Vert e \Vert_{1,\Mesh}. 
\]
Using again Young's inequality, there exists $C_5$ only depending on $u$, $\vbf$, $b$
and $\O$ such that

\be
\Vert e_\Mesh \Vert_{1, \Mesh}^2  + b\Vert e_\Mesh \Vert_{L^2(\O)}^2
\le  C_5 (\size(\Mesh))^2. 
\label{ellraslebol4}
\ee

This inequality yields Estimate \eqref{eqesterrh1} and, in the case $b >0$,
Estimate \eqref{eqesterr}.  In the case where $b  = 0$, one uses the discrete
Poincar\'e inequality  \eqref{ellinpoin} and the inequality \eqref{ellraslebol4}
to obtain 

\[
\Vert e_\Mesh \Vert_{L^2(\O)}^2  \le  \diamO^2 C_5 (\size(\Mesh))^2,
\]

which yields \eqref{eqesterr}.

\medskip

Remark now  that \eqref{eqesterrh1} can be written

\be
\bac
\dsp{ \sum_{\substack{\edge \in \edgesint\\\edge = \edgecvcvv} }   \medge   
\dedge  \Bigl( \frac {u_\cvv  - u_\cv} \dedge -
 \frac {u(\centercvv) - u(\centercv)} \dedge \Bigr)^2   + }
\\
\dsp{\sum_{\substack{\edge \in \edgesext \\ \edge \in \cvbar \cap \dr \O} }  \medge   
\dedge  \Bigl( \frac {g(\centeredge) - u_\cv} \dedge  -
 \frac {u(\centeredge) - u(\centercv)} \dedge \Bigr)^2\le (C 
\size(\Mesh))^2}.
\ean
\label{caspie3}
\ee

From Definition \eqref{eller2} and the consistency of the fluxes, one has

\be
\bac
\dsp{\sum_{\substack{\edge \in \edgesint\\\edge = \edgecvcvv} }   \medge    \dedge 
\Bigl(\frac {u(\centercvv) - u(\centercv)}
{\dedge} - \frac 1 {\medge } \int_{\edge } \grad
u(x)\cdot \ncvedge \dfrontiere(x) \Bigr)^2 }  +
\\
\dsp{\sum_{\substack{\edge \in \edgesext\\ \edge \in \cvbar \cap \dr \O} }   \medge   
\dedge  \Bigl(
 \frac {u(\centeredge) - u(\centercv)} \dedge - 
\frac 1 {\medge } \int_{\edge }
\grad u(x) \cdot \ncvedge \dfrontiere(x) \Bigr)^2 = } 
\\
\dsp{\sum_{\edge \in \edges}} \medge  \dedge 
\Redge^2 \le d \meas(\O) C_1^2 (\size(\Mesh))^2.
\ean
\label{caspie4}
\ee

Then \eqref{caspie3} and \eqref{caspie4} give \eqref{eqesterr1}. 
\eof

\subsection{$H^2$ error estimate} \label{sh2ee}

In Theorem \ref{ellesterr}, the hypothesis $u \in
C^2(\Obar)$ was used. In the following theorem (Theorem \ref{testh2}), one
obtains Estimates \eqref{eqesterrh1} and \eqref{eqesterr}, in the case $b=\vbf=0$ and
assuming some additional assumption on the mesh (see Definition \ref{ram}
below), under the weaker assumption
$u \in H^2(\O)$. This additional assumption on the mesh is not completely necessary
(see Remark \ref{h2plus} and \citet{ghv}). It is also possible to obtain Estimates \eqref{eqesterrh1}
and \eqref{eqesterr} in the cases $b \ne 0$ or $\vbf\ne 0$ assuming $u \in 
H^2(\O)$ (see Remark \ref{h2plus} and \citet{ghv}). Some similar
results are also in \citet{lazarovsiam} and \citet{coudiere}.

\begin{defi}[Restricted admissible meshes]\label{ram}  
\normalfont\index{Mesh!Admissible --!Restricted -- for Dirichlet boundary conditions}
Let $\O$ be an open bounded polygonal subset of $\R^d$, $d = 2$ or $3$. A
restricted admissible finite volume mesh of $\O$, denoted by $\mesh$, is an
admissible mesh in the sense of Definition \ref{meshdirichlet} such that, for
some $\zeta>0$, one has  $\dcvedge \ge
\zeta\diam(\cv)$ for all control volumes $\cv$ and for all $\edge \in \edgescv$.
\end{defi}

\begin{theo}[$H^2$ regularity] \label{testh2}\index{Error estimate!-- for an elliptic 
equation!-- in the general case}
Under Assumption \refp{ellH} with $b=\vbf=0$, let $\mesh$ be a restricted admissible
mesh in the sense of Definition \ref{ram} and $u_\mesh \in
X(\Mesh)$ (see Definition \refp{Xmesh}) be the approximate solution defined in
$\O$ by $u_\mesh (x)  = u_\cv$ for a.e. $x \in \cv$, for all $\cv \in
\mesh$, where $(u_\cv)_{\cv \in \mesh}$ is the (unique) solution to
 \eqref{ellschema}-\eqref{ellschemad} (existence and uniqueness of $(u_\cv)_{\cv
\in \mesh}$ are given by Lemma \ref{ellexistu}). 
Assume that the unique solution, $u$, of \eqref{ellsolvar} (with $b=\vbf=0$)
belongs to $H^2({\O})$. For each control volume $\cv$, let  $e_\cv =
u(\centercv) - u_\cv$, and $e_\Mesh \in X(\Mesh)$ defined by $e_\Mesh (x) = e_\cv$
for a.e.
$x \in \cv$, for all $\cv \in \mesh$. 

Then, there exists $C$, only depending  on $u$, $\zeta$ and $\O$, such that
\eqref{eqesterrh1}, \eqref{eqesterr} and \eqref{eqesterr1} hold.
\end{theo}

\begin{remark} \normalfont\label{uish2}~

\benum 
\item  In  Theorem \ref{testh2}, the function
$e_\mesh$ is still well defined, and  so is the quantity ``$\grad u \cdot \nedge$" 
on $\edge$, for all $\edge \in \edges$. Indeed, since $u
\in H^2(\O)$ (and $d \le 3$), one has $u \in C(\Obar)$ (and then $u(\centercv)$ is well
defined for all control volumes $\cv$) and $\grad u \cdot \nedge$ belongs to
$L^2(\edge)$ (for the
$(d-1)$-dimensional Lebesgue measure on $\edge$) for all $\edge \in \edges$.
\item Note that, under Assumption \ref{ellH} with
$b=\vbf=g=0$ the (unique) solution of
\eqref{ellsolvar} is necessarily in $H^2(\O)$ provided that $\O$ is convex.
\eenum
\end{remark}

\proof of Theorem \ref{testh2}

\smallskip

Let $\cv$ be a control volume and $\edge \in \edgescv$. Define $\V_{\cv,\edge}=$
$\{t\centercv +(1-t)x,$ $x \in \edge,$ $t \in [0,1]\}$. 
For $\edge \in \edgesint$,
let $\V_\edge=\V_{\cv,\edge}\cup \V_{\cvv,\edge}$, if $\cv$ and $\cvv$ 
are the control volumes such that $\edge=\edgecvcvv$. For $\edge \in \edgesext \cap
\edgescv$, let $\V_\edge=\V_{\cv,\edge}$.

\smallskip

The main part of the proof consists in proving the existence of some $C$,
only depending on the space dimension $d$ and $\zeta$ (given in Definition
\ref{ram}), such that, for all control volumes $\cv$ and for all $\edge \in
\edgescv$,

\be
\vert \Rcvedge \vert^2 \le C \frac {(\size(\mesh))^2} {\medge \dedge}
\int_{\V_\edge} \vert H(u)(z) \vert^2 dz,
\label{consh2}
\ee

where $H$ is the Hessian matrix of $u$ and

\[
\vert H(u)(z)\vert^2= \sum_{i,j=1}^{d} \vert D_iD_j u (z) \vert^2,
\]

and $D_i$ denotes the (weak) derivative with respect to the component $z_i$ of
$z=(z_1,\cdots,z_d)^t \in \R^d$.

Recall that
$\Rcvedge$ is the consistency error \index{Consistency!-- error} on the diffusion flux (see \eqref{eller2}), that
is:
\[
\Rcvedge=\frac {u(\centercvv)-u(\centercv)} {\dedge} - \frac 1 {\medge}\int_\edge
\grad u (x) \cdot \ncvedge \dfrontiere(x), \text{ if } \edge \in 
\edgesint \text{ and }
\edge=\edgecvcvv,
\]

\[
\Rcvedge=\frac {u(\centeredge)-u(\centercv)} {\dedge} - \frac 1 {\medge}\int_\edge
\grad u (x) \cdot \ncvedge \dfrontiere(x),  \text{ if } \edge \in \edgesext \cap \edgescv.
\]

Note that $\Rcvedge$ is well defined, thanks to $u \in H^2(\O)$, see 
Remark \ref{uish2}.

\smallskip

In Step $1$, one proves \eqref{consh2}, and, in Step 2, we conclude the proof of
Estimates \eqref{eqesterrh1}
and \eqref{eqesterr}.

\medskip

\textbf{Step $1$}. Proof of \eqref{consh2}.

Let $\edge \in \edges$. Since $u \in H^2(\O)$, the restriction of $u$ to $\V_\edge$
belongs to $H^2(\V_\edge)$. The space
$C^2(\overline {\V_\edge})$ is dense in $H^2(\V_\edge)$ (see, for 
instance, \citet{necas}, this can be 
proved quite easily be a regularization technique). Then, by a density
argument, one needs only to prove
\eqref{consh2} for $u \in C^2(\overline {\V_\edge})$. Therefore, in the remainder of
Step 1, it is assumed
$u \in C^2(\overline {\V_\edge})$.

\smallskip

First, one proves \eqref{consh2} if $\edge \in \edgesint$. Let $\cv$ and $\cvv$ be
the 2 control volumes such that $\edge=\edgecvcvv$.

It is possible to assume, for simplicity of notations and 
without loss of generality, that $\edge={0}\times
\tilde{\edge}$, with some $\tilde{\edge} \subset \R^{d-1}$, and
$\centercv=(-\a,0)^t$, $\centercvv=(\b,0)^t$, with some $\a > \zeta \diam(\cv)$,
$\b>\zeta \diam(\cvv)$ ($\zeta$ is defined in Definition \refp{ram}).
 
Let $x=(0,\tilde{x})^t \in \edge$.
In order to obtain a suitable integral remainder for the consistency error, as suggested in Remark \ref{rem:integral-remainder}, we introduce the function $\varphi: [0,1] \to \R$, defined by $\varphi(t) = u(t x_K +(1-t) x)$, which is twice continuously differentiable and we have:  
\begin{multline*}
\varphi (1) = x_K, \varphi(0) = u(x), \varphi'(t) = \nabla  u(t x_K +(1-t) x) \cdot  (\centercv-x)  \\ \text{ and } \varphi''(t) =H u(t x_K +(1-t) x) (\centercv-x)  \cdot  (\centercv-x), 
\end{multline*}
where $H(u)(z)$ denotes the Hessian matrix of $u$ at point $z$.
Therefore, writing that 
\[
\varphi(1) - \varphi(0)  = \int_0^1 \varphi'(t) \, dt  = \int_0^1 \varphi'(t)  \, dt  = \varphi'(0) -  \int_0^1(t-1)  \varphi''(t)   \, dt 
\]
 yields that
\[
u(\centercv)-u(x)= \int_0^1 H(u)(tx+(1-t)\centercv)(\centercv -x) \cdot (\centercv -x) t dt \text{ for a.e. }x=(0,\tilde{x})^t \in \edge
\]
 (for the $(d-1)$-dimensional Lebesgue measure on $\edge$).
Similarly, we have
\[
u(\centercvv)-u(x)=\grad u(x) \cdot (\centercvv-x) + \int_0^1
H(u)(tx+(1-t)\centercvv) (\centercvv -x) \cdot (\centercvv -x) t dt.
\]
Subtracting one equation to the other and integrating over $\edge$ yields (note that $\centercvv-\centercv =\ncvedge \dedge$) $\vert \Rcvedge \vert \le B_{\cv,\edge} + B_{\cvv,\edge}$, with
\be
B_{\cv,\edge} = \frac {C_1} {\medge \dedge} \int_\edge\int_0^1 \vert H(u)(tx+(1-t)\centercv) \vert \vert \centercv - x \vert^2 tdt\dfrontiere(x),
\label{donotpanic}
\ee
for some $C_1$ only depending  on $d$,
The quantity $B_{\cvv,\edge}$ is obtained with $B_{\cv,\edge}$ by changing $\cv$ in $\cvv$.

Let us perform the change of variables 
\begin{align*}
  h : & ]0,1[ \times \edge \to \V_{\cv,\edge}  \\
      & (t,x) \mapsto h(t,x) = tx+(1-t)\centercv,
\end{align*}
 in \eqref{donotpanic}. 
 let $z_1$ denote the first component of $z$ and $\bar z$ the $d-1$ last components of $z$;  thus  $z=(z_1,\overline{z})^t$ and $z_1=(t-1)\a$, so that 
\[
 dz=t^{d-1}\a \d t \dfrontiere(x).
\]
Since $\vert \centercv - x \vert \le \diam (\cv)$ we obtain\[
B_{\cv,\edge} \le \frac {C_1 (\diam(\cv))^2}{\medge \dedge} \int_{\V_{\cv,\edge}} \vert H(u)(z) \vert \frac {\a^{d-2}} {\a (z_1+\a)^{d-2}} dz.
\]
This gives, with the famous Cauchy-Schwarz inequality,
\[
B_{\cv,\edge} \le \frac {C_1 \a^{d-3}(\diam(\cv))^2}{\medge \dedge}\big(\int_{\V_{\cv,\edge}}\vert H(u)(z) \vert^2 dz \big)^{\frac 1 2}\big(\int_{\V_{\cv,\edge}}
\frac 1 {(z_1+\a)^{(d-2)2}} dz \big)^{\frac 1 2}.
\label{onrespirefort}
\]
For $d=2$, \eqref{onrespirefort} gives
\[
B_{\cv,\edge} \le \frac {C_1 (\diam(\cv))^2}{\a \medge \dedge} \big(\frac {\a \medge} 2 \big)^{\frac 1 2} \big(\int_{\V_{\cv,\edge}} \vert H(u)(z) \vert^2 dz \big)^{\frac 1 2},
\]
and therefore
\[
B_{\cv,\edge} \le \frac {C_1 (\diam(\cv))^2}{2^{\frac 1 2}(\medge \dedge)^{\frac 1 2} (\dedge \a)^{\frac 1 2}} \big(\int_{\V_{\cv,\edge}}
\vert H(u)(z) \vert^2 dz \big)^{\frac 1 2}.
\]

A similar estimate holds on $B_{\cvv,\edge}$ by changing $\cv$ in $\cvv$ and $\a$
in $\b$. Since $\a, \b \ge \zeta \diam(\cv)$ and $\dedge= \a + \b \ge 
\zeta \diam(\cv)$, these estimates on $B_{\cv,\edge}$ and $B_{\cvv,\edge}$ yield \eqref{consh2}
for some $C$ only depending  on $d$ and $\zeta$.



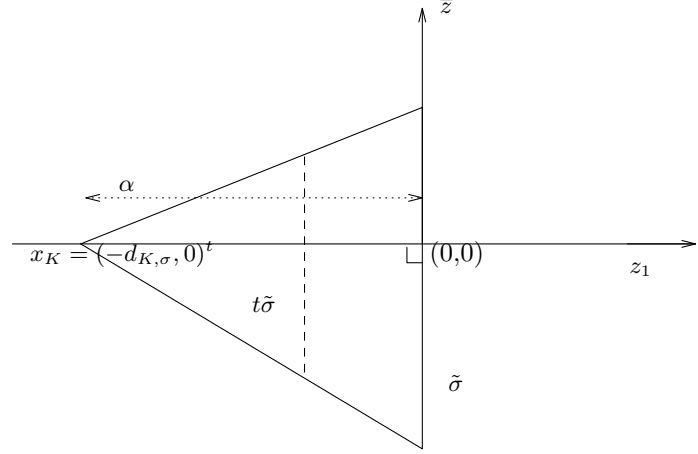
\begin{figure}[htb]
\begin{center}
\begin{tikzpicture}[x=0.00063333in, y=0.00063333in, line cap=round, line join=round]
  \draw (5007,1677) -- (5562,1677);
  \draw (5451.01,1649.255) -- (5562,1677) -- (5451.01,1704.745);
  \draw[dashed] (2385,2384) -- (2385,608);
  \draw (719.99,2078.745) -- (609,2051) -- (719.99,2023.255);
  \draw[dotted] (609,2051) -- (3342,2051);
  \draw (3231.01,2023.255) -- (3342,2051) -- (3231.01,2078.745);
  \draw (3217,1649) -- (3217,1524) -- (3342,1524);
  \draw (3342,1677) -- (3342,3591);
  \draw (3369.745,3480.01) -- (3342,3591) -- (3314.255,3480.01);
  \draw (567,1677) -- (3342,2786) -- (3342,12) -- (567,1677);
  \draw (12,1677) -- (5562,1677);
  \node[inner sep=0pt, anchor=base west, font=\small] at (3550,484) {$\tilde \sigma$};
  \node[inner sep=0pt, anchor=base west, font=\small] at (5021,1441) {$z_1$};
  \node[inner sep=0pt, anchor=base west, font=\small] at (872,2107) {$\alpha$};
  \node[inner sep=0pt, anchor=base west] at (3397,1538) {(0,0)};
  \node[inner sep=0pt, anchor=base west, font=\small] at (3481,3549) {$\overline z$};
  \node[inner sep=0pt, anchor=base west, font=\small] at (1954,1122) {$t \tilde \sigma$};
  \node[inner sep=0pt, anchor=base west, font=\small] at (151,1538) {$x_K = (-d_{K,\sigma},0)^t$};
\end{tikzpicture}

\end{center}
\caption{Consistency error, $d= 3$}
\label{fig-consist-h2}
\end{figure}

For $d=3$, the computation of the integral $A = \int_{\V_{\cv,\edge}}
\frac 1 {(z_1+\a)^{2}} dz$ by the following change of variable 
(see Figure~\ref{fig-consist-h2}):
\[
A = \int_{-d}^0 \frac {1} {(z_1 + \alpha)^2} (\int_{\zbar \in t \tilde
\edge} d \zbar) dz_1, \mbox{ where } t = \frac {z_1 + \alpha} \dcvedge.
\]
Now, 
\[\int_{\zbar \in t \tilde
\edge} d \zbar = \int_{y \in   \tilde
\edge} t^2   dy = \frac {(z_1 + \alpha)^2} {\alpha^2} \medge,\]
and therefore $A = \frac \medge \alpha$, and \eqref{onrespirefort}  
yields that:
\[
B_{\cv,\edge} \le \frac {C_3\, (\diam(\cv))^2}{(\medge\, \dedge^2\, 
\dcvedge)^{1/2}} \left(\int_{\V_{\cv,\edge}}
\vert H(u)(z) \vert^2\, dz \right)^{1/2}\le 
\frac {C_3\, \size(\mesh)}{\sqrt{2}\,\zeta\,(\medge\, \dedge)^{1/2}}\, 
\|H(u)\|_{L^2(\V_{\cv,\edge})}.
\]


and therefore  \eqref{onrespirefort} gives:

\[
B_{\cv,\edge} \le \frac {C_1 (\diam(\cv))^2}{\medge \dedge}
\big(\int_{-\a}^0 \frac {\medge} {\a^2} dz_1 \big)^{\frac 1 2}
\big(\int_{\V_{\cv,\edge}}
\vert H(u)(z) \vert^2 dz \big)^{\frac 1 2},
\]

and then

\[
B_{\cv,\edge} \le \frac {C_1 (\diam(\cv))^2}{(\medge \dedge)^{\frac 1 2} (\dedge
\a)^{\frac 1 2}} \big(\int_{\V_{\cv,\edge}}
\vert H(u)(z) \vert^2 dz \big)^{\frac 1 2}.
\]

With a similar estimate on $B_{\cvv,\edge}$, this yields \eqref{consh2} for some $C$
only depending  on $d$ and $\zeta$.

\medskip

Now, one proves \eqref{consh2} if $\edge \in \edgesext$. Let $\cv$ be the
control volume such that $\edge \in \edgescv$. One can assume, without loss of
generality, that $\centercv=0$ and $\edge=\{2\a\}\times \tilde{\edge}$ with
$\tilde{\edge} \subset \R^{d-1}$ and some $\a \ge \frac 1 2 \zeta \diam(\cv)$. The
above proof gives (see Definition \refp{meshdirichlet} for the definition of
$\centeredge$), with some $C_2$ only depending  on $d$,

\be
\vert \frac {u(\centeredge)-u(\centercv)} {2\a} - \frac {1} {\meas (\hat{\edge})}
\int_{\hat{\edge}} \grad u(x)\cdot \ncvedge \dfrontiere(x) \vert^2 \le C_2 \frac
{(\size(\mesh))^2} {\medge \dedge} \int_{\V_{\hat{\edge}}} \vert H(u)(z) \vert^2 dz,
\label{consbord1}
\ee
with $\hat{\edge}=\{(\a\,\frac x 2), x \in  \tilde{\edge}\}$, and
$\V_{\hat{\edge}}=\{t\centeredge + (1-t)x,$ $x \in \hat{\edge}$, $t \in [0,1]\}$
$\cup$ $\{t\centercv + (1-t)x,$ $x \in \hat{\edge}$, $t \in [0,1]\}$. Note that
$\meas (\hat{\edge}) = \frac {\medge} {2^{d-1}}$ and that $\V_{\hat{\edge}} \subset
\V_\edge$.

One has now to compare $I_\edge=\frac 1 {\medge}\int_\edge
\grad u (x) \cdot \ncvedge \dfrontiere(x)$ with $I_{\hat{\edge}}=\frac 1 {\meas
(\hat{\edge})} \int_{\hat{\edge}} \grad u(x)\cdot \ncvedge \dfrontiere(x)$.

A Taylor expansion gives

\[
I_\edge - I_{\hat{\edge}}= \frac 1 {\medge} \int_\edge \int_{\frac 1 2}^1
H(u)(\centercv + t (x-\centercv))(x-\centercv) \cdot \ncvedge dt \dfrontiere(x).
\]

The change of variables in this last integral $z=\centercv + t
(x-\centercv)$, which gives $dz=2\a t^{d-1} dt \dfrontiere(x)$, yields, with
$E_\edge=\{tx +(1-t)\centercv,$ $x \in \edge$, $t \in [\frac 1 2,1]\}$ and some
$C_3$ only depending  on $d$ (note that $t \ge \frac 1 2$),

\[
\vert I_\edge - I_{\hat{\edge}} \vert \le \frac {C_3} {\medge \a} \int_{E_\edge} \vert
H(u)(z) \vert \vert x-\centercv \vert dz.
\]

Then, from  the Cauchy-Schwarz 
inequality and since $\vert x-\centercv \vert \le \diam (\cv)$,

\be
\vert I_\edge - I_{\hat{\edge}} \vert^2 \le \frac {C_4 (\diam (\cv))^2}
{\medge \dedge}\int_{E_\edge} \vert H(u)(z) \vert^2 dz,
\label{consbord2}
\ee

with some $C_4$ only depending  on $d$ and $\zeta$.

Inequalities \eqref{consbord1} and \eqref{consbord2} yield \eqref{consh2} for some
$C$ only depending  on $d$ and $\zeta$.

\smallskip

One may therefore choose $C \in \R_+$ such that  \eqref{consh2} 
holds for $\edge \in \edgesint$ or  $\edge \in \edgesext$. This concludes Step $1$.

\medskip

\textbf{Step $2$}. Proof of Estimates \eqref{eqesterrh1}, \eqref{eqesterr} and
\eqref{eqesterr1}.

In order to obtain Estimate \eqref{eqesterrh1} (and therefore \eqref{eqesterr} from the discrete Poincar\'e inequality \eqref{ellinpoin}), one proceeds as in Theorem \ref{ellesterr}. 
Inequality \eqref{ellraslebol1} reads here, since $\Rcvedge=-\Rcvvedge$, if $\edge=\edgecvcvv$,
\[
\norm{e_{\mesh}}{1,\mesh}^2 \le \sumedges R_\edge \vert D_\edge e \vert \medge,
\]
with $R_\edge=\vert \Rcvedge \vert$, if $\edge \in \edgescv$. Recall also that
$\vert D_\edge e\vert =\vert e_\cv - e_\cvv \vert$ if $\edge \in
\edgesint$, $\edge=\edgecvcvv$ and
$\vert D_\edge e \vert =\vert e_\cv \vert$, if $\edge \in \edgesext \cap \edgescv$.
Cauchy and Schwarz strike  again:

\[
\norm{e_{\mesh}}{1,\mesh}^2 \le \big(\sumedges R_\edge^2 \medge
\dedge\big)^{\frac 1 2} \big(\sumedges \vert D_\edge e \vert^2 \frac {\medge}
{\dedge} \big)^{\frac 1 2}.
\]

The main consequence of \eqref{consh2} is that

\be
\sumedges \medge \dedge R_\edge^2 \le C(\sizemesh)^2 \sumedges\int_{\V_\edge}
\vert H(u)(z)\vert^2 dz =C(\sizemesh)^2 \int_{\O}\vert H(u)(z)\vert^2 dz.
\label{bouh}
\ee

Then, one obtains

\[
\norm{e_{\mesh}}{1,\mesh} \le \sqrt {C} \size(\mesh) \big( \int_{\O} \vert
H(u)(z)\vert^2 dz \big)^{\frac 1 2}.
\]

This concludes the proof of \eqref{eqesterrh1} since
$u \in H^2(\O)$ implies $\int_{\O} \vert H(u)(z)\vert^2 dz < \infty$.

Estimate \eqref{eqesterr1} follows from \eqref{bouh} in a similar manner as 
in 
the proof  of Theorem
\ref{ellesterr}.
This concludes the proof of Theorem \ref{testh2}.
\eof

\begin{remark}[Generalizations]\label{h2plus} \normalfont
~

\begin{enumerate}

\item By developing the method used to bound the consistency error 
\index{Consistency!-- error} on the flux
on the elements of $\edgesext$, it is possible to replace, in Theorem
\ref{testh2}, the hypothesis $\dcvedge
\ge \zeta \diam(\cv)$ in Definition \refp{ram} by the weaker hypothesis $\dedge
\ge \zeta \diam(\edge)$ provided that $\V_\edge$ is convex. Note also that, in this
case, the hypothesis $\centercv \in \cv$ is not necessary, it suffices that
$\centercvv-\centercv=\dedge\ncvedge$, for all $\edge \in \edgesint$,
$\edge=\edgecvcvv$ (for $\edge \in \edgesext$, one always needs
$\centeredge-\centercv=\dedge \ncvedge$).

\item It is also possible to prove Theorem \ref{testh2} if $b \ne 0$ or $\vbf \ne 0$ (or, of course, $b \ne 0$ and $\vbf \ne 0$). Indeed, if the solution, $u$, to
\eqref{ellsolvar} is not only in $H^2(\O)$ but is also Lipschitz continuous on
$\Obar$ (this is the case if, for instance, there exists $p > d$ such 
that $u \in W^{2,p}(\O)$), the treatment of the 
consistency\index{Consistency!-- error}
error terms due to the terms involving  $b$ and $\vbf$ are exactly as in Theorem
\ref{ellesterr}. If $u$ is not Lipschitz continuous on
$\Obar$, one has to deal with the consistency
error terms due to $b$ and $\vbf$ similarly as in the proof of Theorem 
\ref{testh2} (see also \citet{cvnl} or \citet{ghv}).

\end{enumerate}

\end{remark}

It is also possible, essentially under Assumption \refp{ellH}, to
obtain an $L^q$ estimate of the error, for $ 2 \le q < + \infty $ if $d = 2$,
and for  $ 1 \le q \le 6 $ if $d = 3$, see \cite{sobolev}.  The
error estimate for the $L^q$ norm is a consequence of the following lemma: 

\begin{lemma}[Discrete Sobolev Inequality]\index{Discrete!-- Sobolev 
inequality|(}\index{Sobolev!Discrete --  inequality} \label{lq}
Let $\O$ be an open bounded polygonal subset of $\R^d$ and 
$\mesh$ be a general finite volume mesh of $\O$ in the sense of 
definition \refp{meshneuman}, and let  $\zeta >0$ be such that
\be 
\forall K \in \mesh, \forall \edge \in \edges_{\cv}, \qquad
\dcvedge \geq \zeta \dedge, 
\label{regularitebis}\ee

Let be $u \in X(\Mesh)$  (see definition \refp{Xmesh}), then, 
there exists $C > 0$ only depending on $\O$ and $\zeta$,  such
that for all $q \in [2, + \infty)$, if $d=2$, and $q \in [2, 6]$, 
if $d=3$,
\be
\Vert u \Vert_{L^q(\O)} \le C  q \Vert u  \Vert_{1,\Mesh},
\label{ellinlq}
\ee
where $\Vert \cdot  \Vert_{1,\Mesh}$ is the discrete $\Hunz$ norm 
defined in definition \refp{defh10d}. 

\end{lemma}

\proof of Lemma \ref{lq}

\smallskip

Let us first prove the two-dimensional case. Assume $d=2$ and let $q \in
[2, + \infty)$. Let $\bfd_1 = (1,0)^t$ and  $\bfd_2 = (0,1)^t$; for $x \in \O$, let   $\D^1_x$ and
$\D^2_x$  be the  straight lines going through $x$ and defined by the  vectors
$\bfd_1$ and  $\bfd_2$. 

Let $v \in \Xmesh$. For all control volume $\cv$, one
denotes by $v_\cv$ the value of $v$ on $\cv$. For any control volume $\cv$ and
 a.e. $x \in \cv$, one has

\be  
 v_\cv^2 \le  \sumedges D_\edge v \, \chi^{(1)}_\edge (x) \sumedges D_\edge v \,
\chi^{(2)}_\edge (x), 
\label{lq1}
\ee

where  $\chi^{(1)}_\edge$ and  $\chi^{(2)}_\edge$ are defined by 

\[ 
\chi^{(i)}_\edge(x) = \ba 1 & \text{ if } \edge \cap \D^i_x \not = \emptyset \\
                          0 & \text{ if } \edge \cap \D^i_x   = 
                          \emptyset   \ea   \text{ for  } i = 1,2. 
\]

Recall that $D_\edge v=\vert
v_\cv-v_\cvv \vert$, if $\edge
\in \edgesint$, $\edge=\edgecvcvv$ and $D_\edge v=\vert v_\cv \vert$, if $\edge \in
\edgesext \cap \edgescv$.
Integrating \eqref{lq1} over $\cv$ and summing over $\cv \in \mesh$ yields 

\[
\int_\O v^2(x) dx \le \int_\O 
\Bigl( \sumedges D_\edge v \, \chi^{(1)}_\edge (x) \sumedges D_\edge v \, 
\chi^{(2)}_\edge (x) \Bigr)dx.
\]

Note that 
$\chi^{(1)}_\edge$  (resp.  $\chi^{(2)}_\edge$ ) only depends on the second
component
$x_2$ (resp. the first component $x_1$) of $x$ and that both functions are non zero
on a region the width of which is less than $\medge$; hence

\be
\int_\O v^2(x) dx \le \Bigl( 
 \sumedges \medge D_\edge v \Bigr)^2.
\label{poinbis}
\ee

Applying the inequality \eqref{poinbis} to $v = \vert u \vert^\a \sgn(u)$, where 
$u \in \Xmesh$ and $\a > 1$  yields 

\[
\int_\O \vert u(x) \vert^{2\a} dx \le \Bigl( 
 \sumedges \medge D_\edge v \Bigr)^2.
\]

Now, since
$\vert v_\cv-v_\cvv \vert \le
\a  (\abs{u_\cv}^{\a-1} + \abs{u_\cvv}^{\a-1}) \abs{u_\cv  - u_\cvv}$, if $\edge
\in \edgesint$, $\edge=\edgecvcvv$ and $\vert v_\cv \vert \le
\a  (\abs{u_\cv}^{\a-1}) \abs{u_\cv}$, if $\edge \in
\edgesext \cap \edgescv$,

\[
\Bigl(\int_\O \vert u(x) \vert^{2\a} dx\Bigr)^\half  \le \a \summesh \sumedgecv \medge \vert u_\cv
\vert^{\a-1} D_\edge u.
\]

Using H\"older's inequality with $p,{p'} \in \R_+$ such that $\frac 1 p + \frac 1
{p'} = 1$ yields that

\[
\Bigl(\int_\O \vert u(x) \vert^{2\a} dx\Bigr)^\half \le \a \big(\summesh \sumedgecv \vert u_\cv \vert^{p(\a-1)} \medge \dcvedge \big)^{\frac 1 p} \big(\summesh \sumedgecv \frac {\vert D_\edge u \vert^{p'}}{\dcvedge^{p'}} \medge \dcvedge \big)^{\frac 1 {p'}}.
\]
Since  $\sumedgecv \medge \dcvedge  =2\mcv$, this gives

\[
\Bigl(\int_\O \vert u(x) \vert^{2\a} dx\Bigr)^\half \le \a 2^{\frac 1 p} \big(\int_\O
\vert u(x)
\vert^{p(\a-1)}dx\big)^{\frac 1 p}  \big(\summesh \sumedgecv 
\frac {\vert D_\edge u \vert^{p'}}{\dcvedge^{p'}} \medge \dcvedge \big)^{\frac 1
{p'}},
\]

which yields, choosing $p$ such that $p(\a -1) = 2\a $, i.e. $p = \frac {2\a } {\a - 1}$ and  $p' = \frac {2\a }
{\a + 1}$,

\be
 \norm{u}{L^q(\O)} = \Bigl(\int_\O \vert u(x) \vert^{2\a} dx \Bigr)^{\frac 1 {2 \a}} \le 
\a  2^{\frac 1 p} \big(\summesh \sumedgecv 
\frac {\vert D_\edge u \vert^{p'}}{\dcvedge^{p'}} \medge \dcvedge \big)^{\frac 1
{p'}},
\label{lq2}
\ee

where $q = 2 \a$. Let $r = \frac 2 {p'}$ and $r' = \frac 2 {2-p'}$, 
 H\"older's inequality yields

\[
\summesh \sumedgecv 
\frac {\vert D_\edge u \vert^{p'}}{\dcvedge^{p'}} \medge \dcvedge \le
\big(\summesh \sumedgecv 
\frac {\vert D_\edge u \vert^{2}}{\dcvedge^2} \medge \dcvedge \big)^{\frac {p'} 2}
\big(\summesh \sumedgecv 
\medge \dcvedge \big)^{\frac {1} {r'}},
\]

replacing in \eqref{lq2} gives

\[
 \norm{u}{L^q(\O)}  \le 
\a 2^{\frac 1 p} (\frac 2 {\zeta})^{\frac 1 2} (2 \meas(\O))^{\frac 1 {p'r'}}
\Vert u \Vert_{1,\mesh}
\]

and then \eqref{ellinlq} with, for instance, $C=(\frac 2 {\zeta})^{\frac 1 2}((2
\meas(\O))^{\frac 1 2}+1)$.  
\smallskip 

Let us now prove the three-dimensional case. Let $d= 3$. Using the same 
notations as in the two-dimensional case, 
let $\bfd_1 = (1,0,0)^t$,   $\bfd_2 = (0,1,0)^t$ and  $\bfd_3 = (0,0,1)^t$ ; 
for $x \in \O$, let  $\D^1_x$, 
$\D^2_x$  and $\D^3_x$ be the  straight lines going through $x$ and defined by the  vectors
$\bfd_1$, $\bfd_2$ and  $\bfd_3$. Let us again define the functions
$\chi^{(1)}_\edge$,   $\chi^{(2)}_\edge$ and  $\chi^{(3)}_\edge$  by 

\[ 
\chi^{(i)}_\edge(x) = \ba 1 & \text{ if } \edge \cap \D^i_x \not = \emptyset \\
                          0 & \text{ if } \edge \cap \D^i_x   = 
                          \emptyset   \ea   \text{ for  } i = 1,2,3. 
\]

Let $v \in \Xmesh$ and let $A \in\R_+$ such that $\Omega \subset
[-A,A]^3$; we also denote by $v$ the function defined on $[-A,A]^3$ which 
equals $v$ on $\O$ and 0 on $[-A,A]^3 \setminus \O$.  By the Cauchy-Schwarz inequality, one has:
\be
\bal
\dsp{
\int_{-A}^A \int_{-A}^A   \vert v(x_1,x_2,x_3) \vert^{\frac 3 2} dx_1 
dx_2} \\
\hfill \dsp{ \le \Bigl(\int_{-A}^A \int_{-A}^A   \vert v(x_1,x_2,x_3) \vert  dx_1 
dx_2 \Bigr)^{\frac 1 2}   \Bigl(\int_{-A}^A \int_{-A}^A   \vert 
v(x_1,x_2,x_3) \vert^2  dx_1 
dx_2 \Bigr)^{\frac 1 2}. }
\ean
\label{sob3d1}
\ee
Now remark that 
\[  
\int_{-A}^A \int_{-A}^A   \vert v(x_1,x_2,x_3) \vert  dx_1 
dx_2  \le \sumedges D_\edge v  \int_{-A}^A \int_{-A}^A \chi^{(3)}_\edge(x) dx_1 dx_2 \le 
\sumedges \medge D_\edge v.\]

Moreover,   computations which were already performed in the two-dimensional case give 
that
\[  
\int_{-A}^A \int_{-A}^A   \vert v(x_1,x_2,x_3) \vert^2  dx_1 
dx_2  \le \int_{-A}^A \int_{-A}^A   \sumedges D_\edge v  \chi^{(1)}_\edge(x)   
\sumedges D_\edge v  \chi^{(2)}_\edge(x) dx_1 dx_2 \le 
\Bigr(\sumedges \meas(\edge_{x_3}) D_\edge v\Bigl)^2,\]
where $\edge_{x_3}$ denotes the intersection of $ \edge$ with  the plane 
which contains the point $(0,0,x_3)$ and is orthogonal to $\bfd_3$.
Therefore, integrating \eqref{sob3d1} in the third direction yields:

\be 
\int_\O \vert v(x) \vert^{\frac 3 2} dx \le \Bigr(\sumedges \medge D_\edge v 
\Bigl)^{\frac 3 2}.
\label{sob3d2}
\ee
Now let  $v = \vert u \vert^4 \sgn(u)$, since
$\vert v_\cv-v_\cvv \vert \le
4  (\abs{u_\cv}^{3} + \abs{u_\cvv}^{3}) \abs{u_\cv  - u_\cvv}$, 
Inequality \eqref{sob3d2} yields:
\[
\int_\O  \vert u(x) \vert^6 dx \le\Bigr[ 4 \summesh \sumedgecv \vert u_\cv 
\vert^3 D_\edge u \medge \Bigl]^{\frac 3 2 }.
\]
By Cauchy-Schwarz' inequality and since $\sumedgecv \medge \dcvedge = 3 
\mcv$, this yields
\[
\Vert u \Vert_{L^6} \le 
4 \sqrt{3} \summesh \sumedgecv  (D_\edge u )^2 \frac \medge \dcvedge, 
\] 
and since $\dcvedge \ge \zeta \dedge$, this yields 
  \eqref{ellinlq} with, for instance, $C=\frac {4\sqrt{3}} {\sqrt{\zeta}} $.  

\eof
\index{Discrete!-- Sobolev inequality|)}\index{Sobolev!Discrete --  inequality}
\begin{remark}[Discrete Poincar\'e Inequality] \label{poin2}  \normalfont
\index{Discrete!-- Poincar\'e inequality}\index{Poincar\'e!Discrete -- inequality}
 In the above proof, Inequality \eqref{poinbis} leads to another proof of some discrete
Poincar\'e inequality (as in Lemma \refp{thinpoin}) in the two-dimensional case. Indeed, let $\O$ be an open bounded
polygonal subset of $\R^2$. Let $\mesh$ be an admissible finite volume mesh of $\O$ in the
sense of Definition \refp{meshdirichlet} (but more general meshes are possible). Let $v \in
X(\Mesh)$. Then,   \eqref{poinbis}, the Cauchy-Schwarz  inequality and the fact that
$\sumedges \medge \dedge = 2 \, \meas(\O)$ yield 

\[
\norm{v}{L^2(\O)}^2 \le 2 \, \meas(\O) \norm{v}{1,\mesh}^2.
\]
A similar result holds in the three-dimensional case.\end{remark}

\begin{corol}[Error estimate]
	\label{norminf}\index{Error estimate!-- for an elliptic equation!-- in the general case}
	Under the same assumptions and with the same notations as in Theorem \refp{ellesterr}, or as in Theorem \refp{testh2}, and assuming that the mesh satisfies, for some $\zeta >0$, $\dcvedge \ge \zeta\dedge$, for all $\edge \in \edgescv$ and for all control volume \cv, there exists $C>0$ only depending  on $u$, $\zeta$ and $\O$ such that
	\begin{equation}
		\Vert  e_\Mesh  \Vert_{L^q(\O)} \le C q  \size(\Mesh); \text{ for any } q \in 
					\begin{cases}  
						[1,6]     & \mbox{ if  } \dim = 3,   \\
						[1,+\infty) & \mbox{ if  } \dim = 2;   
					\end{cases}
		\label{eqesterrq}
	\end{equation}
	furthermore, there exists $C \in \R_+$ only depending on $u$, $\zeta$, $\zetamesh = \min \{\frac {\mcv} {\sizemesh^d}, \cv \in \mesh\} $,  and  $\O$, such that
 	\begin{subequations}
 		\begin{align}
		&\Vert  e_\Mesh  \Vert_{L^\infty(\O)} \le  C \size(\Mesh) (\vert \ln( \size(\Mesh)) \vert +1)  ,  {\qquad\mathrm{if}\;} d=2. \label{eqesterrinfdeux}\\
		&\Vert  e_\Mesh  \Vert_{L^\infty(\O)} \le  C \size(\Mesh) ^{2/3}, {\qquad\mathrm{if}\;} d=3.\label{eqesterrinftrois}
		\end{align}
		\label{eqesterrinf}
 	\end{subequations}
\end{corol}

\proof  of Corollary \ref{norminf}
\smallskip

	Estimate \eqref{eqesterrh1} of Theorem \ref{ellesterr} (or Theorem \ref{testh2}) and Inequality \eqref{ellinlq} of Lemma  \ref{lq} immediately yield Estimate \eqref{eqesterrq} in the case $\dim = 2$.  
	Let us now prove \eqref{eqesterrinf}. 
	Remark that 
	\be
		\Vert  e_\Mesh  \Vert_{L^\infty(\O)} = \max \{ \vert \ecv \vert, \cv \in \mesh\} \le \Bigr(\frac 1 {\zetamesh \sizemesh^2}\Bigl)^{\frac 1 q} \norm{e_\mesh}{L^q}.
		\label{inf1}
	\ee
	For $d=2$, a study of the real function defined, for $q \ge 2$, by $q \mapsto \ln q + (1 -\frac 2 q )\ln h$ (with $h=\sizemesh$) shows that its minimum is attained for $q = -2 \ln h$, if $\ln h \le -\frac 1 2$. 
	Therefore \eqref{eqesterrq} and \eqref{inf1} yield\eqref{eqesterrinf}.  

	The 3 dimensional case is  an immediate consequence of\eqref{eqesterrq} with $q = 6$.
\eof

\section{Neumann boundary conditions}\label{nll}
\index{Neumann boundary conditions|(}
\hpage{Elliptic multid}

This section is devoted to the proof of convergence of the finite volume scheme when Neumann boundary conditions are imposed.
The discretization of a general convection-diffusion equation
\index{Equation!Convection-diffusion --}  with Dirichlet, Neumann and Fourier boundary conditions is considered in section \ref{ell2Dgen} below, and the convection term\index{Convection term} is largely studied in the previous section. 
Hence we shall limit here the presentation to the pure diffusion operator. 
Consider the following elliptic problem:
\be        
	- \Delta u(x) = f(x), \; x \in \O,        
	\label{nlldifstaf}
\ee  
with Neumann boundary conditions:
\be        
	\grad u(x) \cdot \bfn(x) = g(x), \;  x \in \dr\O,     
	\label{nlldifstab}
\ee 
where $\dr \O$ denotes the boundary of $\O$ and $\bfn$ its unit normal
vector outward to $\O$. 

The following assumptions are made on the data:

\begin{hyp}\label{nllH}~

\benum
\item $\O$ is an open bounded polygonal  connected 
subset of $\R^d$, $d=2$ or $3$, 
\item $g \in L^2(\dr\O)$, $f \in L^2(\O)$ and $\int_{\dr\O} g(x) \dfrontiere(x) +
\int_{\O} f(x) dx = 0$.
\eenum
\end{hyp}

Under Assumption \ref{nllH}, Problem \eqref{nlldifstaf}, \eqref{nlldifstab} has
a unique (variational) solution, $u$, belonging to $H^1(\O)$ and such that
$\int_\O u(x) dx = 0$. It is the unique solution of the following problem:

\be
u \in H^1(\O),  \;  \int_\O u(x)dx=0,
\label{nllspace}
\ee

\be
\int_\O \grad u(x) \grad \psi(x)  =  \int_\O f(x)  \psi(x) dx +\int_{\dr\O} g(x) 
\gammat(\psi)(x) \dfrontiere(x), \;    
\forall \psi\in H^1(\O).
\label{nllsolvar}
\ee

Recall that $\gammat$ is the ``trace" operator from $H^1(\O)$ to $L^2(\dr \O)$ (or to
$H^{\half}(\dr \O)$).

\subsection{Meshes and schemes}\label{msn}

\subsubsection{Admissible meshes}

The definition of the scheme in the case of Neumann boundary conditions is easier, since the finite volume scheme naturally introduces the fluxes on the boundaries in its formulation.  
Hence  the class of admissible meshes considered here is  somewhat wider than the one considered in Definition \refp{meshdirichlet}, thanks to the Neumann boundary conditions and the absence of convection term.  

\begin{defi}[Admissible meshes] \label{meshneuman} \index{Mesh!Admissible --!--for Neumann boundary conditions}
 \normalfont Let $\O$ be an open  bounded polygonal connected subset of $\R^d$, $d =  2,$ or $3$. 
An admissible finite volume mesh of $\O$ for the discretization of Problem \eqref{nlldifstaf}, \eqref{nlldifstab}, denoted by
$\mesh$,   is given by  a family of ``control volumes", which are open disjoint
 polygonal convex subsets of $\O$, a family of
 subsets  of $\Obar$ contained in hyperplanes of $\R^d$, denoted by  $\edges$ (these are the
``sides" of the control volumes), with strictly positive ($d-1$)-dimensional Lebesgue
measure, and a family of points of $\O$ denoted by 
$\points$ satisfying   properties $(i)$, $(ii)$, $(iii)$ and $(iv)$ of Definition
\refp{meshdirichlet}.

\medskip

The same notations as in Definition \refp{meshdirichlet}  are used in the sequel.

\end{defi}

 One defines the  set $ X(\Mesh)$ of piecewise constant functions on the control
volumes of an admissible mesh as in Definition \refp{Xmesh}.

\begin{defi}[Discrete $H^1$ seminorm] \label{defh1d}  \normalfont Let $\O$ be an open bounded
polygonal subset of $\R^d$, $d  = 2$ or $3$, and $\Mesh$ an admissible finite
volume mesh in the sense of Definition \ref{meshneuman}. 

For  $u \in X(\Mesh)$, the discrete $H^1$ seminorm of $u$ is defined by

\[
\vert u \vert_{1,\Mesh} = \Bigl(\sumedgesint \transedge (D_\edge u)^2
\Bigr)^{\frac 1 2},
\]

where $\transedge = \frac \medge \dedge$ and $\edgesint$ are defined
in Definition
\refp{meshdirichlet}, $u_\cv$ is the value of $u$ in the control volume $\cv$ and $D_\edge u = \vert
u_\cv - u_\cvv \vert$ if $\edge \in
\edgesint$, $\edge = \edgecvcvv$.

\end{defi}

\subsubsection{The finite volume scheme}     
        
Let  $\Mesh$ be an admissible mesh in the sense of Definition
\ref{meshneuman} \index{Finite volume scheme!-- for 
elliptic problems!-- in two or three 
space dimensions}. For $\cv \in
\mesh$, let us define:
        
\be        
f_\cv = \frac 1 {\mcv} \int_\cv f(x) dx,
\label{nllcondsou}
\ee  

\be
\bal      
\dsp{g_\cv = \frac 1 {\mcvbord} \int_{\dr\cv\cap\dr\O} g(x) \dfrontiere(x) \text{ if }
\mcvbord
\not = 0,}
\\
\dsp{g_\cv = 0   \text{ if } \mcvbord
 = 0.}
\ean
\label{nllcondbor}
\ee  

Recall that, in formula \eqref{nllcondsou}, $\mcv$ denotes the $d$-dimensional
Lebesgue measure of $\cv$, and, in \eqref{nllcondbor}, $\mcvbord$ denotes the
$(d-1)$-dimensional Lebesgue measure of $\dr \cv\cap\dr\O$. Note that $g_\cv=0$ if
the dimension of $\dr\cv\cap\dr\O$ is less than $d-1$.
Let $(u_\cv)_{\cv \in \mesh}$ denote the discrete unknowns; the numerical scheme
is defined by  \eqref{ellschema}-\refep{ellschemac}, with $b = 0$ and $\vbf 
= 0$. This yields:

\be        
- \sumv \transcvcvv  \Bigl(\ucvv - \ucv\Bigr) =  \mcv f_\cv + \mcvbord g_\cv, \
\forall \cv\in\mesh,  
\label{nllschema}\ee 

(see the notations in Definitions \refp{meshdirichlet} and 
\refp{meshneuman}). The condition \eqref{nllspace} is discretized by:  

\be        
\summesh  \mcv \ucv = 0.
\label{nllschemaO}\ee 

Then, the  approximate solution, $u_{\Mesh}$, belongs to $X(\Mesh)$ (see
Definition \refp{Xmesh}) and is defined by

\[
u_\Mesh(x)=u_\cv, \mathrm{~for~} \ae x \in \cv, \; \forall \cv \in \mesh.
\]

The following lemma gives existence and uniqueness of the solution of
\eqref{nllschema} and \eqref{nllschemaO}.

\begin{lemma}\label{nllexistu}
Under Assumption \ref{nllH}. let $\Mesh$ be an admissible mesh  (see Definition
\ref{meshneuman}) and $\{f_\cv$, $\cv \in \mesh\}$, $\{g_\cv$,
$\cv \in \mesh\}$ defined by \eqref{nllcondsou}, \eqref{nllcondbor}. Then, there
exists a unique solution
$(u_\cv)_{\cv\in \mesh}$ to \eqref{nllschema}-\eqref{nllschemaO}.
\end{lemma}
    
\proof  of lemma \ref{nllexistu}

\smallskip

Let $N=\card(\mesh)$. The equations \eqref{nllschema} are a system of $N$
equations with $N$ unknowns, namely $(u_\cv)_{\cv\in \mesh}$. Ordering the
unknowns (and the equations), this system can be written under a matrix form with  a
$N \times N$ matrix $A$. 
Using the connexity of $\O$, the null space of this matrix is the set of ``constant"
vectors (that is $u_\cv=u_\cvv$, for all $\cv,\cvv \in \mesh$). Indeed, if
$f_\cv=g_\cv=0$ for all $\cv \in \mesh$ and $\{u_\cv$, $\cv \in \mesh\}$ is
solution of \eqref{nllschema}, multiplying \eqref{nllschema} (for $\cv \in \mesh$)
by $u_\cv$ and summing over $\cv \in \mesh$ yields

\[
\sumedgesint \transedge (D_\edge u)^2 =0,
\]

where $D_\edge u = \vert
u_\cv - u_\cvv \vert$ if $\edge \in
\edgesint$, $\edge = \edgecvcvv$. This gives, thanks to the positivity of $\transedge$ and the
connexity of $\O$, $u_\cv=u_\cvv$, for all $\cv,\cvv \in \mesh$.

\smallskip

For general $(f_\cv)_{\cv\in \mesh}$ and $(g_\cv)_{\cv\in \mesh}$, a necessary
condition, in order that
\eqref{nllschema} has a solution, is that
\be
\sum_{\cv \in \mesh} (\mcv f_\cv + \mcvbord  g_\cv)=0. \label{finetfin}
\ee
Since the dimension of the null
space of $A$ is one, this condition is also a sufficient condition. Therefore,
System
\eqref{nllschema} has a solution if and only if \eqref{finetfin} holds, 
and this solution is unique up to an additive constant.
Adding condition \eqref{nllschemaO} yields uniqueness. 
Note that   \eqref{finetfin} holds thanks to the second item of Assumption
\ref{nllH}; this  concludes the proof
of Lemma \ref{nllexistu}.
\eof

\subsection{Discrete Poincar\'e inequality}\index{Discrete!-- Poincar\'e inequality|(}\index{Poincar\'e!Discrete -- inequality|(}
\index{Poincar\'e  inequality!-- Neumann boundary conditions}
The proof of an error estimate, under a regularity assumption on the exact solution, and of a convergence result, in the general case (under Assumption
\ref{nllH}), requires a ``discrete Poincar\'e" inequality as in the case of the Dirichlet problem.

\begin{lemma}[Discrete mean Poincar\'e inequality]\label{thinpoin1}
Let $\O$ be an open bounded polygonal connected subset of $\R^d$, $d=2$ or $3$. Then, there exists $C \in \Rp$, only depending  on $\O$, such that for all admissible meshes (in the sense of Definition \refp{meshneuman}), $\Mesh$, and for all $u \in X(\Mesh)$ (see Definition \refp{Xmesh}), the following inequality holds:
\be
\Vert u \Vert_{L^2(\O)}^2 \le C \vert u  \vert_{1,\Mesh}^2  +  2
(\meas(\O))^{-1} (\int_\O u(x) dx)^2,
\label{ellinpoin1}
\ee

where $\vert \cdot \vert_{1,\Mesh}$ is the discrete $H^1$ seminorm defined
in Definition \ref{defh1d}.
\end{lemma}

\proof of Lemma \ref{thinpoin1}

\smallskip

The proof given here is a ``direct proof"; another proof,
by contradiction, is possible (see Remark \ref{poincontr}).
Let $\Mesh$ be an admissible mesh and $u
\in X(\Mesh)$. 
Let $m_\O(u)$ be the mean value of $u$ over $\O$, that is
\[
m_\O(u)=\frac 1 {\meas(\O)} \int_\O u(x)dx.
\]
Since 
\[\norm{u}{L^2(\O)}^2\le
2\norm{u-m_\O(u)}{L^2(\O)}^2 + 2 (m_\O(u))^2\meas(\O),
\] 
proving Lemma \ref{thinpoin1} amounts
to proving the existence of $D \ge 0 $, only depending  on $\O$, such that

\be
\norm{u-m_\O(u)}{L^2(\O)}^2 \le D \vert u \vert_{1,\Mesh}^2.
\label{poineff}
\ee

The proof of \eqref{poineff} may be decomposed into three steps (indeed, if $\O$ is convex, the 
first step is sufficient).

\smallskip


\textit{Step 1 (Estimate on a convex part of $\O$)}

Let $\omega$ be an open 
convex subset of $\O$, $\omega \ne \emptyset$ and $m_\omega (u)$ be the mean value of $u$ on 
$\omega$. In this step, one proves that there exists $C_0$, depending 
only on $\O$, such that

\be
\norm {u-m_\omega(u)}{L^2(\omega)}^2  \le \frac 1 {\meas(\omega)} 
C_0 \vert u \vert_{1,\Mesh}^2.
\label{pstep2}
\ee

(Taking $\o=\O$, this proves \eqref{poineff} and Lemma 
\ref{thinpoin1} in the case where $\O$ is convex.)

~

Noting that

\[
\int_\o(u(x)-m_\o(u))^2dx \le \frac 1 {\meas(\o)}\int_\o \big(
\int_\o( u(x)-u(y))^2dy \big)dx,
\]

\eqref{pstep2} is proved provided that
there exists $C_0 \in \Rp$, only depending on $\O$, such that

\be
\int_\o\int_\o (u(x)-u(y))^2dxdy \le C_0 \vert u  \vert_{1,\Mesh}^2.
\label{poinconv}
\ee

~

For $\edge \in \edgesint$, let the function
$\chi_{\edge}$ from $\R^\dim \times \R^\dim$ to $\{0,1\}$ be defined by

\[
\bac
\chi_{\edge}(x,y)=1, \mathrm{~if~} x, y \in \Obar, \, [x,y]\cap\edge \neq
\emptyset, 
\\
\chi_{\edge}(x,y)=0, \mathrm{~if~} x \notin \Obar \, \mathrm{~or~} y \notin \Obar \, 
\mathrm{~or~} [x,y]\cap\edge = \emptyset.
\ean
\]

(Recall that $[x,y]=\{tx+(1-t)y,$ $t\in [0,1]\}$.)
For a.e. $x,y \in \o$, one has, with $D_\edge u = \vert
u_\cv - u_\cvv \vert$ if $\edge \in
\edgesint$, $\edge = \edgecvcvv$,

\[
(u(x)-u(y))^2 \le \big(\sumedgesint \vert D_{\edge}u \vert \chi_{\edge}(x,y)
\big)^2,
\]

(note that the convexity of $\o$ is used here) which yields,  thanks to   the Cauchy-Schwarz  inequality,

\be
(u(x)-u(y))^2 \le \sumedgesint \frac {\vert D_{\edge}u
\vert^2}{\dedge c_{\edge,y-x}}\chi_{\edge}(x,y) \sumedgesint \dedge
c_{\edge,y-x}\chi_{\edge}(x,y),
\label{astcospoin}
\ee

with 
\[
c_{\edge,y-x}=\vert \frac {y-x}{\vert y-x \vert}\cdot \nedge \vert,
\]

recall
that $\nedge$ is a unit normal vector to $\edge$, and that $x_\cv-x_\cvv=\pm
\dedge \nedge$ if $\edge \in \edgesint$, $\edge = \edgecvcvv$.
For a.e. $x,y \in \o$, one has 

\[
\sumedgesints \dedge
c_{\edge,y-x}\chi_{\edge}(x,y)= \vert (x_\cv - x_\cvv) \cdot \frac {y-x}{\vert
y-x \vert} \vert,
\]

for some convenient control volumes $\cv$ and $\cvv$,
depending on $x$, $y$ and $\edge$ (the convexity of $\o$ is used  again 
here). Therefore,

\[
\sumedgesints \dedge
c_{\edge,y-x}\chi_{\edge}(x,y) \le \diam(\O).
\]

Thus, integrating
\eqref{astcospoin} with respect to $x$ and
$y$ in $\o$,

\[
\int_\o \int_\o(u(x)-u(y))^2 dxdy \le \diamO \int_\o \int_\o \sumedgesint \frac
{\vert D_{\edge}u
\vert^2}{\dedge c_{\edge,y-x}}\chi_{\edge}(x,y)dxdy,
\]

which gives, by a change of variables,

\be
\int_\o \int_\o(u(x)-u(y))^2 dxdy \le \diamO \int_{\R^d} \big(\sumedgesint \frac
{\vert D_{\edge}u
\vert^2}{\dedge c_{\edge,z}} \int_\o \chi_{\edge}(x,x+z)dx \big) dz.
\label{astcospoin1}
\ee

Noting that, if $\vert z \vert
> \diamO$, $\chi_{\edge}(x,x+z)=0$, for a.e. $x \in \O$,  and 
\[
\int_\O
\chi_{\edge}(x,x+z)dx \le \medge \vert z \cdot \nedge \vert = \medge \vert z
\vert c_{\edge,z} \text{ for a.e. }z \in \R^d,
\]
therefore, with \eqref{astcospoin1}:

\[
\int_\o \int_\o (u(x)-u(y))^2 dxdy \le (\diamO)^2 \mathrm{m}(B_{\O}) \sumedgesint
\frac {\medge\vert D_{\edge}u
\vert^2}{\dedge},
\]
where $B_{\O}$ denotes the ball of $\R^d$ of center $0$ and radius $\diamO$.

This inequality proves \eqref{poinconv} and then \eqref{pstep2} with $C_0=(\diamO)^2 \mathrm{m}(B_{\O})$ (which only depends on $\O$). Taking $\o=\O$, it 
concludes the proof of Lemma \ref{thinpoin1} in the
case where $\O$ is convex.

~

In order to consider in Step 2 the case $d=3$, it is useful to remark that if $\tilde \omega$ is an open convex subset of $\omega$,
we also have 
\[
\norm {u-m_{\tilde \omega}(u)}{L^2(\tilde \omega)}^2  \le \frac 1 {\meas(\tilde \omega)} 
C_0 \vert u \vert_{1,\Mesh}^2,
\]

\textit{Step 2 (Estimate with respect to the mean value on a part of the boundary)}

In this step, one proves the same inequality than \eqref{pstep2} but with 
the mean value of $u$ on a
(arbitrary) part $I$ of the boundary of $\o$ instead of $m_\o(u)$ and with 
a convenient $C_1$ depending on $I$, $\O$ and $\o$ instead of $C_0$.

More precisely, let $\o$ be a polygonal open 
convex subset of $\O$ and let $I \subset \dr \o$, with $m(I)>0$ ($m(I)$ 
is the ($\dim-1$)-Lebesgue measure
of $I$). Assume that $I$ is included
in a hyperplane of $\R^\dim$.
Let $\ug$ be the ``trace" of $u$ on the boundary of $\o$, that
is $\ug(x)=u_\cv$ if $x \in \dr \o \cap \cvbar$, for $\cv \in \mesh$. (If 
$x \in \cvbar \cap \cvvbar$, the choice of $\ug(x)$ between 
$u_\cv$ and $u_\cvv$ does not matter). Let $m_I(u)$ be the mean value of $\ug$ on $I$.
This step is devoted to the proof that there exists $C_1$, only depending on 
$\O$, $\o$ and $I$, such that

\be
\norm {u-m_I(u)}{L^2(\omega)}^2  \le C_1 \vert u \vert_{1,\Mesh}^2.
\label{pstep3}
\ee

We first consider  the case $\dim=2$. Since $I$ is 
included in a hyperplane, it may be assumed, without loss of generality, 
that
$I=\{0\}\times J$, with $J \subset \R$ and $\o \subset \R_+\times
\R$ (one uses here the convexity of $\o$). 

\smallskip


Let $\a =\max\{x_1$, $x=(x_1,x_2)^t \in \obar\}$ and
$a=(\a,\b)^t \in \obar$. In the following, $a$ is fixed.
For a.e. $x=(x_1,x_2)^t \in \o$ and  for a.e. (for the $1$-Lebesgue 
measure)
$y=(0,\ybar)^t
\in I$ (with $\ybar \in J$), one sets
$z(x,y)=ta+(1-t)y$ with $t= x_1/\a$. Note that, thanks to the 
convexity of $\o$,  $z(x,y)=(z_1,z_2)^t \in
\obar$, with $z_1=x_1$. The following inequality holds:

\[
\pm (u(x)-\ug(y)) \le \vert u(x)-u(z(x,y))\vert + \vert u(z(x,y)-\ug(y))\vert.
\]

In the following, the notation  $C_i$, $i \in \N^\star$, will be used for 
quantities
only depending  on $\O$, $\o$ and $I$.

Let us integrate the above inequality over $y \in I$, take the power 
$2$, from  the Cauchy-Schwarz 
inequality,  an  integration over $x \in \o$ leads to

\[
\bal
\dsp{\int_\o(u(x)-m_I(u))^2dx \le \frac 2 {m(I)} \int_\o\int_I
(u(x)-u(z(x,y)))^2\dfrontiere(y)dx}
\\
\dsp{ + \frac 2 {m(I)} \int_\o\int_I
(u(z(x,y))-u(y))^2\dfrontiere(y)dx.}
\ean
\]

Then,

\[
\int_\o(u(x)-m_I(u))^2dx \le \frac 2 {m(I)} (A+B),
\]

with, since $\o$ is  convex, 

\[
A=\int_\o\int_I \big( \sumedgesint \vert D_\edge u \vert
\chi_{\edge}(x,z(x,y))\big)^2\dfrontiere(y)dx,
\]
and
\[
B=\int_\o\int_I \big( \sumedgesint \vert D_\edge u \vert
\chi_{\edge}(z(x,y),y)\big)^2\dfrontiere(y)dx.
\]

Recall that, for $\xi, \eta \in \Obar$, $\chi_{\edge}(\xi, \eta)=1$ if 
$[\xi, \eta]\cap\edge \neq \emptyset$ and $\chi_{\edge}(\xi, \eta)=0$ if 
$[\xi, \eta]\cap\edge=\emptyset$.
Let us now look for some bounds of $A$ and $B$ of the form
 $C\vert u \vert_{1,\mesh}^2$.

\smallskip

The bound for $A$ is easy. Using the Cauchy-Schwarz  inequality and the fact that

\[
\sumedgesint c_{\edge,x-z(x,y)} \dedge \chi_{\edge}(x,z(x,y)) \le \diamO
\]

(recall that $c_{\edge,\eta}=\vert \frac {\eta} {\vert \eta \vert}\cdot 
\nedge \vert$ for $\eta \in \R^2 \setminus {0}$) 
gives 

\[
A \le C_2 \int_\o\int_I \sumedgesint \frac {\vert D_\edge u
\vert^2\chi_{\edge}(x,z(x,y))} {c_{\edge,x-z(x,y)} \dedge}dx\dfrontiere(y).
\]

Since $z_1=x_1$, 
one has $c_{\edge,x-z(x,y)}=c_{\edge,e}$, with $e=(0,1)^t$. Let us perform the
integration of the right hand side  of the previous inequality, 
 with respect to the first component of $x$, denoted by $x_1$, first. 
  The result of the integration with respect to
$x_1$ is bounded by 
$\vert u \vert_{1,\mesh}^2$. Then, integrating with respect to $x_2$ and
$y \in I$ gives $A \le C_3 \vert u \vert_{1,\mesh}^2$.

\smallskip

In order to obtain a bound $B$, one remarks, as for $A$, that 

\[
B \le C_4 \int_\o\int_I \sumedgesint \frac {\vert D_\edge u
\vert^2\chi_{\edge}(z(x,y),y)} {c_{\edge,y-z(x,y)} \dedge}dx\dfrontiere(y).
\]

In the right hand side of 
this inequality, the integration with respect to $y \in I$ is transformed 
into 
an integration with respect to $\xi=(\xi_1,\xi_2)^t \in \edge$, 
this yields (note that
$c_{\edge,y-z(x,y)}=c_{\edge,a-y}$)

\[
B \le C_4 \sumedgesint \frac {\vert D_\edge u
\vert^2} {\dedge} \int_\o \int_\edge \frac{\psi_\edge(x,\xi)}
{c_{I,a-y(\xi)}} \frac {\vert a-y(\xi) \vert} {\vert a-\xi \vert} dx
\dfrontiere(\xi),
\]

where $y(\xi)=s\xi + (1-s)a$, with $s\xi_1 + (1-s)\a=0$, and where $\psi_\edge$ is
defined by 

\[
\bac
\psi_\edge(x,\xi)=1, \text{~if~} y(\xi) \in I \text{~and~} \xi_1 \le x_1 \\
\psi_\edge(x,\xi)=0, \text{~if~} y(\xi) \not\in I \text{~or~} \xi_1 > x_1.
\ean
\]

Noting that $c_{I,a-y(\xi)} \ge C_5 >0$, one
deduces that

\[
B \le  C_6 \sumedgesint \frac {\vert D_\edge u
\vert^2} {\dedge} \int_\edge \big( \int_\o  \psi_\edge(x,\xi) \frac 
{\vert a-y(\xi) \vert}
{\vert a-\xi \vert} dx \big) \dfrontiere( \xi) \le C_7 \vert u \vert_{1,\mesh}^2,
\]

with, for instance, $C_7=C_6 (\diam(\o))^2$.
The bounds on $A$ and $B$ yield \eqref{pstep3}.

\bigskip

We now consider  the case $\dim=3$, which is a little more tricky. Since $I$ is 
included in a hyperplane, it may be assumed, without loss of generality, 
that
$I=\{0\}\times J$, with $J \subset \R^2$ and $\o \subset \R_+\times
\R^{2}$ (one uses here the convexity of $\o$). 

\smallskip


Let $\a =\max\{x_1$, $x=(x_1,x_2)^t \in \obar\}$ and
$a=(\a,\b)^t \in \obar$. In the following, $a$ is fixed.

We set $\tilde \o=\{(x_1,x_2)^t \in \o,$ $0< x_1 < \frac 1 2 \a\}$.
Note that $\tilde\o$ is an open convex subset of $\o$.

For a.e. $x=(x_1,x_2)^t \in \tilde\o$ and  for a.e. (for the $2$-Lebesgue 
measure)
$y=(0,\ybar)^t
\in I$ (with $\ybar \in J$), one sets
$z(x,y)=ta+(1-t)y$ with $t= x_1/\a$. Note that, thanks to the 
convexity of $\o$,  $z(x,y)=(z_1,z_2)^t \in
\obar$, with $z_1=x_1$ (so that $0<t < \frac 1 2$). The following inequality holds:

\[
\pm (u(x)-\ug(y)) \le \vert u(x)-u(z(x,y))\vert + \vert u(z(x,y)-\ug(y))\vert.
\]

In the following, the notation  $C_i$, $i \in \N^\star$, will be used for 
quantities
only depending  on $\O$, $\o$ and $I$.

Let us integrate the above inequality over $y \in I$ and over $x \in \tilde \o$, 
\[
m(I)m(\tilde w) |m_{\tilde\o}(u)-m_I(u)| \le \int_{\tilde \o }\int_I
|u(x)-u(z(x,y))|\dfrontiere(y)dx +
\int_{\tilde \o }\int_I
|u(z(x,y))-u(y)|\dfrontiere(y)dx.
\]

Then, using the Cauchy-Schwarz inequality,  
\[
\frac {m(I)m(\tilde w)} 2 |m_{\tilde\o}(u)-m_I(u)|^2 \le \int_{\tilde \o }\int_I
|u(x)-u(z(x,y))|^2\dfrontiere(y)dx +
\int_{\tilde \o }\int_I
|u(z(x,y))-u(y)|^2\dfrontiere(y)dx.
\]

%
Then,

\[
|m_{\tilde \o}(u)-m_I(u)]^2 \le \frac 2 {m(I)m(\tilde w)} (A+B),
\]

with, since $\tilde \o$ is  convex, 

\[
A=\int_{\tilde \o}\int_I \big( \sumedgesint \vert D_\edge u \vert
\chi_{\edge}(x,z(x,y))\big)^2\dfrontiere(y)dx,
\]
and
\[
B=\int_{\tilde\o}\int_I \big( \sumedgesint \vert D_\edge u \vert
\chi_{\edge}(z(x,y),y)\big)^2\dfrontiere(y)dx.
\]

Recall that, for $\xi, \eta \in \Obar$, $\chi_{\edge}(\xi, \eta)=1$ if 
$[\xi, \eta]\cap\edge \neq \emptyset$ and $\chi_{\edge}(\xi, \eta)=0$ if 
$[\xi, \eta]\cap\edge=\emptyset$.
Let us now look for some bounds of $A$ and $B$ of the form
 $C\vert u \vert_{1,\mesh}^2$.

\smallskip

Using the Cauchy-Schwarz  inequality and the fact that

\[
\sumedgesint c_{\edge,x-z(x,y)} \dedge \chi_{\edge}(x,z(x,y)) \le \diamO
\]

(recall that $c_{\edge,\eta}=\vert \frac {\eta} {\vert \eta \vert}\cdot 
\nedge \vert$ for $\eta \in \R^2 \setminus {0}$) 
gives 

\[
A \le C_2 \int_{\tilde \o}\int_I \sumedgesint \frac {\vert D_\edge u
\vert^2\chi_{\edge}(x,z(x,y))} {c_{\edge,x-z(x,y)} \dedge}dx\dfrontiere(y).
\]

Using an integration with respect to the two first components of $x$, denoted by $x_1$ and $x_2$,
and then an integration with respect to the third component, $x_3$, we obtain, for all $\edge$,
\[
\int_{\tilde \o}\frac {\chi_{\edge}(x,z(x,y))} {c_{\edge,x-z(x,y)} }dx\le C_3= \medge \diamO.
\]
Then, integrating with respect to
$y \in I$ gives $A \le C_4 \vert u \vert_{1,\mesh}^2$.

\smallskip

In order to obtain a bound $B$, one remarks, as for $A$, that 

\[
B \le C_5 \int_{\tilde \o}\int_I \sumedgesint \frac {\vert D_\edge u
\vert^2\chi_{\edge}(z(x,y),y)} {c_{\edge,y-z(x,y)} \dedge}dx\dfrontiere(y).
\]

In the right hand side of 
this inequality, the integration with respect to $y \in I$ is transformed 
into 
an integration with respect to $\xi=(\xi_1,\xi_2)^t \in \edge$, 
this yields (note that
$c_{\edge,y-z(x,y)}=c_{\edge,a-y}$)

\[
B \le C_5 \sumedgesint \frac {\vert D_\edge u
\vert^2} {\dedge} \int_{\tilde \o} \int_\edge \frac{\psi_\edge(x,\xi)}
{c_{I,a-y(\xi)}} \frac {\vert a-y(\xi) \vert^2} {\vert a-\xi \vert^2} dx
\dfrontiere(\xi),
\]

where $y(\xi)=s\xi + (1-s)a$, with $s\xi_1 + (1-s)\a=0$, and where $\psi_\edge$ is
defined by 

\[
\bac
\psi_\edge(x,\xi)=1, \text{~if~} y(\xi) \in I \text{~and~} \xi_1 \le x_1 \\
\psi_\edge(x,\xi)=0, \text{~if~} y(\xi) \not\in I \text{~or~} \xi_1 > x_1.
\ean
\]

Noting that $c_{I,a-y(\xi)} \ge C_6 >0$ and $s \ge \frac 1 2$ (since $y(\xi)$ is between $y$ and $z(x,y)$), one
deduces that

\[
B \le  C_7 \sumedgesint \frac {\vert D_\edge u
\vert^2} {\dedge} \int_\edge \big( \int_{\tilde \o}  \psi_\edge(x,\xi) \frac 
{\vert a-y(\xi) \vert^2}
{\vert a-\xi \vert^2} dx \big) \dfrontiere( \xi) \le C_8 \vert u \vert_{1,\mesh}^2.
\]
The bounds on $A$ and $B$ yield
\[
|m_{\tilde \o}(u)-m_I(u)] \le C_9 \vert u \vert_{1,\mesh}.
\]
We now recall that at the end of Step 1 we proved 
\[
\norm {u-m_{\tilde \omega}(u)}{L^2(\tilde \omega)} \le C_{10} \vert u \vert_{1,\Mesh},
\]
These two inequalities give \eqref{pstep3}.

\smallskip

\textit{Step 3 (proof of \eqref{poineff})}

Let us now prove that there exists $D \in
\Rp$, only depending  on $\O$ such that \eqref{poineff} hold.
Since $\O$ is a polygonal set ($d=2$ or $3$), there
exists a finite number of disjoint convex polygonal sets, 
denoted by $\{\O_1,\ldots,\O_n\}$, such
that $\overline \O = \cup_{i=1}^n \overline {\O_i}$. 
Let $I_{i,j}=\overline {\O_i}\cap\overline
{\O_j}$, and $B$ be the set of couples $(i,j) \in \{1,\ldots,n\}^2$ such that
$i \neq j$ and the $(d-1)$-dimensional Lebesgue measure of $I_{i,j}$, denoted by
$m(I_{i,j})$, is positive.

Let $m_i$ denote the mean value of $u$ on $\O_i$, $i \in \{1,\ldots,n\}$,
and $m_{i,j}$ denote the mean value of $u$ on $I_{i,j}$,
$(i,j) \in B$. (For $\edge \in \edgesint$, in order that $u$ be defined 
on $\edge$, a.e. for the ($d-1$)-dimensional Lebesgue measure, let $\cv \in 
\mesh$  be  a control 
volume such that $\edge \in \edgescv$,  
one sets $u=u_\cv$ on $\edge$.) Note that $m_{i,j}=m_{j,i}$ for all
$(i,j)\in B$. 

~

Step 1 gives the 
existence of $C_i$, $i \in \{1,\ldots,n\}$, only depending  on $\O$ (since 
the $\O_i$  only depend on $\O$), such that

\be
\norm{u-m_i}{L^2(\O_i)}^2 \le C_i \snorm{u}{1,\Mesh}^2, \; \forall i \in
\{1,\ldots,n\},
\label{morc}
\ee

Step 2 gives the existence of $C_{i,j}$, ${i,j} \in B$, only depending  on $\O$, such that

\[
\norm{u-m_{i,j}}{L^2(\O_i)}^2 \le C_{i,j} \snorm{u}{1,\Mesh}^2, \; \forall (i,j)
\in B.
\]

Then, one has $(m_i-m_{i,j})^2 \mathrm{m}(\O_i)  \le 2 (C_i+C_{i,j})
\snorm{u}{1,\Mesh}^2$, for all $(i,j) \in B$.
Since $\O$ is connected, the above inequality  yields the existence of 
 $M$, only depending  on $\O$, such that $\vert m_i-m_j\vert \le M \snorm{u}{1,\Mesh}$
for all $(i,j) \in \{1,\ldots,n\}^2$, and therefore $\vert m_\O(u) -m_i \vert \le M
\snorm{u}{1,\Mesh}$ for all $i \in \{1,\ldots,n\}$. Then,   \eqref{morc}
yields the existence of $D$, only depending  on $\O$, such that \eqref{poineff}
holds. This completes the proof of Lemma \ref{thinpoin1}.
\eof

An easy consequence of the proof of Lemma \ref{thinpoin1} is the 
following lemma. Although this lemma is not  used in the sequel, it is 
interesting in its own sake.

\begin{lemma}[Mean boundary Poincar\'e inequality]\label{thinpoin2}
\index{Poincar\'e inequality ! - Neumann boundary conditions}
Let $\O$ be an open bounded polygonal connected subset of $\R^d$, $d=2$ or $3$. Let $I \subset \dr \O$ such that the ($\dim-1$)- Lebesgue measure of $I$ is positive. 
Then, there exists $C \in \Rp$, only depending  on $\O$ and $I$, such that for all admissible mesh (in the sense of Definition \refp{meshneuman}) $\Mesh$ and for all $u \in X(\Mesh)$ (see Definition \refp{Xmesh}), the following inequality holds:
\[
 \Vert u - m_I(u) \Vert_{L^2(\O)}^2 \le C \vert u  \vert_{1,\Mesh}^2 
\]
where $\vert \cdot \vert_{1,\Mesh}$ is the discrete $H^1$ seminorm defined in Definition \ref{defh1d} and $m_I(u)$ is the mean value of $\ug$ on $I$ with $\ug$ defined a.e. on $\dr \O$ by $\ug(x)=u_\cv$ if $x \in \edge$, $\edge \in \edgesext \cap  \edgescv$, $\cv \in \mesh$.
\end{lemma}

Note that this last lemma also gives as a by-product a discrete Poincar\'e inequality in the case of a Dirichlet boundary condition on a part of the boundary if 
the domain is assumed to be connex, see Remark \ref{dirich-part}. 
\index{Poincar\'e inequality!-- Dirichlet boundary conditions}

Finally, let us point out that a continuous version of lemmata \ref{thinpoin1} (known as the Poincar\'e-Wirtinger inequality)
and \ref{thinpoin2} holds and that the proof is similar and rather easier. 
Let us state this continuous version which can be proved by contradiction or with a technique similar to Lemma \refp{thuh1}. 
The advantage of the latter is that it gives a more explicit bound.

\begin{lemma}\label{thinpoincont}
Let $\O$ be an open bounded polygonal connected
subset of $\R^d$, $d=2$ or $3$. Let $I \subset \dr \O$ such that the 
($\dim-1$)- Lebesgue measure of $I$ is positive. 

Then, there exist $C \in \Rp$, only depending  on $\O$,  and  $\tilde C \in \Rp$, only depending  on $\O$ and $I$, such that, for all $u \in H^1(\O)$, the following inequalities hold:

\[
\Vert u  \Vert_{L^2(\O)}^2 \le C \vert u \vert_{H^1(\O)}^2  +  2 (\meas(\O))^{-1} (\int_\O u (x) dx)^2
\]

and 

\[
\Vert u  - m_I(u) \Vert_{L^2(\O)}^2 \le \tilde C \vert u     \vert_{H^1(\O)}^2,
\]

where $\vert \cdot  \vert_{H^1(\O)}$ is the $H^1$ seminorm defined by 
$\vert v  \vert_{H^1(\O)}^2 = \Vert \grad u  \Vert_{(L^2(\O))^d}^2 = \int_\O \vert \grad v(x) \vert^2 dx$ for all 
$v \in H^1(\O)$,  and $m_I(u)$ is the mean value of $\ug$ on $I$.
Recall that $\gammat$ is the trace operator from $H^1(\O)$ to
$H^\halft(\dr \O)$.
\end{lemma}

\index{Discrete!-- Poincar\'e inequality|)}\index{Poincar\'e!Discrete -- inequality|)}

\subsection{Error estimate}

Under Assumption \ref{nllH}, let $\Mesh$ be an  admissible mesh (see Definition
\ref{meshneuman}) and $\{f_\cv$, $\cv \in \mesh\}$, $\{g_\cv$,
$\cv \in \mesh\}$ defined by \eqref{nllcondsou}, \eqref{nllcondbor}. By 
Lemma \ref{nllexistu}, there
exists a unique solution
$(u_\cv)_{\cv\in \mesh}$ to \eqref{nllschema}-\eqref{nllschemaO}.
Under an additional regularity assumption on the exact solution, the 
following error
estimate holds:  

\begin{theo}\label{nllesterr}\index{Error estimate!-- for an elliptic 
equation!-- in the general case}
Under Assumption \refp{nllH}, let $\Mesh$ be an
admissible mesh (see Definition \refp{meshneuman}) and $h=\size(\Mesh)$. Let
$(u_\cv)_{\cv \in
\mesh}$ be the unique solution to \eqref{nllschema} and \eqref{nllschemaO}
(thanks to  \eqref{nllcondsou} and \eqref{nllcondbor}, existence and uniqueness of
$(u_\cv)_{\cv \in \mesh}$ is given in  Lemma \ref{nllexistu}). Let $u_\mesh \in
X(\Mesh)$ (see Definition \refp{Xmesh}) be defined by $u_\mesh (x) = u_\cv$ for
a.e. $x \in \cv$, for all $\cv \in \mesh$. 
Assume that the unique solution, $u$, to Problem 
(\ref{nllspace}), (\ref{nllsolvar}) satisfies $u\in C^2({\overline \O})$.

Then there exists $C \in \Rp$ which only depends on $u$ and $\O$ 
such that

\be
\Vert u_\mesh - u \Vert_{L^2(\O)} \le Ch,
\label{nqesterr}
\ee
 
\be
\sumedgesintkl \medge \dedge (\frac {u_\cvv - u_\cv} {\dedge} -
\frac 1 {\medge} \int_{\edge} \grad u(x) \cdot \ncvedge \dfrontiere(x))^2  
\le Ch^2.
\label{nqesterr1}\ee

\end{theo}

Recall that, in the above theorem, $\edgecvcvv$ denotes 
the element $\edge$ of
$\edgesint$ such that $\overline \edge=\dr \cv \cap \dr \cvv$, 
with $\cv$, $\cvv \in \mesh$.

\medskip

\proof of Theorem \ref{nllesterr}

\smallskip

Let $C_\mesh \in \R$ be such that 

\[
\summesh \ubar(\centercv) \mcv=0,
\]

where $\ubar=u+C_\mesh$.

Let, for each $\cv \in \mesh$,  $e_\cv = \ubar(\centercv) -
u_\cv$, and $e_\Mesh \in X(\Mesh)$ defined by $e_\Mesh (x) = e_\cv$ for a.e.
$x \in \cv$, for all $\cv \in \mesh$. Let us first prove
 the existence of $C$ only depending  on $u$ and $\O$ such that

\be
\snorm{e_\mesh}{1,\mesh} \le Ch \text{~~and~~} \norm{e_\mesh}{L^2(\O)} \le Ch.
\label{nllestpre}
\ee

Integrating \refep{nlldifstaf} over  $\cv\in\mesh$, and taking 
 \refep{nlldifstab}  into account yields:

\be
\sumedgecv \int_{\edge} \grad u(x) \cdot \ncvedge \dfrontiere(x) = \int_\cv f(x) dx
+ 
\int_{\dr\cv\cap\dr\O} g(x) \dfrontiere(x).
\label{nller1}\ee

For $\edge \in \edgesint$ such that $\edge=\edgecvcvv$, let us define the 
consistency error \index{Consistency!-- error} on the flux from $\cv$ through $\edge$ by:

\be
\left.\begin{array}{ll}
\Rcvedge = \dsp{\frac 1 {\medge} \int_{\edge}} \grad
u (x) \cdot \ncvedge \dfrontiere (x) - \frac {u(\centercvv) - u(\centercv)}
{\dedge}.
\end{array}\right.
\label{nller2}
\ee

Note that the definition of $\Rcvedge$ remains with $\ubar$ instead of $u$
in \eqref{nller2}.

Thanks to the regularity of the solution $u$,  there exists $C_1 \in \Rp$,
only depending on $u$, such that $\vert R_{\cv,\cvv}\vert \le C_1 h$. Using
\eqref{nller1}, \eqref{nller2} and \eqref{nllschema} yields

\[
\sumE \transcvcvv (e_\cvv - e_\cv)^2  \le  d \meas(\O) (C_1 h)^2,
\]

which gives the first part of \eqref{nllestpre}.

Thanks to the discrete
Poincar\'e inequality \eqref{ellinpoin1}
 applied to the function $e_\mesh$, 
and since
\[\summeshs\mcv e_\cv=0\]  (which is the reason
why  $e_\mesh$ was defined with $\ubar$ instead of $u$) 
one obtains  the second part of \eqref{nllestpre}, that is the
existence of $C_2$ only depending  on $u$ and $\O$ such that

\[
\summesh \mcv (e_\cv)^2  \le  C_2 h^2.
\]

From \eqref{nllestpre}, one deduces \eqref{nqesterr} from the fact that $u \in
C^1(\Obar)$. Indeed, let $C_2$ be the maximum value of $\vert\grad u\vert$ in
$\O$. One has $\vert u(x) - u(y) \vert \le C_2 h$, for all $x$, $y \in \cv$, for
all $\cv \in \mesh$. Then, from $\int_\O u(x) dx = 0$, one deduces $C_\mesh \le
C_2h$. Furthermore, one has

\[
\summesh  \int_\cv (u(x_\cv) - u(x))^2 dx \le \summesh \mcv
(C_2 h)^2 = \meas(\O) (C_2 h)^2.
\]

Then, noting that

\[
\bac
\Vert u_\mesh - u \Vert_{L^2(\O)}^2  = \summesh  \int_\cv (u_\cv - u(x))^2 dx   \\
 \le 3 \summesh \mcv (e_\cv)^2 + 3 (C_\mesh)^2\mesO + 3 \summesh  \int_\cv (u(x_\cv)
- u(x))^2 dx
\ean
\]

yields \eqref{nqesterr}.

\separ

The proof of Estimate \eqref{nqesterr1} is exactly the same as in the Dirichlet
case. This property will be useful in the study of 
the convergence of finite volume
methods in the case of a system consisting 
of an elliptic equation and a hyperbolic equation (see
Section \ref{casimplec}).
\eof

As for the Dirichlet problem, the hypothesis $u \in C^2(\Obar)$ is not necessary
to obtain error estimates. Assuming an additional assumption on the mesh (see
Definition \ref{nram}), Estimates \eqref{nllestpre} and \eqref{nqesterr1} hold under
the weaker assumption $u \in H^2(\O)$ (see Theorem \ref{ntesth2} below). It 
is therefore also possible to obtain \eqref{nqesterr} under the additional 
assumption  that $u$ is
Lipschitz continuous.

\begin{defi}[Neumann restricted admissible meshes]\label{nram}  \normalfont
\index{Mesh!Admissible --!Restricted -- for Neumann boundary conditions}
Let $\O$ be an open bounded polygonal connected subset of $\R^d$, $d = 2$ or $3$. 
A restricted admissible mesh for the Neumann problem, denoted by $\mesh$, is an admissible mesh in the sense of Definition \ref{meshneuman} such that, for some $\zeta>0$, one has  $\dcvedge \ge
\zeta\diam(\cv)$ for all control volume $\cv$ and for all $\edge \in \edgescv \cap \edgesint$.
\end{defi}

\begin{theo}[$H^2$ regularity, Neumann problem] 
\label{ntesth2}\index{Error estimate!-- for an elliptic 
equation!-- in the general case}
Under Assumption \refp{nllH}, let $\mesh$ be an admissible
mesh in the sense of Definition \ref{nram} and $h=\size(\Mesh)$. Let $u_\mesh \in
X(\Mesh)$ (see Definition \refp{Xmesh}) be the approximated solution defined in
$\O$ by $u_\mesh (x)  = u_\cv$ for a.e. $x \in \cv$, for all $\cv \in
\mesh$, where $(u_\cv)_{\cv \in \mesh}$ is the (unique) solution to
\eqref{nllschema} and \eqref{nllschemaO} (thanks to
 \eqref{nllcondsou} and \eqref{nllcondbor}, existence and uniqueness of
$(u_\cv)_{\cv \in \mesh}$ is given in Lemma \ref{nllexistu}). 
Assume that the unique solution, $u$, of (\ref{nllspace}), (\ref{nllsolvar})
belongs to $H^2({\O})$. Let $C_\mesh \in \R$ be such that 
\[\summeshs
\ubar(\centercv) \mcv=0\text{ where } \ubar=u+C_\mesh.\]
 Let, for each control volume 
$\cv \in \mesh$,  $e_\cv = \ubar(\centercv) - u_\cv$, and $e_\Mesh \in X(\Mesh)$ defined by
$e_\Mesh (x) = e_\cv$ for a.e.
$x \in \cv$, for all $\cv \in \mesh$. 

Then there exists $C$, only depending  on $u$, $\zeta$ and $\O$, such that
\eqref{nllestpre} and \eqref{nqesterr1} hold.
\end{theo}

Note that, in  Theorem \ref{ntesth2}, the function
$e_\mesh$ is well defined, and the quantity ``$\grad u \cdot \nedge$" is
well defined on $\edge$, for all $\edge \in \edges$ (see Remark \ref{uish2}).

\medskip

\proof of Theorem \ref{ntesth2}

\smallskip

The proof is very similar to that of Theorem \refp{testh2}, from which the same
notations are used.

There exists some $C$, depending only
on the space dimension ($d$) and
$\zeta$ (given in Definition
\ref{nram}), such that, for all $\edge \in
\edgesint$,

\be
\vert R_\edge \vert^2 \le C \frac {h^2} {\medge \dedge}
\int_{\V_\edge} \vert (H(u)(z) \vert^2 dz,
\label{nconsh2}
\ee

and therefore 

\be
\sumedgesint \medge \dedge R_\edge^2 \le Ch^2 \int_{\O}\vert
H(u)(z)\vert^2 dz.
\label{nbouh}
\ee

The proof of \eqref{nconsh2} (from which \eqref{nbouh} is an easy 
consequence) was already done in
the proof of Theorem \ref{testh2} (note that, here, 
there is no need to consider the
case of $\edge \in
\edgesext$).
In order to obtain Estimate \eqref{nllestpre}, one proceeds as in
Theorem \ref{testh2}. Recall 

\[
\snorm{e_{\mesh}}{1,\mesh}^2 \le \sumedgesint R_\edge \vert D_\edge e \vert 
\medge,
\]
where
$\vert D_\edge e\vert =\vert e_\cv - e_\cvv \vert$ if $\edge \in
\edgesint$ is such that 
$\edge=\edgecvcvv$; hence, from the Cauchy-Schwarz  inequality, one obtains that

\[
\snorm{e_{\mesh}}{1,\mesh}^2 \le \big(\sumedgesint R_\edge^2 \medge
\dedge\big)^{\frac 1 2} \big(\sumedgesint \vert D_\edge e \vert^2 \frac {\medge}
{\dedge} \big)^{\frac 1 2}.
\]

Then, one obtains, with \eqref{nbouh},

\[
\snorm{e_{\mesh}}{1,\mesh} \le \sqrt {C} h \big( \int_{\O} \vert
H(u)(z)\vert^2 dz \big)^{\frac 1 2}.
\]

This concludes the proof of the first part of \eqref{nllestpre}. The second part of
\eqref{nllestpre} is a consequence of the discrete
Poincar\'e inequality \eqref{ellinpoin1}. Using \eqref{nbouh} also 
easily leads \eqref{nqesterr1}.

Note also that, if $u$ is Lipschitz continuous, Inequality \eqref{nqesterr} follows
from the second part of
\eqref{nllestpre} and the definition of $\ubar$ as in Theorem \ref{nllesterr}.

This
concludes the proof of Theorem
\ref{ntesth2}.
\eof

Some generalizations of Theorem \ref{ntesth2} are possible, as for the Dirichlet
case, see Remark \refp{h2plus}.

\subsection{Convergence}

A convergence result, under Assumption \ref{nllH}, may be proved
without any  regularity assumption on the exact solution.

The  proof of  convergence uses the following   preliminary inequality
on the ``trace" of an element of $X(\Mesh)$  on the boundary:

\begin{lemma}[Trace inequality]\label{nlltrace}
Let $\O$ be an open bounded polygonal connected subset of $\R^d$, $d=2$ or 
$3$ (indeed, the connexity of $\O$ is not used in this lemma). Let 
$\mesh$ be an admissible mesh, in the sense of Definition 
\refp{meshneuman}, and $u \in X(\Mesh)$ (see Definition \refp{Xmesh}). 
Let $u_\cv$ be the value of $u$ in the control volume $\cv$. Let 
$\gammat(u)$ be defined by $\gammat(u)=u_\cv$ a.e. (for the 
$(d-1)$-dimensional Lebesgue measure) on $\edge$, if $\edge \in \edgesext$ and $\edge 
\in \edgescv$. Then, there exists $C$, only depending on $\O$, such that

\be
\norm{\ug}{L^2(\dr \O)} \le C(\snorm{u}{1,\Mesh}+\norm{u}{L^2(\O)}).
\label{nlltraceineq}
\ee

\end{lemma}
\begin{remark} \normalfont
The result stated in this lemma still holds if $\O$ is not assumed 
connected. Indeed, one needs only modify (in an obvious way) the definition 
of admissible meshes (Definition \refp{meshneuman}) so as to take into 
account non connected subsets.
\end{remark}

\proof of Lemma \ref{nlltrace}

\smallskip

By compactness of the boundary of $\dr \O$, there exists a finite
number of open hyper-rectangles ($d=2$ or $3$), $\{R_i, i = 1,\ldots,N\}$, and
normalized vectors of $\R^\dim$, $\{\eta_i, i = 1,\ldots,N\}$, such that 
\[
\ba \dsp{ \dr \O
\subset \cup_{i=1}^N R_i,} \\
\dsp{ \eta_i \cdot \n(x) \ge \a > 0 }\text{ for all } x \in R_i \cap \dr \O, 
i \in
\{1,\ldots,N\}, \\
\{x + t \eta_i, x
\in R_i \cap \dr
\O, t \in \Rp\} \cap R_i \subset \O, 
\ea
\]
where $\alpha$ is some positive
number and $\n(x)$ is the normal vector  to
$\dr \Omega$ at $x$,  inward to $\Omega$.
Let $\{\alpha _i, i = 1,\ldots,N\}$ be a family of functions such
that $\sum_{i=1}^N \a_i(x) = 1$, for all $x \in \dr \Omega$,
$\a_i \in C^\infty_c(\R^\dim,\R_+)$ and $\a_i = 0$ outside of
$R_i$, for all $i = 1, \ldots,N$. 
Let $\Gamma_i = R_i \cap \dr \Omega$; let us prove that there
exists $C_i$ only depending  on $\alpha$ and $\a_i$ such that

\be 
\Vert \a_i \ug \Vert_{L^2(\Gamma_i)} \le C_i\bigl(\vert u
\vert_{1,\mesh} + \Vert u \Vert_{L^2(\O)}\bigr). \label{trace1}
\ee

The existence of $C$, only depending on $\O$, such that
\eqref{nlltraceineq} holds, follows easily (taking $C = \sum_{i=1}^N C_i$, and using 
$\sum_{i=1}^N \a_i(x) = 1$, note that $\alpha$ and $\a_i$ depend only on 
$\O$). It remains to prove \eqref{trace1}.

~

Let us introduce some notations. 
For $\edge \in \edges$ and $\cv \in \mesh$, define $\chi_\edge$ and 
$\chi_\cv$ from $\R^\dim \times \R^\dim$ to $\{0,1\}$ by 
$\chi_\edge(x,y)=1$, if $[x,y] \cap \edge \ne \emptyset$, 
$\chi_\edge(x,y)=0$, if $[x,y] \cap \edge = \emptyset$, and 
$\chi_\cv(x,y)=1$, if $[x,y] \cap \cv \ne \emptyset$, 
$\chi_\cv(x,y)=0$, if $[x,y] \cap \cv = \emptyset$.

~

Let $ i \in
\{1,\ldots,N\}$ and let $x \in \Gamma_i$. There exists a 
unique $t>0$ such that $x+t\eta_i 
\in \dr R_i$, let $y(x) = x + t \eta_i$.
For $\edge \in \edges$, let $z_\edge(x)=[x,y(x)] \cap \edge$ if $[x,y(x)] 
\cap \edge \ne \emptyset$ and is reduced to one point. For $\cv \in 
\mesh$, let $\xi_\cv(x),\eta_\cv(x)$ be such that $[x,y(x)] \cap 
\cv=[\xi_\cv(x),\eta_\cv(x)]$ if $[x,y(x)] \cap 
\cv \ne \emptyset$.

One has, for a.e. (for the 
$(d-1)$-dimensional Lebesgue measure) $x \in \Gamma_i$,

\[
\vert \a_i \ug(x) \vert \le \sum_{\edge=\edgecvcvv \in \edgesint} 
\vert  \a_i(z_{\edge}(x))
 (u_{\cv} -  u_\cvv) \vert
\chi_\edge(x,y(x)) + \sum_{\cv \in \mesh} \vert 
(\a_i(\xi_\cv(x)-\a_i(\eta_\cv(x))u_\cv \vert \chi_\cv(x,y(x)),
\]

that is,

\be
\vert \a_i \ug (x) \vert^2  \le  A(x) + B(x)
\label{trace2}
\ee

with

\[
A(x) =  2(\sum_{\edge=\edgecvcvv \in \edgesint} 
\vert  \a_i(z_{\edge}(x))
 (u_{\cv} -  u_\cvv) \vert
\chi_\edge(x,y(x)))^2 ,
\]

\[
B(x) = 2(\sum_{\cv \in \mesh} \vert 
(\a_i(\xi_\cv(x))-\a_i(\eta_\cv(x)))u_\cv \vert \chi_\cv(x,y(x)))^2.
\]

A bound  on $A(x)$ is obtained  for a.e. $x \in \Gamma_i$, by remarking
 that, from the Cauchy-Schwarz   inequality:

\[
A(x) \le D_1 \sumedgesint \frac {\vert D_\edge u \vert^2} {\dedge 
c_\edge} \chi_\edge(x,y(x)) \sumedgesint \dedge c_\edge \chi_\edge(x,y(x)),
\]

where 
$D_1$ only depends on $\a_i$ and  $c_\edge = \vert \eta_i \cdot \nedge \vert$. (Recall that $D_\edge 
u=\vert u_\cv-u_\cvv \vert$.)
Since 
\[\sumedgesint \dedge \cedge \chi_\edge(x,y(x)) \le \diamO,\]
this yields:
\[ 
A(x)  \le \diamO D_1 \sumedgesint \frac {\vert D_\edge u \vert^2} {\dedge 
c_\edge} \chi_\edge(x,y(x)).
\]

Then, since
\[\int_{\Gamma_i} \chi_\edge (x,y(x)) \dfrontiere (x) \le \frac 1 \a 
\cedge \medge,\]
 there exists $D_2$, only depending  on $\O$, such that

\[
A= \int_{\Gamma_i}A(x) \dfrontiere (x) \le D_2 \snorm{u}{1,\Mesh}^2.
\]

A  bound $B(x)$ for a.e. $x \in \Gamma_i$ is obtained with the Cauchy-Schwarz  
inequality:

\[
B(x) \le D_3 \sum_{\cv \in \mesh} u_\cv^2 \chi_\cv(x,y(x)) \vert 
\xi_\cv(x)-\eta_\cv(x) \vert \sum_{\cv \in \mesh}\vert 
\xi_\cv(x)-\eta_\cv(x) \vert \chi_\cv(x,y(x)),
\]

 where
$D_3$ only depends on $\a_i$. Since 
\[\sum_{\cv \in \mesh}\vert 
\xi_\cv(x)-\eta_\cv(x) \vert \chi_\cv(x,y(x)) \le \diamO
\text{ and }
\int_{\Gamma_i} 
\chi_\cv (x,y(x)) \vert \xi_\cv(x)-\eta_\cv(x) \vert \dfrontiere (x)
 \le \frac 1 \a 
\mcv,
\]
 there exists $D_4$, only depending  on $\O$, such that

\[
B=\int_{\Gamma_i}B(x) \dfrontiere (x) \le D_4 \norm{u}{L^2(\O)}^2.
\]

Integrating \eqref{trace2} over $\Gamma_i$, the bounds on $A$ and $B$ lead 
\eqref{trace1} for some convenient $C_i$ and it concludes the proof of
Lemma \ref{nlltrace}.
\eof

\begin{remark} \normalfont\label{poincontr}
Using this ``trace inequality" \eqref{nlltraceineq} and the Kolmogorov theorem (see 
Theorem \refp{Kolm}, it
is possible to prove Lemma \refp{thinpoin1} (Discrete Poincar\'e 
inequality)\index{Discrete!-- Poincar\'e inequality} \index{Poincar\'e!Discrete -- inequality} by way
of contradiction. Indeed, assume that there exists a sequence $(u_n)_{n\in \N}$ such
that, for all
$n\in\N$, $\Vert u_n\Vert_{L^2(\O)} = 1$,
$\int_\O u_n(x) dx = 0$, $u_n \in X(\Mesh_n)$ (where $\Mesh_n$ is an admissible
mesh in the sense of Definition
\ref{meshneuman}) and $\vert u_n \vert_{1,\Mesh_n}
\le \frac 1 n$.
 Using the trace inequality, one proves that $(u_n)_{n\in \N}$ is relatively compact in $L^2(\O)$, as in Theorem \refp{nllconv}.  \index{Compactness!-- of the approximate solutions in $L^2$}
Then, one can assume that $u_n \tends u$ in $L^2(\O)$  as $n \tends \infty$. 
The function $u$ satisfies $\Vert u\Vert_{L^2(\O)} = 1$, since $\Vert u_n\Vert_{L^2(\O)} = 1$, and $\int_\O u(x) dx = 0$, since $\int_\O u_n(x) dx = 0$. 
Using $\vert u_n \vert_{1,\Mesh_n} \le \frac 1 n$, a proof similar to that of Theorem \refp{kolmh1}, yields that $D_iu=0$, for all $i
\in \{1,\ldots,n\}$ (even if $\size(\mesh_n) \not \tends 0$, as $n \tends  \infty$), where $D_iu$ is the derivative in the distribution sense
with respect to $x_i$ of $u$. 
Since $\O$ is connected, one deduces that $u$ is constant on $\O$, but this is impossible since $\Vert u\Vert_{L^2(\O)} = 1$ and $\int_\O u(x) dx = 0$.
\end{remark}

\separ

Let us now prove that the scheme  \eqref{nllschema} and \eqref{nllschemaO}, where 
$(f_\cv)_{\cv \in \mesh}$ and $(g_\cv)_{\cv \in \mesh}$ are given by 
\eqref{nllcondsou} and \eqref{nllcondbor} is stable\index{Stability}:  
the approximate solution given by the scheme is bounded independently of 
the
mesh, as we proceed to show.

\begin{lemma}[Estimate for the Neumann problem]\label{nllest}
Under Assumption \refp{nllH}, let $\mesh$ be an admissible
mesh (in the sense of Definition \refp{meshneuman}). Let
$(u_\cv)_{\cv \in
\mesh}$ be the unique solution to \eqref{nllschema} and \eqref{nllschemaO}, where 
$(f_\cv)_{\cv \in \mesh}$ and $(g_\cv)_{\cv \in \mesh}$ are given by 
\eqref{nllcondsou} and \eqref{nllcondbor};  the  
existence and uniqueness of
$(u_\cv)_{\cv \in \mesh}$ is given in Lemma \ref{nllexistu}. Let $u_\mesh \in
X(\Mesh)$ (see Definition \ref{Xmesh}) be defined by $u_\mesh (x) = u_\cv$ for
a.e. $x \in \cv$, for all $\cv \in \mesh$. Then, there exists $C \in \Rp$, only
depending on $\O$, $g$ and $f$, such that

\be
\snorm{u_\Mesh}{1,\Mesh} \le C,
\label{nllestimate}
\ee

\index{Estimate on the approximate solution!-- for an elliptic equation}
where $\snorm{\cdot}{1,\Mesh}$ is defined in Definition \refp{defh1d}.

\end{lemma}

\proof of Lemma \ref{nllest}

\smallskip

Multiplying \eqref{nllschema} by $u_\cv$ and summing  over $\cv\in\mesh$ yields

\be
\sumE\transcvcvv (u_\cvv - u_\cv)^2 = \summesh \mcv f_\cv u_\cv
+ \sumedgesext  u_\cvedge g_\cvedge \medge,
\label{nqesterr2a}\ee

where, for $\edge \in \edgesext$, $\cvedge \in \mesh$ is such that $\edge \in
\edgescvedge$.

 We get \eqref{nllestimate} from \eqref{nqesterr2a} using
\eqref{nlltraceineq}, 
\eqref{ellinpoin1} and  the Cauchy-Schwarz   inequality.
\eof

Using the estimate \eqref{nllestimate} on the approximate solution, a convergence
result is given in the following theorem.

\begin{theo}[Convergence in the case of the Neumann problem]\label{nllconv}~

\index{Convergence of the finite volume method!-- for an elliptic equation}
Under Assumption \refp{nllH}, let $u$ be the unique solution
to \eqref{nllspace},\eqref{nllsolvar}. For an admissible mesh (in the sense of
Definition \refp{meshneuman}) $\mesh$, let
$(u_\cv)_{\cv \in
\mesh}$ be the unique solution to \eqref{nllschema} and \eqref{nllschemaO} (where 
$(f_\cv)_{\cv \in \mesh}$ and $(g_\cv)_{\cv \in \mesh}$ are given by 
\eqref{nllcondsou} and \eqref{nllcondbor}, the  
existence and uniqueness of
$(u_\cv)_{\cv \in \mesh}$ is given in Lemma \ref{nllexistu}) and define $u_\mesh
\in X(\Mesh)$ (see Definition \ref{Xmesh}) by $u_\mesh (x) = u_\cv$ for
a.e. $x \in \cv$, for all $\cv \in \mesh$. Then,

\[
u_\Mesh \tends u \text{ in } L^2(\O) \text{ as } \size(\Mesh) \tends 0,
\]

\[
\snorm{u_\mesh}{1,\mesh}^2 \tends \int_\O \vert \grad u (x) \vert^2 dx \text{ as }
\size(\Mesh) \tends 0 
\]

and
\[
\gammat (u_\Mesh) \tends \gammat (u) \text{ in } L^2(\dr \O) \text{ for the 
weak topology}  \text{ as } 
\size(\Mesh) \tends 0,
\]
where the function $\gammat(u)$ stands for the trace of $u$
on $\dr \O$ in the sense given in Lemma \ref{nlltrace} when $u \in X(\mesh)$ and in the
sense of the classical trace operator from $H^1(\O)$ to $L^2(\dr \O)$ (or $H^{\half}(\dr
\O)$) when $u \in H^1(\O)$.
\end{theo}

\proof of Theorem \ref{nllconv}

\smallskip

\textit{Step 1 (Compactness)}

 \index{Compactness!-- of the approximate solutions in $L^2$}
 Denote by $Y$ the set of approximate solutions $u_{\mesh}$ for all 
admisible meshes
$\mesh$. Thanks to Lemma \ref{nllest} and to
the discrete Poincar\'e inequality  \eqref{ellinpoin1}, the set $Y$ is bounded in
$L^2(\O)$.
Let us prove that $Y$ is relatively compact in $L^2(\O)$, and that, if $(\mesh_n)_\nnn$ is a  sequence of admissible meshes such that $\size(\mesh_n)$ tends to $0$ and $u_{\mesh_n}$ tends to $u$, in $L^2(\O)$, as $n$ tends to infinity, then $u$ belongs to $H^1(\O)$. Indeed, these results follow from theorems \ref{Kolm} and \refp{kolmh1}, provided that there exists a  real positive number $C$  only depending  on $\O$, $f$ and $g$ such that
\be
\norm{\tilde u_{\mesh}(\cdot+\eta)-\tilde u_{\mesh}}{L^2(\R^d)}^2 \le C \vert \eta \vert,
\text{ for any  admissible mesh  } \mesh   \text{ and  for any  }  \eta \in \R^d,
\; \vert \eta \vert \le 1,
\label{cvneum1}
\ee
and that, for any compact subset $\comp$ of $\O$, 
\be
\ban 
\norm{{u_{\mesh}}(\cdot+\eta)-{u_{\mesh}}}{L^2(\comp)}^2 \le C \vert \eta \vert
(\vert \eta \vert + 2 \, \size(\mesh)),
\text{ for any  admissible mesh  } \mesh   & & \\
\text{ and  for any  }  \eta \in \R^d \text{ 
such that  }  \vert \eta \vert < d(\comp, \O^c). & & 
\ean
\label{cvneum2}
\ee
Recall that $\tilde u_{\mesh}$ is defined by $\tilde u_{\mesh}(x)  =  u_{\mesh}(x)$ if $x \in \O$ and
$\tilde u_{\mesh}(x)  = 0$ otherwise. 
In order to prove
\eqref{cvneum1} and \eqref{cvneum2}, define
$\chi_{\edge}$ from $\R^\dim \times \R^\dim$ to $\{0,1\}$ by 
$\chi_{\edge}(x,y)=1$ if
$[x,y] \cap \edge \neq \emptyset$ and $\chi_{\edge}(x,y)=0$ if
$[x,y] \cap \edge = \emptyset$.
Let $\eta \in  \R^d \setminus
\{0\}$. Then:

\be
\vert \tilde u(x+\eta)-\tilde u(x) \vert \le \sumedgesint \chi_{\edge}(x,x+\eta)
\vert D_{\edge}u \vert + \sumedgesext \chi_{\edge}(x,x+\eta)
\vert u_\edge \vert, \; \text{ for } \ae x \in \O,
\label{cvneum3}
\ee

where, for $\edge \in \edgesext$, $u_\edge = u_{\cv}$, and $\cv$ is the
control volume such that $\edge \in \edgescv$. Recall also that $D_\edge u=\vert u_\cv-u_\cvv \vert$,
if $\edge=\edgecvcvv$.
Let us first prove Inequality \eqref{cvneum2}. Let $\comp$ be a compact subset of $\O$. 
If  $x  \in \comp$ and $\vert \eta  \vert < d(\comp,\O^c)$, the second term of the right
hand side of  \eqref{cvneum3} is 0, and the same proof as in Lemma 
\refp{h1dtot} gives, from an integration over $\comp$ instead of $\O$ and from 
\eqref{aie} with $C=2$ since $[x,x+\eta] \subset \O$ for $x \in \comp$,
 
\be
\norm{{u_{\mesh}}(\cdot+\eta)-{u_{\mesh}}}{L^2(\comp)}^2 \le \vert u \vert^2_{1, \mesh}
\vert \eta \vert
(\vert \eta \vert + 2 \, \size(\mesh)).
\label{preppar}
\ee

\separ

In order to prove \eqref{cvneum1}, remark  that the number of non zero terms in the 
second term of the right hand side of \eqref{cvneum3} is, for a.e. $x \in \O$, bounded by 
some real positive number, which only depends on $\O$, which can be taken, for 
instance, as the number of sides of $\O$, denoted by $N$.  
Hence, with $C_1=(N+1)^2$ (which only depends on $\O$. Indeed, if $\O$ is convex, $N=2$ is 
also convenient), one has

\be
\vert \tilde u(x+\eta)-\tilde u(x) \vert ^2 \le C_1 (\sumedgesint 
\chi_{\edge}(x,x+\eta)
\vert D_{\edge}u \vert )^2+ C_1\sumedgesext \chi_{\edge}(x,x+\eta)
u_\edge^2, \; \text{ for } \ae x \in \O.
\label{cvneum4}
\ee

Let us integrate this inequality over $\R^d$. As seen in the proof of 
Lemma \refp{h1dtot},

\[
\int_{\R^d} \bigl( \sumedgesint \chi_{\edge}(x,x+\eta)
\vert D_{\edge}u \vert \bigr)^2   dx \le \vert u \vert_{1,\mesh}^2 \vert 
\eta\vert (\vert 
\eta\vert  + 2(N-1) \size(\mesh));
\]

hence, by Lemma \refp{nllest}, there exists a real positive
number $C_2$, only depending on $\O$, $f$ and $g$, such that (if $\vert \eta \vert
\le 1$)

\[
\int_{\R^d} \bigl( \sumedgesint \chi_{\edge}(x,x+\eta)
\vert D_{\edge}u \vert \bigr)^2   dx \le C_2 \vert  \eta\vert.
\]

Let us now turn to the second term of the right hand side of
\eqref{cvneum4} integrated over $\R^d$; 

\[ 
\ban
\dsp{ \int_{\R^d} \bigl(  \sumedgesext \chi_{\edge}(x,x+\eta)
 u_\edge ^2  \bigr)  dx  }& \le  & 
 \sumedgesext \medge  \vert  \eta \vert  u_\edge ^2 \\
  & \le  & \norm{\gammat(u_{\mesh})}{L^2(\dr \O)}^2  \vert \eta \vert; \ean
\]

therefore, thanks to Lemma \ref{nlltrace}, Lemma \ref{nllest} and to the discrete
Poincar\'e inequality \eqref{ellinpoin1}, 
there exists a real positive number
$C_3$, only depending on 
$\O$, $f$ and $g$, such that

\[
\int_{\R^d} \bigl( \sumedgesext \chi_{\edge}(x,x+\eta)
 u_\edge ^2  \bigr)  dx   \le   C_3  \vert  \eta \vert.
\]

Hence \eqref{cvneum1} is proved for some real positive number $C$ only depending 
on $\O$, $f$ and $g$.

\separ

\textit{Step 2 (Passage to the limit)}

In this step, the convergence of $u_\mesh$ to the solution of \eqref{nllspace},
\eqref{nllsolvar} (in $L^2(\O)$ as $\size(\mesh) \tends 0$) is first proved.

Since the solution to \eqref{nllspace}, \eqref{nllsolvar} is unique, and thanks to
the compactness of the set $Y$ described in Step 1, it is sufficient to prove that, if 
$u_{\mesh_n} \tends u$ in $L^2(\O)$ and $\size(\mesh_n) \tends 0$ as $n \tends 0$, 
then $u$ is a solution to  \eqref{nllspace}-\eqref{nllsolvar}.

Let $(\mesh_n)_\nnn$ be a sequence of admissible meshes and $(u_{\mesh_n })_\nnn$ 
be the corresponding solutions to \eqref{nllschema}-\refep{nllschemaO} with $\mesh =
\mesh_n$. Assume $u_{\mesh_n} \tends u $ in 
$L^2(\O)$ and $\size (\mesh_n) \tends 0$ as $n \tends 0$. By Step 1,  one has $u \in 
H^1(\O)$ and since the mean value of $u_{\mesh_n}$ is zero, one also has $\int_\O 
u(x) dx = 0$. Therefore, $u$ is a solution of \eqref{nllspace}. It remains to show that $u$ 
satisfies \eqref{nllsolvar}.  Since $(\gammat(u_{\mesh_n}))_\nnn$ is bounded in 
$L^2(\partial \O)$, one may assume (up to a subsequence) that it converges to some $v$ weakly in
$L^2(\dr \O)$.  Let us first prove that

\be\bari
\dsp{- \int_\O u  (x) \Delta \varphi(x) dx + \int_{\dr \O}  \nabla \varphi (x)
\cdot \mathbf{n} (x)  v(x) \dfrontiere(x)  = \int_\O    f(x) \varphi (x) dx}
\\
\dsp{+ \int_{\dr\O} g(x)  \varphi (x) \dfrontiere(x), \, \forall \varphi \in
C^2(\Obar),}
\ean
\label{cvneum5}
\ee
and then that $u$ satisfies  \eqref{nllsolvar}.

\separ 

Let $\mesh$ be an admissible mesh,  $u_{\mesh}$ the corresponding approximate solution to 
the Neumann problem, given by  \eqref{nllschema} and \eqref{nllschemaO}, 
where $(f_\cv)_{\cv \in \mesh}$ and $(g_\cv)_{\cv \in \mesh}$ are given by \eqref{nllcondsou} and \eqref{nllcondbor} and let $\varphi \in 
C^2(\Obar)$. Let $\varphi_\cv = \varphi(x_\cv)$, define $\varphi_\mesh$ by 
$\varphi_\mesh(x) = \varphi_\cv$,  for a.e. $x \in \cv$ and for any control volume
$\cv$, and $\gammat (\varphi_\mesh)(x) = \varphi_\cv$ for a.e. $x \in \edge$ (for
the $(d-1)$-dimensional Lebegue measure), for any $\edge \in \edgesext$ and  control volume
$\cv$ such that 
$\edge \in \edgescv$.

Multiplying  \eqref{nllschema}  by $\varphi_\cv$, summing over $\cv \in \mesh$ and 
reordering the terms yields

\be
\summesh u_\cv \sumv \transcvcvv (\varphi_\cvv - \varphi_\cv)  = \int_\O   f(x) 
\varphi_\mesh(x) dx 
+ \int_{\dr\O} \gammat(\varphi_\mesh)(x) g(x) \dfrontiere(x).
\label{cvneum6}\ee

Using the consistency  of the fluxes and the fact that $\varphi \in C^2(\Obar)$, there exists 
$C$ only depending on $\varphi$ such that 

\[
 \sumv \transcvcvv (\varphi_\cvv - \varphi_\cv)  =  \int_\cv \Delta \varphi (x) dx - 
\int_{\dr \O \cap \dr \cv} \nabla \varphi(x) \cdot \mathbf{n} (x)  \dfrontiere (x) +  
\sumv  R_{\cv,\cvv}(\varphi),
\]

with  $ R_{\cv,\cvv} = - R_{\cvv,\cv}$, for all $\cvv \in \Nv$ and $\cv \in \mesh$, and 
$\vert  R_{\cv,\cvv}  \vert \le C_4 \medgecvcvv \size(\mesh)$, where   $C_4$ 
only depends on $\varphi$. Hence \eqref{cvneum6} may be rewritten 
as

\be
\bari
\dsp{- \int_\O u_{\mesh} (x) \Delta \varphi (x) dx + \int_{\dr \O}  \nabla
 \varphi (x) \cdot \mathbf{n} (x)  
\gammat(u_{\mesh})(x) \dfrontiere (x) + r(\varphi, \mesh)  =} \\
\dsp{\int_\O   
f(x) \varphi_\mesh(x) dx 
+ \int_{\dr\O} \gammat(\varphi_\mesh)(x) g (x) \dfrontiere(x),}
\ean
\label{cvneum7}
\ee

where
\[
\ban
\vert r(\varphi, \mesh)  \vert  & = & C_4  \sumedgesint  \vert D_\edge u \vert 
\medge \size(\mesh) 
\\
& \le & \dsp{ C_4  \bigl(\sumedgesint  \vert D_\edge u \vert^2 \frac 
{\medge}{\dedge }\bigr)^{\half}   \bigl(\sumedgesint   {\medge}{\dedge
}\bigr)^{\half} \size(\mesh) }
\\ 
& \le & C_5  \size(\mesh),
\ean
\]

where $C_5$  is a real positive number  only depending on $f$, $g$, $\O$ and $\varphi$ 
(thanks to Lemma \ref{nllest}).

Writing \eqref{cvneum7} with $\mesh = \mesh_n$ and passing to the limit as $n$ tends to 
infinity yields \eqref{cvneum5}.

\separ

Let us now prove that $u$ satifies  
\eqref{nllsolvar}. Since $u \in H^1(\O)$, an integration by parts in \eqref{cvneum5} yields

\be
\bari
\dsp{\int_\O \nabla u  (x) \cdot \nabla \varphi(x) dx + \int_{\dr \O} 
 \nabla \varphi(x) \cdot \mathbf{n }(x)  (v(x)-  \gammat(u) (x)) \dfrontiere(x)}
 \\
 = \dsp{\int_\O f(x) \varphi(x) dx 
+ \int_{\dr\O}g (x)  \varphi (x) \dfrontiere(x),} \, \forall \varphi \in C^2(\Obar),
\ean
\label{cvneum9}
\ee
where $\gammat(u)$ denotes the trace of $u$ on $\dr \O$ (which belongs to $L^2(\dr
\O)$). In order to prove that $u$ is solution to 
\eqref{nllsolvar} (this will conclude the proof of Theorem \ref{nllconv}), 
it is sufficient, thanks to the density of $C^2(\Obar)$ in 
$H^1(\O)$, to prove that $v = \gammat(u)$ a.e. on $\dr \O$ (for the 
$(d-1)$ dimensional Lebesgue measure on $\dr \O$). Let us now  prove that $v = \gammat(u)$
 a.e. on $\dr \O$ by first 
remarking that \eqref{cvneum9} yields

\[
\int_\O \nabla u  (x) \cdot \nabla \varphi(x) dx   = \int_\O   f(x) 
\varphi (x) dx, \, \forall \varphi \in C_c^\infty(\O),
\]

and therefore, by density of $C_c^\infty(\O)$ in $H^1_0(\O)$, 

\[
\int_\O \nabla u  (x) \cdot \nabla \varphi(x) dx   = \int_\O   f(x) 
\varphi (x) dx, \, \forall \varphi \in H^1_0(\O).
\]

With \eqref{cvneum9}, this yields

\be
- \int_{\dr \O} \nabla \varphi(x) \cdot \n (x)  (v(x)- \gammat(u) (x))
\dfrontiere (x) = 0, \, \forall \varphi \in C^2(\Obar) \text{ such that }  \varphi = 0 \text{ on }  \dr \O.
\label{cvneum10}
\ee

There  remains to show that the wide choice of $\varphi$ in \eqref{cvneum10}
 allows to conclude $v =
\gammat(u)$  a.e. on $\dr \O$ (for the $(d-1)$-dimensional Lebesgue measure of $\dr 
\O$). Indeed, let 
$I$ be a part of the boundary $\dr \O$, such that $I$ is included in a hyperplane of $\R^\dim$.
Assume that $I = \{0\} \times J$, where $J$ is an open ball of
$\R^{d-1}$ centered on the origin. Let $z = (a,\tilde z) \in \R^d$ with 
$a \in \R^\star_+$, $\tilde z \in \R^{\dim -
1}$ and $B = \{(t,\frac {a - \vert t \vert} a y + \frac{\vert t \vert}a 
\tilde z); t \in
(-a,a), y \in J\}$; 
assume that, for a convenient $a$, one has

\[
B \cap \O =  \{(t,\frac {a - \vert t \vert} a y + \frac{\vert t \vert}a 
\tilde z); t \in
(0,a), y \in J\}.
\] 

Let $\psi \in C^\infty_c(J)$, and for $x = (x_1,y) \in \R\times J$, 
define $\varphi_1(x) = -x_1 \psi(y)$. Then, 
\[\varphi_1 \in C^\infty(\R^\dim) \text{ and }
\dsurdn{\varphi_1} = \psi \text{ on }I.\]
(Recall that $\n$ is the normal unit vector to
$\dO$, outward to $\O$.) Let $\varphi_2 \in C_c^\infty(B)$ such that 
$\varphi_2 = 1$
on a neighborhood of $\{0\} \times \{\psi \ne 0\}$, where $\{\psi \ne 
0\}=\{x \in J$; $\psi(x) \ne 0\}$, and set
$\varphi = \varphi_1
\varphi_2$;  $\varphi$ is an admissible test function in \eqref{cvneum10}, and
therefore 
\[\int_J \psi(y) \bigl(\gammat(u)(0,y)-v(0,y)\bigr) dy = 0,\]  
which yields, since $\psi$ is arbitrary in $C^\infty_c(J)$, 
$v=\gammat(u)$ a.e. on
$I$. Since $J$ is arbitrary, this implies  that
 $v=\gammat(u)$ a.e. on $\dO$.  

\smallskip

This conclude the proof of $u_\mesh \tends u$ in $L^2(\O)$ as $\sizemesh \tends 0$, where $u$
is the solution to \eqref{nllspace},\eqref{nllsolvar}.

\smallskip

Note also that the above proof gives (by way of contradiction) that
$\gammat(u_\mesh) \tends \gammat(u)$ weakly in $L^2(\dr \O)$, as $\size(\mesh) \tends 0$.

Then, a passage to the limit in \eqref{nqesterr2a} together with \eqref{nllsolvar} yields

\[
\snorm{u_\mesh}{1,\mesh}^2 \tends \norm{\vert \grad u
\vert}{L^2(\O)}^2, \text{ as } \size(\mesh) \tends 0.
\]

This concludes
the proof of Theorem  \ref{nllconv}.
\eof

Note that, with some discrete Sobolev inequality (similar to 
\eqref{ellinlq}), the hypothesis ``$f \in L^2(\O)$ $g \in L^2(\dr\O)$" may 
be relaxed in some way similar to that of Item 2 of Remark \ref{convplus}.

\index{Neumann boundary conditions|)}

\section{General elliptic operators} \label{ell2Dgen}

\hpage{Elliptic multid}

\subsection{Discontinuous matrix diffusion coefficients}

\index{Discontinuous coefficients}
\subsubsection{Meshes and schemes}

Let $\O$ be an open bounded polygonal subset of
$\R^\dim$, $\dim=2$ or $3$.  We are interested here in the discretization 
of an elliptic operator with discontinuous matrix diffusion coefficients,
which may appear in real case problems such as electrical or thermal transfer
problems or, more generally, diffusion problems in heterogeneous media. 
In this case, the mesh is adapted to fit the discontinuities of the data.
Hence the definition of an admissible mesh given in Definition \ref{meshdirichlet}
must be adapted. As an illustration, let us consider here the following
problem, which was studied in Section \refp{ell1D} in the one-dimensional 
case\index{Equation!Convection-diffusion --}\index{Discontinuous 
coefficients}: 

\be
 -\div (\La \grad  u) (x) + \div (\vbf u) (x) + b u (x)
= f(x) , \; x \in \O, \label{roueneq}
\ee
\be
 u(x)=g(x), \;  x \in \dr \O,
\label{rouencl}
\ee

with the following assumptions on the data (one denotes by $\R^{\dim 
\times \dim}$ the set of $\dim 
\times \dim$ matrices with real coefficients):

\begin{hyp} \label{hypdiscont}~

\begin{enumerate}
\item $\La$ is a bounded measurable function from $\O$ to $\R^{\dim 
\times \dim}$ such that for any  $x \in \O$, $\La(x)$
is symmetric,  and that there exists $\lund$ and $\lbar \in \R^\star_+$ such
that  $\lund  \xi \cdot \xi \le \La(x) \xi \cdot \xi \le \lbar \xi \cdot \xi$ for
any $x \in \O$ and any $\xi \in \R^\dim$.
\item $\vbf \in C^1(\Obar,\R^\dim)$, $\div \vbf \ge 0$ on $\O$, $b \in \R_+$.
\item $f$ is a bounded piecewise continuous function from $\O$ to $\R$.
\item $g$ is such that there exists $\tilde g \in
H^1(\O)$ such that $\gammat(\tilde g) = g$ (a.e. on $\dr \O$) and is a 
bounded piecewise continuous function from $\partial \O$ to $\R$.
\end{enumerate}
\end{hyp}

(Recall that $\gammat$ denotes 
the trace operator from $H^1(\O)$ into $L^2(\dr \O)$.) 
As in Section \ref{sdbc}, under Assumption \ref{hypdiscont}, there exists 
a unique variational solution $u\in H^1(\O)$ of Problem  \eqref{roueneq}, 
\eqref{rouencl}. This solution satisfies $u = w + \tilde g$, where $\tilde 
g \in H^1(\O)$ is such that $\gammat(\tilde g) = g$,  a.e. on $\dr \O$, 
and $w$ is the unique function of $H^1_0(\O)$ satisfying

\[
\bac
\dsp{\int_\O \Bigl( \La(x) \grad w(x) \cdot \grad \psi(x) + \div (\vbf w) (x) 
\psi(x)  +b w(x) \psi(x) \Bigr) dx = }\\

\dsp{ \int_\O \Bigl(- \La(x) \grad \tilde  g(x) \cdot \grad \psi(x)
- \div (\vbf  \tilde  g) (x)  \psi(x) -b \tilde g(x) \psi(x)  + f(x)
\psi(x)\Bigr) dx, \;     \forall \psi\in H^1_0(\O).}
\ean
\]

Let us now define an admissible mesh  for the discretization  of Problem
\eqref{roueneq}-\eqref{rouencl}. 

\begin{defi}[Admissible mesh for a general diffusion operator]  
\label{rouenmesh}  \normalfont
\index{Mesh!Admissible --!--for a general diffusion operator}
 
Let $\O$ be an open bounded  polygonal subset of $\R^d$, $d = 2$ or $3$. An
admissible finite volume mesh for the discretization of Problem 
\eqref{roueneq}-\eqref{rouencl} is an admissible mesh $\Mesh$ of
$\O$ in the sense of Definition \refp{meshdirichlet} where items $(iv)$ and $(v)$
are 
replaced by the two following conditions:

\begin{enumerate}

\item[(iv)'] The set $\Mesh$ is such that 

  the restriction of $g$
 to each edge $\edge \in \edgesext$ is continuous.

  For any $\cv \in \mesh$, let $\La_\cv$ denote the mean value of
$\La$ on $\cv$, that is 
\[
\La_\cv = \frac 1 {\mcv} \int_\cv \La(x)dx.
\] 
There exists a family of points 
\[
\points = (\centercv)_{\cv  \in \mesh} \text{
such that  } \centercv =
\cap_{\edge \in \edgescv} \D_{\cv,\edge} \in \cvbar,
\]
where
$\D_{\cv,\edge}$ is a straigth line perpendicular to $\edge$  with
respect to the scalar product induced by $\La_\cv^{-1}$   such that 
$\D_{\cv,\edge} \cap \edge = \D_{\cvv,\edge}  \cap \edge \ne \emptyset$  if $\edge =
\edgecvcvv$. Furthermore, if $\edge =
\edgecvcvv$, let $\centeredge=\D_{\cv,\edge} \cap \edge (= 
\D_{\cvv,\edge}  \cap \edge)$ and   assume that $\centercv \ne 
\centercvv$.

\item[(v)'] For any $\edge \in \edgesext$, let 
$\cv$ be the control volume  such
that $\edge \in \edgescv$ and let  $\Dcvedge$ be the straight line going through
$\centercv$  and orthogonal to $\edge$ with respect to the scalar product induced by
$\La_\cv^{-1}$;  then, there exists $\centeredge \in \edge \cap \Dcvedge$; let $g_\edge =
g(\centeredge)$.
\end{enumerate}


The  notations are are the same as those introduced in 
Definition  \refp{meshdirichlet}.

\end{defi}

We shall now define the discrete unknowns of  the numerical scheme, with the
same notations as in Section \ref{vfquatre}. As in the case of the
Dirichlet problem, the  primary unknowns
$(u_\cv)_{\cv \in \mesh}$ will be used, which aim to be approximations of the values
$u(\centercv)$, and some auxiliary unknowns,
 namely the fluxes $F_{\cv,\edge}$,
for all $\cv \in \mesh$ and $\edge \in \edgescv$, and some (expected)
approximation of $u$ in $\edge$, say $u_{\edge}$, for all $\edge \in \edges$.
Again,  these auxiliary unknowns are helpful to write the
scheme, but they can be eliminated locally so that the discrete 
equations will only  be written with respect to the primary unknowns
$(u_\cv)_{\cv \in \mesh}$. For any $\edge \in \edgesext$, set
$u_\edge =   g(\centeredge)$. 
The finite volume scheme for the numerical approximation of the solution to
Problem  (\ref{roueneq})-(\ref{rouencl}) is obtained by integrating Equation
(\ref{roueneq}) over each control volume $\cv$, and  approximating the fluxes 
over each edge $\edge$ of $\cv$. This yields\index{Finite volume scheme!-- for 
elliptic problems!-- in two or three 
space dimensions|(}  
\be
  \sumedgecv \Fcvedge  + \sumedgecv \vcvedge \uedgep + \mcv b
u_\cv =  f_\cv, \, \forall \cv \in \mesh,  \label{rouenvfn} 
\ee
 where

  $\vcvedge = \int_\edge \vbf(x) \cdot 
\n_{\cv,\edge}\dfrontiere(x)$  (where $\n_{\cv,\edge}$ denotes the normal unit vector to
$\edge$ outward to $\cv$);  if  $\edge =  {\cvedgep\vert\cvedgem}$, $\uedgep  =  u_{\cvedgep}$, where $\cvedgep$ is the 
upstream\index{Upstream}
control volume, i.e. $\vcvedge \ge 0$, with $\cv= \cvedgep$; if   $\edge
\in \edgesext$, then  $\uedgep  =  u_{\cv}$ if  $\vcvedge \ge 0$ (i.e. 
$\cv$ is upstream to $\edge$ with respect to $v$), and   $\uedgep  =  u_{\edge}$ otherwise.

 $\Fcvedge$ is an approximation of 
$\int_\edge - \La_\cv \nabla u(x)\cdot \n_{\cv,\edge} \dfrontiere(x)$; 
the approximation $\Fcvedge$ is  written with respect to the discrete
unknowns  $(u_\cv)_{\cv \in \mesh}$ and  $(u_\edge)_{\edge \in \edges}$.  For
$\cv \in \mesh$ and $\edge \in \edgescv$, let  $\l_{\cv,\edge} = \vert
\La_\cv  \n_{\cv,\edge} \vert$ (recall that $\vert \cdot \vert$  denote the Euclidean norm).

\begin{itemize}
\item If $\centercv \not \in  \edge$, a natural expression
for $\Fcvedge $ is then  
\[
\Fcvedge  = - \medge  \l_{\cv,\edge} \frac {{u_\edge - u_\cv}} {\dcvedge}.
\] 
Writing the
conservativity  of the scheme, i.e.   $F_{\cvv,\edge} = -
F_{\cv,\edge}$ if $\edge = \edgecvcvv \subset \O$, yields the value 
of  $u_\edge$, if $\centercvv \notin \edge$, with respect to  $(u_\cv)_{\cv \in \mesh}$; 
\[ u_\edge = \dsp{ 
\frac  1 
{ 
\frac { {  \l_{\cv,\edge}  }  }  {\dcvedge} +  \frac {  
{\l_{\cvv,\edge} }  }  {\dcvvedge}   
}
\bigl(  \frac {{\l_{\cv,\edge}}} {\dcvedge}   u_\cv + \frac
{{\l_{\cvv,\edge}}}  {\dcvvedge} u_\cvv  \bigr).
}   
\] 
Note that this expression is similar to that of \refep{uipud} in the 1D 
case.
\item If $\centercv  \in \edge$, one sets $u_\edge = u_\cv.$  
\end{itemize}

Hence the value of $\Fcvedge$;

\begin{itemize}
\item internal edges:
\be
\Fcvedge  = - \transedge (u_{\cvv} - u_\cv), 
\text{ if } \edge \in \edgesint, \, \edge=\edgecvcvv,
\label{rouenfluxni}
\ee
where 
\[
 \transedge  = \medge \frac { {\l_{\cv,\edge} \l_{\cvv,\edge}} }  {
{\l_{\cv,\edge} \dcvvedge + \l_{\cvv,\edge} \dcvedge} } \text{ if } \centeredge \ne
\centercv \text{ and } \centeredge \ne \centercvv  \]
and 
\[
 \transedge  = \medge \frac  {{ \l_{\cv,\edge}  }} {\dcvedge}   \text{
if } \centeredge \ne \centercv \text{ and } \centeredge = \centercvv;  \]

\item boundary edges:

\be
\Fcvedge  = - \transedge (g_\edge -
u_\cv),  \text{ if } \edge \in \edgesext \text{ and } \centercv \not \in 
\edge,
\label{rouenfluxnf} 
\ee

where

\[
\transedge= \medge  \frac {\l_{\cv,\edge}}{\dcvedge};
\]

if $\centercv \in \edge$, then the equation associated to $u_\cv$ is $u_\cv =
g_\edge$ (instead of that given by \eqref{rouenvfn}) and the numerical 
flux $\Fcvedge$ is an unknown which may be deduced from 
\eqref{rouenvfn}\index{Finite volume scheme!-- for 
elliptic problems!-- in two or three 
space dimensions|)}.
\end{itemize}

\begin{remark} \normalfont  Note that if  
$\La =
Id$, then  the scheme  (\ref{rouenvfn})-(\ref{rouenfluxnf}) is the same scheme
than the one described in Section \ref{vfquatre}.  
 \end{remark} 

\subsubsection{Error estimate}

\begin{theo}~  \label{rouenesterr}\index{Error estimate!-- for an elliptic 
equation!-- in the general case}

Let $\O$ be an open bounded polygonal 
subset of $\R^\dim$, $\dim = 2$ or $3$. Under Assumption \ref{hypdiscont}, 
let $u$ be the unique variational solution 
to Problem \eqref{roueneq}-\eqref{rouencl}. Let $\Mesh$ be an
admissible mesh for the discretization of Problem
\eqref{roueneq}-\eqref{rouencl}, in the sense of Definition \ref{rouenmesh}. 
Let $\zeta_1$ and $\zeta_2 \in \R_+$ such that 

\[
\ban
 \zeta_1(\size(\mesh))^2 \le \mcv
\le \zeta_2 (\size(\mesh))^2,  \\
\zeta_1 \size(\Mesh)\le \medge 
\le \zeta_2 \size(\mesh), \\
 \label{vfquatreH2bis}
\zeta_1 \size(\Mesh)\le \dedge 
\le \zeta_2 \size(\mesh).
\ean
\]

Assuming moreover that

   the restriction of $f$ to $\cv$ belongs to $C(\cvbar)$, 
for any $\cv \in \mesh$;

  the restriction of $\La$ to $\cv$ belongs to 
$C^1(\cvbar,\R^{\dim \times \dim})$, for any $\cv \in \mesh$;  

  the restriction of $u$ (unique variational solution of
Problem \eqref{roueneq}-\eqref{rouencl}) to $\cv$ belongs to $C^2(\cvbar)$, for any $\cv \in \mesh$.

(Recall that $C^m(\cvbar,\R^N)=\{v_{\vert_\cv},$ $v \in 
C^m(\R^\dim,\R^N)\}$ and $C^m(\cdot)=C^m(\cdot,\R)$.)

\smallskip

Then, there exists a unique family $(u_\cv)_{\cv \in \mesh}$ satisfying
(\ref{rouenvfn})-(\ref{rouenfluxnf});  furthermore, denoting by $\errorcv =
u(\centercv) - u_\cv$, there exists $C \in \R_+$ only depending on $\zeta_1,
\zeta_2,$ $\gamma = \sup_{\cv \in \mesh}  (\Vert D^2 u \Vert_{L^\infty(\cv)})$  and 
$\delta = \sup_{\cv \in \mesh}$  $ (\Vert D \La \Vert_{L^\infty(\cv)})$ such
that 

\be 
\sumedges \frac {(D_\edge \error)^2} {\dedge }
\medge  \le  C (\size(\Mesh))^2 \label{rouenesthuz}
\ee

and

\be 
\summesh \errorcv^2 \mcv  \le C (\size(\Mesh))^2.
\label{rouenestdeux}
\ee

Recall that $D_\edge \error = \vert \error_\cvv - \error_\cv \vert$ for $\edge \in
\edgesint$, $\edge=\edgecvcvv$ and  
 $D_\edge \error = \vert  \error_\cv \vert$ for
 $\edge \in \edgesext \cap \edgescv$.

\end{theo}

\proof of Theorem \ref{rouenesterr}

\smallskip

First, one may use Taylor expansions and the same technique as in the 1D case (see step 2 of the proof of Theorem \ref{ee1dell}, Section \ref{ell1D}) to show that the expressions (\ref{rouenfluxni}) and (\ref{rouenfluxnf}) are consistent approximations of th exact diffusion flux $\int_\edge - \La(x)  \grad u(x) \cdot \ncvedge \dfrontiere(x)$, i.e. there exists $C_1$ only depending on $u$ and $\La$ such that, for all $\edge \in \edges$, with $\Fastcvedge = \transedge (u(\centercvv) - u(\centercv))$, if $\edge = \edgecvcvv$, and $\Fastcvedge = \transedge (u(\centeredge)-u(\centercv))$, if $\edge \in \edgesext \cap \edgescv$,

\[
\bal
\Fastcvedge - \int_\edge -\La(x) \grad u(x) \cdot \ncvedge \dfrontiere(x)  =  \Rcvedge, \\
\text{ with } \vert \Rcvedge  \vert  \le  C_1 \size(\Mesh)\medge.
\ean 
\]

There also exists $C_2$ only depending on $u$ and $\vbf$ such that, for all 
$\edge \in \edges$,

\[
\bal
\vcvedge u(x_{\cv_{\edge,+}}) - \int_\edge \vbf \cdot \ncvedge u =  \rcvedge, \\
\text{ with } \vert \rcvedge  \vert   \le  C_2 \sizemesh \medge.
\ean 
\]

Let us then integrate 
Equation (\ref{roueneq}) over each control volume,  subtract  to (\ref{rouenvfn}) and 
use the consistency of the fluxes  to obtain  the following equation on the error:

\[
\ba
- \sumedgecv \Gcvedge + \sumedgecv \vcvedge \eedgep + \mcv b
e_\cv =  & & \\ 
  \sumedgecv (\Rcvedge +  \rcvedge) + S_\cv, \,
\forall \cv \in \mesh, 
\ea
\]

where  $\Gcvedge = \transedge  (e_\cvv - e_\cv)$, if $\edge =
\edgecvcvv$, and $\Gcvedge = \transedge  (- e_\cv)$, if $\edge \in 
\edgesext \cap \edgescv$,  $\eedgep = e_{\cv_{\edge,+}}$ 
is the error associated to the
upstream control volume to $\edge$ and $S_\cv = b (\mcv u(\centercv)  - 
\int_\cv u(x) dx)$ is such that $\vert S_\cv \vert \le \mcv C_3 h$, where $C_3 \in \R_+$ only depends on $u$ and $b$. 
Then, similarly to the proof of Theorem \refp{ellesterr}, let us multiply  by  $\error_\cv$, sum  over $\cv \in \mesh$, and use the conservativity  of the scheme, which yields that if $\edge = \edgecvcvv$ then  $\Rcvedge = - \Rcvvedge$. 
A reordering of  the summation over $\edge \in \edges$ yields the ``discrete  $H^1_0$ estimate" (\ref{rouenesthuz}). 
Then, following \citet{VF4}, one shows the following discrete Poincar\'e inequality:\index{Discrete!-- Poincar\'e inequality}\index{Poincar\'e!Discrete -- inequality}
\be
\summesh  e_\cv^2 \mcv \le C_4 \sumedges  \frac {(D_\edge \error)^2}
{\dedge } \medge, 
\label{ellinpoindisc}
\ee
where $C_4$ only depends on $\O$, $\zeta_1$ and $\zeta_2$, which in turn  yields the $L^2$ estimate (\ref{rouenestdeux}).
\eof

\begin{remark} \normalfont In the case where $\La$ is constant, or more generally, in the case where $\La(x) = \l(x) Id$, where $\l(x) >0$, the proof of Lemma \ref{thinpoin} is
easily extended. However, for a general matrix $\Lambda$, the generalization of this proof is not so clear; this is the reason of the dependency of the estimates \eqref{rouenesthuz} and \eqref{rouenestdeux} on $\zeta_1$ and $\zeta_2$,   which arises when proving   \eqref{ellinpoindisc} as in \citet{VF4}.
\end{remark}

\subsection{Other boundary conditions} \label{obt}
The finite volume scheme may be used to discretize elliptic
 problems with Dirichlet or Neumann boundary conditions, as  we saw in the
previous sections. It is also easily implemented in the case of Fourier (or
Robin) and periodic boundary conditions. The case of interface conditions
between two geometrical regions is also generally easy to implement; the
purpose here is to present the treatment of some of these boundary and
interface conditions. One  may also refer to \citet{Angot} and references
therein, \citet{FH1} for the treatment of more complex boundary conditions and
coupling terms in a system of elliptic equations.

Let $\O$ be (for the sake of simplicity) the open rectangular subset of $\R^2$
defined by $\O= (0,1) \times (0,2)$,  let   $\O_1= (0,1) \times (0,1)$, $\O_2=
(0,1) \times (1,2)$, $\Gamma_1 = [0,1] \times \{0\}$,  $\Gamma_2 =  \{1\}  \times
[0,2]$, $\Gamma_3 =  [0,1] \times  \{2\}$, $\Gamma_4 = \{0\} \times 
[0,2]$ and  $I = [0,1] \times \{1\}$. Let $\l_1$ and $\l_2 >0$, $f 
\in C(\Obar)$,  $\a >0$,   $\ubar \in \R$,
$g \in C(\Gamma_4)$, $\theta$ and $\Phi \in C(I)$. Consider here the following
problem (with some ``natural" notations):  

\be
 -\div (\l_i \grad u) (x)  
= f(x), \; x \in \O_i, \, i = 1,2,  \label{pbeq}
\ee
\be
- \l_i \nabla u(x) \cdot \n(x)  = \a (u(x) - \ubar), \, x \in \Gamma_1 \cup \Gamma_3, 
\label{pbFourier}
\ee
\be
 \nabla u(x) \cdot \n(x)  = 0, \, x \in \Gamma_2,
\label{pbNeumann}
\ee
\be
u(x) = g(x), \, x \in \Gamma_4, 
\label{pbDirichlet}
\ee
\be
(\l_2 \nabla u(x) \cdot \n_I(x))_{\vert 2}   = (\l_1 \nabla u(x) 
\cdot \n_I(x))_{\vert 1} +
\theta(x), \, x \in I,  \label{pbjumpflux}
\ee
\be
u_{\vert 2} (x)   -  u_{\vert 1}(x)   =
\Phi(x), \, x \in I,  \label{pbjumpun}
\ee

where $\n$ denotes the unit normal vector to $\dr \O$ outward to $\O$ and 
$\n_I=(0,1)^t$ (it is a unit normal vector to $I$).

Let $\Mesh$ be an admissible mesh for the
discretization of \eqref{pbeq}-\eqref{pbjumpun} in the sense of
Definition \ref{rouenmesh}. For the sake of simplicity, let us assume
here that $\dcvedge >0$ for all $\cv \in \mesh, \edge \in \edgescv$. 
Integrating
Equation (\ref{pbeq}) over each control volume $\cv$, and  approximating the
fluxes  over each edge $\edge$ of $\cv$  yields the following finite volume
scheme\index{Finite volume scheme!-- for 
elliptic problems!-- in two or three space dimensions|(}:    

\be 
\sumedgecv \Fcvedge  
  =  f_\cv, \, \forall \cv \in \mesh, 
\label{pbvfn}  
\ee
where 
$\Fcvedge$ is an approximation of 
$\int_\edge - \l_i \nabla u(x)\cdot \n_{\cv,\edge} \dfrontiere(x)$, with $i$ such that $\cv
\subset \Omega_i$.  

\smallskip

Let $N_\mesh = \card(\mesh)$,  $N_\edges = \card(\edges)$, $N^0_\edges = 
\card (\{\edge \in \edges; \edge \not \subset
\dr \O \cup I\})$, $N^i_\edges = \card (\{\edge \in \edges; \edge  \subset
 \Gamma_i\})$,  and  $N^I_\edges = \card (\{\edge \in \edges; \edge  \subset
 I\})$ (note that $N_\edges = N^0_\edges + \sum_{i=1}^4  N^i_\edges + 
N^I_\edges$).  Introduce the $N_\mesh$ (primary) discrete unknowns 
$(u_\cv)_{\cv \in \mesh}$; note that the number of (auxiliary) unknowns 
 of the
type  $ \Fcvedge $ is $2 (N^0_\edges + N^I_\edges) + \sum_{i=1}^4  N^i_\edges$;
let us introduce the discrete unknowns $(u_\edge)_{\edge \in \edges}$, which
aim to be approximations of $u$ on $\edge$. In order to take into account the
jump condition \eqref{pbjumpun}, two unknowns of this type  are necessary on
the edges $\edge \subset I$, namely $u_{\edge,1}$ and  $u_{\edge,2}$. Hence
the number of (auxiliary) unknowns of the type $u_\edge$ is $ N^0_\edges  +
\sum_{i=1}^4  N^i_\edges +  2 N^I_\edges$. 
Therefore, the total number of discrete unknowns is

\[
N_{tot} = N_\mesh +  3
N^0_\edges + 4 N^I_\edges  + 2 \sum_{i=1}^4  N^i_\edges.
\]

Hence, it is convenient, in order to obtain a well-posed system, to write $N_{tot}$ 
discrete equations. We already have $N_\mesh$ equations from \eqref{pbvfn}.
The expression of $\Fcvedge$ with respect to the unknowns $u_{\cv}$ and $u_\edge$ is

\be
\Fcvedge  = - \medge  \l_i \frac {{u_\edge - u_\cv}} {\dcvedge}, \, \forall \,
\cv \in \mesh; \cv \subset \O_i \, (i = 1, 2), \, \forall \, \edge \in \edgescv;
 \label{pbvfnflux} 
\ee 

which yields $2(N^0_\edges + N^I_\edges)  +   \sum_{i=1}^4  N^i_\edges $. 
(In \eqref{pbvfnflux}, $u_\edge$ stands for $u_{\edge,i}$ if $\edge 
\subset I$.)

Let us now take into
account  the various boundary and
interface conditions: 

\begin{itemize}
\item \underline {Fourier boundary conditions.} Discretizing
condition \eqref{pbFourier} yields
\be
\Fcvedge  = \a \medge ({{u_\edge - \ubar}}),  \forall \, \cv \in \mesh, 
\forall \, \edge \in \edgescv; \, \edge \subset \Gamma_1 \cup \Gamma_3,
\label{pbvfnFourier} 
\ee  
that is $ N^1_\edges +  N^3_\edges $ equations. 

\item \underline {Neumann boundary conditions.} Discretizing condition
\eqref{pbNeumann} yields 

\be
\Fcvedge  = 0,  \forall \, \cv \in \mesh,
\forall \, \edge \in \edgescv; \, \edge \subset \Gamma_2 ,
\label{pbvfnNeumann} 
\ee  
that is $N^2_\edges$ equations. 

\item \underline {Dirichlet boundary conditions.} Discretizing
condition \eqref{pbDirichlet} yields
\be
u_\edge  = g(\centeredge),  \, \forall \, \edge \in \edges; \edge \subset
\Gamma_4,  
\label{pbvfnDirichlet} 
\ee  

that is $ N^4_\edges$ equations.

\item \underline {Conservativity \index{Conservativity} of the flux.}
Except at interface $I$, the flux is continuous, and therefore
\be
\Fcvedge  = -\Fcvvedge ,  \forall \, \edge \in \edges; \edge \not \subset
(\bigcup_{i = 1}^4 \Gamma_i \cup I) \text{ and } \edge = \edgecvcvv,
\label{pbvfnconserv}
\ee  
that is $ N^0_\edges$ equations.

\item \underline {Jump condition on  the flux.}
At interface $I$, condition \eqref{pbjumpflux} is discretized into

\be
\Fcvedge  + \Fcvvedge = \int_\edge \theta (x) ds ,  \forall \, \edge \in \edges; 
\edge
 \subset  I  \text{ and } \edge = \edgecvcvv; \, \cv \subset \O_2,
\label{pbvfnjumpflux} 
\ee  
that is $N^I_\edges$ equations.

\item \underline {Jump condition on the unknown.}
At interface $I$, condition \eqref{pbjumpun} is discretized into
\be
 u_{\edge,2} =  u_{\edge,1}  + \Phi(\centeredge),  \forall \, \edge \in
\edges;  \edge
 \subset  I  \text{ and } \edge = \edgecvcvv.
\label{pbvfnjumpun} \ee  
that is another $N^I_\edges$ equations. 

\end{itemize}

Hence the total number of equations  from \eqref{pbvfn} to \eqref{pbvfnjumpun} is
$N_{tot}$, so that the numerical scheme can be expected to be well 
posed\index{Finite volume scheme!-- for 
elliptic problems!-- in two or three space dimensions|)}.

\medskip

The finite volume scheme for the discretization of equations 
\eqref{pbeq}-\eqref{pbjumpun} is therefore completely defined by 
\eqref{pbvfn}-\eqref{pbvfnjumpun}. Particular cases of this scheme are
the schemes \eqref{ellschema}-\refep{ellschemad} (written for Dirichlet
boundary conditions)  and  
\eqref{nllschema}-\refep{nllschemaO} (written for Neumann boundary 
conditions and no convection term)  which were thoroughly studied in the 
two previous sections.

\section{Dual meshes and unknowns located at vertices} \label{scvfe}
\index{Mesh!Dual --}
\hpage{Elliptic multid}

One of the principles of the classical finite volume method is to
associate the discrete unknowns to the grid cells. However, it is sometimes 
useful to associate  the discrete unknowns with the vertices of the mesh; for instance, 
  the finite
volume method may be used  for the discretization of a hyperbolic equation coupled
with an elliptic equation (see Chapter \ref{systemes}). Suppose that  an existing
finite element code is implemented for the elliptic equation and    yields the discrete
values of the unknown at the vertices of the mesh. One might then want to implement a
finite volume method  for the hyperbolic equation with the values of the unknowns at 
the vertices of the mesh. Note
also that for some physical problems, e.g. the modelling of two phase flow in porous
media,  the conservativity  principle is easier to respect if the discrete
unknowns have the same location.  For these various reasons, we introduce here some
finite volume methods where the discrete unknowns are located at the vertices of an
existing mesh.
 
For the sake of simplicity, the treatment of the  boundary conditions will be omitted
here. Recall that the construction of a finite volume method is carried 
out (in particular) along the
following principles:

\begin{enumerate}
\item Divide the spatial domain in control volumes,
\item Associate to each control volume and, for time dependent problems,  to each
discrete time, one discrete unknown, 
\item Obtain the discrete equations (at each discrete time) by integration of the
equation over the control volume and the definition of one exchange term between two
(adjacent) control volumes.
\end{enumerate}

Recall, in particular, that the definition of one (and one only) exchange term
between two control volumes is important; this is called the property of conservativity
of a finite volume method. 
The aim here is to present finite volume methods for which the discrete 
unknowns are located at the vertices of the mesh. Hence, 
to each vertex must correspond a control volume. Note that these control volumes may be
somehow ``fictive" (see the next section); 
the important issue is to respect the principles given 
above in the construction of the finite volume scheme.  In the three 
following sections, we shall deal with the two dimensional case;
 the generalization to the three-dimensional case is the purpose of section
\ref{cvfe3D}. 

\subsection{The piecewise linear finite element method 
\index{Finite element method|(}
 viewed as a finite volume
method}\label{lfefv}

Links between the finite volume method and the finite element method have already been explored \cite{Sel,ide-94-fin,hug-00-con}.
Here we prove that for the Dirichlet problem, the piecewise linear finite element method may indeed be viewed as a finite volume method.
Let $\O$ be a bounded open polygonal subset of $\R^2$,$f$ and $g$ be some ``regular" functions (from $\O$ or $\partial \O$ to $\R$).
Consider the following problem:

\be 
\left\{\begin{array}{ll}  
- \Delta u (x)= f  (x), & \; x \in  \Omega, \\ 
u(x) = g(x), & \; x \in  \partial \Omega. 
\end{array}\right. 
\label{rouene2D} 
\ee

Let us show that the  ``piecewise linear" finite
element method for the discretization of \eqref{rouene2D} may be viewed as 
a kind of finite volume
method.  Let
$\meshprim$ be a finite element mesh of $\O$, consisting of triangles 
(see e.g. \citet{ciarlet} for the conditions on the triangles), and  let
$\Vertices \subset \Obar $ be the set of vertices of
$\meshprim$. For $\cv \in \Vertices$ (note that here $\cv$ denotes a 
point of $\Obar$),  let $\varphi_\cv$ be the shape function associated
to $\cv$ in the piecewise linear finite element method for 
the mesh $\meshprim$. 
We remark that 

\[
\sum_{\cv \in \Vertices}\varphi_\cv(x)=1, \; \forall x \in  \O,
\]

and therefore

\be
\sum_{\cv \in \Vertices} \int_{\O} \varphi_\cv(x)dx = \mesO
\label{consglo}
\ee

and

\be
\sum_{\cv \in \Vertices}\grad \varphi_\cv(x)=0, \; \text{ for a.e.} x \in  \O.
\label{rouencruu}
\ee

Using the latter equality, the discrete finite element
equation associated to the unknown $u_\cv$, if $\cv \in \O$, 
can therefore be written as

\[
\sum_{\cvv \in \Vertices} \int_{\O}(u_{\cvv}-u_\cv)\grad
\varphi_{\cvv}(x) \cdot \grad \varphi_\cv(x)dx =
 \int_\O f(x) \varphi_\cv (x) dx.
\]

Then the finite element method may be written as

\[
\sum_{\cvv \in \Vertices} - \transcvcvv (u_{\cvv}-u_\cv) =
 \int_\O f(x) \varphi_\cv (x) dx, \text{ if } \cv \in \Vertices \cap \O,
\]
\[
u_\cv=g(\cv), \text{ if } \cv \in \Vertices \cap \dr \O,
\]

with 

\[
\transcvcvv = - \int_{\O}\grad
\varphi_{\cvv}(x) \cdot \grad \varphi_\cv(x)dx.
\]

Under this form, the finite element method may be viewed as a finite 
volume method, except that there are no ``real" control volumes associated to 
the vertices of $\meshprim$. Indeed, thanks to \eqref{consglo}, the 
control volume associated to $\cv$ may be viewed as the support of 
$\varphi_\cv$ ``weighted" by $\varphi_\cv$. This interpretation of 
the finite element method as a finite volume method was also used in 
\citet{F1}, \citet{F2} and \citet{EG} in order to design a numerical 
scheme for a
 transport equation for which the velocity field is the gradient of the
 pressure, which is itself the solution to an elliptic equation (see also 
  \citet{HL} for numerical tests). 
 This method  is often referred to as the "control 
volume finite element" method\index{Control 
volume finite element method}.

\smallskip

In this finite volume interpretation of the finite element scheme, the notion of 
``consistency of the fluxes" does not appear. This notion of consistency, 
however, seems to be an interesting tool in the study of the  ``classical" finite volume schemes.

\smallskip

Note that the (discrete) maximum principle is satisfied with this scheme if
only if the transmissibilities $\transcvcvv$ are nonnegative (for all 
$\cv, \cvv \in \Vertices$ with $\cv \in \O$)\index{Transmissibility} ; this is 
the case  under the classical Delaunay condition\index{Delaunay condition}; this
condition states that the (interior of the) circumscribed circle (or 
sphere  in the three 
dimensional case) of any triangle (tetrahedron in the three 
dimensional case) of the mesh does not contain 
any element of $\Vertices$. This is equivalent, in the case of 
two dimensional triangular meshes, to the fact that the sum of two 
opposite  angles 
facing a common edge is less or equal $\pi$. 
\label{Delaunay}

\subsection{Classical finite volumes on a dual mesh}

Let $\meshprim$ be a mesh of $\O$ ($\meshprim$ may consist of triangles, 
but it is not necessary) and $\Vertices$ be the set of vertices of $\meshprim$. 
In order to
associate to each vertex (of $\meshprim$) a control volume (such that the whole spatial
domain is the ``disjoint union" of the control volumes), a possibility is to
construct a ``dual mesh"  which will be denoted by $\mesh$. In order for 
this mesh to be admissible in the sense of Definition 
\refp{meshdirichlet}, a simple way is to use the 
Vorono\"{\i}\index{Mesh!Vorono\"{\i} --}\index{Vorono\"{\i}} mesh defined with $\Vertices$ (see Example \refp{evoronoi}).  
For a description  of the Delaunay-Vorono\"{\i} discretization  and its use for covolume methods, we refer to \cite{nicolbook} (and references therein).
In order to write the ``classical" finite volume scheme with this mesh (see \eqref{ellschema}-\refep{ellschemad}), a slight modification is 
necessary at the boundary for some particular $\meshprim$ (see Example \ref{evoronoi}); this method is denoted CFV/DM (classical 
finite volume on dual mesh); it is conservative, the numerical fluxes are  consistent, and the transmissibilities \index{Transmissibility}  are nonnegative.  
Hence, the convergence results and error estimates which were studied in previous sections hold (see, in particular, theorems \refp{ellcvgce} and \refp{ellesterr}).

\separ

A case of particular interest is found when the primal mesh (that is 
$\meshprim$) consists in triangles with acute angles. One uses, as 
dual mesh, the Vorono\"{\i} mesh defined with $\Vertices$. Then, the dual 
mesh is admissible in the sense of Definition 
\refp{meshdirichlet} and is constructed with the orthogonal bisectors of the edges of the elements
of $\meshprim$, parts of these orthogonal bisectors (and parts of $\dr 
\O$) give the boundaries to 
the control volumes of the dual mesh.
In this case, the CFV/DM 
scheme is ``close" to the piecewise linear finite element  scheme on the primal 
mesh. Let us elaborate on this point. 

For $\cv \in \Vertices$, let  $\cv$   also denote  the control volume (of the dual mesh) associated to $\cv$ (in the sequel, the   notation ``$\cv$", which denotes  either the vertex or the control volume will be used in such a context that it does not yield any confusion) and let  $\varphi_\cv$  be the shape function associated to the vertex $\cv$ (in the  piecewise linear finite element associated to $\meshprim$).
The term $\transcvcvv$ (ratio between the length of the edge $\edgecvcvv$ and the distance between vertices), which is used in the finite volume scheme, verifies

\[
\transcvcvv = - \int_\O  \grad \varphi_{\cv}(x) \cdot  \grad \varphi_{\cvv}(x)dx. 
\]

\begin{figure}[htb]
\begin{center}
\begin{tikzpicture}[
    scale=0.9,
    tri/.style={line width=0.8pt},
    vor/.style={line width=2.2pt},
    ext/.style={dashed, line width=0.5pt},
    dim/.style={densely dotted, {Straight Barb[open]}-{Straight Barb[open]}},
  ]
  \coordinate (K)  at (0,4);
  \coordinate (L)  at (0.8,-3.6);
  \coordinate (M)  at (5.4,0);
  \coordinate (Mt) at (-5,0.4);
  \coordinate (xT)  at (1.4484,0.3104);
  \coordinate (xTt) at (-0.9571,0.0571);
  \coordinate (mKL)  at ($(K)!.5!(L)$);
  \coordinate (mKM)  at ($(K)!.5!(M)$);
  \coordinate (mLM)  at ($(L)!.5!(M)$);
  \coordinate (mKMt) at ($(K)!.5!(Mt)$);
  \coordinate (mLMt) at ($(L)!.5!(Mt)$);

  \draw[tri] (K) -- (M) -- (L) -- (Mt) -- cycle;
  \draw[tri] (K) -- (L);

  \draw[vor] (xTt) -- (xT);
  \foreach \c/\m in {xT/mKM, xT/mLM, xTt/mKMt, xTt/mLMt}{
    \draw[vor] (\c) -- ($(\c)!1.45!(\m)$);
    \draw[ext] ($(\c)!1.45!(\m)$) -- ($(\c)!2!(\m)$);
  }
  \draw[ext] (xTt) -- ($(xTt)!-1!(mKL)$);   

  \draw ($(mKL)!.35cm!(K)$) -- ++($(xT)!.35cm!(xTt)-(xT)$) -- ($(mKL)!.35cm!(xTt)$);

  \draw[dashed]         (K) -- (xT);
  \draw[dash dot]       (L) -- (xT);
  \draw[dash dot]       (M) -- (xT);

  \coordinate (H1) at ($(L)!(K)!(M)$);
  \draw[dotted, -{Straight Barb[open]}] (K) -- (H1);
  \node at ($(K)!.42!(H1)+(0.3,0.1)$) {$h_1$};
  \draw ($(H1)!.25cm!(M)$) -- ++($(H1)!.25cm!(K)-(H1)$) -- ($(H1)!.25cm!(K)$);
  \draw ($(mLM)!.25cm!(M)$) -- ++($(mLM)!.25cm!(xT)-(mLM)$) -- ($(mLM)!.25cm!(xT)$);

  \pic[draw, "$\theta$", angle radius=0.75cm, angle eccentricity=1.35] {angle = K--xT--xTt};
  \pic[draw, "$\theta_3$", angle radius=1cm, angle eccentricity=1.45] {angle = K--M--L};
  \pic[draw, "$\tilde\theta$", angle radius=0.9cm, angle eccentricity=1.4] {angle = L--Mt--K};

  \coordinate (a) at ($(mKL)!.45!(L)$);
  \coordinate (b) at ($(a)+(xT)-(mKL)$);
  \coordinate (c) at ($(a)+(xTt)-(mKL)$);
  \draw[densely dotted] (xT) -- (b);
  \draw[densely dotted] (xTt) -- (c);
  \draw[dim] (a) -- (b) node[midway, below right=2pt] {$m_{K,T}$};
  \draw[dim] (a) -- (c) node[midway, below left=2pt] {$m_{K,\tilde T}$};

  \node[above]      at (K)   {$K = x_1$};
  \node[below]      at (L)   {$L = x_2$};
  \node[right]      at (M)   {$M = x_3$};
  \node[left]       at (Mt)  {$\tilde M = \tilde x_3$};
  \node[fill=white, inner sep=1pt] at ($(xT)+(0.55,0.05)$) {$x_T$};
  \node at (2.95,-0.75)  {$T$};
  \node at (-2.6,0.6)  {$\tilde T$};
  \node[below right=1pt] at ($(mKL)!.5!(xTt)$) {$\sigma = K|L$};
\end{tikzpicture}
\end{center}
\caption{Triangular finite element mesh and associated Vorono\"{\i}  cells}
\label{trianglef}
\end{figure}
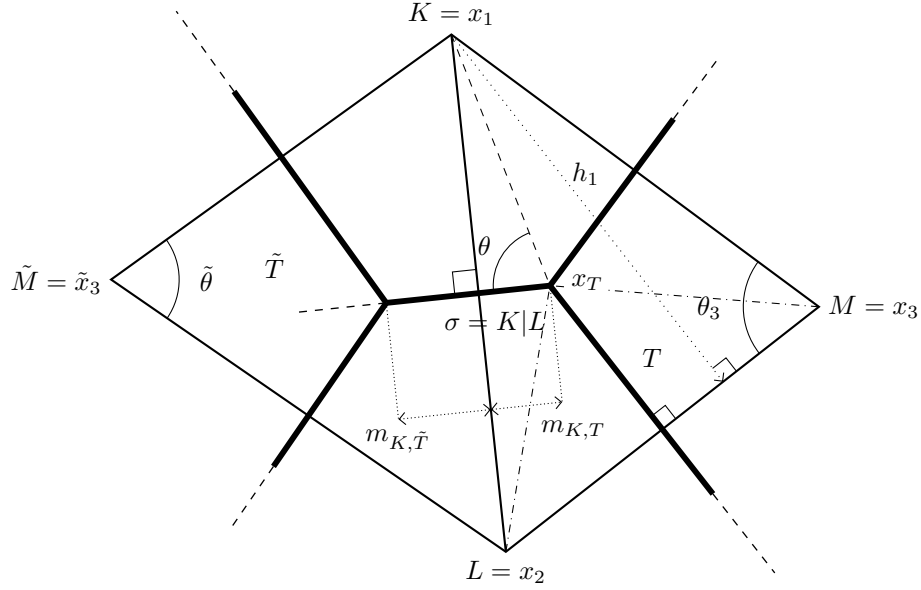

This well-known fact may be proven by considering  two nodes of the mesh, denoted by  $K = x_1$ and $L = x_2$ and the two triangles   $T$ and $\tilde T$  which share  the line segment $x_1 x_2$ as a common edge.
Let  $\phi_1$ and  $\phi_2$  be the two piecewise linear finite element  shape functions  respectively associated  to  the vertices $x_1$ and  $x_2$.  
Let us compute 
\[\int_\Omega \grad \phi_1 (x) \cdot \grad \phi_2(x) dx =    \int_{T \cup \tilde T}  \grad \phi_1 (x) \cdot \grad \phi_2(x) dx.\]
 Now let   $\theta_1$, $\theta_2$ and $\theta_3$ be the angles at vertices $x_1, x_2 $ and $x_3 $ of $T$
 (see Figure \ref{trianglef}). 
For $i = 1,2,3$, let $\bfn_i$ denote the outward unit normal vector to the side opposite to $x_i$.
Since $\phi_1$ (resp. $\phi_2$) is a linear function on $T$, its gradient is constant, and since $\phi_i(x_j) = \delta_{i,j}$, we obtain that
$\grad \phi_i(x) = - \frac 1 {h_i} \bfn_i $ for $i =1,2$, where $h_i$ is the height of the triangle with respect to the vertex $x_i$. 
Since $\bfn_1 \cdot \bfn_2 = -\cos \theta$, 
\[- \int_T \grad \phi_1 (x) \cdot \grad \phi_2(x) dx  = \vert T\vert\frac 1 {h_1 h_2} \cos \theta,  \]
where $\theta$ is the angle of the triangle at its vertex $M = x_3$. 
Now  the area of  the triangle $T$ is equal to $ \vert T  \vert =  \frac 1 2 \sin \theta d(L,M) d( M,K)  = \dsp\frac{h_1 h_2}{2  \sin \theta}$.
Thus, 
\[- \int_T \grad \phi_1 (x) \cdot \grad \phi_2(x) dx  = \frac 1 {2 \tan \theta}. \]
Let $x_T$ denote the circumcenter of $T$. 
Since the triangles $x_1 x_T x_3$, $x_2 x_T x_3$ and  $x_1 x_T x_2$ are isosceles, one gets that the angle between the line segment $x_1 x_T$ and the orthogonal bisector of $[x_1 x_2]$ is also $\theta$ and therefore, $\tan \theta =   \frac {d_{K,L}}  { 2 m_{K,T}}$
where $m_{K,T}$ denotes the distance between $\edgecvcvv$ and the circumcenter of $T$, and $d_{K,L}$ denotes the distance between $K$ and $L$, which is also the length of $\edge$.
Therefore, 
\[- \int_T \grad \phi_1 (x) \cdot \grad \phi_2(x) dx  = \frac {  m_{K,T}} {d_{K,L}}\]
The same computation on $\tilde T$ yields that:
 \[-\int_{\tilde T} \grad \phi_1 (x) \cdot \grad \phi_2(x) dx =  
  \frac {m_{K,\tilde T}} {d_{K,L}}\] if the angles 
 $\theta$ and $\tilde \theta$ are such that $\theta+ \tilde \theta< \pi$  (weak Delaunay condition), which, in turn yields the expected result. 

The CFV/DM scheme (finite volume scheme on the dual mesh) reads
\[
- \sumv \transcvcvv (\ucvv - \ucv) = \int_\cv f(x) dx, \text{ if } \cv 
\in \Vertices \cap \O,
\]
\[
u_\cv = g(\cv), \text{ if } \cv \in \Vertices \cap \dr \O,
\]

where $\cv$ stands for an element of $\Vertices$ or for the control 
volume (of the dual mesh) associated to this point.

The finite element scheme \index{Finite element method!-- on the primal mesh} (on the primal mesh) reads

\[
- \sumv \transcvcvv (\ucvv - \ucv) = \int_\O f(x) \varphi_\cv(x) dx, \text{ if } \cv 
\in \Vertices \cap \O,
\]
\[
u_\cv = g(\cv), \text{ if } \cv \in \Vertices \cap \dr \O.
\]

\begin{figure}[ht]
\begin{center}
\begin{tikzpicture}[x=0.00063333in, y=0.00063333in, line cap=round, line join=round]
  \draw (12,3012) -- (12,12) -- (3012,12) -- (12,3012) -- (3012,3012) -- (12,12) -- (3012,12) -- (3012,3012);
  \draw[dashed] (12,1512) -- (1512,3012) -- (3012,1512) -- (1512,12) -- (12,1512);
  \node[inner sep=0pt, anchor=base west] at (1137,1512) {$K$};
\end{tikzpicture}

\caption{An example of a triangular primal mesh (solid line) and
a dual Vorono\"{\i}  control volume (dashed line)}
\label{cfvdm}
\end{center}
\end{figure}
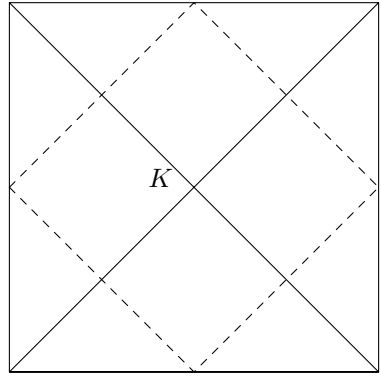

Therefore, the only difference between the finite element and the finite volume schemes is in the definition of the right hand sides. Note that these right hand sides may be quite  different. 
Consider for example a node $\cv$ which is the vertex of four identical triangles featuring  an angle of $\frac \pi 2$ at the vertex $\cv$, as depicted in Figure \ref{cfvdm},  and denote by $a$ the area of each of these triangles. 

Then, for $f \equiv 
1$, the right hand side computed for the discrete equation associated to 
the node $\cv$ is equal to $a$ in the case of the finite element 
(piecewise linear finite element) scheme, 
and equal to $2a$ for the dual mesh finite volume (CFV/DM) scheme. 
Both schemes may be shown to converge, by 
using  finite volume techniques for the CFV/DM scheme (see previous 
sections), and finite element techniques for the piecewise linear finite element (see 
e.g.\citet{ciarlet}).
\smallskip

Let us now weaken  the   hypothesis 
that all angles of the triangles of the primal mesh   
$\meshprim$ are acute  to   the so called
Delaunay condition \index{Delaunay condition} and the additional assumption that an angle of an 
element of $\meshprim$ is less or equal $\pi / 2$ if its opposite edge 
lies on $\dr \O$ (see e.g. \citet{Vanselow}).
Under this new assumption  the 
schemes  (piecewise linear finite element
finite element and CFV/DM with the Vorono\"{\i} mesh defined 
with $\Vertices$) still 
lead to the same transmissibilities \index{Transmissibility}  and still differ in the definition
of the right hand sides.

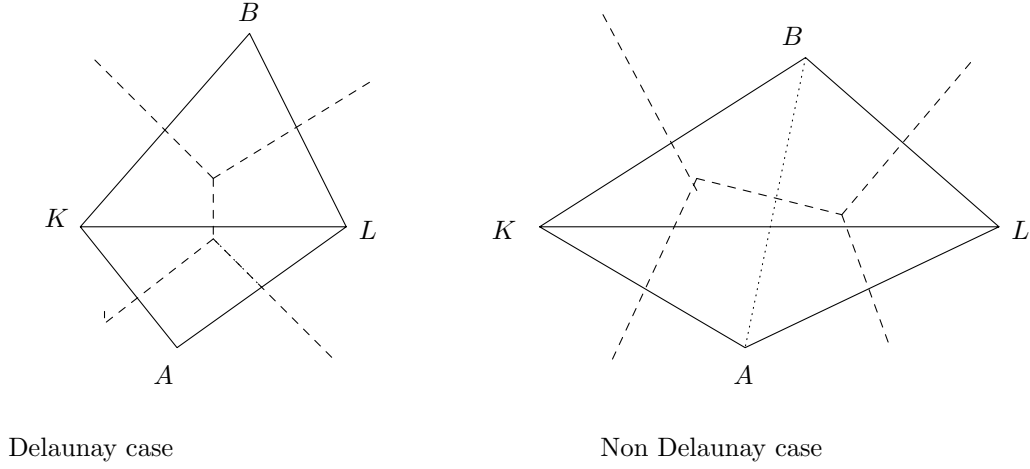
\begin{figure}[ht]
\begin{center}
\begin{tikzpicture}[x=0.00083333in, y=0.00083333in, line cap=round, line join=round]
  \node[inner sep=0pt, anchor=base west] at (3675,58) {Non Delaunay case};
  \node[inner sep=0pt, anchor=base west] at (0,58) {Delaunay case};
  \node[inner sep=0pt, anchor=base west] at (1425,2758) {$B$};
  \node[inner sep=0pt, anchor=base west] at (900,508) {$A$};
  \node[inner sep=0pt, anchor=base west] at (4800,2608) {$B$};
  \node[inner sep=0pt, anchor=base west] at (4500,508) {$A$};
  \node[inner sep=0pt, anchor=base west] at (6225,1408) {$L$};
  \node[inner sep=0pt, anchor=base west] at (2175,1408) {$L$};
  \node[inner sep=0pt, anchor=base west] at (3000,1408) {$K$};
  \node[inner sep=0pt, anchor=base west] at (225,1483) {$K$};
  \draw[dotted] (4950,2533) -- (4575,733);
  \draw[dashed] (1275,1783) -- (1275,1408);
  \draw[dashed] (1275,1783) -- (2250,2383);
  \draw[dashed] (5175,1558) -- (5475,733);
  \draw[dashed] (5175,1558) -- (6000,2533);
  \draw[dashed] (4275,1783) -- (5175,1558);
  \draw[dashed] (4275,1783) -- (3750,658);
  \draw[dashed] (4275,1708) -- (3675,2833);
  \draw[dashed] (1275,1408) -- (2025,658);
  \draw[dashed] (1275,1408) -- (600,883) -- (600,958);
  \draw[dashed] (1275,1783) -- (525,2533);
  \draw (3300,1483) -- (4950,2533) -- (6150,1483) -- (4575,733) -- (3300,1483) -- (6150,1483);
  \draw[dotted] (1575,1108) -- (1275,1408);
  \draw (450,1483) -- (1500,2683) -- (2100,1483) -- (1050,733) -- (450,1483) -- (2100,1483);
\end{tikzpicture}

\caption{Construction of the Vorono\"{\i} 
dual cells (dashed line) in the case
of  a triangular primal mesh (solid line) with and without the Delaunay condition}
\label{fig.delaunay}
\end{center}
\end{figure}

Recall that the Delaunay condition states that no 
neighboring element (of $\meshprim$) is included in the circumscribed circle 
of an 
arbitrary element of $\meshprim$. This is equivalent to saying that the 
sum of two opposite angles to an edge is less or equal $\pi$. 
As shown in  Figure \ref{fig.delaunay}, 
the dual 
mesh is still admissible in the sense of Definition 
\refp{meshdirichlet} and is still constructed with the orthogonal 
bisectors of the edges of the elements
of $\meshprim$, parts of these orthogonal bisectors (and parts of $\dr 
\O$) give the boundaries to 
the control volumes of the dual mesh (see Figure~\ref{fig.delaunay})
is not the 
case when $\meshprim$ does not satisfy the Delaunay condition.

\medskip

Consider now a primal mesh, $\meshprim$, consisting of triangles, but 
which does not satisfy the Delaunay \index{Delaunay condition} condition and let the 
dual mesh be the Vorono\"{\i} mesh defined with $\Vertices$. Then, the two 
schemes, piecewise linear finite element and CFV/DM are quite different. If the 
Delaunay condition does not hold say between the
angles $\dsp{\widehat{\cv A\cvv}}$ and $\dsp{\widehat{\cv B\cvv}}$ (the triplets 
($\cv,A,\cvv$) and ($\cv,B,\cvv$) defining two elements of $\meshprim$), 
the sum of 
these two angles is greater than $\pi$ and the transmissibility \index{Transmissibility}  
$\transcvcvv= - \int_\O  \grad \varphi_{\cv}(x) \cdot  \grad
\varphi_{\cvv}(x)dx$ between the two control volumes associated 
respectively to $\cv$ and $\cvv$ becomes negative with the piecewise linear finite element
 scheme; 
 there is no transmissibility between 
$A$ and $B$ (since $A$ and $B$ do not belong to a common element of 
$\meshprim$). Hence the maximum 
principle is no longer respected for the finite element scheme, while it remains
valid for the CFV/DM finite volume scheme.  This is due to the fact that 
the CFV/DM scheme allows an exchange term between $A$ and $B$, with a 
positive transmissibility (and leads to no exchange term between $\cv$ and 
$\cvv$), while the finite element scheme does not. Also note also that the common edge to 
the control volumes (of the dual mesh) associated to $A$ and $B$ is not a 
part of an orthogonal bisector of an edge of an element
of $\meshprim$ (it is a part of the orthogonal bisector of the segment 
$[A,B]$).

~


To conclude this section, note that an admissible mesh for the classical 
finite volume is  generally not a dual mesh of a primal triangular mesh 
consisting of triangles (for instance, the general triangular
meshes which are considered in \citet{VF4} are not dual meshes to triangular meshes).

\subsection{``Finite Volume Finite Element" methods}\index{Finite volume 
finite element method|(}

The ``finite volume finite element" method for
elliptic problems also  uses a dual mesh $\mesh$
constructed from a finite element primal mesh, such that each cell of $\mesh$ is
associated with a vertex of the primal mesh $\meshprim$. Let $\Vertices$ 
again denote  the set of vertices of $\meshprim$. 
As in the classical finite volume method, the conservation law is integrated over each
cell of the (dual) mesh. Indeed, this integration is performed only if the cell 
is associated to a vertex (of the primal mesh) belonging to $\O$. 

Let us consider Problem \eqref{rouene2D}. Integrating
the conservation law over $\cv_P$, where $P \in \Vertices \cap \O$ and 
$\cv_P$ is the control volume (of the dual mesh) associated to $P$ yields

\[
- \int_{\dr \cv_P} \grad u(x) \cdot \n_P(x) \dfrontiere(x) = \int_{\cv_P} 
f(x) dx.
\]

(Recall that $\n_P$ is the unit normal vector to $\dr \cv_P$ outward to
$\cv_P$.) Now, following the idea of finite element methods, the function  $u$ is 
approximated by a Galerkin expansion\index{Galerkin expansion}
$\sum_{M\in \Vertices} u_M \varphi_M$, 
where the functions $\varphi_M$ are the shape
functions of the piecewise linear finite element method. Hence, the 
discrete unknowns are $\{u_P,$ $P\in \Vertices\}$ and the scheme reads

\be
 -\sum_{M \in \Vertices}  \Bigl(\int_{\dr \cv_P}  \grad \varphi_M(x) 
 \cdot \n_P(x) \dfrontiere(x) \Bigr)
u_M =  \int_{\cv_P} f(x) dx, \; \forall \, P \in \Vertices \cap \O,
\label{cvfe}
\ee

\[
u_P=g(P), \; \forall \, P \in \Vertices \cap \dr \O.
\]

Equations \eqref{cvfe} may also be written under the conservative form

\be
\sum_{Q \in \Vertices}  E_{P,Q} =  \int_{\cv_P} f(x) dx,\; \forall \, P \in \Vertices \cap \O,
\label{cvfe1}
\ee

\be
u_P=g(P), \; \forall \, P \in \Vertices \cap \dr \O,
\label{cvfe2}
\ee

where 
\be
 E_{P,Q} = -\sum_{M \in \Vertices} \int_{\dr \cv_P \cap \dr \cv_Q}  \grad
\varphi_M(x) \cdot \n_P(x) \dfrontiere(x).
\label{cvfe3}
\ee

Note that   $ E_{Q,P} = - E_{P,Q}$.
Unfortunately, the exchange term $E_{P,Q}$  between $P$ and $Q$ 
 is not, in general, a function of the only unknowns $u_P$ and $u_Q$ 
(this property was used, in the previous sections, to obtain convergence 
results of finite volume schemes).
Another way to write \eqref{cvfe} is, thanks to  \eqref{rouencruu},

\[
- \sum_{Q \in \Vertices}  \Bigl(\int_{\dr \cv_P}  \grad \varphi_Q(x) 
\cdot \n_P(x) \dfrontiere(x) \Bigr)
(u_Q - u_P) =  \int_{\cv_P} f(x) dx, \; \forall \, P \in \Vertices \cap \O.
\]

Hence a new exchange term from $P$ to $Q$ might be $\bar{E}_{P,Q} = - \Bigl(\int_{\dr \cv_P} 
\grad \varphi_Q(x) \cdot \n_P(x) \dfrontiere(x) \Bigr) (u_Q - u_P)$ and the scheme is therefore
conservative  if $\bar{E}_{P,Q} = - \bar{E}_{Q,P}$. Unfortunately, 
this is not the case for
a general  dual mesh.

~

There are several ways of constructing a dual mesh from a primal mesh. 
A common way (see e.g.  \citet{La})is to take a primal mesh ($\meshprim$) consisting of triangles and to construct the dual mesh with the medians (of the triangles of $\meshprim$), joining the centers of gravity of the triangles to the midpoints of the edges of the primal mesh. 
The main interest of this way is that the resulting scheme (called FVFE/M below, Finite Volume Finite Element with Medians) is very close to the piecewise linear finite element scheme associated to $\meshprim$. Indeed the FVFE/M scheme is defined by \eqref{cvfe1}-\eqref{cvfe3} while the  piecewise linear finite element scheme reads
\[
 \sum_{Q \in \Vertices}  E_{P,Q} =  \int_{\O} f(x)\varphi_P(x) dx,\; \forall \, P \in \Vertices \cap \O,
\]

\[
u_P=g(P), \; \forall \, P \in \Vertices \cap \dr \O,
\]

where  $E_{P,Q}$ is defined by \eqref{cvfe3}.

These two schemes only differ by the right hand sides and, in fact, 
these right hand sides are ``close" since 

\[
\meas(K_P) = \int_\O  \varphi_{P}(x)dx, \; \forall \, P \in \Vertices.
\]

This is due to the fact that $\int_T  \varphi_{P}(x)dx=\meas(T)/3$ and 
$\meas(K_P \cap T)=\meas(T)/3$, for 
all $T \in \meshprim$ and all vertex $P$ of $T$.

Thus, convergence properties of the FVFE/M scheme can be proved by using 
the finite element techniques. Recall however that the piecewise linear finite element scheme (and the 
FVFE/M scheme) does not satisfy the (discrete) maximum principle if 
$\meshprim$ does not satisfy the Delaunay \index{Delaunay condition} condition.

~
 
There are other means to construct a dual mesh starting from a primal 
triangular mesh. One of them is the Vorono\"{\i} 
mesh associated to the vertices of the primal mesh, another possibility 
is to join the centers of gravity; in the latter case, the control
volume associated to a vertex, say $S$, 
of the primal mesh is then limited by the
lines joining the centers of gravity of the neighboring triangles of which 
$S$ is a
vertex (with some convenient modification for the vertices which are on the boundary of
$\O$). See also \citet{Barth} for descriptions of dual meshes.

\smallskip

Note that the proof of convergence which we designed for finite volume 
with admissible meshes 
does not generalize to any ``FVFE" (Finite Volume Finite
Element) method for several reasons. In particular, since the exchange 
term between $P$ and $Q$ (denoted by 
$E_{P,Q}$) is not, in general, a function of the only unknowns $u_P$ and 
$u_Q$ (and even if it is the transmissibilities \index{Transmissibility}  may become 
negative) and also since, 
as in the case of the finite element method, the concept of
consistency  of the fluxes is 
not clear with the FVFE schemes. 
\index{Finite element method|)}

\subsection{Generalization to the three dimensional case} \label{cvfe3D}

The methods described in the three above sections
generalize to the three-dimensional case, in particular when the primal mesh is a tetrahedral
mesh. With such a mesh, the Delaunay\index{Delaunay condition} condition no longer ensures
the non negativity of the transmissibilities in the case of the piecewise linear 
finite element
method. It is however possible to construct a dual mesh (the ``three-dimensional 
Vorono\"{\i}" mesh) to a
Delaunay triangulation such that the FVFE scheme leads to positive 
transmissibilities, and
therefore such that the maximum principle holds, see \citet{Putti}. 

\smallskip

Note that   the theoretical results
(convergence and error estimate) which were shown for the classical finite volume
method on an admissible mesh (sections \refp{vfquatre} 
and \refp{nll}) still hold for CFV/DM in three-dimensional, since the dual mesh is admissible.

\section{Mesh refinement and singularities} \label{meshref}

\hpage{Elliptic multid}

Some problems involve singular source terms. In the case of petroleum engineering 
for instance, one may model (in two space dimensions) the well with a Dirac measure.
Other problems may  require a better precision of some unknown in certain areas.
This section is devoted to the treatment of this kind of problem, either with an
adequate treatment of the singularity or by mesh refinement. 

\index{Finite volume finite element method|)}
\subsection{Singular source terms and finite volumes}
\index{Singular source terms in elliptic equations|(}
It is possible to take into account, in the discretization with the finite volume
method, the singularities of the solution of an elliptic problem. A common example is
the study of wells in petroleum engineering. As a model example we can consider the
following problem, which appears, for instance, in the study of a two phase flow in a
porous medium. Let $B$ be the ball of $\R^2$ of center $0$ and radius $r_p$ ($B$
represents a well of radius $r_p$). Let
$\Omega=(-R,R)^2$ be the whole domain of simulation;  $r_p$
is of the order of 10 cm while $R$ can be  of the order of 1 km for 
instance. An approximation to 
the solution of the following problem is sought: 

\be
\begin{array}{l}
- \div (\grad u)(x) = 0, \ \ x \in \Omega \setminus B, \\
u(x)=P_p, \ \  x \in \partial B, \\
\mathrm{``BC"  on} \ \  \partial \Omega,
\end{array}
\label{sing}
\ee

where ``BC" stands for some ``smooth" boundary conditions on $\partial \Omega$ (for
instance, Dirichlet or Neumann condition). This system is a mathematical model (under
convenient assumptions\ldots) of the two phase flow problem, with $u$ representing the
pressure of the fluid and $P_p$  an imposed pressure at the well. In order to
discretize (\ref{sing}) with the finite volume method, a mesh $\mesh$ of
$\Omega$ is introduced. For the sake of simplicity, the elements of
$\mesh$ are assumed to be squares of length $h$ (the method is easily generalized to
other meshes). It is assumed that the well, represented by $B$, is located in the middle
of one cell, denoted by $K_0$, so that the origin  $0$ is  the center of $K_0$. It is
also assumed that the mesh size, $h$, is large with respect to the radius of the well,
$r_p$  (which is the case in real applications, where, for instance, $h$ 
ranges between 10 and 100 m). Following the principle of the finite volume method,  one discrete unknown
$u_\cv$ per cell $\cv$ ($\cv \in \mesh$) is introduced in order to discretize the following
system:

\be
\begin{array}{l}
\dsp{ \int_{\partial \cv} \grad u(x) \cdot \bfn_\cv(x) \dfrontiere(x)= 0, \ \  \cv \in \mesh, \ \ \cv \ne K_0,} \\
\dsp{ \int_{\partial K_0} \grad u(x) \cdot \bfn_{K_0}(x) \dfrontiere (x) = 
\int_{\partial B} \grad u(x) \cdot \bfn_{B}(x)
\dfrontiere(x),}
\end{array}
\label{singd}
\ee

where $\bfn_P$ denotes the normal to $\partial P$, outward to $P$ (with $P= \cv, K_0$ or $B$).

Hence, we have to discretize $\grad u \cdot \bfn_\cv$ on $\partial \cv$ (and $\grad u \cdot \bfn_B$ on $\partial B$) in terms of 
$\{u_\cvv,\cvv \in \mesh\}$ (and ``BC" and $P_p$).

The problems arise in the discretization of $\grad u \cdot \bfn_{K_0}$ and $\grad u \cdot n_B$. 
Indeed, if $\edge=\edgecvcvv$ is the common edge to $\cv$ and 
$\cvv$ (elements of $\mesh$), with $\cv \ne K_0$ and $\cvv \ne K_0$, 
since the solution of (\ref{sing}) is ``smooth" enough with respect to 
the mesh size, except ``near" the well, $\grad u \cdot n_{\cv}$ can be 
discretized by $\frac 1 h (u_{\cvv}-u_{\cv})$ on $\edge$.

In order to 
discretize $\grad u$ near the well, it is assumed that $\grad u \cdot n_B$ 
is constant on $\partial B$. Let $q(x) =-2 \pi r_p \grad u \cdot n_B$  for  
$x \in \partial B$ (recall that $\bfn_B$ is the normal to $\partial B$, outward 
to $B$). Then $q \in \R$ is a new unknown, which satisfies

\[
\int_{\partial B} - \grad u \cdot n_{B}\dfrontiere (x) = q.
\]

Denoting by $\vert \cdot \vert$ the euclidian norm in $\R^2$, and $u$ the solution to
(\ref{sing}), let $v$ be defined by

\be
v(x)= \frac q {2 \pi} \ln(\vert x \vert) + u(x), \ \ x \in \Omega \setminus B,
\label{reg1}
\ee
\be
v(x)  =  \frac q {2 \pi} \ln(r_p) + P_p, \ \ x \in B. 
\label{reg2}
\ee

Thanks to the boundary
conditions satisfied by $u$ on $\partial B$, the function $v$ satisfies
 $- \div (\grad v)
= 0$ on the whole domain $\Omega$, and therefore $v$ is regular on the whole domain
$\Omega$. Note that, if we set 
\[u(x)=-\frac q {2 \pi} \ln(\vert x \vert) 
+ v(x), \text{  a.e. } x \in \Omega,
\]
 then 
 \[- \div (\grad u) = q \delta_0\text{ on } \Omega,\]
 where $\delta_0$ is the Dirac mass at $0$. 
A discretization of $\grad u \cdot n_{K_0}$ is now obtained in the following way. 
Let $\edge$ be the common edge to $K_1 \in \mesh$ and  $K_0$, since $v$ is smooth, it is possible to approximate $\grad v \cdot \bfn_{K_0}$ on $\edge$ by $\frac 1 h (v_{K_1}-v_{K_0})$, where $v_{K_i}$
is some approximation of $v$ in $K_i$ 
(e.g. the value of $v$ at the center of $K_i$).
Then, by (\ref{reg2}), it is natural to set 
\[ 
v_{K_0} = \frac q {2 \pi} \ln(r_p) + P_p,
\]
 and by 
 (\ref{reg1}),
 \[
 v_{K_1}=\frac q {2 \pi} \ln(h) + u_{K_1}.
 \]
 By (\ref{reg1}) and 
from the fact that the integral over $\edge$ of $\grad (\frac q {2 \pi} \ln(\vert 
x \vert)) \cdot
\bfn_{K_0}$ is equal to $\frac q 4$, we find the following approximation for $\int_\edge \grad
u \cdot n_{K_0} \dfrontiere$:

\[
-\frac q  4 + \frac q {2 \pi} \ln(\frac h {r_p}) + u_{K_1} - P_p.
\]

The discretization is now complete, there are as many equations as unknowns. 
The discrete unknowns
appearing in the discretized problem are $\{u_\cv, \cv \in \mesh, \cv \ne 
K_0\}$ and $q$. Note
that, up to now, the unknown $u_{K_0}$ has not been used. The discrete 
equations are given by \eqref{singd} where each term of \eqref{singd} is 
replaced by its approximation in terms of $\{u_\cv, \cv \in \mesh, \cv \ne 
K_0\}$ and $q$. In particular, the
discrete equation ``associated" to the unknown $q$ is the discretization of the second
equation of (\ref{singd}), which is

\be
\sum_{i=1}^4 (\frac q {2 \pi} \ln(\frac h {r_p}) + u_{K_i} - P_p)=0,
\label{do1}
\ee

where $\{K_i, i=1,2,3,4\}$ are the four neighbouring cells   to $K_0$.

It is possible to replace the unknown $q$ by the unknown $u_{K_0}$ (as it is done in
petroleum engineering) by setting

\be
u_{K_0}= \frac q  4 - \frac q {2 \pi} \ln(\frac h {r_p}) + P_p,
\label{uko}
\ee

the interest of which is that it yields the usual formula for the 
discretization of $\grad u \cdot n_{K_0}$ on $\edge$ if $\edge$ is the 
common edge to $K_1$ and $K_0$, namely 
$\frac 1 h$ $(u_{K_1}-u_{K_0})$; the discrete equation associated to the 
unknown $u_{K_0}$ is then (from (\ref{do1}))

\[
\sum_{i=1}^4(u_{K_i} -  u_{K_0})=-q
\]

and  (\ref{uko}) may be written as:

\[
 q=i_p(P_p-u_{K_0}), \text{ with }
i_p=\frac 1 {-\frac 1 4 +\frac1 {2 \pi} \ln (\frac h {r_p})}.
\]

This last equation defines $i_p$, the so called ``well-index" in petroleum engineering.
With this formula for $i_p$, the discrete unknowns are now $\{u_\cv, \cv 
\in \mesh\}$. The discrete 
equations associated to $\{u_\cv, \cv \in \mesh, \cv \ne 
K_0\}$ are given by the first part of \eqref{singd} where each terms of \eqref{singd} is 
replaced by its approximation in terms of $\{u_\cv, \cv \in \mesh\}$ 
(using also ``BC" on $\dr \O$). The
discrete equation associated to the unknown $u_{K_0}$ is 

\[
\sum_{i=1}^4(u_{K_i} -  u_{K_0})=-i_p(P_p-u_{K_0}),
\]

where $\{K_i, i=1,2,3,4\}$ are the four neighbouring cells  to $K_0$.

Note that the discrete unknown $u_{K_0}$ is somewhat artificial, it does not really
represent the value of $u$ in $K_0$. In fact, if $x \in K_0$, the ``approximate value"
of $u(x)$ is $ - \frac q {2 \pi} \ln(\frac {\vert x \vert} {r_p}) + P_p$
and $u_{K_0}=\frac q  4 - \frac q {2 \pi} \ln(\frac h {r_p}) + P_p$.

\index{Singular source terms in elliptic equations|)}
\subsection{Mesh refinement}\index{Mesh!-- Refinement} 

Mesh refinement consists in using, in certain areas of the domain,  control volumes
of smaller size than elsewhere. In the case of triangular grids, a refinement may
be performed for instance by dividing each triangle in the refined area into four
subtriangles, and those at the boundary of the refined area in two triangles.  
Then, with some additional technique (e.g. change of diagonal), one may obtain an
 admissible mesh in the
sense of definitions \refp{meshdirichlet},  \refp{meshneuman} and 
\refp{rouenmesh}; therefore the error estimates \index{Error estimate!-- for an elliptic 
equation!-- in the general case} \refp{ellesterr}, \refp{nllesterr}
and \refp{rouenesterr}  hold under the same assumptions. 

\medskip

In the case of rectangular grids, the same refining procedure leads to ``atypical"
nodes and edges, i.e. an edge $\sigma$ of a given control volume $\cv$  may be
common to two other control volumes, denoted by $\cvv$ and $\cvw$. 
This is also true in the triangular case if the triangles of the boundary 
of the refined area are left untouched. 


Let us consider for instance the same problem as in section \refp{smesh}, 
with the same assumptions and notations, namely the discretization of

\[
\bal
-\Delta u(x,y)  =  f(x,y), \, (x,y)  \in \O = (0,1)\times(0,1),  \\
u (x,y)  =  0, (x,y) \in \partial \O.              
\ean
\]

 It is easily seen that, in this case, if the  approximation of the fluxes
is performed using differential quotients such as in \refep{fluxrect},  the fluxes
on the
``atypical" edge  $\edge$ cannot be consistent, since the lines joining the centers
of $\cv$ and $\cvv$ and the centers
of $\cv$ and $\cvw$  are not orthogonal to $\edge$.  However, the error which
results from this lack of consistency can be controlled if the number of atypical
edges is not too large. 

\smallskip

In the case of rectangular grids (with a refining procedure), denoting by 
$\edgesatyp$ the set of ``atypical" edges of a given
mesh $\mesh$, i.e. edges  with separate more than two control volumes,  and
$\meshatyp$ the set of ``atypical" control volumes, i.e. the control volumes
containing an atypical edge in their boundaries; let $e_\cv$ denote  the 
error between $u(x_\cv)$ and $\ucv$ for
each control volume $\cv$, and $e_\mesh$ denote the piecewise constant function
defined by $e(x) = e_\cv$ for any $x \in \cv$, then one has 

\[
\Vert e
\Vert_{L^2(\O)} \le C(\sizemesh + (\sum_{\cv \in \meshatyp} 
\meas(\cv))^{\frac 1 2}).
\]
The proof is
similar to that of Theorem \refp{ellesterr}. 
It is detailed in \citet{belmouhoub}.





\section{Compactness results}\label{rellichd}

\hpage{Elliptic multid}

\index{Compactness!-- results in $L^2$|(}
This section is devoted to some functional analysis results 
which were 
used in the previous section.  Let $\O$ be a bounded open set of $\R^d$, $d \ge 1$. 
 Two
relative compactness results in $L^2(\O)$ for sequences ``almost" bounded in
$H^1(\O)$ which were used in the proof of convergence of the schemes are 
presented here. Indeed, they are variations of the Rellich theorem (relative
compactness in $L^2(\O)$ of a bounded sequence in $H^1(\O)$ or $H^1_0(\O)$). The
originality of these results is not the fact that the sequences are relatively compact in $L^2(\O)$, which is an immediate consequence of the Kolmogorov theorem (see below), but the fact that the eventual limit, in $L^2(\O)$, of the sequence (or of a subsequence) is necessarily in $H^1(\O)$ (or in $H^1_0(\O)$ for Theorem \ref{kolmh10}), a space which does not contain the elements of the sequence.

We shall make use in this section of the Kolmogorov compactness theorem in $L^2(\O)$ which we now recall.
The essential part of the proof of this theorem may be found in \citet{bre}.

\begin{theo}[Kolmogorov compactness lemma] \label{Kolm} \index{Compactness!Kolmogorov's -- theorem}\index{Kolmogorov's theorem}
Let $\Omega$ be an open  bounded  set of $\R^N$, $N \ge 1$, $1 \le q < \infty$ and $A \subset L^q(\o)$. 
Then, $A$ is relatively compact in $L^q(\o)$ if and only if there exists $\{p(u)$, $u \in A\}$ $\subset L^q(\R^N)$ such that
\benum
\item $p(u)=u$ a.e. on $\o$, for all $u \in A$,
\item $\{p(u)$, $u \in A\}$ is bounded in $L^q(\R^N)$,
\item $\norm{p(u)(\cdot+\eta)-p(u)}{L^q(\R^N)} \tends 0$ as $\eta \tends 0$, uniformly with respect to $u \in A$.
\eenum
\end{theo}
Let us now state the compactness results used in this chapter.
\begin{theo}[Compactness of a bounded sequence and regularity of the limit] \label{kolmh10}
Let $\O$ be an open bounded set of $\R^d$ with a Lipschitz continuous boundary, $d \ge 1$, and $\{u_n,$ $n \in \N\}$ a bounded sequence of $L^2(\O)$. 
For $n \in \N$, one defines $\tilde{u}_n$ by $\tilde{u}_n=u_n$ a.e. on $\O$ and $\tilde{u}_n=0$ a.e. on $\R^d \setminus \O$. 
Assume that there exist $C \in \R$ and $\{h_n,$ $n \in \N\} \subset \Rp$ such that $h_n \tends 0$ as $n \tends \infty$
and 
\be
\norm{\tilde{u}_n(\cdot+\eta)-\tilde{u}_n}{L^2(\R^d)}^2 \le C \vert \eta \vert
(\vert \eta \vert + h_n),
\forall n \in \N, \; \forall \eta \in \R^d.
\label{transh10}
\ee
Then, $\{u_n,$ $n \in \N\}$ is relatively compact in $L^2(\O)$. 
Furthermore, if $u_n \tends u$ in $L^2(\O)$ as $n \tends \infty$, then $u \in H^1_0(\O)$.
\end{theo}

\proof of Theorem \ref{kolmh10}

\smallskip

Since $\{h_n,$ $n \in \N\}$ is bounded, the fact that $\{u_n,$ $n \in \N\}$ is relatively compact in $L^2(\O)$ is an immediate consequence of Theorem
\ref{Kolm}, taking $N=\dim$, $\o=\O$, $q=2$ and $p(u_n)=\tilde u_n$. 
Then, assuming that $u_n \tends u$ in
$L^2(\O)$ as $n \tends \infty$, it is only necessary to prove that 
$u \in H^1_0(\O)$. Let us first remark that $\tilde {u}_n \tends \tilde u$ in $L^2(\R^d)$, as $n
\tends \infty$, with $\tilde u=u$ a.e. on $\O$ and $\tilde u=0$ a.e. on
$\R^d \setminus \O$.

Then, for $\varphi \in C^\infty_c(\R^d)$, one has, for all $\eta \in \R^d$, $\eta\neq0$ and $\nnn$, using the Cauchy-Schwarz inequality and thanks to
\eqref{transh10},

\[
\int_{\R^d} \frac {(\tilde u_n(x+\eta)-\tilde u_n(x))} {\vert \eta
\vert}\varphi(x)dx
\le \frac {\sqrt{C \vert \eta \vert
(\vert \eta \vert + h_n)}}{\vert \eta \vert} \norm{\varphi}{L^2(\R^d)},
\]

which gives, letting  $n \tends \infty$, since $h_n \tends 0$,

\[
\int_{\R^d} \frac {(\tilde u(x+\eta)-\tilde u(x))} {\vert \eta \vert}\varphi(x)dx
\le \sqrt{C}\norm{\varphi}{L^2(\R^d)},
\]

and therefore, with a trivial change of variables in the integration,

\be
\int_{\R^d} \frac {(\varphi(x-\eta)-\varphi(x))} {\vert \eta
\vert}\tilde u(x)dx
\le \sqrt{C}\norm{\varphi}{L^2(\R^d)}.
\label{ueh1}
\ee

Let $\{e_i$, $i=1,\ldots,d\}$ be the canonical basis of $\R^d$. For $i \in
\{1,\ldots,d\}$ fixed, taking $\eta = he_i$ in \eqref{ueh1} and letting $h \tends
0$ (with $h>0$, for instance) leads to 

\[
-\int_{\R^d}\frac {\pa \varphi(x)} {\pa x_i} \tilde u(x)dx \le \sqrt{C}\norm{\varphi}{L^2(\R^d)},
\]

for all $\varphi \in C^\infty_c(\R^d)$.

This proves that $D_i\tilde u$ (the derivative of $\tilde u$ with respect to
$x_i$ in the sense of distributions) belongs to $L^2(\R^d)$, and therefore that
$\tilde u \in H^1(\R^d)$. Since $u$ is the restriction of $\tilde u$ on $\O$ and
since $\tilde u=0$ a.e. on $\R^d \setminus \O$, therefore $u \in
H^1_0(\O)$. This completes the proof of Theorem \ref{kolmh10}.
\eof

\begin{remark}[Direct proof of the regularity of the limit]In fact, the proof that a possible limit is in $H^1_0$ can be directly drawn from  the boundedness of the sequence of approximate solutions in  the $H^1$ discrete norm.
Let us detail this point.

Let $\O$ be an open bounded set of $\R^d$ with a Lipschitz continuous boundary, $d \ge 1$.
Let $(\mesh_n)_\nnn$ be a  sequence of admissible meshes (in fact, the orthogonality condition is not required for this direct proof) such that $\size(\mesh_n)$ tends to $0$, and let  $\{u_n \in X_{\Mesh_n},$ $n \in \N\}$  be a  sequence of functions  of $L^2(\O)$ weakly converging to $u$.
Let us assume that there exists a real number $C$ not depending on $n$ such that  $\Vert u_n \Vert_{1,\Mesh_n} \le C$.
Then $u \in H^1_0 (\Omega)$.

Indeed,  let $\bfvarphi  \in (C_c^\infty(\R^d))^d$ (note that $\bfvarphi $ does not vanish on the boundary of $\Omega$);
as in Step 4 of the proof of Theorem \ref{kolmh10}, let us define  for $n \in \N$, $\tilde{u}_n$ by $\tilde{u}_n=u_n$ a.e. on $\O$ and $\tilde{u}_n=0$ a.e. on $\R^d \setminus \O$. 

Then 
\begin{align*}
 \int_{\R^d} \tilde u_n  \dive \bfvarphi  dx  & =    \summesh u_\cv \int_{\partial \cv} \bfvarphi \cdot \bfn \dfrontiere  \\
								  &  \le   \sumedgesintkl \medge \vert u_K - u_L \vert \vert  \varphi_\edge \vert  +  \sum_{\substack{\edge\in \edgesext\\ \sigma \in \edgescv}} \medge \vert u_{K}   \vert \vert  \varphi_\edge \vert 
\end{align*}
where $\varphi_\edge$ is the mean value of $\bfvarphi \cdot \bfn $ over $\edge$.
By the Cauchy-Schwarz inequality, we obtain
\begin{equation*}
 \int_{\R^d} \tilde u_n  \dive \bfvarphi  dx  \le  \Vert u \Vert_{1,\mesh_n}  \left(\sumedges \medge \dedge \varphi_\edge^2 \right)^{1/2}.
\end{equation*}
Defining by $\widetilde\varphi_\edge$ the mean value of $\bfvarphi \cdot \bfn $ over $D_\edge$, we get
\[
 |\widetilde\varphi_\edge - \varphi_\edge| \le \size(\mesh_n) \Vert \bfvarphi\Vert_{1,\infty},
\]
where $\Vert \bfvarphi\Vert_{1,\infty}$ is a bound of the first derivatives of $\bfvarphi$. Therefore we get
\[
\left(\sumedges \medge \dedge \varphi_\edge^2 \right)^{1/2} \le  \left(\sumedges \medge \dedge \widetilde\varphi_\edge^2 \right)^{1/2} + \sqrt{d \meas(\O)} \size(\mesh_n) \Vert \bfvarphi\Vert_{1,\infty}.
\]
From the above results, we may write
\begin{equation*}
 \int_{\R^d} \tilde u_n  \dive \bfvarphi  dx  \le \sqrt{d} \ \Vert u \Vert_{1,\mesh_n} \Vert \vert \bfvarphi \vert \Vert_{L^2(\Omega)} + \sqrt{d \meas(\O)} \size(\mesh_n) \Vert \bfvarphi\Vert_{1,\infty}.
\end{equation*}

Passing to the limit, we then get that 
\[
 \lim_{n\to +\infty} \int_{\R^d} \tilde u_n \dive \bfvarphi dx \le C \Vert \vert \bfvarphi \vert \Vert_{L^2(\Omega)} .
\]
This in turn shows that that $\nabla u \in (L^2 (\R^d) )^d$ and $u \in H^1_0(\Omega)$. 
\end{remark}

\begin{theo}\label{kolmh1} Let $\O$ be an open bounded set of $\R^d$, $d \ge
1$, and $\{u_n,$ $n \in \N\}$ a bounded sequence of $L^2(\O)$. For $n \in \N$,
one defines $\tilde{u}_n$ by $\tilde{u}_n=u_n$ a.e. on $\O$ and $\tilde{u}_n=0$
a.e. on $\R^d \setminus \O$. Assume that
there exist $C \in \R$ and $\{h_n,$ $n
\in \N\} \subset \Rp$ such that $h_n \tends 0$ as $n \tends \infty$ and such
that

\be
\norm{\tilde{u}_n(\cdot+\eta)-\tilde{u}_n}{L^2(\R^d)}^2 \le C \vert \eta \vert,
\forall n \in \N, \; \forall \eta \in \R^d,
\label{transh11}
\ee

and, for all compact $\comp \subset \O$,

\be
\norm{{u_n}(\cdot+\eta)-{u_n}}{L^2(\comp)}^2 \le C \vert \eta
\vert(\vert \eta \vert + h_n),
\forall n \in \N, \; \forall \eta \in \R^d, \; \vert \eta \vert < d(\comp,\O^c).
\label{transh12}
\ee

(The distance between $\comp$ and
$\R^d \setminus \O$ is denoted by $d(\comp,\O^c)$.)

Then $\{u_n,$ $n \in \N\}$ is relatively compact in $L^2(\O)$. Furthermore, if
$u_n \tends u$  in $L^2(\O)$  as $n \tends \infty$, then $u \in H^1(\O)$.
\end{theo}

\proof of Theorem \ref{kolmh1}

\smallskip

The proof is very similar to that of Theorem \ref{kolmh10}. Using assumption 
\ref{transh11}, Theorem \ref{Kolm} yields that $\{u_n,$ $n \in \N\}$
is relatively compact in $L^2(\O)$. Assuming now that $u_n \tends u$ in
$L^2(\O)$, as $n \tends \infty$, one has to prove that $u \in H^1(\O)$.

Let $\varphi \in C^\infty_c(\O)$ and $\eps >0$ such that $\varphi(x)=0$
if the distance from $x$ to $\R^d \setminus \O$ is less than $\eps$. 
Assumption \ref{transh12} yields

\[
\int_{\O} \frac {({u_n}(x+\eta)-{u_n}(x))} {\vert \eta
\vert}\varphi(x)dx
\le \frac {\sqrt{C \vert \eta \vert
(\vert \eta \vert + h_n)}}{\vert \eta \vert} \norm{\varphi}{L^2(\O)},
\]

for all $\eta \in \R^d$ such that $0 < \vert \eta \vert <\eps$.

From this inequality, it may be proved, as in the proof of Theorem \ref{kolmh10}
(letting $n \tends \infty$ and using a change of variables in the integration),

\[
\int_{\O} \frac {(\varphi(x-\eta)-\varphi(x))} {\vert \eta
\vert}u(x)dx
\le \sqrt{C}\norm{\varphi}{L^2(\O)},
\]

for all $\eta \in \R^d$ such that $0 < \vert \eta \vert <\eps$.

Then, taking $\eta = he_i$ and letting $h \tends
0$ (with $h>0$, for instance) one obtains, for all $i \in \{1,\ldots,d\}$,

\[
-\int_{\O}\frac {\pa \varphi(x)} {\pa x_i} u(x)dx
\le \sqrt{C}\norm{\varphi}{L^2(\O)},
\]

for all $\varphi \in C^\infty_c(\O)$.

This proves that $D_i u$ (the derivative of $u$ with respect to
$x_i$ in the sense of distributions) belongs to $L^2(\O)$, and therefore that
$ u \in H^1(\O)$. This completes the proof of Theorem \ref{kolmh1}.
\eof

 \index{Compactness!-- results in $L^2$|)}

\chapter{Parabolic equations}\label{parabolic}

\section{Introduction}

The aim of this chapter is the study of finite volume schemes applied to a class 
of linear or nonlinear parabolic problems. We consider the following
transient diffusion-convection 
equation\index{Equation!Convection-diffusion --}:

\be
u_t(x,t) - \lapl \varphi(u)(x,t) + \div (\vbf u)(x,t) + b u (x,t) = f(x,t), \,  \, x \in
\O, \, t \in (0,T), \label{parabch}
\ee

where  $\O$ is an open polygonal bounded subset of $\R^d$, with $d=2$ or $d=3$,
$T >0$, $b \ge 0$, $\vbf\in\R^d$ is, for the sake of simplicity, 
a constant velocity field,
$f$ is a function defined on  $\O\times\R_+$ which represents a volumetric source term.
The function $\varphi$ is a nondecreasing Lipschitz continuous function,
 which arises in the modelling of 
general diffusion processes. A simplified version of Stefan's problem may 
be expressed with the formulation
\eqref{parabch} where $\varphi$  is  a continuous piecewise linear function,
which is constant on an interval. The porous medium
equation is also included in equation \eqref{parabch}, with $\varphi(u) = 
u^m$, $m>1$.
However, the linear case, i.e. $\varphi(u) = u$, is of full interest and
the error estimate of section \ref{parabest} will be given in such a case.
In section \refp{paranl}, we study the convergence of the explicit and of 
the implicit Euler scheme for the  nonlinear case with $\vbf = 0$ and $b = 0$.

\begin{remark} \normalfont
One could also consider a nonlinear convection term\index{Convection term}
 of the form 
 $\div (\vbf \psi(u))(x,t)$ where $\psi \in C^1(\R,\R)$. Such a nonlinear 
 convection term 
will be largely studied in the framework of nonlinear hyperbolic equations
(chapters \ref{hyper1d} and \ref{hypmd}) and we restrain here to a linear 
convection term for the sake of simplicity.
\end{remark}

An initial condition is given by

\be
u(x,0) = u_0(x),\, x \in \O. \label{parabcich}
\ee

Let  $\dr\O$ denote the boundary of $\O$, and let $\dr\O_d\subset\dr\O$ and 
$\dr\O_n\subset\dr\O$ such that $\dr\O_d \cup \dr\O_n = \dr\O$ and 
$\dr\O_d \cap \dr\O_n = \emptyset$. A Dirichlet boundary condition  is specified on 
$\dr\O_d\subset\dr\O$. Let   $g$ be  a real function defined on 
$\dr\O_d\times\R_+$, the Dirichlet boundary condition states that

\be
u(x,t) = g(x,t), \, x \in \dr \O_d, \, t \in (0,T).   \label{parabclch}
\ee

A Neumann boundary condition is given with a function $\tilde g$  defined
on $\dr\O_n\times\R_+$:

\be
-\grad \varphi(u)(x,t) \cdot \n(x) = \tilde g(x,t), \, x \in \dr \O_n, \, t \in 
(0,T), \label{parabc2ch}
\ee

where $\n$ is the unit normal vector to $\dr\O$, outward to $\O$.

\begin{remark} \normalfont Note that, formally, $\lapl \varphi(u)=\div (\varphi'(u) 
\grad u)$. Then, if $\varphi'(u)(x,t) = 0$ for some $(x,t) \in 
\O\times(0,T)$, the diffusion coefficient vanishes, so that 
Equation \eqref{parabch} is a ``degenerate" parabolic equation. In this 
case of degeneracy, the choice of the boundary conditions is important in 
order for the problem to be well-posed. In the case where 
$\varphi' $ is positive, the problem is always parabolic.
\end{remark}

In the next section, a finite volume scheme for the discretization of 
\eqref{parabch}-\eqref{parabc2ch} is presented. An error estimate in the 
linear case (that is $\varphi(u) = u$) is given in section \ref{parabest}.
Finally, a nonlinear (and degenerate) case is studied in section
\ref{paranl}; a convergence result is given for subsequences of sequences 
of approximate solutions, and, when the weak solution is unique, for the 
whole set of approximate solutions. A uniqueness result is therefore
proved for the case of a smooth boundary.

\section{Meshes and schemes} \label{parasch}

\hpage{Parabolic}

In order to perform a finite volume discretization of system 
(\ref{parabch})-(\ref{parabc2ch}), admissible meshes are used in a similar way to
the elliptic cases.
Let $\mesh$ be an admissible mesh of $\O$ in the sense of Definition 
\refp{meshdirichlet} with the additional assumption that any $\edge \in 
\edgesext$ is included in the closure of $\dr\O_d$ or included in the closure of 
$\dr \O_n$. The time discretization may be performed with a 
variable time step; in order to simplify the notations, we shall choose a 
constant time step $k \in (0,T)$. Let $\Nk \in \N^\star$ such 
that $\Nk = \max\{n \in \N, nk < T\}$, and we shall denote  $t_n 
= nk$, for $n \in \NNparpun$. Note that with a variable time step, error 
estimates and convergence results similar to that which are given 
in the next sections hold.

\medskip 
 Denote by $\{u_\cv^n,$ $\cv \in \mesh,$ $n \in 
\NNparpun\}$ the discrete unknowns; the value $u_\cv^n $ is an expected 
approximation of $u(\centercv,n k)$.

In order to obtain the numerical scheme, let us integrate formally Equation \eqref{parabch} over each control volume $\cv$ of $\mesh, 
$ and time interval $(nk, (n+1)k)$, for $n \in \NNpar$:

\be \bal
\dsp{ \int_\cv (u(x,t_{n+1}) - u(x,t_n))dx -
\int_{nk}^{(n+1)k} \int_{\dr\cv}  \grad \varphi(u)(x,t) \cdot \ncv(x)  
\dfrontiere(x) dt + } \\
\dsp{\int_{nk}^{(n+1)k} \int_{\dr\cv} \vbf \cdot \ncv(x) u(x,t) \dfrontiere(x) dt +
b \int_{nk}^{(n+1)k}\int_\cv u(x,t)dx dt =} 
\dsp\int_{nk}^{(n+1)k}\int_\cv f(x,t) dx dt.
\ean
\label{eqint}
\ee

where $\ncv$ is the unit normal vector to $\dr \cv$, outward to $\cv$.

\medskip

Recall that, as usual, the
stability \index{Stability} condition for an explicit discretization of a parabolic equation requires the time
step to be limited by a power two of the space step, which is generally too strong a
condition in terms of computational cost. Hence the choice of an implicit 
formulation in the left hand side of \eqref{eqint} which yields

\be \bal
\dsp{ \frac 1 k \int_\cv (u(x,t_{n+1}) - u(x,t_n))dx -
\int_{\dr\cv}  \grad \varphi(u)(x,t_{n+1}) \cdot \ncv(x)  
\dfrontiere(x) + } \\
\dsp{ \int_{\dr\cv} \vbf \cdot \ncv(x) u(x,t_{n+1}) \dfrontiere(x)  +
b \int_\cv u(x,t_{n+1})dx dt =} 
\dsp {\frac 1 k \int_{nk}^{(n+1)k}\int_\cv f(x,t) dx dt,} 
\ean
\label{eqinti}
\ee

There now remains to replace in Equation \eqref{eqint} each term by its 
approximation with respect to the discrete unknowns (and  the data). 
Before doing so, let us remark that another way to obtain \eqref{eqinti} is to 
integrate (in space) formally Equation \eqref{parabch} over each control volume $\cv$ 
of $\mesh$, at time $t \in (0,T)$. This gives  

\be \bal
\dsp{\int_\cv u_t(x,t)dx -
\int_{\dr\cv}  \grad \varphi(u)(x,t) \cdot \ncv(x)  
\dfrontiere(x) + } \\
\dsp{ \int_{\dr\cv} \vbf \cdot \ncv(x) u(x,t) \dfrontiere(x)  +
b \int_\cv u(x,t)dx =} 
\dsp {\int_\cv f(x,t) dx.} 
\ean
\label{eqintrt}
\ee

An implicit time discretization is then obtained by taking  $t=t_{n+1}$  
in the left hand side 
of \eqref{eqintrt},  and replacing $u_t(x,t_{n+1})$ by 
$(u(x,t_{n+1}) - u(x,t_n)) / k$. For the right hand side of \eqref{eqintrt} a mean 
value of $f$ between $t_n$ and $t_{n+1}$ may be used. This gives 
\eqref{eqinti}. It is also possible to take $f(x,t_{n+1})$ in the right hand side of 
\eqref{eqintrt}. This latter choice is simpler for the  proof of some 
error estimates 
(see Section \ref{parabest}).
 
~

Writing the
approximation of the various terms in Equation \eqref{eqinti} with respect to the discrete unknowns 
(namely, $\{u_\cv^n,$ $\cv \in \mesh,$ $n \in \NNparpun\}$) and taking into account the initial and
boundary conditions yields the  following implicit finite volume 
scheme\index{Finite volume scheme!-- for 
parabolic problems|(}\index{Scheme!Implicit -- for parabolic 
equations} for the
discretization  of (\ref{parabch})-(\ref{parabc2ch}), using the same 
notations and introducing some auxiliary unknowns  as in Chapter 
\ref{ellmd} (see equations \eqref{ellschema}-\refep{ellschemad}):

\be\ban
\dsp \mcv \frac {u_\cv^{n+1}- u_\cv ^{n}} {k} +  \sumedgecv
\Fcvedge^{n+1}  + \sumedgecv \vcvedge \uedgep^{n+1} + \mcv b
u_\cv^{n+1} =  \mcv f^n_\cv, \\ \forall \cv \in \mesh, \,  \forall
n \in \NNpar,
\ean  \label{parab} \ee

with

\be
{\dcvedge} \Fcvedge^{n}  = - {\medge} \Bigl(\varphi(u_\edge^{n}) -
\varphi(u_\cv^{n})\Bigr)  \text{ for } \edge \in \edgescv, \text{ for } n 
\in \NNparun,
\label{parabflux} \ee

\be
F_{\cv,\edge}^{n}= - F_{\cvv,\edge}^{n}
\text{ for all }\edge \in \edgesint\text{ such that }\edge = \edgecvcvv, \text{ for } n 
\in \NNparun,
\label{parabcons}\ee

\be \Fcvedge^{n} =  \frac 1 {k} \int_{(n-1)k}^{nk} \int_\edge \tilde g(x,t)
\dfrontiere(x)dt\  \text{ for } \edge \in \edgescv \text{ such that  } \edge \subset \dr \O_n,  \text{ for } n 
\in \NNparun, \label{parabc2} \ee

and

\be u_\edge^n =  g(\centeredge, n k)\   \text{ for } \ \edge \subset \dr \O_d,
\text{ for } n 
\in \NNparun, \label{parabcl}\ee

The upstream\index{Upstream} choice for the convection term
\index{Convection term} is performed as in the 
elliptic case (see page \pageref{easyvf1}, recall that $\vcvedge = 
\medge \vbf. \ncvedge$),

\be
\uedgep^{n} =\ba u_\cv^{n}, \text{ if } \vbf \cdot \ncvedge \ge 0, \\
                 u_\cvv^{n}, \text{ if } \vbf \cdot \ncvedge < 0,
\ea \text{ for all }\edge \in \edgesint\text{ such that }\edge = 
\edgecvcvv, \label{parabconv} \ee

\be
\uedgep^{n}= \ba u_\cv^{n}, &\text{ if } \vbf \cdot \ncvedge \ge 0, \\
                u_\edge^n 
&\text{ if } \vbf \cdot \ncvedge < 0,
\ea \text{ for all }\edge \in \edgescv \text{ such that }\edge \subset
\dr\O. \label{parabconvb} \ee

Note that, in the same way as in the elliptic case, the unknowns $u_\edge^{n+1}$
may be  eliminated using \eqref{parabflux}-\eqref{parabcl}.
There remains to define the right hand side, which may be defined by:

\be
f^n_\cv = \frac 1 {k\ \mcv} \int_{nk}^{(n+1)k} \int_\cv f(x,t) dx dt, \\ \forall \cv \in \mesh, \,  \forall
n \in \NNpar,
\label{parabvfnp}
\ee

or by:

\be
f^n_\cv = \frac 1 {\mcv}  \int_\cv f(x,t_{n+1}) dx, \\ \forall \cv \in \mesh, \,  \forall
n \in \NNpar.
\label{parabvfnq}
\ee

Initial conditions can be taken into account by different ways, 
depending on the regularity of the
data $u_0$. For example, it is possible to take

\be
u_\cv^0 = \frac {1} {\mcv} \int_\cv u_0(x) dx, \,  \cv \in \mesh, \label{parnlcondini}
\ee

or

\be
u_\cv^0 = u_0(\centercv),  \cv \in \mesh. \label{parabvf3}
\ee

\begin{remark} \normalfont
It is not obvious to prove that the implicit finite volume scheme 
\eqref{parab}-\eqref{parabconvb}  (with \eqref{parabvfnp} or \eqref{parabvfnq} 
and \eqref{parnlcondini} or \eqref{parabvf3}) has a solution. Once
the unknowns $\Fcvedge^{n+1}$ are eliminated, 
a nonlinear system of equations has to be solved\index{Finite volume 
scheme!-- for 
parabolic problems|)}.
A proof of the existence and uniqueness of a solution to this system is proved in the next section
for the linear case, and is sketched in Remark \ref{parnleximp} for the 
nonlinear case.
\end{remark}

\begin{remark}[Comparison with finite differences]\index{Finite difference method|(} \normalfont
Let us  consider the case of the heat equation, that is the case where $\vbf = 0$, $b=0$, $\varphi(s)=$ $s$ for all $s \in \R$, with 
Dirichlet condition on the whole boundary ($\dr\O_d =$ $\dr\O$). 
If the mesh consists in rectangular control volumes with constant space step in each direction, then the discretization obtained with the finite volume method gives (as in the case of the Laplace operator), the same scheme than the one obtained with the finite difference method (for which the discretization points are the centers of the elements of $\mesh$) except at the boundary. 
In the general  nonlinear case, finite difference methods\index{Finite difference method|)} have been used in \citet{ATT}, \citet{KAM} and \citet{MEY}, for example.

\end{remark}

\begin{remark}[Comparison with mass-lumped finite element]\index{Finite element method|(} \label{rem:mass-lumping} 
     \index{Mass lumping}\normalfont
Finite element  methods are classically used for elliptic or parabolic problems, see for instance   \citet{AG} or \citet{CIA}.
 Let $\meshprim$ be a finite element mesh of $\O$, consisting of triangles (see e.g. \citet{ciarlet} for the conditions on the triangles), with $N$ internal nodes.
A finite element formulation for (\ref{parabch}), with the implicit Euler scheme\index{Scheme!Implicit Euler --} in time, yields
\[
	\frac 1 k \Bigl(\int_\O (u^{n+1} (x) - u^n (x)) \phi_i(x) dx \Bigr) + \int_\O \grad u^{n+1}(x) \cdot \grad \phi_i(x) dx = \int_\O f(x,t_{n+1}) \phi_i(x) dx,
\]
where $\phi_i$ is the shape function of the finite element basis, associated with node $i$, for $i = 1, \ldots N$. 
Let us approximate $u^n$ by the following  Galerkin expansion\index{Galerkin expansion}:  
\[
u^{n+1} = \sum_{j =1}^{\overline{N}}  u_j^{n+1} \phi_j \text{ and } u^{n} = \sum_{j =1}^{\overline{N}}  u_j^{n} \phi_j,
\]
where $\overline{N}$ is the total number of nodes, and $u_j^n$ is expected to be an approximation of $u$ at time  $t_n$ and node $j$, for all $j$ and $n$; replacing in the above  equation, this yields:
\be
\frac 1 k \sum_{j = 1}^{\overline{N}} \int_\O (u_j^{n+1}-u_j^n) \phi_j(x) \phi_i(x) dx + \sum_{j = 1}^{\overline{N}} \int_\O u_j^{n+1} \grad \phi_j(x) \cdot \grad \phi_i(x) dx =  \int_\O f(x,t_{n+1}) \phi_i(x) dx.
\label{efp}
\ee
Hence, the finite element formulation yields, at each time step, a linear system
of the form $ C U^{n+1} + A U^{n+1} = B$ (where $U^{n+1}=(u_1,\ldots,u_N)^t$,
and $A$ and $C$ are $N \times N$ matrices); 
this scheme, however, is generally used after a mass-lumping\index{Mass lumping}, i.e. by assigning to the diagonal term of $C$ the sum of the coefficients of the corresponding line and setting the extra-diagonal terms to 0, thereby transforming $C$ into a diagonal matrix; we already saw in section \ref{lfefv} that the part $A U^{n+1}$ may be seen as a linear system derived from a finite volume formulation over the associated Vorono\"{\i} cells.  
With the mass lumping technique, the term $C_{i,i}$ corresponding to the $i$-th node is in fact equal to the integral $\int_K \phi_i$ of the $i$-th shape function $\phi_i$; since for an element $K$ whose vertices contain the $i$-th node, one has $\int_K \phi_i = \frac 1 3 | K |$,  therefore the integral $\int_K \phi_i$ is also equal to the area of the so-called Donald dual cell \index{Donald dual cell}, which is the dual cell around the $i$-th node obtained by joining the barycenter of each cell around the node to the middle of its edges.  
The term $C U^{n+1}$ may thus be  interpreted as a discretization by  a finite volume scheme over this Donald dual cell. 
\index{Finite element method|)} 
\end{remark}

\section{Error estimate for the linear case}\label{parabest}

\hpage{Parabolic}

We consider, in this section, the linear case, $\varphi (s) = s$ for all $s 
\in \R$, and assume   $\dr \O_d = \dr \O$, i.e. that a Dirichlet boundary condition
is given on the whole boundary, in which case Problem 
\eqref{parabch}-\eqref{parabc2ch} becomes
\[
u_t(x,t) - \lapl  u(x,t) + \div (\vbf u)(x,t) + b u (x,t) = f(x,t), \,  \, x \in
\O, \, t \in (0,T), 
\]
\index{Equation!Convection-diffusion --}
\[
u(x,0) = u_0(x),\, x \in \O,
\]

\[
u(x,t) = g(x,t), \, x \in \dr \O, \, t \in (0,T);   
\]

the finite volume scheme \eqref{parab}-\eqref{parabconvb} then becomes, 
assuming, for the sake of simplicity, that $\centercv \in \cv$ for all 
$\cv \in \mesh$,

\begin{multline}
 \mcv \frac {u_\cv^{n+1}- u_\cv ^{n}} {k} +  \sumedgecv \Fcvedge^{n+1}  + \sumedgecv \vcvedge \uedgep^{n+1} + \mcv b
u_\cv^{n+1} =  \mcv f^n_\cv,
\\
\forall \cv \in \mesh, \,  \forall n \in \NNpar,
 \label{parablin}
 \end{multline}

with

\be
\Fcvedge^{n}  = - \transcvcvv (u_\cvv^{n} -
u_\cv^{n}) \text{ for all }\edge \in \edgesint \text{ such that }\edge = \edgecvcvv,
 \text{ for } n 
\in \NNparun,
\label{parablinflux} \ee

\be \Fcvedge^{n} =  - \transedge (g(\centeredge, n k) -u_\cv^{n}) \text{ for all  }\   
\edge  \in \edgescv  \text{ such that } \edge \subset \dr \O,  \text{ for } n 
\in \NNparun, \label{parablincl} \ee

and

\be\ba
\uedgep^{n} = u_\cv^{n}, \text{ if } \vbf \cdot \ncvedge \ge 0, \\
\uedgep^{n} = u_\cvv^{n}, \text{ if } \vbf \cdot \ncvedge < 0,
\ea \text{ for all }\edge \in \edgesint\text{ such that }\edge = 
\edgecvcvv, \label{parablinconv} \ee

\be
\ba
\uedgep^{n} = u_\cv^{n}, &\text{ if } \vbf \cdot \ncvedge \ge 0, \\
\uedgep^{n} = g(\centeredge, n k),
&\text{ if } \vbf \cdot \ncvedge < 0,
\ea \text{ for all }\edge \in \edgescv \text{ such that }\edge \subset
\dr\O.
\label{parablinconvb}
\ee

The source term and initial condition  $f$ and $u_0$ are  discretized by \eqref{parabvfnq} and  \eqref{parabvf3}. 

\smallskip

A convergence analysis of a one-dimensional vertex-centered scheme was performed in \citet{guo} by writing the scheme in a finite element framework. 
Here we shall use direct finite volume techniques which also handle the multi-dimensional case. 

The following theorem gives an $L^\infty$ estimate (on the approximate 
solution) and an error estimate. Some easy generalizations are possible 
(for instance, the same theorem holds with $b<0$, the only difference is 
that in the $L^\infty$ estimate \eqref{parinfty} the bound $c$ also depends 
on $b$). 

\begin{theo}\label{parabestim}\index{Estimate on the approximate solution!-- for a parabolic 
equation}\index{Error estimate!-- for a parabolic 
equation}~
Let $\Omega$ be an open polygonal bounded subset of $\R^\dim$, $T >0$,
$u\in C^2(\Obar \times \R_+, \R)$,  $b \ge 0$ and  $\vbf \in \R^\dim$.
Let $u_0 \in C^2(\Obar,\R)$ be defined by $u_0 = u(\cdot,0)$,
let $f \in C^0(\Obar \times \R_+, \R)$ be defined by
$f = u_t - \div (\grad u) + \div (\vbf u) + b u$
and $g\in C^0(\dr \O \times \R_+, \R)$ defined by $g = u$ on $\dr\O \times \R_+$.
Let $\mesh$ be an admissible mesh in the sense of Definition 
\refp{meshdirichlet} and $k \in (0,T)$. Then there exists a unique
vector $(u_\cv)_{\cv \in \mesh}$ satisfying 
\eqref{parablin}-\eqref{parablinconvb} (or \eqref{parab}-\eqref{parabconvb}) 
with \eqref{parabvfnq} and 
\eqref{parabvf3}. There exists $c$ only depending  on 
$u_0$, $T$, $f$ and $g$ such that

\be
\sup\{\vert \ucvn \vert, \cv \in \mesh, n \in \NNparun\} \le c.
\label{parinfty}
\ee 

Furthermore, let $\errorcv^n =
u(\centercv,t_n) - u^n_\cv$, for $\cv \in \mesh$ and $n \in \NNparun$, 
and  $h  = \size(\mesh)$. Then there exists $C \in \R_+$ only depending 
on $b$, $u$, $\vbf$, $\Omega$ and $T$ such that

\be
(\summesh (\errorcv^n)^2 \mcv)^{\frac 1 2}\le C (h+k), \; \forall \, n \in \NNparun.
\label{estdeuxp}
\ee

\end{theo}

\proof of Theorem \ref{parabestim}

\smallskip

For simplicity,  let us  assume   that $\centercv \in \cv$ for all $\cv \in \mesh$. Generalization without this condition is straightforward.

\smallskip

\textit{(i) Existence,  uniqueness, and $L^\infty$  estimate}

For a given $n \in \NNpar$, set $f_\cv^n = 0$ and $u_\cv^{n} = 0$ in \eqref{parablin}, and $g(\centeredge, (n+1)k) = 0$ for all $\edge \in 
\edges$ such that $\edge \subset \dr\O$. 
Multiplying \eqref{parablin}  by $\ucvnp$ and using the same technique as in the proof of Lemma \refp {ellexistu} yields that $ \ucvnp 
= 0 $ for all $\cv \in \mesh$.
 This yields the uniqueness of the solution $\{\ucvnp$, $\cv \in \mesh\}$  to \eqref{parablin}-\eqref{parablinconvb} 
for given $\{\ucvn$, $\cv \in \mesh\}$, $\{f_\cv^n$, $\cv \in \mesh\}$ and $\{g(\centeredge, (n+1)k)$, $\edge \in \edges, \edge \subset \dr\O_d\}$. 
The existence follows immediately, since \eqref{parablin}-\eqref{parablinconvb} is a 
finite dimensional linear system with respect to the unknown
$\{\ucvnp$, $\cv \in \mesh\}$ (with as many unknowns as equations). 

\smallskip

Let us now prove the estimate \eqref{parinfty}. 
Let $m_f=\min\{f(x,t)$, $x \in \Obar$, $t \in [0,2T]\}$ and  $m_g=\min\{g(x,t)$, $x \in \dr \O$, $t \in [0,2T]\}$. 
Let $n \in \NNpar$. Then, we claim that

\be
\min\{u^{n+1}_\cv, \, \cv \in \mesh\} \ge \min\{\min\{u^{n}_\cv, \, \cv \in \mesh\}+k m_f,0,m_g\}.
\label{claim1}
\ee

Indeed, if $\min\{u^{n+1}_\cv, \, \cv \in \mesh\} < \min\{0,m_g\}$, let 
$\cv_0 \in \mesh$ such that $u^{n+1}_{\cv_0}=\min\{u^{n+1}_\cv, \, \cv \in \mesh\}$. 
Let us write \eqref{parablin} with $\cv = \cv_0$.
 Since $u^{n+1}_{\cv_0} < 0$ and $u^{n+1}_{\cv_0}<m_g$, we get that $F_{K_0,\edge}^{n+1}\le 0.$
Moreover, since $\bfv$ is constant, we have $  \sumedgecv \vcvedge = 0,$ so that
\[
  \sumedgecv \vcvedge \uedgep^{n+1} =   \sumedgecv \vcvedge \left( \uedgep^{n+1} -u_{K_0}^{n+1} \right) \le 0;
\]
therefore 
\[
u^{n+1}_{\cv_0} \ge u^{n}_{\cv_0} + kf^n_{\cv_0} \ge \min\{u^{n}_\cv, \, \cv  \in \mesh\}+km_f,
\]
this proves \eqref{claim1}, which yields, by induction, that:

\[
\min\{u^{n}_\cv, \, \cv \in \mesh\} \ge \min\{\min\{u^{0}_\cv, \, \cv  \in \mesh\},0,m_g\} + nk\min\{m_f,0\}, \forall \, n \in \NNparpun.
\]

Similarly,
\[
\max\{u^{n}_\cv, \, \cv \in \mesh\} \le \max\{\max\{u^{0}_\cv, \, \cv  \in \mesh\},0,M_g\} + nk\max\{M_f,0\}, \forall \, n \in \NNparpun,
\]

with $M_f=\max\{f(x,t)$, $x \in \Obar$, $t \in [0,2T]\}$ and 
$M_g=\max\{g(x,t)$, $x \in \dr \O$, $t \in [0,2T]\}$.

This proves \eqref{parinfty} with $c= \norm{u_0}{L^\infty(\O)} + 
\norm{g}{L^\infty(\dr \O \times (0,2T))} + 2T \norm{f}{L^\infty(\O \times 
(0,2T))}$.

\medskip

\textit{(ii) Error estimate}

As in the stationary case (see the proof of Theorem \refp{ellesterr}), 
one uses the regularity of the
data and the solution to write an equation for the error
$e_\cv^n = u(\centercv,t_n) - u_\cv^n$, defined for $\cv \in \mesh$ and 
$n \in \NNparpun$. Note that $e_\cv^0=0$ for $\cv \in \mesh$.
Let $n \in \NNpar$. Integrating (in space) Equation \eqref{parabch} over each control volume $\cv$ 
of $\mesh$, at time $t=t_{n+1}$, gives, thanks to the choice of 
$f^n_\cv$ (see \eqref{parabvfnq}),

\be
\int_\cv u_t(x,t_{n+1})dx -
\int_{\dr\cv}  \Bigl(\grad u(x,t_{n+1}) - \vbf u(x,t_{n+1})\Bigr) \cdot \ncv(x) \dfrontiere(x)  +
b \int_\cv u(x,t_{n+1})dx =
\mcv f^n_\cv. 
\label{eqintrtl}
\ee

Note that, for all $x \in \cv$ and all $\cv \in \mesh$, a Taylor expansion yields, thanks to the regularity of $u$:
\[
u_t(x,t_{n+1})=(1/k)(u(\centercv,t_{n+1})-u(\centercv,t_{n})) + s^n_\cv(x)\text{ with } \vert s^n_\cv(x) \vert \le C_1(h+k)
\]
with some $C_1$ only depending  
on $u$ and $T$. Therefore, defining $ S^n_\cv  =\dsp{\int_\cv 
s^n_\cv(x) dx}$, one has: $\vert S^n_\cv \vert \le C_1 \mcv (h+k).
$

One follows now the lines of the proof of Theorem \refp{ellesterr}, 
adding the terms due to the time derivative $u_t$.
Substracting
\eqref{parablin} to  \eqref{eqintrtl} yields

\be\bal
\dsp{\mcv \frac {e^{n+1}_\cv - e^{n}_\cv} k + \sumedgecv \Bigl(\Gcvedge^{n+1} 
+ \Wcvedge^{n+1} \Bigr) + b \mcv e^{n+1}_\cv = }
\\
\dsp{ b \mcv
\rho^{n}_\cv  - \sumedgecv \medge   (\Rcvedge^{n} + \rcvedge^{n}) - 
S^n_\cv,\, \forall \cv \in \mesh,}
\ean
\label{diffea}
\ee

where (with the notations of Definition \refp{meshdirichlet}),

\smallskip

\[
\bal
\Gcvedge^{n+1} = -\transedge (e^{n+1}_\cvv -e^{n+1}_\cv), \, \forall \cv 
\in \mesh, \, \forall  
\edge \in \edgescv \cap \edgesint, \, \edge=\edgecvcvv,
\\
\Gcvedge^{n+1} = \transedge e^{n+1}_\cv, \, \forall \cv \in \mesh, \, \forall  \edge \in \edgescv 
\cap \edgesext,
\\
\Wcvedge^{n+1} = \medge \vbf \cdot \ncvedge (u(x_{\edge,+},t_{n+1}) 
-u^{n+1}_{\edge,+}),
\ean
\]

where $x_{\edge,+} = \centercv$ (resp. $\centercvv$) if $\edge \in \edgesint$, 
$\edge = \edgecvcvv$ and $\vbf \cdot \ncvedge \ge 0$ (resp. $\vbf \cdot 
\ncvedge < 0$) and   
$x_{\edge,+} = \centercv$ (resp. 
$\centeredge$) if $\edge = \edgescv \cap 
\edgesext$ and $\vbf \cdot \ncvedge \ge 0$ (resp. $\vbf \cdot \ncvedge < 
0$),

\[ \rho_\cv^n = u (x_\cv, t^{n+1}) - \frac 1 \mcv \int_\cv  u(x,t_{n+1}) dx,\]
\[ \medge \Rcvedge^{n} = \transedge(u (x_\cv, t^{n+1}) - u (x_\cvv, t^{n+1})) + 
\int_\edge \grad u (x,t_n). \ncvedge \dfrontiere(x)
\text{ if } \edge = \edgecvcvv \in \edgesint,\]
\[  \medge \Rcvedge^{n} = \transedge(u (x_\cv, t^{n+1}) - g(\centeredge, t^{n+1}) + 
\int_\edge \grad u (x,t_n). \ncvedge \dfrontiere(x)
\text{ if } \edge  \in \edgescv \cap \edgesint,\]
and 
\[\medge \rcvedge^{n} =  \vbf \cdot \ncvedge (\medge u(x_{\edge,+},t_{n+1}) - 
\int_\medge u(x,t_{n+1}) \dfrontiere(x), \text{ for any } \edge \in \edges.\]

~

As in Theorem \ref{ellesterr}, thanks to the regularity of $u$, 
there exists $C_2$, only depending on $u$, 
$\vbf$ and $T$, such that 
$\vert \Rcvedge^{n} \vert  + \vert \rcvedge^{n} \vert  \le C_2h$ and 
$\abs {\rho^n_\cv} \le C_2h$, for 
any $\cv \in \mesh$ and $\edge \in \edgescv$.

~

Multiplying \eqref{diffea} by $e^{n+1}_\cv$, summing for $\cv\in\mesh$, 
and performing the same computations as in  the proof of Theorem \ref{ellesterr} 
between \eqref{ellraslebol1} to \refep{ellraslebol4} yields, with some $C_3$ 
only depending  on $u$, $\vbf$, $b$, $\O$ and $T$,

\be
\bal
\dsp{ \frac 1 k }
\summesh \mcv (e_\cv^{n+1})^2 + \frac 1 2 \Vert e^{n+1}_{\mesh} \Vert_{1,\mesh}^2 
+ \frac 1 2 b \Vert e^{n+1}_{\mesh} \Vert_{L^2(\O)}^2 \le
\\
C_3h^2 + C_1 (h+k) 
\summesh \mcv  \vert e^{n+1}_\cv  \vert+ \frac 1 k \summesh \mcv 
e_\cv^{n+1}e_\cv^{n},
\ean
\label{finesterr}
\ee

where the second term of the right hand side is due to the bound on 
$S^n_\cv$ and  where $e^{n+1}_{\mesh}$ is a piecewise constant function defined by
 
\[
e^{n+1}_{\mesh}(x) = e^{n+1}_{\cv}, \text{ for }  x \in \cv,  \cv \in 
\mesh.
\]
Inequality \eqref{finesterr} yields

\[
\Vert e^{n+1}_{\mesh} \Vert_{L^2(\O)}^2 \le 2kC_3h^2 + 2k C_1 \mesO (k+h)
\Vert e^{n+1}_{\mesh} \Vert_{L^2(\O)} + \Vert e^{n}_{\mesh} \Vert_{L^2(\O)}^2,
\]

which gives

\be
\Vert e^{n+1}_{\mesh} \Vert_{L^2(\O)}^2 \le 
\Vert e^{n}_{\mesh} \Vert_{L^2(\O)}^2
  + C_4 \bigl( kh^2 + k(k+h) \Vert e^{n+1}_{\mesh} \Vert_{L^2(\O)}
   \bigr), \label{An}
\ee

where $C_4\in \R_+$
only depends on $u$, $\vbf$, $b$, $\O$ and $T$.
Remarking that for $\eps >0$, the following inequality holds:

\[
C_4k(k+h)\Vert e^{n+1}_{\mesh} \Vert_{L^2(\O)} 
\le \eps^2\Vert e^{n+1}_{\mesh} \Vert_{L^2(\O)}^2 + (1/\eps^2) C_4^2 k^2 
(k+h)^2,\]

taking  $\eps^2=k/(k+1)$,  \eqref{An} yields

\be
\Vert e^{n+1}_{\mesh} \Vert_{L^2(\O)}^2  \le (1+k) 
\Vert e^{n}_{\mesh} \Vert_{L^2(\O)}^2  + C_4 kh^2 (1+k) +  
(1+k)^2 C_4^2 k (k+h)^2.
\label{Anp}
\ee

Then, if $\Vert e^{n}_{\mesh} \Vert_{L^2(\O)}^2
\le c_n (h+k)^2$, with  $c_n \in \R_+$, one deduces from \eqref{Anp}, 
using $h\le h+k$ and $k < T$, that 
\[\Vert e^{n+1}_{\mesh} \Vert_{L^2(\O)}^2
\le c_{n+1} (h+k)^2 \text{ with } c_{n+1}=(1+k)c_n + C_5 k \text{ and }
C_5=C_4(1+T)+C_4^2(1+T)^2.\] (Note that $C_5$ only depends on $u$, $\vbf$, 
$b$, $\O$ and $T$).

Choosing $c_0=0$ (since $\Vert e^{0}_{\mesh} \Vert_{L^2(\O)}=0$), 
the relation between $c_n$ and 
$c_{n+1}$ yields (by induction) $c_n \le C_5 e^{2kn}$. Estimate 
(\ref{estdeuxp}) follows with $C^2=C_5e^{4T}$.
\eof

\begin{remark} \normalfont

The error estimate given in Theorem \ref{parabestim} 
may be generalized to the case of discontinuous
coefficients. The admissibility of the mesh is then redefined 
so that the data and
the solution are piecewise regular on the control volumes as in 
Definition \refp{rouenmesh}, see also \citet{rouenrh}.

\end{remark}

\medskip

\section{Convergence in the nonlinear case}\label{paranl}

\hpage{Parabolic}

\subsection{Solutions to the continuous problem}\label{parnlpbcont}

We consider Problem \eqref{parabch}-\eqref{parabc2ch} with $\vbf=0$, $b = 0$, 
$\dr\O_n = \dr\O$ and $\tilde g=0$, that is a homogeneous Neumann condition on the 
whole boundary, in which case the problem\index{Equation!Diffusion --} becomes

\be
u_t(x,t) - \lapl \varphi(u)(x,t) = f(x,t), \ \ \text{ for } (x,t)\in \O\times (0,T),
\label{paranlanf}\ee

with

\be
\grad \varphi(u)(x,t)\cdot \bfn(x) = 0, \ \ \text{ for } (x,t)\in \dr\O\times(0,T),
\label{paranlanb}\ee

and the initial condition

\be
u(x,0) = u_0(x), \ \ \text{ for all } x\in \O.
\label{paranlani}\ee

\separ

We suppose that the following hypotheses are satisfied:

\begin{hyp}\label{parnlH} ~

\begin{enumerate}
\item[(i)] $\O$ is an open bounded polygonal subset of $\R^\dim$ and $T > 
0$.
\item[(ii)] The function $\varphi\in C(\R,\R)$ is a nondecreasing locally 
Lipschitz continuous function.
\item[(iii)] The initial data $u_0$ satisfies $u_0\in L^\infty(\O)$.
\item[(iv)] The right hand side $f$ satisfies $f\in L^\infty(\O\times \Rpast)$.
\end{enumerate}
\end{hyp}

Equation (\ref{paranlanf}) is a degenerate parabolic equation. Formally, 
$\lapl \varphi (u) = \div (\varphi'(u) \grad u)$, so that, if 
$\varphi'(u) = 0$, the diffusion coefficient vanishes.
Let us give a definition of a weak solution
$u$ to Problem (\ref{paranlanf})-(\ref{paranlani}) (the proof of the existence
of such a solution is given in \citet{KAM}, \citet{LSU}, \citet{MEI}, 
\citet{Ol2}).

\begin{defi}\label{parnldefess}  \normalfont
Under Assumption \ref{parnlH}, a measurable function $u$ is a weak solution of
(\ref{paranlanf})-(\ref{paranlani}) if

\be\ban
u \in L^\infty(\O\times(0,T)),
\\
\dsp{\int_0^T\int_\O} \Bigl( u(x,t) \psi_t(x,t)
+ \varphi(u(x,t)) \Delta \psi(x,t) + f(x,t) \psi(x,t) \Bigr) \ dx \ dt \ + 
\\
\dsp{\int_\O} u_0(x) \psi(x,0) dx = 0, \text{ for all } \psi\in \mathcal{A}_T,
\ean\label{paranlanw}\ee

where $\AT = \{ \psi\in C^{2,1}(\overline{\O}\times[0,T])$, 
$\grad \psi \cdot \n= 0$ on $\dr\O\times [0,T]$,
and $\psi(\cdot,T)=0 \}$, and $C^{2,1}(\overline{\O}\times[0,T])$ denotes 
the set of functions which are  restrictions on 
$\overline{\O}\times[0,T]$ of functions from $\R^{\dim} \times \R$ into $\R$ 
which are twice (resp. once) continuously differentiable with 
respect to the first (resp. second) variable. (Recall that, as usual, $\n$ is the 
unit normal vector to $\dr \O$, outward to $\O$.)
\end{defi}

\begin{remark} \normalfont
It is possible to use a solution in a stronger sense, using only one integration
by parts for the space term. It then leads to a larger test function 
space than $\AT$.
\end{remark}
\begin{remark} \normalfont
Note that the function $u$ formally satisfies the conservation law
\be
\int_\O u(x,t) dx = \int_\O u_0(x) dx + \int_0^t\int_\O f(x,t) dx dt,
\label{parnlcsvt}\ee
for all $t\in [0,T]$. This property is also satisfied by the finite volume approximation.
\end{remark}

\subsection{Definition of the finite volume approximate solutions}
As in sections \refp{vfquatre} and \refp{msn}, an admissible mesh of $\O$ is defined, 
with respect to  which a functional 
space is introduced: this space contains the approximate solutions obtained from the finite 
volume discretization over the admissible mesh.

\begin{defi}\label{Xmeshk}  \normalfont 
Let $\O$ be an open bounded polygonal subset of $\R^\dim$, $\mesh$ 
be an admissible mesh in the sense of Definition \refp{meshneuman}, $T > 0$, $k 
\in (0,T)$ and $\Nk = \max \{n \in \N ; nk < T\}$. 
Let $X(\mesh,k)$ be the set of functions $u$
from $\O\times(0,(\Nk +1) k)$ to $\R$ such that there exists a family of real values
$\{u^n_\cv$, $\cv \in \mesh$, $n\in\NNpar\}$, with $u(x,t) = \ucvn$ for a.e. 
$x \in \cv$, $\cv \in \mesh$ and for
a.e. $t\in [nk,(n+1)k)$, $n\in\NNpar$.
\end{defi}

\separ

Since we only consider, for the sake of simplicity, a Neumann boundary condition,
we can easily eliminate the unknowns  $\Fcvedge^n$ located at the edges in 
equation \eqref{parab} using the equations \eqref{parabflux}, \eqref{parabcons},
and \eqref{parabc2}.
An explicit version of the scheme \index{Finite volume scheme!-- for 
parabolic problems|(} can then be written in the following way:

\be\bal
\dsp{ \mcv \frac {u_\cv^{n+1}- u_\cv ^{n}} {k} -
\sumv \transcvcvv  \Bigl(\varphi(\ucvvn) - \varphi(\ucvn)\Bigr)
=  \mcv f^n_\cv,} 
\\
\forall \cv \in \mesh, \,  \forall
n \in \NNpar.
\ean
\label{parnlschema}
\ee

\be
 u_\cv^0 = \frac {1} {\mcv} \int_\cv u_0(x) dx, \,  \forall \cv \in \mesh,
\label{parnlcondinj}\ee

\be
f^n_\cv = \frac 1 {k\  \mcv} \int_{n k}^{(n+1)k} \int_\cv f(x,t) dx dt,
\forall \cv \in \mesh, \, \forall n\in\NNpar.
\label{parnlcondsou}\ee

\separ

(Recall that $\dsp {\transcvcvv = \dsp \frac \medgecvcvv \dedgecvcvv}$, see 
Definition \refp{meshneuman}.)

\begin{remark} \normalfont
The definition using the mean value in \eqref{parnlcondinj} is motivated by 
the lack of  regularity assumed on the data $u_0$.
\end{remark}

The scheme (\ref{parnlschema})-(\ref{parnlcondsou})
is then used to build an approximate solution,
$\uapp \in X(\Mesh,k)$ by

\be 
\uapp(x,t)  = \ucvn, \forall x\in \cv,  \forall t\in [nk, 
(n+1)k),  
 \forall \cv\in\mesh,  \forall n\in\NNpar.
 \label{parnldefuapp}\ee

\begin{remark} \normalfont\label{parnleximp}
The implicit finite volume scheme \index{Scheme!Implicit -- for parabolic 
equations}  is defined by
\be
\ban
\dsp \mcv \frac {u_\cv^{n+1}- u_\cv ^{n}} {k} -
\sumv \transcvcvv \Bigl(\varphi(\ucvvnp) - \varphi(\ucvnp)\Bigr) =
\mcv f^n_\cv, \\ \forall \cv \in \mesh, \,  \forall
n \in \NNpar.
\ean
\label{parnlschemaimp}
\ee
\index{Finite volume scheme!-- for 
parabolic problems|)}

The proof of the existence of $\ucvnp$, for any $n\in\NNpar$, can be obtained
using the following fixed point method:

\be
u_\cv^{n+1,0} = \ucvn, \ \text{ for all } \cv\in\mesh,
\label{parnlpfixei}\ee
and

\be\ban
\dsp{\mcv \frac {u_\cv^{n+1,m+1}-\ucvn} k} - \sumv \transcvcvv
\Bigl(\varphi(u_\cvv^{n+1,m}) - \varphi(u_\cv^{n+1,m+1})\Bigr) =
\mcv f^n_\cv,\\ \ \forall \cv\in\mesh, \ \forall m\in\N.
\ean\label{parnlpfixem}\ee

Equation $(\ref{parnlpfixem})$ gives a contraction property,
which leads first to
prove that for all $\cv\in\mesh$, the sequence
$(\varphi(u_\cv^{n+1,m}))_{m\in\N}$ converges.
Then we deduce that $(u_\cv^{n+1,m})_{m\in\N}$ also converges.

We shall see further that all results obtained for the explicit
scheme are also true, with convenient adaptations, for the implicit scheme. 
The function $\uapp$ is then defined
by $\uapp(x,t) = \ucvnp, \text{ for all } x\in \cv, \text{ for all } t\in [nk, (n+1)k)$.
\end{remark}

The mathematical problem is to study, under Assumption
\ref{parnlH} and with a mesh in the sense of Definition \ref{meshneuman}, the convergence of
$\uapp$ to a weak solution of Problem
(\ref{paranlanf})-(\ref{paranlani}),
when $h=\sizemesh \tends 0$ and $k\tends 0$.
Exactly in the same manner as for the elliptic case, we shall use
estimates on the approximate solutions 
which are discrete versions of the  estimates  
which hold on the solution of the continous problem and which ensure the 
stability of the scheme.
We present the proofs in the case of the explicit
scheme and show in several remarks how they can be extended to the
case of the implicit scheme  (which is significantly easier to study).
The proof of convergence of the scheme uses 
a weak-$\star$  convergence property, as in \citet{CIA},
which is proved  in a general setting  in section \refp{nlparaws}.
For the sake of completeness,   the proof
of uniqueness of the weak solution
of Problem \eqref{paranlanf}-\eqref{paranlani}  is given for the
case of a regular boundary; this allows to prove that the whole sequence 
of approximate solutions converges to the weak solution of problem 
\eqref{paranlanf}-\eqref{paranlani}, in which case an admissible mesh for a 
smooth domain can easily be defined (see Definition \refp{parnlHT}).

\subsection{Estimates on the approximate solution}
\index{Estimate on the approximate solution!-- for a parabolic equation|(}
\index{Stability}

\subsubsection{Maximum principle}

\index{Discrete!-- maximum principle}
\begin{lemma}\label{parnlestinf}
Under Assumption \ref{parnlH}, let $\mesh$ be an admissible  mesh in 
the sense of Definition \refp{meshneuman} and $k \in (0,T)$. Let
$U = \normuzr + T \normef$,
$B = \dsp{ \sup_{-U \le x < y \le U} \frac
{\varphi(x) - \varphi(y)} {x - y}}$. Assume that the condition

\be
k \le \frac {\mcv} {\dsp{B \sumv} \transcvcvv}, \ \
\text{ for all } \cv\in\mesh,
\label{parnlcflt}\ee
is satisfied. Then the function $\uapp$ defined by 
(\ref{parnlschema})-(\ref{parnldefuapp}) verifies

\be
\norm{ \uapp } {L^\infty(\O\times(0,T))} \le U.
\label{parnlinfsta}
\ee
\end{lemma}

\proof of Lemma \ref{parnlestinf}

\smallskip

Let $n \in \NNparmun$ and assume $u^n_\cv \in [-U,+U]$ for 
all $\cv \in \mesh$.

Let $\cv\in\mesh$, Equation (\ref{parnlschema}) can be written as

\[
\left.\begin{array}{ll}
\ucvnp =  &\dsp{ \Bigl( 1 - \frac k {\mcv} \sumv \transcvcvv
\frac {\varphi(\ucvvn) - \varphi(\ucvn)} {\ucvvn - \ucvn} \Bigr) \ucvn}\ + \\
&\dsp{ \frac k {\mcv} \sumv \Bigl( \transcvcvv
\frac {\varphi(\ucvvn) - \varphi(\ucvn)} {\ucvvn - \ucvn} \Bigr) \ucvvn +
k \fcvn},
\end{array}\right .
\]

with the convention that $\dsp{\frac {\varphi(\ucvvn) - \varphi(\ucvn)} 
{\ucvvn - \ucvn}}=0$ if $\ucvvn - \ucvn=0$.

Thanks to the condition (\ref{parnlcflt}) and since $\varphi$ is 
nondecreasing, the following inequality can be deduced:

\[
\vert \ucvnp  \vert \le \sup_{\cvv\in\mesh} \vert \ucvvn \vert + k \normef.
\]

Then, since $\cv$ is arbitrary in $\mesh$,

\be
\sup_{\cv \in \mesh} \vert \ucvnp  \vert \le \sup_{\cvv \in \mesh} \vert \ucvvn \vert + k \normef.
\label{parnlinfsta1}
\ee

Using (\ref{parnlinfsta1}), an induction on $n$ yields, for $n\in\NNpar$, 
$\sup_{\cv \in \mesh} \vert \ucvn \vert \le \normuzr + nk \normef$,
which leads to Inequality (\ref{parnlinfsta}) since $\Nk k \le T$.
\eof

\begin{remark} \normalfont Assume that there exist $\a$, $\beta$, $\gamma \in 
\R_+^\star$ such that $\mcv \ge \a h^{\dim}$, $\meas(\dr \cv) \le \beta 
h^{\dim-1}$, for all $\cv \in \mesh$, and $\dcvcvv \ge \gamma h$, for all 
$\edgecvcvv \in \edgesint$ (recall that $h=\sizemesh$). Then, $k \le 
Ch^2$ with $C= (\a \gamma)/(B\beta)$ yields \eqref{parnlcflt}.
\end{remark}

\begin{remark} \normalfont
Let $(\mesh_n,k_n)_\nnn$ be a sequence of admissible meshes and time 
steps, and  $(u_{\mesh_n,k_n})_\nnn$ the associated sequence of
approximate  finite  volume solutions; then , thanks to   
(\ref{parnlinfsta}),
there exists a function
$u\in L^\infty(\O\times(0,T))$ and a subsequence of $(u_{\mesh_n,k_n})_\nnn$
which converges
to $u$ for the weak-$\star$ topology of $L^\infty(\O\times(0,T))$.
\end{remark}

\begin{remark} \normalfont
Estimate $(\ref{parnlinfsta})$ is also true, with $U = \normuzr + 2T 
\normefd$, for the implicit scheme,
because the fixed point method guarantees  $(\ref{parnlinfsta1})$ (with 
$\normefd$ instead of $\normef$ and until $n=\Nk$),
without any condition on $k$.
\end{remark}

\subsubsection{Space translates of approximate solutions}

Let us now define a seminorm, which is the discrete version of the seminorm
in the space $L^2(0,T;H^1(\O))$.

\begin{defi}[Discrete $L^2(0,T;H^1(\O))$ seminorm]
\index{Discrete!--  $L^2(0,T;H^1(\O))$ seminorm}
 \label{defh1dk}  \normalfont
Let $\O$ be an open bounded polygonal subset of $\R^\dim$, 
$\mesh$ an admissible finite volume mesh in the sense of Definition 
\refp{meshneuman}, $T>0$,  $k \in (0,T)$ and $\Nk = \max \{n \in \N ; nk < T\}$. 
For  $u \in X(\mesh,k)$, let the following seminorms be defined by:

\be
\snorm{u(\cdot,t)}{1,\Mesh}^2 =
\sum_{\edgecvcvv\in\edgesint} \transcvcvv (\ucvvn - \ucvn)^2, \text{ for 
a.e. } t \in (0,T) \text{  and } n= \max\{n \in \N; nk \le t\},
\label{normh1dsk}\ee

 and
\be
\snorm{u}{1,\Mesh,k}^2 =
\sumparTsk k \sum_{\edgecvcvv\in\edgesint} \transcvcvv (\ucvvn - \ucvn)^2.
\label{normh1dk}\ee
\end{defi}

Let us now state some preliminary lemmata to the 
use of Kolmogorov's theorem (compactness properties
in $L^2(\O\times(0,T))$) in the proof of convergence of the approximate 
solutions.

\begin{lemma}\label{lsteftrx}
Let $\O$ be an open bounded polygonal subset of $\R^d$, $\mesh$ an 
admissible mesh in the sense of Definition 
\refp{meshneuman}, $T>0$, $k \in (0,T)$ and
$u \in X(\mesh,k)$. For all $\eta \in \R^d$, let $\O_\eta$ be defined by
$\O_\eta = \{x\in\O,\,[x,x+\eta]\subset\O\}$. Then:

\be
\norm{u(\cdot+\eta,\cdot)-u(\cdot,\cdot)}{L^2(\O_\eta\times(0,T))}^2 \le
\vert u \vert_{1,\Mesh,k}^2 \vert \eta \vert
(\vert \eta \vert + 2 \, \size(\mesh)), \forall \eta \in \R^d,
\label{parnltrx}
\ee
\end{lemma}

\proof of Lemma \ref{lsteftrx}

\smallskip

Reproducing the proof of Lemma \refp{h1dtot} (see also the proof of 
\refep{cvneum2}), we get,
for a.e. $t\in(0,T)$:

\be
\norm{u(\cdot+\eta,t)-u(\cdot,t)}{L^2(\O_\eta)}^2 \le
\snorm{u(\cdot,t)} {1,\Mesh}^2 \vert \eta \vert
(\vert \eta \vert + 2 \, \size(\mesh)), \forall \eta \in \R^d.
\label{parnltrx1}
\ee

Integrating \eqref{parnltrx1} on $t\in(0,T)$ gives \eqref{parnltrx}. \eof

The set $\O_\eta$ defined in Lemma \ref{lsteftrx} verifies
$\O\setminus\O_\eta \subset \cup_{\edge \in \edgesext} \bar{\omega_{\eta,\edge}}$, 
with $\o_{\eta,\edge} = \{y-t\eta$, $y \in 
\edge$, $t \in [0,1]\}$. Then, $\meas(\O\setminus\O_\eta)\le  \vert \eta 
\vert \  \meas(\dr\O)$, since $\meas(\bar{\o}_\eta) \le \eta \medge$. 
Then, an immediate corollary of Lemma \ref{lsteftrx} is the following:

\begin{lemma}\label{lsteftrxt}
Let $\O$ be an open bounded polygonal subset of $\R^d$, $\mesh$ an 
admissible mesh in the sense of Definition 
\refp{meshneuman}, $T>0$, $k \in (0,T)$ and
$u\in X(\mesh,k)$. Let $\tilde{u}$ be defined by $\tilde{u}=u$
a.e. on $\O\times(0,T)$, and $\tilde{u}=0$  a.e. on
$\R^{d+1} \setminus \O\times(0,T)$. Then:

\be\ba
\norm{\tilde{u}(\cdot+\eta,\cdot)-\tilde{u}(\cdot,\cdot)}{L^2(\R^{\dim+1})}^2 \le
\vert \eta \vert \Bigl( \vert u \vert_{1,\Mesh,k}^2
(\vert \eta \vert + 2 \, \size(\mesh)) +
2 \meas(\dr\O) \norm{u}{L^\infty(\O\times(0,T))}^2 \Bigr),\\
\forall \eta \in \R^d.
\ea
\label{parnltrxt}
\ee
\end{lemma}

\begin{remark} \normalfont
Estimate \eqref{parnltrxt} makes use of the 
$L^\infty(\O\times(0,T))$-norm of $u\in X(\mesh,k)$. 
A similar estimate may be proved with the $L^2(\O\times(0,T))$-norm of 
$u$ (instead of the $L^\infty(\O\times(0,T))$-norm). 
Indeed, the right hand side  
of \eqref{parnltrxt} may be replaced by $C \eta (\vert u \vert_{1,\Mesh,k}^2 + 
\norm{u}{L^2(\O\times(0,T))}^2)$, where $C$  only depends on 
$\O$. This estimate is proved in Theorem \refp{nllconv} where it is
used for the convergence of numerical 
schemes for the Neumann problem (for which no  $L^\infty$ estimate on the 
approximate solutions is available).
The key  to its proof  is the ``trace lemma" 
\refp{nlltrace}.
\end{remark}

Let us now state the following lemma, which
gives  an estimate of the discrete 
$L^2(0,T;H^1(\O))$ seminorm of the nonlinearity.

\begin{lemma}\label{parnlestimeg}
Under Assumption \ref{parnlH}, let $\mesh$ be an admissible mesh in the sense of Definition 
\refp{meshneuman}. Let $\xi \in (0,1)$ and $k \in (0,T)$ such that
\be k \le (1-\xi) \frac {\mcv}
{B \sumv \transcvcvv}, \ \ \text{ for all } \cv\in\mesh.
\label{parnlcfl}
\ee

Let $\uapp\in X(\Mesh,k)$ be given by (\ref{parnlschema})-(\ref{parnldefuapp}). 

Let $U = \normuzr + T \normef$ and $B$ be the Lipschitz constant
of $\varphi$ on $[-U,U]$.
Then there exists  $F_1 \ge 0 $, which only depends on
$\O$, $T$, $\varphi$, $u_0$, $f$ and $\xi$ such that

\be
\snorm{\varphi(\uapp)} {1,\Mesh,k}^2 \le F_1.
\label{parnlestimx}\ee
\end{lemma}

\proof of lemma \ref{parnlestimeg}

\smallskip

Let us first remark that the condition \eqref{parnlcfl} is stronger than
(\ref{parnlcflt}). Therefore, the result of lemma \ref{parnlestinf}
holds, i.e. $\vert \ucvn \vert \le U$, for all $\cv\in\mesh$, $n\in\NNpar$.
Multiplying   equation (\ref{parnlschema}) by $k \ucvn$, and
summing  the result over $n\in\NNpar$ and $\cv\in\mesh$ yields:

\be
\left.\begin{array}{ll}
\sumparTsk \summesh \mcv (\ucvnp-\ucvn) \ucvn - &\\
\sumparTsk k \summesh \sumv \transcvcvv
\Bigl(\varphi(\ucvvn) - \varphi(\ucvn)\Bigr) \ucvn & =
\sumparTsk k \summesh \mcv \ucvn \fcvn.
\end{array}\right.
\label{parnl1e}\ee

In order to obtain a lower bound on the first term on the left hand side
of (\ref{parnl1e}), let us first remark that:

\be
(\ucvnp-\ucvn) \ucvn = \half (\ucvnp)^2 -  \half (\ucvn)^2 -
\half (\ucvnp - \ucvn)^2.
\label{parnl2e}
\ee

Now, applying (\ref{parnlschema}), using  Young's inequality, the following 
inequality is obtained:

\be
(\ucvnp - \ucvn)^2 \le k^2 (1+\xi)
\Bigl[ \Bigl( \frac {1} {\mcv} \sumv \transcvcvv
(\varphi(\ucvvn) - \varphi(\ucvn))\Bigr)^2 + \frac {(\fcvn)^2} {\xi} \Bigr].
\label{parnl3eb}\ee

which yields in turn, using the Cauchy-Schwarz  inequality:

\be
\bac
\dsp{(\ucvnp - \ucvn)^2 \le \frac {k^2} {\mcv^2} (1+\xi)
\Bigl[ \sumv \transcvcvv\Bigr]
\Bigl[ \sumv \transcvcvv
\Bigl(\varphi(\ucvvn) - \varphi(\ucvn)\Bigr)^2 \Bigr]} \\
\dsp{+ \frac {(1+\xi)(k \ \fcvn)^2} {\xi}}.
\ean
\label{parnl3e}
\ee

Taking   condition  \eqref{parnlcfl} into account gives:

\be
(\ucvnp - \ucvn)^2 \le
(1 - \xi^2) \frac {k} {B \mcv}
\Bigl[ \sumv \transcvcvv \Bigl(\varphi(\ucvvn) - \varphi(\ucvn)\Bigr)^2\Bigr] +
\frac {(1+\xi)(k\ \fcvn)^2} {\xi}.
\label{parnl4e}\ee

Using  (\ref{parnl2e})   and (\ref{parnl4e}) 
leads to the following lower bound on the first term of the left hand 
side of \eqref{parnl1e}:

\be
\ban 
\sumparTsk \summesh \mcv (\ucvnp-\ucvn) \ucvn \ge &
 \dsp{\half\summesh \mcv \Bigl((u_\cv^{\Nk+1})^2 - (\ucvz)^2\Bigr) }\\
& - \dsp{\frac {1 - \xi^2} {2 B} \sumparTsk k \summesh
\Bigl[ \sumv} \transcvcvv \Bigl(\varphi(\ucvvn) - \varphi(\ucvn)\Bigr)^2 \Bigr] \\
& - \dsp{\frac {k (1+\xi)} {2 \xi}} \sumparTsk k \summesh \mcv (\fcvn)^2.
\ean\label{parnl5e}
\ee

Let us now handle the second term on the left hand side of (\ref{parnl1e}). 
Let $\phi\in C(\R,\R)$ be defined by $\phi(x) = x \varphi(x) -
\dsp{\int_{x_0}^x \varphi(y) dy}$, where
$x_0\in\R$ is an arbitrary given real value. Then the following
equality holds:

\be
\phi(\ucvvn) - \phi(\ucvn) = \ucvn ( \varphi(\ucvvn) - \varphi(\ucvn) ) -
\int_\ucvn^\ucvvn (\varphi(x) - \varphi(\ucvvn)) dx.
\label{parnl6e}\ee

The following technical lemma  is used here and several times in the 
sequel:

\begin{lemma}\label{petitlemme}
Let $g: \R \tends \R $ be a monotone Lipschitz continuous function,
with a Lipschitz constant $G>0$. Then:

\be
\vert \int_c^d (g(x)-g(c)) dx \vert \ge
\frac 1 {2 G} (g(d) - g(c))^2, \ \ \forall c, d \in \R.
\label{estplemme}\ee
\end{lemma}

\proof of Lemma \ref{petitlemme}

\smallskip

In order to prove Lemma \ref{petitlemme}, we assume, for instance,
that $g$ is nondecreasing and $c < d$ (the
other cases are similar). Then, one has $g (s) \ge h(s)$, for all $s \in
[c,d]$, where $h(s)=g(c)$ for $s \in [c,d-l]$ and $h(s)=g(c)+(s-d+l)G$
for $s \in [d-l,d]$, with $lG=g(d)-g(c)$, and therefore:

\[
\int_c^d (g(s)-g(c))ds \ge \int_c^d (h(s)-g(c))ds=\frac l
2(g(d)-g(c))=\frac 1 {2G} (g(d)-g(c))^2,
\]

this completes the proof of Lemma \ref{petitlemme}.

It is interesting to notice that, for this proof, the fact that $g$ is Lipschitz continuous is not
necessary. We only use the fact that $g(s) \ge g(c)$  and $g(d)-g(s) \le G(d-s)$ for all $s \in [c,d]$
(indeed we use
$g(s) \ge g(c)$ for $s \in [c,d-l]$)  and $g(d)-g(s) \le G(d-s)$
for $s\in [d-l,d]$). \eof

Using Lemma \ref{petitlemme}, (\ref{parnl6e}) and the equality
$\summesh \sumv \transcvcvv(\phi(\ucvvn) - \phi(\ucvn)) = 0$ yields:

\be
- \sumparTsk k \summesh \sumv \transcvcvv \Bigl(\varphi(\ucvvn) - \varphi(\ucvn)\Bigr) \ucvn \ge
 \frac 1 {2 B}
\sumparTsk k \summesh \sumv \transcvcvv (\varphi(\ucvvn) - \varphi(\ucvn))^2.
\label{parnl7e}\ee

Since $k<T$ we deduce from  (\ref{parnlinfsta}) that the right hand side
of equation (\ref{parnl1e}) satisfies

\be
\vert \sumparTsk k \summesh  \mcv\ucvn\fcvn \vert
\le 2 T \meas(\O) U \normefd.
\label{parnl8e}\ee

Relations $k<T$, (\ref{parnl1e}), (\ref{parnl5e}), (\ref{parnl7e}) and (\ref{parnl8e}) lead to

\be\bac
\dsp{\frac {\xi^2} {2 B}}
\sumparTsk k \summesh \sumv
\transcvcvv (\varphi(\ucvvn) - \varphi(\ucvn))^2 \le 
\\
2 T \meas(\O) \normefd \Bigl(U +  \dsp{\frac {1+\xi} {2 \xi}} \normefd T 
\Bigr)+
\dsp{\frac 1 2 \meas(\O) \norm{u_0}{L^\infty(\O)}^2}
\ean\label{parnl9e}\ee

which concludes the proof of the lemma. \eof

\begin{remark} \normalfont
Estimate (\ref{parnlestimx}) also holds for the implicit scheme\index{Scheme!Implicit --
for parabolic equations} ,
without any condition on $k$. One multiplies (\ref{parnlschemaimp})
by $\ucvnp$: the last term on the right hand side of (\ref{parnl2e})
appears with the opposite sign,
which considerably simplifies the previous proof.
\end{remark}

\subsubsection{Time translates of approximate solutions}

In order to  fulfill the hypotheses of Kolmogorov's theorem, 
  the study of time translates must now be performed. The following 
  estimate holds:

\begin{lemma}\label{lsteftrt}
Under Assumption \refp{parnlH}, let $\mesh$ be an admissible mesh 
in the sense of Definition \refp{meshneuman} and $k \in (0,T)$. Let 
$\uapp\in X(\Mesh,k)$ be given by (\ref{parnlschema})-(\ref{parnldefuapp}). Let
$U= \norm{\uapp}{L^\infty(\O\times(0,T))}$ and $B$ be the Lipschitz constant
of $\varphi$ on $[-U,U]$. Then:

\be\ba
\norm{\varphi(\uapp(\cdot,\cdot+\tau))-\varphi(\uapp(\cdot,\cdot))}{L^2(\O\times(0,T-\tau))}^2 \le
\\
2B\tau\Bigl(\vert \varphi(\uapp) \vert_{1,\Mesh,k}^2 +
B T \meas(\O) U \normef\Bigr), \forall \tau\in(0,T).
\ea\label{parnltrt}
\ee
\end{lemma}

\proof of Lemma \ref{lsteftrt}

\smallskip

Let $\tau \in (0,T)$.
Since $B$ is the Lipschitz constant
of $\varphi$ on $[-U,U]$, $U= \norm{\uapp}{L^\infty(\O\times(0,T))}$ and 
$\varphi$ is nondecreasing, the following inequality holds:

\be
\int_{\O\times (0,T-\tau)} \Bigl(\varphi(\uapp(x,t+\tau))-\varphi(\uapp(x,t))\Bigr)^2 dx dt
\le B \int_0^{T-\tau} A(t) dt,
\label{parnltrt1}\ee

where, for almost every $t\in (0,T-\tau)$,

\[
A(t) =  \int_\O  \Bigl(\varphi(\uapp(x,t+\tau))-\varphi(\uapp(x,t))\Bigr)
\Bigl(\uapp(x,t+\tau)-\uapp(x,t)\Bigr) dx.
\]

Let $t\in(0,T-\tau)$. Using the definition of   $\uapp$ \eqref{parnldefuapp},
this may also be written:

\be
A(t) =  \summesh \mcv \Bigl(\varphi(u_\cv^{n_1(t)})-\varphi(u_\cv^{n_0(t)})\Bigr)
\Bigl(u_\cv^{n_1(t)}-u_\cv^{n_0(t)}\Bigr),
\label{parnltrt3}\ee

with $n_0(t)$, $n_1(t) \in \NNpar$ such that $n_0(t)k \le t <
(n_0(t)+1)k$ and $n_1(t)k \le t+\tau < (n_1(t)+1)k$.

Equality \eqref{parnltrt3} may be written as

\[
A(t) =  \summesh (\varphi(u_\cv^{n_1(t)})-\varphi(u_\cv^{n_0(t)}))
\Bigl( \sum_{n = n_0(t)+1}^{n_1(t)}\mcv (\ucvn - \ucvnm) \Bigr),
\]

which also reads 

\be
A(t) =  \summesh  (\varphi(u_\cv^{n_1(t)})-\varphi(u_\cv^{n_0(t)}))
\Bigl( \sum_{n=1}^{\Nk} \chi_n(t,t+\tau) \mcv (\ucvn - \ucvnm) \Bigr),
\label{parnltrt4}\ee

with $\chi_n(t,t+\tau)=1$ if $nk \in (t,t+\tau]$ and $\chi_n(t,t+\tau)=0$ 
if $nk \notin (t,t+\tau]$.

In \eqref{parnltrt4}, the order of 
summation between $n$ and $\cv$ is changed and  the scheme 
\eqref{parnlschema} is used. Hence,

\[
\bal
\dsp{A(t) = k \sum_{n=1}^{\Nk} \chi_n(t,t+\tau) \Bigr[ \summesh  
(\varphi(u_\cv^{n_1(t)})-
\varphi(u_\cv^{n_0(t)}))}
\\
\dsp{\Bigl( \sumv \transcvcvv (\varphi(\ucvvnm)-\varphi(\ucvnm)) + \mcv 
\fcvnm \Bigr)\Bigr]}.
\ean
\]

Gathering by edges, this yields:

\[\bal
\dsp{A(t) = k \sum_{n=1}^{\Nk}
\Bigl[  \sumE \transcvcvv (\varphi(u_\cv^{n_1(t)})-\varphi(u_\cvv^{n_1(t)})
-\varphi(u_\cv^{n_0(t)})+\varphi(u_\cvv^{n_0(t)}))}
\\
\dsp{(\varphi(\ucvvnm)-\varphi(\ucvnm)) +
\summesh (\varphi(u_\cv^{n_1(t)})-\varphi(u_\cv^{n_0(t)})) \mcv \fcvnm 
\Bigr] \chi_n(t,t+\tau)}.
\ean\]

Using the inequality $2 a b \le a^2 + b^2$, this yields:
\be
A(t) \le \half A_0(t) + \half A_1(t) + A_2(t) + A_3(t),
\label{parnltrt7}\ee

with

\[
A_0(t) = k \sum_{n=1}^{\Nk} \chi_n(t,t+\tau)
\big( \sumE \transcvcvv 
(\varphi(u_\cvv^{n_0(t)})-\varphi(u_\cv^{n_0(t)}))^2 \big),
\]

\[
A_1(t) = k \sum_{n=1}^{\Nk} \chi_n(t,t+\tau)
\big( \sumE \transcvcvv 
(\varphi(u_\cvv^{n_1(t)})-\varphi(u_\cv^{n_1(t)}))^2 \big),
\]

\[
A_2(t) = k \sum_{n=1}^{\Nk} \chi_n(t,t+\tau) \big( \sumE
\transcvcvv (\varphi(\ucvvnm)-\varphi(\ucvnm))^2 \big),
\]

and

\[
A_3(t) = k \sum_{n=1}^{\Nk} \chi_n(t,t+\tau)
\big( \summesh (\varphi(u_\cv^{n_1(t)})-\varphi(u_\cv^{n_0(t)})) \mcv
\fcvnm \big).
\]

Note that, since $t \in (0,T-\tau)$, $n_0(t) \in \{0,\ldots,\Nk\}$, 
and, 
for $m \in \{0,\ldots,\Nk\}$, $n_0(t)=m$ if and only if $t \in 
[mk,(m+1)k)$. Therefore,

\[
\int_0^{T-\tau} A_0(t) dt  \le \sum_{m=0}^{\Nk} \int_{mk}^{(m+1)k} k
\sum_{n=1}^{\Nk} \chi_n(t,t+\tau)
\big( \sumE \transcvcvv (\varphi(u_\cvv^{m})-\varphi(u_\cv^{m}))^2
\big)dt,
\]

which also reads 

\be
\int_0^{T-\tau} A_0(t) dt  \le \sum_{m=0}^{\Nk} k \int_{mk}^{(m+1)k}
\Bigl( \sum_{n=1}^{\Nk} \chi_n(t,t+\tau) \Bigr) dt \sumE \transcvcvv 
(\varphi(u_\cvv^{m})-\varphi(u_\cv^{m}))^2.
\label{parnltrt12}\ee

The change of variable $t=s+(n-m)k$ yields

\[
\int_{mk}^{(m+1)k} \chi_n(t,t+\tau)dt=\int_{2mk-nk}^{2mk-nk+k}
\chi_n(s+(n-m)k,s+(n-m)k+\tau)ds=\int_{2mk-nk}^{2mk-nk+k}
\chi_m(s,s+\tau)ds,
\]

then, for all $m \in \{0,\ldots, \Nk\}$,

\[
\int_{mk}^{(m+1)k} 
 \Big( \sum_{n=1}^{\Nk} \chi_n(t,t+\tau) \Big) dt \le \int_\R 
 \chi_m(s,s+\tau)ds = \tau,
\]

since $\chi_m(s,s+\tau)=1$ if and only if $mk \in (s,s+\tau]$ which is 
equivalent to $s \in [mk-\tau,mk)$.

Therefore \eqref{parnltrt12} yields

\be
\int_0^{T-\tau} A_0(t) dt \le \tau \vert \varphi(\uapp) \vert_{1,\Mesh,k}^2.
\label{parnltrt14}
\ee

Similarly:

\be
\int_0^{T-\tau} A_1(t) dt \le \tau \vert \varphi(\uapp) \vert_{1,\Mesh,k}^2.
\label{parnltrt14b}\ee

Let us now study the term $\int_0^{T-\tau} A_2(t) dt$: 

\be
\int_0^{T-\tau} A_2(t) dt  \le \sum_{n=1}^{\Nk} k  \sumE
\transcvcvv (\varphi(\ucvvnm)-\varphi(\ucvnm))^2
\int_0^{T-\tau} \chi_n(t,t+\tau)  dt.
\label{parnltrt15}\ee

Since $\int_0^{T-\tau} \chi_n(t,t+\tau)  \le \tau$ (recall that 
$\chi_n(t,t+\tau)=1$ if and only if $t \in [nk-\tau, nk)$),
the following inequality holds:

\be
\int_0^{T-\tau} A_2(t) dt  \le \tau \vert \varphi(\uapp) \vert_{1,\Mesh,k}^2.
\label{parnltrt16}\ee

In the same way:

\be\ban
\int_0^{T-\tau}  A_3(t) dt  &\le \dsp{\sum_{n=1}^{\Nk}} k  \big( \summesh \mcv 2 
B U \normef \big) 
\int_0^{T-\tau} \chi_n(t,t+\tau)  dt
\\ &\le \tau T \meas(\O) 2 B U \normef.
\ean\label{parnltrt17}\ee

Using   inequalities \eqref{parnltrt1}, \eqref{parnltrt7} and 
\eqref{parnltrt14}-\eqref{parnltrt17},  
\eqref{parnltrt} is proved. \eof

\begin{remark} \normalfont
Estimate (\ref{parnltrt}) is again true for the implicit 
scheme\index{Scheme!Implicit -- for parabolic equations} ,
with $\normefd$ instead of $\normef$.
\end{remark}

An immediate corollary of Lemma \ref{lsteftrt} is the following.

\begin{lemma}\label{lsteftrtt}
Under Assumption \refp{parnlH}, let $\Mesh$ be an admissible mesh in the sense of
Definition \refp{meshneuman} and $k \in (0,T)$. Let $\uapp\in X(\Mesh,k)$ be given by 
(\ref{parnlschema})-(\ref{parnldefuapp}). Let
$U= \norm{\uapp}{L^\infty(\O\times(0,T)}$ and $B$ be the Lipschitz constant
of $\varphi$ on $[-U,U]$. One defines $\tilde{u}$ by $\tilde{u}=\uapp$
a.e. on $\O\times(0,T)$, and $\tilde{u}=0$  a.e. on
$\R^{d+1} \setminus \O\times(0,T)$. Then:

\[\ban
\norm{\varphi(\tilde{u}(\cdot,\cdot+\tau))-\varphi(\tilde{u}(\cdot,\cdot))}
{L^2(\R^{\dim+1})}^2 \le
2 \vert\tau\vert B  \Bigl( &\vert \varphi(\uapp) \vert_{1,\Mesh,k}^2 + \\&
B T \meas(\O) U \normef  + B \meas(\O) U^2\Bigr),\\ \forall \tau\in\R.
\ean
\]

\end{lemma}
\index{Estimate on the approximate solution!-- for a parabolic equation|)}

\subsection{Convergence}

\begin{theo}\label{parnlcvgce}
\index{Convergence of the finite volume method!-- for a nonlinear parabolic 
equation}
Under Assumption \refp{parnlH}, let $U = \normuzr + T \normef$ and

\[
B = \dsp{ \sup_{-U \le x < y \le U} \frac
{\varphi(x) - \varphi(y)} {x - y}}.
\]

Let $\xi\in (0,1)$ be a given real value. For $m \in \N$, let $\Mesh_m$ be an 
admissible mesh in the sense of
Definition \refp{meshneuman} and $k_m \in (0,T)$ satisfying the condition 
\eqref{parnlcfl} with $\mesh=\mesh_m$ and $k=k_m$. Let $\uappm$ be given by
(\ref{parnlschema})-(\ref{parnldefuapp}) with $\mesh=\mesh_m$ and 
$k=k_m$. Assume that $\size(\mesh_m) \tends 0$ as $m 
\tends \infty$.

Then, there exists a subsequence of  the sequence of approximate 
solutions, still denoted by $(\uappm)_{m \in \N}$, which converges to a weak solution $u$ of
Problem (\ref{paranlanf})-(\ref{paranlani}), as
$m \tends \infty$, in the following sense:

(i) $\uappm$ converges to $u$ in 
$L^\infty(\O\times(0,T))$, for the weak-$\star$ topology as $m$ tends to $+ \infty$,

(ii) $(\varphi(\uappm))$ converges to $\varphi(u)$ 
in $L^1(\O\times(0,T))$ as $m$ tends to $+ \infty$,

where  $\uappm$ and $\varphi(\uappm)$  also denote the restrictions  of 
these functions to $\O \times 
(0,T)$.

\end{theo}

\proof of Theorem \ref{parnlcvgce}

\smallskip

Let us set $u_m=\uappm$ and assume, without loss of generality, that 
$\varphi(0)=0$. First remark that, by \eqref{parnlcfl}, $k_m \tends 0$ as $m \tends 0$.
Thanks to Lemma \refp{parnlestinf}, the sequence $(u_m)_{m \in \N}$ is 
bounded in $L^\infty(\O \times (0,T)$. Then, there exists a subsequence, 
still denoted by $(u_m)_{m \in \N}$, such that $u_m$ converges, as $m 
\tends \infty$, to $u$ in 
$L^\infty(\O\times(0,T))$, for the weak-$\star$ topology.

For the study of the sequence $(\varphi(u_m))_{m \in \N}$,
we shall apply Theorem \refp{Kolm} with $N=\dim+1$, $q=2$, $\o=\O 
\times (0,T)$ and $p(v)=\tilde v$ with $\tilde v$ defined, as usual, by $\tilde v=v$ on 
$\O \times (0,T)$ and $\tilde v=0$ on $\R^{\dim+1} \setminus \O \times 
(0,T)$.

The first and second items of Theorem 
\ref{Kolm} are clearly satisfied; let us prove hereafter that the third 
is also satisfied.
By Lemma \ref{parnlestimeg}, the sequence 
$(\snorm{\varphi(u_m)}{1,\mesh_m,k_m})_{m \in \N}$ is bounded. 
Let $\eta \in \R^\dim$ and $\tau \in \R$, since

\[
\bal
\norm{\varphi(\tilde{u}_m(\cdot+\eta,\cdot+\tau))-\varphi(\tilde{u}_m(\cdot,\cdot))}
{L^2(\R^{\dim+1})} \le
\\
\norm{\varphi(\tilde{u}_m(\cdot+\eta,\cdot))-\varphi(\tilde{u}_m(\cdot,\cdot))}
{L^2(\R^{\dim+1})}+
\norm{\varphi(\tilde{u}_m(\cdot,\cdot+\tau))-\varphi(\tilde{u}_m(\cdot,\cdot))}
{L^2(\R^{\dim+1})},
\ean
\]

lemmata \ref{lsteftrxt} and \ref{lsteftrtt} give the third item of 
Theorem \ref{Kolm} and this yields the compactness of the sequence 
$(\varphi(u_m))_{m \in \N}$ in $L^2(\O \times (0,T))$.

\smallskip

Therefore, there exists a subsequence, still denoted by $(\varphi(u_m))_{m 
\in \N}$, and there exists $\chi \in L^2(\O \times (0,T))$ such that  
$\varphi(\uappm)$ converges, as $m \tends \infty$, to $\chi$ 
in $L^2(\O\times(0,T))$. Indeed, since $(\varphi(u_m))_{m \in \N}$ is bounded in 
$L^\infty(\O\times(0,T))$, this convergence holds in $L^q(\O\times(0,T))$ 
for all $1 \le q < \infty$.
Furthermore, since $\varphi$ is nondecreasing, Theorem \refp{parnlAF} gives that $\chi=\varphi(u)$.

\medskip

Up to now, the following properties have been shown to be satisfied by
a convenient subsequence:

\begin{enumerate}

\item[(i)] $(u_{m})_{m\in\N}$ converges to $u$, as $m\tends\infty$, in 
$L^\infty(\O\times(0,T))$ for the weak-$\star$ topology,

\item[(ii)] $(\varphi(u_{m}))_{m\in\N}$ converges to $\varphi(u)$ in  $L^1(\O\times(0,T))$ (and even in $L^p(\O\times(0,T))$ for all $p \in [1,\infty)$).
\end{enumerate}

There remains to show that $u$ is a weak solution of
Problem (\ref{paranlanf})-(\ref{paranlani}), which concludes the 
proof of Theorem \ref{parnlcvgce}.

\smallskip

Let $m\in\N$. For the sake of simplicity, we shall use the notations $\mesh = \mesh_{m}$, $h=\size(\mesh)$
and $k=k_{m}$. Let $\psi\in\AT$.
We multiply \refep{parnlschema} by $k \psi(\centercv,n k) $, and sum
the result on $n\in\NNpar$ and $\cv\in\mesh$. We obtain

\be
T_{1m} + T_{2m} = T_{3m},
\label{parnl1g}\ee

with

\[
T_{1m} = \sumparTsk \summesh \mcv (\ucvnp-\ucvn) \psi(\centercv,n k),
\]

\[
T_{2m} = - \sumparTsk k \summesh \sumv \transcvcvv
\Bigl(\varphi(\ucvvn) - \varphi(\ucvn)\Bigr) \psi(\centercv,n k),
\]

and

\[
T_{3m} = \sumparTsk k \summesh \psi(\centercv,n k) \mcv\fcvn.
\]

We first consider $T_{1m}$. 

\[
\left.\begin{array}{lll}
T_{1m} = &
\dsp{\sum_{n=1}^{\Nk}} \summesh \mcv \ucvn
\Bigl(\psi(\centercv,(n-1) k) - \psi(\centercv,n k) \Bigr) &+ \\
& \summesh \mcv \Bigl(u_\cv^{\Nk+1} \psi(\centercv,k \Nk) - u_\cv^0
\psi(\centercv,0)\Bigr).
\end{array}\right.
\]

Performing one more step of the induction in Lemma 
\ref{parnlestinf}, it is clear that
 $\vert u_\cv^{\Nk+1} \vert < U + 2T \normefd$, for all $\cv \in \mesh$.

Since $0 < T- \Nk k \le k$, there exists $C_{1,\psi}$ which only depends 
on $\psi$, $T$ and $\O$, such that  $\vert \psi(\centercv, \Nk k) \vert \le 
kC_{1,\psi}$. Hence,

\[
\summesh \mcv u_\cv^{\Nk+1} \psi(\centercv,k \Nk) \tends 0  \text{ as }m\tends \infty.
\]

Since  

\[ \Vert \summesh \ucvz 1_\cv - u_0 \Vert_{L^1(\O)} \tends  0,  \text{ as }m\tends \infty,
\] 

(where $1_\cv (x) = 1$ if $x \in \cv$, $0$ otherwise), one has 

\[
\summesh \mcv  u_\cv^0 \psi(\centercv,0) \tends \int_\O u_0(x) \psi(x,0) 
dx \text{ as } m\tends \infty.
\]

Since  $(u_{m})_{m\in\N}$ converges, as $m \tends + \infty$, to $u$ in 
$L^\infty(\O\times(0,T))$, for the weak-$\star$ topology, and 
since $\vert u_\cv^{\Nk} \vert < U + T \normef$, for all $\cv \in 
\mesh$, the following property also holds:

\[
\sum_{n=1}^{\Nk} \summesh \mcv \ucvn
\Bigl(\psi(\centercv,(n-1) k) - \psi(\centercv,n k) \Bigr) \tends 
-\int_0^T\int_\O u(x,t) \psi_t(x,t) dx dt \text{ as }m\tends \infty.
\]

Therefore,

\[
T_{1m}  \tends -\int_0^T\int_\O u(x,t) \psi_t(x,t) dx dt
- \int_\O u_0(x) \psi(x,0) dx, \text{ as } m \rightarrow \infty.
\]

\separ

We now study $T_{2m}$. This term can be rewritten as

\[
\left.\begin{array}{ll}
T_{2m} = & - \sumparTsk k \sum_{\edgecvcvv \in \edgesint} m(\edgecvcvv)
(\varphi(\ucvvn) - \varphi(\ucvn)) \frac {\psi(\centercv,n k) - \psi(\centercvv,n k)}
{\dcvcvv}.
\end{array}\right.
\]

It is useful to introduce the following expression:

\[
\left.\begin{array}{ll}
T'_{2m} & =
\sumparTsk \dsp{\int_{nk}^{(n+1)k} \int_\O} \varphi(\uapp(x,t)) \Delta \psi(x,n k) dx dt\\
& = \sumparTsk k \summesh \varphi(\ucvn) \dsp{\int_\cv} \Delta \psi(x,n k) dx \\
& = \sumparTsk k \sum_{\edgecvcvv \in \edgesint} (\varphi(\ucvn) - \varphi(\ucvvn))
\dsp{\int_{\edgecvcvv}} \grad\psi(x ,n k)\cdot\bfn_{\cv, \cvv} d\g(x).
\end{array}\right.
\]

The sequence  $(\varphi(u_{m}))_{m\in\N}$ converges to
$\varphi(u)$  in $L^1(\O\times(0,T))$; furthermore, it is bounded in 
$L^\infty$ so that the 
integral between $T$ and $(\Nk+1)k$ tends to 0. Therefore:

\[
T'_{2m}  \tends 
 \int_0^T \int_\O  \varphi(u(x,t)) \Delta \psi(x,t) dx dt, \text{ 
as }m\rightarrow\infty.
\]

The term $T_{2m}+T'_{2m}$ can be written as

\[
\left.\begin{array}{ll}
T_{2m}+T'_{2m} & =
\sumparTsk k \sum_{\edgecvcvv \in \edges} m(\edgecvcvv) (\varphi(\ucvn) - \varphi(\ucvvn))
R_{\cv ,\cvv}^n,
\end{array}\right.
\]

with

\[
\left.\begin{array}{ll}
R_{\cv ,\cvv}^n = \dsp{\frac 1 {m(\edgecvcvv)} \int_{\edgecvcvv}}
\grad\psi(x ,n k)\cdot\bfn_{\cv, \cvv} d\g(x)
- \frac {\psi(\centercvv,n k) - \psi(\centercv,n k)} {\dcvcvv}.
\end{array}\right.
\]

Thanks to the regularity properties of $\psi$ there exists $C_{\psi}$, 
which only depends on $\psi$, such that
$\vert R_{\cv ,\cvv}^n \vert \le C_{\psi} h$. Then,
using the estimate (\ref{parnlestimx}), we conclude that
$ T_{2m}+T'_{2m} \tends 0$ as $m\rightarrow\infty$. Therefore,

\[
T_{2m} \tends -\int_0^T \int_\O  \varphi(u(x,t)) \Delta \psi(x,t) dx dt, \text{ 
as }m\rightarrow\infty.
\]

Let us now study $T_{3m}$.

Define $f_{\mesh,k} \in X(\mesh,k)$ by $f_{\mesh,k}(x,t)=f^n_\cv$ if 
$(x,t)\in \cv \times (nk, nk+k)$. Since $f_{\mesh,k} \tends f$ in $L^1(\O 
\times (0,T)$ and since $f \in L^\infty(\O \times (0,2T)$,

\[
T_{3m} \tends \int_\O\int_0^T f(x,t) \psi(x,t) dt dx, \text{ 
as }m\rightarrow\infty. 
\]

Passing to the limit in Equation \eqref{parnl1g} gives that $u$ is a weak solution of
Problem (\ref{paranlanf})-(\ref{paranlani}). This concludes the 
proof of Theorem \ref{parnlcvgce}.\eof

\begin{remark} \normalfont
This convergence proof is quite similar in the case of
the implicit scheme, with the additional condition that
$(k_m)_{m\in\N}$ converges to zero, since condition \eqref{parnlcfl}
does not have to be satisfied.
\end{remark}

\begin{remark} \normalfont\label{subsequence}

The above convergence result was shown for a subsequence only.
A convergence theorem is obtained for the full set of approximate solutions, if a uniqueness
result is valid. Such a result can be easily  obtained in the case of a smooth boundary
and is given in section \ref{up} below. 
For this case, an extension to the definition \refp{meshneuman} of
admissible meshes is given hereafter.
\end{remark}

\begin{defi}[Admissible meshes for regular domains]\label{parnlHT}  \normalfont
\index{Mesh!Admissible --!--for regular domains}
Let $\O$ be an open bounded connected subset of $\R^d$, $d = 2$ or $3$ 
with a $C^2$ boundary $\dr\O$. An
admissible finite volume mesh of $\O$ is given by an open bounded polygonal set $\O'$ containing $\O$,
and an admissible mesh  $\Mesh'$ of $\O'$ in the sense 
of Definition \refp{meshneuman}. The set of control volumes of the mesh of $\O$ 
are $\{ \cv'\cap\O$, $\cv' \in \mesh'$ such that
$m_d(\cv'\cap\O) > 0\}$ and the set of edges of the mesh is $\edges = 
\{\edge \cap \Obar$, $\edge \in \edges'$ such that $m_{d-1}(\edge \cap 
\Obar) >0\}$, where $\edges'$ denotes the set of edges of $\mesh'$ and 
$m_N$ denotes the $N$-dimensional Lebesgue measure.
\end{defi}

\begin{remark} \normalfont
For smooth domains $\O$, the set of edges $\edges$ of an admissible mesh
of $\O$ does not contain the parts of the boundaries of the control volumes 
which are  included in the boundary $\dr\O$ of $\O$.
\end{remark}

\subsection{Weak convergence and nonlinearities}\label{nlparaws}

We show here a property which was used in the proof of Theorem \ref{parnlcvgce}.

\begin{theo}\label{parnlAF}\index{Convergence!-- for the weak star topology}
Let $U>0$ and $\varphi \in C([-U,U])$ be a nondecreasing function. Let $\o$ 
be an open bounded subset of $\R^N$, $N \ge 1$. Let $(u_n)_{\nnn} \subset 
L^\infty(\o)$ such that

(i) $-U \le u_n \le U$ a.e. in $\o$, for all $n\in \N$;

(ii) there exists  $u \in L^\infty(\o)$ such that $(u_n)_{n\in\N}$
converges to $u$ in $L^\infty(\o)$ for the weak-$\star$ topology;

(iii) there exists a function $\chi\in L^1(\o)$ such that
$(\varphi(u_n))_{n\in\N}$ converges to $\chi$ in $L^1(\o)$.

Then $\chi(x) = \varphi(u(x))$, for a.e. $x\in \o$.
\end{theo}

\proof of Theorem \ref{parnlAF} 

\smallskip

First we extend the definition
of $\varphi$ by $\varphi(v) = \varphi(-U) + v + U$ for all
$v<-U$ and $\varphi(v) = \varphi(U) + v - U$ for all $v>U$, and denote
again by $\varphi$ this extension of $\varphi$ which now maps $\R$
into $\R$, is continuous and nondecreasing.
Let us define $\a_\pm$ from $\R$ to $ \R$ by
$\a_-(t) = \inf \{v\in \R, \varphi(v) = t \}$ and $\a_+(t) = \sup \{v\in \R, \varphi(v) = t \}$,
for all $t\in \R$.

Note that the functions $\a_\pm$ are  increasing and that

(i) $\a_-$ is left continuous and therefore
lower semi-continuous, that is

\[
t=\lim_{n \tends \infty} t_n \Longrightarrow \a_-(t)\le \liminf_{n \rightarrow 
\infty} \a_-(t_n),
\]

(ii) $\a_+$ is right continuous  and therefore
upper semi-continuous, that is
\[
t=\lim_{n \tends \infty} t_n  \Longrightarrow \a_+(t)\ge \limsup_{n \rightarrow 
\infty} \a_+(t_n).
\]

Thus, since we may assume, up to a subsequence, that $\varphi(u_n) \tends 
\chi$ a.e. in $\o$,

\be
\a_-(\chi(x)) \le \liminf_{n\rightarrow \infty} \a_-\Bigl(\varphi(u_n(x))\Bigr) 
\le \limsup_{n\rightarrow \infty} \a_+\Bigl(\varphi(u_n(x))\Bigr)
\le \a_+(\chi(x)),
\label{parnl1}
\ee
for a.e. $x \in \o$.

~

A direct application of the definition of the functions $\a_-$ and $\a_+$ 
gives

\be
\a_-\Bigl(\varphi(u_n(x))\Bigr) \le u_n(x) \le 
\a_+\Bigl(\varphi(u_n(x))\Bigr).
\label{parnl5}
\ee

Let $L^1_+=\{\psi \in L^1(\o)$, $\psi \ge 0$ a.e.$\}$. Let $\psi \in 
L^1_+$. We multiply \eqref{parnl5} by $\psi(x)$ and integrate over $\o$, 
it yields

\be
\int_\o \a_-\Bigl(\varphi(u_n(x))\Bigr) \psi(x)dx \le \int_\o 
u_n(x)\psi(x)dx \le \int_\o\a_+\Bigl(\varphi(u_n(x))\Bigr)\psi(x)dx.
\label{plus}
\ee

Applying Fatou's lemma to the sequences of $L^1$ positive
functions $\a_-(\varphi(u_n)) \psi - \a_-(\varphi(-U)) \psi$ and 
$\a_+(\varphi(U)) \psi-\a_+(\varphi(u_n)) \psi$
yields, with \eqref{parnl1},

\[
\int_\o \a_-(\chi(x)) \psi(x) dx \le
\liminf_{n\rightarrow\infty}  \int_\o \a_-\Bigl(\varphi(u_n(x))\Bigr) 
\psi(x) dx,
\]

and

\[
\limsup_{n\rightarrow\infty}  \int_\o \a_+\Bigl(\varphi(u_n(x))\Bigr)
\psi(x) dx \le \int_\o \a_+(\chi(x)) \psi(x) dx.
\]

Then, passing to the  $\liminf$ and $\limsup$ in \eqref{plus} and using the  
convergence of $(u_n)_{n\in\N}$ to $u$ in $L^\infty(\o)$ for the 
weak-$\star$ topology gives

\[
\int_\o \a_-(\chi(x)) \psi(x) dx \le
\int_\o u(x) \psi(x) dx \le
\int_\o \a_+(\chi(x)) \psi(x) dx.
\]

Thus, since $\psi$ is arbitrary in $L^1_+$, the following inequality 
holds for  a.e. $x\in \o$:
\[\a_-(\chi(x)) \le u(x) \le \a_+(\chi(x)),\]
which implies in turn 
that $\chi(x) = \varphi(u(x))$ for a.e. $x\in \o$. This
completes the proof of Theorem \ref{parnlAF}.
\eof

\begin{remark} \normalfont
Another proof of Theorem \ref{parnlAF} is possible by passing to the limit 
in the inequality 
\[
0 \le \int_\omega (\varphi(u_n)(x) - \varphi(v(x)))(u_n(x) - v(x)) dx, 
\forall v \in L^\infty(\omega),
\]
which leads to
\[
0 \le \int_\omega (\chi(x) - \varphi(v(x)))(u(x) - v(x)) dx, 
\forall v \in L^\infty(\omega).
\]
From this inequality, one deduces that $\chi = \varphi(u)$ a.e. on 
$\omega$.

\medskip

A third proof is  possible by using  the concept of
nonlinear weak-$\star$ convergence, see Definition \refp{nlwsdef}.
\end{remark}

\subsection{A uniqueness result for nonlinear diffusion equations} \label{up}

The uniqueness of the weak solution to  variations of Problem 
(\ref{paranlanf})-(\ref{paranlani}) has been proved by several
authors. For precise references we refer to \citet{MEI}.
Also rather similar proofs have been given in \citet{BKP} and \citet{GHP}.
Recall that this uniqueness result allows to obtain a convergence result 
on the whole set of finite volume approximate solutions to Problem \eqref{parabch}-\eqref{parabc2ch}
(see Remark \ref{subsequence}).

\smallskip

The uniqueness of the weak solution to
Problem (\ref{paranlanf})-(\ref{paranlani})
immediately results from the following property.
\index{Uniqueness!-- of the solution to a nonlinear diffusion equation}
\begin{theo}\label{parnlthunicite}
Let $\O$ be an open bounded subset of $\R^\dim$ with a $C^2$ boundary, 
and suppose that items (ii), (iii) and (iv) of  Assumption \ref{parnlH} are satisfied. 
Let $u_1$ and $u_2$ be two solutions of
Problem (\ref{paranlanf})-(\ref{paranlani}) in the sense of Definition 
\refp {parnldefess},
with initial conditions $u_{0,1}$ and $u_{0,2}$ and
source terms $v_1$ and $v_2$ respectively, that is, for $u_1$ (resp. 
$u_2$), $u_0=u_{0,1}$ (resp. 
$u_0=u_{0,2}$) in \eqref{paranlani} and $f=v_1$ 
(resp. $v_2$) in \eqref{paranlanf}.

Then for all $T>0$,

\[
\int_0^T \int_\O\vert u_1(x,t) -u_2(x,t) \vert dx dt
\leq T \int_\O \vert u_{0,1}(x) -u_{0,2}(x) \vert dx +
\int_0^T \int_\O (T-t) \ \vert v_1(x,t) - v_2(x,t) \vert dx dt.
\]
\end{theo}

Before proving Theorem \ref{parnlthunicite},
let us first show the following auxiliary result.

\subsubsection{The existence of regular solutions to the adjoint problem}
\begin{lemma}\label{parnlauxi}
Let $\O$ be an open bounded subset of $\R^\dim$ with a $C^2$ boundary, 
and suppose that 
$\varphi$ is a nondecreasing locally Lipschitz-continuous function. 
Let $T>0$,
$w\in C^\infty_ c(\O \times (0,T))$ such that $\vert w \vert \le 1$,
and  $g\in C^\infty(\Obar \times [0,T])$
such that there exists $r\in\R$ with $0 < r \le g(x,t)$, for all
$(x,t)\in \O\times (0,T)$.

\separ

Then there exists a unique function $\psi \in C^{2,1}(\overline{\O}\times [0,T])$
such that
\be
\psi_t(x,t) + g(x,t) \Delta \psi(x,t) = w(x,t), \ \
\ \ \text{ for all } (x,t)\in \O\times (0,T),
\label{parnllemme1}\ee
\be
\grad \psi \cdot \n (x,t) = 0, \ \ \text{ for all } (x,t)\in \dr\O\times (0,T),
\label{parnllemme3}\ee
\be
\psi(x,T) = 0, \ \ \text{ for all } x \in \O.
\label{parnllemme2}\ee
Moreover the function $\psi$ satisfies
\be
\vert \psi(x,t) \vert \le T-t,
\ \ \text{ for all } (x,t)\in \O\times (0,T),
\label{parnllemme4}\ee

and

\be
 \int_0^T \int_\O g(x,t) \Bigl(\Delta\psi(x,t)\Bigr)^2 dx dt
\le 4 T \int_0^T \int_\O  \vert \grad w(x,t) \vert^2 dx dt.
\label{parnllemme5}\ee
\end {lemma}

\proof of Lemma \ref{parnlauxi}

\smallskip

It will be useful in the following to
point out that the right hand side of (\ref{parnllemme5}) does not depend on $g$.
Since the function $g$ is bounded away from zero,
equations (\ref{parnllemme1})-(\ref{parnllemme2}) define
a boundary value problem for a usual heat equation with an initial 
condition, in which the time variable is reversed.  Since
$\O$, $g$ and $w$ are sufficiently smooth,
this problem has a unique solution
$\psi\in \AT$, see \citet{LSU}.
Since $\vert w\vert \le 1$, the functions $T-t$ and $-(T-t)$ are
respectively upper and lower solutions of Problem 
(\ref{parnllemme1})-(\ref{parnllemme3}). Hence we get (\ref{parnllemme4}) 
(see \citet{LSU}).

\separ

In order to show (\ref{parnllemme5}), multiply (\ref{parnllemme1}) by
$\Delta\psi(x,t)$, integrate by parts on $\O\times(0,\tau)$, for $\tau 
\in (0,T]$. This gives

\be\ban
\dsp{\half\int_\O \vert \grad \psi(x,0) \vert^2 dx \ - \
\half\int_\O \vert \grad \psi(x,\tau) \vert ^2 dx \ + \
\int_0^\tau \int_\O  g(x,t) \Bigl(\Delta\psi(x,t)\Bigr)^2 dx dt } = \\
\dsp{- \int_0^\tau \int_\O \grad w(x,t) \cdot \grad \psi(x,t) dx dt}.
\ean\label{parnllemme6}\ee

Since  $\grad \psi(\cdot,T) = 0$, letting $\tau = T$ in
(\ref{parnllemme6}) leads to

\be\ban
\dsp{\half\int_\O \vert \grad \psi(x,0) \vert^2 dx +
\int_0^T \int_\O g(x,t) \Bigl(\Delta\psi(x,t)\Bigr)^2 dx dt } = \\
\dsp{- \int_0^T \int_\O \grad w(x,t) \cdot \grad \psi(x,t) dx dt}.
\ean\label{parnllemme7}\ee

Integrating (\ref{parnllemme6}) with respect to $\tau\in (0,T)$ leads to

\be\ban
\dsp{\half\int_0^T \int_\O \vert \grad \psi(x,\tau) \vert^2 dx d\tau } \le &
\dsp{ \frac {T} {2} \int_\O \vert \grad \psi(x,0) \vert^2 dx} &+ \\&
\dsp{T \int_0^T \int_\O g(x,t) \Bigl(\Delta\psi(x,t)\Bigr)^2 dx dt} &+\\&
\dsp{T \int_0^T \int_\O \vert\grad w(x,t) \cdot \grad \psi(x,t)\vert dx dt}.
\ean\label{parnllemme8}\ee

Using (\ref{parnllemme7}) and (\ref{parnllemme8}), we get

\be
\half \int_0^T \int_\O \vert \grad \psi(x,\tau) \vert^2 dx d\tau \le
2 T \int_0^T \int_\O \vert\grad w(x,t) \cdot \grad \psi(x,t)\vert dx dt.
\label{parnllemme9}\ee

Thanks to   the Cauchy-Schwarz  inequality, the right hand side of \eqref{parnllemme9}
may be estimated as follows:

\[\ban
\dsp{\Bigl[ \int_0^T \int_\O \vert\grad w(x,t) \cdot \grad \psi(x,t)\vert dx dt
\Bigr]^2 } \le &
\dsp{ \int_0^T \int_\O \vert \grad \psi(x,t) \vert^2 dx dt}\\ &
\times \dsp{\int_0^T \int_\O \vert \grad w(x,t) \vert^2 dx dt}.
\ean\]

With (\ref{parnllemme9}), this implies

\[\ban
\dsp{\Bigl[ \int_0^T \int_\O \vert\grad w(x,t) \cdot \grad \psi(x,t)\vert dx dt
\Bigr]^2} \le &
\dsp{4 T \int_0^T \int_\O \vert\grad w(x,t) \cdot \grad \psi(x,t)\vert dx 
dt} \\ &
\times \dsp{\int_0^T \int_\O \vert \grad w(x,t) \vert^2 dx dt}.
\ean\]

Therefore,

\[
\int_0^T \int_\O \vert\grad w(x,t) \cdot \grad \psi(x,t)\vert dx dt
\le
4 T \int_0^T \int_\O \vert \grad w(x,t) \vert^2 dx dt,
\]

which, together with  (\ref{parnllemme7}), yields
 (\ref{parnllemme5}).\eof

\subsubsection{Proof of the uniqueness theorem}

Let $u_1$ and $u_2$ be two solutions of Problem
(\ref{paranlanw}), with initial conditions $u_{0,1}$  and $u_{0,2}$ and
source terms $v_1$ and $v_2$ respectively. We set $u_d = u_1 - u_2$,
$v_d = v_1-v_2$ and $u_{0,d} = u_{0,1}-u_{0,2}$.
Let us also define, for all $(x,t)\in \O\times \R^\star_+$,
$q(x,t) = \dsp{\frac{\varphi(u_1(x,t)) - \varphi(u_2(x,t))} {u_1(x,t)-u_2(x,t)}}$ if
$u_1(x,t) \neq u_2(x,t)$, else $q(x,t)=0$.
For all $T \in \R^\star_+$ and for all
$\psi \in\ \AT$, we deduce from (\ref{paranlanw}) that

\be
\left.\begin{array}{ll}
\dsp{\int_0^T\int_\O} \Bigl[ u_d(x,t) \Bigl( \psi_t(x,t)
+ q(x,t) \Delta \psi(x,t)\Bigr) + v_d(x,t) \psi(x,t) \Bigr] dx dt
\ & + \\
\dsp{\int_\O} u_{0,d}(x) \psi(x,0) dx & = 0.
\end{array}\right.
\label{parnlun1}\ee

Let $w\in C^\infty_ c(\O\times (0,T))$, such that
$\vert w \vert \le 1$.
Since $\varphi$ is locally Lipschitz continuous, we can define its
Lipschitz constant, say $B_M$, on $[-M,M]$, where
$M=\max\{\|u_1\|_{L^\infty (\O\times(0,T))},
\|u_2\|_{L^\infty (\O\times(0,T))}\}$ so that $0 \le q \le B_M$ a.e. on 
$\O\times(0,T)$.

\separ
Using mollifiers, functions
$q_{1,n} \in C^\infty_c(\O\times(0,T))$ may be constructed 
such that $ \| q_{1,n} - q \|_{L^2 (\O\times(0,T))} \le \frac 1 n$ and
$0 \le q_{1,n} \le B_M$, for $n\in \N^\star$.
Let $q_n = q_{1,n} + \frac 1 n$. Then

\[
\frac 1 n \leq q_n(x,t) \leq B_M + \frac 1 n,
\text{ for all } (x,t)\in \O\times(0,T),
\]

and

\[
\ban
\dsp{\int_0^T \int_\O \frac{(q_n(x,t)-q(x,t))^2} {q_n(x,t)}} dx dt \le
2 \Bigl(
& \dsp{\int_0^T \int_\O \frac{(q_n(x,t)-q_{1,n}(x,t))^2} {q_n(x,t)}}
dx dt \ + \\
& \dsp{\int_0^T \int_\O \frac{(q_{1,n}(x,t)-q(x,t))^2} {q_n(x,t)}}
dx dt \Bigr),
\ean
\]

which shows that
\[
\int_0^T \int_\O \frac{(q_n(x,t)-q(x,t))^2} {q_n(x,t)} dx dt \le 2 n \Bigl( \frac{T \meas(\O)}{n^2} + \frac 1 {n^2} \Bigr).
\]

It leads to
\be
\| \frac{q_n-q}{\sqrt {q_n}} \|_{L^2 (\O\times (0,T))}  \tends  0 \ \
\text{as} \ \ n \tends \infty.
\label{parnlun3}\ee

Let $\psi_n\in \AT$ be given by lemma \ref{parnlauxi}, with $g=q_n$.
Substituting $\psi$ by $\psi_n$ in (\ref{parnlun1}), using (with $g=q_n$ 
and $\psi=\psi_n$) (\ref{parnllemme1}) and (\ref{parnllemme4}) give

\be\ban
\dsp{\vert \int_0^T\int_\O} u_d(x,t) \Bigl( w(x,t)
+ (q(x,t) - q_n(x,t)) \Delta \psi_n(x,t)\Bigr) dx dt \vert & \le \\
\dsp{\int_0^T\int_\O} \vert v_d(x,t)\vert (T-t) dx dt +
T \dsp{\int_\O} \vert u_{0,d}(x) \vert dx.
\ean\label{parnlun10}\ee

The Cauchy-Schwarz   inequality yields

\be\ban
\Bigl[ \dsp{\int_0^T\int_\O} \vert u_d(x,t) \vert
\vert (q(x,t) - q_n(x,t)) \Delta \psi_n(x,t)\vert dx dt \Bigr]^2 \le
\dsp{4 M^2 } \\
\dsp{\int_0^T\int_\O} \Bigl(\frac {q(x,t) - q_n(x,t)}
{\sqrt{q_n(x,t)}} \Bigr)^2 dx dt
\int_0^T \int_\O q_n(x,t) \Bigl(\Delta\psi_n(x,t)\Bigr)^2 dx dt.
\ean\label{parnlun11}\ee

We deduce from (\ref{parnllemme5}) and  (\ref{parnlun3}) that the right hand side
of (\ref{parnlun11}) tends to zero as $n\tends\infty$. Hence the
left hand side of (\ref{parnlun11}) also tends to zero as $n\tends\infty$.
Therefore letting $n\tends\infty$ in (\ref{parnlun10}) gives

\be\ban
\dsp{\vert \int_0^T\int_\O} u_d(x,t) w(x,t) dx dt \vert \le &
\dsp{\int_0^T\int_\O} \vert v_d(x,t)\vert (T-t)  dx dt \ + \\ &
T \dsp{\int_\O} \vert u_{0,d}(x) \vert dx.
\ean\label{parnlun12}\ee

Inequality (\ref{parnlun12}) holds for any function
$w\in C^\infty_ c(\O\times (0,T))$, with
$\vert w\vert \le 1$. Let us take as functions $w$ the elements of
a sequence $(w_m)_{m\in\N}$ such that
$w_m\in C^\infty_ c(\O\times (0,T))$ and $\vert w_m\vert \le 1$  for all $m\in\N$,  and
the sequence $(w_m)_{m\in\N}$
converges to sign$(u_d(\cdot,\cdot))$ in $L^1(\O\times(0,T))$.
Letting $m\tends \infty$ yields

\[\ban
\dsp{\int_0^T\int_\O} \vert u_d(x,t) \vert dx dt \le
\dsp{\int_0^T\int_\O} \vert v_d(x,t)\vert (T-t) dx dt +
T \dsp{\int_\O} \vert u_{0,d}(x) \vert dx,
\ean\]

which concludes the proof of Theorem \ref{parnlthunicite}.
\eof

\chapter{Hyperbolic equations in  the one  dimensional case} \label{hyper1d}

This chapter is devoted to the numerical schemes for one-dimensional hyperbolic conservation laws.
Some basics on the solution to linear or nonlinear hyperbolic equations with initial data and without boundary conditions are first recalled.
We  refer to \citet{godlewski-ellipses}, \citet{godlewski-springer}, \citet{kronerbook}, \citet{leveque} and  \citet{serre} for extensive studies of  theoreticaland/or numerical aspects;  we shall highlight here the finite volume point of view for several
well known schemes,  comparing them with finite difference schemes, either for the linear and the nonlinear case. 
Convergence results for numerical schemes are presented, using a ``weak $BV$ inequality"  which is not really necessary in the 1D case (at least for $BV$ initial data), but is crucial in the multidimensional case. 
We also recall the classical proof of convergence, which uses a ``strong $BV$ estimate" and the Lax-Wendroff  theorem, and does not seem to extend to the mutidimentional case.
Error estimates  which also hold are not adressed in this chapter:  they are  given in the chapter in the multidimensional case (Chapter \ref{hypmd}).

\medskip

Throughout this chapter, we shall  focus on explicit schemes. 
However, all the results which are presented here can be extended  to implicit schemes;  this requires a bit of work and is  detailed  in the multidimensional case (see \refep{estschemai} for the scheme).

\section{The continuous problem}\label{hyppb}

\hpage{Hyperbolic 1d}

Consider the nonlinear hyperbolic equation with initial data:
\index{Equation!Hyperbolic --} 
\be
\ba  \dsp
  u_t(x,t) +  (f(u))_x (x,t)
   & = 0 & x \in \R, \ \ t \in \R_+, 
\\
 u(x,0)  & = u_0 (x), &  x \in \R, 
\ea \label{hyperbolic1D}
\ee

where $f$ is a given function from $\R$ to $\R$, of class $C^1$, $u_0 \in L^\infty
(\R)$ and where  the partial derivatives of 
$u$ with respect to time and space are denoted by $u_t$ and $u_x$. 

\begin{example}[B\"urgers equation] \normalfont
A simple flow model was introduced by B\"urgers and yields the following equation:
\be
u_t(x,t) + u(x,t) u_x(x,t) - \vareps u_{xx}(x,t) = 0
 \label{Burgeps}
\ee
B\"urgers studied the limit case which is obtained when $\eps$ tends to 0; the resulting  equation is
(\ref{hyperbolic1D})
with $\dsp f(s) = \frac {s^2} 2$, i.e. 

\[
u_t(x,t) +\frac 1 2  (u^2)_x (x,t) = 0
\]

\end{example}

\begin{defi}[Classical solution] \label{defclass1} \normalfont Let $f \in C^1(\R,\R)$ and $u_0 \in C^1(\R,\R)$; a classical solution to Problem (\ref{hyperbolic1D}) is a function $ u \in C^1(\R \times 
\R_+,\R)$ such that  

\[
\left\{\begin{array}{rll} 
 u_t  (x,t) + f'(u(x,t))  u_x  (x,t)  & = 0, & \forall \  x \in \R, \, \forall \ t \in \R_+, 
\\
 u(x,0)  & = u_0 (x), &  \forall \,  x \in \R.
\end{array}\right. 
\]
\end{defi}
Recall that in the linear case, i.e.
$f(s) = c s$ for all $s \in \R$, for some  $c \in \R$, there exists (for 
$u_0 \in C^1(\R,\R)$) a unique classical
solution. It is  $u(x,t) = u_0(x-ct)$, for all $x \in \R$ and for all $t \in \R_+$. In the
nonlinear case, the existence of such a solution depends on the initial data
$u_0$; in fact, the following result holds: 

\begin{prop} \label{nocr} Let $f \in C^1(\R,\R)$ be a nonlinear function, i.e. such that there
exist $s_1$, $s_2 \in \R$ with $f'(s_1) \ne f'(s_2)$;  then there exists $u_0 \in
C^\infty_c(\R,\R) $ such that Problem (\ref{hyperbolic1D}) has no 
classical solution.
\end{prop}

Proposition \ref{nocr} is an easy consequence of the following remark.

\begin{remark} \normalfont
If $u$ is a classical solution to \eqref{hyperbolic1D}, then $u$ is constant along the
characteristic lines which are defined by
\[ x(t) = f'(u_0(x_0)) t + x_0, \; t \in \R_+,\]
where $x_0 \in \R$ is the origin of the characteristic. This is the equation of a straight
line issued from the point  $(x_0,0)$ (in the $(x,t)$ coordinates). Note that if
$f$ depends on $x$ and $u$ (rather than only on $u$), the characteristics are no
longer straight lines.
\end{remark}

The concept of weak solution is introduced in order to define solutions of 
\eqref{hyperbolic1D} when classical solutions do not exist.

\begin{defi}[Weak
solution] \label{defweak1}  \normalfont Let $f \in C^1(\R,\R)$ and $u_0 \in L^\infty(\R)$; a weak
solution to Problem (\ref{hyperbolic1D})  is a function $ u $ such
that 

\be \ba  u \in L^\infty (\R \times \Rpast), \\
\dsp{ \dint u(x,t) \varphi_t(x,t) dt dx + \dint f(u(x,t))
\varphi_x(x,t) dt dx + \int_\R u_0(x) \varphi(x,0) dx = 0, } \\
\hfill \forall \varphi \in C^1_c(\R \times \R_+,\R).
\ea
\label{h1dweak}
\ee\end{defi}

\begin{remark} \normalfont ~

1. If $u \in C^1(\R \times \R_+,\R)\cap  L^\infty(\R \times \Rpast)$ then $u$ is a weak solution if and only if $u$ is a
classical solution. 

2. Note that in the above definition,  we require  the test function 
$\varphi$ to belong to $C^1_c(\R\times \R_+,\R)$, 
so that $\varphi$ may be non zero at time $t=0$. 
\end{remark}

One may show that there exists at least one weak solution to 
\eqref{hyperbolic1D}. In the linear case, i.e.
$f(s) = cs $, for all $s \in \R$, for some $c \in \R$, this solution is unique (it is
$u(x,t)=u_0(x-ct)$ for a.e. $(x,t) \in \R \times \Rp$). However, the uniqueness of
this weak solution in the general  nonlinear case is no longer true. Hence  the concept of entropy 
weak solution, for which an existence and uniqueness result is known.

\begin{defi}[Entropy weak solution] \label{defent1} \index{Entropy!-- weak 
solution}\index{Entropy!Continuous --  inequalities} \normalfont
Let $f \in C^1(\R,\R)$ and $u_0 \in L^\infty(\R)$; 
  the entropy weak solution to Problem (\ref{hyperbolic1D}) is a function  $u$ such that 
 
\be\left\{\ban
u \in L^\infty(\R\times \Rpast), \\
\dsp{ \dint \eta(u(x,t)) \varphi_t(x,t) dt dx + 
\dint \Phi(u(x,t)) \varphi_x(x,t) dt dx + \int_\R
\eta(u_0(x)) \varphi(x,0) dx \ge 0},
\\
\hfill \forall \varphi \in C^1_c(\R \times \R_+,\R_+),
\\
\hfill \text{for all convex  function } \eta \in C^1(\R,\R) \text{~and } \Phi
\in C^1(\R,\R) \text{~such that }  \Phi' = \eta' f'. 
\ean\right.\label{h1dent}
\ee
\end{defi}

\begin{remark} \normalfont The solutions of \eqref{h1dent} are necessarily 
solutions of \eqref{h1dweak}. This can be shown by taking in \eqref{h1dent} $\eta 
(s) = s$  for 
all $s\in\R$, $\eta (s) = -s$, 
for all $s\in\R$,  and regularizations of the positive and negative parts 
of the test functions of the weak formulation.
\end{remark}

\begin{theo}\label{hyp1dexistu} Let $f \in C^1(\R,\R)$, $u_0 \in L^\infty(\R)$, then 
there exists a unique
entropy weak solution \index{Entropy!-- weak 
solution} to Problem (\ref{hyperbolic1D}).
\end{theo}

The proof of this result was first given by Vol'pert in \citet{Vo},
introducing the space $BV(\R)$\index{BV!-- Space} which is defined hereafter and assuming
$u_0\in BV(\R)$, see also \citet{Ol} for the convex case.
In \citet{krushkov}, Krushkov  proved the theorem of existence
and uniqueness in the general case $u_0 \in L^\infty(\R)$, 
using a regularization of $u_0$ in $BV(\R)$, 
under the slightly stronger assumption $f\in C^3(\R,\R)$. 
Krushkov also proved that the solution
is in the space $C(\Rp,L^1_{loc}(\R))$. Krushkov's proof uses
particular entropies, namely the functions $\vert \cdot - \sk \vert$ 
for all $\sk \in\R$, which are generally referred to as ``Krushkov's 
entropies". The ``entropy flux" associated to $\vert \cdot - 
\sk \vert$ may be taken as $f(\cdot\top\sk) -
f(\cdot\bot\sk)$, where $a\top b$ denotes the maximum
of $a$ and $b$ and $a\bot b$ denotes the minimum of $a$ and $b$, 
for all real values $a,b$ (recall that $f(a\top b) -f(a \bot b) = \sgn(a-b)(f(a)- f(b))$).

\begin{defi}[$BV(\R)$]\index{BV!-- Space|(} \label{bv} \normalfont A function $v \in L^1_{loc}(\R)$ is of bounded variation,
that is $v \in  BV(\R)$, if

\be
\vert v \vert_{BV(\R)} = \sup \{\int_\R v(x) \varphi_x(x) dx, \varphi \in 
C^1_c(\R,\R), \vert \varphi (x) \vert \le 1\ \forall x \in \R\} < + 
\infty.
\label{defBV}
\ee
\end{defi}

\begin{remark} \normalfont\label{divbv1d} ~
\benum

\item If $v: \R \tends \R$ is piecewise constant, that is if
there exists an increasing sequence $(x_i)_{i\in\Z}$ with $\R = \cup_{i \in \Z} [x_i, 
x_{i+1}]$ and a sequence $(v_i)_{i\in\Z}$  such that $v\vert_{(x_i, x_{i+1})} =
v_i$, then $\vert  v \vert_{BV(\R)} = \sum_{i \in \Z}
\vert v_{i+1} - v_i \vert$. 

\item  If $v \in C^1(\R,\R)$ then $\vert v \vert_{BV(\R)} = \Vert v_x 
\Vert_{L^1(\R)}$.  

\item The space  $BV(\R)$ is included in the space $L ^\infty(\R)$; 
furthermore,  if $u \in BV(\R)  \cap L^1(\R)$ 
then $\norm{u}{L^\infty(\R)} \le \vert  u \vert_{BV(\R)}$.

\item Let $u \in BV(\R)$ and let $(x_{i+1/2})_{i \in \Z}$ be an 
increasing sequence of real values such that $\R = \cup_{i\in\Z} 
[x_{i-1/2}, x_{i+1/2}]$. For $i \in \Z$, let $\cv_i = (x_{i-1/2}, 
x_{i+1/2})$ and $u_i$ be the mean value of $u$ over $\cv_i$. Then, 
choosing conveniently $\vphi$ in the definition of $\snorm{u}{BV(\R)}$, it 
is easy to show that

\be
\sum_{i \in \Z}\vert u_{i+1}-u_i \vert \le \snorm{u}{BV(\R)}.
\label{bvz}
\ee

Inequality \eqref{bvz} is used for the classical proof of ``$BV$ estimates" 
for the approximate solutions given by finite volume schemes (see Lemma 
\refp{nlstabv} and Corollary \refp{nlstabvco}).

Note that \eqref{bvz} is also true when $u_i$ is the mean value of $u$ 
over a subinterval of $\cv_i$ instead of the mean value of $u$ 
over $\cv_i$.
\eenum
\index{BV!-- Space|)}\end{remark}

Krushkov used a
characterization of entropy weak solutions which is given in the following proposition.

\begin{prop}[Entropy weak solution 
using ``Krushkov's entropies"]\label{prentkr1}~
\index{Entropy!-- weak solution}\index{Krushkov's entropies}
Let $f \in C^1(\R,\R)$ and $u_0 \in L^\infty(\R)$, 
$u$ is the unique entropy weak solution to Problem (\ref{hyperbolic1D}) 
if and only if $u $ is such that 
 
\be\ba 
u \in L^\infty(\R\times \Rpast),
\\
\dsp \dint \vert u(x,t) - \sk \vert \varphi_t(x,t) dt dx + 
\\
\dsp \dint \Bigl( f(u(x,t)\top\sk) - f(u(x,t)\bot\sk) \Bigr) 
\varphi_x(x,t) dt dx + \int_\R
\vert u_0(x) - \sk \vert \varphi(x,0) dx \ge 0,
\\
\hfill \forall \varphi \in C^1_c(\R \times \R_+,\R_+), \forall \sk \in \R. 
\ea\label{h1dentkr}
\ee
\end{prop}

The result of existence of an entropy weak solution  defined by 
\eqref{h1dent} was already  
proved by passing to the limit  on the solutions of an appropriate 
numerical scheme, see e.g. \citet{Ol}, and may  also be  obtained by 
passing to the limit
on finite volume approximations of the solution (see Theorem \refp{nlcv1d} 
in the one-dimensional case and Theorem \refp{eghbcewsol} in the 
multidimensional case).

\begin{remark} \normalfont An entropy weak solution \index{Entropy!-- weak 
solution} 
is sometimes  defined as a function $u$ satisfying:

\be\left\{\ban
\dsp{ \dint u(x,t) \varphi_t(x,t) dt dx + \dint f(u(x,t))
\varphi_x(x,t) dt dx + \int_\R u_0(x) \varphi(x,0) dx = 0, } \\
\hfill \forall \varphi \in C^1_c(\R \times \R_+,\R).
\\
\dsp{ \dint \eta(u(x,t)) \varphi_t(x,t) dt dx + 
\dint \Phi(u(x,t)) \varphi_x(x,t) dt dx \ge 0,} 
\\
\hfill \forall \varphi \in C^1_c(\R \times \R_+^\star,\R_+),
\\
\hfill \text{for all convex  function } \eta \in C^1(\R,\R) \text{~and } \Phi
\in C^1(\R,\R) \text{~such that }  \Phi' = \eta' f'. 
\ean\right.\label{h1dentf}
\ee

The uniqueness of an entropy 
weak solution thus defined depends on the functional space to which $u$ is 
chosen to belong. Indeed, the  uniqueness 
result given in Theorem \ref{hyp1dexistu} is no longer true with 
$u$ defined by \eqref{h1dentf} such that

\be
u, \, f(u) \in L^1_{loc}(\R\times \R_+), \, u \in L^\infty(\R \times (\eps,\infty)), \,
\forall \eps \in \R_+^\star. \label{h1space}
\ee

Under Assumption \eqref{h1space}, every term in \eqref{h1dentf} makes 
sense. Note that \eqref{h1space}-\eqref{h1dentf} is weaker than \eqref{h1dent}. 
An easy counterexample to a uniqueness result of the solution to 
\eqref{h1dentf}-\eqref{h1space} is obtained
with $f(s)=s^2$ for all $s \in \R$ and $u_0(x)=0$ for a.e. $x \in \R$. In 
this case, a first solution to
\eqref{h1dentf}-\eqref{h1space}  is $u(x,t)=0$ for a.e. $(x,t) \in \R\times \Rp$ (it is the entropy weak
solution). A second solution to \eqref{h1dentf}-\eqref{h1space} 
is defined for a.e. $(x,t) \in \R\times \Rp$
by

\[\bal
u(x,t)=0, \text{ if } x < -\sqrt t \text{ or } x > \sqrt t,
\\
u(x,t)=\frac x {2t}, \text{ if } -\sqrt t <x< \sqrt t.
\ean\]

This second solution is not an entropy weak
solution: it does not satisfy \eqref{h1dent}. Also note that this second solution is not in
the space $C(\Rp,L^1_{loc}(\R))$ nor in the space $L^\infty(\R \times 
\Rp)$ (it belongs to $L^\infty(\Rp,L^1(\R)$)). 
Indeed, under the assumption $u \in L^\infty(\R \times \Rp) \cap C(\Rp,L^1_{loc}(\R))$, the solution of
\eqref{h1dentf} is unique.
\end{remark}

The entropy weak solution to
\eqref{hyperbolic1D} satisfies the following $L^\infty$ and $BV$ 
stability\index{Stability}  properties:  

\begin{prop} Let $f \in C^1(\R,\R)$ and $u_0 \in L^\infty(\R)$.
Let $u$ be the entropy weak solution  to
\eqref{hyperbolic1D}. Then, $u \in C(\Rp,L^1_{loc}(\R))$; furthermore, the following estimates
hold:

\begin{enumerate}
\item $\Vert u(\cdot,t) \Vert_{L^\infty(\R)} \le \Vert u_0 
\Vert_{L^\infty(\R)}$, for all $t \in \R_+$.
\item  If $u_0 \in BV(\R)$, then $\vert u(\cdot,t) 
\vert_{BV(\R)} \le  \vert u_0 \vert_{BV(\R)}$, for all $t \in \R_+$.
\end{enumerate}
\end{prop}

\section{Numerical schemes in the linear case}

\hpage{Hyperbolic 1d}

We shall first introduce the numerical schemes in the linear case $f(u) = u$ in
(\ref{hyperbolic1D}). The problem considered in this section is therefore 
\index{Equation!Transport --}
\be
\left\{\begin{array}{rll}
u_t(x,t)
+ u_x(x,t)
   & = 0 &  x \in \R, \ \ t \in \R_+, 
\\
 u(x,0)  & = u_0 (x), &  x \in \R. 
\end{array}\right. \label{lhy1D}
\ee
Assume that  $u_0 \in C^1(\R,\R)$; Problem \eqref{lhy1D}  has a unique classical
solution,  as defined in
Definition \ref{defclass1}, which is $u(x,t)= u_0(x-t)$ for all $(x,t) \in \R \times \Rp$. If
$u_0 \in L^\infty(\R)$, then Problem \eqref{lhy1D} has a unique weak solution, as defined in
Definition \ref{defweak1}, which  is again $u(x,t)= u_0(x-t)$ for a.e. $(x,t) \in \R \times
\Rp$. Therefore, if $u_0 \ge 0$,   the  solution $u$ is also nonnegative.  Hence, it is
advisable for many problems that the solution given  by the numerical scheme should preserve
the nonnegativity of the solution. 

\subsection{The centered finite difference scheme} 

Assume $u_0 \in C(\R,\R)$. Let $h\in\Rpast$ and $x_i = i h$ for
all $i  \in \Z$. Let $k\in\Rpast$ be the time step. With 
the explicit Euler scheme for the
time discretization (the implicit Euler scheme could also be used), the centered finite
difference scheme associated to points $x_i$ and $k$ is 

\be
\left\{\begin{array}{rl}  
\dsp{  \frac  {u_i^{n+1} - u_i^n}  k + \frac  {u_{i+1}^n - u_{i-1}^n} {2 h}  }& =
0, \ \  \forall n
\in \N, \ \ \forall i \in \Z,
 \\
 u_i^0  & = u_0 (x_i) , \ \ \forall i \in \Z.
\end{array}\right. \label{DFCl}
\ee

The discrete unknown $u_i^n$ is expected to be an approximation of $u(x_i,nk)$ where $u$ is
the solution to \eqref{lhy1D}.

It is well known that this scheme should be avoided. In
particular, for the following reasons:

\begin{enumerate}

\item it does not preserve positivity, i.e. $u_i^0 \ge 0$ for all $i \in \Z$ does
not imply  $u_i^1 \ge 0$ for all $i \in  \Z$; take for instance $u_i^0  = 0$ for $i
\le 0$ and  $u_i^0  = 1$ for $i > 0$, then  $u_0^1 = -  k /(2h) < 0$;  

\item it is not ``$L^\infty$-diminishing", i.e.     
$\max\{\vert u_i^0 \vert$, $i \in \Z\} = 1$ does not imply that $\max\{\vert u_i^1 \vert$,
$i \in \Z\} \le 1$; take for instance $u_i^0  =1$ for $i
\le 0$ and  $u_i^0  = 0$ for $i > 0$, then $\max\{\vert u_i^0 \vert$, $i \in
\Z\} = 1$ and $\max\{\vert u_i^1 \vert$, $i \in \Z\}= 1 + k /(2h)$;

\item it is not ``$L^2$-diminishing", i.e. $\sum_{i \in \Z} (u_i^0)^2 = 1$ does not
imply that  $\sum_{i \in \Z} (u_i^1)^2 \le 1$;  take for instance $u_i^0  
= 0$ for $i \neq 0$ and 
$u_i^0  =  1$ for $i = 0$, then  $u_0^1 = 1, u_1^1 = k /(2h), u_{-1}^1 = - k /(2h)$,  so
that $\sum_{i \in \Z} (u_i^1)^2  = 1 +  k^2/(2h^2) > 1$;

\item it is unstable in the von Neumann sense: if the initial condition is taken
under the  form $u_0(x) = \exp(ipx)$, where $p$ is given in $\Z$, then 
$u(x,t) = \exp(-ipt)
\exp(ipx)$ ($i$ is, here, the usual complex number, $u_0$ and $u$ take values in $\C$). Hence
$\exp(-ipt)$ can  be seen as an amplification factor, and its modulus is 1. 
The numerical scheme is stable in the von Neumann sense if the 
amplification factor for the discrete
solution is less than or equal to  1. For the scheme (\ref{DFCl}), we 
have $u^1_j=u^0_j-(u^0_{j+1}-u^0_{j-1}) k/(2h)$ $=\exp(ipjh) \xi_{p,h,k}$, 
with $\xi_{p,h,k} =1-(\exp(iph)-\exp(-iph)) k/(2h)$. 
Hence $\vert \xi_{p,h,k} \vert^2 = 1 +  (k^2/h^2) \sin^2 ph > 1$ if 
$ph \neq q \pi$ for any $q$ in $\Z$. 

\end{enumerate}

In fact, one can also show that there exists $u_0 \in C_c^1(\R,\R)$ such that the solution
given by the numerical scheme does not tend to the solution of the continuous problem 
when $h$ and $k$ tend to $0$ (whatever the relation between $h$ and $k$). 

\begin{remark} \normalfont The scheme \eqref{DFCl} is also a finite volume scheme with 
the (spatial) mesh
$\mesh$ given by $x_{i+1/2}=(i+1/2)h$ in Definition \ref{meshhyp1d} below and with a
centered choice for the approximation of $u(x_{i+1/2},nk)$: 
the value of 
$u(x_{i+1/2},nk)$ is approximated by $(u^n_i+u^n_{i+1})/2$, see \eqref{VFDA} where an 
upstream
choice for $u(x_{i+1/2},nk)$ is performed. In fact, the choice of
$u_i^0$ is different in \eqref{VFDA} and in \eqref{DFCl} but this  does not
change the unstability of the centered scheme.
\end{remark}

\subsection{The upstream\index{Upstream} finite difference scheme}

Consider now a nonuniform distribution of points
$x_i$, i.e. an increasing sequence of real values  $(x_i)_{i\in\Z}$ such that
$\lim_{i\tends\pm\infty} x_i = \pm\infty$. For all $i\in\Z$, we set $h_{i-1/2} = x_i -
x_{i-1}$. 
The time discretization is performed with the explicit Euler scheme with time step $k>0$. 
Still assuming  $u_0 \in C(\R,\R)$,  consider  the upwind (or upstream)  finite difference scheme defined by  

\be
\left\{\begin{array}{rl}  
\dsp{  \frac  {u_i^{n+1} - u_i^n}  k + \frac  {u_{i}^n - u_{i-1}^n} 
{h_\imud}  }& = 0, \ \  \forall n \in \N, \ \ \forall i \in \Z,
 \\
 u_i^0  & = u_0 (x_i) , \ \ \forall i \in \Z.
\end{array}\right. \label{DFDA}
\ee

Rewriting the scheme as
\[
\dsp{  u_i^{n+1} = (1 - \frac  k {h_\imud})  u_{i}^n +  \frac  k {h_\imud} u_{i-1}^n},
\]
it appears that if $\inf_{i\in \Z} h_{i-1/2} > 0$ and if $k$ is such that 
\begin{equation}
	k \le \inf_{i\in \Z} h_{i-1/2} 
\label{cfl-df}
\end{equation}
then $u_i^{n+1}$ is a convex combination of $u_{i}^n$ and $ u_{i-1}^n$.  
By induction, this proves that the scheme \eqref{DFDA} is stable, in the sense that if $u_0$ is such that $U_m \le u_0(x) \le U_M$  for a.e. $x \in \R$, where $U_m$, $U_M \in \R$, then  $U_m \le u_i^n\le U_M$  for any $i \in \Z$ and $n \in \N$. 

Moreover, if $u_0 \in C^2(\R,\R) \cap L^\infty(\R)$ and $u_0'$ and $u_0''$ belong to $L^\infty(\R)$, it is easily shown that the scheme is consistent in the finite  difference sense\index{Consistency!-- in the finite difference sense}; indeed, the consistency error \index{Consistency!-- error} defined by
\be 
	R_i^n  =  \frac {\uexinp- \uexin} k +  \frac {\uexin- \ueximn} \himud 
	\label{consisthyplin}
\ee
is such that, if the CFL condition \eqref{cfl-df} holds, $\vert R_i^n \vert \le C  h$, where $h = \sup_{i \in \Z} h_i$ and $C \ge 0$ only depends on $u_0$ (recall that $u$ is the solution to problem \eqref{lhy1D}). 
Hence the following error estimate holds.

\begin{prop}[Error estimate for the upwind finite difference 
scheme]\label{eefd}\index{Error estimate!-- for a hyperbolic 
equation!-- in the one dimensional case}~

Let $u_0 \in C^2(\R,\R) \cap L^\infty(\R)$,
such that $u_0'$ and $u_0''$ $\in$  $L^\infty(\R)$.
Let $(x_i)_{i\in\Z}$ be an increasing sequence of real values 
such that $\lim_{i\tends\pm\infty} x_i = \pm\infty$. Let
$h = \sup_{i\in\Z} h_\imud$, and  assume that $h <\infty$ and $\inf_{i\in \Z} h_{i-1/2} >
0$. Let $k>0$ such that $k \le \inf_{i\in \Z} h_{i-1/2}$. 
Let $u$ denote the unique solution to \eqref{lhy1D} and 
$\{\uin$, $i \in \Z$, $n \in \N\}$ be given by \eqref{DFDA}; 
let $\ein = \uexin - \uin$, for any $n \in \N$
and $i \in \Z$, and  let $T \in ]0,+\infty[$ (note that $\uexin$ is well 
defined since $u \in C^2(\R \times \Rp,\R)$).

Then there exists $C \in \R_+$, only depending on $u_0$, such that 
$ \abs{\ein} \le C h T,  $ for any $n \in \N$  such that $nk \le T$, and for any $i
\in \Z$.
\end{prop}

\proof of Proposition \ref{eefd}

\smallskip

Let $i\in\Z$ and $n\in\N$. By definition of the consistency error 
$\Rin$ in \eqref{consisthyplin}, the error $\ein$ satisfies

\[
\frac {\einp- \ein} k +  \frac
{\ein- \eimn} \himud = R_i^n.
\]

Hence
\[
\einp = \ein(1 - \frac  k \himud) +    \frac  k \himud \eimn  + k \Rin.
\]

Using $\vert R_i^n \vert  \le C  h$ (for some $C$ only depending on 
$u_0$) and the assumption 
$k \le \inf_{i\in \Z} h_\imudt$, this yields

\[
\vert \einp \vert \le \sup_{j \in \Z} \vert e_j^n \vert + C k h.
\]

Since $e_i^0=0$ for any $i \in \Z$, an induction yields 

\[
\sup_{i \in \Z} \vert \ein \vert \le Cnkh
\] 

and the result follows. 
\eof

Note that in the above proof, the linearity of the equation and the
regularity of  $u_0$ are used. The next questions to arise are what to do in the case of 
a nonlinear equation and in the case $u_0 \in L^\infty(\R)$.

\subsection{The upwind finite volume scheme} 

 Let us
first give a definition of the admissible meshes  for the finite volume schemes.

\begin{defi}[One-dimensional admissible mesh]  \normalfont
\index{Mesh!Admissible --!--for hyperbolic equations}
An admissible mesh $\mesh$ of $\R$ is given
by an increasing sequence of real values $(x_{i+1/2})_{i \in \Z}$, such that 
$\R = \cup_{i\in\Z} [x_{i-1/2}, x_{i+1/2}]$.
The mesh $\mesh$ is the set $\mesh = \{K_i,\ i\in\Z\}$ of subsets of 
$\R$ defined by $\cv_i = (x_{i-1/2}, x_{i+1/2})$ for all $i \in \Z$.
The length of $\cv_i$ is denoted by $h_i$, so that $h_i = x_{i+1/2} - x_{i-1/2}$ 
for all $i\in \Z$. It is assumed that $h=\sizemesh=\sup\{h_i,$ $i\in \Z\} 
< +\infty$ and that, for some $\a \in \Rpast$, $\a h \le \inf\{h_i,$ $i\in \Z\}$. 
\label{meshhyp1d}
\end{defi}

Consider an admissible mesh in the sense of Definition \ref{meshhyp1d}.
Let $k\in\Rpast$ be the time step.
Assume $u_0 \in L^\infty(\R)$ (this is a natural hypothesis for the finite volume framework).
Integrating \eqref{lhy1D} on each control volume of the
mesh,  approximating the time derivatives by differential quotients and using
an upwind choice for $u(x_\ipud,nk)$ yields the following (time explicit) 
scheme\index{Finite volume scheme!-- for hyperbolic problems!-- in one space 
dimension}: 

\be
\left\{\begin{array}{rl}  
\dsp{  h_i \frac  {u_i^{n+1} - u_i^n}  k + u_{i}^n - u_{i-1}^n }&
= 0, \ \ \forall n
\in \N, \ \ \forall i \in \Z,
\\
u_i^0  & = \dsp{ \frac 1 {h_i} \int_{\cv_i} u_0(x) dx,} \ \ \forall i \in \Z.
\end{array}\right. 
\label{VFDA}
\ee

The value $u_i^n$ is expected to be  an approximation of $u$ (solution to \eqref{lhy1D}) in
$\cv_i$ at time $nk$. It is easily shown that 
this scheme is not consistent \index{Consistency!-- in the finite difference 
sense} in the finite difference sense
if $u_i^n$ is considered to be an approximation of $u(x_i,nk)$ with, for instance,
$x_i=(x_{i-1/2}+x_{i+1/2})/2$ for all $i \in \Z$. Even if $u_0 \in C^\infty_c(\R,\R)$, the
quantity $R_i^n$ defined by \eqref{consisthyplin} does not satisfy (except 
in particular cases)
$\vert R_i^n \vert \le Ch$, with some $C$ only depending on $u_0$.

\smallskip

It is  however  possible to interpret this scheme as another 
expression of the upwind finite difference scheme \eqref{DFDA} (except for the minor
modification of $u_i^0$, $i \in \Z$). One simply needs to consider $u_i^n$ as an approximation of
$u(x_{i+1/2},nk)$ which leads to a consistency property in the finite difference sense.
Indeed, taking $x_j=x_{j+1/2}$ (for $j=i$ and $i-1$) in the definition \eqref{consisthyplin} of
$R_i^n$ yields $\vert R_i^n \vert \le Ch$, where $C$ only depends on $u_0$. Therefore, a
convergence result for this scheme is given by the proposition 
\ref{eefd}. This analogy cannot be extended to the general case 
of ``monotone flux schemes" (see Definition \refp{monflux} below)
for a nonlinear equation for which there may be no value of $x_i$ (independent of $u$) leading to such a consistency property, see Remark \refp{vfordf} for a 
counterexample (the analogy holds however for the scheme \eqref{up1D}, 
convenient for a nondecreasing function $f$, see Remark \ref{vfdf}).

\medskip

The approximate finite volume solution $\uapp$ may be defined 
 on $\R \times \Rp$ from 
the discrete unknowns $u_i^n$, $i \in \Z$, $\nnn$ which  are computed in 
\eqref{VFDA}:
 
\be 
\uapp(x,t) = \uin \text{ for } x \in \cv_i \text{  and  }  t \in
[nk,(n+1)k). \label{uapp}
\ee 

The following $L^\infty$ estimate holds:

\begin{lemma}[$L^\infty$ estimate in the linear case] \label{linstab}
\index{Estimate on the approximate solution!-- for a hyperbolic equation}
Let $u_0 \in L^\infty(\R)$ and $U_m$, $U_M\in\R$ such that $U_m \le 
u_0(x) \le U_M$ for a.e. $x\in\R$.  Let $\mesh$  be an admissible
mesh in the sense of Definition \ref{meshhyp1d} and let $k\in\Rpast$
satisfying the Courant-Friedrichs-Levy (CFL) condition\index{CFL condition} 

\[
k \le \inf_{i \in \Z} h_i.
\]

(note that taking  $k\le\a h$ implies the above condition). Let 
$\uapp$ be the finite volume approximate solution defined by
\eqref{VFDA} and \eqref{uapp}. 

Then,

\[
U_m \le  \uapp(x,t) \le U_M \, \text{ for a.e. } x\in\R \text{ and a.e. } 
t\in\Rp.
\]

\end{lemma}

\proof of Lemma \ref{linstab}

\smallskip

The proof that $U_m \le u_i^n \le U_M$, for all $i\in\Z$ and $n\in\N$,
as in the case of the 
upwind finite difference scheme (see \refep{DFDA}), consists in remarking that
equation \eqref{VFDA} gives, under the CFL condition, an expression of
$u_i^{n+1}$ as a linear convex combination of $u_i^n$ and $u_{i-1}^n$,
for all $i\in\Z$ and $n\in\N$. \eof

The following  inequality  will be crucial for the proof of convergence.
 
\begin{lemma}[Weak $BV$ estimate, linear case]\label{bvwhyplin}
\index{BV!Weak -- estimate}
Let $\mesh$ be an admissible mesh in the sense of Definition \refp{meshhyp1d} and let $k\in\Rpast$ satisfying the CFL condition 
\be
k \le (1-\xi) \inf_{i \in \Z} h_i,
\label{cflhypl}
\ee 
for some $\xi\in(0,1)$ (taking $k \le (1-\xi)\a h$ implies this condition).

Let $\{u_i^n$, $i\in\Z$, $n\in\N\}$ be given by the finite volume scheme \eqref{VFDA}. 
Let $R\in\Rpast$ and $T\in\Rpast$ and assume $h=\sizemesh < R$, $k < T$. Let $i_0 \in \Z$, $i_1 \in \Z$ and $N \in \N$ be such that $-R\in \cvbar_{i_0}$,  $R\in \cvbar_{i_1}$ and $T \in (Nk,(N+1)k]$ (note that $i_0 < i_1$).

Then there exists $C\in\Rpast$, only depending on $R$, $T$, $u_0$, $\a$ and $\xi$, such that

\be
\sum_{i=i_0}^{i_1} \sumN k \vert u_i^n - u_{i-1}^n\vert \le C h^{-1/2}.
\label{bvf}\ee
\end{lemma}

\proof of Lemma \ref{bvwhyplin}

\smallskip

Multiplying the first equation of \eqref{VFDA} by $k u_i^n$ and summing on
$i=i_0,\ldots,i_1$ and $n=0,\ldots N$ yields $A+B = 0$ with

\[
A = \sum_{i=i_0}^{i_1} \sumN h_i (u_i^{n+1} - u_i^n) u_i^n
\]

and 

\[
B = \sum_{i=i_0}^{i_1} \sumN k (u_i^n - u_{i-1}^n) u_i^n.
\]

Noting that

\[
A = - \half \sum_{i=i_0}^{i_1} \sumN h_i (u_i^{n+1} - u_i^n)^2
+ \half \sum_{i=i_0}^{i_1} h_i [(u_i^{N+1})^2 - (u_i^0)^2]
\]

and using the scheme \eqref{VFDA} gives

\[
A = - \half \sum_{i=i_0}^{i_1} \sumN \frac{k^2} {h_i} (u_i^n - u_{i-1}^n)^2
+ \half \sum_{i=i_0}^{i_1} h_i [(u_i^{N+1})^2 - (u_i^0)^2];
\]

therefore, using the CFL condition \eqref{cflhypl},

\[
A \ge - (1-\xi) \half \sum_{i=i_0}^{i_1} \sumN k (u_i^n - u_{i-1}^n)^2
- \half \sum_{i=i_0}^{i_1} h_i (u_i^0)^2.
\]

We now study the term $B$, which  may be rewritten as

\[
B = \half \sum_{i=i_0}^{i_1} \sumN k (u_i^n - u_{i-1}^n)^2
+ \half \sumN k [(u_{i_1}^n)^2 - (u_{i_0-1}^n)^2].
\]

Thanks to the $L^\infty$ estimate of  Lemma \refp{linstab}, this last equality implies that

\[
B \ge \half \sum_{i=i_0}^{i_1} \sumN k (u_i^n - u_{i-1}^n)^2
- T \max\{-U_m,U_M\}^2.
\]

Therefore, since $A+B=0$ and $\sum_{i=i_0}^{i_1}h_i \le 4R$, the following inequality holds:

\[
0 \ge \xi \sum_{i=i_0}^{i_1} \sumN k (u_i^n - u_{i-1}^n)^2
- (4R + 2T)  \max\{-U_m,U_M\}^2,
\]

which, in turn, gives the existence of $C_1\in\Rpast$, only depending on $R$, $T$, $u_0$ and
$\xi$ such that

\be
\sum_{i=i_0}^{i_1} \sumN k (u_i^n - u_{i-1}^n)^2 \le C_1.
\label{bvfsq}
\ee

Finally, using  
\[
\sum_{i=i_0}^{i_1} 1 \le \sum_{i=i_0}^{i_1} \frac {h_i} 
{\a h} \le \frac {4R} {\a h},
\] 

the Cauchy-Schwarz inequality leads to

\[
[\sum_{i=i_0}^{i_1} \sumN k \vert u_i^n - u_{i-1}^n \vert ]^2 \le C_1 2 T 
\frac {4R} {\a h},
\]

which concludes the proof of the lemma.
\eof

Contrary to the discrete $H^1_0$ estimates which were obtained on the approximate finite volume solutions of elliptic equations, see e.g. 
\eqref{ellestimx}, the  weak $BV$ estimate\index{BV!Weak -- estimate} \eqref{bvf} is not related to an \textit{a priori} estimate on the
solution to the continuous problem \eqref{lhy1D}. 
It does not give any compactness property  \index{Compactness!-- of the approximate solutions in $L^1_{loc}$} in
the space $L^1_{loc}(\R)$ (there are some counterexamples); such a compactness property is obtained thanks to a ``strong $BV$ estimate" (with, for instance, an $L^\infty$ estimate) as it is recalled below (see Lemma \ref{helly}). 
In the one-dimensional case which is studied here such a ``strong $BV$ estimate" can be obtained if $u_0 \in BV(\R)$, see Corollary \ref{nlstabvco}; this is no longer true in the multidimensional case with general meshes, for which only the above weak $BV$ estimate is available.

\begin{remark} \normalfont\label{rbvweff} The weak $BV$ estimate is a crucial point for the proof of
convergence. Indeed, the property which is used in the proof of 
convergence (see Proposition \ref{convhyplin} below) is, with the notations of Lemma \ref{bvwhyplin},

\be
h\sum_{i=i_0}^{i_1} \sumN k \vert u_i^n - u_{i-1}^n \vert \tends 0, \text{ as } h \tends 0,
\label{bvweff}
\ee

for $R$, $T$, $u_0$, $\a$ and $\xi$ fixed.

\smallskip

If a  piecewise constant function $u_{\mesh,k}$, such as given by \eqref{uapp} 
(with some $u^n_i$ in $\R$, not
necessarily given by \eqref{VFDA}), is bounded in (for instance) $L^\infty(\R \times \Rp)$
and converges in $L^1_{loc}(\R \times \Rp)$ as $h \tends 0$ and $k \tends 0$ (with a
possible relation between $k$ and $h$) then \eqref{bvweff} holds. This proves that the
hypothesis \eqref{bvweff} is included in the hypotheses of the classical Lax-Wendroff theorem of
convergence (see Theorem \refp{laxw}); note that   \eqref{bvweff} is implied by 
\eqref{bvf}  and that it is weaker 
than \eqref{bvf}).
\end{remark}

 We show in the  following remark how the `` weak" and `` strong" $BV$  estimates   may  ``formally" be obtained on the ``continuous equation";  this gives a hint of the reason why this estimate may be obtained even if the exact solution does not belong to the space $BV(\R \times
\Rp)$. A similar remark also holds in the nonlinear case (i.e. for Problem \eqref{hyperbolic1D}).

\begin{remark}[Formal derivations of the strong and weak $BV$ estimates] \normalfont~
\index{Estimate on the approximate solution!-- for a hyperbolic equation} 
When approximating the solution to \eqref{lhy1D} by the finite volume 
scheme \eqref{VFDA} (with $h_i=h$ for all $i$, for the sake of simplicity), 
the equation to which an approximation of a solution is sought 
is ``close" to the equation

\be u_t + u_x - \eps u_{xx}  = 0 \label{rem} \ee
where $\eps = \frac{h-k} 2$ is positive under the CFL condition \eqref{cflhypl}, which ensures 
that the scheme is diffusive.

We assume that $u$ is regular enough, with null limits for $u(x,t)$ and 
its derivatives as $x\tends\pm\infty$.

\medskip

\textit{(i)} ``Strong" $BV$ estimate\index{BV!Strong -- estimate}.

 Derivating the equation \eqref{rem} with respect to the variable $x$, multiplying by
$\sgn_r(u_x(x,t))$, where $\sgn_r$ denotes a nondecreasing 
 regularization of the function $\sgn$,  and integrating over $\R$  yields
 
 \[\bigl(\int_\R \phi_r(u_x(x,t))  dx \bigr)_t + 
 \int_\R  u_{xx}(x,t) \sgn_r(u_x(x,t))dx = - \eps 
 \int_\R \sgn_r'(u_x(x,t)) (u_{xx}(x,t))^2 dx \le 0,\]

where $\phi_r'=\sgn_r$ and $\phi_r(0)=0$. Since

\[
\int_\R  u_{xx}(x,t) \sgn_r(u_x(x,t))dx=\int_\R  (\phi_r(u_x(x,t)))_x dx=0,
\]

this yields, passing to the limit on the regularization, 
that $\norm{u_x(\cdot,t)}{L^1(\R)}$ is 
nonincreasing with respect to $t$. Copying this formal proof on the 
numerical scheme yields a strong $BV$ estimate, which is an \textit{a priori}
estimate giving compactness properties in $L^1_{loc}(\R \times \Rp)$, see
Lemma \ref{nlstabv}, Corollary \ref{nlstabvco} and Lemma \refp{helly}.

\medskip

\textit{(ii)} ``Weak" $BV$ estimate\index{BV!Weak -- estimate}
  
Multiplying  \eqref{rem} by $u$ and summing over $\R\times (0,T)$ yields

\[
\frac 1 2 \int_\R u^2(x,T) dx- \frac 1 2 \int_\R u^2(x,0) dx + 
\int_0^T \int_{\R} \eps u_x^2(x,t) dx dt = 0,
\]

which yields in turn

\[
\eps \int_0^T \int_{\R} u_x^2(x,t) dx dt \le 
\frac 1 2 \Vert u_0\Vert^2_{L^2(\R)}.
\]

This is the continuous analogous of \eqref{bvfsq}. 
Hence if $h - k = \eps \ge \xi h$ (this is Condition \eqref{cflhypl}, note that this 
condition is more restrictive than the usual CFL condition required for 
the $L^\infty$ stability), the discrete equivalent of this formal proof 
yields \eqref{bvfsq} (and then \eqref{bvf}). 

\medskip

In the first case, we derivate the equation and we use some regularity on 
$u_0$ (namely $u_0 \in BV(\R)$). In the second case, it is sufficient to 
have $u_0 \in L^\infty(\R)$ but we need the 
diffusion term to be large enough in order to obtain 
the estimate which, by the way, does not yield any estimate on the 
solution of \eqref{rem} with $\eps=0$. This formal derivation may be carried out similarly in 
the nonlinear case.

\end{remark}

Let us now  give a convergence result for the scheme \eqref{VFDA} in $L^\infty(\R\times\Rpast)$
for the weak-$\star$ topology. 
Recall that a sequence $(v_n)_\nnn \subset L^\infty(\R\times\Rpast)$ 
converges to $v \in L^\infty(\R\times\Rpast)$
in $L^\infty(\R\times\Rpast)$ for the
weak-$\star$ topology\index{Convergence!-- for the weak star topology}
if

\[
\int_{\Rp} \int_\R (v_n(x,t) - v(x,t)) 
\varphi(x,t) dx dt \tends 0 \text{ as } n \tends \infty, \, \forall 
\varphi\in L^1(\R\times\Rpast).
\]

A stronger convergence result is available, and comes from the 
nonlinear study given in Section \ref{hyp1dnl}.

\begin{prop}[Convergence in the linear case]\label{convhyplin}
\index{Convergence of the finite volume method!-- for a linear hyperbolic 
equation}

Let $u_0  \in L^\infty(\R)$ and $u$ be the unique weak solution to Problem
\refep{lhy1D} in the sense of
Definition \refp{defweak1}, with $f(s) = s$ for all $s \in \R$. Let $\xi 
\in (0,1)$ and $\a >0$ be given. Let  $\mesh$ be an
admissible mesh in the sense of
Definition \refp{meshhyp1d} and let $k\in\Rpast$
satisfying the CFL condition \refep{cflhypl} (taking
$k\le (1-\xi)\a h$ implies this condition, note that $\xi$ and $\a$ do 
not depend on $\mesh$).

Let $u_{\mesh,k}$ be the finite volume 
approximate solution defined by
\eqref{VFDA} and \eqref{uapp}. 
Then $\uapp \tends u$ in $L^\infty(\R\times\Rpast)$ for the weak-$\star$
topology as $h=\sizemesh \tends 0$.
\end{prop}

\proof of Proposition \ref{convhyplin}

\smallskip

Let $(\mesh_m,k_m)_{m\in\N}$ be a sequence of meshes and
time steps satisfying the hypotheses of Proposition \ref{convhyplin} and 
such that $\size(\mesh_m) \tends 0$ as $m \tends \infty$.

Lemma \ref{linstab} gives the existence of a subsequence, still denoted 
by $(\mesh_m,k_m)_{m\in\N}$, and
of a function $u\in L^\infty(\R\times\Rpast)$ such that
$u_{\mesh_m,k_m} \tends u$ in $L^\infty(\R\times\Rpast)$ for the weak-$\star$
topology, as $m\tends +\infty$. There remains to show that $u$ is the
solution of \eqref{h1dweak} (with $f(s) = s$ for all $s \in \R$). The uniqueness of the weak
solution to Problem \eqref{lhy1D} will then imply that the full sequence
converges to $u$.

\smallskip

Let $\varphi \in C^1_c(\R \times \R_+,\R)$. Let $m\in\N$ and 
$\mesh = \mesh_m$, $k=k_m$ and $h=\sizemesh$.
Let us multiply the first equation of \eqref{VFDA} by 
$(k/h_i)\varphi(x,n k)$,   integrate over $x\in\cv_i$ and  sum for
all $i\in\Z$ and $n\in\N$. This yields

\[
A_m + B_m = 0
\] 
with

\[
A_m = \sum_{i\in\Z} \sumIN (u_i^{n+1} - u_i^n) \int_{\cv_i} \varphi(x,n k) dx
\]

and

\[
B_m = \sum_{i\in\Z} \sumIN k (u_i^n - u_{i-1}^n) 
 \frac{1} {h_i} \int_{\cv_i} \varphi(x,n k) dx.
\]

Let us remark that  $A_m = A_{1,m} -  A'_{1,m}$ with
 
\[
A_{1,m} = - \int_{k}^\infty \int_\R \uapp(x,t) \varphi_t(x,t-k) dx dt - \int_\R 
u_0(x) \varphi(x,0) dx
\]

and 

\[
A'_{1,m} = \sum_{i\in\Z} u_i^{0} \int_{\cv_i} \varphi(x,0) dx - 
\int_\R u_0(x) \varphi(x,0) dx.
\]

Using the fact that $\sum_{i\in\Z} u_i^0 1_{\cv_i} \tends u_0$ in
$L^1_{loc}(\R)$ as $m\tends\infty$, we get that $A'_{1,m}\tends 0$
as $m\tends\infty$. (Recall that $1_{\cv_i}(x)=1$ if $x \in \cv_i$ and $1_{\cv_i}(x)=0$ if $x
\notin \cv_i$.)

Therefore, since $\uapp \tends u$ in $L^\infty(\R \times \Rpast)$ for the 
weak-$\star$ topology   as $m \tends \infty$, and $\vphi_t(\cdot,\cdot-k)1_{\R \times (k,\infty)} 
\tends \vphi_t$ in $L^1(\R \times \Rpast)$ (note that $k \tends 0$ thanks 
to \eqref{cflhypl}),

\[
\lim_{m \tends + \infty} A_m = \lim_{m \tends + \infty} A_{1,m} = 
- \int_{\Rp} \int_\R u(x,t) \varphi_t(x,t) dx dt - \int_\R 
u_0(x) \varphi(x,0) dx.
\]

Let us now turn  to the study of $B_m$. We compare $B_m$ with

\[
B_{1,m} = - \sumIN \int_{nk}^{(n+1)k} \int_{\R} \uapp(x,t) 
\varphi_x(x,nk) dx dt,
\]

which tends to 
$- \int_{\Rp} \int_\R u(x,t) \varphi_x(x,t) dx dt$ as
$m\tends\infty$. The term $B_{1,m}$ can be rewritten as

\[
B_{1,m} =   \sum_{i\in\Z} \sumIN k (u_i^n - u_{i-1}^n)  
 \varphi(x_\imud,n k).
\]

Let $R>0$ and $T>0$ be such that $\varphi(x,t)=0$ if $\vert x \vert \ge 
R$ or $t \ge T$. Then, there exists $C\in\Rpast$, only depending  on $\varphi$, such that, if $h
< R$ and $k < T$ (which is true for $h$ small enough, thanks to \eqref{cflhypl}),

\be
\vert B_m - B_{1,m} \vert \le C h 
\sum_{i=i_0}^{i_1} \sumN k \vert u_i^n - u_{i-1}^n\vert,
\label{eqcvhl2}\ee

where $i_0 \in \Z$, $i_1 \in \Z$ and $N \in \N$ are such that $-R \in \cvbar_{i_0}$, $R
\in \cvbar_{i_1}$ and $T \in (Nk,(N+1)k]$. 

Using \eqref{eqcvhl2} and Lemma \ref{bvwhyplin}, we get that $B_m - B_{1,m} \tends 0$ and then 

\[
B_m \tends - \int_{\Rp} \int_{\R} 
u(x,t) \varphi_x(x,t) dx dt \text{ as } m\tends\infty,
\]

which completes the proof that
$u$ is the weak solution to Problem \refep{lhy1D} (note that here  the useful consequence of lemma
\ref{bvwhyplin} is \eqref{bvweff}). \eof

\begin{remark} \normalfont In Proposition \ref{convhyplin}, a simpler proof of convergence could be achieved,
with  $\xi = 0$, using a multiplication of the first equation of \eqref{VFDA} by
$(k/h_i)\varphi(x_{i-1/2}, nk)$. However, this proof does not
generalize to the general case of nonlinear hyperbolic problems.
\end{remark}

\begin{remark} \normalfont
Proving the convergence of the finite difference method (with the scheme \eqref{DFDA}) with $u_0\in
L^\infty(\R)$ can be done using the same technique as the proof
of the finite volume method (that is considering the finite difference scheme as a finite
volume scheme on a convenient mesh).
\end{remark}

\section{The nonlinear case} \label{hyp1dnl}

\hpage{Hyperbolic 1d}

In this section, finite volume schemes for the discretization of Problem \eqref{hyperbolic1D}
are presented and a theorem of convergence is given (Theorem 
\ref{nlcv1d}) which will be generalized to the multidimensional case in 
the next chapter. We also recall the classical proof of convergence which 
uses a ``strong $BV$ estimate" and the Lax-Wendroff theorem. This proof, 
however, does not seem to extend to the 
multidimensional case for general meshes.
The following properties are assumed to be satisfied by the data of problem \eqref{hyperbolic1D}.

\begin{hyp}\label{datah1d}
The flux function $f$ belongs to $C^1(\R,\R)$, the initial data $u_0$ belongs
to $L^\infty(\R)$ and $U_m$, $U_M \in \R$ are such that $U_m \le u_0 \le U_M$ a.e. on $\R$.
\end{hyp}

\subsection{Meshes and schemes}\label{ms1d}

Let $\mesh$ be an admissible mesh in the sense of Definition
\refp{meshhyp1d} and $k\in\Rpast$ be the time step. 
In the general nonlinear case, 
the finite volume scheme \index{Finite volume scheme!-- for hyperbolic 
problems!-- in one space dimension}
 for the discretization of Problem
\refep{hyperbolic1D} reads

\be
\left\{\begin{array}{rl}  
\dsp{  \frac {h_i} k  (u_i^{n+1} - u_i^n) + f_{i+1/2}^n - f_{i-1/2}^n}  
& = 0, \ \ \forall n \in \N, \ \ \forall i \in \Z,
 \\
 u_i^0  & = 
\dsp{  \frac 1 {h_i} \int_{x_{i-1/2}}^{x_{i+1/2}} u_0 (x) dx} , 
\ \ \forall i \in \Z,
\end{array}\right.
\label{h1D}
\ee

where  $u_i^n$ is expected to be an approximation of $u$ at time $t_n = nk$ in cell $K_i$.
The quantity  $f_{i+1/2}^n$ is often called the numerical flux at point $x_{i+1/2}$ 
and time $t_n$ (it is expected to be an approximation of $f(u)$ at point
$x_{i+1/2}$ and time $t_n$).  
Note that a common expression of $f_{i+1/2}^n$ is used for both equations $i$ and
$i+1$ in (\ref{h1D}); therefore the scheme (\ref{h1D}) satisfies the
property of conservativity, common to all finite volume schemes.
%
%
In the case of a so-called $2p +1$ point scheme with  ($p \in \N^\star$), the numerical flux may be written   
\be
f_{i+1/2}^n = g_\ipud^n (u_{i-p+1}^n,\ldots,u_{i+p}^n),
\label{nf} 
\ee
where $g_\ipud^n$ is the numerical flux function at point $x_{\ipud}$ and time $t_n$, which determines the scheme. 
Note that the numerical flux may thus depend on the interface and the time. 
This is important in applications, for instance in the case of boundary faces (see section \ref{hyppb}) or  interfaces coupling different domains.
For $p=1$, the flux reads
\be
f_{i+1/2}^n = g_\ipud^n (u_{i}^n,\ldots,u_{i+1}^n),
\label{nf-twopoints} 
\ee
and yields a 3-point scheme. 

As in the linear case \refep{uapp}, 
the approximate finite volume solution  is defined by
 
\be 
\uapp(x,t) = \uin \text{ for } x \in \cv_i \text{  and  }  t \in
[nk,(n+1)k). \label{nluapp}
\ee 

The property of consistency 
 for the finite volume scheme 
\eqref{h1D}, \eqref{nf} with $2p+1$ points, is 
ensured by writing the following condition: 

\be
g(s,\ldots,s) = f(s), \ \ \forall s \in \R.
\label{consflux}
\ee

This condition is equivalent to writing the consistency of the 
approximation of the flux (as in
the elliptic and parabolic cases, which were described in the previous 
chapters, see e.g. Section \ref{Diric1D}). 

\begin{remark}[Finite volumes and finite differences] \label{vfordf}  \normalfont
We can remark that, as in the elliptic case, the condition
(\ref{consflux}) does not generally give the consistency  of the scheme 
(\ref{h1D}) when it is considered as a finite difference scheme. 
For instance, assume $f(s)=s^2$ for all $s \in \R$, $p=1$ and $g(a,b) = 
f_1(a)+f_2(b)$ for all $a$, $b\in\R$ with $f_1(s)=\max\{s,0\}^2$, 
$f_2(s)=\min\{s,0\}^2$ (which is shown below to be a  good choice, see 
Example \ref{flumono}). Assume also $h_{2i} = h$ and $h_{2i+1} = h/2$ for all $i\in\Z$. 
In this case, there is no choice
of points $x_i \in \R$ such that the quantity
$(f_{i+ 1/2}^n-f_{i-1/2}^n) / h_i$ is an approximation 
of order $1$ of $(f(u))_x(x_i, n k)$, for any regular function $u$, when 
$u_i^n = u(x_i,n k)$ for all $i \in \Z$. Indeed, up to second order terms, this property of consistency is achieved if and only if 
$f_2'(a)\vert x_{i+1} -x_i \vert + f_1'(a)\vert x_{i-1} -x_i \vert=f'(a) 
h_i$ for all $i \in \Z$ and for all $a \in \R$. Choosing $a>0$ and $a 
<0$, this condition leads to 
$\vert x_{i+1} -x_i \vert=h_i$ and $\vert x_{i+1} -x_i \vert=h_{i+1}$ for 
all $i \in \Z$, which is impossible.
\end{remark}

Examples of convenient choices
for the function $g$ will now be given. 
An interesting class of schemes is the class of $3$-points schemes with a monotone flux, 
which we now define. 


\begin{defi}[Monotone flux schemes] \label{monflux}  \normalfont

\index{Scheme!Monotone flux --|(}
Let $p=1$. 
Under Assumption 
\ref{datah1d}, the finite volume scheme \eqref{h1D}\--\-\eqref{nf-twopoints} is said 
to be a ``monotone flux scheme" if the function $g$, 
only depending on $f$, $U_m$ and $U_M$,  satisfies the following assumptions:

\begin{itemize} 
\item $g$ is locally Lipschitz continuous from $\R^2$ to $\R$,
\item $g(s,s) = f(s)$, for all $s \in [U_m,U_M]$,
\item $(a,b) \mapsto g(a,b)$, from $[U_m,U_M]^2$ to $\R$, is nondecreasing with respect to $a$
and nonincreasing with respect to $b$.
\end{itemize}
\end{defi}

The monotone flux property seems to be remarkable; indeed, monotone flux schemes are consistent in the finite volume sense, they are $L^\infty$-stable under a condition (the so called Courant-Friedrichs-Levy condition) of the type $k \le C_1 h$, where $C_1$ only depends on $g$ and $u_0$ (see Section \refp{linfstah1d} below),  and they are ``consistent with the entropy inequalities" also under a condition of the type $k \le C_2h$, where $C_2$ only depends on $g$ and $u_0$ (but $C_2$ may be different of $C_1$, see Section \refp{discenth1d}). 

\begin{remark} \normalfont\label{schmono} A monotone flux scheme is a monotone scheme, under a Courant-Friedrichs-Levy condition, which means that the scheme can be written under the form 
\[
\uinp = H(\uimn,\uin, \uipn), 
\]
with $H$ nondecreasing with respect to its three arguments.
\end{remark}

\begin{example}[Examples of monotone flux schemes] \label{flumono}\normalfont (see also
\citet{godlewski-springer}, \citet{leveque} and references therein). 
Under Assumption \ref{datah1d},  here are some numerical flux functions $g$ for which
the  finite volume scheme \eqref{h1D}-\eqref{nf} is a monotone flux scheme (in the sense
of Definition \ref{monflux}):

\begin{itemize}

\item the flux splitting scheme: assume $f = f_1 +f_2$, with $f_1$, $f_2 \in C^1(\R,\R)$, $f_1'(s) \ge 0$ and  $f_2'(s) \le 0$
for all $s \in [U_m,U_M]$ (such a decomposition for $f$ is always possible, 
see the modified Lax-Friedrichs  scheme below), and take

\[
g(a,b) = f_1(a) +f_2(b).
\] 

Note that if $f' \ge 0$, taking $f_1=f$ and $f_2=0$, the flux splitting scheme boils down to
the upwind scheme, i.e. $g(a,b) = f(a)$.

\item the Godunov scheme: the Godunov scheme, which was introduced in \citet{Go}, may  be summarized by
the following expression.

\be
g(a,b) = 
\left\{
\begin{array}{l} 
\min \{f(\xi), \xi \in [a,b]\} \text{ if } a \le b, \\
\max\{f(\xi), \xi \in [b,a]\} \text{ if } b \le a. 
\end{array}\right.
\label{fluxgodunov}
\ee

\item the modified Lax-Friedrichs scheme\index{Lax-Friedrichs scheme}
\index{Scheme!Lax-Friedrichs}:  take

\be g(a,b) = \frac {f(a) +f(b)} 2 + D(a-b),
\label{fluxLF}
\ee
with $D \in \R$ such that $ 2 D \ge \max \{\abs{f'(s)}, s \in [U_m,U_M]\}$. Note that in this
modified version of the Lax-Friedrichs  scheme, the coefficient $D$ only depends 
on $f$, $U_m$ and $U_M$, while the original Lax-Friedrichs  scheme consists
in taking $D = h/(2k)$, in the case $h_i = h$ for all $i \in \N$, and therefore satisfies
the three items of Definition \ref{monflux} under the condition $  h/k \ge
\max \{\abs{f'(s)}, s \in [U_m,U_M]\}$. However, an  inverse CFL condition
appears to be necessary for the convergence of the original Lax-Friedrichs  
scheme (see remark \refp{estext}); such a condition is not necessary for the modified 
version.

Note also that the modified Lax-Friedrichs  scheme consists in a 
particular flux splitting scheme with  $f_1(s)=(1/2)f(s)+Ds$
and $f_2(s)=(1/2)f(s)-Ds$ for $s \in [U_m,U_M]$.
\end{itemize}
\end{example}

\begin{remark} \normalfont \label{vfdf}
In the case of a nondecreasing (resp. nonincreasing) function $f$, the Godunov
monotone flux scheme \eqref{fluxgodunov} reduces to
$g(a,b) = f(a)$ (resp. $ f(b)$). Then,  in the case of a nondecreasing
function $f$, the scheme \eqref{h1D}, \eqref{nf} reduces to

\be
 h_i \frac {u_i^{n+1} - u_i^n} k + f(u_i^n) -
f(u_{i-1}^n) = 0, 
\label{up1D} 
\ee

i.e. the upstream\index{Upstream} (or upwind) finite volume scheme. The scheme (\ref{up1D}) is sometimes
called ``upstream finite difference" scheme. In that particular case ($f$ 
monotone and 1D) it is possible to find points $x_i$ in order to obtain a consistent
scheme  in the finite difference sense 
(if $f$ is nondecreasing, take $x_i=x_{{i+1/2}}$ as for
the scheme \refep{VFDA}).
\end{remark}

\index{Scheme!Monotone flux --|)}

\subsection{$L^\infty$-stability for monotone flux 
schemes}\index{Stability|(}
\label{linfstah1d}

\begin{lemma}[$L^\infty$ estimate in the nonlinear case] \label{nlh1dstab}
\index{Estimate on the approximate solution!-- for a hyperbolic equation}
Under Assumption \ref{datah1d}, let $\mesh$ be an admissible
mesh in the sense of definition \refp{meshhyp1d} and let $k\in\Rpast$ be 
the time step. 

Let $\uapp$ be the finite volume approximate solution defined by
\eqref{h1D}-\eqref{nluapp} and assume that the scheme is a monotone flux 
scheme in the sense of definition \refp{monflux}. Let $g_1$ and $g_2$ be 
the Lipschitz constants of $g$ on $[U_m,U_M]^2$ with respect to its two 
arguments.

Under the Courant-Friedrichs-Levy (CFL) condition\index{CFL condition} 

\be
k \le \frac {\inf_{i \in \Z} h_i} {g_1 + g_2},
\label{cflhypnl}
\ee

(note that taking $k \le \a h/ (g_1 + g_2)$ implies \eqref{cflhypnl}), 

the approximate solution  $\uapp$ satisfies

\[
U_m \le  \uapp(x,t) \le U_M \, \text{ for a.e. } x\in\R \text{ and a.e. } 
t \in \Rp.
\]

\end{lemma}

\proof of Lemma \ref{nlh1dstab}

\smallskip

Let us prove that 
\be
U_m \le  u_i^n \le U_M, \, \forall i \in \Z, \, \forall n \in \N,
\label{induc}  
\ee
by induction on $n$, which proves the lemma.
Assertion \eqref{induc} holds for $n=0$ thanks to the definition of $u_i^0$  in \refep{h1D}. 
Suppose that  it holds  for $n\in\N$.

For all $i\in\Z$,  scheme \eqref{h1D}, \eqref{nf} (with $p=1$) gives
\begin{equation}
u_i^{n+1} =  (1 - b_{\ipud}^n - a_{\imud}^n) u_i^n  +  b_{\ipud}^n u_{i+1}^n + a_{\imud}^n u_{i-1}^n,
\label{convex-comb}
\end{equation}
with
\[
b_{\ipud}^n  = 
\ba \dsp{ \frac k {h_i} \frac  {g(u_i^n,u_{i+1}^n) - f(u_i^n)} {u_i^n - u_{i+1}^n}}  \text{ if } u_i^n \ne  u_{i+1}^n,
\\
0 \text{ if }  u_i^n  = u_{i+1}^n,
\ea
\]
and

\[ 
a_{\imud}^n  = 
\ba \dsp{ \frac k {h_i} \frac 
{g(u_{i-1}^n,u_i^n) - f(u_i^n)} {u_{i-1}^n - u_i^n} 
}  \text{ if } u_i^n  \ne  u_{i-1}^n,
\\
0 \text{ if }  u_i^n  = u_{i-1}^n.
\ea
\]

Since $f(u_i^n)=g(u_i^n,u_i^n)$ and thanks to the monotonicity of $g$, $0 
\le  b_{\ipud}^n \le g_2 k/h_i$ and $0 \le  a_{\imud}^n \le g_1 k/h_i$, for all $i\in\Z$.
Therefore, under condition \eqref{cflhypnl}, the value $u_i^{n+1}$
may be written as  a convex linear combination of the values
$u_i^n$ and $u_{i-1}^n$. Assertion  \eqref{induc}
is thus proved for $n+1$, which concludes the proof of the lemma.
\eof
\index{Stability|)}
\subsection{Discrete entropy inequalities}\label{discenth1d}
\index{Discrete!-- entropy inequalities}\index{Entropy!Discrete -- 
inequalities}

\begin{lemma}[Discrete entropy inequalities] \label{nlh1dentdis}
Under Assumption \ref{datah1d}, let $\mesh$ be an admissible
mesh in the sense of definition \refp{meshhyp1d} and let $k\in\Rpast$ be 
the time step. 

Let $\uapp$ be the finite volume approximate solution defined by
\eqref{h1D}-\eqref{nluapp} and assume that the scheme is a monotone flux 
scheme in the sense of definition \refp{monflux}. Let $g_1$ and $g_2$ be 
the Lipschitz constants of $g$ on $[U_m,U_M]^2$ with respect to its two 
arguments.
Under the CFL condition \eqref{cflhypnl}, the following inequation holds:

\be\ban
\dsp \frac {h_i} k  \Bigl( \vert u_i^{n+1} - \sk\vert - 
\vert u_i^n - \sk\vert \Bigr) +
\\
\dsp  
g(u_i^n\top\sk,u_{i+1}^n\top\sk)
- g(u_i^n\bot\sk,u_{i+1}^n\bot\sk) -   
g(u_{i-1}^n\top\sk,u_i^n\top\sk)
+ g(u_{i-1}^n\bot\sk,u_i^n\bot\sk)
\le 0,
\\
\hfill \forall \ n \in \N,\ \forall \ i \in \Z,\ \forall \ \sk \in \R.
\ean\label{eqnlh1de}\ee

Recall that $a\top b$ (resp. $a\bot b$) denotes the maximum (resp. the 
minimum) of the  two real numbers $a$ and $b$.
\end{lemma}

\proof of Lemma \ref{nlh1dentdis}

\smallskip

Thanks to the monotonicity properties of $g$ and to the condition 
\eqref{cflhypnl} (see remark \ref{schmono}),

\[
u_i^{n+1} = H(u_{i-1}^n,u_i^n,u_{i+1}^n), \, \forall i \in \Z, \, \forall 
\nnn,
\]

where $H$ is a function from $\R^3$ to $\R$ which is nondecreasing with respect to 
all its arguments and such that $\sk = H(\sk,\sk,\sk)$ for all $\sk \in \R$.

Hence, for all $\sk \in \R$,

\[
u_i^{n+1} \le  H(u_{i-1}^n\top\sk,u_i^n\top\sk,u_{i+1}^n\top\sk),
\]

and

\[
\sk \le H(u_{i-1}^n\top\sk,u_i^n\top\sk,u_{i+1}^n\top\sk),
\]

which yields

\[
u_i^{n+1}\top\sk \le 
H(u_{i-1}^n\top\sk,u_i^n\top\sk,u_{i+1}^n\top\sk).
\]

In the same manner, we get

\[
u_i^{n+1}\bot\sk \ge 
H(u_{i-1}^n\bot\sk,u_i^n\bot\sk,u_{i+1}^n\bot\sk),
\]

and therefore, by substracting the last two  equations, 

\[
\vert u_i^{n+1} - \sk \vert \le 
H(u_{i-1}^n\top\sk,u_i^n\top\sk,u_{i+1}^n\top\sk) - 
H(u_{i-1}^n\bot\sk,u_i^n\bot\sk,u_{i+1}^n\bot\sk),
\]

that is \eqref{eqnlh1de}.
\eof

In the two next sections, we study the convergence of the schemes
defined by \eqref{h1D}, \eqref{nf} with $p=1$ (see the remarks 
\ref{5pts1} and \ref{5pts2} and Section \ref{hyperhos} for the schemes 
with $2p+1$ points). 

We first develop a proof of convergence for the 
monotone flux schemes; this proof is based 
on a weak $BV$ estimate  similar to \eqref{bvf} like  the proof of 
proposition \refp{convhyplin} in the linear case. It will be generalized 
in the multidimensional case studied in Chapter \ref{hypmd}.
We then  briefly  describe the $BV$ framework which gave the first
convergence results; its generalization to the 
multidimensional case is not so easy, except in the case of Cartesian  meshes.

\subsection{Convergence of the upstream\index{Upstream} scheme in the general case}

A proof of convergence similar to the proof of convergence given in 
the linear case can be developed. For the sake of simplicity, we shall consider only 
the case of a nondecreasing function $f$ and of the classical 
upstream\index{Upstream} scheme 
(the general case for $f$ and for the monotone flux schemes being handled
in Chapter \ref{hypmd}). We shall first prove a ``weak $BV$" 
estimate.

\index{BV!Weak -- estimate|(}
\begin{lemma}[Weak $BV$ estimate for the nonlinear case] 
\label{nlwb1d}
Under Assumption \ref{datah1d}, assume that   $f$ is nondecreasing. Let $\xi\in(0,1)$ 
be a given value. Let $\mesh$ be an admissible
mesh in the sense of definition \refp{meshhyp1d}, let $M$
be the Lipschitz constant of $f$ in $[U_m,U_M]$ and let $k\in\Rpast$
satisfying the CFL condition 

\be
k \le (1 - \xi) \frac {\inf_{i \in \Z} h_i} M.
\label{cflhypnlc}\ee

(The condition $k\le(1 - \xi)\a h/ M$ implies the above condition.)
Let $\{u_i^n$, $i\in\Z$, $n\in\N\}$ be given by the finite volume
scheme \eqref{h1D}, \eqref{nf} with $p=1$ and $g(a,b)=f(a)$.
Let $R\in\Rpast$ and $T\in\Rpast$ and assume $h < R$ and $k < T$.
Let $i_0 \in \Z$, $i_1 \in \Z$ and $N \in \N$ be such that $-R \in 
\cvbar_{i_0}$,  $R \in \cvbar_{i_1}$,and $T \in (Nk,(N+1)k]$.
Then there exists $C\in\Rpast$, only depending on
$R$, $T$, $u_0$, $\a$, $f$ and $\xi$, such that

\be
\sum_{i=i_0}^{i_1} \sumN k 
\vert f(u_i^n) - f(u_{i-1}^n) \vert \le C h^{-1/2}.
\label{bvfnl}\ee

\end{lemma}

\proof of Lemma \ref{nlwb1d}

\smallskip

We multiply the first equation of \eqref{h1D} by $k u_i^n$, and we sum on
$i=i_0,\ldots,i_1$ and $n=0,\ldots,N$. We get $A+B = 0$,
with

\[
A = \sum_{i=i_0}^{i_1} \sumN h_i (u_i^{n+1} - u_i^n) u_i^n,
\]

and 

\[
B = \sum_{i=i_0}^{i_1} \sumN k \Bigl(f(u_i^n) - f(u_{i-1}^n)\Bigr) u_i^n.
\]

We have

\[
A = - \half \sum_{i=i_0}^{i_1} \sumN h_i (u_i^{n+1} - u_i^n)^2
+ \half \sum_{i=i_0}^{i_1} h_i [(u_i^{N+1})^2 - (u_i^0)^2].
\]

Using the scheme \eqref{h1D}, we get

\[
A = - \half \sum_{i=i_0}^{i_1} \sumN \frac{k^2} {h_i} 
\Bigl(f(u_i^n) - f(u_{i-1}^n)\Bigr)^2
+ \half \sum_{i=i_0}^{i_1} h_i [(u_i^{N+1})^2 - (u_i^0)^2],
\]

and therefore, using the CFL condition \eqref{cflhypnlc},

\be
A \ge - \frac 1 {2M} (1-\xi) \sum_{i=i_0}^{i_1} \sumN k 
\Bigl(f(u_i^n) - f(u_{i-1}^n)\Bigr)^2
- \half \sum_{i=i_0}^{i_1} h_i (u_i^0)^2.
\label{eqbvwhnl}\ee

We now study the term $B$. 

Denoting by $\Phi$ the function 
$\Phi(a) = \int_{U_m}^{a} s f'(s) ds$, for all $a\in\R$,
an integration by parts yields, 
for all $(a,b)\in\R^2$,

\[
\Phi(b) - \Phi(a) = b(f(b)-f(a)) - \int_a^b(f(s)-f(a))dx.
\]

Using the technical lemma \refp{petitlemme} which states
$\int_a^b(f(s)-f(a))dx \ge \frac 1 {2M} (f(b)-f(a))^2$, we obtain

\[
b(f(b)-f(a)) \ge \frac 1 {2M} (f(b)-f(a))^2 + \Phi(b) - \Phi(a).
\]

The above inequality with $a=u_{i-1}^n$ and $b=u_i^n$
yields

\[
B \ge \frac 1 {2M} \sum_{i=i_0}^{i_1} \sumN k 
\Bigl(f(u_i^n) - f(u_{i-1}^n)\Bigr)^2
+ \sumN k [\Phi(u_{i_1}^n) - \Phi(u_{i_0-1}^n)].
\]

Thanks to the $L^\infty$ estimate of Lemma \refp{linstab}, there exists
$C_1>0$, only depending on $u_0$ and $f$  such that 

\[
B \ge \frac 1 {2M} \sum_{i=i_0}^{i_1} \sumN k 
\Bigl(f(u_i^n) - f(u_{i-1}^n)\Bigr)^2 -  T  C_1.
\]

Therefore, since $A+B=0$ and $\sum_{i=i_0}^{i_1}h_i \le 4 R$, the following inequality holds:

\[
0 \ge \xi \sum_{i=i_0}^{i_1} \sumN k 
\Bigl(f(u_i^n) - f(u_{i-1}^n)\Bigr)^2
- 4 R M \max\{-U_m,U_M\}^2 - 2 M T  C_1,
\]

which gives the existence of $C_2\in\Rpast$, only depending 
on $R$, $T$, $u_0$, $f$ and $\xi$ such that

\[
\sum_{i=i_0}^{i_1} \sumN k \Bigl(f(u_i^n) - f(u_{i-1}^n)\Bigr)^2 
\le C_2.
\]

The Cauchy-Schwarz inequality yields

\[
\Bigl[\sum_{i=i_0}^{i_1} \sumN k 
\vert f(u_i^n) - f(u_{i-1}^n) \vert \Bigr]^2 \le 
C_2 2 T \frac {4 R} {\a h},
\]

which concludes the proof of the lemma.
\eof
\index{BV!Weak -- estimate|)}

We can now state the convergence theorem.

\begin{theo}[Convergence in the nonlinear case] 
\index{Convergence of the finite volume method!-- for a nonlinear 
hyperbolic equation}\label{nlcv1d}
Assume Assumption \ref{datah1d} and $f$ nondecreasing. Let $\xi\in(0,1)$ 
and $\a >0$ be given. Let $M$ be the Lipschitz constant of $f$ in 
$[U_m,U_M]$. For an admissible mesh $\mesh$ in the sense of Definition 
\refp{meshhyp1d} and for a time step  $k\in\Rpast$
satisfying the CFL condition \eqref{cflhypnlc} (taking 
$k\le(1 - \xi)\a h/ M$ is a sufficient condition, note that $\xi$ and $\a$ do not depend of 
$\mesh$), let $\uapp$ be the finite volume approximate solution defined by
\eqref{h1D}-\eqref{nluapp} with $p=1$ and $g(a,b)=f(a)$.

Then the function $\uapp$ converges to the unique entropy
weak solution  $u$  of \refep{hyperbolic1D} in $L^1_{loc}(\R\times\Rp)$
as $\sizemesh$ tends to $0$.
\end{theo}

\proof

\smallskip

Let $Y$ be the set of approximate solutions, that is the set of $\uapp$, defined by
\eqref{h1D}-\eqref{nluapp} with $p=1$ and $g(a,b)=f(a)$, for all 
$(\mesh,k)$ where $\mesh$ is an admissible 
mesh in the sense of Definition 
\refp{meshhyp1d} and $k\in\Rpast$
satisfies the CFL condition \eqref{cflhypnlc}.
Thanks to Lemma \ref{nlh1dstab}, the set $Y$ is bounded in $L^\infty(\R 
\times \Rp)$.

\smallskip

The proof of Theorem \ref{nlcv1d} is performed in three steps. 
In the first step, a compactness 
result is given for $Y$, 
only using the boundeness 
of $Y$ in $L^\infty(\R \times \Rp)$. In the second 
step, it is proved that the eventual limit (in a convenient sense) of a sequence of approximate 
solutions is a solution (in a convenient sense) of problem 
\eqref{hyperbolic1D}. In the third step a uniqueness result yields the 
conclusion. For steps 1 and 3, we refer to chapter \ref{hypmd} for a 
complete proof.

\smallskip

\textit{Step 1 (compactness result)}

Let us first use a compactness result in 
 $L^\infty(\R \times \Rp)$ which 
is stated in Proposition \refp{nlwsprop}.
Since $Y$ is bounded in 
$L^\infty(\R \times \Rp)$, for any sequence $(u_m)_{\mnn}$ of $Y$ there 
exists a subsequence, still denoted by $(u_m)_{\mnn}$, and there exists 
$\mu \in L^\infty(\R \times \Rp \times (0,1))$ such that $(u_m)_{\mnn}$ 
converges to $\mu$ in the ``nonlinear weak-$\star$ sense", that is

\[
\int_{\R}\int_{\Rp}\theta(u_m(x,t))\varphi(x,t)dtdx \tends 
\int_{\R}\int_{\Rp}\int_0^1\theta(\mu(x,t,\a)) \vphi(x,t) d\a dt dx, \text{ as } m \tends 
\infty,
\]

for all $\varphi \in L^1(\R \times \Rp)$ and all $\theta \in C(\R,\R)$.
In other words, for any $\theta \in C(\R,\R)$,

\be
\theta(u_m) \tends \mu_\theta \text{ in } L^\infty(\R \times \Rp) \text{ for the 
weak-$\star$ topology as } m \tends \infty,
\label{cvmug}
\ee

where $\mu_\theta$ is defined by

\[
\mu_\theta(x,t)=\int_0^1\theta(\mu(x,t,\a)) d\a, \text{ for a.e. } (x,t) \in \R 
\times \Rp.
\]

\smallskip

\textit{Step 2 (passage to the limit)}

Let $(u_m)_{\mnn}$ be a sequence of $Y$. Assume that $(u_m)_{\mnn}$ converges to 
$\mu$ in the nonlinear weak-$\star$ sense and that $u_m=u_{\mesh_m,k_m}$ (for all 
$\mnn$) with $\size(\mesh_m) \tends 0$ as $m \tends \infty$ (note that 
$k_m \tends 0$ as $m \tends \infty$, thanks to \eqref{cflhypnlc}).

Let us prove that $\mu$ is a ``solution" to problem \eqref{hyperbolic1D} 
in the following sense (we shall
say that $\mu$ is ``an entropy process solution"\index{Entropy!-- process 
solution} to problem \eqref{hyperbolic1D}):

\be\left\{
\bal
\mu \in L^\infty(\R \times \Rp \times (0,1)),
\\
\dsp \dint \int_0^1 \Bigl( \vert \mu(x,t,\a) - \sk \vert \varphi_t(x,t) + 
( f(\mu(x,t,\a)\top\sk) - f (\mu(x,t,\a)\bot\sk) ) 
\varphi_x(x,t) \Bigr) d\a dt dx
\\
\hfill \dsp + \int_\R
\vert u_0(x) - \sk \vert \varphi(x,0) dx \ge 0, \, \forall \varphi \in C^1_c(\R \times \R_+,\R_+), \forall \sk \in \R. 
\ean
\right.
\label{h1deps}
\ee

Let $\sk \in \R$. Setting 

\[
v(x,t)=\int_0^1 \vert \mu(x,t,\a) - \sk \vert d\a , \text{ for a.e. } 
(x,t) \in \R \times \Rp
\]

and

\[
w(x,t)=\int_0^1 \Bigl( f(\mu(x,t,\a)\top\sk) - f(\mu(x,t,\a)\bot\sk) \Bigr) d\a, \text{ for a.e. } 
(x,t) \in \R \times \Rp,
\]

the inequality in \eqref{h1deps} reads

\be\bari
\dsp \dint [ v(x,t) \varphi_t(x,t) + 
 w(x,t) \varphi_x(x,t) ] dtdx + \int_\R
\vert u_0(x) - \sk \vert\varphi(x,0) dx \ge 0,
\\
\forall \varphi \in C^1_c(\R \times \R_+,\R_+).
\ean\label{cvh1dnl1}
\ee

Let us  prove that \eqref{cvh1dnl1} holds; for  $m\in\N$  we shall denote by 
$\mesh = \mesh_m$ and $k= k_m$. We use
the result of Lemma \ref{nlh1dentdis}, which reads in the present
particular case $f'\ge 0$, 

\[
h_i \frac {v_i^{n+1} - v_i^n} k + w_i^n - w_{i-1}^n \le 0, \, \forall 
i\in\Z, \, \forall n\in\N,
\]

where  $v_i^n = \vert u_i^n - \sk \vert$ and
$w_i^n =$ $f(u_i^n\top\sk) - f(u_i^n\bot\sk) =$ 
$\vert f(u_i^n) - f(\sk) \vert$. 

The functions $v_{\mesh_m,k_m}$ and $w_{\mesh_m,k_m}$ are
defined in the same way as the function $u_{\mesh_m,k_m}$, i. e. with
constant values $v_i^n$ and $w_i^n$ 
in each control volume $\cv_i$ during each time step $(nk,(n+1)k)$. 
Choosing $\theta$ equal to the continuous functions $\vert \cdot - \sk\vert$ and
$\vert f(\cdot) - f(\sk)\vert$ in \eqref{cvmug} yields that the sequences
$(v_{\mesh_m,k_m})_{m\in\N}$ and $(w_{\mesh_m,k_m})_{m\in\N}$ 
converge to $v$ and $w$ in $L^\infty(\R\times\Rpast)$ for the weak-$\star$ topology.

Applying the method
which was used  in the proof of Proposition \refp{convhyplin},
taking $v_i^l$ instead of $u_i^l$ in the
definition of $A_m$ (for $l=n$ and $n+1$) and $w_j^n$ instead of $u_j^n$ in the
definition of $B_m$ (for $j=i$ and $i-1$), we conclude that \eqref{cvh1dnl1} holds.

Indeed, a weak $BV$ inequality holds on the values $w_i^n$ (that is 
\refep{bvf} holds with $w^n_j$ instead of $u^n_j$ for $j=i$ and $i-1$), thanks to 
Lemma \refp{nlwb1d} and the relation 

\[ 
\left\vert \vert f(u_i^n) - \sk\vert - 
\vert f(u_{i-1}^n) - \sk\vert \right\vert
\le \vert f(u_i^n) - f(u_{i-1}^n) \vert, \, \forall i\in\Z, \, \forall n\in\N.
\]

(Note that here, as in the linear case, the useful consequence of the weak $BV$ 
inequality, is \refep{bvweff} with $w^n_j$ instead of $u^n_j$ for $j=i$ 
and $i-1$.)

This concludes Step 2.

\smallskip

\textit{Step 3 (uniqueness result for \eqref{h1deps} and conclusion)}

Theorem \refp{ewsunicite} states that there exists at most one 
solution to \eqref{h1deps} and that there exists $u \in L^\infty(\R \times 
\Rp)$ such that $\mu$ solution to \eqref{h1deps} implies 
$\mu(x,t,\a)=u(x,t)$ for a.e. $(x,t,\a) \in \R \times  \Rp \times (0,1)$. 
Then, $u$ is necessarily the entropy weak solution to  \eqref{hyperbolic1D}.

Furthermore, if $(u_m)_{\mnn}$ converges to $u$ in 
the nonlinear weak-$\star$ sense, an easy argument shows that $(u_m)_{\mnn}$ 
converges to $u$ in $L^1_{loc}(\R \times \Rp)$ (and even in $L^p_{loc}(\R 
\times \Rp)$ for all $1 \le p < \infty$), see Remark \refep{cvforte}.

Then, the conclusion of Theorem \ref{nlcv1d} follows easily from Step 2 
and Step 1 by way of contradiction (in order to prove the convergence of 
a sequence $u_{\mesh_m,k_m} \subset Y$ to $u$, if $\size(\mesh_m) \tends 
0$ as $m \tends \infty$, without any extraction of 
a ``subsequence").
\eof

\begin{remark} \normalfont\label{5pts1}
In Theorem \ref{nlcv1d}, we only consider the case $f'\ge 0$ and the 
so called ``upstream\index{Upstream} scheme". It is quite easy to generalize the result 
for any $f \in C^1(\R,\R)$ and any monotone flux scheme (see the 
following chapter). It is also 
possible to consider other schemes (for instance, some 5-points schemes, 
as in Section \ref{hyperhos}). For a given scheme, the proof of 
convergence of the approximate solution towards the entropy weak solution 
contains 2 steps:

\benum
\item prove an $L^\infty$ estimate on the approximate solutions, which 
allows to use the compactness result of Step 1 of the proof of Theorem 
\ref{nlcv1d},
\item prove a ``weak $BV$" estimate and some ``discrete entropy 
inequality" in order to have the following property:

If $(u_m)_{\mnn}$ is a sequence of approximate solutions which converges 
in the nonlinear weak-$\star$ sense, then

\[
\bari
\dsp{\lim_{m \tends \N} \dint \Bigl( \vert u_m(x,t) - \sk \vert \varphi_t(x,t) + 
(f(u_m(x,t)\top\sk) - f(u_m(x,t)\bot\sk)) 
\varphi_x(x,t)\Bigr) dt dx}
\\
\dsp{+ \int_\R
\vert u_0(x)- \sk \vert \varphi(x,0) dx \ge 0, \, \forall \varphi \in 
C^1_c(\R \times \R_+,\R_+), \, \forall \sk \in \R.}
\ean
\] 
\eenum
\end{remark}

\begin{remark} \normalfont An interesting consequence of the proof of the uniqueness theorem (Theorem \refp{ewsunicite}) is the principle of ``finite speed of propagation'' for the entropy weak solution  of  \refep{hyperbolic1D} and
the ``continuity in time'' of this entropy weak solution, see Proposition \ref{prop-pvf} and remark \ref{rem-pvf}.
\end{remark}


\subsection{Convergence proof using $BV$}\label{bvmethod}\index{BV!Strong -- 
estimate}

We   now give the details of the classical proof of convergence (considering only $3$ 
points schemes), which requires regularizations of $u_0$ in $BV(\R)$. It consists
 in using
Helly's compactness theorem 
(which may also be used in the linear case to 
obtain a strong convergence of $\uapp$ to $u$ in $L^1_{loc}(\R \times \Rp)$).
This theorem is a direct consequence of Kolmogorov's 
theorem
(theorem \refp{Kolm}). We give below the definition of $BV(\O)$ where 
$\O$ is an open subset of $\R^p(\O)$, $p \ge 1$ (already given in 
Definition \refp{defBV} for $\O=\R$) and we give a straightforward consequence of
Helly's theorem for the case of interest here.

\begin{defi}[$BV(\O)$]\index{BV!-- Space} \label{bvrd}  \normalfont
Let $p\in\N^\star$ and let $\O$ be an open subset of $\R^p$.
A function $v \in L^1_{loc}(\O)$ has a bounded variation, that is $v \in 
BV(\O)$, if $\vert v \vert_{BV(\O)} < \infty$ where

\be
\vert v \vert_{BV(\O)} = \sup \{\int_\O v(x) \div \varphi(x) dx, 
\varphi \in C^1_c(\O,\R^p), \vert \varphi(x) \vert \le 1, \, \forall x 
\in \O \}.  
\label{defBVRp}
\ee
\end{defi}

\begin{lemma}[Consequence of Helly's theorem]\label{helly}
\index{Helly's theorem}
\index{Compactness!Helly's -- theorem}
Let $\mathcal{A } \subset L^\infty(\R^2)$. Assume that there exists $C \in \Rp$
and, for all $T>0$, there exists $C_T \in \Rp$ such that

\[
\Vert v \Vert_{L^\infty(\R^2)} \le C, \, \forall v \in \mathcal{A },
\]

and 
\[
\vert v \vert_{BV(\R\times(-T,T))} \le C_T, \, \forall v \in \mathcal{A }, \, 
\forall T > 0.
\]

Then  for any sequence
$(v_n)_{\nnn}$  of elements of $\mathcal{A}$, 
there exists a subsequence, still denoted by $(v_{n})_{n \in \N}$, and there 
exists $v \in L^\infty(\R^2)$, with 
$\Vert v \Vert_{L^\infty(\R^2)} \le C$
and $\vert v \vert_{BV(\R\times(-T,T))} \le C_T$ for all $T>0$, such that 
$v_{n} \tends v$ in  $L^1_{loc}(\R^2)$ as $n\tends\infty$, 
that is $\int_\comp \vert v_{n}(x) - v(x) \vert dx \tends 0$,
as $n\tends\infty$ for any compact set $\comp$ of $\R^2$.

\end{lemma}
 
In order to use Lemma \ref{helly},  one first shows the following 
$BV$  stability\index{Stability} estimate for the approximate solution:

\begin{lemma}[Discrete space $BV$ estimate]\index{BV!Strong -- estimate}
\index{Estimate on the approximate solution!-- for a hyperbolic equation}
\label{nlstabv}
Under Assumption \ref{datah1d}, assume that $u_0 \in BV(\R)$; let $\mesh$ be an admissible
mesh in the sense of Definition \refp{meshhyp1d} and let $k\in\Rpast$ be 
the time step.
Let $\{u_i^n$, $i\in\Z$, $n\in\N\}$ be given by \eqref{h1D}, 
\eqref{nf} and assume that the scheme is a monotone flux 
scheme<\index{Scheme!Monotone flux --} in the 
sense of Definition \refp{monflux}. Let $g_1$ and $g_2$ be the Lipschitz 
constants of $g$ on $[U_m,U_M]^2$ with respect to its two arguments.
Then, under the CFL condition \eqref{cflhypnl}, the following inequality holds:

\be
\sum_{i\in\Z} \vert u_{i+1}^{n+1} - u_i^{n+1} \vert
\le \sum_{i\in\Z} \vert u_{i+1}^{n} - u_i^{n} \vert, \, 
\forall \ n \in \N.
\label{eqnlstabv}\ee

\end{lemma}

\proof of Lemma \ref{nlstabv}

\smallskip

First remark that, for $n=0$, 
$\sum_{i\in\Z} \vert u_{i+1}^0 - u_i^0 \vert \le 
\vert u_0 \vert_{BV(\R)}$ (see Remark \refp{divbv1d}).

For all $i\in\Z$, the  scheme \eqref{h1D}, \eqref{nf} (with $p=1$) leads to

\[
u_i^{n+1}  = u_i^n + b_{\ipud}^n (u_{i+1}^n - u_i^n) +
a_{\imud}^n (u_{i-1}^n - u_i^n),
\]

and

\[
u_{i+1}^{n+1}  = u_{i+1}^n + b_{\iptd}^n (u_{i+2}^n - u_{i+1}^n) +
a_{\ipud}^n (u_i^n - u_{i+1}^n),
\]

where $a_{{i+1/2}}$ and $b_{{i+1/2}}$ are defined (for all $i \in \Z$) in 
Lemma \refp{nlh1dstab}. Substracting one equation to the other leads to

\[
u_{i+1}^{n+1}-u_i^{n+1}= (u_{i+1}^{n}-u_i^{n})(1 
- b_{\ipud}^n - a_{\ipud}^n) 
+ b_{\iptd}^n (u_{i+2}^n - u_{i+1}^n) + a_{\imud}^n (u_i^n - u_{i-1}^n).
\]

Under the condition \eqref{cflhypnl}, we get

\[
\vert u_{i+1}^{n+1}-u_i^{n+1} \vert \le \vert u_{i+1}^{n}-u_i^{n} \vert (1 
- b_{\ipud}^n - a_{\ipud}^n) 
+ b_{\iptd}^n \vert u_{i+2}^n - u_{i+1}^n \vert 
+ a_{\imud}^n \vert u_i^n - u_{i-1}^n \vert.
\]

Summing the previous equation over $i\in\Z$ gives \eqref{eqnlstabv}.
\eof

\begin{corol}[Discrete $BV$ estimate]\index{BV!Strong -- estimate} \label{nlstabvco}
\index{Estimate on the approximate solution!-- for a hyperbolic equation}
Under assumption \ref{datah1d}, let $u_0 \in BV(\R)$; let $\mesh$ be an admissible
mesh in the sense of Definition \refp{meshhyp1d} and let $k\in\Rpast$
be the time step.
Let $\uapp$ be the finite volume approximate solution defined by
\eqref{h1D}-\eqref{nluapp} and assume that the scheme is a monotone flux scheme in the 
sense of Definition \refp{monflux}. Let $g_1$ and $g_2$ be the Lipschitz 
constants of $g$ on $[U_m,U_M]^2$ with respect to its two arguments and 
assume that $k$ satisfies the CFL condition \eqref{cflhypnl}. Let
$\uapp(x,t)=u^0_i$ for a.e. $(x,t) \in \cv_i \times \R_{-}$, for all $i \in 
\Z$ (hence $\uapp$ is defined a.e. on $\R^2$).
Then, for any $T>0$, there exists $C\in\Rpast$, only depending  on
$u_0$, $g$ and $T$ such that:

\be
\vert \uapp \vert_{BV(\R\times(-T,T))} \le C.
\label{conlstabv}\ee
\end{corol}

\proof of Corollary \ref{nlstabvco}

\smallskip
 
As in Lemma \ref{nlstabv}, remark that $\sum_{i \in \Z}\vert u^0_{i+1}-u^0_i \vert \le 
\snorm{u_0}{BV(\R)}$.

\smallskip

Let us first assume that $T \le k$. Then, the $BV$ semi-norm of $\uapp$ 
satisfies
 
\[
\vert \uapp \vert_{BV(\R\times(-T,T))} \le 2T\sum _{i 
\in \Z}\vert u^0_{i+1}-u^0_i \vert.
\]

Hence the estimate  \eqref{conlstabv} is true 
for $C = 2T \snorm{u_0}{BV(\R)}$.

\smallskip

Let us now assume that $k<T$. Let $N \in \N^\star$ such that $Nk < T \le (N+1)k$.
The definition of $\vert \cdot \vert_{BV(\R\times(-T,T))}$ yields

\be\bal
\vert \uapp \vert_{BV(\R\times(-T,T))} \le T\sum _{i \in \Z}\vert 
u^0_{i+1}-u^0_i \vert +
\\
\dsp{\sum_{n=0}^{N-1}\sum_{i \in \Z} k 
\vert u_{i+1}^n-u_i^n \vert + (T-Nk)\sum_{i \in \Z} 
\vert u_{i+1}^N-u_i^N \vert +\sum_{n=0}^{N-1}\sum_{i \in \Z} h_i 
\vert u_{i}^{n+1}-u_i^n \vert.}
\ean
\label{bvtot}
\ee

Lemma \ref{nlstabv} gives $\sum_{i \in \Z} \vert u_{i+1}^n-u_i^n \vert 
\le \snorm{u_0}{BV(\R)}$ for all $\nnn$, and therefore,  

\be
\sum_{n=0}^{N-1}\sum_{i \in \Z} k 
\vert u_{i+1}^n-u_i^n \vert + (T-Nk)\sum_{i \in \Z} 
\vert u_{i+1}^N-u_i^N \vert \le T \snorm{u_0}{BV(\R)}.
\label{bvtot1}
\ee

In order to bound the last term of \eqref{bvtot}, using the scheme 
\eqref{h1D} yields, for all $i \in \Z$ and all $\nnn$,

\[
\vert u_i^{n+1}-u_i^n \vert \le \frac k {h_i} g_1 \vert u_i^n-u_{i-1}^n \vert 
+ \frac k {h_i} g_2 \vert u_{i}^n-u_{i+1}^n \vert.
\]

Therefore, 

\[
\sum_{i \in \Z} h_i \vert u_i^{n+1}-u_i^n \vert \le k (g_1+g_2) \sum_{i 
\in \Z} \vert u_{i}^n-u_{i+1}^n \vert, \text{ for all }\nnn,
\]

which yields, since $Nk < T$,

\be
\sum_{n=0}^{N-1}\sum_{i \in \Z} h_i 
\vert u_{i}^{n+1}-u_i^n \vert \le T(g_1+g_2)\snorm{u_0}{BV(\R)}.
\label{bvtot2}
\ee

Therefore Inequality \eqref{conlstabv} follows from \eqref{bvtot}, 
\eqref{bvtot1} and \eqref{bvtot2} with $C=T(2+g_1+g_2)\snorm{u_0}{BV(\R)}$.\eof

Consider a sequence of admissible meshes and time steps
verifying the CFL condition, and the associated sequence of approximate
solutions (prolonged on $\R \times \R_-$ as in Corollary 
\ref{nlstabvco}). By Lemma \refp{nlh1dstab} and Corollary \ref{nlstabvco}, 
the sequence of approximate
solutions satisfies  the hypotheses
of Lemma \refp{helly}. It is therefore  possible to extract a
subsequence which converges in $L^1_{loc}(\R\times\Rp)$
to a function $u\in L^\infty(\R\times\Rpast)$. It must still be shown that
the function $u$ is the unique weak entropy solution of
Problem \eqref{hyperbolic1D}. This may be proven by 
using the discrete entropy inequalities
\eqref{eqnlh1de} and the strong $BV$ estimate \eqref{eqnlstabv} or the 
classical Lax-Wendroff
theorem recalled below.

\begin{theo}[Lax-Wendroff]\label{laxw}\index{Lax-Wendroff theorem} 
Under Assumption \ref{datah1d}, let $\a > 0$ be given  and let $(\mesh_m)_{m 
\in \N}$ be a sequence of admissible meshes in the sense of Definition 
\refp{meshhyp1d} (note that, for all $m \in \N$, the mesh $\mesh_m$ satisfies the 
hypotheses of Definition \ref{meshhyp1d} where $\mesh=\mesh_m$ and $\a$ 
is independent of $m$). Let $(k_m)_{m \in \N}$ be a sequence of 
(positive) time steps. Assume that $\size(\mesh_m) \tends 0$ and $k_m \tends 
0$ as $m \tends \infty$.

\smallskip

For $m \in \N$, setting $\mesh=\mesh_m$ and $k=k_m$, let $u_m=\uapp$ be 
the solution of \eqref{h1D}-\eqref{nluapp} with $p=1$ and some $g$ from 
$\R^2$ to $\R$, only depending  on $f$ and $u_0$, locally Lipschitz 
continuous and such that $g(s,s)=f(s)$ for all $s \in \R$.

\smallskip

Assume that $(u_m)_{m \in \N}$ is bounded in $L^\infty(\R \times \Rp)$ 
and that $u_m \tends u$ a.e. on $\R \times \Rp$.
Then, $u$ is a weak solution to problem \eqref{hyperbolic1D} (that is $u$ 
satisfies 
\eqref{h1dweak}).

\smallskip

Furthermore, assume that for any $\sk \in \R$ there exists some 
locally Lipschitz continuous function $G_{\sk}$ from $\R^2$ to $\R$, 
only depending on $f$, $u_0$ and $\sk$, such that $G_{\sk}(s,s)=f(s \top 
\sk)-f(s \bot \sk)$ for all $s \in \R$ and such that for all $m \in \N$ 

\be
\frac 1 k (\vert u_i^{n+1}- \sk \vert -\vert u_i^{n}- \sk \vert) + \frac 
1 {h_i} (G_{\sk}(u^n_i,u^n_{i+1})-G_{\sk}(u^n_{i-1},u^n_{i})) \le 0, \, \forall i \in 
\Z, \, \forall \nnn,
\label{entlw}
\ee

where $\{u^n_i$, $i \in \Z$, $\nnn\}$ is the solution to 
\eqref{h1D}-\eqref{nf} for $\mesh=\mesh_m$ and $k=k_m$.
Then, $u$ is the entropy weak solution  to Problem \eqref{hyperbolic1D} 
(that is $u$ is the unique solution of \eqref{h1dent}).

\end{theo}

\proof of Theorem \ref{laxw}

\smallskip

Since $(u_m)_{m \in \N}$ is bounded in $L^\infty(\R \times \Rp)$ 
and $u_m \tends u$ a.e. on $\R \times \Rp$, the sequence $(u_m)_{m \in 
\N}$ converges to $u$ in $L^1_{loc}(\R \times \Rp)$. 
This implies in particular (from Kolmogorov's theorem, see Theorem 
\ref{Kolm}) that, for all $R>0$ and all $T>0$,

\[
\sup_{m \in \N} \int_0^{2T}\int_{-2R}^{2R} \vert u_m(x,t)-u_m(x-\eta,t) \vert 
dxdt \tends 0 \text{ as } \eta \tends  0.
\]

Then, taking $\eta = \a \, \size(\mesh_m)$ (for $\mnn$) and letting $m \tends \infty$ 
yields, in particular,

\be
\int_0^{2T}\int_{-2R}^{2R} \vert u_m(x,t)-u_m(x-\a \, \size(\mesh_m),t) \vert dxdt \tends 
0 \text{ as } m \tends \infty.
\label{astuce}
\ee

For $\mnn$, let $\{u^n_i$, $i \in \Z$, $\nnn\}$ be the solution to 
\eqref{h1D}-\eqref{nf} for $\mesh=\mesh_m$ and $k=k_m$ (note that $u^n_i$ 
depends on $m$, even though this dependency is not so clear in the 
notation). 
We also set $k_m=k$ and $\size(\mesh_m)=h$, so that $k$ and $h$ depend on $m$ (but recall that $\alpha$ does not depend on $m$).

\smallskip

Let $R>0$ and $T>0$. Let $i_0 \in \Z$, $i_1 \in \Z$ and $N \in \N$ be such that $-R \in 
\cvbar_{i_0}$, $R \in \cvbar_{i_1}$ and $T \in (Nk,(N+1)k]$.
Then, for $h < R$ and $k < T$ (which is true for $m$ large enough),

\[
\a h\sum_{i=i_0}^{i_1} \sumN k \vert u_i^n - u_{i-1}^n \vert \le 
\int_0^{2T}\int_{-2R}^{2R} \vert u_m(x,t)-u_m(x-\a h,t) \vert dx dt.
\]

Therefore, Inequality \eqref{astuce} leads to \eqref{bvweff}, that is

\be
h\sum_{i=i_0}^{i_1} \sumN k \vert u_i^n - u_{i-1}^n \vert \tends 0 
\text{ as } m \tends \infty.
\label{bvwelw}
\ee

Using \eqref{bvwelw}, the remainder of the proof of Theorem \ref{laxw} is very similar to 
the proof of Proposition \refp{convhyplin} and to Step 2 in the proof of 
Theorem \refp{nlcv1d} (Inequality \eqref{bvwelw} replaces the weak $BV$ 
inequality).

\smallskip

In order to prove that $u$ is solution to \eqref{h1dweak}, let us multiply 
the first equation of \eqref{h1D} by 
$(k/h_i)\varphi(x,nk)$, integrate over $x\in\cv_i$ and sum for
all $i\in\Z$ and $n\in\N$. This yields

\[
A_m + B_m = 0
\]

with

\[
A_m = \sum_{i\in\Z} \sumIN (u_i^{n+1} - u_i^n) \int_{\cv_i} \varphi(x,nk) dx
\]

and

\[
B_m = \sum_{i\in\Z} \sumIN k (g(u^n_i,u^n_{i+1}) - g(u_{i-1}^n,u^n_i))
 \frac{1} {h_i} \int_{\cv_i} \varphi(x,nk) dx.
\]

As in the proof of Proposition \ref{convhyplin}, one has

\[
\lim_{m \tends + \infty} A_m = 
- \int_{\Rp} \int_{\R} u(x,t) \varphi_t(x,t) dx dt - \int_\R 
u_0(x) \varphi(x,0) dx.
\]

Let us now turn  to the study of $B_m$. We compare $B_m$ with

\[
B_{1,m} = - \sumIN \int_{nk}^{(n+1)k} \int_{\R} f(\uapp(x,t)) 
\varphi_x(x,nk) dx dt,
\]

which tends to 
$- \int_{\Rp} \int_{\R} f(u(x,t)) \varphi_x(x,t) dx dt$ as
$m\tends\infty$ since $f(\uapp) \tends f(u)$ in $L^1_{loc}(\R \times 
\Rp)$ as $m\tends\infty$.

The term $B_{1,m}$ can be rewritten as

\[
B_{1,m} =   \sum_{i\in\Z} \sumIN k (f(u_i^n) - f(u_{i-1}^n))  
 \varphi(x_\imudt,n k),
\]

which yields, introducing $g(u_{i-1}^n,u_i^n)$,

\[
\bari
\dsp{B_{1,m} =  \sum_{i\in\Z} \sumIN k (f(u_i^n) - g(u_{i-1}^n,u_i^n)) 
\varphi(x_\imudt,n k)}
\\
\dsp{ + \sum_{i\in\Z} \sumIN k (g(u_{i-1}^n,u_i^n) - 
f(u_{i-1}^n))  \varphi(x_\imudt,n k).}
\ean
\]

Similarly, introducing $f(u_i^n)$ in $B_m$,

\[
\bari
\dsp{B_m =\sum_{i\in\Z} \sumIN k 
(f(u_i^n) - g(u_{i-1}^n,u^n_i)) \frac{1} {h_i} \int_{\cv_i} \varphi(x,nk) 
dx}
\\
\dsp{+ \sum_{i\in\Z} \sumIN k (g(u^n_i,u^n_{i+1}) - f(u_i^n)) 
\frac{1} {h_i} \int_{\cv_i} \varphi(x,nk) dx.}
\ean
\]

In order to compare $B_m$ and $B_{1,m}$, let $R>0$ and $T>0$ be such that $\varphi(x,t)=0$ if $\vert x \vert \ge 
R$ or $t \ge T$. Let $A>0$ be such that $\norm{u_m}{L^\infty(\R \times 
\Rp)} \le A$ for all $\mnn$. Then there exists $C >0$, only depending  on 
$\varphi$ and the Lipschitz constants on $g$ on $[-A,A]^2$, such that, if 
$h < R$ and $k < T$ (which is true for $m$ large enough),

\be
\vert B_m - B_{1,m} \vert \le C h 
\sum_{i=i_0}^{i_1} \sumN k \vert u_i^n - u_{i-1}^n\vert,
\label{bazar}\ee

where $i_0 \in \Z$, $i_1 \in \Z$ and $N \in \N$ are such that $-R \in \cvbar_{i_0}$, $R
\in \cvbar_{i_1}$ and $T \in (Nk,(N+1)k]$. 

Using \eqref{bazar} and \eqref{bvwelw}, we get $\vert B_m - B_{1,m} \vert \tends 
0$ and then

\[
B_m \tends - \int_{\R}\int_{\Rp} f(u(x,t)) \varphi_x(x,t) dt dx \text{ as 
} m\tends\infty,
\]

which completes the proof that $u$ is a solution to problem \eqref{h1dweak}. 

\smallskip

Under the additional assumption that $u_m$ satisfies \eqref{entlw}, one 
proves that $u$ satisfies \refep{h1dentkr} (and therefore that $u$ satisfies 
\eqref{h1dent}) and is the entropy weak solution to Problem 
\eqref{hyperbolic1D} by a similar method.

Indeed, let $\sk \in \R$. One replaces $u^l_i$ by $\vert u^l_i - 
\sk \vert$ in $A_m$ (for $l=n$ and $n+1$) and one replaces $g$ by $G_{\sk}$ in $B_m$. Then, 
passing to the limit in $A_m+B_m \le 0$ (which is a consequence of the 
inequation \eqref{entlw}) leads the desired result.

This concludes the proof of Theorem \ref{laxw} 
\eof

\begin{remark} \normalfont\label{5pts2}
Theorem \ref{laxw} still holds with $(2p+1)$-points schemes ($p > 1$). The 
generalization of the first part of Theorem  \ref{laxw} (the proof that $u$ 
is a solution to \eqref{h1dweak}) is quite easy. For the second part of Theorem  
\ref{laxw} (entropy inequalities) the discrete entropy inequalities 
may 
be replaced by some weaker ones (in order to handle interesting schemes 
such as those which are described in the following section). 

However, the use of Theorem  \ref{laxw} needs a compactness property 
of sequences of approximate solutions \index{Compactness!-- of the approximate
solutions in $L^1$}
 in the space $L^1_{loc}(\R \times 
\Rp)$. Such a compactness property is generally achieved with a ``strong 
$BV$ estimate" (similar to \eqref{eqnlstabv}). Hence an extensive 
literature on ``TVD schemes" (see \citet{HA}), ``ENO schemes"\ldots (see \citet{godlewski-ellipses}, 
\citet{godlewski-springer} and references therein). The generalization of 
this method in the multidimensional case (studied in the following 
chapter) does not seem so clear except in the case of  Cartesian meshes.
\end{remark}

\section{Higher order schemes}\label{hyperhos}\index{Scheme!Higher order --|(}

\hpage{Hyperbolic 1d}

Consider a  monotone flux 
scheme in the sense of Definition \refp{monflux}. By definition, the 
considered scheme  is a $3$ 
points scheme; recall that the numerical flux function is denoted by $g$. The
approximate solution obtained with this scheme converges to the entropy 
weak solution   of Problem
\refep{hyperbolic1D} as the mesh size tends to $0$ and under a so called 
CFL condition (it is proved in Theorem \ref{nlcv1d} for a 
particular case and in the next chapter for the general case). However, 
$3$-points schemes are known to be diffusive, so that the
approximate solution is not very precise near the discontinuities. 
An idea to reduce the diffusion is to go to a $5$-points scheme by introducing ``slopes" on each discretization cell and limiting the slopes in order for the scheme to remain stable.
A classical way to do this is the ``MUSCL" (Monotonic Upwind Scheme for Conservation  Laws, see \citet{vanleer3}) technique\index{Scheme!MUSCL --}
\index{Scheme!Van Leer --}\index{Van Leer scheme}.

Reconstructing a slope on each cell enables to compute interface values on each side of an interface $x_\ipud$.
These values are then used in the computation of the fluxes 
\medskip

We briefly describe, with the notations of Section \ref{ms1d},  an example of such a scheme, see e.g. 
\citet{godlewski-ellipses} and \citet{godlewski-springer} for further 
details. Let $n \in \N$.

\begin{itemize}
\item{Computation of the slopes}

\[ 
\tilde p_i^n = \frac {\uipn - \uimn} {h_i + \frac {h_{i-1}} 2 + \frac 
{h_{i+1}} 2}, i \in \Z.
\] 

\item{Limitation of the slopes}

$ p_i^n = \a_i^n \tilde p_i^n$, $i \in \Z$, where $\a_i^n $ is the largest 
number in 
$[0,1]$ such that

\[
\uin + \frac {h_i} 2 \a_i^n \tilde p_i^n  \in [\uin \bot \uipn,\uin 
\top \uipn] \text{ and } \uin - \frac {h_i} 2 \a_i^n \tilde p_i^n  \in [\uin \bot 
u^n_{i-1},\uin \top u^n_{i-1}].
\]

In practice, other formulas giving smaller values of $\a_i^n$ are 
sometimes needed for stability reasons.

\item{Computation of  $u^{n+1}_i$ for $i \in \Z$}
One replaces $g(u^n_i,u^{n}_{i+1})$ in \eqref{nf} by :
\[ \overline g (\uimn,\uin,\uipn, u_{i+2}^n) = g(\uin + \frac {h_i} 2 
p_i^n, \uipn - \frac {h_{i+1}} 2 p_{i+1}^{n}).
\]
\end{itemize}

The scheme thus constructed is less diffusive than the original one and 
it remains stable thanks to the limitation of the slope. Indeed, if the 
limitation of the slopes is not 
active (that is $\a_i^n = 1$), the space diffusion term disappears from 
this new scheme, while the time ``antidiffusion" term
remains. 
Hence it seems appropriate to use a higher order scheme for the time discretization. 
This may be done by using, for instance, an 
RK2 (Runge Kutta order 2, or Heun) method for the discretization of the
time derivative. The MUSCL scheme may be written as

\[ \frac {U^{n+1} - U^n} k = \overline H (U^n) \text{ for } n \in \N,\]
where $U^n = (u_i^n)_{i \in \Z}$; hence it may be seen as the explicit 
Euler discretization \index{Scheme!Explicit Euler --} of 

\[ U_t = \overline H (U);\]

therefore, the RK2 time discretization yields to the following scheme:
\[ \frac {U^{n+1} - U^n} k = \half \overline H (U^n) + \half 
\Bigl(\overline H(U^n + k\overline H(U^n)) \Bigr) \text{ for } n \in \N.\]

Going to a second order discretization in time allows larger time steps, 
without loss of stability.

\smallskip

Results of convergence are possible with these new schemes (with 
eventually some  adaptation of the slope limitations to obtain convenient discrete 
entropy inequalities, 
see \citet{vilathese}. It is also possible to obtain 
error estimates in the 
spirit of those given in the following chapter, in the multidimensional 
case, see e.g. \citet{Chainais1}, \citet{Noelle}, \citet{kroner2}. However these error 
estimates are somewhat unsatisfactory  since they are of a similar order  to that of 
the original $3$-points scheme (although these schemes are  numerically  
more precise that the original $3$-points schemes).

\medskip

The higher order  schemes are  nonlinear even if  
Problem \refep{hyperbolic1D} is linear, because of the limitation of the slopes.
 
\medskip

Implicit versions  of these higher order schemes are more or less straightforward. 
However, the numerical implementation of these implicit versions requires the solution of nonlinear systems. 
In many cases, the solutions to these nonlinear systems seem impossible to reach for large $k$; in fact, the  existence of the solutions is not so clear, see \citet{pfertzel}. Since the advantage of implicit schemes is essentially the possibility to use large values of $k$, the above flaw considerably reduces the opportunity of their use.  
Therefore, although  implicit 3-points schemes are very diffusive, they remain the basic schemes in several industrial environments; see also Section \refp{turbo} for some clues on implicit schemes  applied to complex industrial  applications.

\index{Scheme!Higher order --|)}

\section{Boundary conditions}\label{hyperbc1d}\index{Boundary conditions for
hyperbolic equations}

A general convergence result is presented here in the case of a scalar equation. Then, this
result will be applied to understand the sense of the boundary condition, described at $x=1$
in the previous section, in a simplified scalar case.

\subsection{A general convergence result}\label{conv}

The unknown is now a function $u~:~(0,1) \times \mathbf{R}+ \rightarrow \mathbf{R}$. The flux
is a function $f \in C^1(\mathbf{R},\mathbf{R}$)  (or $f~:~\mathbf{R} \rightarrow \mathbf{R}$
Lipschitz continuous) and the initial datum is $u_0 \in L^\infty((0,1))$.
Let $A, B \in \mathbf{R}$ be such that $A \le u_0 \le B$ a.e..
The problem to solve is:

\begin{equation}
u_t + (f(u))_x =0, \; x \in (0,1), \; t \in  \mathbf{R}_+,
\label{pipeeq}
\end{equation}

with the initial condition : 

\begin{equation}
u(x,0)=u_0(x), \; x \in (0,1),
\label{pipeeqic}
\end{equation}

and some boundary conditions which will be prescribed later.

\vspace{.1cm}

As in the previous section, let $h=\frac 1 N$ (with $N \in \mathbf{N}^\star$) be the mesh size and $k>0$ be the time step (assumed to be constant, for the sake of simplicity). The discrete unknowns are now
the values $u_i^n \in \mathbf{R}$ for $i \in \{1,\ldots,N\}$ and $n \in \mathbf{N}$. 
In order to define the approximate solution a.e. in $(0,1) \times \mathbf{R}$, one sets $u_{h,k} (x,t)=u_i^n$ for
$x \in ((i-1)h,ih)$, $t \in (nk,(n+1)k)$, $i \in \{1,\ldots,N\}$, $n \in \mathbf{N}$.

\vspace{.1cm}

The discretization of the initial condition leads to

\begin{equation}
u_i^0=\frac 1 h \int_{(i-1)h}^{ih} u_0(x)dx, \; i \in \{1,\ldots,N\}.
\label{pipeeqicd}
\end{equation}

For the computation of $u_i^n$ for $n >0$, one uses, as before, an explicit, $3$-points scheme:

\begin{equation}
\frac h k (u_i^{n+1}-u_i^n) + f_{i+\frac 1 2}^n-f_{i-\frac 1 2}^n= 0, \; i \in \{1,\ldots,N\}, \; n \in \mathbf{N}.
\label{pipeeqd}
\end{equation}

For $i \in 1,\ldots,N-1$, one takes

\begin{equation}
f_{i+\frac 1 2}^n=g(u_i^n, u_{i+1}^n),
\label{numf}
\end{equation}

where $g$ is the
numerical flux. Sufficient conditions on $g~:~[A,B]^2 \rightarrow \mathbf{R}$, in order to have a convergent scheme if $x \in \mathbf{R}$ instead of $(0,1)$, are:

\begin{enumerate}
\item[C1:] $g$ is non decreasing with respect to its first argument and
nonincreasing with respect to its
second argument,
\item[C2:] $g(s,s)=f(s)$, for all $s \in [A,B]$,
\item[C3:] $g$ is Lipschitz continuous.
\end{enumerate}




Let $L$ be a Lipschitz constant for $g$ (on $[A,B]^2$) and $\zeta >0$. If $(0,1)$ is replaced
by $\mathbf{R}$, It is well known (see e.g. \cite{egh}) that, if 
$k \le (1-\zeta) \frac h L$, the approximate solution $u_{h,k}$, that is the solution defined
by (\ref{pipeeqicd})-(\ref{numf})  (with $i \in \mathbf{Z}$), takes its values in $[A,B]$
and converges towards the unique entropy weak solution of (\ref{pipeeq})-(\ref{pipeeqic}) in
$L^p_{loc}(\mathbf{R} \times \mathbf{R}_+)$ as $h \rightarrow 0$.

\vspace{.1cm}

In the case $x \in (0,1)$ instead of $x \in \mathbf{R}$, one assumes the same conditions on $g$, namely
(C1)-(C3).
In order to complete the scheme, one has to define $f_{\frac 1 2}^n$ and  $f_{N+ \frac 1 2}^n$.

\vspace{.1cm}

Let $\overline{u}, \overline{\overline{u}} \in L^\infty(\mathbf{R}_+)$ be
such that $A \le \overline{u}, \overline{\overline{u}} \le B$, a.e. on $\mathbf{R}_+$,
let $g_0, g_1~: ~[A,B]^2
\rightarrow \mathbf{R}$, satisfying (C1)-(C3), and define:

\begin{equation}
\begin{array}{l}
 f_{\frac 1 2}^n=g_0(\overline{u}^n,u_1^n); \; \; \overline{u}^n=\frac 1 k \int_{nk}^{(n+1)k} \overline{u}(t)dt
\\
 f_{N+\frac 1 2}^n=g_1(u_N^n, \overline{\overline{u}}^n,); \; \; \overline{\overline{u}}
 =\frac 1 k \int_{nk}^{(n+1)k} \overline{\overline{u}}(t)dt,
\end{array}
\label{eqbc}
\end{equation}

Then, a convergence theorem can be proven as in the case $x \in \mathbf{R}$, see \cite{v}:

\begin{theo}\label{thm}Let $f \in C^1(\mathbf{R},\mathbf{R}$)  (or $f~:~\mathbf{R} \rightarrow \mathbf{R}$
Lipschitz continuous). Let  $u_0 \in L^\infty((0,1))$, $\overline{u}, \overline{\overline{u}} \in L^\infty(\mathbf{R}_+)$ and $A, B \in \mathbf{R}$ be such that $A \le u_0 \le B$ a.e. on $(0,1)$, 
$A \le \overline{u}, \overline{\overline{u}} \le B$ a.e. on $\mathbf{R}_+$. Let $g_0, g_1~: ~[A,B]^2
\rightarrow \mathbf{R}$, satisfying (C1)-(C3).
Let $L$ be a common Lipschitz constant for $g$, $g_0$ and $g_1$ (on $[A,B]^2$) and let $\zeta >0$.
Then, if $k \le (1-\zeta) \frac h L$, the equations (\ref{pipeeqicd})-(\ref{eqbc}) define
an approximate solution $u_{h,k}$ which takes its values in $[A,B]$
and converges towards the unique solution of (\ref{ewsbc}) in
$L^p_{loc}([0,1] \times \mathbf{R}_+)$ for any $1\le p < \infty$, as $h \rightarrow 0$:

\begin{equation}
\begin{array}{l}
u \in L^\infty((0,1) \times (0,\infty)),
\\
\displaystyle \int_0^\infty \int_0^1 [(u-\kappa)^{\pm}\varphi_t + 
\text{sign}_{\pm}(u-\kappa)(f(u)-f(\kappa)) \varphi_x ] dx dt
\\
\hfill \displaystyle + M \int_0^\infty (\overline u (t) - 
 \kappa)^{\pm} \varphi(0,t) dt + M \int_0^\infty (\overline{\overline{u}}(t) - 
 \kappa)^{\pm} \varphi(1,t) dt
 \\
\hfill  \displaystyle+ \int_0^1 (u_0-\kappa)^{\pm}\varphi(x,0)dx \ge 0,
 \\
\hfill  \forall \kappa \in [A,B], \; \forall \varphi \in C^1_c([0,1] \times 
 [0,\infty), \mathbf{R}_+).
\end{array}
\label{ewsbc}
\end{equation}

In (\ref{ewsbc}), $M$ is any bound for  $\vert 
f'\vert$ on $[A,B]$ (and the solution of  (\ref{ewsbc}) does not depends on the choice of $M$). The
definition of $\text{sign}_{\pm}$ is: $\text{sign}_{+}(s)=1$ if $s>0$, $\text{sign}_{+}(s)=0$ if $s<0$,
$\text{sign}_{-}(s)=0$ if $s>0$, $\text{sign}_{-}(s)=-1$ if $s<0$.
\end{theo}

\begin{remark}~

\begin{enumerate}
\item It is interesting to remark that this convergence result
is also true if the function $g$ depends on  $i$ and $n$, provided that $L$
is a common Lipschitz constant for all these functions.
\item The definition (\ref{ewsbc})
of solution of (\ref{pipeeq})-(\ref{pipeeqic}) with the  ``weak"  boundary conditions $\overline u$
and $\overline{\overline{u}}$ at $x=0$ and $x=1$ is essentially due to F. Otto, see \cite{o}.
\item It is interesting also to remark that if one replaces, in (\ref{ewsbc}), the two entropies 
$(u-\kappa)^{\pm}$ by the sole entropy $\vert u-\kappa \vert$, one has an existence result (since
$\vert u-\kappa \vert=(u-\kappa)^{+} +(u-\kappa)^{-}$) but no uniqueness result,
see \cite{v} for a counter-example to uniqueness. 
\item This convergence result can be generalized to the multidimensional case, see Sect. 
\ref{hyperbcmd} and \cite{v}.
\end{enumerate}
\end{remark}

If $u$, solution of (\ref{ewsbc}), is regular enough (say $u \in C^1([0,1] \times \mathbf{R}_+)$, for
instance), $u$ satisfies $u(0,t)=\overline{u}(t)$ and $u(1,t)=\overline{\overline{u}}(t)$ in the weak
sense given in \cite{bln}.
This condition is very simple if $f$ is monotone:

If $f^\prime >0$, then $u(0,\cdot) = \overline{u}$ and $u$ does not depend on $\overline{\overline{u}}$.

If $f^\prime<0$, then $u(1,\cdot) = \overline{\overline{u}}$ and $u$ does not depend on $\overline{u}$.

\subsection{A very simple example}

One considers here Equation (\ref{pipeeq}), with initial condition (\ref{pipeeqic}) and
weak boundary condition $\overline u$
and $\overline{\overline{u}}$ at $x=0$ and $x=1$, that is in the sense of (\ref{ewsbc}),
in the particular case $f^\prime >0$.
In this case, the main example of numerical flux is $g=g_0=g_1$, $g(a,b)=f(a)$, which leads to
the well known upstream scheme. With this choice of $g_0$ and $g_1$, using the notations of  Sect.
 \ref{conv},
the boundary conditions are taken into account in the form:

\begin{equation}
f_{\frac 1 2}^n=f(\overline{u}^n), \; f_{N+ \frac 1 2}^n=f(u_N^n),
\label{up1}
\end{equation}
with $\overline{u}^n=\frac 1 k \int_{nk}^{(n+1)k} \overline{u}(t)dt$. One may apply the general convergence theorem. The approximate solutions converge (as $h \rightarrow 0$)
towards the solution of (\ref{ewsbc}). In this case, the approximate solutions, as well as
the solution of (\ref{ewsbc}), do not depends on $\overline{\overline{u}}$. 

\vspace{.1cm}

In the case $f^\prime <0$ the main example is  $g=g_0=g_1$, $g(a,b)=f(b)$, which also leads to
the upstream scheme. The boundary conditions are taken into account in the following way:

\begin{equation}
f_{\frac 1 2}^n=f(u_1^n), \; f_{N+ \frac 1 2}^n=f(\overline{\overline{u}}^n),
\label{upbis}
\end{equation}
with  $\overline{\overline{u}}^n=\frac 1 k \int_{nk}^{(n+1)k} \overline{\overline{u}}(t)dt$.

\vspace{.2cm}

These simple cases suggest the following scheme for any $f$, which is
the scalar version of the scheme described in Sect. \ref{discsys} (note that $f'(u)$
is the Jacobian matrix at point $u \in \mathbf{R}$):

\begin{itemize}
\item Boundary condition at $x=0$: 
\begin{equation}
\left\{
\begin{array}{ccc}
 f_{\frac 1 2}^n=f(\overline{u}^n), & \text{ if }  &   f^\prime(u_1^n)>0,
\\
f_{\frac 1 2}^n=f(u_1^n),  & \text{ if }   &  f^\prime(u_1^n)<0.
\end{array}
\right.
\label{sbcl}
\end{equation}
\item Boundary condition at $x=1$: 
\begin{equation}
\left\{
\begin{array}{ccc}
 f_{N+\frac 1 2}^n=f(\overline{\overline{u}}^n), & \text{ if }  &   f^\prime(u_N^n)<0,
\\
f_{N+\frac 1 2}^n=f(u_N^n),  & \text{ if }   &  f^\prime(u_N^n)>0. 
\end{array}
\right.
\label{sbcr}
\end{equation}
\end{itemize}

This solution is not always satisfactory as can be shown on the following
simple example with the B\"{u}rgers equation:

Let $f(s)=s^2$, $u_0=1$ a.e. on $(0,1)$, $\overline {u}=1$ a.e. on $\mathbf{R}_+$ and
$\overline{\overline {u}}=-2$ a.e. on $\mathbf{R}_+$. The exact solution which has to be
approached by the numerical scheme is the unique solution of (\ref{ewsbc}) with
these values of $f$, $u_0$,  $\overline {u}$ and $\overline{\overline {u}}$.
Computing the approximate solution with (\ref{pipeeqicd})-(\ref{numf}), the function
$g$ satisfying $(C2)$, and with (\ref{sbcl})-(\ref{sbcr}), leads to an approximate
solution which is constant and equal to $1$ for any $h$ and $k$. Then, it does not converge (as
$h$ and $k$ go to $0$)
towards the exact solution which is not constant and equal to $1$ since, for the exact
solution, a shock wave with a negative speed starts from the point $x=1$ at
time $t=0$.
Indeed, one can also remark that this approximate solution is the exact solution
of (\ref{ewsbc}) with the same values of $f$, $u_0$,  $\overline {u}$ and with any $\overline{\overline {u}}$ satisfying $\overline{\overline {u}} \ge -1$ a.e. on $\mathbf{R}_+$.
In order to obtain a convergent approximation of the exact solution corresponding to
$\overline{\overline {u}}=-2$, a good choice is, instead of (\ref{sbcr}), $f_{N+\frac 1 2}^n=
g_1(u_N^n, -2)$ with $g_1$ satisfying (C1)-(C3).

\subsection{A simplifed model for two phase flows in pipelines}\label{simple}

It is now possible to understand the treatment of the boundary described in Sect. \ref{pipe} on a simplified model.
This simplifed model for two phase flows in pipelines is given in \cite{tran}. For this model,
the densities are constant so that there are no longer pressure waves but only
the void fraction wave, corresponding to the second eigenvalue of the
original system (\ref{pipesys}). It is also easy to see that for this model, the total flux
(that is the sum of the fluxes of the two phases) is constant in space. One also assumes
that this total flux is constant in time (and positive). System (\ref{pipesys}) is then
reduced to a scalar equation, Equation (\ref{pipeeq}), where the unknown, $u~:~(0,1)
\times \mathbf{R} \rightarrow \mathbf{R}$, is the gas fraction which takes its values between $0$ and $1$.

\vspace{.1cm}

The function $f$ can be taken as $f(s)=a s - b s^2$, where $a, b \in \mathbf{R}$ are given and such that
$0 < b < a <2b$. In (\ref{pipeeq}), the quantity $f(u)$ is the flux of gas. Then, $f(1)-f(u)$ is the
flux of liquid. The function $f$ is increasing between $0$ and $u_M=a/(2b)$ and
decreasing between $u_M$ and $1$. An important value is $u_m \in [0,u_M]$ such that
$f(u_m)=f(1)$.

\vspace{.1cm}

One takes $u_0=0$ a.e. on $[0,1]$ as an initial condition. At $x=0$, the gas flux is given (as
in the complete model, see Sect. \ref{pipe}), one takes $f(u(0,\cdot))=\overline f$ with
$\overline f(t)=c$ for $t \le T$ and $\overline f(t)=0$ for $t \ge T$, where $c$ and $T$ are given
with $c > f(1)$ and $T$ large enough so that $f^\prime$ changes sign at $x=1$ during the simulation. Indeed, in this simplified model, it
is also necessary to take $T$ not too large in order to avoid a problem at $x=0$ (for
$T$ too large, $f^\prime$ will also changes sign at $x=0$). The boundary
condition at $x=1$ will be described on the discrete problem below.

\vspace{.1cm}

The discretization of the problem is performed as before with (\ref{pipeeqicd})-(\ref{numf}),
with $g$ satisfying (C1)-(C3).

\vspace{.1cm}

For the discretization of the boundary condition at $x=0$, the method described in Sect. \ref{pipe}
leads here to

\begin{equation}
f_{\frac1 2}^n=\overline f (nk),
\label{bcl}
\end{equation}

which is indeed  in accordance with the fact that $f^\prime(u_1^n) >0$ for all $n$, at least if
$T$ is not too large.

\vspace{.1cm}

For the discretization of the boundary condition at $x=1$, the first method described in Sect. \ref{pipe}
and given in (\ref{sbcr}), using the sign of $f^\prime(u_N^n)$ leads to
\begin{equation}
\left\{
\begin{array}{ccc}
 f_{N+\frac 1 2}^n=f(u_N^n), & \text{ if }  &   u_N^n<u_M,
\\
f_{N+\frac 1 2}^n=f(1),  & \text{ if }   &   u_N^n>u_M,
\end{array}
\right.
\label{sbcrs}
\end{equation}
and does not lead to the desired results. Note also that $f_{N+\frac 1 2}^n$, given by
(\ref{sbcrs}), is a discontinuous
function of $u_N^n$.

\vspace{.1cm}

The second method, described in Sect. \ref{pipe}, uses the fact that the liquid flux cannot
be negative at $x=1$. Since the liquid flux at $x=1$ is $f(1)-f_{N+\frac 1 2}$ and since
$f(u_m)=f(1)$, this method leads to
\begin{equation}
\left\{
\begin{array}{ccc}
 f_{N+\frac 1 2}^n=f(u_N^n), & \text{ if }  &   u_N^n \le u_m,
\\
f_{N+\frac 1 2}^n=f(u_m),  & \text{ if }   &   u_N^n>u_m,
\end{array}
\right.
\label{bcr}
\end{equation}

Note  that $f_{N+\frac 1 2}^n$, given by
(\ref{bcr}), is a continuous
function of $u_N^n$.
We shall apply the convergence theorem, Theorem \ref{thm} given in Sect. \ref{conv}, for the boundary conditions
(\ref{bcl}) and (\ref{bcr}), and understand the boundary
conditions satisfied by the limit of the approximate solutions. In order to do so, we need
to find $g_0$ and $g_1$,
satisfying (C1)-(C3), and $\overline{u}, \overline{\overline{u}} \in L^\infty(\mathbf{R}_+)$ such that  $f_{\frac 1 2}^n$ and $ f_{N+\frac 1 2}^n$, respectively defined by (\ref{bcl}) and (\ref{bcr}),
satisfy (\ref{eqbc}). Indeed, it is shown in \cite{egv} that both boundary fluxes
$f_{\frac 1 2}^n$ and $ f_{N+\frac 1 2}^n$ may be expressed with the Godunov flux in
the following way:

\begin{itemize}
\item Boundary flux at $x=1$. One takes $\overline{\overline{u}}=1$ a.e. on $\mathbf{R}_+$ and
$g_0$ equal to the Godunov flux, that is $g_0=g_G$ with
\[
g_G(\alpha,\beta) =\left\{\begin{array}{ccc}
\min\{f(s), \; s \in [\alpha, \beta]\}  & \text{ if }  & \alpha \le \beta, \\
\max\{f(s), \; s \in [\beta, \alpha]\}  & \text{ if }  & \alpha > \beta.
\end{array}\right.
\]

The formula (\ref{bcr}) reads
\begin{equation} 
f_{N+\frac 1 2}^n=g_G(u_N^n, 1)=\left\{\begin{array}{ccc}
f(u_N^n)  & \text{ if }  & u_N^n \le u_m, \\
f(1)  & \text{ if }  & u_N^n > u_m.
\end{array}\right.
\label{bcrf}
\end{equation}

\item Boundary flux at $x=0$. One assumes (for simplicity) that $\frac T k \in \mathbf{N}$. let
$\alpha, \beta \in (0,1)$ such that $\alpha < \beta$ and $f(\alpha)=f(\beta)=c$. One takes

\begin{equation}
\overline{u} (t) =\left\{\begin{array}{ccc}
\alpha  & \text{ if }  & t <T, \\
0  & \text{ if }  & t >T,
\end{array}\right.
\label{ubar}
\end{equation}

so that, recalling that $\overline{u}^n=\frac 1 k \int_{nk}^{(n+1)k} \overline{u}(t)dt$,

\[
f(\overline{u}^n)=\left\{\begin{array}{ccc}
c  & \text{ if }  & nk <T, \\
0  & \text{ if }  & nk  \ge T,
\end{array}\right.
\]

Then, if $u_1^n \le  \beta$, the  formula (\ref{bcl}) reads
\begin{equation} 
f_{\frac 1 2}^n=g_G(\overline{u}^n, u_1^n),
\label{bclf}
\end{equation}
since, in this case, $g_G(\overline{u}^n, u_1^n)=f(\overline{u}^n)$.
The fact that $u_1^n \le  \beta$ is true for all $n$ if $T$ is not too large. If $T$ is too
large, the convergence result can be applied with (\ref{bclf}) instead of (\ref{bcl}).

\end{itemize}

It is now possible to apply Theorem \ref{thm}.
Let $L$ be a common Lipschitz constant for $g$ and $g_G$ (on $[0,1]^2$) and let $\zeta >0$.
If $k \le (1-\zeta) \frac h L$, the approximate solution $u_{h,k}$, that is the solution defined
by (\ref{pipeeqicd})-(\ref{numf}), with the boundary fluxes (\ref{bcrf})-(\ref{bclf}) (and $u_0=0$, $\overline{\overline{u}}=1$ and $\overline u$ given by (\ref{ubar})), takes its values in $[0,1]$
and converges towards the unique solution of (\ref{ewsbcpipe}) in
$L^p_{loc}([0,1] \times \mathbf{R}_+)$ for any $1\le p < \infty$, as $h \rightarrow 0$:

\begin{equation}
\begin{array}{l}
u \in L^\infty((0,1) \times (0,\infty)),
\\
\displaystyle \int_0^\infty \int_0^1 [(u-\kappa)^{\pm}\varphi_t + 
\text{sign}_{\pm}(u-\kappa)(f(u)-f(\kappa)) \varphi_x ] dx dt
\\
\hfill \displaystyle + M \int_0^\infty (\overline u (t) - 
 \kappa)^{\pm} \varphi(0,t) dt + M \int_0^\infty (1 - 
 \kappa)^{\pm} \varphi(1,t) dt
 \\
\hfill  \displaystyle+ \int_0^1 (0-\kappa)^{\pm}\varphi(x,0)dx \ge 0,
 \\
\hfill  \forall \kappa \in [0,1], \; \forall \varphi \in C^1_c([0,1] \times 
 [0,\infty), \mathbf{R}_+).
\end{array}
\label{ewsbcpipe}
\end{equation}

If $u$, solution of (\ref{ewsbcpipe}), is regular enough on $[0,1] \times (0,T)$, then, it is possible to prove that $u$ satisfies the boundary conditions, for $0<t < T$, in the following sense (see \cite{v} and \cite{egv}):

\begin{itemize}
\item Boundary condition at $x=0$ (recall that $\overline{u}$ is given by (\ref{ubar})): $u(0,t)=\alpha$ or $u(0,t) \ge  \beta$. In fact, if $T$ is not too
large, one has $u(0,t)=\alpha$.
\item Boundary condition at $x=1$: $u(1,t) \le u_m$ or $u(1,t) =1$,
\end{itemize}

Thanks to Theorem \ref{thm}, it is possible to give other choices for $f_{N+\frac 1 2}^n$ for which the approximate solutions
obtained with this new choice of $f_{N+\frac 1 2}^n$ converge towards the same function $u$, which
is the unique solution
of (\ref{ewsbcpipe}). Indeed, let $h~:~[0,1] \rightarrow \mathbf{R}$ be a nondecreasing function
such that $h \le f$ and $h(1)=f(1)$ and take: 

\begin{equation} 
f_{N+\frac 1 2}^n=h(u_N^n).
\label{bcrnew}
\end{equation}

One may construct a function $g_1$ satisfying (C1)-(C3) such that 
$h(s)=g_1(s,1)$, for all $s \in [0,1]$, and then use Theorem \ref{thm}.
Let $L$ be a common Lipschitz constant for $g$ and $g_G$ and $g_1$ (on $[0,1]^2$) and let $\zeta >0$.
If $k \le (1-\zeta) \frac h L$, the approximate solution $u_{h,k}$, that is the solution defined
by (\ref{pipeeqicd})-(\ref{numf}), with the boundary fluxes (\ref{bcrnew}) and (\ref{bclf}) (and $u_0=0$, $\overline{\overline{u}}=1$ and $\overline u$ given by (\ref{ubar})), takes its values in $[0,1]$
and converges towards the unique solution of (\ref{ewsbcpipe}) in
$L^p_{loc}([0,1] \times \mathbf{R}_+)$ for any $1\le p < \infty$, as $h \rightarrow 0$.

\vspace{.3cm}

Turning back to the complete system described in Sect. \ref{pipe}, the analysis of this
simplified model for two phase flows in pipelines may also suggest another way to take into account
the boundary condition at $x=1$ (with a given numerical flux $g_1$):

\begin{enumerate}

\item Compute $DF(w_N^n)$, its eigenvalues $\{\lambda_1, \lambda_2,\lambda_3\}$  and
a basis of $\mathbf{R}^3$, $\{\varphi_1, \varphi_2,\varphi_3\}$, such that $DF(w_N^n) \varphi_i=
\lambda_i \varphi_i$, $i=1,2,3$,

\item write $w_N^n$ on the basis $\{\varphi_1, \varphi_2,\varphi_3\}$, namely
$w_N^n=\alpha_1 \varphi_1 +\alpha_2 \varphi_2 + \alpha_3 \varphi_3$,

\item Since $\lambda_3 <0$ and since one wants $Q_l \ge 0$, compute 
$w_{N+1}^n=\beta_1 \varphi_1 +\beta_2 \varphi_2 + \beta_3 \varphi_3$ and
$F_{N+\frac 1 2}^n=g_1(w_N^n,w_{N+1}^n)$ with the following $3$ conditions on the components
of $w_{N+1}^n$:
usual condition on the pressure, $\beta_3= \alpha_3$ and $R_{N+1}^n=1$ where
$R_{N+1}^n$ is the gas fraction computed with $w_{N+1}^n$.

\end{enumerate}


\chapter{Multidimensional nonlinear hyperbolic equations} \label{hypmd}

The aim of this chapter is to define and study finite volume schemes for the 
approximation of the solution to a nonlinear scalar hyperbolic problem in 
several space dimensions. 
Explicit and implicit time discretizations are considered. We   
prove the 
convergence of the approximate solution towards the entropy weak solution 
of the problem and  give an 
error estimate between the approximate solution and the entropy weak solution  
with respect to the discretization mesh size.

\section{The continuous problem}

\hpage{Hyperbolic multid}

We consider here the following nonlinear hyperbolic equation in $\dim$
space dimensions ($\dim \ge 1$), with initial 
condition\index{Equation!Hyperbolic --} 

\be
u_t (x,t) + \div (\vbf f(u))(x,t)  = 0, \,  x  \in  \R^\dim, 
\,   t \in \R_+, \label{estpbcont}
\ee
\be
u(x,0)  = u_0(x), \,  x  \in \R^\dim,
\label{estcondini}
\ee

where $u_t$ denotes the time derivative of $u$ ($t \in \R_+$), and $\div$
the divergence operator with respect to the space variable (which belongs to
$\R^\dim$). Recall that $\vert x \vert$ denotes the euclidean
norm of $x$ in $\R^\dim$, and $x \cdot y$ the usual scalar
product of $x$ and $y$ in $\R^\dim$.

\medskip

The following hypotheses are made on the data:

\begin{hyp}\label{estdonnees} 

\[
\bal
(i) \; u_0\in L^\infty(\R^\dim), \, U_m, U_M \in \R, \; U_m \le  u_0 \le 
U_M \; a.e.,
\\
(ii) \; \vbf \in C^1(\RNRp,\R^\dim),
\\
(iii) \; \div\vbf(x,t) = 0, \, \forall (x,t) \in \RNRp,
\\
(iv) \; \exists V < \infty \text{ such that } \vert\vbf(x,t)\vert \le V, \hskip 
.1 cm \forall (x,t) \in \RNRp,
\\
(v) \;  f \in C^1(\R,\R).
\ean
\]
\end{hyp}

\begin{remark} \normalfont \label{baptisee}
Note that part (iv) of Assumption \ref{estdonnees} is crucial. It ensures 
the property of ``propagation in finite time" which is 
needed for the uniqueness of the solution of \eqref{estkruzkov} and for the stability (under a 
``Courant-Friedrichs-Levy" (CFL) condition) of the time  explicit  numerical scheme. Part (iii) 
of Assumption \ref{estdonnees}, on the other hand, is only considered for 
the sake 
of simplicity; the results of existence and uniqueness of the entropy weak 
solution   and convergence (including error estimates as in the theorems 
\refp{esth14} and \refp{esth14i}) of the numerical 
schemes presented below may be extended  to the case $\div \vbf \not = 0$. 
However, part (iii) 
of Assumption \ref{estdonnees} is natural in many ``applications" and avoids 
several technical complications. Note, in particular, that, for instance, 
if $\div \vbf \ne 0$, the $L^\infty$-bound on the solution of 
\eqref{estkruzkov} and the $L^\infty$ estimate (in Lemma \ref{estinftyest} and 
Proposition \ref{esteu}) on  the approximate solution depend on  $\vbf$ and 
$T$.
The case $F(x,t,u)$ instead of 
$\vbf(x,t)f(u)$ is also feasible, but somewhat more technical, see 
\citet{Chainais1} and \citet{Chainais2}.
\end{remark}

Problem \eqref{estpbcont}-\eqref{estcondini} has a unique 
entropy weak solution,
which is the solution to  the following equation (which is the multidimensional extension
of the one-dimensional definition \refp{defent1}).

\begin{equation}\label{def-ewsmd}
\ba 
u \in L^\infty(\R^\dim \times \Rpast), \\
\dsp{\int_{\Rp}\int_{\R^\dim}} \Bigl[ \eta(u(x,t)) \varphi_t(x,t) +
\Phi(u(x,t)) 
\vbf(x,t)\cdot\grad\varphi(x,t) \Bigr] dx dt \, + \\
\dsp{\int_{\R^\dim}}  \eta (u_0(x)) \varphi(x,0) dx \ge 0,
\, \forall\varphi\in C^\infty_c(\RNRp,\R_+), \\
\forall \eta \in C^1(\R,\R), \text{ convex function, and } \Phi \in C^1(\R,\R) 
\text{ such that } \Phi' = f' \eta',
\ea
\end{equation}

where $\grad\varphi$ denotes the gradient of the function $\varphi$
with respect to the space variable (which belongs to $\R^\dim$). 
Recall that $C^m_c(E,F)$ denotes the set of functions $C^m$  
from $E$ to $F$, with compact support in $E$.

The characterization of the entropy weak solution \index{Entropy!-- weak 
solution} by the Krushkov 
entropies (proposition \refp{prentkr1}) still holds in the multidimensional case.
Let us define again, for all $\sk \in\R$, the Krushkov entropies $(\vert \cdot - \sk
\vert)$ for which the entropy flux is $f(\cdot\top\sk) - f(\cdot\bot\sk)$ 
(for any pair of real values $a,b$, we denote again by $a\top b$ the maximum
of $a$ and $b$, and by $a\bot b$ the minimum of $a$ and $b$).
The unique entropy weak solution \index{Entropy!-- weak 
solution} 
is also the unique solution to the following problem: 

\be\ba 
u \in L^\infty(\R^\dim \times \Rpast), \\
\dsp{\int_{\Rp}\int_{\R^\dim}} \Bigl[ \vert u(x,t) - \sk \vert \varphi_t(x,t) +
\Bigl( f(u(x,t)\top\sk) - f(u(x,t)\bot\sk) \Bigr) 
\vbf(x,t)\cdot\grad\varphi(x,t) \Bigr] dx dt \, + \\
\dsp{\int_{\R^\dim}} \vert u_0(x) - \sk \vert \varphi(x,0) dx \ge 0,
\, \forall \sk \in \R, \, \forall\varphi\in C^\infty_c(\RNRp,\R_+).
\ea\label{estkruzkov}\ee

As in  the one-dimensional  case (Theorem \refp{hyp1dexistu}), 
existence and uniqueness results are also known for the 
entropy weak solution  to Problem  (\ref{estpbcont})-(\ref{estcondini}) under assumptions which 
differ slightly from 
assumption \ref{estdonnees}  (see e.g. \citet{krushkov}, \citet{Vo}). 
In particular, these results are
obtained with a nonlinearity  $F$ (in our case $F = \vbf f$) of class $C^3$. 
We recall that the
methods which were used in \citet{krushkov} require a regularization in
$BV(\Rn)$ of the function $u_0$,
in order to take advantage, for any $T>0$, of compactness properties 
which are similar to those given in Lemma \refp{helly} for the case 
$\dim=1$. Recall that the space   $BV(\O)$ where $\O$ is an 
open subset of $\R^p$,  $p \ge 1$, was defined in Definition 
\refp{bvrd}; it will be  used later with  $\O=\Rn$ or $\O=\Rn \times (-T,T)$.

\medskip

The existence of solutions to  similar  problems to  (\ref{estpbcont})-(\ref{estcondini}) 
was already proved by passing to the limit on solutions of an appropriate 
numerical scheme, see \citet{CS}. 
The work of \citet{CS} uses a finite difference scheme on a uniform 
rectangular grid, in two space dimensions, and requires that the initial condition $u_0$ belongs 
to $BV(\R^\dim)$ (and thus, the
solution to Problem  (\ref{estpbcont})-(\ref{estcondini}) also  has a locally 
bounded variation). These assumptions (on meshes and on $u_0$) yield, as 
in Lemma \refp{helly}, a (strong) compactness property in 
$L^1_{loc}(\R^\dim \times \R_+)$ on 
a family of approximate solutions. In the following, 
however, we shall only require that  $u_0 \in L^\infty(\R^\dim)$ and we 
shall be able to deal with more 
general meshes. We may use, for instance, a triangular mesh in the case 
of two space dimensions. For each of these reasons, the $BV$  framework 
may not be used and a (strong) compactness property in $L^1_{loc}$ on 
a family of approximate solutions is not easy to obtain (although this 
compactness property does hold and results from this chapter). 
In order to prove the existence of a solution to 
(\ref{estpbcont})-(\ref{estcondini}) by passing to the limit on the 
approximate solutions given by  finite 
volume schemes on general meshes (in the sense used below) in two or 
three space dimensions, 
we shall work with some ``weak" compactness result in $L^\infty$, namely 
Proposition \ref{nlwsprop}, which yields the ``nonlinear 
weak-$\star$ convergence"
\index{Convergence!Nonlinear weak star --} 
(see Definition \refp{nlwsdef}) of a family of 
approximate solutions. 
When doing so, passing to the limit with the 
approximate solutions will give the existence of an
``entropy process solution"\index{Entropy!-- process 
solution}  to Problem 
(\ref{estpbcont})-(\ref{estcondini}), see Definition \refp{estepsol}.
A uniqueness result for the entropy process solution \index{Entropy!-- process
solution} to Problem 
(\ref{estpbcont})-(\ref{estcondini}) is then  proven. 
This uniqueness result proves that the entropy process solution is 
indeed the entropy weak solution, hence the existence and  uniqueness of 
the entropy weak solution. This uniqueness result  also  allows  us to 
conclude to the convergence of the approximate solution
given by the numerical scheme (that is (\ref{estschema}), (\ref{estschemaz})) towards the 
entropy weak solution  to (\ref{estpbcont})-(\ref{estcondini}) (this 
convergence holds in $L^p_{loc}(\R^\dim \times \Rp)$ for any $1 \le p < 
\infty$).

\medskip

Note that uniqueness results for ``generalized" solutions (namely measure valued 
solutions)  to (\ref{estpbcont})-(\ref{estcondini}) have recently been 
proved (see \citet{DP},  \citet{Sz}, \citet{GH}). The proofs of these 
results rely on the one hand on the concept of measure valued solutions and on the
other hand on the existence of an entropy weak solution. The direct proof 
of the uniqueness of a measure valued solution (i.e. without assuming any 
existence result of entropy weak solutions) leads to a difficult problem 
involving the application of the theorem of continuity in mean.
This difficulty is easier to deal within the framework of  
entropy process solutions  (but in fact,  measure valued 
solutions and entropy process solutions are two presentations of the same concept). 

\medskip

Developing the above analysis gives a (strong) convergence result of approximate 
solutions towards the entropy weak solution. But moreover,  we  also derive some   
error estimates depending on the regularity of 
$u_0$.

\medskip

In the case of a Cartesian grid, the convergence and error 
analysis reduces essentially to a one-dimensional discretization problem  for which 
results were proved some time ago, see e.g. \citet{Kuznetsov}, \citet{CM}, 
\citet{S}.
In the case of general meshes, the numerical schemes are not 
generally ``TVD" (Total Variation Diminushing) and therefore the classical framework of the 1D case (see 
Section \refp{bvmethod}) may 
not be used. More recent works deal with several 
convergence results and error estimates for time explicit finite volume 
schemes, see e.g. \citet{Coquel}, \citet{cgh}, \citet{Vi}, 
\citet{kroner1}, \citet{kroner2}, \citet{kronerbook}: 
following Szepessy's work on the convergence of the streamline diffusion 
method (see \citet{Sz}), most of these works use DiPerna's uniqueness
theorem, see \citet{DP} (or an adaptation of it, see \citet{GH} and \citet{EGH2}), 
and the error 
estimates generalize the work by  \citet{Kuznetsov}.
Here we  use the framework of \citet{cgh}, \citet{eggh}; 
we   prove directly that any monotone flux scheme  (defined below)
satisfies a ``weak $BV$" estimate (see 
lemmata \refp{estwbvest} and \refp{estwbvesti}). 
This inequality appears to be a key 
for the proof of convergence and for the error estimate. Some convergence 
results and error estimates are also possible with some so called ``higher 
order schemes" which are not monotone flux schemes (briefly presented 
for the 1D case in section \refp{hyperhos}). These results are 
not presented here, see \citet{Noelle} and \citet{Chainais1}  for some of them.

\medskip

Note that the nonlinearity considered here is of the form $\vbf(x,t) f(u)$. This 
kind of flux is often encountered in porous medium modelling, where the 
hyperbolic equation may then be coupled with an elliptic or parabolic 
equation (see e.g. \citet{EG}, \citet{Vignal1}, \citet{Vignal2},
\citet{HL}). It adds an extra difficulty to 
the case $F(u)$ because of the dependency on $x$ and $t$. Note again (see 
Remark \ref{baptisee}) that the 
method which we present here for a nonlinearity of the form  $\vbf(x,t) 
f(u)$  also yields the same results in  the case of a nonlinearity of the
form $F(x,t,u)$, see the recent work  of \citet{Chainais2}.  

\medskip

The time implicit\index{Scheme!Implicit -- for hyperbolic equations} 
discretization  adds the extra difficulties of proving the existence of the
approximate solution  (see Lemma \refp{esteu}) and proving a so called 
``strong time $BV$ estimate" (see Lemma 
\refp{estlbvt}) 
in order to show that the error estimate for the implicit scheme may still 
be of order $h^{1/4}$ even if the time step 
$k$ is of order $\sqrt h$, at least in 
particular cases.

\medskip

We first describe  in section \ref{eghbns}  finite volume schemes using a
``general" mesh for the discretization of  
(\ref{estpbcont})-(\ref{estcondini}). In sections \ref{stabe} and 
\ref{stabi} some estimates on the approximate solution given by the 
numerical schemes are shown and in Section \ref{entineq} some entropy 
inequalities are proven. We then prove  in section \ref{eghcvsol}
the convergence of convenient subsequences of sequences of
approximate solutions towards an entropy process solution\index{Entropy!-- 
process
solution}, by  
passing to the limit when the mesh size
and the time step go to $0$. 
A byproduct of this result is the existence of an entropy process solution
to (\ref{estpbcont})-(\ref{estcondini}) (see Definition \refp{estepsol}).
The uniqueness of the \epsol  to problem  
(\ref{estpbcont})-(\ref{estcondini}) is then proved; we can therefore 
conclude to the existence and uniqueness of the entropy weak solution  and also
 to the $L^p_{loc}$ convergence for any finite $p$ of the 
approximate solution towards the entropy weak solution (Section 
\ref{eghcvsol}). Using the existence of the entropy weak solution, an 
error estimate result is given in Section \ref{esterr} (which also yields 
the convergence result). Therefore the 
main interest of this convergence result is precisely to prove the 
existence of the entropy weak solution to 
(\ref{estpbcont})-(\ref{estcondini}) without any regularity assumption  on
the initial data.
Section \ref{nlwscv} describes the notion of nonlinear weak-$\star$ 
convergence, which is widely used in the proof of convergence of section 
\ref{eghcvsol}.

\smallskip

Section \ref{fve} is not related to the previous sections. It describes 
a finite volume approach which may be used to stabilize finite element schemes for 
the discretization of a hyperbolic equation (or system).

\medskip

\section{Meshes and schemes}\label{eghbns}

\hpage{Hyperbolic multid}

Let us first define an admissible mesh of $\R^\dim$ as a generalization
of the notion of admissible mesh of $\R$ as defined in definition \refp{meshhyp1d}.

\begin{defi}[Admissible
meshes]\index{Mesh!Admissible --!--for hyperbolic equations}\label{meshypmd}
  \normalfont An admissible finite volume mesh of $\R^d$, with 
$d = 1, 2$ or $3$ (for the
discretization of Problem  (\ref{estpbcont})-(\ref{estcondini})),  denoted by
$\mesh$, is given by  a family of disjoint polygonal connected subsets of 
$\R^d$ such that $\R^\dim$ is the union of the closure of the elements 
of $\mesh$ (which are called
control volumes in the following) and such that the common 
``interface" of any two control volumes is
included in a  hyperplane of $\R^\dim$  (this is not 
necessary but  is introduced to simplify the formulation). Denoting by 
$h = \sizemesh = \sup \{\diam (\cv), \cv \in \mesh\}$, it is assumed that $h < + 
\infty$ and that, for some $\a>0$,

\be \bal
\a h^\dim  \le  \mcv, \\
\meas(\dr\cv)  \le   \frac 1 {\a} h^{\dim-1}, \, \forall \cv \in \mesh,
\ean
\label{estmaillage} 
\ee
\medskip

where $\mcv$ denotes the $\dim$-dimensional Lebesgue measure of $\cv$, 
 $\meas(\dr\cv)$ denotes the $(\dim-1)$-dimensional Lebesgue measure 
of $\dr\cv$ ($\dr\cv$ is the boundary of $\cv$) and  
$\Nv$ denotes the set of neighbours
of the control volume $\cv$; for $\cvv \in \Nv$, 
we denote by $\edgecvcvv$ the
common interface between $\cv$ and $\cvv$, and by $\bfn_{\cv,\cvv}$ the unit normal
vector to $\edgecvcvv$ oriented from $\cv$ to $\cvv$.
The set of all the interfaces is denoted by $\edges$.

\end{defi}

Note that, in this definition, the terminology is ``mixed". 
For $d=3$, ``polygonal" stands for ``polyhedral" and, for $d=2$, 
``interface" stands for ``edge". For $d=1$ definition \ref{meshypmd} is 
equivalent to 
definition \refp{meshhyp1d}.

\medskip

In order to define the numerical flux, we consider functions  
$g \in C(\R^2,\R)$ satisfying the following assumptions: 

\begin{hyp}\label{estcondflu}
Under Assumption \ref{estdonnees} the function $g$, only depending  on $f$, 
$\vbf$, $U_m$ and $U_M$, satisfies

\begin{itemize} 
\item $g$ is locally Lipschitz continuous from $\R^2$ to $\R$,
\item $g(s,s) = f(s)$, for all $s \in [U_m,U_M]$,
\item $(a,b) \mapsto g(a,b)$, from $[U_m,U_M]^2$ to $\R$, is nondecreasing with respect to $a$
and nonincreasing with respect to $b$.
\end{itemize}
Let us denote by $g_1$ and $g_2$ 
the Lipschitz constants of $g$ on $[U_m,U_M]^2$ with respect to its 
two arguments.
\end{hyp}

The hypotheses on $g$ are the same as those presented 
for monotone flux schemes in the one-dimensional case (see definition 
\refp{monflux}); the function $g$ allows the
construction of a numerical flux, see Remark \ref{deflunum} below.

\begin{remark} \normalfont  In Assumption \ref{estcondflu}, the third item will 
ensure some stability  properties of the schemes defined below. In 
particular, in the case of the ``explicit scheme" (see \eqref{estschema}), it yields 
the monotonicity of the scheme under a CFL condition (namely, condition
\eqref{estcfl} with $\xi=0$). The second item is essential since it ensures
the consistency  of the fluxes. All the examples of functions
$g$ given in Examples \refp{flumono} satisfy these assumptions. We again give
the important example of the ``generalized 1D Godunov scheme"
obtained with a one-dimensional Godunov scheme for each interface (see e.g.,
for the explicit scheme, see \citet{Coquel}, \citet{Vi}),

\[
g(a,b) =\left\{
\begin{array}{ll}
\max\{f(s), \, b\le s\le a\} \ \text{if} \ b\le a \\
\min\{f(s), \, a\le s\le b\} \ \text{if} \ a\le b,
\end{array}\right.
\]

and also the framework of some ``flux splitting" schemes:

\[
g(a,b) = f_1(a) + f_2(b),
\]

with $f_1$, $f_2 \in 
C^1(\R,\R)$, $f=f_1+f_2$, $f_1$ nondecreasing and $f_2$ nonincreasing (this framework 
is considerably more simple that 
the general framework, because it reduces the study to the particular case
of two monotone nonlinearities). 

\medskip

Besides, it is possible to replace Assumption \ref{estcondflu} on $g$ by some slightly more general
assumption, in order to handle, in particular, the case of some ``Lax-Friedrichs 
type" schemes (see Remark \ref{estext} below).
\end{remark}

In order to describe the numerical schemes considered here, let $\mesh$ 
be an admissible mesh in the sense of Definition  \ref{meshypmd} and $k>0$ 
be the time step.
The discrete unknowns are $\ucvn$, $n \in \N^\star$, $\cv\in\mesh$.
The set $\{\ucvz,$ $\cv\in \mesh\}$ is given by the initial condition,

\be 
\ucvz  = \frac 1 {\mcv} \int_\cv u_0(x) dx, \forall \, \cv \in \mesh.
\label{estschemaz}\ee

The equations satisfied by the discrete unknowns, $u^n_\cv$, $n \in 
\N^\star$, $\cv \in \mesh$, are obtained by discretizing  equation 
\eqref{estpbcont}.
We now describe the explicit and implicit schemes.

\subsection{Explicit schemes}

We   present here the ``explicit scheme" associated to a function $g$ 
satisfying Assumption~\ref{estcondflu}. In this case, for 
stability reasons \index{Stability} (see lemmata \ref{estinftyest} and \ref{estwbvest}), 
the time step $k\in\Rpast$ is chosen such that

\be
k \le (1 - \xi) \frac {\a^2 h} { V (g_1 + g_2)},
\label{estcfl}\ee
where  $\xi\in (0,1)$ is a given  real value; recall that $g_1$ and $g_2$
are the Lipschitz constants of $g$ with respect to the first and second
variables on $[U_m,U_M]^2$ and  that $U_m \le u_0 \le
U_M$ a.e. and $\vert \vbf(x,t) \vert \le V< + \infty$, for all 
$(x,t) \in \R^d\times \R_+$.  
Consider the following explicit numerical 
scheme\index{Finite volume scheme!-- for hyperbolic problems!-- in two or 
three space dimensions}:

\be
\mcv \frac {\ucvnp - \ucvn} k + \sumv 
\left(\vcvcvvn \ g(\ucvn,\ucvvn) - \vcvvcvn \ g(\ucvvn,\ucvn)\right) = 0, \, \forall
\cv\in\mesh, \, \forall n\in\N,  \label{estschema}\ee

where

\[
\vcvcvvn = \dsp{ \frac 1 k \int_{n k}^{(n+1)k} 
\int_{\edgecvcvv} (\vbf(x,t)\cdot\bfn_{\cv,\cvv})^+ \dfrontiere(x) dt}
\] 

and

\[
\bari
\vcvvcvn = \dsp{\frac 1 k \int_{n k}^{(n+1)k}
\int_{\edgecvcvv} (\vbf(x,t)\cdot\bfn_{\cvv,\cv})^+
\dfrontiere(x) dt}
\\
= \dsp{\frac 1 k \int_{n k}^{(n+1)k} 
\int_{\edgecvcvv} (\vbf(x,t)\cdot\bfn_{\cv,\cvv})^- 
\dfrontiere(x) dt}.
\ean
\]

Recall that $a^+=a \top 0$ and $a^-=- (a \bot 0)$ for all $a \in \R$ and 
 that $\dfrontiere$ is the 
integration symbol for the $(\dim-1)$-dimensional Lebesgue measure on the 
considered hyperplane.

\begin{remark}[Numerical fluxes]\label{deflunum}  \normalfont
The numerical flux at the interface between  the control volume $\cv$ and  the control volume
$\cvv\in\Nv$ is then equal to 
$\vcvcvvn \ g(\ucvn,\ucvvn) - \vcvvcvn \ g(\ucvvn,\ucvn)$;  this 
expression yields a monotone flux such as defined  in definition \refp{monflux}, given 
in the one-dimensional case. However, in the multidimensional case,
 the expression of the numerical flux depends on the considered 
interface; this was not so in the one-dimensional case for which the 
numerical flux is completely
defined by the function $g$.
\end{remark}

The approximate solution, denoted by $\uapp$, is defined a.e. 
from $\RNRp$ to $\R$ by 
\be 
\uapp(x,t) = \ucvn, \text{ if } x \in \cv, \, t \in \, [n k ,(n+1) k), \,  \cv
\in \mesh, \,  n 
\in \N.
\label{estuapp} 
\ee

\medskip

\subsection{Implicit schemes}

The use of implicit schemes is steadily increasing in 
industrial codes  for  reasons  such as robustness and 
computational cost. Hence we    consider in our analysis 
the following implicit numerical scheme \index{Finite volume scheme!--
for hyperbolic problems!-- in two or three space dimensions} 
(for which condition \eqref{estcfl} is
no longer needed) associated to a function $g$ 
satisfying Assumption~\ref{estcondflu}:

\be
\mcv \frac {\ucvnp - \ucvn} k + \sumv 
(\vcvcvvn \ g(\ucvnp,\ucvvnp) - \vcvvcvn \ g(\ucvvnp,\ucvnp)) = 0, \, \forall
\cv\in\mesh, \, \forall n\in\N.  \label{estschemai}\ee

where $\{u^0_\cv$, $\cv \in \mesh\}$ is still determined by 
\eqref{estschemaz}.
The implicit approximate solution $\uapp$, is defined now a.e. from $\RNRp$ to $\R$ by
\be 
\uapp(x,t) = \ucvnp, \text{ if } x \in \cv, \, t \in \, (n k ,(n+1) k], \, \cv \in
\mesh, \,  n 
\in \N.
\label{estuappimp} 
\ee

\subsection{Passing to the limit}

We    show in section \refp{eghcvsol}  the convergence of the
approximate solutions $u_{\mesh,k}$ (given by the numerical schemes 
above described)
towards the unique entropy weak solution    $u$  to  (\ref{estpbcont})-(\ref{estcondini})
in an adequate sense, 
when
$\sizemesh \rightarrow 0$ and $k \rightarrow 0$ (with, possibly, a stability condition).
In order to describe the general line of thought leading to this convergence result, 
we shall simply
consider  the explicit scheme (that is \eqref{estschemaz}, 
\eqref{estschema} and \eqref{estuapp}) (the implicit scheme will also be 
fully investigated later).

First, in section \ref{stabe}, 
by writing $u^{n+1}_\cv$ as a convex combination of $u^n_\cv$ and  
$(u^n_{\cvv})_{\cvv \in \Nv}$,
the $L^\infty$ stability is easily shown under the  CFL condition \eqref{estcfl}
($u_{\mesh,k}$ is proved to  be bounded in 
$L^\infty(\R^\dim \times \Rpast)$, 
independently of  $\sizemesh$ and $k$).

\smallskip

Let $\mesh$ be an admissible mesh in the sense of Definition \refp{meshypmd} and let $k$ satisfy (\ref{estcfl})); by a classical argument, if any possible limit of a family of approximate solutions  $u_{\mesh,k}$ is the entropy weak solution  to problem (\ref{estpbcont})-(\ref{estcondini}) then  $u_{\mesh,k}$ converges (in  $L^\infty(\R^\dim \times \Rpast)$ for the weak-$\star$ topology, for instance), as  $h = \size(\mesh)\rightarrow 0$  (and $k$  satisfies \eqref{estcfl}),  towards the unique entropy weak solution to problem  (\ref{estpbcont})-(\ref{estcondini}). Unfortunately, the $L^\infty$ estimate of section \ref{stabe} does not yield  that any possible limit of a family  of approximate solutions is solution to problem (\ref{estpbcont})-(\ref{estcondini}), even in the linear case ($f(u) = u)$ (see the proofs of convergence of Chapter \ref{hyper1d}). 
The ``$BV$  stability" can be used (combined with the  $L^\infty$ stability) to show the convergence in the case of one space dimension (see section \refp{bvmethod})  and in the case of Cartesian meshes in two or three space dimensions. 
Indeed, in the case of Cartesian meshes, assuming $u_0 \in BV(\Rn)$ and assuming (for simplicity) $\vbf$ to be constant (a generalization is possiblefor $\vbf$ regular enough), the following estimate holds, for all $T \ge k$:

\[ 
k \sumTsk \sum_{\edgecvcvv \in \edges} \medgecvcvv \vert
u^n_{\cv} - u^n_{\cvv} \vert \le T \vert u_0 \vert_{BV(\Rn)}, 
\]

where $\NTk \in \N$ is such that $(\NTk+1) k
\le T < (\NTk+2)k$, and  the values $u^n_\cv$ are given by 
\eqref{estschemaz} and \eqref{estschema}. Such an estimate is wrong in the general case
of admissible meshes in the sense of Definition \refp{meshypmd}, as it
can be shown with easy counterexamples. 
It is, however, not necessary for the proof of
convergence. A weaker inequality, which is called ``weak $BV$" as in the  
one-dimensional case (see lemma \refp{nlwb1d}) will be shown in the 
multidimensional case for both 
explicit and implicit schemes (see lemmata \refp{estwbvest} and 
\refp{estwbvesti});
the weak $BV$ estimate  yields the 
convergence of the scheme in the general case. As an illustration, consider 
the case $f' \ge 0$; using an upwind scheme, i.e. $g(a,b) = f(a)$, the 
weak $BV$ inequality  \refep{estbvfx}, which is very close to that of 
the 1D case (lemma \refp{nlwb1d}), reads
 
\be
\begin{array}{r}
\sumTsk k \sumER 
(\vcvcvvn + \vcvvcvn) \abs{f(\ucvn) - f(\ucvvn)} \le \ \ \dsp{\frac {C} {\sqrt{h}}},
\end{array} 
\label{BVfaible}\ee

where $\edges_R^n = \{(\cv,\cvv)\in \mesh^2, \cvv\in \Nv, 
\edgecvcvv\subset B(0,R)$ and $\ucvn > \ucvvn \}$ and $C$ only depends  
on $\vbf$, $g$, $u_0$, $\a$, $\xi$, $R$ and $T$ (see Lemma \ref{estwbvest}).

We say that Inequality \eqref{BVfaible} is ``weak", but it is in 
fact ``three times weak" for the following reasons:

\begin{enumerate}
\item the inequality is of order $\frac  1 {\sqrt h}$, and not of order 1. 

\item In the left hand side of \eqref{BVfaible}, the quantity which is 
associated to the $\edgecvcvv\in\edges_R^n$ interface is zero if $f$ 
is constant on the interval to which the values
$u^n_\cv$ and $u^n_\cvv$ belong; variations of the discrete unknowns
in this interval are therefore not taken into account.

\item The left hand side of  \eqref{BVfaible} involves  terms $(\vcvcvvn + 
\vcvvcvn)$ which are not uniformly bounded from below
by $C \ \medgecvcvv$ with some $C>0$ only depending on the data (that is 
$\vbf$, $u_0$ and $g$) and not on $\mesh$ (note that, for instance, 
$\vcvcvvn=\vcvvcvn=0$ if 
$\vbf \cdot \ncvcvv =0$).

\end{enumerate}

For the convergence result (namely Theorem \refp{eghbcewsol}) the useful 
consequence of \eqref{BVfaible}  is 

\[
h \sumTsk k \sumER 
(\vcvcvvn + \vcvvcvn) \abs{f(\ucvn) - f(\ucvvn)} \tends 0 \text{ as } h 
\tends 0,
\]

as in the 1D case, see Theorem \refp{nlcv1d}. For the error estimate in 
Theorem \refp{esth14}, the bound $C/\sqrt h$ in \eqref{BVfaible} is 
crucial. Note that a ``twice weak $BV$"  inequality in the sense (ii) and (iii), but 
of order 1 (that is $C$ instead of $C/\sqrt h$ in the right hand side of 
\eqref{BVfaible}), would yield a sharp error estimate, i.e. $C_e h^{1/2}$ 
instead of $C_e h^{1/4}$ in \refep{esteqcorol}. 

Note that, in order to obtain \eqref{BVfaible}, $\xi >0$ is crucial 
in the CFL condition \eqref{estcfl}.

Recall also that \eqref{BVfaible} together with the $L^\infty(\R^\dim 
\times \Rpast)$ bound does not yield any (strong) compactness property in 
$L^1_{loc}(\R^\dim \times \R_+)$ on a family of 
approximated solutions.

\medskip

In the linear case (that 
is $f(s)=cs$ for all $s \in \R$, for some $c$ in $\R$), 
the inequality  
 (\ref{BVfaible}) is used
in the same manner as in the previous chapter;
one proves that the approximate
solution satisfies the weak formulation to 
(\ref{estpbcont})-(\ref{estcondini}) (which is equivalent to 
\eqref{estkruzkov}) with an error 
which goes to $0$ as $h \rightarrow 0$, under condition \eqref{estcfl}. 
We deduce from this  the
convergence of $u_{\mesh,k}$ (as $h \tends 0$ and under condition 
\eqref{estcfl}) towards the unique weak solution of (\ref{estpbcont})-(\ref{estcondini})
in $L^\infty(\R^\dim \times \Rpast)$ for the weak-$\star$ topology. In fact, the convergence
holds in $L^p_{loc}(\R^\dim \times \R_+)$ (strongly) for any $1 \le p < 
\infty$, thanks to the argument developped for the study of the nonlinear case.

\smallskip

The nonlinear case adds an extra  difficulty,  as in the 1D case; it will 
be   handled in detail in
the present chapter. This difficulty arises from the fact that,
if $u_{\mesh,k}$ converges to $u$ (as $h \tends 0$, under condition 
\eqref{estcfl}) and $f(u_{\mesh,k})$ to $\mu_f$,
in $L^\infty(\R^\dim \times \Rpast)$ for the weak-$\star$ topology, there remains to show
that $\mu_f = f(u) $ and that $u$ is the entropy weak solution  to problem
(\ref{estpbcont})-(\ref{estcondini}). 
The weak $BV$ inequality 
(\ref{BVfaible}) is  used to show that, for any  ``entropy" function 
\index{Entropy!-- function} $\eta$, i.e. 
convex function of class $C^1$ from $\R$ to $\R$, with associated entropy 
flux $\phi$, i.e. $\phi$
such that $ \phi' = f' \eta'$, the following entropy inequality is satisfied: 

\be
\left\{\begin{array} {l} 
\displaystyle{ 
\int_{\R_+}
\int_{\R^\dim} \Bigl( \mu_\eta (x,t) \varphi_t (x,t) + \mu_\phi 
(x,t) \vbf(x,t) \cdot \grad  \varphi (x,t) \Bigr) dx dt + \int_{\R^\dim} \eta (u_0 (x)) \varphi (x,0) dx
\ge 0,  } \\  \forall \, \varphi \in C^\infty_c(\R^\dim \times
\R_+,\R_+),
\end{array}\right.\label{entropy}
\ee

where $\mu_\eta$ (resp. $\mu_\phi$) is the limit of 
$\eta(u_{\mesh,k})$ (resp.  $\phi(u_{\mesh,k})$) in $L^\infty(\R^\dim \times
\Rpast)$ for the weak-$\star$ topology (the existence of these limits can indeed be
assumed). From (\ref{entropy}), it is shown that $u_{\mesh,k}$
converges to $u$ in $L^1_{loc}(\R^\dim \times \R_+)$  (as $h \rightarrow 
0$, $k$ satisfying (\ref{estcfl})), and that $u$ 
is the entropy weak solution  to 
problem (\ref{estpbcont})-(\ref{estcondini}). This last result 
uses a generalization of a result on  measure valued  solutions of  DiPerna (see 
\citet{DP}, \citet{GH}), and is developped in section \refp{eghcvsol}.

\section{Stability results for the explicit scheme}\label{stabe}

\hpage{Hyperbolic multid}

\subsection{$L^\infty$ stability}\index{Stability}
\index{Estimate on the approximate solution!-- for a hyperbolic equation}
\begin{lemma}\label{estinftyest}  Under
Assumption \ref{estdonnees}, let $\mesh$ be an admissible mesh in the 
sense of Definition \ref{meshypmd} and $k>0$, let $g \in C(\R^2,\R)$ 
satisfy Assumption \ref{estcondflu} and assume that \eqref{estcfl} 
holds; let $\uapp$ be given by \eqref{estuapp}, \eqref{estschema}, \eqref{estschemaz}; then,  
\be   
 U_m \le \ucvn \le U_M, \, \forall n \in \N, \, \forall \cv \in
\mesh, \label{estinduction} 
\ee

and

\be
\Vert \uapp\Vert_{L^\infty(\R^\dim \times \Rpast)} \le \normuz.
\label{eststabinf}\ee

\end{lemma} 

\proof of Lemma \ref{estinftyest}

\smallskip

Note that \eqref{eststabinf} is a straightforward  consequence of
\eqref{estinduction}, which will be proved by induction. For $n =
0$,  since $U_m \le u_0 \le U_M$ a.e., \eqref{estinduction} follows from
\eqref{estschemaz}. 

Let $n \in \N$, assume that $U_m \le \ucvn \le U_M$ for all $\cv\in \mesh$.
Using the fact that  $\div\vbf=0$, which yields $\sumv
(\vcvcvvn - \vcvvcvn) = 0$,   we can rewrite \eqref{estschema} as

\be
 \mcv \frac {\ucvnp - \ucvn} k + \sumv  \Bigl( \vcvcvvn
(g(\ucvn,\ucvvn) - f(\ucvn)) -  \vcvvcvn (g(\ucvvn,\ucvn) - f(\ucvn))\Bigr) =
0.  \label{estschemab}
\ee

Set, for $\ucvn \not=\ucvvn$, 

\[
\tau_{\cv,\cvv}^n =  
\vcvcvvn \frac {g(\ucvn,\ucvvn) - f(\ucvn)} {\ucvn -\ucvvn} - 
\vcvvcvn \frac {g(\ucvvn,\ucvn) - f(\ucvn)} {\ucvn -\ucvvn},
\]

and $\tau_{\cv,\cvv}^n = 0$ if $\ucvn = \ucvvn$.

Assumption \ref{estcondflu} on  $g$ and Assumption \ref{estdonnees} yields
$0 \le \tau_{\cv,\cvv}^n \le V m(\edgecvcvv) (g_1 + g_2)$. Using \eqref{estschemab},
we can write

\[
\ucvnp = \Bigl(1 - \frac k {\mcv} \sumv \tau_{\cv,\cvv}^n \Bigr) \ucvn  
 + \frac k {\mcv} \sumv \tau_{\cv,\cvv}^n \ucvvn, 
\]

which gives, under condition \eqref{estcfl}, 
$ \dsp \inf_{\cvv\in\mesh} u_\cvv^n \le \ucvnp  \le \sup_{\cvv\in\mesh} u_\cvv^n$,
for all $\cv \in \mesh$.  This concludes the proof of \eqref{estinduction},
which, in turn, yields \eqref{eststabinf}. 

\begin{remark} \normalfont Note that the stability result
\eqref{eststabinf}  holds even if $\xi = 0 $ in  \eqref{estcfl}. However, we
shall need $\xi > 0$ for the following ``weak $BV$" inequality.
\end{remark}

\subsection{A ``weak $BV$" estimate}

\index{BV!Weak -- estimate|(}In the following 
lemma, $B(0,R)$ denotes the ball of $\R^\dim$ of center $0$ and radius $R$ 
($\R^\dim$ is always endowed with its usual scalar product).

\begin{lemma}\label{estwbvest}  Under
Assumption \ref{estdonnees}, let $\mesh$ be an admissible mesh in the 
sense of Definition \ref{meshypmd} and $k>0$. Let $g \in C(\R^2,\R)$ satisfy 
Assumption \ref{estcondflu} and assume that \eqref{estcfl} holds. Let 
$\uapp$ be given by \eqref{estuapp}, \eqref{estschema}, 
\eqref{estschemaz}.

Let $T>0$, $R>0$, $\NTk =
\max\{n\in\N, n < T/k\}$, $\mesh_R =
\{\cv\in\mesh, \cv \subset B(0,R)\}$ and 
$\edges_R^n = \{(\cv,\cvv)\in \mesh^2, \cvv\in \Nv,
 \edgecvcvv\subset B(0,R)$ and $\ucvn > \ucvvn \}$.
 
Then there exists 
$C \in \R$, only depending on $\vbf$, $g$, $u_0$, $\a$, $\xi$, $R$,
$T$ such that, for $h < R$ and $k < T$,

\be\begin{array}{r}
\sumTsk k \sumER \Bigl[
\vcvcvvn \Bigl(
\dsp{\max_{\ucvvn\le p \le q  \le\ucvn}} ( g(q,p) -  f(q)) + 
\dsp{\max_{\ucvvn\le p \le q  \le\ucvn}} ( g(q,p) -  f(p)) \Bigr) +\\ 
\vcvvcvn \Bigl(
\dsp{\max_{\ucvvn\le p \le q  \le\ucvn}} ( f(q) - g(p,q) ) + 
\dsp{\max_{\ucvvn\le p \le q  \le\ucvn}} ( f(p) - g(p,q) ) \Bigr) 
\Bigr] \\
\le \ \ \dsp{\frac {C} {\sqrt{h}}},
\end{array}\label{estbvfx}\ee

and 

\be
\sumTsk \summeshR \mcv \vert \ucvnp - \ucvn \vert \le \frac {C}
{\sqrt{h}}, \label{estbvft}\ee
\end{lemma}

\proof of Lemma \ref{estwbvest}

\smallskip

In this proof, we shall denote by $C_i$ ($i \in \N$) various quantities 
only depending on  $\vbf$, $g$, $u_0$, $\a$, $\xi$, $R$,
$T$. 

Multiplying \eqref{estschemab} by $k \ucvn$ and  summing the result over
$\cv\in \mesh_R$, $n\in\NN$  yields

\be
B_1 + B_2 = 0,
\label{est1bv}\ee

with

\[
B_1 = \sumTsk \summeshR \mcv \ucvn (\ucvnp - \ucvn),
\]

and

\[
B_2 = \sumTsk k \summeshR \sumv  \Bigl(
\vcvcvvn (g(\ucvn,\ucvvn) - f(\ucvn)) \ucvn - 
\vcvvcvn (g(\ucvvn,\ucvn) - f(\ucvn)) \ucvn \Bigr). 
\]

Gathering the last two summations by edges in $B_2$ leads to 
the definition of  $B_3$:

\[\ban
B_3 = \sumTsk k \sumER \Bigl[ &
\vcvcvvn \Bigl(\ucvn (g(\ucvn,\ucvvn) - f(\ucvn)) -
            \ucvvn (g(\ucvn,\ucvvn) - f(\ucvvn)) \Bigr) - \\& 
\vcvvcvn \Bigl(\ucvn (g(\ucvvn,\ucvn) - f(\ucvn)) -
            \ucvvn (g(\ucvvn,\ucvn) - f(\ucvvn)) \Bigr) \Bigr]. 
\ean\]

The expression $\vert B_3 - B_2 \vert$ can be reduced to a sum of
terms each of which corresponds to the boundary of a control volume which is included in $B(0,R+h)
\setminus B(0,R-h)$; since the measure of  $B(0,R+h)
\setminus B(0,R-h)$ is less than $C_2 h$,  
  the number of such terms is, for $n$ fixed, lower than $(C_2h)/(\a h^\dim)=C_3 
h^{1-\dim}$.  Thanks to
\eqref{eststabinf}, using the fact that $\meas(\dr\cv)  \le (1/\a) h^{\dim-1}$, 
that $\vert \vbf(x,t) \vert \le V$, that $g$ is bounded on $[U_m,U_M]^2$, 
and that $g(s,s)=f(s)$, one may   show that each of the non zero term in 
$\vert B_3 - B_2 \vert$  is bounded by $C_1  h^{\dim-1}$.  Furthermore, 
since $(\NTk +1)k \le 2k$, we deduce that 

\be
\vert B_3 - B_2 \vert \le C_{4}.
\label{est5bv}\ee

Denoting by $\Phi$ a primitive of the function $(\cdot) f'(\cdot)$, an 
integration by parts yields, for all $(a,b)\in\R^2$,

\be
\Phi(b) - \Phi(a) = \int_a^b s f'(s) ds = 
b(f(b)-g(a,b))-a(f(a)-g(a,b)) - \int_a^b(f(s)-g(a,b))ds.
\label{est6bv}\ee

Using \eqref{est6bv}, the term $B_3$ may be decomposed as
\[ 
B_3 = B_4 -B_5,
\]

where

\[
B_4 = \sumTsk k \sumER \left(\vcvcvvn \int_\ucvn^\ucvvn (f(s) -
g(\ucvn,\ucvvn)) ds + \vcvvcvn \int_\ucvvn^\ucvn (f(s) -
g(\ucvvn,\ucvn)) ds \right)
\]

and 

\[
B_5 =  \sumTsk k \sumER (\vcvcvvn -
\vcvvcvn)  \Bigl(\Phi(\ucvn) - \Phi(\ucvvn)\Bigr).
\]

The term $B_5$ is again reduced to a sum of terms corresponding to 
control volumes included in
$B(0,R+h) \setminus B(0,R-h)$, thanks to  $\div\vbf = 0$; therefore,
as for \eqref{est5bv}, there exists $C_5\in \R$ such
that 

\[
B_5 \le C_5.
\]

Let us now turn to an estimate of $B_4$. To this purpose, let $a,b
\in \R$,  define $\mathcal{C}(a,b) = \{(p,q) \in [a\bot
b,a\top b]^2$; $(q-p)(b-a)\ge 0\}$. Thanks to the
monotonicity properties of $g$ (and using the fact that 
$g(s,s)=f(s)$), the following inequality holds, for
any $(p,q) \in \mathcal{C}(a,b)$:

\be
\int_a^b(f(s)-g(a,b))ds \ge \int_c^d(f(s)-g(a,b))ds 
\ge \int_p^q(f(s)-g(p,q) )ds \ge 0.
\label{est7bv}\ee

The technical lemma \refp{petitlemme} can then be applied. 
It states that

\[
\vert \int_p^q (\theta(s)-\theta(p)) ds \vert \ge 
\frac 1 {2 G} (\theta(q) - \theta(p))^2, \ \ \forall p,q \in \R,
\]

for all monotone, Lipschitz continuous function  $\theta: \R \tends \R $, 
with a Lipschitz constant $G>0$.

From Lemma \ref{petitlemme}, we can notice that

\be
\int_p^q(f(s)-g(p,q) )ds \ge \int_p^q (g(p,s)-g(p,q) )ds \ge
\frac 1 {2 g_2} ( f(p) - g(p,q) )^2,
\label{est8bv}\ee

and

\be
\int_p^q (f(s)-g(p,q) )ds \ge \int_p^q (g(s,q)-g(p,q) )ds \ge
\frac 1 {2 g_1} ( f(q) - g(p,q) )^2.
\label{est9bv}\ee

Multiplying \eqref{est8bv} (resp. \eqref{est9bv}) by $g_2/(g_1 + g_2)$ (resp.
$g_1/(g_1 + g_2)$),
taking the maximum for $(p,q)  \in \mathcal{C}(a,b)$,  and adding the two
equations yields, with \eqref{est7bv},

\be
\int_a^b (f(s)-g(a,b))ds \ge 
\frac 1 {2 (g_1+g_2) } \Bigl( 
\dsp{\max_{(p,q) \in \mathcal{C}(a,b)}} ( f(p) - g(p,q) )^2 +
\dsp{\max_{(p,q) \in \mathcal{C}(a,b)}} ( f(q) - g(p,q) )^2 \Bigr).
\label{est9bvbis}\ee

We can then deduce, from \eqref{est9bvbis}:

\be\ban
B_4 \ge & \dsp{\frac 1 {2(g_1 + g_2)}} \sumTsk k \sumER \Bigl[ \\
&
\vcvcvvn \Bigl(
\dsp{\max_{\ucvvn\le p \le q  \le\ucvn}} ( g(q,p) -  f(q))^2 + 
\dsp{\max_{\ucvvn\le p \le q  \le\ucvn}} ( g(q,p) -  f(p))^2 \Bigr) +\\& 
\vcvvcvn \Bigl(
\dsp{\max_{\ucvvn\le p \le q  \le\ucvn}} ( f(q) - g(p,q) )^2 + 
\dsp{\max_{\ucvvn\le p \le q  \le\ucvn}} ( f(p) - g(p,q) )^2 \Bigr) 
\Bigr].
\ean\label{est10bv}\ee

This gives a bound on $B_2$, since (with $C_6=C_4+C_5$):

\be
B_2 \ge B_4 - C_6.
\label{est10bva}
\ee

Let us now turn to $B_1$. We have 

\be
B_1 = - \half \sumTsk \summeshR \mcv (\ucvnp - \ucvn)^2 + 
\half \summeshR \mcv \Bigl(u_\cv^{\NTk+1}\Bigr)^2 - 
\half \summeshR \mcv \Bigl(u_\cv^{0}\Bigr)^2.
\label{est11bv}\ee

Using \eqref{estschemab} and  the Cauchy-Schwarz inequality yields the
following inequality:

\[
\bal
(\ucvnp - \ucvn)^2 \le 
\\
{\dsp \frac {k^2} {\mcv^2}} \sumv (\vcvcvvn+\vcvvcvn)
\sumv
\Bigl[ \vcvcvvn \Bigl(g(\ucvn,\ucvvn) - f(\ucvn)\Bigr)^2 + 
\vcvvcvn \Bigl(g(\ucvvn,\ucvn) - f(\ucvn)\Bigr)^2\Bigr]. 
\ean
\]

Then, using the CFL condition \eqref{estcfl}, Definition \ref{meshypmd} and 
part (iv) of Assumption \ref{estdonnees} gives

\be\bal
\mcv (\ucvnp - \ucvn)^2 \le \\
k \dsp{\frac {1-\xi} {g_1 + g_2}}\sumv
\Bigl[ \vcvcvvn \Bigl(g(\ucvn,\ucvvn) - f(\ucvn)\Bigr)^2 + 
\vcvvcvn \Bigl(g(\ucvvn,\ucvn) - f(\ucvn)\Bigr)^2\Bigr]. 
\ean\label{est13bv}\ee

Summing equation \eqref{est13bv} over $\cv\in \mesh_R$ and over $n
=0,\ldots,\NTk$, and reordering the summation leads to 

\be\bac
\dsp{\half} \sumTsk \summeshR \mcv (\ucvnp - \ucvn)^2 \le   
  \dsp{\frac {1-\xi} {2(g_1 + g_2)}} \sumTsk k \sumER \Bigl[ \\
\vcvcvvn \Bigl(
(g(\ucvn,\ucvvn) - f(\ucvn))^2 + 
(g(\ucvn,\ucvvn) - f(\ucvvn))^2 \Bigr) +\\ 
\vcvvcvn \Bigl(
(f(\ucvn) - g(\ucvvn,\ucvn))^2 + 
(f(\ucvvn) - g(\ucvvn,\ucvn))^2 \Bigr) 
\Bigr]+ C_7,
\ean\label{est13bvbis}\ee

where $C_7$ accounts for the interfaces $\edgecvcvv \subset B(0,R)$ such 
that $\cv \notin \mesh_R$ and/or $\cvv \notin \mesh_R$ (these control 
volumes are included in
$B(0,R+h) \setminus B(0,R-h)$).

Note that the right hand side of \eqref{est13bvbis} is bounded by $(1-\xi) B_4 + C_7$
(from \eqref{est10bv}).  Using \eqref{est1bv}, \eqref{est10bva} and  
\eqref{est11bv} gives

\be\ban
\dsp{\frac \xi {2(g_1 + g_2)}} \sumTsk k \sumER \Bigl[& 
\vcvcvvn \Bigl(
\dsp{\max_{\ucvvn\le p \le q  \le\ucvn}} ( g(q,p) -  f(q))^2 + 
\dsp{\max_{\ucvvn\le p \le q  \le\ucvn}} ( g(q,p) -  f(p))^2 \Bigr) +\\ 
&\vcvvcvn \Bigl(
\dsp{\max_{\ucvvn\le p \le q  \le\ucvn}} ( f(q) - g(p,q) )^2 + 
\dsp{\max_{\ucvvn\le p \le q  \le\ucvn}} ( f(p) - g(p,q) )^2 \Bigr) 
\Bigr] \\&\le \ \ 
\dsp{\half} \summeshR \mcv \Bigl(u_\cv^{0}\Bigr)^2 + C_6 + C_7= C_8.       
\ean\label{est14bv}\ee

Applying  the Cauchy-Schwarz inequality to the left hand side of
\eqref{estbvfx} and using \eqref{est14bv} yields

\be\ban
\sumTsk k \sumER \Bigl[&
\vcvcvvn \Bigl(
\dsp{\max_{\ucvvn\le p \le q  \le\ucvn}} ( g(q,p) -  f(q)) + 
\dsp{\max_{\ucvvn\le p \le q  \le\ucvn}} ( g(q,p) -  f(p)) \Bigr) +\\ 
&\vcvvcvn \Bigl(
\dsp{\max_{\ucvvn\le p \le q  \le\ucvn}} ( f(q) - g(p,q) ) + 
\dsp{\max_{\ucvvn\le p \le q  \le\ucvn}} ( f(p) - g(p,q) ) \Bigr) 
\Bigr] \\
&\le \ \ C_9 \Bigr(\sumTsk k \sumER (\vcvcvvn + \vcvvcvn)\Bigl)^{\frac 1 2}.
\ean\label{est15bv}\ee

Noting that

\[
\sumER (\vcvcvvn + \vcvvcvn) \le \sum_{\cv \in \mesh_{R+h}} V \meas(\dr \cv) \le V \frac 1 
{\a} h^{\dim -1} \frac {\meas(B(0,R+h))} {\a h^\dim}=\frac {C_{10}} h
\]

and $(\NTk+1)k \le 2T$, one obtains \eqref{estbvfx} from
\eqref{est15bv}. 

\medskip

Finally, since \eqref{estschemab} yields

\[
\mcv \vert\ucvnp - \ucvn\vert \le k \sumv 
\Bigl( \vcvcvvn \vert g(\ucvn,\ucvvn) - f(\ucvn) \vert + 
\vcvvcvn \vert g(\ucvvn,\ucvn) - f(\ucvn) \vert \Bigr), 
\]

Inequality  \eqref{estbvft} immediately  follows from \eqref{estbvfx}.
This completes the proof of Lemma \ref{estwbvest}.
\eof

\section{Existence of the solution and stability results for the
implicit scheme}\label{stabi} \index{Scheme!Implicit -- for 
hyperbolic equations|(}

\hpage{Hyperbolic multid}

This section is devoted to the time implicit scheme (given by
\eqref{estschemai} and \eqref{estschemaz}). 
We first prove the existence and uniqueness of the solution $\{\ucvn, n \in \N,
\cv \in \mesh\}$ of \eqref{estschemaz}, \eqref{estschemai} and 
such that $\ucvn \in [U_m,U_M]$ for all $\cv\in \mesh$ and all $n \in \N$. 
Then, one gives a ``weak space $BV$" inequality (this is equivalent to 
the inequality \eqref{estbvfx} for the explicit scheme)  and a
``(strong) time $BV$" estimate (Estimate \eqref{estbvt} below). 
This last estimate requires that $\vbf$ does not depend on $t$
(and it leads to the term ``$k$" in the right hand side of 
\eqref{esteqcoroli} in Theorem \ref{esth14i}). The error estimate, in the case where
$\vbf$ depends on $t$, is given in Remark \ref{estexti}.

\index{BV!Weak -- estimate|)}
\subsection{Existence, uniqueness and $L^{\infty}$ stability}\index{Stability}
\index{Estimate on the approximate solution!-- for a hyperbolic equation}
The following proposition gives an existence and uniqueness result of the
solution to \eqref{estschemaz}, \eqref{estschemai}. In this proposition, $\vbf$
may depend on $t$ and one does not need to assume 
$u_0 \in BV(\R^\dim)$.

\begin{prop}\label{esteu}
Under Assumption \ref{estdonnees}, let $\mesh$ be an admissible mesh in the 
sense of Definition \ref{meshypmd} and $k>0$. Let $g \in C(\R^2,\R)$ satisfy 
Assumption \ref{estcondflu}.

Then there
exists a unique solution $\{\ucvn,$ $n \in \N$, 
$\cv\in \mesh\}$ $\subset [U_m,U_M]$ to 
\eqref{estschemaz},\eqref{estschemai}.
\end{prop}

\proof  of Proposition \ref{esteu}

\smallskip

One proves Proposition \ref{esteu} by induction. Indeed, $\{\ucvz,$ $\cv\in
\mesh\}$ is uniquely defined by \eqref{estschemaz} and one has $u^0_\cv \in
[U_m,U_M]$, for all $\cv\in \mesh$, since $U_m \le u_0 \le U_M$ a.e..
Assuming that, for some $n \in \N$, the set $\{\ucvn,$ $\cv\in \mesh\}$
is given and that $u^n_\cv \in
[U_m,U_M]$, for all $\cv\in \mesh$,  the existence and
uniqueness of $\{\ucvnp, \, \cv \in \mesh\}$, such that $\ucvnp \in 
[U_m,U_M]$ for all $\cv \in \mesh$, solution
of \eqref{estschemai}, must be shown.

\medskip

\textit{Step 1 (uniqueness of $\{\ucvnp, \, \cv \in \mesh \}$, such that $\ucvnp \in 
[U_m,U_M]$ for all $\cv \in \mesh$, solution of \eqref{estschemai})}

Recall that $n\in\N$ and $\{\ucvn,$ $\cv\in \mesh\}$
are given. Let us consider two solutions of \eqref{estschemai}, 
respectively denoted by $\{\ucv,$ $\cv\in \mesh\}$ and 
$\{\wcv,$ $\cv\in \mesh\}$; therefore,  $\{\ucv,$ $\cv\in \mesh\}$ and 
$\{\wcv,$ $\cv\in \mesh\}$ satisfy $\{\ucv,$ $\cv\in \mesh\} \subset 
[U_m,U_M]$, $\{\wcv,$ $\cv\in \mesh\} \subset 
[U_m,U_M]$, 

\be
\mcv \frac {\ucv - \ucvn} k + \sumv 
(\vcvcvvn \ g(\ucv,\ucvv) - \vcvvcvn \ g(\ucvv,\ucv)) = 0, \, \forall
\cv \in\mesh,  \label{estschemaiut}\ee

and

\be
\mcv \frac {\wcv - \ucvn} k + \sumv 
(\vcvcvvn \ g(\wcv,\wcvv) - \vcvvcvn \ g(\wcvv,\wcv)) = 0, \, \forall
\cv \in\mesh.  \label{estschemaiwt}\ee

Then, substracting (\ref{estschemaiwt}) to  (\ref{estschemaiut}), for all
$\cv \in \mesh$,

\be\bal
{\dsp \frac {\mcv} k }( \ucv - \wcv ) + \sumv \vcvcvvn (g(\ucv,\ucvv) - g(\wcv,\ucvv)) 
\\
+ \sumv \vcvcvvn ( g(\wcv,\ucvv) - g(\wcv,\wcvv)) 
- \sumv \vcvvcvn ( g(\ucvv,\ucv) - g(\wcvv,\ucv)) \\
- \sumv \vcvvcvn ( g(\wcvv,\ucv) - g(\wcvv,\wcv))=0
\ean
\label{estschemaiu1t}
\ee

thanks to the monotonicity properties of $g$, \eqref{estschemaiu1t} leads to

\be\bal
{\dsp \frac {\mcv} k} \vert \ucv - \wcv \vert + \sumv \vcvcvvn \vert g(\ucv,\ucvv) -
g(\wcv,\ucvv) \vert  \\
+ \sumv \vcvvcvn \vert g(\wcvv,\ucv) - g(\wcvv,\wcv) \vert 
\le \sumv \vcvcvvn \vert g(\wcv,\ucvv) - g(\wcv,\wcvv) \vert \\
+ \sumv \vcvvcvn \vert
g(\ucvv,\ucv) - g(\wcvv,\ucv) \vert. 
\ean
\label{estschemaiu2t}
\ee

Let $\vphi: \R^\dim \mapsto \R_+^\star$ be defined by $\vphi(x)=\exp(-\gamma \vert x
\vert)$, for some positive $\gamma$ which will be specified later. 
For $\cv\in \mesh$, let $\vphi_\cv$ be the
mean value of $\vphi$ on $\cv$. Since
$\vphi$ is integrable over
$\R^\dim$ (and thanks to \eqref{estmaillage}), one has $\sum_{\cv \in 
\mesh}\vphi_\cv$ $\le (1/(\a h^\dim)) \norm{\varphi}{L^1(\R^\dim)}$ $<
\infty$. Therefore the series

\[
\sum_{\cv \in \mesh}\vphi_\cv(\sum_{\cvv \in
\Nv} \vcvcvvn \vert g(\wcv,\ucvv) - g(\wcv,\wcvv) \vert ) \text{ and } \sum_{\cv \in
\mesh}\vphi_\cv (\sum_{\cvv \in
\Nv} \vcvvcvn \vert g(\ucvv,\ucv) - g(\wcvv,\ucv)
\vert)
\]

are convergent (thanks to \eqref{estmaillage} and the boundedness of
$\vbf$ on $\R^\dim$ and $g$ on $[U_m,U_M]^2$).

Multiplying \eqref{estschemaiu2t} by $\vphi_\cv$ and summing for $\cv \in \mesh$
yields five convergent series which can be reordered in order to give

\[
\begin{array}{l}
\dsp{\sum_{\cv \in \mesh}} \frac {\mcv} k \vert \ucv - \wcv \vert \vphi_\cv
\le 
 \dsp{\sum_{\cv \in \mesh}} \sumv \vcvcvvn \vert g(\ucv,\ucvv) - g(\wcv,\ucvv)
\vert \vert
\vphi_\cv-\vphi_\cvv
\vert \\
+\dsp{\sum_{\cv \in \mesh}} \sumv \vcvvcvn
\vert g(\wcvv,\ucv) - g(\wcvv,\wcv) \vert \vert \vphi_\cv-\vphi_\cvv \vert, 
\end{array}
\]

from which one deduces

\be
\sum_{\cv \in \mesh} a_\cv \vert \ucv - \wcv \vert \le 
\sum_{\cv \in \mesh} b_\cv \vert \ucv - \wcv \vert,
\label{estschemaiu4t}
\ee

with, for all $\cv\in \mesh$, $a_\cv= \frac {\mcv} k \vphi_\cv$ and 
$b_\cv= \sumv(\vcvcvvn g_1 + \vcvvcvn g_2)\vert \vphi_\cv-\vphi_\cvv \vert$.

For $\cv \in \mesh$, let $\centercv$ be an arbitrary point of $\cv$.
Then,

\[
a_\cv \ge \frac 1 k \a h^\dim \inf \{\vphi(x), \, x \in B(\centercv,h)\}
\]

and

\[
b_\cv 
\le \frac {2V(g_1+g_2)}{\a} h^\dim \sup \{ \vert 
\grad \vphi(x) \vert, \, x \in B(\centercv,2h)\}.
\]

Therefore, taking $\gamma>0$ small enough in order to have

\be
\inf \{\vphi(y), \, y \in B(x,h)\} > C \sup \{\vert 
\grad \vphi(y) \vert, \,  y \in B(x,2h)\}, \, \forall x \in
\R^\dim
\label{choixg}
\ee

with $C= (2 k V (g_1 + g_2))/\a^2$, yields $a_\cv > b_\cv$ for all $\cv
\in \mesh$. Hence \eqref{estschemaiu4t} gives $\ucv=\wcv$, for all $\cv\in 
\mesh$.

A choice of $\gamma >0$ verifying \eqref{choixg} is always possible. 
Indeed, since $\vert \grad \vphi(z) \vert = \gamma \exp(-\gamma \vert z 
\vert)$, taking $\gamma >0$ such that $\gamma 
\exp(3\gamma h) < 1/C$ is convenient.

This concludes Step 1.

\medskip

\textit{Step 2 (existence of $\{\ucvnp, \, \cv \in \mesh\}$, such that 
$\ucvnp \in [U_m,U_M]$ for all $\cv \in \mesh$, solution of 
\eqref{estschemai})}.

Recall that $n\in\N$ and $\{\ucvn,$ $\cv\in \mesh\}$
are given. For $r \in \N^\star$, let $B_r=B(0,r)=\{x \in \R^\dim$, $\vert 
x \vert < r\}$ and $\mesh_r=\{\cv \in \mesh,$ $\cv\subset B_r\}$ (as in 
Lemma \ref{estwbvest}). Let us 
assume that $r$
is large enough, say $r \ge r_0$, in order to have $\mesh_r \ne \emptyset$. 

\smallskip

If $\cv\in
\mesh \setminus \mesh_r$, set $\ucvr=\ucvn$. Let us first prove that 
there exists $\{\ucvr$, $\cv \in \mesh_r\} \subset [U_m,U_M]$,
solution to 

\be
\mcv \frac {\ucvr - \ucvn} k + \sumv 
(\vcvcvvn \ g(\ucvr,\ucvvr) - \vcvvcvn \ g(\ucvvr,\ucvr)) = 0, \, \forall
\cv \in \mesh_r.
\label{estschemair}
\ee

Then, we will prove that passing to the limit as $r \tends
\infty$ (up to a subsequence) leads to a solution 
$\{\ucvnp, \, \cv \in\mesh \}$ to \eqref{estschemai} such that $\ucvnp \in [U_m,U_M]$ for all 
$\cv \in \mesh$.

\smallskip

For a fixed $r \ge r_0$, in order to prove the existence of  $\{\ucvr$, 
$\cv \in \mesh_r\} \subset [U_m,U_M]$
solution to \eqref{estschemair}, a ``topological degree" argument is used 
(see, for instance, \citet{deimling} for a presentation of the degree).

Let $U_r^n=\{ \ucvn,$ $\cv\in \mesh_r \}$ and assume that $U_r= \{\ucvr,$ 
$\cv\in \mesh_r\}$ is a solution of
\eqref{estschemair}. The families $U_r$
and $U_r^n$ may be viewed as vectors of $\R^N$, with $N=\card(\mesh_r)$.
Equation \eqref{estschemair} gives

\[
  \ucvr + \frac k {\mcv} \sumv 
(\vcvcvvn \ g(\ucvr,\ucvvr) - \vcvvcvn \ g(\ucvvr,\ucvr)) =  \ucvn , \, \forall
\cv \in \mesh_r,
\]

which can be written on the form

\be
U_r-G_r(U_r)=U_r^n,
\label{estschemairb}
\ee

where $G_r$ is a continuous map from $\R^N$ into $\R^N$. 

One may assume that $g$ is nondecreasing with respect to its first 
argument and nonincreasing with respect to its second argument on $\R^2$ 
(indeed, thanks to the monotonicity properties of $g$ given by 
Assumption \ref{estcondflu}, it is sufficient to change, if necessary, 
$g$ on $\R^2 \setminus [U_m,U_M]^2$, setting, for instance, $g(a,b)=g(U_m 
\top (U_M \bot a), U_m \top (U_M \bot b))$). Then, since $\ucvn \in
[U_m,U_M]$, for all $\cv\in \mesh$, and $\ucvr = \ucvn \in [U_m,U_M]$, for all $\cv\in
\mesh \setminus \mesh_r$, it is easy to show (using $\div(\vbf)=0$) that
if $U_r$ satisfies \eqref{estschemairb}, then one has $\ucvr \in [U_m,U_M]$, for
all $\cv\in \mesh_r$. Therefore, if $\mathcal{C}_r$ is a ball of $\R^N$ of center
$0$ and of radius great enough, Equation \eqref{estschemairb} has no solution on the
boundary of $\mathcal{C}_r$, and one can define the topological degree of the
application
$Id-G_r$ associated to the set $\mathcal{C}_r$ and to the point $U_r^n$, that is 
$deg(Id-G_r,\mathcal{C}_r,U_r^n)$.
Furthermore, if $\l \in [0,1]$, the same argument allows us to define 
$deg(Id-\l G_r,\mathcal{C}_r,U_r^n)$. Then, the property of invariance of the 
degree by continuous transformation asserts that $deg(Id-\l G_r,\mathcal{C}_r,U_r^n)$ does not depend on $\l \in [0,1]$. This gives

\[deg(Id-G_r,\mathcal{C}_r,U_r^n)=deg(Id,\mathcal{C}_r,U_r^n).\] 

But, since $U_r^n \in \mathcal{C}_r$, 

\[
deg(Id,\mathcal{C}_r,U_r^n)=1.
\] 

Hence

\[
deg(Id-G_r,\mathcal{C}_r,U_r^n)
\ne 0.
\]

This proves that there exists a solution $U_r\in \mathcal{C}_r$ to
\eqref{estschemairb}. Recall also that we already proved that the 
components of $U_r$ are necessarily 
in $[U_m,U_M]$.

\medskip

In order to prove the existence of $\{\ucvnp,$ $\cv\in \mesh\} \subset 
[U_m,U_M]$ solution to \eqref{estschemai}, let us pass to the limit as $r 
\tends \infty$. For $r \ge r_0$, let $\{\ucvr,$ $\cv\in \mesh\}$ be a solution
of \eqref{estschemair} (given by the previous proof). Since
$\{\ucvr,$ $r \in \N\}$ is included in $[U_m,U_M]$, for all $\cv
\in \mesh$, one can find (using a ``diagonal process") a
sequence $(r_l)_{l \in \N}$, with $r_l \tends \infty$, as $l \tends
\infty$, such that $(u^{r_l}_\cv)_{l \in \N}$ converges (in
$[U_m,U_M]$) for all $\cv
\in \mesh$. One sets $\ucvnp=\lim_{l \tends \infty}u^{r_l}_\cv$. Passing
to the limit in \eqref{estschemair} (this is possible since for all $\cv 
\in \mesh$, this equation is satisfied
for all $l\in\N$ large enough) shows that $\{\ucvnp,$ $\cv\in \mesh\}$
is solution to \eqref{estschemai}.

Indeed, using the uniqueness of the solution of \eqref{estschemai}, one can
show that $\ucvr \tends \ucvnp$, as $r \tends \infty$, for all $\cv\in
\mesh$.

This completes the proof of Proposition \ref{esteu}. \eof

\subsection{``Weak space $BV$" inequality}

\index{BV!Weak -- estimate|(}One gives here an inequality similar to Inequality
\eqref{estbvfx} (proved for the explicit scheme). This inequality does not make use of 
$u_0 \in BV(\R^\dim)$ and
$\vbf$ can depend on $t$. Inequality \eqref{estbvft} also holds but is 
improved in Lemma \ref{estlbvt} when $u_0 \in BV(\R^\dim)$ and
$\vbf$ does not depend on $t$.

\begin{lemma}\label{estwbvesti}  Under
Assumption \ref{estdonnees}, let $\mesh$ be an admissible mesh in the 
sense of Definition \ref{meshypmd} and $k>0$. Let $g \in C(\R^2,\R)$ satisfy 
Assumption \ref{estcondflu} and let $\{\ucvn, n \in \N, \cv \in \mesh \}$ be the solution of
\eqref{estschemai}, \eqref{estschemaz} such that
$\ucvnp \in [U_m,U_M]$ for all $\cv \in \mesh$ and all $\nnn$ (existence and uniqueness of such a solution is
given by Proposition
\ref{esteu}).

Let $T>0$, $R>0$, $\NTk =
\max\{n\in\N, n < T/k\}$, $\mesh_R =
\{\cv\in\mesh, \cv \subset B(0,R)\}$ and 
$\edges_R^n = \{(\cv,\cvv)\in \mesh^2, \cvv\in \Nv, \edgecvcvv \subset B(0,R)$ and $\ucvn > \ucvvn \}$.
 
Then there exists 
$C_{v} \in \R$, only depending on $\vbf$, $g$, $u_0$, $\a$, $R$,
$T$ such that, for $h < R$ and $k < T$,

\be\begin{array}{r}
\sumTsk k \sumERi \Bigl[
\vcvcvvn \Bigl(
\dsp{\max_{\ucvvnp\le p \le q  \le\ucvnp}} ( g(q,p) -  f(q)) + 
\dsp{\max_{\ucvvnp\le p \le q  \le\ucvnp}} ( g(q,p) -  f(p)) \Bigr) +\\
\vcvvcvn \Bigl(
\dsp{\max_{\ucvvnp\le p \le q  \le\ucvnp}} ( f(q) - g(p,q) ) + 
\dsp{\max_{\ucvvnp\le p \le q  \le\ucvnp}} ( f(p) - g(p,q) ) \Bigr) 
\Bigr] \\
\le \ \ \dsp{\frac {C_{v}} {\sqrt{h}}}.
\end{array}\label{estbvfxi}\ee

Furthermore, Inequality \refp{estbvft} holds.

\end{lemma}

\proof of Lemma \ref{estwbvesti}

\smallskip

We multiply \eqref{estschemai} by $k\ucvnp$, and sum the result
over $\cv \in \mesh_R$ and $n\in\NN$. We can then follow,
step by step, the proof of Lemma \refp{estwbvest} until 
Equation \eqref{est11bv} in which the first term of the right hand side
appears with the opposite sign. We can then directly conclude
an inequality similar to \eqref{est14bv}, which is sufficient to conclude 
the proof of Inequality \eqref{estbvfxi}. Inequality \refp{estbvft} 
follows easily from \eqref{estbvfxi}.
\eof

\index{BV!Weak -- estimate|)}

\subsection{``Time $BV$" estimate}
\index{BV!Strong time -- estimate|(}
\index{Estimate on the approximate solution!-- for a hyperbolic equation}
This section gives a so called ``strong time $BV$ estimate" (estimate 
\eqref{estbvt}). For this estimate, the fact that
$u_0 \in BV(\R^\dim)$ and that $\vbf$ does not depend on $t$ is required. 
Let us begin this section with a preliminary lemma on the space $BV(\R^d)$.

\begin{lemma}\label{divbv}
Let $\mesh$ be an admissible mesh in the 
sense of Definition \refp{meshypmd} and let $u \in BV(\R^\dim)$ (see 
Definition \refp{defBVRp}). For $\cv \in \mesh$, let $u_\cv$ be the mean 
value of $u$ over $\cv$. Then,

\be
\sum_{\edgecvcvv \in \edges} \medgecvcvv \vert u_\cv - u_\cvv \vert \le 
\frac C {\a^4} \snorm{u}{BV(\R^d)},
\label{bvzd}
\ee

where $C$ only depends on the space dimension ($d=1$, $2$ or $3$).
\end{lemma}

\proof of Lemma \ref{divbv}

\smallskip

Lemma \ref{divbv} is proven in two steps. In the first step, it is proved 
that if \eqref{bvzd} holds for all $u \in BV(\R^d)\cap C^1(\R^d,\R)$ then 
\eqref{bvzd} holds  
for all $u \in BV(\R^d)$. In Step 2, \eqref{bvzd} is 
proved to hold for $u \in BV(\R^d) \cap C^1(\R^d,\R)$.

\smallskip

\textit{Step 1 (passing from $BV(\R^d)\cap C^1(\R^d,\R)$ to $BV(\R^d)$)}

Recall that $BV(\R^d) \subset L^1_{loc}(\R^d)$. Let $u \in BV(\R^d)$, let 
us regularize $u$ by a sequence of mollifiers.

Let $\rho \in C^\infty_c(\R^d,\R_+)$ such that $\int_{\R^d} \rho(x)dx=1$. 
Define, for all $n \in \N^\star$, $\rho_n$ by $\rho_n(x)=n^\dim \rho 
(nx)$ for all $x \in \R^d$ and $u_n=u \star \rho_n$, that is

\[
u_n(x)=\int_{\R^d} u(y) \rho_n(x-y)dy, \, \forall x \in \R^d.
\]

It is well known that $(u_n)_{n \in \N^\star}$ is included in $C^\infty(\R^d,\R)$ 
and converges to $u$ in $L^1_{loc}(\R^d)$ as $n \tends \infty$. Then, the mean value of $u_n$ over $\cv$ 
converges, as $n \tends \infty$, to $u_\cv$, for all $\cv \in \mesh$.
Hence, if \eqref{bvzd} holds with $u_n$ instead of $u$ (this will be proven 
in Step 2) and if $\snorm{u_n}{BV(\R^d)} \le 
\snorm{u}{BV(\R^d)}$ for all $n \in \N^\star$, Inequality \eqref{bvzd} is 
proved by passing to the limit as $n \tends \infty$.

\smallskip

In order to prove $\snorm{u_n}{BV(\R^d)} \le 
\snorm{u}{BV(\R^d)}$ for all $n \in \N^\star$ (this will conclude step 1), 
let $n \in \N^\star$ and  $\vphi \in 
C^\infty_c(\R^d,\R^d)$ such that $\vert \vphi(x) \vert \le 1$ for all $x 
\in \R^d$. A simple computation gives, using Fubini's theorem,

\be
\int_{\R^d} u_n(x)\div \vphi(x)dx =\int_{\R^d} \bigl(\int_{\R^d} u(x-y)\div \vphi 
(x) dx \bigr) \rho_n(y) dy \le \snorm{u}{BV(\R^d)},
\label{bba}
\ee

since, setting $\psi_y=\vphi(y+\cdot) \in C^\infty_c(\R^d,\R^d)$ (for all $y 
\in \R^d$),

\[
\int_{\R^d} u(x-y)\div \vphi 
(x) dx=\int_{\R^d} u(z)\div \psi_y
(z) dz \le \snorm{u}{BV(\R^d)}, \, \forall y \in \R^d,
\]

and

\[
\int_{\R^d} \rho_n(y)dy =1.
\]

Then, taking in \eqref{bba} the supremum over $\vphi \in 
C^\infty_c(\R^d,\R^d)$ such that $\vert \vphi(x) \vert \le 1$ for all $x 
\in \R^d$ leads to $\snorm{u_n}{BV(\R^d)} \le 
\snorm{u}{BV(\R^d)}$.

\smallskip

\textit{Step 2 (proving \eqref{bvzd} if $u \in BV(\R^d)\cap C^1(\R^d,\R)$)}

Recall that $B(x,R)$ denotes the ball of $\R^d$ of center $x$ and radius 
$R$.
Since $u \in C^1(\R^d,\R)$,

\[
\int_{\R^d} u(x)\div \vphi(x)dx= - \int_{\R^d} \grad u(x) \cdot \vphi(x)dx.
\]

Then $\snorm{u}{BV(\R^d)}=\norm{ (\vert \grad u \vert) }{L^1(\R^d)}$ and we 
will prove \eqref{bvzd} with $\norm{ (\vert \grad u \vert) }{L^1(\R^d)}$ 
instead of $\snorm{u}{BV(\R^d)}$.

\smallskip

Let $\edgecvcvv \in \edges$, then $\cv \in \mesh$, $\cvv \in \Nv$ and

\[
u_\cv - u_\cvv = \frac 1 {\mcv \mcvv} \int_\cvv \int_\cv (u(x) - 
u(y))dxdy.
\]

For all $x \in \cv$ and all $y \in \cvv$,

\[
u(x)-u(y)=\int_0^1 \grad u (y + t(x-y)) \cdot (x-y) dt.
\]

Then,

\[
\bari
\dsp{\mcv \mcvv \vert u_\cv - u_\cvv \vert \le \int_\cvv (\int_\cv \int_0^1 \vert \grad u 
(y + t(x-y)) \vert \vert x-y \vert dt dx ) dy}
\\
\dsp{\le \int_\cvv (\int_0^1 \int_\cv \vert \grad u 
(y + t(x-y)) \vert \vert x-y \vert dx dt ) dy.}
\ean
\]

Using $\vert x-y \vert \le 2h$ and changing the variable $x$ in $z=x-y$ 
(for all fixed $y \in \cvv$ and $t \in (0,1)$) yields 

\[
\mcv \mcvv \vert u_\cv - u_\cvv \vert \le 2h \int_\cvv (\int_0^1 
\int_{B(0,2h)} \vert \grad u 
(y + tz) \vert dz dt) dy,
\]

which may also  be written (using Fubini's theorem)

\be
\mcv \mcvv \vert u_\cv - u_\cvv \vert \le 2h \int_{B(0,2h)} (\int_0^1 \int_\cvv \vert \grad u 
(y + tz) \vert dy dt ) dz.
\label{oufa}
\ee

For all $\cv \in \mesh$, let $\centercv$ be an arbitrary point of $\cv$.

Then, changing the variable $y$ in $\xi=y+tz$ (for all 
fixed $z \in \cvv$ and $t \in (0,1)$) in \eqref{oufa},

\[
\mcv \mcvv \vert u_\cv - u_\cvv \vert \le 2h \int_{B(0,2h)} (\int_0^1 
\int_{B(x_\cvv,3h)} \vert \grad u 
(\xi) \vert d\xi dt ) dz,
\]

which yields, since $\mesh$ is an adimissible mesh in the sense of 
 Definition \refp{estmaillage},

\[
\medgecvcvv \vert u_\cv - u_\cvv \vert \le \frac {2h^{\dim}} {\a^3 
h^{2d}} \meas(B(0,2h)) \int_{B(x_\cvv,3h)} \vert \grad u 
(\xi) \vert d\xi.
\]

Therefore there exists $C_1$, only depending on the space dimension, such 
that

\be
\medgecvcvv \vert u_\cv - u_\cvv \vert \le \frac {C_1} {\a^3}\int_{B(x_\cvv,3h)} \vert \grad u 
(\xi) \vert d\xi, \, \forall \edgecvcvv \in \edges.
\label{oufb}
\ee

Let us now remark that, if $\cvw \in \mesh$ and $\cvv \in \mesh$, $\cvw \cap B(x_\cvv,3h) 
\ne \emptyset$ implies $\cvv \subset B(\centercvw,5h)$. Then, for a fixed 
$\cvw \in \mesh$, the number of $\cvv \in \mesh$ such that $\cvw \cap B(x_\cvv,3h) 
\ne \emptyset$ is less or equal to $\meas(B(0,5h))/(\a h^d)$ that is less or 
equal $C_2/\a$ where $C_2$ only depends on the space dimension.

Then, summing \eqref{oufb} over $\edgecvcvv \in \edges$ leads to

\[
\sum_{\edgecvcvv \in \edges} \medgecvcvv \vert u_\cv - u_\cvv \vert \le 
\frac {C_1C_2} {\a^4} \sum_{\cvw \in \mesh} \int_\cvw \vert \grad u 
(\xi) \vert d\xi= \frac {C_1C_2} {\a^4} \norm{(\vert \grad u 
\vert)}{L^1(\R^d)},
\]

that is \eqref{bvzd} with $C=C_1C_2$.
\eof

Note that, in Lemma \ref{divbv} the estimate \eqref{bvzd} depends on $\a$. 
This dependency on $\a$ is not necessary in the one dimensinal case (see \eqref{bvz} in Remark 
\ref{divbv1d}) and for particular meshes in the two and three dimensional 
cases. Recall also that, 
except if $d=1$, the space $BV(\R^d)$ is not included in 
$L^\infty(\R^d)$. In particular, it is then quite easy to 
prove that, contrary to the 1D case given in Remark \ref{divbv1d}, it is 
not possible, for $d=2$ or $3$, to replace, in \eqref{bvzd}, $u_\cv$ by the mean value of $u$ over an 
arbitrary ball (for instance) included in $\cv$. 

\smallskip

Let us  now give the ``strong time $BV$ estimate".

\begin{lemma}\label{estlbvt}
Under Assumption \ref{estdonnees}, let $\mesh$ be an admissible mesh in the 
sense of Definition \ref{meshypmd} and $k>0$. Let $g \in C(\R^2,\R)$ satisfy 
Assumption \ref{estcondflu}. Assume that
$u_0 \in BV(\R^\dim)$\index{BV!Initial condition in --} and that $\vbf$ does not depend on $t$.

Let $\{ u^n_\cv, n \in \N, \cv \in \mesh\}$ be the solution of
\eqref{estschemai}, \eqref{estschemaz} such that $\ucvn \in
[U_m,U_M]$ for all $\cv \in \mesh$ and all $\nnn$ (existence and uniqueness of such a solution is
given by Proposition \refp{esteu}).

Then, there exists $C_b$, only depending on  $\vbf$, $g$, $u_0$ and $\a$ such that

\be
\sum_ {\cv \in \mesh}\frac {\mcv} k \vert \ucvnp- \ucvn \vert \le C_b, \,
\forall n \in \N.
\label{estbvt}
\ee
\end{lemma}

\proof  of lemma \ref{estlbvt}

\smallskip

Since $\vbf$ does not depend on $t$, one denotes $\vcvcvv=\vcvcvvn$, for 
all $\cv\in \mesh$ and all $\cvv \in \Nv$.

For $n \in \N$, let

\[
A_n=\sum_ {\cv \in \mesh} \mcv \frac {\vert \ucvnp -
\ucvn \vert} k
\]

and

\[
B_n=\sum_ {\cv \in \mesh} \vert \sumv 
[\vcvcvv \ g(\ucvn,\ucvvn) - \vcvvcv \ g(\ucvvn,\ucvn)] \vert.
\]

Since $u_0 \in BV(\R^\dim)$ and $\div \vbf=0$, there exists $C_b>0$, 
only depending on $\vbf$, $g$, $u_0$ and $\a$, such that
$B_0 \le C_b$. Indeed, 

\[
B_0 \le \summesh \sumv V(g_1+g_2)\medgecvcvv \vert u_\cv^0-u_\cvv^0 \vert.
\]

Thanks to lemma \ref{divbv}, $B_0 \le C_b$ with $C_b= 2V(g_1+g_2)C 
(1/\a^4)\snorm{u_0}{BV(\R^\dim)}$, where $C$ only depends on the space 
dimension ($\dim=1$, $2$ or $3$).

\smallskip

From \eqref{estschemai}, one deduces that $B_{n+1} \le A_n$, for all $n \in
\N$. In order to prove Lemma \ref{estlbvt}, there  only remains to
prove that $A_n \le B_n$ for all $n \in
\N$ (and to conclude by induction).

\smallskip

Let $n \in \N$, in order to prove that $A_n \le B_n$, recall that the implicit scheme \eqref{estschemai} reads

\be
\mcv \frac {\ucvnp - \ucvn} k + \sumv 
\Bigl(\vcvcvv \ g(\ucvnp,\ucvvnp) - \vcvvcv \ g(\ucvvnp,\ucvnp)\Bigr) = 0.
\label{estschemir}
\ee

From \eqref{estschemir}, one deduces, for all $\cv\in \mesh$,

\[
\bal
\mcv \dsp{\frac {\ucvnp - \ucvn} k } 
+\sumv \vcvcvv \ (g(\ucvnp,\ucvvnp)-g(\ucvn,\ucvvnp))
\\
+\sumv \vcvcvv \ \Bigl(g(\ucvn,\ucvvnp)-g(\ucvn,\ucvvn)\Bigr)
-\sumv \vcvvcv \ \Bigl(g(\ucvvnp,\ucvnp)-g(\ucvvn,\ucvnp)\Bigr)
\\
-\sumv \vcvvcv \ \Bigl( g(\ucvvn,\ucvnp)-g(\ucvvn,\ucvn)\Bigr)
\\
= -\sumv \vcvcvv   \, g(\ucvn,\ucvvn) + \sumv \vcvvcv \ g(\ucvvn,\ucvn).
\ean
\]

Using the monotonicity properties of $g$, one obtains for all $\cv\in
\mesh$,

\be
\bal
\mcv \dsp{\frac {\vert \ucvnp - \ucvn \vert} k} 
+\sumv \vcvcvv \ \vert g(\ucvnp,\ucvvnp)-g(\ucvn,\ucvvnp) \vert
\\
+\sumv \vcvvcv \ \vert g(\ucvvn,\ucvnp)-g(\ucvvn,\ucvn) \vert
\\
\le \vert -\sumv \vcvcvv g(\ucvn,\ucvvn) + 
\sumv \vcvvcv \ g(\ucvvn,\ucvn) \vert
\\
+\sumv \vcvcvv \ \vert g(\ucvn,\ucvvnp)-g(\ucvn,\ucvvn)\vert
+\sumv \vcvvcv \ \vert g(\ucvvnp,\ucvnp)-g(\ucvvn,\ucvnp)\vert.
\ean
\label{estschemirb}
\ee

In order to deal with convergent series, let us proceed as in the proof 
of proposition 
\ref{esteu}. For $0 < \gamma < 1$, let $\vphig: \R^\dim \mapsto \R_+^\star$ be 
defined by $\vphig(x)=\exp(-\gamma \vert x
\vert)$. 

For $\cv\in \mesh$, let $\vphi_{\gamma,\cv}$ be the
mean value of $\vphig$ on $\cv$. As in Proposition \ref{esteu}, since
$\vphig$ is integrable over
$\R^\dim$, $\sum_{\cv \in 
\mesh}\vphi_{\gamma,\cv} < \infty$. Therefore, multiplying \eqref{estschemirb} by 
$\vphi_{\gamma,\cv}$ (for a fixed $\gamma$) and summing over $\cv \in 
\mesh$ yields six convergent series which can be reordered to give

\[
\bal
\summesh \mcv \dsp{\frac {\vert \ucvnp - \ucvn \vert} k} \vphi_{\gamma,\cv}
\\
\le \summesh \vert -\sumv \vcvcvv g(\ucvn,\ucvvn) + 
\sumv \vcvvcv \ g(\ucvvn,\ucvn) \vert \vphi_{\gamma,\cv}
\\
+ \summesh \sumv \vcvcvv \ \vert g(\ucvnp,\ucvvnp)-g(\ucvn,\ucvvnp) \vert 
\vert \vphi_{\gamma,\cv} - \vphi_{\gamma,\cvv}\vert
\\
+\summesh \sumv \vcvvcv \ \vert g(\ucvvn,\ucvnp)-g(\ucvvn,\ucvn) \vert
\vert \vphi_{\gamma,\cv} - \vphi_{\gamma,\cvv}\vert.
\ean
\]

For $\cv \in \mesh$, let $\centercv \in \cvbar$ be such that 
$\vphi_{\gamma,\cv}=\vphig(\centercv)$. Let $\cv \in \mesh$ and 
$\cvv \in \Nv$. Then there exists $s \in (0,1)$ such that  
$\vphi_{\gamma,\cvv} - \vphi_{\gamma,\cv}=\grad \vphig(\centercv 
+ s (\centercvv-\centercv)) \cdot (\centercvv-\centercv)$. Using $\vert 
\grad \vphig(x) \vert = \gamma \exp(-\gamma \vert x \vert)$, this yields 
$\vert \vphi_{\gamma,\cvv} - \vphi_{\gamma,\cv} \vert \le 2 h \gamma 
\exp(2h \gamma) \vphi_{\gamma,\cv} \le 2 h  \gamma
\exp(2h) \vphi_{\gamma,\cv}$.

Then, using the assumptions \ref{estdonnees} and \ref{estcondflu}, 
there exists some $a$ only depending on $k$, $V$, $h$, $\a$, $g_1$ and 
$g_2$ such that

\[
\bal
\summesh \mcv \dsp{\frac {\vert \ucvnp - \ucvn \vert} k} 
\vphi_{\gamma,\cv}(1-\gamma a)
\\
\le \summesh \vert -\sumv \vcvcvv g(\ucvn,\ucvvn) + 
\sumv \vcvvcv \ g(\ucvvn,\ucvn) \vert \vphi_{\gamma,\cv} \le B_n.
\ean
\]

Passing to the limit in the latter inequality as $\gamma \tends 0$ 
yields $A_n \le B_n$. This
completes the proof of Lemma \ref{estlbvt}. \eof

\index{Scheme!Implicit -- for hyperbolic equations|)}
\section{Entropy inequalities for the approximate solution}\label{entineq}

\hpage{Hyperbolic multid}

In this section, an entropy estimate on the approximate solution is
proved (Theorem \ref{estentest}), which will be used in the proofs of 
convergence and error estimate of the numerical scheme. In order to obtain this entropy 
estimate, some discrete entropy inequalities satisfied by the 
approximate solution are first derived.

\index{BV!Strong time -- estimate|)}


\subsection{Discrete entropy inequalities}

\index{Discrete!-- entropy inequalities|(}\index{Entropy! 
Discrete -- inequalities|(}
In the case of the explicit scheme, the following lemma
asserts that the scheme \eqref{estschema} satisfies a discrete entropy condition (this is 
classical in the study of 1D schemes, see e.g. \citet{godlewski-ellipses}, 
\citet{godlewski-springer}).

\begin{lemma}\label{estldei} 
Under assumption \refp {estdonnees}, let $\mesh$ be an admissible mesh in the 
sense of Definition \ref{meshypmd} and $k>0$. Let $g \in C(\R^2,\R)$ 
satisfying 
assumption \ref{estcondflu} and assume that \eqref{estcfl} holds.

Let $\uapp$ be given by \eqref{estuapp}, \eqref{estschema}, 
\eqref{estschemaz}; then, for all $\sk \in \R$, $\cv \in \mesh$ and $n\in\N$,
the following inequality holds:

\be\ban
\mcv \dsp{\frac {\vert\ucvnp -  \sk\vert - \vert\ucvn -  \sk\vert} k} + 
\sumv \Bigl[ & \vcvcvvn \Bigl( 
g(\ucvn\top\sk,\ucvvn\top\sk) - g(\ucvn\bot\sk,\ucvvn\bot\sk) \Bigl) - \\
& \vcvvcvn \Bigl( 
g(\ucvvn\top\sk,\ucvn\top\sk) - g(\ucvvn\bot\sk,\ucvn\bot\sk) \Bigl) 
\Bigr] &\le 0. 
\ean\label{estschemak}\ee

\end{lemma}

\proof of lemma \ref{estldei}

\smallskip

From relation \eqref{estschema}, we express $\ucvnp$ as a
function of $\ucvn$ and $\ucvvn$, $\cvv \in \Nv$,

\[
\ucvnp = \ucvn + \frac k {\mcv} \sumv 
(\vcvvcvn \ g(\ucvvn,\ucvn) -   \vcvcvvn \ g(\ucvn,\ucvvn)).
\]

The right hand side is nondecreasing with respect to $\ucvvn$, $\cvv \in 
\Nv$. It is also nondecreasing with respect to $\ucvn$, thanks to the Courant-Friedrichs-Levy 
condition \eqref{estcfl}, and the Lipschitz continuity of $g$.

Therefore, for all $\sk\in\R$, using $\div \vbf = 0$, we have:

\be
\ucvnp\top\sk \le \ucvn\top\sk + \frac k {\mcv} \sumv 
\Bigl[\vcvvcvn \ g(\ucvvn\top\sk,\ucvn\top\sk)  - \vcvcvvn \
g(\ucvn\top\sk,\ucvvn\top\sk) \Bigr] 
\label{est2ent}\ee

and

\be
\ucvnp\bot\sk \ge \ucvn\bot\sk + \frac k {\mcv} \sumv 
(\vcvvcvn \ g(\ucvvn\bot\sk,\ucvn\bot\sk) - \vcvcvvn \ 
g(\ucvn\bot\sk,\ucvvn\bot\sk)).
\label{est3ent}\ee

The difference between \eqref{est2ent} and \eqref{est3ent} leads directly to
\eqref{estschemak}. Note that using $\div\vbf = 0$ leads to

\be\ban
&\mcv \dsp{\frac {\vert\ucvnp -  \sk\vert - \vert\ucvn -  \sk\vert} k} +\\ 
\sumv \Bigl[ &\vcvcvvn \Bigl( 
g(\ucvn\top\sk,\ucvvn\top\sk) - f(\ucvn\top\sk)
- g(\ucvn\bot\sk,\ucvvn\bot\sk) + f(\ucvn\bot\sk)\Bigl) -\\
&\vcvvcvn \Bigl( 
g(\ucvvn\top\sk,\ucvn\top\sk) - f(\ucvn\top\sk)
- g(\ucvvn\bot\sk,\ucvn\bot\sk) + f(\ucvn\bot\sk)\Bigl) \Bigr]
\le 0. 
\ean\label{estschemakb}\ee \eof

For the implicit scheme, one obtains the same kind of discrete entropy
inequalities.

\begin{lemma}\label{estldeii} 
Under assumption \refp{estdonnees}, let $\mesh$ be an admissible mesh in the 
sense of Definition \refp{meshypmd} and $k>0$. Let $g \in C(\R^2,\R)$ 
satisfying 
assumption \ref{estcondflu}.

Let $\{u^n_\cv,$
$n \in \N,$ $\cv\in \mesh\} \subset [U_m,U_M]$ be the solution of
\eqref{estschemai},\eqref{estschemaz} (the 
existence and uniqueness of such a solution is
given by Proposition
\ref{esteu}). Then, for all $\sk\in\R$, $\cv \in \mesh$ and $n\in\N$,
the following inequality holds:

\be\bac
\mcv \dsp{\frac {\vert\ucvnp -  \sk\vert - \vert\ucvn -  \sk\vert} k} + 
\sumv  \Bigl[ \vcvcvvn \Bigl( 
g(\ucvnp\top\sk,\ucvvnp\top\sk) - g(\ucvnp\bot\sk,\ucvvnp\bot\sk) \Bigl) \\
- \vcvvcvn \Bigl( 
g(\ucvvnp\top\sk,\ucvnp\top\sk) - g(\ucvvnp\bot\sk,\ucvnp\bot\sk) \Bigl) \Bigr]
\le 0. 
\ean
\label{estschemaki}\ee

\end{lemma}

\proof of lemma \ref{estldeii}

\smallskip

Let $\sk\in\R$, $\cv \in \mesh$ and $n\in\N$.
Equation \eqref{estschemai} may be written as

\[
\ucvnp=\ucvn- \frac k {\mcv} \sumv 
(\vcvcvvn \ g(\ucvnp,\ucvvnp) - \vcvvcvn \ g(\ucvvnp,\ucvnp)).
\]

The right hand side of this last equation is nondecreasing with 
respect to $\ucvn$ and with respect to $\ucvvnp$ for all $\cvv \in \Nv$. 
Thus,

\[
\ucvnp \le \ucvn \top \sk - \frac k {\mcv} \sumv 
(\vcvcvvn \ g(\ucvnp,\ucvvnp\top \sk) - \vcvvcvn \ g(\ucvvnp\top \sk,\ucvnp)).
\]

Writing $\sk=\sk- \dsp{\frac k {\mcv}} \sumv 
(\vcvcvvn \ g(\sk,\sk) - \vcvvcvn \ g(\sk,\sk))$, one may remark that

\[
\sk \le \ucvn \top \sk - \frac k {\mcv} \sumv 
(\vcvcvvn \ g(\sk,\ucvvnp\top \sk) - \vcvvcvn \ g(\ucvvnp\top \sk,\sk)).
\]

Therefore, since $\ucvnp \top \sk=\ucvnp$ or $\sk$,

\be
\ucvnp \top \sk \le \ucvn \top \sk - \frac k {\mcv} \sumv 
(\vcvcvvn \ g(\ucvnp \top \sk,\ucvvnp\top \sk) - \vcvvcvn \ g(\ucvvnp\top \sk,\ucvnp \top \sk)).
\label{impt}
\ee

A similar argument yields

\be
\ucvnp \bot \sk \ge \ucvn \bot \sk - \frac k {\mcv} \sumv 
(\vcvcvvn \ g(\ucvnp \bot \sk,\ucvvnp \bot \sk) - \vcvvcvn \ g(\ucvvnp \bot \sk,\ucvnp \bot \sk)).
\label{impb}
\ee

Hence, substracting \eqref{impb} to \eqref{impt} gives \eqref{estschemaki}.
\eof
\index{Discrete!-- entropy inequalities|)}\index{Entropy!
Discrete -- inequalities|)}



\subsection{Continuous entropy estimates for the approximate solution}

For $\Omega = \R^\dim$ or $\RNRp$, we denote by $\mathcal{M}(\Omega)$ the set
of positive measures on $\Omega$, that is of $\s$-additive applications
from the
 Borel $\s$-algebra of $\Omega$ in $\overline \R_+$. If $\mu \in
\mathcal{M}(\Omega)$ and $\psi \in C_c(\O)$, one sets $\act{\mu}{\psi}=\int
\psi d\mu$.

The following
theorems investigate the  entropy inequalities which are satisfied by
the approximate solutions  $\uapp$  in the case of the time
explicit scheme
(Theorem \ref{estentest}) and in the case of the time implicit scheme (Theorem
\ref{estentesti}). 

\begin{theo}\label{estentest}
Under assumption \ref{estdonnees}, let $\mesh$ be an admissible mesh in the 
sense of Definition \refp{meshypmd} and $k>0$. Let $g \in C(\R^2,\R)$ satisfy 
assumption \ref{estcondflu} and assume that \eqref{estcfl} holds.

Let $\uapp$ be given by \eqref{estuapp}, \eqref{estschema},  \eqref{estschemaz}; then
there exist   $\mapp \in \mathcal{M}(\RNRp)$ and
$\mapz \in \mathcal{M}(\R^\dim)$ such that

\be\ba
\dsp{\int_{\Rp}\int_{\R^\dim}} \Bigl( &\vert \uapp(x,t) - \sk \vert \varphi_t(x,t) +\\
&(f(\uapp(x,t)\top\sk) - f(\uapp(x,t)\bot\sk)) 
\vbf(x,t)\cdot\grad\varphi(x,t) \Bigr) dx dt & + \\
\dsp{\int_{\R^\dim}} & \vert u_0(x) - \sk \vert \varphi(x,0) dx & \ge \\ 
- \dsp{\int_\RNRp} & \Bigl( \vert\varphi_t(x,t)\vert + 
\vert\grad\varphi(x,t)\vert\Bigr) d\mapp(x,t) 
- \dsp{\int_{\R^\dim}}   \varphi(x,0)  d\mapz(x), \\ \\ &
 \ \ \forall \sk\in\R,
\ \ \forall\varphi\in C^\infty_c(\R^\dim\times\R_+,\R_+). 
\ea\label{kruzkovap}\ee

The measures $\mapp$ and
$\mapz$ verify the following properties:

\begin{enumerate}

\item For all $R>0$ and $T>0$, there exists $C$ depending
only on
$\vbf$, $g$, $u_0$, $\a$, $\xi$, $R$ and $T$ such that, for $h < R$ and $k 
< T$,
 
\be
\mapp(B(0,R)\times [0,T]) \le C \sqrt{h}.
\label{estmapp}
\ee

\item The measure $\mapz$ is the measure of density $\vert
u_{0}(\cdot)-u_{\mesh,0}(\cdot) \vert$ with respect to the Lebesgue 
measure, where $u_{\mesh,0}$ is defined by $u_{\mesh,0}(x)=u_\cv^0$ for 
a.e. $x \in \cv$, for all $\cv \in \mesh$.

If $u_0 \in BV(\R^\dim)$\index{BV!Initial condition in --}, then there exists
$D$, only depending on $u_0$ and $\a$, such that

\be
\mapz(\R^\dim) \le D h.
\label{estmapz}\ee

\end{enumerate}

\end{theo}

\begin{remark} \normalfont \label{estentropiex}~
\benum
\item Let $u$ be the weak entropy solution to
\eqref{estpbcont}-\eqref{estcondini}. Then \eqref{kruzkovap} is satisfied with $u$
instead of $\uapp$ and $\mapp = 0$ and $\mapz = 0$. 
\item Let $BV_{loc}(\R^\dim)$\index{BV!-- Space} be the set of $v \in L^1_{loc}(\R^\dim)$ 
such that the restriction of $v$ to $\O$ belongs to $BV(\O)$ for all 
open bounded subset $\O$ of $\R^\dim$.

An easy adaptation of the following proof gives that if $u_0 \in 
BV_{loc}(\R^\dim)$ instead of $BV(\R^\dim)$ (in the second item 
of Theorem \ref{estentest}) then, for all
$R>0$, there exists
$D$, only depending on $u_0$, $\a$ and $R$,  such that $\mapz(B(0,R)) \le 
D h$.
\eenum
\end{remark}

\proof of Theorem \ref{estentest}

\smallskip

Let $\varphi\in C^\infty_c(\R^\dim\times\R_+,\R_+)$ and $\sk \in \R$.

Multiplying \eqref{estschemakb} by 
$k\vphi_\cv^n=(1/{\mcv}) \int_{n k}^{(n+1)k} \int_\cv \varphi(x,t) dx dt$ and
summing the result for all $\cv \in \mesh$ and $n\in\N$ yields

\[
T_1 + T_2 \le 0,
\]

with 

\be
T_1 = \sumIN \summesh \frac {\vert\ucvnp -  \sk\vert - \vert\ucvn -  \sk\vert} k
\int_{n k}^{(n+1)k} \int_\cv \varphi(x,t) dx dt, 
\label{est1kap}\ee

and

\be\bal
T_2 = k \sumIN \dsp\sum_{(\cv,\cvv) \in \edges_n} \Bigl[
\\  
\vcvcvvn \vphi_\cv^n \Bigl(g(\ucvn\top\sk,\ucvvn\top\sk) -  f(\ucvn\top\sk) - 
g(\ucvn\bot\sk,\ucvvn\bot\sk) +  f(\ucvn\bot\sk)\Bigr)
\\
- \vcvcvvn \vphi_\cvv^n \Bigl(g(\ucvn\top\sk,\ucvvn\top\sk) -  f(\ucvvn\top\sk) - 
g(\ucvn\bot\sk,\ucvvn\bot\sk) +  f(\ucvvn\bot\sk)\Bigr)
\\
-\vcvvcvn \vphi_\cv^n \Bigl(g(\ucvvn\top\sk,\ucvn\top\sk) -  f(\ucvn\top\sk) - 
g(\ucvvn\bot\sk,\ucvn\bot\sk) +  f(\ucvn\bot\sk)\Bigr)
\\
+ \vcvvcvn \vphi_\cvv^n \Bigl(g(\ucvvn\top\sk,\ucvn\top\sk) -  f(\ucvvn\top\sk) - 
g(\ucvvn\bot\sk,\ucvn\bot\sk) +  f(\ucvvn\bot\sk)\Bigr) \Bigr],
\ean\label{est2kap}\ee

where $\edges_n=\{(\cv,\cvv) \in \mesh^2$, $u^n_\cv > u^n_\cvv\}$.

\smallskip

One has to prove

\be
T_{10} + T_{20} \le \dsp{\int_\RNRp} \Bigl(
\vert\varphi_t(x,t)\vert + 
\vert\grad\varphi(x,t)\vert\Bigr) d\mapp(x,t) 
+ \dsp{\int_{\R^\dim}}   \varphi(x,0)  d\mapz(x),\label{estentcont}
\ee

for some convenient measures $\mapp$ and $\mapz$, and $T_{10}$,
$T_{20}$ defined as follows

\[
T_{10}=
- \int_{\Rp}\int_{\R^\dim} \vert \uapp(x,t) - \sk\vert\varphi_t(x,t) dx dt
- \int_{\R^\dim} \vert u_{0}(x) - \sk \vert \varphi(x,0) dx,
\]

\be
T_{20} = - \int_{\Rp}\int_{\R^\dim} \Bigl( 
(f(\uapp(x,t)\top\sk) - f(\uapp(x,t)\bot\sk)) 
\vbf(x,t)\cdot\grad\varphi(x,t) \Bigr) dx dt. 
\label{est7kap}\ee

In order to prove \eqref{estentcont}, one compares $T_1$ and $T_{10}$
(this will give $\mapz$, and a part of $\mapp$) and one compares
$T_2$ and
$T_{20}$ (this will give
another part of $\mapp$).

Inequality \eqref{estbvft} (in the comparison of $T_1$ and $T_{10}$) and
Inequality \eqref{estbvfx} (in the comparison of $T_2$ and
$T_{20}$) will be used in order to obtain \eqref{estmapp}.

\textbf{Comparison of $T_1$ and $T_{10}$}

Using the definition of $\uapp$ and introducing the function 
$u_{\mesh,0}$ (defined by $u_{\mesh,0}(x)= u_\cv^0$, for a.e. $x \in 
\cv$, for all $\cv \in \mesh$) yields

\[
\ban
T_{10}= \sumIN \summesh 
\dsp{\frac {\vert\ucvnp -  \sk\vert - \vert\ucvn -  \sk\vert} k
\int_{n k}^{(n+1)k} \int_\cv \varphi(x,(n+1)k) dx dt} \ \ + \\
\dsp{\int_{\R^\dim}} (\vert u_{\mesh,0}(x) - \sk \vert -
\vert u_{0}(x)- \sk \vert) \varphi(x,0) dx.
\ean
\]

The function $\vert \cdot - \sk\vert$ is Lipschitz
continuous with a Lipschitz constant equal to $1$, we then obtain

\[
\ban 
\vert T_1 - T_{10} \vert \le 
& \sumIN \summesh \frac {\vert\ucvnp - \ucvn\vert} k \int_{n k}^{(n+1)k} \int_\cv 
\vert \varphi(x,(n+1)k) - \varphi(x,t) \vert dx dt \ \ + \\
& \dsp{\int_{\R^\dim}} \vert u_{0}(x)-u_{\mesh,0}(x)\vert \varphi(x,0) dx,
\ean
\]

which leads to

\be\ban 
\vert T_1 - T_{10} \vert \le 
& \sumIN \summesh \vert\ucvnp - \ucvn\vert \int_{n k}^{(n+1)k} \int_\cv 
\vert \varphi_t(x,t) \vert dx dt \ \ + \\
& \dsp{\int_{\R^\dim}} \vert u_{0}(x)-u_{\mesh,0}(x)\vert \varphi(x,0)
dx.
\ean\label{est5kapb}\ee

Inequality \eqref{est5kapb} gives

\be 
\vert T_1 - T_{10} \vert \le 
\dsp{\int_\RNRp} \vert\varphi_t(x,t)\vert d\nu_{\mesh,k}(x,t) 
+ \dsp{\int_{\R^\dim}}   \varphi(x,0)  d\mapz(x),\label{estmajtime}\ee

where the measures $\mapz \in \mathcal{M}(\R^\dim)$ and $\nu_{\mesh,k} \in
\mathcal{M}(\RNRp)$ are defined, by their action on $C_c(\R^\dim)$ and
$C_c(\RNRp)$, as follows

\[
\act{\mapz}{\psi} =\dsp{\int_{\R^\dim}} \vert
u_{0}(x)-u_{\mesh,0}(x)\vert \psi(x) dx,
\ \forall \psi\in C_c(\R^\dim),
\]

\[
\bal
\act{\nu_{\mesh,k}}{\psi}= 
\dsp{\sum_{n\in\N} \sum_{\cv\in\mesh} 
\vert\ucvnp - \ucvn\vert \int_{n k}^{(n+1)k}} \int_\cv \psi(x,t) dx dt, \\
\ \forall \psi\in C_c(\RNRp).
\ean
\]

The measures $\mapz$ and $\nu_{\mesh,k}$ are absolutely
continuous with respect to the Lebesgue measure. Indeed, one has $d\mapz(x)=
\vert u_{0}(x)-u_{\mesh,0}(x) \vert dx$ and
$d\nu_{\mesh,k}(x,t)=(\sum_{n\in\N} \sum_{\cv\in\mesh} 
\vert\ucvnp - \ucvn\vert 1_{\cv \times [nk,(n+1)k)})dxdt$ (where $1_\O$
denotes the characteristic function of $\O$ for any Borel subset $\O$ of
$\R^{d+1}$).

\smallskip

If $u_0 \in  BV(\R^\dim)$, the measure $\mapz$ verifies \eqref{estmapz} with 
some $D$ only depending on $\vert u_0 \vert_{BV(\R^\dim)}$ and $\a$ (this is 
classical  result which is 
given in Lemma \ref{divbvd} below  for the sake of completeness).

The measure $\nu_{\mesh,k}$ satisfies \eqref{estmapp}, with
$\nu_{\mesh,k}$ instead of $\mapp$, thanks to \eqref{estbvft} and
condition \eqref{estcfl}. Indeed, for $R>0$ and $T>0$,

\[  
\nu_{\mesh,k}(B(0,R) \times [0,T]) = \int_0^T \int_{B(0,R)} \sum_{n\in\N} \sum_{\cv\in\mesh} 
\vert\ucvnp - \ucvn\vert 1_{\cv \times [nk,(n+1)k)} dxdt,
\]

which yields, with $\mesh_{2R}=\{\cv \in \mesh$, $\cv \subset B(0,2R)\}$ 
and $\NTk k < T \le (\NTk+1)k$, $h < R$ and $k<T$,

\[  
\nu_{\mesh,k}(B(0,R) \times [0,T]) \le k \sum_{n=0}^{\NTk} \sum_{\cv \in 
\mesh_{2R}} \mcv \vert \ucvnp - \ucvn \vert \le \frac {kC_1} {\sqrt h},
\]

where $C_1$ is given by lemma \ref{estwbvest} and only depends on $\vbf$, $g$, $u_0$, $\a$, $\xi$, $R$,
$T$. Finally, since the condition \eqref{estcfl} gives $k \le C_2 h$, where $C_2$  only  depends 
on $\vbf$, $g$, $u_0$, $\a$, $\xi$, the last inequality 
yields, for $h < R$ and $k<T$,

\be
\nu_{\mesh,k}(B(0,R) \times [0,T]) \le  C_3 \sqrt h,
\label{boundnu}
\ee

with $C_3=C_1C_2$.

\textbf{Comparison of $T_2$ and $T_{20}$}

Using $\div\vbf=0$, and gathering \eqref{est7kap} by interfaces, we get

\be\ban
T_{20} = - \sumIN \dsp\sum_{(\cv,\cvv) \in \edges_n} \Bigl[ & \Bigl( 
(f(\ucvn\top\sk) - f(\ucvn\bot\sk)) - (f(\ucvvn\top\sk) - f(\ucvvn\bot\sk))
\Bigr) \\ & \dsp{\int_{\edgecvcvv} \int_{n k}^{(n+1)k}} 
\Bigl( \vbf(x,t) \cdot \bfn_{\cv,\cvv} \,\varphi(x,t) \Bigr) 
\dfrontiere(x) dt\Bigr]. 
\ean\label{est8kap}\ee

Define, for all $\cv \in \mesh$, all $\cvv \in \Nv$ and all $\nnn$,

\[
(v\vphi)_{\cv,\cvv}^{n,+}=\frac 1 k \int_{n k}^{(n+1)k} \int_{\edgecvcvv}  
(\vbf(x,t)\cdot\bfn_{\cv,\cvv})^+ \varphi(x,t)  
\dfrontiere(x) dt
\]

and

\[
(v\vphi)_{\cv,\cvv}^{n,-}= \frac 1 k \int_{n k}^{(n+1)k} \int_{\edgecvcvv}  
(\vbf(x,t)\cdot\bfn_{\cv,\cvv})^- \varphi(x,t)  
\dfrontiere(x) dt.
\]

Note that $(v\vphi)_{\cv,\cvv}^{n,+}=(v\vphi)_{\cvv,\cv}^{n,-}$. Then, \eqref{est8kap} gives

\be\bal
T_{20} = k \sumIN \dsp\sum_{(\cv,\cvv) \in \edges_n} \Bigl[
\\  
(v\vphi)_{\cv,\cvv}^{n,+} \Bigl(g(\ucvn\top\sk,\ucvvn\top\sk) -  f(\ucvn\top\sk) - 
g(\ucvn\bot\sk,\ucvvn\bot\sk) +  f(\ucvn\bot\sk)\Bigr)
\\
- (v\vphi)_{\cvv,\cv}^{n,-} \Bigl(g(\ucvn\top\sk,\ucvvn\top\sk) -  f(\ucvvn\top\sk) - 
g(\ucvn\bot\sk,\ucvvn\bot\sk) +  f(\ucvvn\bot\sk)\Bigr)
\\
-(v\vphi)_{\cv,\cvv}^{n,-} \Bigl(g(\ucvvn\top\sk,\ucvn\top\sk) -  f(\ucvn\top\sk) - 
g(\ucvvn\bot\sk,\ucvn\bot\sk) +  f(\ucvn\bot\sk)\Bigr)
\\
+ (v\vphi)_{\cvv,\cv}^{n,+} \Bigl(g(\ucvvn\top\sk,\ucvn\top\sk) -  f(\ucvvn\top\sk) - 
g(\ucvvn\bot\sk,\ucvn\bot\sk) +  f(\ucvvn\bot\sk)\Bigr) \Bigr].
\ean\label{est8bkap}\ee

Let us introduce some terms related to the difference between $\vphi$ on 
$\cv \in \mesh$ and $\edgecvcvv \in \edges$,

\[
r_{\cv,\cvv}^{n,+} =\vert \vcvcvvn \vphi_\cv^n - 
(v\vphi)_{\cv,\cvv}^{n,+}\vert
\]

and

\[
r_{\cv,\cvv}^{n,-} = \vert \vcvvcvn \vphi_\cv^n - (v\vphi)_{\cv,\cvv}^{n,-} \vert.
\]

Then, from \eqref{est2kap} and \eqref{est8bkap}, 
 
\be\ban 
\vert T_2 - T_{20} \vert \le 
\sumIN k \dsp\sum_{(\cv,\cvv) \in \edges_n} \Bigl[ \\  
r_{\cv,\cvv}^{n,+} 
\Bigl(  
g(\ucvn\top\sk,\ucvvn\top\sk) -  f(\ucvn\top\sk) + 
g(\ucvn\bot\sk,\ucvvn\bot\sk) -  f(\ucvn\bot\sk) \Bigl) &+ \\
r_{\cvv,\cv}^{n,-} 
\Bigl( 
g(\ucvn\top\sk,\ucvvn\top\sk) -  f(\ucvvn\top\sk) + 
g(\ucvn\bot\sk,\ucvvn\bot\sk) -  f(\ucvvn\bot\sk) \Bigl) &+ \\
r_{\cv,\cvv}^{n,-} 
\Bigl( 
f(\ucvn\top\sk) - g(\ucvvn\top\sk,\ucvn\top\sk)  + 
f(\ucvn\bot\sk) - g(\ucvvn\bot\sk,\ucvn\bot\sk)  \Bigl) &+ \\
r_{\cvv,\cv}^{n,+} 
\Bigl( 
f(\ucvvn\top\sk) - g(\ucvvn\top\sk,\ucvn\top\sk)  + 
f(\ucvvn\bot\sk) - g(\ucvvn\bot\sk,\ucvn\bot\sk)  \Bigl) \Bigr].
\ean\label{est9kap}\ee

For all $(\cv,\cvv) \in \edges_n$, the following inequality holds:

\[
0 \le g(\ucvn\top\sk,\ucvvn\top\sk) -  f(\ucvn\top\sk) \le
\dsp{\max_{\ucvvn\le p \le q  \le\ucvn}} ( g(q,p) -  f(q)),
\]

more precisely, one has $g(\ucvn\top\sk,\ucvvn\top\sk) - 
f(\ucvn\top\sk)=0$, if $\sk\ge\ucvn$, and one has $g(\ucvn\top\sk,\ucvvn\top\sk) - 
f(\ucvn\top\sk)= g(q,p) -  f(q)$ with $p=\sk$ and
$q=\ucvn$ if $\sk\in [\ucvvn,\ucvn]$, and with $p=\ucvvn$
and $q=\ucvn$ if
$\sk\le\ucvvn$. 

In the same way, we can assert that

\[
0 \le g(\ucvn\bot\sk,\ucvvn\bot\sk) -  f(\ucvn\bot\sk) \le
\dsp{\max_{\ucvvn\le p \le q  \le\ucvn}} ( g(q,p) -  f(q)).
\]

The same analysis can be applied to the six other terms of \eqref{est9kap}.

To conclude the estimate on $\vert T_2-T_{20} \vert$, there remains to
estimate the two  quantities $r_{\cv,\cvv}^{n,\pm}$. This will be done
with convenient measures applied to $\vert\grad\varphi\vert$ and
$\vert\varphi_t\vert$. To estimate $r_{\cv,\cvv}^{n,+}$, for instance, one
remarks that

\[
r_{\cv,\cvv}^{n,+} \le  
\dsp{ \frac 1 {k^2 \mcv}
\int_{n k}^{(n+1)k} \int_{n k}^{(n+1)k} \int_{\cv} \int_{\edgecvcvv}} 
 \vert \varphi(x,t) - \varphi(y,s) \vert 
(\vbf(y,s)\cdot\bfn_{\cv,\cvv})^+ \dfrontiere(y) dx dt ds.
\]

Hence 
 
\[
\ban
r_{\cv,\cvv}^{n,+} \le 
\dsp{\frac 1 {k^2 \mcv}
\int_{n k}^{(n+1)k} \int_{n k}^{(n+1)k} \int_{\cv} \int_{\edgecvcvv} \int_0^1} 
\vert \grad\varphi(x+\theta(y-x),t+\theta(s-t))\cdot(y-x) + \\ 
 \varphi_t(x+\theta(y-x),t+\theta(s-t))(s-t)\vert  
(\vbf(y,s)\cdot\bfn_{\cv,\cvv})^+ d\theta 
\dfrontiere(y) dx dt ds 
\ean
\]
which yields

\[
\ban
r_{\cv,\cvv}^{n,+} \le
\dsp{\frac 1 {k^2 \mcv}
\int_{n k}^{(n+1)k} \int_{n k}^{(n+1)k} \int_{\cv} \int_{\edgecvcvv} \int_0^1} 
\Bigl( h \vert \grad\varphi(x+\theta(y-x),t+\theta(s-t)) \vert + \\ 
 k \vert \varphi_t(x+\theta(y-x),t+\theta(s-t))\vert \Bigr)    
(\vbf(y,s)\cdot\bfn_{\cv,\cvv})^+ d\theta 
\dfrontiere(y) dx dt ds.
\ean
\]

This leads to the definition of a measure $\mu_{\cv,\cvv}^{n,+}$, given by
its action on $C_c(\RNRp)$:

\[\bari
\act{\mu_{\cv,\cvv}^{n,+}}{\psi} =  
\dsp{ \frac 2 {k^2 \mcv} \int_{n k}^{(n+1)k} \int_{n k}^{(n+1)k} 
\int_{\cv} \int_{\edgecvcvv} \int_0^1}  
\Bigl( (h + k) \psi(x+\theta(y-x),t+\theta(s-t)) \Bigr)
\\
(\vbf(y,s)\cdot\bfn_{\cv,\cvv})^+ d\theta 
\dfrontiere(y) dx dt ds, \, \forall \psi\in C_c(\RNRp),
\ean\]

in order to have $2 r_{\cv,\cvv}^{n,+} \le \act{\mu_{\cv,\cvv}^{n,+}}{\vert \grad 
\vphi \vert + \vert \vphi_t \vert}$.

We define in the same way $\mu_{\cv,\cvv}^{n,-}$, changing 
$(\vbf(y,s)\cdot\bfn_{\cv,\cvv})^+$ in $(\vbf(y,s)\cdot\bfn_{\cv,\cvv})^-$. We finally define the
measure $\tilde {\nu}_{\mesh,k}$ by

\be\bari 
\act{\tilde {\nu}_{\mesh,k}}{\psi} = \sumIN k \dsp\sum_{(\cv,\cvv) \in \edges_n} \Bigl[   
 \Bigl(  
\dsp{\max_{\ucvvn\le p \le q  \le\ucvn}} ( g(q,p) -  f(q)) 
\Bigl) \act{\mu_{\cv,\cvv}^{n,+}}{\psi}
\\
+\Bigl(  
\dsp{\max_{\ucvvn\le p \le q  \le\ucvn}} ( g(q,p) -  f(p))
\Bigl) \act{\mu_{\cvv,\cv}^{n,-}}{\psi}
\\
+\Bigl(  
\dsp{\max_{\ucvvn\le p \le q  \le\ucvn}} ( f(q) - g(p,q) ) 
\Bigl) \act{\mu_{\cv,\cvv}^{n,-}}{\psi}
\\
+\Bigl(  
\dsp{\max_{\ucvvn\le p \le q  \le\ucvn}} ( f(p) - g(p,q) ) 
\Bigl) \act{\mu_{\cvv,\cv}^{n,+}}{\psi}\Bigl].
\ean\label{est10kap}\ee

Since $2 r_{\cv,\cvv}^{n,\pm} \le \act{\mu_{\cv,\cvv}^{n,\pm}}{\vert \grad 
\vphi \vert + \vert \vphi_t \vert}$, \eqref{est9kap} and \eqref{est10kap} 
leads to $\vert T_2-T_{20} \vert \le \act{\tilde {\nu}_{\mesh,k}}{\vert \grad \vphi \vert 
+ \vert \vphi_t \vert}$. Therefore, setting $\mapp 
={\nu}_{\mesh,k}+\tilde {\nu}_{\mesh,k}$, using \eqref{estmajtime} and $T_1+T_2 \le 0$,

\[
T_{10} + T_{20} \le \dsp{\int_\RNRp} \Bigl(
\vert\varphi_t(x,t)\vert + 
\vert\grad\varphi(x,t)\vert\Bigr) d\mapp(x,t) 
+ \dsp{\int_{\R^\dim}}   \varphi(x,0)  d\mapz(x),
\]

which is \eqref{estentcont} and yields \eqref{kruzkovap}.

\smallskip

There remains to prove \eqref{estmapp}.

For all $\cv \in \mesh$, let $\centercv$ be an arbitrary point of $\cv$. 
For all $\cv \in \mesh$, all $\cv \in \Nv$ and all $\nnn$, the supports of the 
measures $\mu_{\cv,\cvv}^{n,\pm}$ are included in the closed set 
$\bar{B}(\centercv,h) \cap [nk,(n+1)k]$. Furthermore, 

\[
\mu_{\cv,\cvv}^{n,+}(\R^\dim \times \Rp) \le 2 \vcvcvvn (h+k) \text{ and } 
\mu_{\cv,\cvv}^{n,-}(\R^\dim \times \Rp) \le 2 \vcvvcvn (h+k).
\]

Then, for 
all $R >0$ and $T>0$, the definition of $\mapp$ (i.e. $\mapp 
={\nu}_{\mesh,k}+\tilde {\nu}_{\mesh,k})$) leads to

\[
\bal
\dsp{\mapp(B(0,R) \times [0,T]) \le C_3 \sqrt h}
\\
\dsp{+2(h+k)\sum_{n=0}^{\NTk}k 
\sum_{(\cv,\cvv) \in \edges_{2R}^n} \Bigl[
\vcvcvvn \Bigl(\max_{\ucvvn\le p \le q  \le\ucvn} ( g(q,p) -  f(q)) + 
\max_{\ucvvn\le p \le q  \le\ucvn} ( g(q,p) -  f(p)) \Bigr)}
\\
\dsp{+ \vcvvcvn \Bigl(
\max_{\ucvvn\le p \le q  \le\ucvn} ( f(q) - g(p,q) ) + 
\max_{\ucvvn\le p \le q  \le\ucvn} ( f(p) - g(p,q) ) \Bigr) \Bigr],}
\ean
\]

for $h <R$ and $k <T$, where $C_3 \sqrt h$ is the bound of $\nu_{\mesh,k}(B(0,R) \times [0,T])$ 
given in \eqref{boundnu}.
Therefore, thanks to Lemma \ref{estwbvest},

\[
\mapp(B(0,R) \times [0,T]) \le C_3 \sqrt h + (1+C_2)h\frac {C_4}{\sqrt 
h}=C \sqrt h,
\]

where $C$ only depends on $\vbf$, $g$, $u_0$, $\a$, $\xi$, $R$ and $T$. 
The proof of Theorem \ref{estentest} is complete. \eof

The following theorem investigates the case of the implicit scheme.

\begin{theo}\label{estentesti}
Under Assumption \ref{estdonnees}, let $\mesh$ be an admissible mesh in the 
sense of Definition \ref{meshypmd} and $k>0$. Let $g \in C(\R^2,\R)$ satisfy 
Assumption \ref{estcondflu}.

Let $\{ \ucvn, n \in \N, \cv \in \mesh \}$, such that $\ucvn \in [U_m,U_M]$ 
for all $\cv \in \mesh$ and $\nnn$, be the solution of 
\eqref{estschemai},\eqref{estschemaz} (existence and uniqueness of such a solution are
given by Proposition \ref{esteu}). Let
$\uapp$ be given by \eqref{estuapp}. Assume that $\vbf$ does not depend on
$t$ and that $u_0 \in BV(\R^\dim)$.

Then, there exist $\mapp \in \mathcal{M}(\RNRp)$ and
$\mapz \in \mathcal{M}(\R^\dim)$ such that

\be\ba
\dsp{\int_{\Rp}\int_{\R^\dim}} \Bigl( &\vert \uapp(x,t) - \sk \vert \varphi_t(x,t) +\\
&(f(\uapp(x,t)\top\sk) - f(\uapp(x,t)\bot\sk)) 
\vbf(x,t)\cdot\grad\varphi(x,t) \Bigr) dx dt & + \\
\dsp{\int_{\R^\dim}} & \vert u_0(x) - \sk \vert \varphi(x,0) dx & \ge \\ 
- \dsp{\int_\RNRp} & \Bigl( \vert\varphi_t(x,t)\vert + 
\vert\grad\varphi(x,t)\vert\Bigr) d\mapp(x,t) 
- \dsp{\int_{\R^\dim}}   \varphi(x,0)  d\mapz(x), \\ \\ &
 \ \ \forall \sk\in\R,
\ \ \forall\varphi\in C^\infty_c(\R^\dim\times\R_+,\R_+). 
\ea\label{kruzkovapi}\ee

The measures $\mapp$ and
$\mapz$ verify the following properties:

\begin{enumerate}

\item For all $R>0$ and $T>0$, there exists $C$, only depending on
$\vbf$, $g$, $u_0$,  $\a$, $R$, $T$ such that, for $h < R$ and $k<T$,
\be
\mapp(B(0,R)\times [0,T]) \le C (k+\sqrt{h}).
\label{estmappi} \ee

\item The measure $\mapz$ is the measure of density $\vert
u_{0}(\cdot)-u_{\mesh,0}(\cdot) \vert$ with respect to the Lebesgue 
measure and there exists
$D$, only depending on $u_0$ and $\a$, such that

\be
\mapz(\R^d) \le D h.
\label{estmapzi}
\ee

\end{enumerate}

\end{theo}

\proof of Theorem \ref{estentesti}

\smallskip

Similarly to the proof of Theorem \ref{estentest}, we introduce
a test function $\varphi\in C^\infty_c(\R^\dim\times\R_+,\R_+)$ and
a real number $\sk \in \R$. We multiply \eqref{estschemaki} by 
$(1/{\mcv}) \int_{n k}^{(n+1)k} \int_\cv \varphi(x,t) dx dt$, and
sum the result for all $\cv \in \mesh$ and $n\in\N$. 
We then define $T_1$ and $T_2$ such that $T_1 + T_2 \le 0$ 
using equations \eqref{est1kap} and \eqref{est2kap} in which
we replace $\ucvn$ by $\ucvnp$ and $\ucvvn$ by $\ucvvnp$.
Therefore we get \eqref{estmajtime}, where the measure 
$\nu_{\mesh,k}$  is such that for all $T>0$, 
there exists $C_1$  only depending on
$\vbf$, $g$, $u_0$   $\a$ and $T$, such that, for $k <T$,

\[
\nu_{\mesh,k}(\R^\dim \times [0,T]) \le C_1 k,
\]

using Lemma \refp{estlbvt}, which is available if $\vbf$ does not depend on
$t$ (and for which one needs that $u_0 \in 
BV(\R^\dim)$). 

\smallskip

The treatment of $T_2$ is very similar to that of Theorem 
\ref{estentest}, replacing $\ucvn$ by $\ucvnp$ and $\ucvvn$ by $\ucvvnp$. 
But, since $\vbf$ does not depend on
$t$, the bounds on $r_{\cv,\cvv}^{n,\pm}$ are simpler. Indeed,

\[
r_{\cv,\cvv}^{n,\pm} \le  
\dsp{ \frac 1 {k \mcv}
\int_{n k}^{(n+1)k} \int_{\cv} \int_{\edgecvcvv}} 
 \vert \varphi(x,t) - \varphi(y,t) \vert 
(\vbf(y)\cdot\bfn_{\cv,\cvv})^{\pm} \dfrontiere(y) dx dt.
\]

Now $2 r_{\cv,\cvv}^{n,\pm} \le \act{\mu_{\cv,\cvv}^{n,\pm}}{\vert \grad 
\vphi \vert}$ where ${\mu_{\cv,\cvv}^{n,\pm}}$ is
defined by

\[\bari
\act{\mu_{\cv,\cvv}^{n,\pm}}{\psi} =  
\dsp{ \frac 2 {k \mcv} \int_{n k}^{(n+1)k}  
\int_{\cv} \int_{\edgecvcvv} \int_0^1}  
\Bigl( h \; \psi(x+\theta(y-x),t) \Bigr) \\  
(\vbf(y)\cdot\bfn_{\cv,\cvv})^{\pm} d\theta \dfrontiere(y) dx dt, \, 
\forall \psi \in C_c(\RNRp).
\ean
\]

With this definition of ${\mu_{\cv,\cvv}^{n,\pm}}$, the bound on $\tilde 
{\nu}_{\mesh,k}$ (defined by \eqref{est10kap}, replacing $\ucvn$ by 
$\ucvnp$ and $\ucvvn$ by $\ucvvnp$) becomes, thanks to Lemma
\refp{estwbvesti},

\[
\tilde {\nu}_{\mesh,k}(B(0,R) \times [0,T]) \le C_2 \sqrt h,
\]

for $h<R$ and $k<T$, where $C_2$ only depends on $\vbf$, $g$, $u_0$, 
$\a$, $R$ and $T$.

Hence, defining (as in Theorem \ref{estentest}) 
$\mapp=\nu_{\mesh,k}+\tilde {\nu}_{\mesh,k}$, for all $R>0$ and all $T>0$ 
there exists  $C$, only depending on
$\vbf$, $g$, $u_0$,  $\a$, $R$, $T$ such that, for $h < R$ and $k<T$,

\[
\mapp(B(0,R)\times [0,T]) \le C (k+\sqrt{h}),
\]
which is \eqref{estmappi} and concludes the proof of Theorem 
\ref{estentesti}. \eof

\begin{remark} \normalfont 
In the case where $\vbf$ depends on $t$, Lemma \ref{estlbvt} cannot be 
used. However, it is easy to show (the proof follows that of Theorem 
\ref{estentest}) that Theorem \ref{estentesti} is true if \eqref{estmappi} is replaced
by 

\be
\mapp(B(0,R)\times [0,T]) \le C (\frac k {\sqrt{h}} + 
\sqrt{h}),
\label{estmappit} \ee

which leads to the result given in Remark \ref{estexti}. The estimate 
\eqref{estmappit} may be obtained without assuming that $u_0 \in 
BV(\R^\dim)$ (it is sufficient that $u_0 \in 
L^\infty(\R^\dim)$).
\end{remark}

For the sake of completeness we now prove  a lemma which gives the bound 
on the measure $\mapz$ in the two last theorems.

\begin{lemma}\label{divbvd}
Let $\mesh$ be an admissible mesh in the 
sense of Definition \refp{meshypmd} and let $u \in BV(\R^\dim)$ (see 
Definition \refp{defBVRp}). For $\cv \in \mesh$, let $u_\cv$ be the mean 
value of $u$ over $\cv$. Define $u_{\mesh}$ by $u_{\mesh}(x)=u_\cv$ 
for a.e. $x \in \cv$, for all $\cv \in \mesh$. Then,

\be
\norm{u-u_{\mesh}}{L^1(\R^d)} \le  \frac C {\a^2} h \snorm{u}{BV(\R^d)},
\label{bvzdd}
\ee

where $C$ only depends on the space dimension ($d=1$, $2$ or $3$).
\end{lemma}

\proof of Lemma \ref{divbvd}

\smallskip

The proof is very similar to that of Lemma \ref{divbv} and we will mainly refer 
to the proof of Lemma \ref{divbv}. 

\smallskip

First, remark that if \eqref{bvzdd} holds for all $u \in BV(\R^d)\cap 
C^1(\R^d,\R)$ then  \eqref{bvzdd} holds for all $u \in BV(\R^d)$. Indeed, let $u \in BV(\R^d)$, 
it is proven in Step 1 of the proof of Lemma \ref{divbv} that there 
exists a sequence $(u_n)_{\nnn} \subset C^\infty(\R^d,\R)$ such that $u_n 
\tends u$ in $L^1_{loc}(\R^d)$, as $n \tends \infty$, and $\norm{u_n}{BV(\R^d)} \le 
\norm{u}{BV(\R^d)}$ for all $\nnn$. One may also assume, up to a 
subsequence, that $u_n \tends u$ a.e. on 
$\R^\dim$. Then, if \eqref{bvzdd} is true with 
$u_n$ instead of $u$, passing to the limit in \eqref{bvzdd} (for $u_n$) as $n \tends 
\infty$ leads to \eqref{bvzdd} (for $u$) thanks to Fatou's lemma.

\smallskip

Let us now prove \eqref{bvzdd} if $u \in BV(\R^d)\cap C^1(\R^d,\R)$ (this 
concludes the proof of Lemma \ref{divbvd}). Since $u \in C^1(\R^d,\R)$, 
\[\snorm{u}{BV(\R^d)}=\norm{ 
(\vert \grad u \vert) }{L^1(\R^d)};\] 
hence we shall  prove \eqref{bvzdd} with $\norm{ (\vert \grad u \vert) }{L^1(\R^d)}$ 
instead of $\snorm{u}{BV(\R^d)}$.

For $\cv \in \mesh$,

\[
\int_\cv \vert u(x) - u_\cv \vert dx \le \frac 1 {\mcv} \int_\cv 
(\int_\cv \vert u(x) - 
u(y) \vert dx)dy.
\]

Then, following the lines of Step 2 of Lemma \ref{divbv},

\be
\int_\cv \vert u(x) - u_\cv \vert dx \le \frac 1 {\mcv} h \int_{B(0,h)} (\int_0^1 \int_\cv \vert \grad u 
(y + tz) \vert dy dt ) dz.
\label{oufc}
\ee

For all $\cv \in \mesh$, let $\centercv$ be an arbitrary point of $\cv$.

Then, changing the variable $y$ in $\xi=y+tz$ (for all 
fixed $z \in \cv$ and $t \in (0,1)$) in \eqref{oufc},

\[
\int_\cv \vert u(x) - u_\cv \vert dx \le \frac 1 {\mcv} h \int_{B(0,h)} (\int_0^1 
\int_{B(x_\cv,2h)} \vert \grad u 
(\xi) \vert d\xi dt ) dz,
\]

which yields, since $\mesh$ is an admissible mesh in the sense of 
 Definition \refp{estmaillage},

\[
\int_\cv \vert u(x) - u_\cv \vert dx \le \frac {1} {\a h^{d}} 
\meas(B(0,h)) h \int_{B(x_\cv,2h)} \vert \grad u 
(\xi) \vert d\xi.
\]

Therefore there exists $C_1$, only depending on the space dimension, such 
that

\be
\int_\cv \vert u(x) - u_\cv \vert dx \le \frac {C_1} {\a} h \int_{B(x_\cv,2h)} \vert \grad u 
(\xi) \vert d\xi, \, \forall \cv \in \mesh.
\label{oufd}
\ee

As in Lemma \ref{divbv}, for a fixed 
$\cvw \in \mesh$, the number of $\cv \in \mesh$ such that $\cvw \cap 
B(x_\cv,2h) 
\ne \emptyset$ is less or equal to $\meas(B(0,4h))/(\a h^d)$ that is less or 
equal to $C_2/\a$ where $C_2$ only depends on the space dimension.

Then, summing \eqref{oufd} over $\cv \in \mesh$ leads to

\[
\summesh \int_\cv \vert u(x) - u_\cv \vert dx \le 
\frac {C_1C_2} {\a^2} h \sum_{\cvw \in \mesh} \int_\cvw \vert \grad u 
(\xi) \vert d\xi= \frac {C_1C_2} {\a^2} h \norm{(\vert \grad u 
\vert)}{L^1(\R^d)},
\]

that is \eqref{bvzdd} with $C=C_1C_2$.
\eof

\section{Convergence of the scheme}\label{eghcvsol}

\hpage{Hyperbolic multid}

This section is devoted to the proof of the existence and uniqueness of
the entropy weak solution  and of the convergence of the approximate 
solution towards the entropy weak solution  as the mesh size and time 
step tend to 0. This proof will be performed in two steps. We 
first prove in section \ref{eghcvepsol} the convergence of the 
approximate solution towards an \epsol which is defined in Definition 
\ref{estepsol} below (note that the convergence also yields the 
existence of an entropy process solution).

\begin{defi}\label{estepsol}  \normalfont 
\index{Entropy!-- process solution}A function  $\mu$ 
is an entropy process solution to
problem (\ref{estpbcont})-\eqref{estcondini} if $\mu$ satisfies

\be\left\{\bal
\mu \in L^\infty(\R^\dim \times \Rpast \times (0,1)),
\\
\dsp{ \int_{\R^\dim} \int_0^{+\infty} \int_0^1 
\Bigl( \eta (\mu(x,t,\a)) \varphi_t (x,t) + 
\Phi (\mu(x,t,\a)) \vbf(x,t) \cdot  \grad  \varphi (x,t) \Bigr) d\a dt dx }
\\
+\dsp{\int_{\R^\dim} \eta(u_0 (x)) \varphi (x,0) dx} \ge 0,
\\ 
 \text{ for any }\varphi \in C^1_c(\R^\dim\times\R_+,\R_+),
\\
 \text{ for any  convex function }\eta \in C^1(\R,\R), \text{  and } 
\Phi\in C^1(\R,\R) \text{ such that } \Phi' = f' \eta'.
\ean\right.\label{ews3}
\ee
\end{defi}

\begin{remark} \normalfont\label{ewsews} 
From an entropy weak solution $u$ to problem  (\ref{estpbcont})-\eqref{estcondini}, one may
easily construct an \epsol \index{Entropy!-- process solution} to problem  
(\ref{estpbcont})-\eqref{estcondini} by setting $\mu (x,t,\a) 
= u(x,t)$ for a.e. $(x,t,\a) \in \R^\dim \times \Rpast \times (0,1)$.
Reciprocally, if $\mu$ is an \epsol to problem  
(\ref{estpbcont})-\eqref{estcondini} such that there 
exists $u \in L^\infty(\R^\dim \times \Rpast)$ such that $\mu(x,t,\a) = u(x,t)$, for a.e. $(x,t,\a) \in 
\R^\dim \times \Rpast \times (0,1)$,
then $u$ is an entropy weak solution \index{Entropy!-- weak 
solution} to problem  (\ref{estpbcont})-\eqref{estcondini}.
\end{remark}

In section \ref{ewsp}, we show the uniqueness of the entropy process 
solution, which, 
thanks to remark \ref{ewsews}, also yields the existence and uniqueness 
of the entropy weak solution. This allows us to state and prove, in 
section \ref{eghbscewsol}, the convergence of the approximate solution towards the entropy weak 
solution.

\smallskip

We now give a useful characterization of an \epsol \index{Entropy!-- process
solution} in terms of 
Krushkov's entropies (as for the entropy weak solution).

\begin{prop}\label{sanslabel}
A function $\mu$ is an entropy process solution of 
problem (\ref{estpbcont})-(\ref{estcondini}) 
if and only if,

\be\left\{\bal
\mu \in L^\infty(\R^\dim \times \Rpast \times (0,1)),
\\
\dsp{ \int_{\R^\dim} \int_0^{+\infty}  \int_0^1 
\Bigl( \vert \mu(x,t,\a) -\kappa  \vert \varphi_t (x,t) + \Phi(\mu(x,t,\a),\kappa) \vbf(x,t) \cdot 
 \grad  \varphi (x,t) \Bigr) d\a dt dx }
\\
\hfill +\dsp{\int_{\R^\dim} \vert u_0 (x) - \kappa \vert \varphi (x,0) dx} \ge 0,
\\
\hfill \forall \kappa \in \R, \, \forall \varphi \in C^1_c(\R^\dim\times\R_+,\R_+),
\end{array}\right.\label{ews3k}\ee

where we set $\Phi(a,b) = f(a\top b) - f(a\bot b)$, for all $a$, $b\in\R$.

\end{prop}

\proof of Proposition \ref{sanslabel}

\smallskip

The proof of this result is similar to the case  of classical entropy weak solutions. 
The characterization (\ref{ews3k}) can be obtained from \eqref{ews3}, by using 
regularizations
of the function $\vert \cdot - \kappa \vert$. Conversely, 
(\ref{ews3}) may be obtained from
(\ref{ews3k}) by approximating any convex function  
$\eta\in C^1(\R,\R)$ by functions 
of
the form: 
$\eta_n(\cdot)=\dsp{\sum_{i=1}^n} \a_i^{(n)} \vert \cdot - \kappa_i^{(n)} \vert$, with
$\a_i^{(n)} \ge 0$. \eof

\subsection{Convergence towards an \epsol}\index{Entropy!-- process solution|(} \label{eghcvepsol}
\index{Convergence of the approximate solutions!-- towards the entropy 
process solution|(}
Let $\a>0$ and $0< \xi< 1$. Let $(\mesh_m,k_m)_{m \in
\N}$ be a sequence of admissible meshes  in the sense of Definition 
\refp{meshypmd}  and time steps. Note that 
$\mesh_m$ is admissible with $\a$ independent of $m$. Assume that $k_m$ 
satisfies \eqref{estcfl}, for $\mesh= \mesh_m$ and $k=k_m$, and that 
$\size(\mesh_m) \tends 0$ as $m \tends \infty$.

By Lemma \refp{estinftyest}, the sequence  $(u_{\mesh_m,k_m})_{m \in
\N}$ of approximate solutions defined by the finite volume scheme 
\eqref{estschemaz} and \refep{estschema}, with  $\mesh= \mesh_m$ and $k=k_m$,
is bounded in $L^\infty(\R^\dim\times \Rpast)$; therefore, there exists 
$\mu \in L^\infty(\R^\dim\times \Rpast\times (0,1))$ such that 
$u_{\mesh_m,k_m}$ converges, as $m$ tends to $\infty$, towards $\mu$ in the
nonlinear weak-$\star$ sense (see Definition \refp{nlwsdef} and 
Proposition \refp{nlwsprop}), that is:

\be
\bal
\dsp{\lim_{m\rightarrow\infty} \int_{\R^\dim} \int_{\R_+} 
\theta(u_{\mesh_m,k_m}(x,t)) \varphi(x,t) dt dx =
\int_{\R^\dim} \int_{\R_+} \int_0^1 \theta(\mu(x,t,\a)) \varphi(x,t) d\a 
dt dx,} 
\\
\dsp{\forall \varphi \in L^1(\R\times \Rpast), \, \forall \theta \in 
C(\R,\R)}.
\ean
\label{eghbconvy}
\ee

Taking for $\theta$, in \eqref{eghbconvy}, the Krushkov entropies (namely 
$\theta=\vert \cdot -  \sk \vert$, for all $\sk \in \R$) and the associated functions defining 
the entropy fluxes (namely $\theta=f(\cdot,\sk)=f(\cdot \top \sk)-f(\cdot 
\bot \sk)$) and using Theorem 
\ref{estentest} (that is passing to the limit, as $m \tends \infty$, in 
\eqref{kruzkovap} written with $\uapp=u_{\mesh_m,k_m}$) yields that $\mu$ 
is an entropy process solution.
Hence the following result holds:

\begin{prop}\label{eghbcepsol}
Under assumptions \ref{estdonnees}, let $\a>0$ and $0< \xi< 1$. Let $(\mesh_m,k_m)_{m \in
\N}$ be a sequence of admissible meshes  in the sense of Definition 
\refp{meshypmd}  and time steps. Note that 
$\mesh_m$ is admissible with $\a$ independent of $m$. Assume that $k_m$ 
satisfy \eqref{estcfl}, for $\mesh= \mesh_m$ and $k=k_m$, and that 
$\size(\mesh_m) \tends 0$ as $m \tends \infty$.

Then there exists a subsequence, still denoted by  
$(\mesh_m,k_m)_{m \in
\N}$, and a function $\mu \in L^{\infty}(\R^\dim \times \Rpast \times 
(0,1))$ such that

\begin{enumerate}
\item the approximate solution defined by \eqref{estschema}, \eqref{estschemaz} and 
\eqref{estuapp} with $\mesh=\mesh_m$ and $k=k_m$, that is
$u_{\mesh_m,k_m}$, converges towards $\mu$ in the nonlinear weak-$\star$ 
sense, i.e. \eqref{eghbconvy} holds,
\item $\mu$ is an \epsol of (\ref{estpbcont})-(\ref{estcondini}).
\end{enumerate}
\end{prop}

\begin{remark} \normalfont
The same theorem can be proved for the implicit scheme without condition 
\eqref{estcfl} (and thus without $\xi$).
\end{remark}

\begin{remark} \normalfont\label{eghbexepsol} 

Note that a consequence of Proposition \ref{eghbcepsol} is the existence of an \epsol to 
Problem
(\ref{estpbcont})-(\ref{estcondini}).

\end{remark}
\index{Convergence of the approximate solutions!-- towards the entropy 
process solution|)}

\subsection{Uniqueness of the \epsol} \label{ewsp}

In order to show the uniqueness of an entropy process 
solution, we shall use the  
characterization  of an  entropy process solution given in proposition \ref{sanslabel}. 

\index{Uniqueness!-- of the entropy process solution to a hyperbolic equation}
\begin{theo}\label{ewsunicite}
Under Assumption \ref{estdonnees}, the entropy process solution 
$\mu$ of problem
(\ref{estpbcont}),(\ref{estcondini}), as defined in Definition 
\refp{estepsol}, is unique. 
Moreover, there exists a
function $u\in
L^\infty(\R^\dim\times\Rpast)$ such that $u(x,t) = \mu(x,t,\a)$, for a.e. 
$(x,t,\a) \in \R^\dim \times \Rpast \times(0,1)$. 
(Hence, with Proposition \ref{eghbcepsol} and Remark \ref{ewsews}, there exists a  
unique entropy weak solution  to Problem 
(\ref{estpbcont})-(\ref{estcondini}).) 
\end{theo}

\proof of Theorem \ref{ewsunicite}

\smallskip

Let $\mu$ and $\nu$ be two entropy process solutions to 
Problem (\ref{estpbcont})-(\ref{estcondini}). Then, one has
$\mu \in L^\infty(\R^\dim \times \Rpast \times (0,1))$, $\nu \in 
L^\infty(\R^\dim \times \Rpast \times (0,1))$ and

\be\bal
\dsp{ \int_{\R^\dim} \int_0^{+\infty}  \int_0^1 
\Bigl( \vert \mu(x,t,\a) -\kappa  \vert \varphi_t (x,t)}
\\
\dsp{ + (f(\mu(x,t,\a) \top \kappa)-f(\mu(x,t,\a) \bot \kappa)) \vbf(x,t) \cdot 
 \grad  \varphi (x,t) \Bigr) d\a dt dx }
\\
\dsp{+\int_{\R^\dim} \vert u_0 (x) - \kappa \vert \varphi (x,0) dx \ge 0, 
\,  \forall \kappa \in \R, \, \forall \varphi \in  
C^1_c(\R^\dim\times\R_+,\R_+),}
\end{array}\label{ewsmu}\ee

\be\bal
\dsp{ \int_{\R^\dim} \int_0^{+\infty}  \int_0^1 
\Bigl( \vert \nu(y,s,\b) -\kappa  \vert \varphi_s (y,s)}
\\
\dsp{ + (f(\nu(y,s,\b) \top \kappa)-f(\nu(y,s,\b) \bot \kappa)) \vbf(y,s) \cdot 
 \grad  \varphi (y,s) \Bigr) d\b ds dy }
\\
\dsp{+\int_{\R^\dim} \vert u_0 (y) - \kappa \vert \varphi (y,0) dy \ge 0,
\,  \forall \kappa \in \R, \, \forall \varphi \in 
C^1_c(\R^\dim\times\R_+,\R_+).}
\end{array}\label{ewsnu}\ee

\smallskip

The proof of Theorem \ref{ewsunicite} contains 2 steps. In Step 1, it is 
proven that

\be\bal
\dsp{\int_0^1 \int_0^1 \int_{\Rp}\int_{\R^\dim} 
\Bigl[ \vert \mu(x,t,\a) - \nu(x,t,\b) \vert  \psi_t(x,t)}
\\
\dsp{ +\Bigl(f(\mu(x,t,\a) \top \nu(x,t,\b))-f(\mu(x,t,\a) \bot 
\nu(x,t,\b))\Bigr) \vbf(x,t) \cdot \grad \psi(x,t) \Bigr] 
dx dt d\a d\b \ge 0,}
\\
\dsp{\forall \psi \in C^1_c(\R^\dim \times \Rp,\Rp).}
\ean
\label{step1}
\ee

In Step 2, it is proven that $\mu(x,t,\a)=\nu(x,t,\b)$ for a.e. 
$(x,t,\a,\b) \in \R^\dim \times \Rpast \times (0,1) \times (0,1)$. We 
then deduce that there exists $u \in L^\infty(\R^\dim \times \Rpast)$ such that 
$\mu(x,t,\a)=u(x,t)$ for a.e. 
$(x,t,\a) \in \R^\dim \times \Rpast \times (0,1)$ (therefore $u$ is 
necessarily the unique entropy weak solution   to 
(\ref{estpbcont})-(\ref{estcondini})).

\smallskip

\textit{Step 1 (proof of relation \eqref{step1})}

In order to prove relation \eqref{step1}, a sequence of  mollifiers in 
$\R$ and $\R^\dim$ is introduced .

Let $\rho \in C_c^\infty(\R^\dim,\Rp)$ and $\bar{\r} \in 
C_c^\infty(\R,\Rp)$ be such that

\[
\{x \in \R^\dim; \,  \rho(x) \neq 0\} \subset \{x 
\in \R^\dim; \, \vert x \vert \le 1\},
\]

\be
\{x \in \R; \,  
\bar{\r}(x) \neq 0\} \subset [-1,0]
\label{rhodec}
\ee

and

\[
\int_{\R^\dim} \rho(x) dx = 1, \; \int_{\R} \bar{\r}(x) dx = 1.
\]

For $n\in\N^\star$, define $\rho_{n} = n^\dim \rho(nx)$ for all $x 
\in \R^\dim$ and $\bar{\r}_{n} = n \bar{\r}(nx)$ for all $x 
\in \R$.

Let $\psi \in C^1_c(\R^\dim \times \Rp,\Rp)$. For $(y,s,\b) \in \R^\dim 
\times \Rp \times (0,1)$, let us take, in \eqref{ewsmu},  
$\vphi(x,t)=\psi(x,t) \r_{n}(x-y) \bar{\r}_{n}(t-s)$ and $\sk = 
\nu(y,s,\b)$. Then, integrating the result over $\R^\dim 
\times \Rp \times (0,1)$ leads to

\be
A_1 + A_2 + A_3 + A_4 + A_5 \ge 0,
\label{ewscin}
\ee

where

\[
\bari
\dsp{ A_1 = \int_0^1\int_0^1\int_0^\infty \int_{\R^\dim} \int_0^\infty 
\int_{\R^\dim}
\Bigl[ \vert \mu(x,t,\a) - \nu(y,s,\b) \vert}
\\
\dsp{\psi_t(x,t) \r_{n}(x-y) \bar{\r}_{n}(t-s)  \Bigr] dx dt dy ds d\a 
d\b,}
\ean
\]

\[
\bari
\dsp{ A_2 = \int_0^1\int_0^1\int_0^\infty \int_{\R^\dim} \int_0^\infty 
\int_{\R^\dim}
\Bigl[ \vert \mu(x,t,\a) - \nu(y,s,\b) \vert}
\\
\dsp{\psi(x,t) \r_{n}(x-y) \bar{\r}_{n}'(t-s)  \Bigr] dx dt dy ds d\a 
d\b,}
\ean
\]

\[
\bari
\dsp{ A_3 = \int_0^1\int_0^1\int_0^\infty \int_{\R^\dim} \int_0^\infty 
\int_{\R^\dim} \Bigl[
\Bigl(f(\mu(x,t,\a)\top \nu(y,s,\b)) -  f(\mu(x,t,\a)\bot 
\nu(y,s,\b))\Bigr)}
\\ 
\dsp{\vbf(x,t)\cdot\grad\psi(x,t) \r_{n}(x-y) \bar{\r}_{n}(t-s)  
\Bigr] dx dt dy ds d\a d\b,} 
\ean
\]

\[
\bari
\dsp{ A_4 = \int_0^1\int_0^1\int_0^\infty \int_{\R^\dim} \int_0^\infty 
\int_{\R^\dim} \Bigl[
\Bigl(f(\mu(x,t,\a)\top \nu(y,s,\b)) -  f(\mu(x,t,\a)\bot 
\nu(y,s,\b))\Bigr)}
\\ 
\dsp{\vbf(x,t)\cdot\grad\r_{n}(x-y) \psi(x,t) \bar{\r}_{n}(t-s)  
\Bigr] dx dt dy ds d\a d\b} 
\ean
\]

and

\[
A_5 =\int_0^1 \int_{\R^\dim} \int_0^\infty \int_{\R^\dim} \vert u_0(x) - 
\nu(y,s,\b) \vert
\psi(x,0) \r_{n}(x-y) \bar{\r}_{n}(-s) dydsdxd\b.
\label{a5}
\]

Passing to the limit in \eqref{ewscin} as $n \tends \infty$ (using 
\eqref{ewsnu} for the study of $A_2+A_4$ and $A_5$) will give \eqref{step1}.

\smallskip

Let us first consider $A_1$ and $A_3$. Note that, using \eqref{rhodec},

\[
\int_{\R^\dim} \int_0^\infty 
\rho_{n}(x-y) \bar{\r}_{n}(t-s)dsdy=1, \, 
\forall x \in \R^\dim, \, \forall t \in \R_+.
\]

Then,

\[
\bal
\dsp{ \vert A_1 - \int_0^1\int_0^1\int_{\Rp}\int_{\R^\dim} 
\Bigl[ \vert \mu(x,t,\a) - \nu(x,t,\b) \vert  
\psi_t(x,t) \Bigr] \ dx dt d\a d\b \vert}
\\
\dsp{\le \int_0^1 \int_0^\infty \int_{\R^\dim} \int_0^\infty \int_{\R^\dim} 
\Bigl[ \vert \nu(x,t,\b) - \nu(y,s,\b) \vert  
\vert\psi_t(x,t)\vert \r_{n}(x-y) \bar{\r}_{n}(t-s)  \Bigr] dx dt dy ds d\b}
\\
\le \Vert \psi_t \Vert_{L^\infty(\R^\dim \times \Rpast)} \eps(n,S),
\ean
\]

with $S=\{(x,t) \in \R^\dim \times \Rp$; $\psi(x,t) \ne 0\}$ and

\[
\eps(n,S)=\sup\{ \norm{\nu-\nu(\cdot+\eta,\cdot+\tau,\cdot)}{L^1(S\times 
(0,1))}; \; \vert \eta \vert \le \frac 1
n, \, 0 \le \tau \le \frac 1 n\}.
\]

Since $\nu \in L^1_{loc}(\R^\dim\times \Rp \times [0,1])$ and $S$ is 
bounded, one has $\eps(n,S) \tends 0$ as $n \tends \infty$. Hence,

\be
A_1 \tends \int_0^1 \int_0^1\int_{\Rp}\int_{\R^\dim}
\Bigl[ \vert \mu(x,t,\a) - \nu(x,t,\b) \vert  
\psi_t(x,t) \Bigr] dx dt d\a d\b, \text{ as } n \tends \infty.
\label{lima1}
\ee

Similarly, let $M$ be the Lipschitz constant of $f$ on $[-D,D]$ where 
$D=\max\{\norm{\mu}{\infty},\norm{\nu}{\infty}\}$, with 
$\norm{\cdot}{\infty}=\norm{\cdot}{L^\infty(\R^\dim \times \Rpast \times 
(0,1))}$,

\[
\bari
\dsp{\vert A_3 - \int_0^1\int_0^1\int_{\Rp}\int_{\R^\dim}
\Bigl(f(\mu(x,t,\a)\top \nu(x,t,\b)) -  f(\mu(x,t,\a)\bot 
\nu(x,t,\b))\Bigr)}
\\ 
\dsp{\vbf(x,t)\cdot\grad\psi(x,t) dx dt d\a d\b \vert \le 2MV \norm{(\vert \grad \psi 
\vert)}{L^\infty(\R^\dim \times \Rpast)} \, \eps(n,S),}
\ean
\]

which yields

\be
\bari
\dsp{A_3 \tends \int_0^1\int_0^1\int_{\Rp}\int_{\R^\dim}
\Bigl(f(\mu(x,t,\a)\top \nu(x,t,\b)) -  f(\mu(x,t,\a)\bot 
\nu(x,t,\b))\Bigr)}
\\ 
\dsp{\vbf(x,t)\cdot\grad\psi(x,t) dx dt d\a d\b, \text{ as } n \tends 
\infty.}
\ean
\label{lima3}
\ee

\smallskip

Let us now consider $A_2+A_4$.

For $(x,t,\a) \in \R^\dim \times \Rp \times (0,1)$, let us take 
$\vphi(y,s)=\psi(x,t) \r_{n}(x-y) \bar{\r}_{n}(t-s)$ and 
$\sk = \mu(x,t,\a)$
in \eqref{ewsnu}. Integrating the result over $\R^\dim 
\times \Rp \times (0,1)$ leads to

\be
-A_2-B_4 \ge 0,
\label{a2b4}
\ee

 with 

\[
\bari
\dsp{A_4-B_4 = \int_0^1\int_0^1\int_0^\infty \int_{\R^\dim} \int_0^\infty 
\int_{\R^\dim} \Bigl[
\Bigl(f(\mu(x,t,\a)\top \nu(y,s,\b)) -  f(\mu(x,t,\a)\bot 
\nu(y,s,\b))\Bigr)}
\\ 
\dsp{(\vbf(x,t)-\vbf(y,s))\cdot\grad\r_{n}(x-y) \psi(x,t) \bar{\r}_{n}(t-s)  
\Bigr] dx dt dy ds d\a d\b.} 
\ean
\]

Note that $B_4=A_4$ if $\vbf$ is constant (and one directly obtains \eqref{lima24} below).
In the general case, in order to prove that $A_4-B_4 \tends 0$ as $n 
\tends \infty$ (which then gives \eqref{lima24}), let us remark 
that, using $\div \vbf=0$,

\be
\bari
\dsp{\int_0^1\int_0^1\int_0^\infty \int_{\R^\dim} \int_0^\infty 
\int_{\R^\dim} \Bigl[
\Bigl(f(\mu(x,t,\a)\top \nu(x,t,\b)) -  f(\mu(x,t,\a)\bot 
\nu(x,t,\b))\Bigr)}
\\ 
\dsp{(\vbf(x,t)-\vbf(y,s))\cdot\grad\r_{n}(x-y) \psi(x,t) \bar{\r}_{n}(t-s)  
\Bigr] dx dt dy ds d\a d\b=0.} 
\ean
\label{abc4}
\ee

Indeed, the latter equality follows from an integration by parts for the 
variable $y \in \R^\dim$. Then, substracting the left hand side of 
\eqref{abc4} to $A_4-B_4$ and using the regularity of $\vbf$, there 
exists $C_1$, only depending on $M$, $\vbf$ and $\psi$, such that $\vert A_4-B_4 
\vert \le C_1 \eps(n,S)$. This gives $A_4-B_4 \tends 0$ as $n 
\tends \infty$ and, thanks to \eqref{a2b4},

\be
\limsup_{n \tends \infty} (A_2 + A_4) \le 0.
\label{lima24}
\ee

\smallskip

Finally, let us consider $A_5$.

For $x \in \R^\dim$, let us take $\vphi(y,s) = \psi(x,0) 
\r_{n}(x-y) \int_s^\infty \bar{\r}_{n}(-\tau)d\tau$ and $\sk = u_0(x)$ in 
\eqref{ewsnu}. Integrating the resulting inequality with respect to $x \in 
\R^\dim$ gives

\be
-A_5 + B_{5a} + B_{5b} \ge 0,
\label{a5c}
\ee

with

\[
\bari
\dsp{B_{5a} = -\int_0^1 \int_0^\infty \int_{\R^\dim} \int_{\R^\dim} \int_s^\infty}
(f(\nu(y,s,\b)\top u_0(x)) - f(\nu(y,s,\b)\bot u_0(x)))
\\
\vbf(y,s) \cdot \grad \r_{n}(x-y) \psi(x,0) \bar{\r}_{n}(-\tau)  d\tau 
dy dx ds d\b,
\ean
\]

\[
B_{5b} = \int_{\R^\dim} \int_{\R^\dim} \psi(x,0) \r_{n}(x-y) 
\vert u_0(x) - u_0(y) \vert dy dx.
\]

Let $S_0=\{x \in \R^\dim$; $\psi(x,0) \ne 0\}$ and 

\[
\eps_0(n,S_0)=\sup\{\int_{S_0}\vert
u_0(x)-u_0(x+\eta)\vert dx; \; \vert\eta\vert \le \frac 1 n\},
\]

so that $B_{5b} \le \norm{\psi(\cdot,0)}{L^\infty(\R^\dim)}\eps_0(n,S_0)$.

Since $u_0 \in L^1_{loc}(\R^\dim)$ and since $S_0$ is bounded, one has 
$\eps_0(n,S_0) \tends 0$ as $n \tends \infty$. Then, $B_{5b} \tends 0$ 
as $n \tends \infty$.

\smallskip

Let us now prove that $B_{5a} \tends 0$ as $n \tends \infty$ (then, 
\eqref{a5c} will give 
\eqref{lima5} below). Note that $B_{5a}=-B_{5c} + (B_{5a}+B_{5c})$ with

\[
\bari
\dsp{B_{5c} =  \int_0^1 \int_0^\infty \int_{\R^\dim} \int_{\R^\dim} \int_s^\infty} 
(f(\nu(y,s,\b)\top u_0(y)) - f(\nu(y,s,\b)\bot u_0(y)))
\\
\vbf(y,s) \cdot \grad \r_{n}(x-y) \psi(x,0)  \bar{\r}_{n}(-\tau) d\tau dy 
dx ds d\b.
\ean
\]

Integrating by parts for the $x$ variable yields

\[
\bari
\dsp{B_{5c} =  \int_0^1 \int_0^\infty \int_{\R^\dim} \int_{\R^\dim} \int_s^\infty} 
(f(\nu(y,s,\b)\top u_0(y)) - f(\nu(y,s,\b)\bot u_0(y)))
\\
\vbf(y,s) \cdot \grad \psi(x,0)  \r_{n}(x-y) \bar{\r}_{n}(-\tau) d\tau dy 
dx ds d\b.
\ean
\]

Noting that the integration with 
respect to $s$ is reduced to $[0,1/n]$, $B_{5c} \tends 0$ as $n \tends \infty$.

There remains to study $B_{5a}+B_{5c}$. Noting that $\vert f(a \top b)-f(a \top 
c) \vert \le \bar{M} \vert b-c \vert$ and $\vert f(a \bot b)-f(a \bot 
c) \vert \le \bar{M} \vert b-c \vert$ if  $b$, $c \in 
[-\bar{D},\bar{D}]$, where $\bar{D}=\norm{u_0}{L^\infty(\R^\dim)}$ and 
$\bar{M}$ is the  Lipschitz constant to $f$ on $[-\bar{D},\bar{D}]$,

\[
\vert B_{5a}+B_{5c} \vert \le 2 \bar{M} V
\int_0^\infty \int_{\R^\dim} \int_{\R^\dim} \int_s^\infty  \vert u_0(x) -
u_0(y)\vert  \vert \grad\r_{n}(x-y) \vert \psi(x,0) 
\bar{\r}_{n}(-\tau) d\tau dy dx ds,
\]

which yields the existence of $C_2$, only depending on $\bar{M}$, $V$ and 
$\psi$, such that

\[
\vert B_{5a}+B_{5c} \vert \le C_2
\int_0^{\frac 1 n} \int_{S_0} \int_{B(0,\frac 1 n)} \vert u_0(x) -
u_0(x-z)\vert  n^{\dim + 1}  dz dx ds.
\]

Therefore, $\vert B_{5a}+B_{5c} \vert \le C_3 \eps_0(n,S_0)$, with some 
$C_3$ only depending on $\bar{M}$, $V$ and 
$\psi$. Since $\eps_0(n,S_0) \tends 0$ as $n \tends \infty$,  one deduces 
$\vert B_{5a}+B_{5c} \vert \tends 0$ as $n \tends \infty$. Hence, $B_{5a} 
\tends 0$ as $n \tends \infty$ and \eqref{a5c} yields

\be
\limsup_{n \tends \infty} A_5 \le 0.
\label{lima5}
\ee

It is now possible to conclude Step 1. Passing to the limit as $n 
\tends \infty$ in \eqref{ewscin} and using 
\eqref{lima1}, \eqref{lima3}, \eqref{lima24} and \eqref{lima5} yields 
\eqref{step1}.

\smallskip

\textit{Step 2 (proof of $\mu=\nu$ and conclusion)}

Let $R>0$ and $T>0$. One sets $\o = V M$ (recall that $V$ is given in 
Assumption \ref{estdonnees} and that $M$ is given in Step 1).

Let $\vphi \in C_c^1(\R_+,[0,1])$  be a function such that
$\vphi(r) = 1$ if  $r \in [0,R+\o T]$, $\vphi(r) = 0$ if
$r \in [R+\o T +1,\infty)$ and $\vphi'(r) \le 0$, for all $r \in \R_+$.

One takes, in \eqref{step1}, $\psi$ defined by 

\[
\ba
\psi(x,t) = \vphi(\vert x \vert + \o t) \frac {T-t} T, \text{ for } x 
\in \R^\dim  \text{ and } t\in[0,T],
\\
\psi(x,t)=0,  \text{ for } x \in \R^\dim  \text{ and } t \ge T.
\ea
\]

The function $\psi$ is not in
$C^\infty_c(\RNRp,\R_+)$, but, using a usual regularization technique, it 
may be proved that such a function can be considered in
\eqref{step1}, in which case Inequality \eqref{step1} reads

\[
\bal
\dsp{\int_0^1 \int_0^1 \int_0^T \int_{\R^\dim} 
\Bigl[ \vert \mu(x,t,\a) - \nu(x,t,\b) \vert \Bigl(\frac {T-t} T
 \o \vphi'(\vert x \vert + \o t) 
- \frac 1 T \vphi(\vert x \vert + \o t) \Bigr)
+}
\\
\dsp{\Bigl(f(\mu(x,t,\a)\top \nu(x,t,\b)) - f(\mu(x,t,\a)\bot 
\nu(x,t,\b))\Bigr) \frac {T-t} T \vphi'(\vert x \vert + \o t)
\vbf(x,t) \cdot \frac x {\vert x \vert} \Bigr] dx dt d\a d\b \ge 0.}
\ean
\]

Since $\o=VM$ and $\vphi'\le 0$, one has $(f(a \top b) -  f (a\bot b)) 
\vphi'(\vert x \vert + \o t) \vbf(x,t) \cdot (x/\vert x \vert)  \le \vert a-b \vert \o (-\vphi'(\vert 
x \vert + \o t))$, for a.e. $(x,t) \in \R^\dim \times \Rpast$ and all $a$, 
$b \in [-D,D]$ ($D$ is defined in Step 1).
Therefore, the previous inequality gives
\begin{equation}\label{ine-wphi}
\int_0^1 \int_0^1 \int_0^T \int_{\R^\dim} 
 \vert \mu(x,t,\a) - \nu(x,t,\b) \vert 
 \frac 1 T\vphi(\vert x \vert + \o t) dx dt d\a d\b
\le 0.
\end{equation}
Then, since $\vphi(\vert x \vert + \o t)=1$ if $(x,t) \in B(0,R) 
\times [0,T]$,
\[
\int_0^1 \int_0^1 \int_0^T \int_{B(0,R)}
\vert \mu(x,t,\a) - \nu(x,t,\b) \vert dx dt d\a d\b \le 0,
\]
which yields, since $R$ and $T$ are arbitrary, $\mu(x,t,\a) = 
\nu(x,t,\b)$ for a.e. $(x,t,\a,\b) \in \R^\dim \times \Rpast \times (0,1) 
\times (0,1)$.

\smallskip

Let us now deduce also from this uniqueness result that there exists $u 
\in L^\infty(\R^\dim \times \Rpast)$ such that $\mu(x,t,\a) = u(x,t)$, 
for a.e. $(x,t,\a) \in \R^\dim \times \Rpast \times (0,1)$ (then it is easy 
to see, with Definition \ref{estepsol}, that $u$ is the entropy weak 
solution  to Problem (\ref{estpbcont})-(\ref{estcondini})).


Indeed, it is possible to take, in the preceeding proof, $\mu=\nu$ (recall
that the proposition \ref{eghbcepsol} gives the existence of an entropy 
process solution to Problem (\ref{estpbcont})-(\ref{estcondini}), see Remark \ref{eghbexepsol}).
This
yields $\mu(x,t,\a) = 
\mu(x,t,\b)$ for a.e. $(x,t,\a,\b) \in \R^\dim \times \Rpast \times (0,1) 
\times (0,1)$. Then, for a.e. $(x,t) \in \R^\dim \times \Rpast$, one has 

\[
\mu(x,t,\a)= \mu(x,t,\b) \text{ for a.e. } (\a,\b) \in (0,1) 
\times (0,1)
\]

and, for a.e. $\a \in (0,1)$,

\[
\mu(x,t,\a)= \mu(x,t,\b) \text{ for a.e. } \b \in (0,1).
\]

Thus, defining $u$ from $\R^\dim \times \Rpast$ to $\R$ by

\[
u(x,t)= \int_0^1 \mu(x,t,\b) d\b,
\]

one obtains $\mu(x,t,\a) = u(x,t)$, 
for a.e. $(x,t,\a) \in \R^\dim \times \Rpast \times (0,1)$, and
$u$ is the entropy weak 
solution  to Problem (\ref{estpbcont})-(\ref{estcondini}). 
This completes 
the proof of Theorem \ref{ewsunicite}. \eof


The proof of Theorem  \ref{ewsunicite} may be adapted in order to prove the principle of ``finite speed of propagation" for the solution to a hyperbolic equation,
namely Proposition \ref{prop-pvf}.

\begin{prop}\label{prop-pvf}\index{Finite speed of propagation!}
Under Assumption \ref{estdonnees}, let $u$ and $v$ be the entropy weak solutions 
of (\ref{estpbcont}),(\ref{estcondini}), where $u_0$ is replaced by $v_0$ for $v$, (see \refp{def-ewsmd} for the definition of entropy weak solution). 
Then, for any $x_0 \in \R^d$, $R>0$ and $T>0$,
\begin{equation}\label{ine-pafs}
\int_0^T \int_{B(x_0,R)} |u(x,t)-v(x,t)| dx dt \le T \int_{B(x_0,R+\o T)}|u_0(x)-v_0(x)| dx,
\end{equation}
where, $\o=VM$ and $B(x,a)=\{y \in \R^d, |y-x| <a\}$.
\end{prop}
\proof of Proposition \ref{prop-pvf}
Let $R>0$ and $T>0$.
Taking into account the fact that $v_0 \ne u_0$ in the proof of Theorem \ref{ewsunicite}, Inequality \ref{ine-wphi} becomes
\[
\int_0^T \int_{\R^\dim} 
\vert u(x,t) - v(x,t) \vert 
\frac 1 T\vphi(\vert x \vert + \o t)dx dt - \int_{\R^d} [u_0(x)-v_0(x)| \vphi(x)dx
\le 0.
\]
Then, the choice of $\vphi$ gives
\[
\int_0^T \int_{B_R} 
 \vert u(x,t) - v(x,t) \vert dx dt \le T \int_{B_{R+\o T +1}}  [u_0(x)-v_0(x)| \vphi(x)dx.
\]
In this inequality, it is quite easy to replace $B_{R+\o T +1}$ by $B_{R+\o T +\eps}$ for any $\eps>0$ (changing the choice of $\vphi$ in the proof).
Then, letting $\eps \to 0$,
\[
\int_0^T \int_{B_R} 
 \vert u(x,t) - v(x,t) \vert 
 \vphi(\vert x \vert + \o t)dx dt \le T \int_{B_{R+\o T}}  [u_0(x)-v_0(x)| \vphi(x)dx.
\]
Finally, the point $0$ may be replaced by any point $x_0$ in the proof and we obtain \ref{ine-pafs}.
\eof

\begin{remark}\label{rem-pvf} \index{Finite speed of propagation!}\normalfont
Proposition \ref{prop-pvf} gives the so-called ``Finite speed of propagation'' for a hyperbolic equation.
Namely it gives (under the hypothesis of Proposition \ref{prop-pvf}) that if $u_0=v_0$ a.e. $B_{R+\o T}$, then
$u=v$ a.e. in $B_R \times ]0,T[$.

Another interesting consequence of the proof of Theorem \ref{ewsunicite} (changing conveniently the function $\psi$ in the proof) is the continuity in time of the entropy weak solution of (\ref{estpbcont}),(\ref{estcondini}).
More precisely, if $u$ is the entropy weak solution of (\ref{estpbcont}),(\ref{estcondini}), then, there exists $\tilde u \in C([0,+\infty,$ $L^1_{loc}(\R^d))$ such that $\tilde u = u$ a.e on $R^d \times \R_+$.
Moreover, if $\tilde v$ is the continuous in time entropy weak solution corresponding to the initial condition $v_0$ (as in Proposition \ref{prop-pvf}), the proof of Proposition \ref{prop-pvf} (changing conveniently the choice of $\psi$) gives
\[
 \int_{B(x_0,R)} |\tilde u(x,T)-\tilde v(x,T)| dx dt \le  \int_{B(x_0,R+\o T)}|u_0(x)-v_0(x)| dx,
\]
\end{remark}


\index{Entropy!-- process solution|)}
\subsection{Convergence towards the entropy weak solution}\label{eghbscewsol}

We now know that there exists a unique entropy process solution to problem 
\eqref{estpbcont}-\refep{estcondini}, which is identical to  the entropy weak 
solution of problem \eqref{estpbcont}-\eqref{estcondini}; we may now 
prove the convergence of the approximate solution given by the finite 
volume scheme  \eqref{estschema}, \eqref{estschemaz} and 
\eqref{estuapp} towards the entropy weak solution as the mesh size tends 
to 0.

\begin{theo}\label{eghbcewsol}
\index{Convergence of the approximate solutions!-- towards the entropy weak 
solution}Under Assumptions \refp{estdonnees}, let $\a \in \Rpast$ and $\xi \in 
(0,1)$ be given. For an admissible mesh $\mesh$ in the sense of Definition 
\refp{meshypmd} and for $k>0$ satisfying 
\eqref{estcfl} (note that $\a$ and $\xi$ are fixed), let $u_{\mesh,k}$ be 
the solution to \eqref{estschema}, \eqref{estschemaz} and 
\eqref{estuapp}.

Then, $u_{\mesh,k} \tends u$ in $L^p_{loc}(\R^\dim \times\R_+)$ for all $p \in [1,\infty)$, as 
$h=\sizemesh \tends 0$, where $u$ is the entropy
weak solution  to \eqref{estpbcont}-\refep{estcondini}. 

\end{theo}

\proof of Theorem \ref{eghbcewsol}

\smallskip

In order to prove that $u_{\mesh,k} \tends u$ (in $L^p_{loc}(\R^\dim 
\times\R_+)$ for all $p \in [1,\infty)$, as 
$h=\sizemesh \tends 0$), let us proceed by a classical way of 
contradiction which uses the uniqueness of the entropy process solution  
to Problem \eqref{estpbcont}-\refep{estcondini}. 
Assume that there exists $1 \le p_0 < \infty$, 
$\varepsilon > 0$, $\comp$ a compact subset of $\R^\dim$, $T>0$ and a
sequence $((\mesh_m, k_m))_{m \in \N}$ such that, for any $m \in \N$, 
$\mesh_m$ is an admissible mesh, $k_m$ satisfies  \eqref{estcfl} (with 
$\mesh=\mesh_m$ and $k=k_m$, note that $\a$ and $\xi$ are independent of 
$m$), $\size(\mesh_m) \tends 0$ as $m \rightarrow \infty$ and

\be
\int_0^T \int_\comp \vert u_{\mesh_m,k_m} - u \vert^{p_0} dx dt \ge \varepsilon, 
\, \forall m \in \N,
\label{eghbcontra}
\ee

where $u_{\mesh_m,k_m}$ is the solution to \eqref{estschema}, \eqref{estschemaz} and 
\eqref{estuapp} with $\mesh=\mesh_m$ and $k=k_m$ and $u$ is the entropy
weak solution  to \eqref{estpbcont}-\eqref{estcondini}.

Using Proposition \ref{eghbcepsol}, there exists a subsequence of the sequence 
$((\mesh_m, k_m))_{m
\in \N}$, still denoted by  $((\mesh_m,$ $k_m))_{m \in \N}$, and a function $\mu \in
L^{\infty}(\R^\dim \times \Rpast \times (0,1))$ such that 

\begin{enumerate}
\item $u_{\mesh_m,k_m} \tends \mu$, as $m \tends \infty$, in the nonlinear weak-$\star$ 
sense, that is:

\be\bari
\dsp{\lim_{m\rightarrow\infty} \int_0^\infty \int_{\R^\dim} \theta(u_{\mesh_m,k_m}(x,t)) \varphi(x,t) dx dt =
\int_0^1 \int_0^\infty \int_{\R^\dim} \theta(\mu(x,t,\a)) \varphi(x,t) dx 
dt d\a,}
\\
\forall \varphi \in L^1(\R^\dim \times \Rpast), \, \forall \theta \in 
C(\R,\R),
\ean
\label{eghbww}
\ee
\item $\mu$ is an \epsol to (\ref{estpbcont})-(\ref{estcondini}). 
\end{enumerate}

By Theorem \refp{ewsunicite},
one has $\mu(\cdot,\cdot,\a)=u$, for a.e. $\a \in [0,1]$ (and $u$ is the
entropy weak solution to (\ref{estpbcont})-(\ref{estcondini})). 
Taking first $\theta(s)=s^2$ in (\ref{eghbww}) and then
$\theta(s)=s$ and $\varphi u$ instead of $\varphi$ in (\ref{eghbww}) 
one obtains:

\be
\int_0^\infty \int_{\R^\dim} (u_{\mesh_m,k_m}(x,t)) - u(x,t))^2
\varphi(x,t) dx dt \tends 0, \text{ as } m \tends \infty,
\label{eghbfort}
\ee
for any  function $\varphi \in L^1(\RdR)$. From (\ref{eghbfort}), and thanks to the $L^{\infty}$-bound on
$(u_{\mesh_m,k_m})_{m \in \N}$, one deduces the convergence of
$(u_{\mesh_m,k_m})_{m \in \N}$ towards $u$ in $L^p_{loc}(\R^\dim \times
\R_+)$ for all $p \in [1,\infty)$, which is in contradiction with
(\ref{eghbcontra}).

This completes the proof of our convergence theorem. \eof

\begin{remark} \normalfont~
\benum
\item Theorem \ref{eghbcewsol} is also true with the implicit scheme instead of 
the explicit scheme (that is \eqref{estschemai} and 
\eqref{estuappimp} instead of \eqref{estschema} and 
\eqref{estuapp}) without the condition \eqref{estcfl} (and thus without 
$\xi$).
\item The following section improves this convergence result and 
gives an error estimate.
\eenum
\end{remark}

\section{Error estimate}\label{esterr}

\hpage{Hyperbolic multid}

\subsection{Statement of the results}

This section is devoted to the proof of an
error estimate of time explicit and time implicit finite volume 
approximations to the solution 
$u\in L^\infty(\R^\dim\times\Rpast)$ of Problem 
\eqref{estpbcont}-\refep{estcondini}. Assuming that  $u_0 \in 
BV(\R^\dim)$, a ``$h^{1/4}$" error estimate 
is shown for a large variety of finite volume monotone flux schemes such 
as those which were presented in Section \refp{eghbns}. 

\medskip

Under Assumption \refp{estdonnees}, let $\mesh$ be an admissible mesh in the 
sense of Definition \refp{meshypmd} and $k>0$. Let $g \in C(\R^2,\R)$ 
satisfying
Assumption \ref{estcondflu}. 

Let $u$ be the entropy weak solution   of \eqref{estpbcont}-\eqref{estcondini} and let $\uapp$ be the solution of the 
time explicit scheme \eqref{estschema},
\eqref{estschemaz}, \eqref{estuapp}, assuming that \eqref{estcfl} holds, or 
$\uapp$ be the solution of the time implicit scheme \eqref{estschemai},
\eqref{estschemaz}, \eqref{estuappimp}. Our aim is to give an
error estimate between $u$ and $\uapp$.

\medskip

In the case of the explicit scheme, one proves, in this section,
the following theorem.

\begin{theo}\label{esth14}\index{Error estimate!-- for a hyperbolic 
equation!-- in the general case}
Under Assumption \refp{estdonnees}, let $\mesh$ be an admissible mesh in the 
sense of Definition \refp{meshypmd} and $k>0$. Let $g \in C(\R^2,\R)$ satisfy 
Assumption \ref{estcondflu} and assume that 
condition \eqref{estcfl} holds. Let $u$ be the unique entropy weak 
solution   of
\eqref{estpbcont}-\eqref{estcondini} and $\uapp$ be given by \eqref{estuapp},
\eqref{estschema},  \eqref{estschemaz}. Assume $u_0 \in BV(\R^\dim)$. Then, 
for all $R>0$ and all $T>0$ there exists $C_e \in \R_+$, only depending on 
$R$, $T$, $\vbf$, $g$, $u_0$, $\a$ and
$\xi$,  such that the following  inequality holds:
\be
\int_0^T \int_{B(0,R)} \vert \uapp(x,t) - u(x,t) \vert dx dt \le C_e h^{\frac 1 4}.
\label{esteqcorol}\ee
(Recall that $B(0,R)=\{x \in \R^\dim$, $\vert x \vert < R\}$.)
\end{theo}

In Theorem \ref{esth14}, $u_0$ is assumed to belong to
$BV(\R^\dim)$ (recall that $u_0 \in BV(\R^\dim)$ if $\sup\{ \int u_0(x)\div
\varphi(x)dx$,
$\varphi \in C^\infty_c(\R^\dim,\R^\dim)$;
$\vert \varphi(x) \vert \le 1$, $\forall x \in \R^\dim\}$ $<$ $\infty$).
This assumption allows us to obtain an $h^{1/4}$ estimate in 
\eqref{esteqcorol}. If $u_0 \not \in BV(\R^\dim)$ (but $u_0$ still belongs 
to $L^\infty(\R^\dim)$), one can also give an error estimate which
depends on the functions $\eps(r,S)$ and $\eps_0(r,S)$ defined
in \eqref{estepsrk} and \eqref{estepsork}.

\smallskip

A slight improvement of Theorem \ref{esth14} (and also Theorem 
\ref{esth14i} below) is possible. Using the fact that $u \in 
C(\Rp,L^1_{loc}(\R^\dim))$ and thus $u(\cdot,t)$ is defined for all $t 
\in \Rp$, Theorem \ref{esth14} remains true with

\[
\int_{B(0,R)} \vert \uapp(x,t)  -  u(x,t) \vert dx \le C_e h^{1/4}, \, 
\forall t\in [0,T],
\]

instead of \eqref{esteqcorol}. The proof of such a result 
may be handled with an adaptation of the proof a uniqueness of the entropy 
process solution given for instance in \citet{EGH2}, see \citet{Vi} and 
\citet{Coquel} for some similar results.

\smallskip

In some cases, it is possible to obtain $h^{1/2}$, instead of
$h^{1/4}$, in Theorem \ref{esth14}. This is the case, for instance,
when the mesh $\mesh$ is composed of rectangles ($\dim=2$) and when $\vbf$
does not depend on $(x,t)$, since, in this case, one obtains a ``$BV$
estimate" on $\uapp$. In this case, the right hand sides of inequalities \eqref{estbvfx} and 
\eqref{estbvft}, proven above, are changed from $C/\sqrt h$ to $C$, so that the right
hand side of \eqref{estmapp} becomes $C h$ instead of $C \sqrt h$,
which in turn yields $C_e h^{1/2}$ in \eqref{esteqcorol} instead 
of $C_e h^{1/4}$.
It is, however, still an open problem to know whether it is possible to obtain an error
estimate with $h^{1/2}$, instead of
$h^{1/4}$, in Theorem \ref{esth14} (under the hypotheses of Theorem
\ref{esth14}), even in the case where $\vbf$ does not depend on $(x,t)$ (see
\citet{C} for an attempt in this direction).

\begin{remark} \normalfont\label{estext} 
Theorem \ref{esth14} (and also Theorem \ref{esth14i}) remains true with some slightly more general
assumption on $g$, instead of \ref{estcondflu}, in order to allow $g$
to depend on $\mesh$ and $k$. Indeed, in \eqref{estschema}, one can
replace $g(\ucvn,\ucvvn)$ (and $g(\ucvvn,\ucvn)$) by
$g_{\cv,\cvv}(\ucvn,\ucvvn,\mesh,k)$ (and
$g_{\cvv,\cv}(\ucvvn,\ucvn,\mesh,k)$). Assume that, for all
$\cv \in \mesh$ and all $\cvv \in \Nv$, the
function $(a,b) \mapsto g_{\cv,\cvv}(a,b,\mesh,k)$, from $[U_m,U_M]^2$ to 
$\R$, is nondecreasing with respect to $a$, nonincreasing with respect to 
$b$, Lipschitz continuous uniformly with respect to $\cv$
and $\cvv$ and that $g_{\cv,\cvv}(a,a,\mesh,k)=f(a)$ for all $a \in 
[U_m,U_M]$ (recall that $U_m \le u_0 \le U_M$ a.e. on $\R^\dim$). Then Theorem
\ref{esth14} remains true. 

\smallskip

However, note that condition \eqref{estcfl} and 
$C_e$ in the estimate
\eqref{esteqcorol} of Theorem \ref{esth14} depend on the Lipschitz
constants of $g_{\cv,\cvv}(\cdot,\cdot,\mesh,k)$ on $[U_m,U_M]^2$. 
An interesting form for $g_{\cv,\cvv}$ is
$g_{\cv,\cvv}(a,b,\mesh,k)$ $=$ $c_{\cv,\cvv}(\mesh,k)f(a)$ $+$ $
(1-c_{\cv,\cvv}(\mesh,k))$ $f(b)$ $+$
$D_{\cv,\cvv}(\mesh,k)$ $(a-b)$, with some $c_{\cv,\cvv}(\mesh,k) \in
[0,1]$ and $D_{\cv,\cvv}(\mesh,k)$
$\ge 0$. In order to obtain the desired properties on
$g_{\cv,\cvv}$, it is sufficient to take $\max\{\vert f'(s) \vert,$ $s \in
[U_m,U_M]\} \le D_{\cv,\cvv}(\mesh,k) \le D$ (for all $\cv,\cvv$), with some
$D \in \R$. The Lipschitz constants of $g_{\cv,\cvv}$ on $[U_m,U_M]^2$ 
only depend on $D$, $f$, $U_m$ and $U_M$. 

\smallskip

For instance, a ``Lax-Friedrichs  type" 
scheme\index{Scheme!Lax-Friedrichs --}\index{Lax-Friedrichs scheme} consists, roughly speaking,
in taking $D_{\cv,\cvv}(\mesh,k)$ of order ``$h/k$". The desired properties on $g_{\cv,\cvv}$ are satisfied,
provided that $k/h \le C$, with some $C$ depending on $\max\{\vert f'(s) \vert,$ $s \in
[U_m,U_M]\}$. Note, however, that the condition $k/h \le C$ is not sufficient to give a real ``$h^{1/4}$"
estimate, since the coefficient $C_e$ in \eqref{esteqcorol} depends
on $D$. Taking, for example, $k$ of order ``$h^2$" leads to an
estimate ``$C_eh^{1/4}$" which do not goes to $0$ as $h$ goes to
$0$ (indeed, it is known, in this case, that the approximate solution
does not converge towards the entropy weak solution to 
\eqref{estpbcont}-\eqref{estcondini}). One obtains a real
``$h^{1/4}$" estimate, in the case of that ``Lax-Friedrichs 
type" scheme, by taking
$C_1 \le (k/h) \le C_2$.
In order to avoid the condition $C_1 \le (k/h)$ (note that $(k/h) \le 
C_2$ is imposed by the Courant-Friedrichs-Levy condition \ref{estcfl}), a possibility is to 
take $D_{\cv,\cvv}(\mesh,k)=D=\max\{\vert f'(s) \vert,$ $s \in
[U_m,U_M]\}$ (this is related to the ``modified Lax-Friedrichs " of  
Example \refp{flumono} in the 1D case). Then $D$ only depends on $f$ and 
$u_0$ and, in the estimate ``$C_eh^{1/4}$" of Theorem \ref{esth14}, $C_e$ 
only depends on $R$, $T$, $\vbf$, $f$, $u_0$, $\a$ and
$\xi$, which leads to a convergence result at rate ``$h^{1/4}$" as $h \tends 0$ 
(with fixed $\a$ and $\xi$).
\end{remark}

In the case of the implicit scheme, one proves
the following theorem.

\begin{theo}\label{esth14i}\index{Error estimate!-- for a hyperbolic 
equation!-- in the general case}
Under Assumption \refp{estdonnees}, let $\mesh$ be an admissible mesh in the 
sense of Definition \refp{meshypmd} and $k>0$. Let $g \in C(\R^2,\R)$ satisfy 
Assumption \ref{estcondflu}. Let $u$ be
the unique entropy weak solution of
\eqref{estpbcont}-\eqref{estcondini}. Assume that $u_0 \in BV(\R^\dim)$ and
\index{BV!Initial condition in --} that $\vbf$ does not depend on $t$. 

Let $\{\ucvn$, $n \in \N$,  $\cv\in 
\mesh\}$ be the unique solution to \eqref{estschemai} and \eqref{estschemaz}
such that $\ucvn \in [U_m, U_M]$ 
for all $\cv\in \mesh$ and $n \in \N$ (existence and uniqueness of such a 
solution is given by Proposition \ref{esteu}).  Let $\uapp$ be 
defined by \eqref{estuappimp}. 

Then, for all $R>0$ and $T>0$, there exists 
$C_e$, only depending on $R$, $T$, $\vbf$, $g$, $u_0$ and $\a$, such that 
the following inequality holds:
\be
\int_0^T \int_{B(0,R)} \vert \uapp(x,t) - u(x,t) \vert dx dt \le C_e(k+ h^{\frac 1
2})^{\frac 1 2}.
\label{esteqcoroli}\ee

\end{theo}

\begin{remark} \normalfont\label{estexti} 
Note that, in Theorem \ref{esth14i}, there is no restriction on $k$ (this
is usual for an implicit scheme), and one obtains an ``$h^{1/4}$" error 
estimate for some ``large" $k$, namely if $k \le h^{1/2}$. In Theorem 
\ref{esth14i}, if $\vbf$ depends on $t$ and $u_0 \in L^\infty(\R^\dim)$ 
(but $u_0$ not necessarily in $BV(\R^\dim)$), one can also give
an error estimate. Indeed one obtains

\[
\int_0^T \int_{B(0,R)} \vert \uapp(x,t) - u(x,t) \vert dx dt \le C_e(\frac k {h^{\frac
1 2}}+ h^{\frac 1 2})^{\frac 1 2},
\]

which yields an ``$h^{1/4}$" error estimate if $k$ is of order ``$h$".
\end{remark}

Theorem \ref{esth14} (resp. Theorem \ref{esth14i}) is an easy consequence of 
Theorem \ref{estentest} (resp. \ref{estentesti}) and of a quite general 
theorem of comparison between the entropy weak solution to 
\eqref{estpbcont}-\eqref{estcondini} and an approximate solution. This 
theorem of comparison (Theorem \ref{esterrest}) may be used in other frameworks (for instance, to 
compare the entropy weak solution to 
\eqref{estpbcont}-\eqref{estcondini} and the approximate solution obtained 
with a parabolic regularization of \eqref{estpbcont}). It is stated and proved 
in Section \ref{compthm} where the proofs of  theorems \ref{esth14} and \ref{esth14i} are also given.
First, in Section \ref{sprelim}, two preliminary lemmata are given. 
Indeed, Lemma \ref{estlemma} is the crucial part of the two following 
sections.

\subsection{Preliminary lemmata}\label{sprelim}
\index{BV!-- Space|(}

Let us first give a classical lemma on the space $BV$.

\begin{lemma}\label{bvtrans}
Let $u \in BV_{loc}(\R^p)$, $p \in \N^\star$, that is $u \in 
L^1_{loc}(\R^p)$ and 
the restriction of $u$ to $\O$ belongs to $BV(\O)$ for all open bounded 
subset $\O$ of $\R^p$ (see Definition \refp{defBVRp} for 
the definition of $BV(\O)$).

Then, for all bounded subset $\O$ of $\R^p$ and for all $a>0$,

\be
\norm{u(\cdot+\eta)-u}{L^1(\O)} \le \vert \eta \vert \snorm{u}{BV(\O_a)}, \, \forall 
\eta \in \R^p, \, \vert \eta \vert \le a,
\label{bvdtrans}
\ee

where $\O_a=\{x \in \R^p$; $d(x,\O) < a\}$ and $d(x,\O)=\inf\{\vert x-y 
\vert$, $y \in \O\}$ is the distance from $x$ to  $\O$.

\end{lemma}

\proof of Lemma \ref{bvtrans}

\smallskip

Let $\O$ be a bounded subset of $\R^p$ and $\eta  \in \R^p$.
The following equality classically holds:

\be
\norm{u(\cdot+\eta)-u}{L^1(\O)}=\sup\{\int_\O (u(x+\eta)-u(x))\vphi(x)dx, 
\, \vphi \in C^\infty_c(\O,\R), \, \norm{\vphi}{L^\infty(\O)} \le 1\}.
\label{defl1}
\ee

Let $\vphi \in C^\infty_c(\O,\R)$ such that $\norm{\vphi}{L^\infty(\O)} 
\le 1$.

Since $\vphi(x)=0$ if $x \in 
\O_{\vert \eta \vert} \setminus \O$ (recall that $\O_{\vert \eta \vert}=\{x \in \R^p$; $d(x,\O) < 
\eta\}$),

\[
\int_\O u(x)\vphi(x)dx=\int_{\O_{\vert \eta \vert}}u(x)\vphi(x)dx.
\]

Similarly, using an obvious change of variables,

\[
\int_\O u(x+\eta)\vphi(x)dx=\int_{\O_{\vert \eta \vert}}u(x)\vphi(x-\eta)dx.
\]

Therefore,

\[
\int_\O (u(x+\eta)-u(x))\vphi(x)dx=\int_{\O_{\vert \eta \vert}} 
u(x)(\vphi(x-\eta)-\vphi(x))dx=-\int_{\O_{\vert \eta \vert}} u(x) (\int_0^1 \grad \vphi 
(x-s\eta) 
\cdot \eta ds) dx
\]

and, with Fubini's theorem,

\be
\int_\O (u(x+\eta)-u(x))\vphi(x)dx = \int_0^1 (\int_{\O_{\vert \eta \vert}} u(x) \grad 
\vphi (x-s\eta) 
\cdot \eta dx) ds.
\label{vtt}
\ee

For all $s \in (0,1)$, Define $\psi_s\in C^\infty_c(\O_{\vert \eta 
\vert},\R^p)$ by
$\psi_s(x)=\vphi(x-s\eta)\eta$; since 
$\psi_s \in C^\infty_c(\O_{\vert \eta \vert},\R^p)$ and $\vert \psi_s(x) \vert \le 
\vert \eta \vert$ for all $x \in \R^p$, the definition of 
$\snorm{u}{BV(\O_{\vert \eta \vert})}$ yields

\[
\int_{\O_{\vert \eta \vert}} u(x) \grad \vphi (x-s\eta) \cdot \eta dx=\int_{\O_{\vert \eta \vert}} 
u(x) \div \psi_s (x) dx \le \vert \eta \vert \snorm{u}{BV(\O_{\vert \eta \vert})}.
\]

Then, \eqref{vtt} gives

\be
\int_\O (u(x+\eta)-u(x))\vphi(x)dx \le \vert \eta \vert 
\snorm{u}{BV(\O_{\vert \eta \vert})}.
\label{vttt}
\ee

Taking in \eqref{vttt} the supremum over $\vphi \in C^\infty_c(\O,\R)$ such that $\norm{\vphi}{L^\infty(\O)} 
\le 1$ yields, thanks to \eqref{defl1},

\[
\norm{u(\cdot+\eta)-u}{L^1(\O)} \le \vert \eta \vert \snorm{u}{BV(\O_{\vert \eta \vert})}, \, \forall 
\eta \in \R^p,
\]

and \eqref{bvdtrans} follows, since $\O_{\vert \eta \vert} \subset \O_a$ 
if $\vert \eta \vert \le a$. \eof

\begin{remark} \normalfont \label{finvacances}
Let us give an application of the lemma \ref{bvtrans} which will be quite 
useful further on.
Let $u \in BV_{loc}(\R^p)$, $p \in \N^\star$. Let $\psi$, $\vphi \in 
C_c(\R^p,\R_+)$, $a>0$ and 
$0< \eps <a$ such that $\int_{\R^p} \vphi (x) dx =1$ and $\vphi(x)=0$ for 
all $x \in \R^p$, 
$\vert x \vert > \eps$. Let $S=\{x \in \R^p$, $\psi(x) \ne 0\}$.

Then,

\be
\int_{\R^p} \int_{\R^p} \vert u(x) - u(y) \vert \psi(x) \vphi (x-y) dydx 
\le \eps \norm{\psi}{L^\infty(\R^p)} \snorm{u}{BV(S_a)},
\label{bvdvar}
\ee

where $S_a=\{x \in \R_p$, $d(x,S) < a\}$.

Indeed, Lemma \ref{bvtrans} gives

\be
\norm{u(\cdot+\eta)-u}{L^1(S)} \le \vert \eta \vert \snorm{u}{BV(S_a)}, \, \forall 
\eta \in \R^p, \vert \eta \vert \le a.
\label{ousuisje}
\ee

Using a change of variables in the left hand side of \eqref{bvdvar},

\[
\int_{\R^p} \int_{\R^p} \vert u(x) - u(y) \vert \psi (x) \vphi (x-y) 
dydx \le  \norm{\psi}{L^\infty(\R^p)} \int_{B(0,\eps)} (\int_{S} \vert 
u(x)-u(x-z) \vert dx) \vphi (z) dz.
\]

Then, \eqref{ousuisje} yields

\[
\int_{\R^p} \int_{\R^p} \vert u(x) - u(y) \vert \psi(x) \vphi (x-y) 
dydx \le \eps \norm{\psi}{L^\infty(\R^p)} \snorm{u}{BV(S_a)} \int_{\R^p}\vphi (z)dz,
\]

which gives \eqref{bvdvar}.

\end{remark}
\index{BV!-- Space|)}
\begin{lemma}\label{estlemma}

Under assumption \ref{estdonnees}, let $u_0 \in 
BV(\R^\dim)$\index{BV!Initial condition in --} and $\util \in
L^\infty(\R^\dim \times \Rpast)$ such that
$U_m \le \util \le U_M$ $a.e.$ on $\R^\dim \times \Rpast$. Assume that there exist
$\mu \in \mathcal{M}(\RNRp)$ and $\mu_0 \in \mathcal{M}(\R^\dim)$ such that 

\be\ba
\dsp{\int_{\Rp}\int_{\R^\dim}} \Bigl( &\vert \util(x,t) - \sk \vert \varphi_t(x,t) +\\
&(f(\util(x,t)\top\sk) - f(\util(x,t)\bot\sk)) 
\vbf(x,t)\cdot\grad\varphi(x,t) \Bigr) dx dt & + \\
\dsp{\int_{\R^\dim}} & \vert u_0(x) - \sk \vert \varphi(x,0) dx & \ge \\ 
- \dsp{\int_\RNRp} & \Bigl( \vert\varphi_t(x,t)\vert + 
\vert\grad\varphi(x,t)\vert\Bigr) d\mu(x,t) 
- \dsp{\int_{\R^\dim}}  \vert\varphi(x,0)\vert d\mu_0(x), \\ \\ &
 \ \ \forall \sk\in\R,
\ \ \forall\varphi\in C^\infty_c(\RNRp,\R_+). 
\ea\label{est1est}\ee

Let $u$ be the unique entropy weak solution of 
\eqref{estpbcont}-\eqref{estcondini} (i.e. $u \in L^\infty(\R^\dim \times 
\Rpast)$ is the unique  solution to \eqref{est1est} with
$u$ instead of $\util$ and $\mu = 0$, $\mu_0 = 0$).

Then  for all $\psi\in C^\infty_c(\RNRp,\R_+)$  there exists $C$ 
only depending on $\psi$ (more precisely on $\Vert \psi
\Vert_{\infty}$, $\Vert \psi_t \Vert_{\infty}$, $\Vert \grad\psi \Vert_{\infty}$, and on the
support of $\psi$), $\vbf$,
$f$, and $u_0$, such that

\be
\left\{\begin{array}{ll} 
\dsp{\int_{\Rp}\int_{\R^\dim}} 
&\Bigl[ \vert \util(x,t) - u(x,t) \vert  \psi_t(x,t) \ + \\
&\Bigl(f(\util(x,t)\top u(x,t)) -  f(\util(x,t)\bot u(x,t))\Bigr) \ 
(\vbf(x,t)\cdot\grad\psi(x,t)) \Bigr]  \ dx dt \ \ge \\
& - C (\mu_0(\{\psi(\cdot,0) \not = 0\})+(\mu(\{\psi \not = 0\}))^{\frac 1
2}+\mu(\{\psi \not = 0\})), 
\end{array}\right.
\label{est3est}\ee

where $\{\psi \not = 0\}=\{(x,t) \in \R^\dim \times \Rp$, $\psi(x,t) \not 
= 0\}$ and $\{\psi(\cdot,0) \not = 0\}=\{x \in \R^\dim$, $\psi(x,0) \not = 
0\}$. (Note that $\norm{\cdot}{\infty}=\norm{\cdot}{L^\infty(\R^\dim 
\times \Rpast)}$.)
\end{lemma}

\proof of Lemma \ref{estlemma}

\smallskip

The proof of Lemma \ref{estlemma} is close to that of step 1 in the proof of 
Theorem \ref{ewsunicite}.
Let us first define mollifiers in $\R$ and $\R^\dim$. For $p=1$ and 
$p=\dim$, one defines
$\rho_p \in C_c^\infty(\R^p,\R)$ satisfying the following properties: 

\[
\text{supp}(\rho_p) =  \{x \in \R^p; \, \rho_p(x) \neq 0\} 
\subset \{x \in \R^p;
\vert x \vert \le 1\},
\]

\[
\rho_p(x) \ge 0, \, \forall x \in \R^p,
\]

\[
\dsp{\int_{\R^p}} \rho_p (x) dx = 1
\]

and furthermore, for $p=1$,

\be
\rho_1(x) = 0, \, \forall x \in \R_+.
\label{estrhodec}
\ee

For $r\in\R$, $r \ge 1$, one defines $\rho_{p,r}(x) = r^p \rho_p (r  x )$, for all $x \in \R^p$.

Using the mollifiers $\rho_{p,r}$ will allow to choose convenient test 
functions in \eqref{est1est} (which are the inequalities satisfied by 
$\util$) and in the analogous inequalities satisfied by $u$ which are

\be\ba 
\dsp{\int_{\Rp}\int_{\R^\dim}} \Bigl[ \vert u(y,s) - \sk \vert \varphi_s(y,s) +
\Bigl( f(u(y,s)\top\sk) - f(u(y,s)\bot\sk) \Bigr) 
\vbf(y,s)\cdot\grad\varphi(y,s) \Bigr] dy ds + \\
\dsp{\int_{\R^\dim}} \vert u_0(y) - \sk \vert \varphi(y,0) dy \ge 0,
\ \ \forall \sk\in\R, \ \ \forall\varphi\in C^\infty_c(\RNRp,\R_+).
\ea\label{estkruzkovys}\ee

Indeed, the main tool is to take $\sk = u(y,s)$ in 
\eqref{est1est}, $\sk = \util(x,t)$ in \eqref{estkruzkovys}
and to introduce mollifiers in order to
have  $y$ close to $x$ and $s$ close to $t$. 

\smallskip

Let $\psi\in C^\infty_c(\R^\dim\times\R_+,\R_+)$, and let $\varphi:
(\R^\dim\times\R_+)^2 \tends \R_+$ be defined by:

\[
\varphi(x,t,y,s) = \psi(x,t) \r_{\dim,r}(x-y) \r_{1,r}(t-s).
\]

Note that, for any $(y,s) \in \R^\dim \times \R_+$, one
has $\varphi(\cdot,\cdot,y,s)\in C^\infty_c(\RNRp,\R_+)$ and, for any
$(x,t)
\in \R^\dim \times \R_+$, one has $\varphi(x,t,\cdot,\cdot)\in
C^\infty_c(\RNRp,\R_+)$. 
Let us take $\varphi(\cdot,\cdot,y,s)$ as test function
$\varphi$ in \eqref{est1est} and $\varphi(x,t,\cdot,\cdot)$ as test function
$\varphi$ in \eqref{estkruzkovys}. We take, in \eqref{est1est}, $\sk = u(y,s)$
and we take, in \eqref{estkruzkovys}, $\sk = \util(x,t)$. We then integrate
\eqref{est1est} for $(y,s)\in\RNRp$, and \eqref{estkruzkovys} for
$(x,t)\in\RNRp$. Adding the two  inequalities yields

\be
E_{11} + E_{12} + E_{13} + E_{14}  \ge - E_2,
\label{estcin}
\ee

where

\[ 
E_{11} = \dsp{\int_0^\infty \int_{\R^\dim} \int_0^\infty \int_{\R^\dim}} 
\Bigl[ \vert \util(x,t) - u(y,s) \vert  
\psi_t(x,t) \r_{\dim,r}(x-y) \r_{1,r}(t-s)  \Bigr] dx dt dy ds,
\]

\[\bari
E_{12} = \dsp{\int_0^\infty \int_{\R^\dim} \int_0^\infty \int_{\R^\dim}} \Bigl[
\Bigl(f(\util(x,t)\top u(y,s)) -  f(\util(x,t)\bot u(y,s))\Bigr)
\\ 
\vbf(x,t)\cdot\grad\psi(x,t) \r_{\dim,r}(x-y) \r_{1,r}(t-s)  
\Bigr] dx dt dy ds, 
\ean\]

\[\bari
E_{13} = -\dsp{\int_0^\infty \int_{\R^\dim} \int_0^\infty \int_{\R^\dim}}
\Bigl(f(\util(x,t)\top u(y,s)) -  f(\util(x,t)\bot u(y,s))\Bigr) 
\psi(x,t)
\\
(\vbf(y,s) - \vbf(x,t))\cdot \grad\r_{\dim,r}(x-y) \r_{1,r}(t-s) dx dt dy ds,
\ean\]

\[
E_{14} =\dsp{\int_{\R^\dim} \int_0^\infty \int_{\R^\dim}} \vert u_0(x) - u(y,s) \vert
\psi(x,0) \r_{\dim,r}(x-y) \r_{1,r}(-s) dydsdx
\]

and

\be\bari 
E_2 = \dsp{\int_0^\infty \int_{\R^\dim} \int_\RNRp} \Bigl(\vert 
\r_{\dim,r}(x-y) (\psi_t(x,t) \r_{1,r}(t-s) + 
\psi(x,t) \r'_{1,r}(t-s)) \vert 
\\
+\vert \r_{1,r}(t-s) (\grad \psi(x,t)  \r_{\dim,r}(x-y) + 
\psi(x,t) \grad\r_{\dim,r}(x-y)) \vert \Bigr) d\mu(x,t) dy ds 
\\
+ \dsp{\int_0^\infty \int_{\R^\dim} \int_{\R^\dim}} \vert \psi(x,0) \r_{\dim,r}(x-y) \r_{1,r}(-s) \vert 
d\mu_0(x) dy ds.
\ean\label{est5est}\ee

One may be surprised by the fact that the inequation \eqref{estcin} is 
obtained without using the initial condition which is satisfied by the 
 entropy weak solution  $u$  of \eqref{estpbcont}-\eqref{estcondini}. 
 Indeed, this initial condition appears
only in the third term of the left hand side of
\eqref{estkruzkovys}; since  $\varphi(x,t,\cdot,0)=0$ for all
$(x,t)
\in \R^\dim \times \R_+$, the third
term of the left hand side of
\eqref{estkruzkovys} is zero when   $\varphi(x,t,\cdot,\cdot)$ is chosen as 
a  test
function in \eqref{estkruzkovys}. However, the fact that $u$ satisfies the initial
condition of \eqref{estpbcont}-\eqref{estcondini} will be used later 
in order to get a bound on $E_{14}$.

\medskip

Let us now study  the five terms of
\eqref{estcin}. One sets $S=\{\psi \ne 0\}=\{(x,t)\in \R^\dim\times \R_+;$ $\psi(x,t) \ne
0\}$ and $S_0=\{\psi(\cdot,0) \ne 0\}=\{x \in \R^\dim;$ $\psi(x,0) \ne 
0\}$. In the following, the notation $C_i$ ($i \in \N$) will refer to various 
real quantities only depending on $\Vert
\psi \Vert_{\infty}$, $\Vert \psi_t
\Vert_{\infty}$, $\Vert \grad\psi \Vert_{\infty}$, $S$, $S_0$, $\vbf$, $f$,
and $u_0$. 

Equality \eqref{est5est} leads to

\be
E_2 \le (r+1)C_1\mu(S)+C_2\mu_0(S_0).
\label{est6est}\ee

Let us handle the term $E_{11}$. 
For all $x \in \R^\dim$ and for all $t \in \R_+$, one has, using
\eqref{estrhodec},

\[
\int_{\R^\dim} \int_0^\infty \rho_{\dim,r}(x-y)\rho_{1,r}(t-s)dsdy=1.
\]

Then,

\[
\ban
\vert E_{11} - \dsp{\int_{\Rp}\int_{\R^\dim}} 
\Bigl[ \vert \util(x,t) - u(x,t) \vert  
\psi_t(x,t) \Bigr] \ dx dt \vert \le
\\
\dsp{\int_0^\infty \int_{\R^\dim} \int_0^\infty \int_{\R^\dim}} 
\Bigl[ \vert u(x,t) - u(y,s) \vert  
\vert\psi_t(x,t)\vert \r_{\dim,r}(x-y) \r_{1,r}(t-s)  \Bigr] dx dt dy ds \le 
\Vert \psi_t \Vert_\infty \eps(r,S),
\ean
\]

with 

\be
\eps(r,S)=\sup\{ \norm{u-u(\cdot+\eta,\cdot+\tau)}{L^1(S)}, \; \vert \eta \vert \le \frac 1
r, \, 0 \le \tau \le \frac 1 r\}.
\label{estepsrk}
\ee

Since $u_0 \in BV(\R^\dim)$, the function $u$ (entropy weak solution to 
\eqref{estpbcont}-\eqref{estcondini}) belongs to $BV(\R^\dim \times 
(-T,T))$, for all $T>0$, setting, for instance, $u(.,t)=u_0$ for $t<0$ (see 
\citet{krushkov} or \citet{Chainais2} where this result is proven passing 
to the limit on numerical schemes).

Then, Lemma \ref{bvtrans} gives, since $r \ge 1$, (taking $p=\dim+1$, $\O=S$ and $a=\sqrt 2$ in 
Lemma \ref{bvtrans},)

\be
\eps(r,S) \le \frac {C_3} r.
\label{estubv}
\ee

Hence,

\be
\vert E_{11} - \dsp{\int_{\Rp}\int_{\R^\dim}} 
\Bigl[ \vert \util(x,t) - u(x,t) \vert  
\psi_t(x,t) \Bigr] \ dx dt \vert \le \frac {C_4} r.
\label{est7estb}
\ee

In the same way, using $\vert f(a \top b)-f(a \top c) \vert \le M \vert  
b-c \vert$ and $\vert f(a \bot b)-f(a \bot c) \vert \le M \vert  
b-c \vert$ for all $a$, $b$, $c \in [U_m,U_M]$ where $M$ is the Lipschitz 
constant of $f$ in $[U_m,U_M]$,  

\be\bari
\vert E_{12} - \dsp{\int_{\Rp}\int_{\R^\dim}}
\Bigl(f(\util(x,t)\top u(x,t)) -  f(\util(x,t)\bot u(x,t))\Bigr)
\\ 
(\vbf(x,t)\cdot\grad\psi(x,t)) dx dt\vert \le C_5 \eps(r,S) \le
\frac {C_6} r.
\ean\label{est8est}\ee

Let us now turn to  $E_{13}$. We compare this term with 

\[\bari
E_{13b} = -\dsp{\int_0^\infty \int_{\R^\dim} \int_0^\infty \int_{\R^\dim}}
\Bigl( f(\util(x,t)\top u(x,t)) -  f(\util(x,t)\bot u(x,t))\Bigr) 
\psi(x,t)
\\
(\vbf(y,s) - \vbf(x,t))\cdot
 \grad\r_{\dim,r}(x-y) \r_{1,r}(t-s) \ dx dt dy ds.
\ean
\]

Since $\div(\vbf(\cdot,s)-\vbf(x,t))=0$ (on $\R^\dim$) for all $x \in \R^\dim$,
$t \in \R_+$ and $s \in \R_+$, one has $E_{13b}=0$. Therefore,
substracting $E_{13b}$ from
$E_{13}$ yields

\be\bari
 E_{13} \le C_7 \dsp{\int_0^\infty \int_{\R^\dim} \int_0^\infty \int_{\R^\dim}}
\vert u(x,t) - u(y,s) \vert 
\psi(x,t)
\\
\vert(\vbf(y,s) - \vbf(x,t))\cdot
 \grad\r_{\dim,r}(x-y)\vert \r_{1,r}(t-s) \ dx dt dy ds.
\ean
\label{est8estb}
\ee

The right hand side of \eqref{est8estb} is then smaller than
$C_8 \eps(r,S)$, since $\vert(\vbf(y,s) -
\vbf(x,t))\cdot\grad\r_{\dim,r}(x-y)\vert$ is bounded by $C_9r^\dim$ 
(noting that $\vert x-y \vert \le 1/r$). Then, with
\eqref{estubv}, one has

\be
E_{13} \le \frac {C_{10}} r.
\label{est8estc}
\ee

In order to estimate $E_{14}$, let us take in \eqref{estkruzkovys}, for $x \in
\R^\dim$ fixed, $\vphi = \vphi(x,\cdot,\cdot)$, with

\[
\varphi(x,y,s) =  \psi(x,0) \r_{\dim,r}(x-y) 
\int_s^\infty\r_{1,r}(-\tau)d\tau,
\] 

and $\sk = u_0(x)$. Note that $\vphi(x,\cdot,\cdot) \in 
C^\infty_c(\R^\dim \times \Rp,\R_+)$. We then
integrate the resulting inequality with respect to $x\in\R^\dim$. We get

\[
-E_{14} + E_{15} + E_{16} \ge 0,
\]

with

\[
\bari
E_{15} = -
\dsp{\int_0^\infty \int_{\R^\dim} \int_{\R^\dim} \int_s^\infty}
(f(u(y,s)\top u_0(x)) - f(u(y,s)\bot u_0(x)))
\\
\vbf(y,s)\cdot (\psi(x,0) \grad\r_{\dim,r}(x-y)) \r_{1,r}(-\tau)  d\tau 
dy dx ds,
\ean
\]

\[
\ban
E_{16} = 
 \dsp{\int_{\R^\dim} \int_{\R^\dim} \int_0^\infty} 
\psi(x,0) \r_{\dim,r}(x-y) \r_{1,r}(-\tau)
\vert u_0(x) - u_0(y) \vert d\tau dy dx. 
\ean
\]

To bound $E_{15}$, one introduces $E_{15b}$ defined as

\[\bari
E_{15b} = 
\dsp{\int_0^\infty \int_{\R^\dim} \int_{\R^\dim} \int_s^\infty} 
(f(u(y,s)\top u_0(y)) - f(u(y,s)\bot u_0(y)))
\\
(\vbf(y,s)\cdot\grad\r_{\dim,r}(x-y)) \psi(x,0)  \r_{1,r}(-\tau) 
d\tau dy dx ds.
\ean\]

Integrating by parts for the $x$ variable yields

\[\bari
E_{15b} = -
\dsp{\int_0^\infty \int_{\R^\dim} \int_{\R^\dim} \int_s^\infty} 
(f(u(y,s)\top u_0(y)) - f(u(y,s)\bot u_0(y)))
\\
(\vbf(y,s)\cdot\grad\psi(x,0)) \r_{\dim,r}(x-y)  \r_{1,r}(-\tau) 
d\tau dy dx ds.
\ean\]

Then, noting that the time
support of this integration is reduced to $s \in [0,1/r]$, one has 

\be
E_{15b} \le \frac {C_{11}} r.
\label{est9estd}
\ee

Furthermore, one has

\[
\vert E_{15}+E_{15b} \vert \le C_{12}
\int_0^\infty \int_{\R^\dim} \int_{\R^\dim} \int_s^\infty  \vert u_0(x) -
u_0(y)\vert  \vert \vbf(y,s)\cdot\grad\r_{\dim,r}(x-y) \vert \psi(x,0) 
\r_{1,r}(-\tau) d\tau dy dx ds,
\]

which is bounded by $C_{13}\eps_0(r,S_0)$, since the
time support of the integration is reduced to
$s \in [0,1/r]$, where $\eps_0(r,S_0)$ is defined by

\be
\eps_0(r,S_0)=\sup\{\int_{S_0}\vert
u_0(x)-u_0(x+\eta)\vert dx; \; \vert\eta\vert \le \frac 1 r\}.
\label{estepsork}
\ee

Since $u_0 \in BV(\R^\dim)$, one has (thanks to Lemma \ref{bvtrans}) 
$\eps_0(r,S_0) \le {C_{14}}/r$
and therefore, with \eqref{est9estd}, $E_{15} \le {C_{15}}/ r$.

\smallskip

Since   $u_0 \in
BV(\R^\dim)$, again thanks to Lemma \ref{bvtrans}, see remark 
\ref{finvacances},
the term $E_{16}$ is also bounded by $C_{16}/r$.

\smallskip

Hence, since $E_{14} \le E_{15} + E_{16}$,

\be
E_{14} \le \frac {C_{17}} r.
\label{est9estf}
\ee

Using \eqref{estcin}, \eqref{est6est}, \eqref{est7estb},
\eqref{est8est},\eqref{est8estc}, \eqref{est9estf}, one obtains 

\[
\begin{array}{ll} 
\dsp{\int_{\Rp}\int_{\R^\dim}} 
&\Bigl[ \vert \util(x,t) - u(x,t) \vert  \psi_t(x,t) \ + \\
&\Bigl(f(\util(x,t)\top u(x,t)) -  f(\util(x,t)\bot u(x,t))\Bigr) \ 
(\vbf(x,t)\cdot\grad\psi(x,t)) \Bigr]  \ dx dt \ \ge \\
& - C_1(r+1)\mu(S) -C_2\mu_0(S_0) - \frac {C_{18}} {r}, 
\end{array}
\]

which, taking $r=1/{\sqrt {\mu(S)}}$ if $0 < \mu(S) \le 1$ ($r \tends \infty$ if
$\mu(S)=0$ and $r=1$ if $\mu(S) > 1$), gives
\eqref{est3est}.

This concludes the proof
of the lemma \ref{estlemma}.\eof

\subsection{Proof  of the error estimates}\label{compthm}

Let us now prove a quite general theorem of comparison between the entropy weak solution to 
\eqref{estpbcont}-\eqref{estcondini} and an approximate solution, from 
which  theorems \ref{esth14} and \ref{esth14i} will be  deduced.

\begin{theo}\label{esterrest}

Under assumption \ref{estdonnees}, let $u_0 \in BV(\R^\dim)$ and  $\util \in
L^\infty(\R^\dim \times \Rpast)$ such that $U_m \le \util \le U_M$ a.e. 
on $\R^\dim \times \Rpast$. Assume that there exist $\mu \in \mathcal{M}(\RNRp)$ and $\mu_0 \in \mathcal{M}(\R^\dim)$ such that 
\eqref{est1est} holds. Let $u$ be the unique entropy weak solution of 
\eqref{estpbcont}-\eqref{estcondini} (note that $u \in
L^\infty(\R^\dim \times \Rpast)$ is solution to \eqref{est1est} with
$u$ instead of $\util$ and $\mu = 0$, $\mu_0 = 0$).

Then, for all $R>0$ and all $T>0$ there exists $C_e$ and $\bar{R}$, only depending on $R$, $T$, $\vbf$, 
$f$  and $u_0$, such that the
following  inequality holds:

\[\bari
\int_0^T \int_{B(0,R)} \vert \util(x,t) - u(x,t) \vert dx dt \le C_e ( \mu_0(B(0,\bar{R})) +
[\mu(B(0,\bar{R}) \times [0,T])]^{\frac 1 2}
\\
+\mu(B(0,\bar{R}) \times [0,T]) ).
\ean
\]

Recall that $B(0,R)=\{x \in \R^\dim$; $\vert x \vert < R\}$.
\end{theo}

\proof of Theorem \ref{esterrest}

\smallskip

The proof of Theorem \ref{esterrest} is close to that of Step 2 in the 
proof of Theorem \ref{ewsunicite}. It uses Lemma 
\refp{estlemma}, the proof of which is given in section \ref{sprelim} above.

Let $R>0$ and $T>0$.
One sets $\o = V M$, where $V$ is given in Assumption \ref{estdonnees} and $M$ is
the Lipschitz constant of $f$ in $[U_m,U_M]$ (indeed, since $f \in
C^1(\R,\R)$, one has $M= \sup \{\vert f'(s) \vert;$ $s \in [U_m,U_M]\}$).

Let $\r \in C_c^1(\R_+,[0,1])$  be a function such that
$\r(r) = 1$ if  $r \in [0,R+\o T]$, $\r(r) = 0$ if
$r \in [R+\o T +1,\infty)$ and $\r'(r) \le 0$, for all $r \in \R_+$ ($\r$ 
only depends on $R$, $T$, $\vbf$, $f$ and $u_0$).

One takes, in \eqref{est3est}, $\psi$ defined by 

\[
\ba \psi(x,t) = \r(\vert x
\vert + \o t) \frac {T-t} T, \text{ for } x\in \R^\dim  \text{ 
 and }t\in[0,T], \\
 \psi(x,t)=0,  \text{ for } x \in \R^\dim  \text{  and }t\ge T.\ea
 \]

Note that $\r(\vert x \vert + \o t)=1$, if $(x,t) \in B(0,R) \times [0,T]$.

The function $\psi$ is not in
$C^\infty_c(\RNRp,\R_+)$, but, using a usual regularization  technique, 
it may be proved that such a function can be considered in
\eqref{est3est}, in which case  
 Inequality \eqref{est3est} reads, with $\bar{R}=R+\o T+1$,

\[
\ban
\dsp{\int_0^T \int_{\R^\dim}} 
\Bigl[ \vert \util(x,t) - u(x,t) \vert \Bigl(\frac {T-t} T
 \o \r'(\vert x \vert + \o t) 
- \frac 1 T \r(\vert x \vert + \o t) \Bigr) \ + \\
\Bigl(f(\util(x,t)\top u(x,t)) -  f(\util(x,t)\bot u(x,t))\Bigr) \ 
\frac {T-t} T \r'(\vert x \vert + \o t)
(\vbf(x,t)\cdot\frac x {\vert x\vert}) \Bigr] \ dx dt \ \ge \\
 - C (\mu_0(B(0,\bar{R})) +
(\mu(B(0,\bar{R}) \times [0,T]))^{\frac 1 2}+\mu(B(0,\bar{R}) \times [0,T])), 
\ean
\]

where $C$ only depends on $R$, $T$, $\vbf$, $f$  and $u_0$.

Since $\o=VM$ and $\r'\le 0$, one has

\[\begin{array}{l}
\dsp {\Bigl(f(\util(x,t)\top u(x,t)) -  f(\util(x,t)\bot u(x,t))\Bigr) \ 
\frac {T-t} T \r'(\vert x \vert + \o t)
(\vbf(x,t)\cdot\frac x {\vert x\vert}) \Bigr) \le} \\
\vert \util(x,t) - u(x,t) \vert \frac {T-t} T
 \o (-\r'(\vert x \vert + \o t)),
\end{array}
\]

and therefore, since $\r(\vert x \vert + \o t)=1$, if $(x,t) \in B(0,R) 
\times [0,T]$,

\[
\dsp{\int_0^T \int_{B(0,R)}}
\vert \util(x,t) - u(x,t) \vert  dx dt \le 
C T (\mu_0(B(0,\bar{R})) +
(\mu(B(0,\bar{R}) \times [0,T]))^{\frac 1 2}+\mu(B(0,\bar{R}) \times [0,T])).
\]

This completes the proof of Theorem \ref{esterrest}.\eof

Let us now conclude with the proofs of theorems \refp{esth14} (which gives an error
estimate for the time explicit numerical scheme
\eqref{estschema}, \refep{estschemaz}) and 
\refp{esth14i} (which gives an error estimate for the time implicit numerical scheme
\eqref{estschemai}, \refep{estschemaz}). There are easy consequences of 
theorems \ref{estentest} and \ref{estentesti} and of Theorem \ref{esterrest}.

\medskip

\proof of Theorem \ref{esth14}

\smallskip

Under the assumptions of Theorem \ref{esth14}, let $\util=\uapp$. Thanks 
to the $L^\infty$ estimate on $\uapp$ (Lemma \ref{estinftyest}) and to 
Theorem \ref{estentest}, $\util=\uapp$ satisfies the hypotheses of Theorem 
\ref{esterrest} with $\mu=\mu_{\mesh,k}$ and $\mu_0=\mu_{\mesh}$ (the 
measures $\mu_{\mesh,k}$ and $\mu_{\mesh}$ are given in Theorem 
\ref{estentest}).

Let $R>0$ and $T>0$. Then, Theorem \ref{esterrest} gives the existence of 
$C_1$ and $\bar{R}$, only depending on $R$, $T$, $\vbf$, 
$f$ and $u_0$, such that 

\be\bari
\int_0^T \int_{B(0,R)} \vert \uapp(x,t) - u(x,t) \vert dx dt \le C_1 ( 
\mu_{\mesh}(B(0,\bar{R})) +
[\mu_{\mesh,k}(B(0,\bar{R}) \times [0,T])]^{\frac 1 2}
\\
+\mu_{\mesh,k}(B(0,\bar{R}) \times [0,T]) ).
\ean
\label{ouvaisje}\ee

For $h$ small enough, say $h \le R_0$, one has $h<\bar{R}$ and 
$k<T$ (thanks to condition \ref{estcfl}, note that $R_0$ only depends on  $R$, $T$, $\vbf$, 
$g$, $u_0$, $\a$ and $\xi$).

Then, for $h<R_0$, Theorem \ref{estentest} gives, with \eqref{ouvaisje},

\[
\int_0^T \int_{B(0,R)} \vert \uapp(x,t) - u(x,t) \vert dx dt \le 
C_1(Dh + \sqrt C h^{\frac 1 4} + C \sqrt h) \le C_2 h^{\frac 1 4},
\]

where $C_2$ only depends on $R$, $T$, $\vbf$, 
$g$, $u_0$, $\a$ and $\xi$.

This gives the desired estimate $\eqref{esteqcorol}$ of Theorem 
\ref{esth14} for $h < R_0$.

There remains the case $h \ge R_0$. This case is trivial since, for $h \ge 
R_0$,

\[
\int_0^T \int_{B(0,R)} \vert \uapp(x,t) - u(x,t) \vert dx dt \le 
2\max\{-U_m,U_M\} \meas (B(0,R) \times (0,T)) \le C_3 (R_0)^{\frac 1 4} \le 
C_3 h^{\frac 1 4},
\]

for some $C_3$ only depending on $R$, $T$, $\vbf$, 
$g$, $u_0$, $\a$ and $\xi$.

This completes the proof of Theorem \ref{esth14}.
\eof

\proof of Theorem \ref{esth14i}

\smallskip

The proof of Theorem \ref{esth14i} is very similar to that of Theorem 
\ref{esth14} and we follow the proof of Theorem 
\ref{esth14}.

Under the assumptions of Theorem \ref{esth14i}, using Theorem 
\ref{estentesti} instead of Theorem \ref{estentest} gives that 
$\util=\uapp$ satisfies the hypotheses of Theorem 
\ref{esterrest} with $\mu=\mu_{\mesh,k}$ and $\mu_0=\mu_{\mesh}$ (the 
measures $\mu_{\mesh,k}$ and $\mu_{\mesh}$ are given in Theorem 
\ref{estentesti}).

Let $R>0$ and $T>0$. Theorem \ref{esterrest} gives the existence of 
$C_1$ and $\bar{R}$, only depending on $R$, $T$, $\vbf$, 
$f$  and $u_0$, such that \eqref{ouvaisje} holds.

For $h<\bar{R}$ and $k < T$  Theorem \ref{estentest} gives with 
\eqref{ouvaisje},

\[
\int_0^T \int_{B(0,R)} \vert \uapp(x,t) - u(x,t) \vert dx dt \le 
C_1(Dh +  \sqrt C (k+h^{\frac 1 2})^{\frac 1 2} + C(k+h^{\frac 1 2})) \le 
C_2 (k+h^{\frac 1 2})^{\frac 1 2},
\]

where $C_2$ only depends on $R$, $T$, $\vbf$, 
$g$, $u_0$, $\a$.

This gives the desired estimate $\eqref{esteqcoroli}$ of Theorem 
\ref{esth14i} for $h < \bar{R}$ and $k <T$.

There remains the cases $h \ge \bar{R}$ and $k \ge T$. These cases are trivial since 

\[
\int_0^T \int_{B(0,R)} \vert \uapp(x,t) - u(x,t) \vert dx dt \le 
2\max\{-U_m,U_M\} \meas (B(0,R) \times (0,T)) \le C_3 \inf\{\bar{R}^{\frac 
1 4},T^{\frac 1 2}\}
\]

for some $C_3$ only depending on $R$, $T$, $\vbf$, 
$g$, $u_0$.

This completes the proof of Theorem \ref{esth14i}.
\eof

\subsection{Remarks and open problems}

Theorem \refp{esth14} gives an error estimate of order $h^{1/4}$ for 
the approximate solution of a nonlinear hyperbolic equation of the form 
$ u_t + \div (\vbf f(u)) = 0$\index{Equation!Hyperbolic --}, with initial data in $L^\infty \cap BV$ by 
the explicit finite volume scheme \eqref{estschema} and \refep{estschemaz}, under a 
usual CFL condition $k \le C h$ (see \refep{estcfl}).

\medskip

Note that, in fact, the same estimate holds if $u_0$ is only locally 
$BV$.
More generally, if the initial data $u_0$ is only in $L^\infty$, 
then one still obtains an error 
estimate in terms of the quantities

\[
\eps(r,S)=\sup\{\int_{S}\vert
u(x,t)-u(x+\eta,t+\tau)\vert dxdt; \; \vert \eta \vert \le \frac 1
r, 
\, 0 \le \tau \le \frac 1 r\}
\]

and

\[
\eps_0(r,S_0)=\sup\{\int_{S_0}\vert
u_0(x)-u_0(x+\eta)\vert dx; \; \vert\eta\vert \le \frac 1 r\},
\]

see \refep{estepsrk} and \refep{estepsork}. This is again an obvious consequence of 
Theorem \refp{estentest} and Theorem \refp{esterrest}.

\medskip

We also considered the implicit schemes, which seem to be 
much more widely used in industrial codes in order to ensure their 
robustness.  The implicit case required additional work in order 

(i) to prove the existence of the solution to the finite volume scheme,

(ii) to obtain the ``strong time $BV$" estimate \eqref{estbvt} if 
$\vbf$ does not depend on $t$.

For $\vbf$ depending on $t$, Remark \ref{estexti} yields an estimate of order 
$h^{1/4}$ if $k$ behaves as $h$; however, in the case where $\vbf$ 
does not depend on $t$, then an estimate of order 
$h^{1/4}$ is obtained (in Theorem \ref{esth14i}) for a behaviour of $k$ as $\sqrt h$; Indeed, recent 
numerical experiments suggest that taking $k$ of the order of  $\sqrt h$ yields 
results of the same precision than taking $k$ of the order of $h$, 
with an obvious reduction of the computational cost. 

\medskip

Note that the method described here may also be extended to higher order 
schemes for the same equation, see \citet{Chainais1}; other methods have been 
used for error estimates for higher order schemes with a nonlinearity of 
the form $F(u)$, as in  \citet{Noelle}. However, it is still an 
open problem, to our knowledge, to improve the order of the error estimate 
in the case of higher order schemes.

\section{Boundary conditions}\label{hyperbcmd}\index{Boundary conditions for
hyperbolic equations}

In this section, a generalization of Theorem \ref{thm} is presented
for the multidimensional scalar case together with a rough sketch of proof. For the sake of
simplicity, one considers $d=2$ (the extension to $d=3$ is straightforward) and a flux function
under the form $v(x,t) f(u)$, with $\mathrm{div}(v(\cdot,t))=0$ (see \cite{v} for
the general case of a flux function $f(x,t,u)$). This leads to the following equation:

\begin{equation}
u_t + \mathrm{div}(vf(u)) =0, \text{ in } \Omega \times (0,T),
\label{mdeq}
\end{equation}
where $\Omega$ is a bounded polygonal open set of $\mathbf{R}^2$, $T >0$,
$f \in C^1(\mathbf{R},\mathbf{R}$)  (or $f~:~\mathbf{R} \rightarrow \mathbf{R}$ Lipschitz continuous) and
$v \in C^1(\mathbf{R}^2
\times [0,T]) \rightarrow \mathbf{R}^2$ with $\mathrm{div}(v(\cdot,t))=0$ in $\mathbf{R}^2$ for all
$t \in [0,T]$. The unknown is $u~: ~\Omega \times (0,T) \rightarrow \mathbf{R}$.

\vspace{.1cm}

Let $u_0  \in L^\infty(\Omega)$ and  $\overline u \in L^\infty(\partial \Omega
\times (0,T))$. Let $A, B \in \mathbf{R}$ be such that $A \le u_0 \le B$ a.e. on $\Omega$ and
$A \le \overline u \le B$ a.e. on $\partial \Omega \times (0,T)$.
Following the work of \cite{o}, an entropy weak solution
of (\ref{mdeq}) with the initial condition $u_0$ and the (weak) boundary condition $\overline u$
is a solution of (\ref{ewsmd}):

\begin{equation}
\begin{array}{l}
u \in L^\infty(\Omega \times (0,T)),
\\
\displaystyle \int_0^T \int_{\Omega} [(u-\kappa)^{\pm}\varphi_t + 
\text{sign}_{\pm}(u-\kappa)(f(u)-f(\kappa)) v \cdot \mathrm{grad} \varphi ] dx dt
\\
\hfill \displaystyle + M \int_0^T  \int_{\partial \Omega} (\overline u (t) - 
 \kappa)^{\pm} \varphi(x,t) d\gamma(x) dt  \\
\hfill  \displaystyle+ \int_\Omega (u_0-\kappa)^{\pm}\varphi(x,0)dx \ge 0,
 \\
\hfill  \forall \kappa \in [A,B], \; \forall \varphi \in C^1_c(\overline \Omega \times 
 [0,T), \mathbf{R}_+),
\end{array}
\label{ewsmd}
\end{equation}
where $d \gamma(x)$ stands for the integration with respect to the one dimensional
Lebesgue measure on the boundary of $\Omega$ and $M$ is such that
$\Vert v \Vert_{\infty} \vert  f(s_1)-f(s_2) \vert \le M \vert s_1 - s_2 \vert$ for all $s_1,s_2 \in [A,B]$,
where $\Vert v \Vert_{\infty}=\sup_{(x,t) \in \Omega \times [0,T]} \vert v(x,t) \vert$ (and $\vert \cdot \vert$ denotes here the Euclidean norm in $\mathbf{R}^2$).

\begin{remark}~

\begin{enumerate}
\item If $u$ satisfies the family of inequalities (\ref{ewsmd}), it is possible to prove that $u$ is a solution
of (\ref{mdeq}) (on a weak form), $u$ satisfies some entropy inequalities in $\Omega \times (0,T)$,
namely $\vert u- \kappa \vert_t + \mathrm{div}(v (f(\max(u,\kappa))-f(\min(u,\kappa)))) \le 0$ for all $\kappa \in \mathbf{R}$, but also on the boundary of $\partial \Omega$ and on $t=0$.
$u$ satisfies the initial condition ($u(\cdot,0)=u_0$) and $u$ satisfies partially the boundary
condition. For instance, if $f^\prime >0$ and $u$ is regular enough, then $u(x,t)=\overline u(x,t)$ if $x \in \partial \Omega$,
$t \in (0,T)$ and $v(x,t) \cdot n(x,t) <0$, where $n$ is the outward normal vector to $\partial \Omega$.
\item Let $\overline M \ge 1$. It is interesting to remark that $u$ is solution of (\ref{ewsmd}) if and only if
$u$ is solution of (\ref{ewsmd}) where the term $\int_\Omega (u_0-\kappa)^{\pm}\varphi(x,0)dx$
is replaced by $\overline M \int_\Omega (u_0-\kappa)^{\pm}\varphi(x,0)dx$.
\end{enumerate}
\end{remark}

A sketch of proof of existence and uniqueness of the solution of (\ref{ewsmd}) together
with the convergence of numerical approximations is now given, following \cite{v}.

\vspace{.1cm}

\textsc{Step 1: Approximate solution.} With a quite general mesh of $\Omega$ (with triangles, for instance), denoted by $\mathcal{T}$,
and a time step $k$, it
is possible to define an approximate solution, denoted by $u_{\mathcal{T},k}$, using some numerical fluxes
(on the edges of the mesh) satisfying conditions similar to (C1)-(C3) in Sect. \ref{conv}.
Under a so called CFL condition (like $k \le (1-\zeta) \frac h L$ in Sect. \ref{conv}), it is
easy to prove that $A \le u_{\mathcal{T},k} \le B$ a.e. on $\Omega \times (0,T)$.
Unfortunately, it does not seem easy to obtain directly a strong compactness result on the
familly of approximate solutions (alhough this strong compactness result is true, as we
shall see below). 

\vspace{.1cm}

\textsc{Step 2: Weak compactness.} Using only this $L^\infty$ bound on $u_{\mathcal{T},k}$, one can
assume (for convenient subsequences of sequences of approximate solutions) that
$u_{\mathcal{T},k} \rightarrow u$, as the mesh size goes to zero (with the CFL condition), in a ``non linear weak-$\star$ sense"
(similar to the convergence towards young measures, see \cite{egh} for instance), that is
$u \in L^\infty(\Omega \times (0,T) \times (0,1))$ and
\[
\begin{array}{c}
\displaystyle \int_0^T \int_\Omega g(u_{\mathcal{T},k}(x,t)) \varphi (x,t) dx dt \rightarrow \int_0^1 \int_0^T \int_\Omega
g(u(x,t,\alpha))  \varphi (x,t) dx dt d\alpha,
\\
\hfill \text{for all } \varphi \in L^1(\Omega \times (0,T)).
\end{array}
\]

\vspace{.1cm}

\textsc{Step 3: Passing to the limit.} Using the monotonicity of the numerical fluxes, the approximate
solutions satisfy some discrete entropy inequalities. Passing to the limit in these inequalities
gives that $u$ (defined in Step 2) is solution of some inequalities which are very similar to (\ref{ewsmd}), namely:

\begin{equation}
\begin{array}{l}
u \in L^\infty(\Omega \times (0,T) \times (0,1)),
\\
\displaystyle \int_0^1 \int_0^T \int_{\Omega} [(u-\kappa)^{\pm}\varphi_t + 
\text{sign}_{\pm}(u-\kappa)(f(u)-f(\kappa)) v \cdot \mathrm{grad} \varphi ] dx dt d\alpha
\\
\hfill \displaystyle + M \int_0^T  \int_{\partial \Omega} (\overline u (t) - 
 \kappa)^{\pm} \varphi(x,t) d\gamma(x) dt  \\
\hfill  \displaystyle+ \int_\Omega (u_0-\kappa)^{\pm}\varphi(x,0)dx \ge 0,
 \\
\hfill  \forall \kappa \in [A,B], \; \forall \varphi \in C^1_c(\overline \Omega \times 
 [0,T), \mathbf{R}_+),
\end{array}
\label{epsmd}
\end{equation}

For this step, one chooses $M$ not only greater than the 
Lipschitz constant of $\Vert v \Vert_{\infty} f$ on
$[A,B]$, but also greater than the Lipschitz constant (on $[A,B]^2$) of
the numerical fluxes associated to the edges of the meshes (the equivalent
of $L$ in Theorem \ref{thm}). This choice of $M$ is possible
since the unique solution of (\ref{ewsmd}) does not
depend on $M$ provided that $M$ is greater
than the Lipschitz constant of $\Vert v \Vert_{\infty} f$ on
$[A,B]$ and since it is possible to choose numerical fluxes
(namely, Godunov flux, for instance)
such as the Lipschitz constant of
these numerical fluxes is bounded by the Lipschitz constant of $\Vert v \Vert_{\infty} f$
(then, the present method leads to an existence result with
$M$ only greater than the 
Lipschitz constant of $\Vert v \Vert_{\infty} f$ on $s \in
[A,B]$, passing to the limit on approximate solutions given with these
numerical fluxes).

\vspace{.1cm}

\textsc{Step 4: Uniqueness of the solution of (\ref{epsmd}).} In this step, the ``doubling variables"
method of Krushkov is used to prove the uniqueness of the solution of (\ref{epsmd}). Indeed,
if $u$ and $w$ are two solutions of (\ref{epsmd}), the doubling variables method leads to:

\begin{equation}
\begin{array}{l}
\displaystyle \int_0^1 \int_0^1 \int_0^T \int_{\Omega} \vert u(x,t,\alpha)-w(x,t,\beta)\vert \varphi_t  \, dx dt
d\alpha d\beta 
\\
\hfill + \displaystyle \int_0^1 \int_0^1 \int_0^T \int_{\Omega}(f(\max(u,w))-
f(\min(u,w))) v \cdot \mathrm{grad} \varphi \,  dx dt d\alpha d\beta \ge 0
 \\
\hfill  \forall \varphi \in C^1_c(\overline \Omega \times 
 [0,T), \mathbf{R}_+),
\end{array}
\label{uniq}
\end{equation}

Taking $\varphi (x,t) =(T-t)^+$ in (\ref{uniq}) (which is, indeed, possible) gives
that $u$ does not depend on $\alpha$, $v$ does not depend on $\beta$ and
$u=v$ a.e. on $\Omega \times (0,T)$. As a result, $u$ is also the unique solution of (\ref{ewsmd}).

\vspace{.1cm}

\textsc{Step 5: Conclusion.} Step 4 gives, in particular, the uniqueness of
the solution of (\ref{ewsmd}). It gives also that the non linear weak-$\star$ limit of sequences
of approximate solutions is solution of (\ref{ewsmd}) and, therefore, the existence of the
solution of (\ref{ewsmd}).
Furthermore, since the non linear weak-$\star$
limit of sequences of approximate solution does not depend on $\alpha$, it is quite easy to deduce
that this limit is ``strong" in $L^p(\Omega \times (0,T))$ for any $p \in [1,\infty)$ (see \cite{egh}, for
instance) and, thanks to
the uniqueness of the limit, the convergence
holds without extraction of subsequences.

\section{Nonlinear weak-$\star$ convergence} \label{nlwscv}

\hpage{Hyperbolic multid}

The notion of nonlinear weak-$\star$ convergence was used
in Section \ref{eghbscewsol}. We give here the definition of 
this type of convergence and we prove that a bounded sequence of 
$L^\infty$ converges, up to a subsequence, in the  nonlinear weak-$\star$ sense.

\begin{defi}[Nonlinear weak-$\star$  convergence] \label{nlwsdef}  \normalfont~

\index{Convergence!Nonlinear weak star --}
Let $\O$ be an open subset of $\R^N$ ($N \ge 1$), 
$(u_n)_\nnn \subset L^\infty(\O)$ 
and $u \in L^\infty(\O \times (0,1))$. 
The sequence $(u_n)_\nnn$ converges towards $u$ in the 
``nonlinear weak-$\star$ sense" if

\be
\bari 
\dsp{\int_\O g(u_n(x)) \varphi(x) dx \tends \int_0^1 \int_\O g(u(x,\a)) 
\varphi(x) dx d\a,
\text{ as } n \tends + \infty,}
\\
\forall \varphi \in L^1(\O), \; \forall g \in C(\R,\R).
\ean
\label{ews5}
\ee

\end{defi}

\begin{remark} \normalfont \label{nlwsib} Let $\O$ be an open subset of $\R^N$ ($N \ge 1$), 
$(u_n)_\nnn \subset L^\infty(\O)$ 
and $u \in L^\infty(\O \times (0,1))$ such that $(u_n)_\nnn$ converges 
towards $u$ in the nonlinear weak-$\star$ sense. Then, in particular, the 
sequence $(u_n)_\nnn$ converges towards $v$ in $L^\infty(\O)$, for the 
weak-$\star$ topology, where $v$ is defined by

\[
v(x)=\int_0^1 u(x,\a) d\a, \text{ for a.e. } x \in \O.
\]

Therefore, the sequence $(u_n)_\nnn$ is bounded in $L^\infty(\O)$ (thanks 
to the Banach-Steinhaus theorem). The following proposition gives that, 
up to a subsequence, a bounded sequence of $L^\infty(\O)$ converges in 
the nonlinear weak-$\star$ sense.
\end{remark}

\begin{prop}\label{nlwsprop}
\index{Compactness!Nonlinear weak -- in $L^\infty$}
Let $\O$ be an open subset of 
$\R^N$ ($N \ge 1$) and $(u_n)_{n\in\N}$ 
be a bounded sequence of $L^\infty(\O)$.  Then there exists a 
subsequence of
$(u_n)_{n\in\N}$, which will still be denoted by  $(u_n)_{n\in\N}$, and a 
function  $u  \in L^\infty(\O\times(0,1))$ such that the subsequence 
$(u_n)_{n\in\N}$ converges towards $u$ in the nonlinear weak-$\star$ 
sense.

\end{prop}

\proof

\smallskip

This proposition is classical in the framework of ``Young 
measures" and we only sketch the proof for the sake of completeness.

\smallskip

Let $(u_n)_{n\in\N}$ be a bounded sequence of $L^\infty(\O)$ and $r \ge 0$ such
that $\Vert u_n \Vert_{L^\infty(\O)} \le r, \forall n \in \N$.

\smallskip

\textit{Step 1 (diagonal process)}

Thanks to the separability of
the set of continuous functions defined from $[-r,r]$ into $\R$ (this set 
is endowed with the uniform norm) and the  
sequential weak-$\star$ relative compactness of the bounded
sets of $L^\infty(\O)$ \index{Compactness!Weak -- in $L^\infty$}, 
there exists (using a diagonal process) a 
subsequence, which will still be denoted by  $(u_n)_{n\in\N}$, such that, 
for any function  $g \in C(\R,\R)$, the sequence $(g(u_n))_{n\in\N}$ converges in  $L^\infty(\O)$ for the
weak-$\star$ topology towards a function  $\mu_g \in L^\infty(\O)$.

\smallskip

\textit{Step 2 (Young measure)}

In this step, we prove the existence of a family $(m_x )_{x \in \O}$ such 
that

\begin{enumerate}
\item for all $x \in \O$, $m_x $ is a probability on $\R$ whose support 
is included in $[-r,+r]$ (i.e. $m_x$ is a $\s$-additive application
from the Borel $\s$-algebra of $\R$ in $\R_+$ such that $m_x(\R)=1$ and 
$m_x(\R \setminus [-r,r])=0$),

\item $\mu_g (x) = \int_{\R } g(s)d{m_x }(s)$ for a.e. $x \in
\O$ and for all $g \in C(\R,\R )$. 
\end{enumerate}

The family $m=(m_x )_{x \in \O}$ is called a ``Young measure".

Let us first claim that it is possible to define 
$\mu_g \in L^\infty(\O)$ for $g \in C([-r,r],\R)$ by setting $\mu_g=\mu_f$ 
where $f \in C(\R,\R)$ is such that $f=g$ on $[-r,r]$. Indeed, this 
definition is meaningful since if $f$ and $h$ are two elements 
of $C(\R,\R)$ such that $f=g$ on $[-r,r]$ then $\mu_f$ and $\mu_h$ are the same 
element of $L^\infty(\O)$ (i.e. $\mu_f=\mu_h$ a.e. on $\O$) thanks to the 
fact that $-r \le u_n \le r$ a.e. on $\O$ and for all 
$\nnn$.

For $x \in \O$, let

\[
E_x = \{g \in C([-r,r],\R); \,  \lim_{h \tends 0}
\frac 1 {\meas (B(0,h))} \int_{B(x,h)} \mu_g(z) dz \text{ exists in } \R\},
\]

where $B(x,h)$ is the ball of 
center $x$ and radius $h$ (note that $B(x,h) \subset \O$ for $h$ small 
enough).

If $g \in E_x$, we set

\[
\bar{\mu}_g(x) =  
\lim_{h \tends 0}
\frac 1 {\meas (B(0,h))} \int_{B(x,h)} \mu_g(z) dz.
\]

Then, we define $T_x$ from $E_x$ in $\R$ by $T_x(g) =
\bar{\mu}_g(x)$. It is easily seen that $E_x$ is a vector space which contains the
constant functions, that $T_x$ is a linear application from $E_x$ to $\R$ 
and that $T_x$ is nonnegative (i.e. $g(s)\ge 0$ for all $s \in \R$ 
implies $T_x(g) \ge 0$). Hence, using a modified
version of the Hahn-Banach theorem, one can prolonge
$T_x$ into a linear nonnegative application $\Tbar_x$ defined on the 
whole set $C([-r,r],\R)$. By a classical Riesz theorem, there exists a 
(nonnegative) measure $m_x$ on the Borel sets of $[-r,r]$ 
such that 

\be
\Tbar_x(g)=\int_{-r}^r g(s) d m_x(s), \forall g \in C([-r,r], \R).
\label{ewsTbar}
\ee

If $g(s)=1$ for all $s \in [-r,r]$, the function $g$ belongs to  $E_x$ and 
$\bar{\mu}_g(x) = 1$ (note that $\mu_g=1$ a.e. on $\O$).
Hence, from (\ref{ewsTbar}), $m_x$ is a probability over $[-r,r]$, and therefore a
probability over $\R$ by prolonging it by 0 outside of $[-r,r]$. This 
gives the first item on the family $(m_x)_{x \in \O}$.

\smallskip

Let us prove now the second item on the family $(m_x)_{x \in \O}$. If $g 
\in C([-r,r],\R)$ then $g \in E_x$ for a.e. $x \in \O$ and $\mu_g (x) = 
\bar{\mu}_g(x)$ for a.e. $x \in \O$ (this is a classical 
result, since $\mu_g \in L^1_{loc}(\O)$, see \citet{rudin}). 
Therefore, $\mu_g(x)=T_x(g)=\Tbar_x(g)$ for a.e. $x \in \O$. Hence,

\[
\mu_g(x)=\int_{-r}^r g(s)dm_x(s) \text{ for a.e. } x \in \O,
\]

for all $g \in C([-r,r],\R)$ and therefore for all $g\in C(\R,\R)$. 
Finally, since the support of $m_x$ is included in 
$[-r,r]$,

\[
\mu_g(x)=\int_\R g(s)dm_x(s) \text{ for a.e. } x \in \O, \, \forall g 
\in C(\R,\R).
\]

This completes Step 2.

\smallskip

\textit{Step 3 (construction of $u$)}

It is well known that, if $\bar{m}$ is a 
probability on $\R$, one has

\be
\int_\R g(s)d\bar{m}(s)=\int_0^1g(u(\a))d\a, \, \forall g \in \M_b,
\label{distfct}
\ee

where $\M_b$ is the set of bounded measurable functions from $\R$ to $\R$ and 
with

\[
u(\a)=\sup\{c \in \R; \, \bar{m}((-\infty,c)) < \a\}, \, \forall \a \in 
(0,1).
\]

Note that the function $u$ is measurable, nondecreasing and 
left continuous. Furthermore, if the support of $\bar{m}$ is included in $[a,b]$ 
(for some $a$, $b \in \R$, $a<b$) then $u(\a) \in [a,b]$ for all $\a \in 
(0,1)$ and \eqref{distfct} holds for all $g \in C(\R,\R)$.

Applying this result to the measures $m_x$ leads to the definition of 
$u$ as

\[
u(x,\a)=\sup\{c \in \R; \, m_x((-\infty,c)) < \a\}, \,  \forall \a \in 
(0,1), \, \forall x \in \O.
\]

For all $x \in \O$, the function $u(x,\cdot)$ is measurable (from $(0,1)$ 
to $\R$), nondecreasing, left continuous and takes its values in $[-r,r]$.
Furthermore,

\[
\mu_g(x)=\int_0^1 g(u(x,\a))d\a \text{ for a.e. } x \in \O, \, \forall g 
\in C(\R,\R).
\]

Therefore,

\[
\bari
\dsp{\int_\O g(u_n(x))\vphi(x)dx \tends \int_\O (\int_0^1 
g(u(x,\a))d\a)\vphi(x)dx,  \text{ as } n \tends \infty,}
\\
\forall \vphi \in L^1(\O), \, \forall g \in C(\R,\R).
\ean
\]

In order to conclude the proof of Proposition \ref{nlwsprop}, there 
remains to show that modifying $u$ on a negligible set leads to a 
function (still denoted by $u$) measurable with 
respect to $(x,\a) \in \O \times (0,1)$. Indeed, this mesurability is 
needed in order to assert for instance, applying Fubini's Theorem (see 
\citet{rudin}), that

\[
\int_\O (\int_0^1 g(u(x,\a))d\a)\vphi(x)dx = \int_0^1 (\int_\O 
g(u(x,\a))\vphi(x)dx)d\a,
\]

for all  $\vphi \in L^1(\O)$ and for all $g \in C(\R,\R)$.

\smallskip

For all $g \in C(\R,\R)$, one chooses for $\mu_g$ (which belongs to 
$L^\infty(\O)$) a bounded measurable function from $\O$ to $\R$.

Let us define $\E=\{g_{a,b}$; $a$, 
$b \in \Q$, $a < b\}$ where $g_{a,b} \in C(\R,\R)$ is defined by

\[
\bal
g_{a,b}(x)=1 \text{ if } x \le a,
\\
g_{a,b}(x)= \frac{x-b}{a-b} \text{ if } a < x < b,
\\
g_{a,b}(x)=0 \text{ if } x \ge b.
\ean
\]

Since $\E$ is a countable subset of $C(\R,\R)$, there exists a Borel 
subset $A$ of $\O$ such that $\meas (A)=0$ and 

\be
\mu_g (x) = \int_{\R } g(s)d{m_x}(s), \; \forall x \in
\O \setminus A, \; \forall g \in \E.
\label{partout}
\ee

Define for all $\a \in (0,1)$ $v(.,\a)$ by

\[
\bal
v(x,\a)=0 \text{ if } x \in A,
\\
v(x,\a)=\sup \{c \in \R, \, m_x((-\infty,c)) <\a \} \text{ if } x \in 
\O \setminus A,
\ean
\]

so that $u=v$ on $(\O \setminus A) \times (0,1)$ (and then $u=v$ a.e. on $\O 
\times (0,1)$).

Let us now prove that $v$ is measurable from $\O \times
(0,1)$ to $\R$ (this will conclude the proof of Proposition 
\ref{nlwsprop}).

Since $v(x,.)$ is left continuous on $(0,1)$ for all $x \in \O$, proving that $v(.,\a)$ is measurable (from $\O$ to 
$\R$) for all $\a \in (0,1)$ leads to the mesurability of $v$ on $\O \times
(0,1)$ (this is also classical, see \citet{rudin}).

There remains to show the mesurability of $v(.,\a)$ for all
$\a \in (0,1)$.

\smallskip

Let $\a \in (0,1)$ (in the following, $\a$ is fixed). Let us set 
$w=v(.,\a)$ and define, for $c \in \R$,

\[
f_c(x)=m_x((-\infty,c)) -\a, \, x \in \O \setminus A,
\]

so that $v(x,\a)=w(x)=\sup \{c \in \R,\; f_c(x)<0 \}$ for all $x \in \O 
\setminus A$. 

Using \eqref{partout} leads to

\[
m_x((-\infty,c))=\sup \{\mu_g(x), \, g\le 1_{(-\infty,c)} \text{~and } 
g\in \E \}, \, \forall x \in \O \setminus A.
\]

Then, the function $f_c: \O \setminus A \rightarrow \R$ is
measurable as the supremum of a countable set of measurable functions 
(recall that $\mu_g$ is measurable for all $g \in \E$).

In order to prove the measurability of $w$ (from $\O$
to $\R$), it is sufficient to prove that $\{x\in \O \setminus A$; $w(x) \ge a\}$ is a 
Borel set, for all $a\in \R$ (recall that $w=0$ on $A$).

Let $a\in \R$, since $f_c(x)$ is nondecreasing with respect to $c$, one 
has

\[
\{x\in \O \setminus A;\; w(x) \ge a\}=\cap_{n >0} \{x\in \O \setminus A; 
\, f_{a- \frac 1 n}(x)<0 \}.
\]

Then $\{x\in \O \setminus A;\; w(x) \ge a\}$ is measurable, thanks to the measurability 
of $f_c$ for all $c \in \R$.

This concludes the proof of Proposition \ref{nlwsprop}. \eof

\begin{remark} \normalfont\label{cvforte} Let $\O$ be an open subset of $\R^N$ ($N \ge 1$), 
$(u_n)_\nnn \subset L^\infty(\O)$ 
and $u \in L^\infty(\O \times (0,1))$ such that $(u_n)_\nnn$ converges 
towards $u$ in the nonlinear weak-$\star$ sense. Assume that $u$ does not 
depend on $\a$, i.e. there exists $v \in L^\infty(\O)$ such that 
$u(x,\a)=v(x)$ for a.e. $(x,\a) \in \O \times (0,1)$. Then, it is easy to 
prove that $(u_n)_\nnn$ converges towards $u$ in $L^p(B)$ for all $1 \le p < \infty$ and all 
bounded subset $B$ of $\O$. Indeed, let $B$ be a bounded subset of $\O$. 
Taking, in \eqref{ews5}, $g(s)=s^2$ (for all $s \in \R$) and $\vphi=1_B$ 
and also $g(s)=s$ (for all $s \in \R$) and $\vphi=1_Bv$ leads to

\[
\int_B (u_n(x)-v(x))^2 dx \tends 0, \text{ as } n \tends \infty.
\]

This proves that $(u_n)_\nnn$ converges towards $u$ in $L^2(B)$. The 
convergence of $(u_n)_\nnn$ towards $u$ in $L^p(B)$ for all $1 \le p < 
\infty$ is then an easy consequence of the $L^\infty(\O)$ bound on 
$(u_n)_\nnn$ (see Remark \ref{nlwsib}).
\end{remark}

\section{A stabilized finite element method}\label{fve} 
\index{Stabilization of a finite element method}
\hpage{Hyperbolic multid}

In this section, we shall try to compare the finite element method
\index{Finite element method|(} to the 
finite volume method for the discretization of a nonlinear hyperbolic 
equation. It is well known that the use of the finite element is not 
straightforward in the case of hyperbolic equations, since the lack of 
coerciveness of the operator yields a lack of stability of the finite 
element scheme. There are several techniques to stabilize these schemes, 
which are beyond the scope of this work. 
Here, as in \citet{Sel}, we are interested in viewing the finite element as a finite volume method, by writing it in a conservative form, and using a stabilization as in the third item of Example \refp{flumono}. 

\medskip

Let $F \in C^1(\R,\R^2)$, consider the following scalar conservation law:

\be
u_t(x,t) + \div(F(u))(x,t)  = 0, \, x \in \R^2, \, t \in \R_+, 
\label{claw}
\ee
with an initial condition. Let $\mesh$ be a triangular mesh of $\R^2$, well suited for the finite element method. Let $\mathcal{S}$ denote the set of nodes of this mesh, and let  $(\phi_j)_{j \in \mathcal{S}}$ be the classical piecewise bilinear shape functions. 
Following the finite element principles, let us look for an approximation of $u$ in the space spanned by the shape functions $\phi_j$; hence, at time $t_n = nk$ (where $k$ is the time step), we look for an approximate solution of the form 
\[
	u(.,t_n) = \sum_{j \in \mathcal{S}} u^n_j \phi_j;
\]
multiplying \eqref{claw} by $\phi_i$, integrating over $\R^2$, approximating $F(\sum_{j \in \mathcal{S}} u^n_j \phi_j)$ by $\sum_{j \in \mathcal{S}} F(u^n_j) \phi_j$, using the explicit Euler scheme\index{Scheme!Explicit Euler --} for the time discretization) and the mass lumping technique on the mass matrix as described in Remark \ref{rem:mass-lumping} \index{Mass lumping}  yields the following scheme:
\[
\frac {u_i^{n+1} - u_i^n} k \int_{\R^2} \phi_i(x) dx - \sum_{j \in \mathcal{S}} F(u_j^n) \cdot \int_{\R^2} \phi_j(x) \grad  \phi_i(x) dx = 0,
\]

which reads, noting that $\dsp{ \int \phi_j(x) \grad \phi_i (x) dx =
- \int \phi_i(x) \grad \phi_j (x) dx}$ and that $\dsp{ \sum_{j \in \mathcal{S}}  \grad \phi_j(x) = 0}$,
 
\[
\frac {u_i^{n+1} - u_i^n} k \int_{\R^2} \phi_i(x) dx  +
\sum_{j \in \mathcal{S}} (F(u_i^n) +  F(u_j^n)) \cdot \int_{\R^2} \phi_i(x) \grad 
\phi_j(x) dx = 0.\]

This last equality may also be written
\[
\frac {u_i^{n+1} - u_i^n} k \int_{\R^2} \phi_i(x) dx  +
\sum_{j \in \mathcal{S}} E_{i,j} = 0, \]
where 
\[
 E_{i,j} = \half (F(u_i^n) +  F(u_j^n) ) \cdot \int_{\R^2} (\phi_i(x) \grad 
\phi_j(x) - \phi_j(x) \grad  \phi_i(x)) dx.
\] 
Note that $E_{j,i} = -  E_{i,j}$.

\medskip

This is a centered and therefore unstable scheme. One way to stabilize 
it is to replace $E_{i,j}^n$ by 

\[
\tilde E_{i,j}^n = E_{i,j}^n + D_{i,j} (u_i^n - u_j^n),
\]
where $D_{i,j} =  D_{j,i}$ (in order for the scheme to remain 
``conservative") and 
$D_{i,j} \ge 0$ is chosen large enough so that $E_{i,j}^n$ is a 
nondecreasing function of $u_i^n$ and a nonincreasing function of 
$u_j^n$, which ensure the stability of the scheme, under a so 
called CFL condition, and does not change the ``consistency" (see \refep{fluxLF} 
and Remark \refp{estext})\index{Finite element method|)}.

\section{Moving meshes}\index{Mesh!Moving --|(}
For some evolution problems the use of time variable control volumes is advisable, e.g. 
when the domain of study changes with time. This is the case, for instance, for the simulation 
of a flow in a porous medium, when the porous medium is heterogeneous and its geometry 
changes with time. In this case, the mesh is required to move with the medium. 
The influence of the moving mesh on  the finite volume formulation can be 
explained  by 
considering the following simple transport equation:

\be
u_t (x,t) + \div (u \vbf)(x,t) =0, \ \  x \in \R^2, \ \ t
\in \R_+,
\label{cetv}
\ee

where $\vbf$ depends on the unknown $u$ (and possibly on other unknowns).
Let $k$ be the time step, and set $t_n=nk$, $n \in \N$. Let
$\mesh(t)$ be the mesh at time $t$. Since the mesh moves, the elements of the mesh 
vary in
time. For a fixed $n \in \N$, let $R(\cv,t)$ be the domain of $\R^2$ occupied by the element  
$\cv$ ($\cv \in \mesh(t_n)$) at time
$t$, $t \in [t_n,t_{n+1}]$, that is $R(\cv,t_n)=\cv$. Let $\vbf_s(x,t)$ be the velocity of the
displacement of the mesh  at point $x \in \R^2$ and for all $t \in [t_n,t_{n+1}]$ (note
that $\vbf_s(x,t) \in \R^2$).  Let
$u_\cv^n$ and $u_\cv^{n+1}$ be the discrete unknowns associated to element $\cv$ at 
times
$t_n$ and $t_{n+1}$ (they can be considered as the approximations of the mean values of
 $u(\cdot,t_n)$ and 
$u(\cdot,t_{n+1})$
over $R(\cv,t_n)$ and  $R(\cv,t_{n+1})$ respectively). The discretization of (\ref{cetv}) 
must take
into account the evolution of the mesh in time. In order to do so, let us first consider the
following  differential equation with initial condition: 

\be
\begin{array}{rcl}
\dsp{\frac {\partial y} {\partial t}} (x,t) & = & - \vbf_s(y(x,t),t), \ \ t \in
[t_n,t_{n+1}], \\ y(x,t_n) & = & x.
\end{array}
\label{edy}
\ee

Under suitable assumptions on $\vbf_s$ (assume for instance that $\vbf_s$ is 
continuous, Lipschitz continuous with respect to its first variable and 
that the 
Lipschitz
constant is integrable with respect to its second variable), the problem (\ref{edy}) has, 
for
all  $x \in \R^2$, a unique (global) solution.
 For $x \in \R^2$, define the function $y(x,\cdot)$ from $[t_n,t_{n+1}]$ to
$\R^2$ as the solution of problem (\ref{edy}).  Let $(\varphi_p)_{p \in\N} \subset
C_c^1(\R^2,\R_+)$ such that $0 \le  \varphi_p (x) \le 1$ for $x \in 
\R^2$ and for all $p \in \N$, and such that $
\varphi_p \tends 1_\cv$ a.e. as $p \tends + \infty$. Multiplying (\ref{cetv}) by $\psi_p(x,t) =
\varphi_p(y(x,t))$ and integrating over $\R^2$ yields
 
\be
\int_{\R^2} \Bigr( \frac {\partial (u \psi_p)} {\partial t}(x,t) + u(x,t) 
\grad \varphi_p
(y(x,t)) \cdot \vbf_s(y(x,t),t) -  (u \vbf)(x,t) \cdot \grad \psi_p (x,t) 
\Bigl) dx = 0. 
\label{upsit}
\ee
Using the explicit Euler discretization  
in time on Equation (\ref{upsit}) and denoting by
$u^n(x)$ a (regular) approximate value of $u(x,t_n)$ yields  

\[
\bac
{\dsp \int_{\R^2} \frac 1 k \Bigr( u^{n+1}(x)  \psi_p (x,t_{n+1}) -  u^n(x) 
\psi_p (x,t_{n}) \Bigl) dx   +  } 
\\   
{\dsp\int_{\R^2}  u^n(x)  (\vbf_s(x,t_n) -
\vbf(x,t_n)) \cdot \grad \varphi_p (x)  dx  =  0},  
\ean
\]

which also gives (noting that $\psi_p(x,t) = \varphi_p(y(x,t))$)

\be
\bac
{\dsp \int_{\R^2} \frac 1 k \Bigr( u^{n+1}(x)  \varphi_p(y(x,t_{n+1})) -  u^n(x) 
\varphi_p(y(x,t_{n})) \Bigl) dx -}
\\
{\dsp \int_{\R^2}  \div (u^n  (\vbf_s -
\vbf)) (x,t_n) \cdot \varphi_p (x)  dx  =  0}.
\ean
\label{upsinb}
\ee

Letting $p$ tend to infinity and noting that  $1_\cv (y(x,t_{n})) = 1_{R(\cv,t_{n})} (x) $ 
and  $1_\cv
(y(x,t_{n+1})) = 1_{R(\cv,t_{n+1})}(x)$, (\ref{upsinb}) becomes 

\[
\frac 1 k \Bigl(\int_{R(\cv,t_{n+1})} u^{n+1}(x) dx - \int_{R(\cv,t_{n})} u^{n}(x) 
dx\Bigr)
 +  \int_{R(\cv,t_{n})} \div((\vbf-\vbf_s) u^n) (x,t_n)dx=0,  
\]

which can also be written 

\[
\bac
{\dsp \frac 1 k (u_\cv^{n+1}\meas(R(\cv,t_{n+1}))-u_\cv^n \meas(R(\cv,t_{n})))
+ }
\\
{\dsp \int_{\partial R(\cv,t_{n})}(\vbf-\vbf_s)(x,t_n)\cdot \bfn_\cv (x,t_n) u^n(x) 
\dfrontiere(x)=0},  
\ean
\]

where $u_\cv^n = [1 / \meas(R(\cv,t_n))] \int_{R(\cv,t_{n})} u^n(x) dx$ and  
$u_\cv^{n+1} =  [1 / \meas(R(\cv,t_{n+1})) ]
\int_{R(\cv,t_{n+1})} u^{n+1}(x) dx$.  Recall that $\bfn_\cv$ denotes the normal to
$\partial  \cv$,
outward to $\cv$. The complete   discretization of the problem uses some additional 
equations (on $\vbf$,
$\vbf_s$\ldots).
 
\begin{remark} \normalfont
The above considerations concern a pure convection equation. 
In the case of a convection-diffusion 
equation\index{Equation!Convection-diffusion --},  such a moving
mesh may become non-admissible in the sense of definitions 
\refp{meshdirichlet} or \refp{meshneuman}. It is an interesting open 
problem to understand what should be done in that case. 
\end{remark}\index{Mesh!Moving --|)}

\chapter{Systems} \label{systemes}

In chapters \ref{ellund} to \ref{hypmd}, the finite volume was successively investigated for
the discretization of elliptic, parabolic, and hyperbolic equations.
In most scientific models, however, systems of equations have to be discretized. 
These may be partial differential equations of the same type or of different
types, and they may also be coupled to ordinary differential equations or algebraic
equations.

\smallskip

The discretization of systems of elliptic equations by the finite
volume method is straightforward, following the principles which were introduced in chapters
\ref{ellund} and \ref{ellmd}. Examples of the performance of the finite volume method  for
systems of elliptic equations on  rectangular  meshes, with ``unusual" source
terms (in particular, with source terms located on the edges or interfaces of the mesh) may be found in
e.g. \citet{Angot} (see also references therein), 
\citet{FH1} (where a comparison to a mixed finite element formulation 
\index{Finite element method} is also performed).
Parabolic systems are treated similarly as elliptic systems, with the addition of a
convenient time discretization.

\smallskip 

A huge literature is devoted to the discretization of
hyperbolic systems of equations, in particular to systems related to the
 compressible Euler equations, using structured or
unstructured meshes. We shall give only a short insight on this subject 
in Section \ref{systemh}, without any convergence result. Indeed, very 
few theoretical results of convergence of 
numerical schemes are known on this subject. We refer to 
\citet{godlewski-springer} and references therein for a more complete 
description of the numerical schemes for hyperbolic systems.

\smallskip 

Finite volume methods are also well adapted to the discretization of systems of 
equations of different types (for instance, an elliptic or parabolic equation
coupled with hyperbolic equations). Some examples are considered in 
sections \refp{nseq} and \refp{flowpm}. The classical case of 
incompressible Navier-Stokes\index{Navier-Stokes equations}   (for which,
generally, staggered grids\index{Staggered grid}  are used) and examples which arise in the simulation
of a multiphase flow in a porous medium are described. The latter example also serves as
an illustration of how to deal with  algebraic equations and inequalities.

\section{Hyperbolic systems of equations}  \label{systemh}

\hpage{Systems}

Let us consider a hyperbolic system consisting of $m$  equations (with $m \ge 1$). 
The unknown of the system is a function
$u=(u_1,\ldots,u_m)^t$, from $\Obar \times [0,T]$ to $\R^m$, where $\O$ is an open set of 
$\R^\dim$
(i.e. $\dim \ge 1$ is the space dimension), and $u$ is a solution of the following system:

\be
\bari
\dsp{\frac {\pa u_i}{\pa t} (x,t)
+ \sum_{j=1}^\dim \frac {\pa G_{i,j}}{\pa
x_j} (x,t)= g_i(x,t,u(x,t)),}
\\
 x=(x_1,\ldots,x_\dim) \in \O, \, t \in (0,T), \, i=1,\ldots,m,
\ean\label{rtgeq}
\ee

where
\[ G_{i,j} (x,t) = F_{i,j} (x,t,u(x,t)), \]
and the functions $F_j=(F_{1,j},\ldots,F_{m,j})^t$
($j=1,\ldots,\dim$) and $g=(g_1,\ldots,g_m)^t$ are given functions from $\Obar \times [0,T]
\times \R^m$ (indeed, generally, a part of $\R^m$, instead of $\R^m$) to  
$\R^m$. The function  $F=(F_1,\ldots,F_\dim)$ is assumed to  satisfy the usual hyperbolicity
condition, that is, for any (unit) vector of $\R^\dim$, $\n$, the
derivative of $F \cdot \n$ with respect to its third
argument (which can be considered as an $m \times m$ matrix) has only real eigenvalues 
and is diagonalizable.

Note that in real applications, diffusion terms may also be present in the
equations, we shall omit them here.
In order to complete System \eqref{rtgeq}, an initial condition for $t=0$ and
adequate boundary conditions for $x \in \pa \O$ must be specified.

\medskip

In the first section (Section \ref{godroe}), we shall only
briefly describe the general method of discretization by finite volume and some
classical schemes. In the subsequent sections, some possible treatments of difficulties
appearing in real simulations will be given.

\subsection{Classical schemes}\label{godroe}

Let us first describe some classical finite volume schemes for the discretization of
\eqref{rtgeq} with initial and boundary conditions,  using the concepts and notations which
were introduced in chapter \ref{hypmd}.  Let $\mesh$ be an admissible mesh in the sense of
Definition \refp{meshypmd} and $k$ be the time step, which is
assumed to be constant (the generalization to a variable time step is easy). 
We recall that the interface,
$\edgecvcvv$, between any two elements 
$\cv$ and $\cvv$ of $\T$ is assumed to be included in a hyperplane of $\Rn$. 
The discrete 
unknowns are the $\ucvn$, $\cv \in \T$, $n \in \{0,\ldots,\Nk+1\}$, 
with $N_k \in \N$, $(\Nk+1)k=T$. 
For $\cv\in\T$, let $\Nv$ be the set of its neighbours, that is the set
 of elements $\cvv$  of $\T$ such that the $(\dim-1)$ Lebesgue measure
of $\edgecvcvv$ is positive.  For $\cvv\in\Nv$, let $\ncvcvv$ be the unit normal
vector to
$\edgecvcvv$ oriented from $\cv$ to $\cvv$. Let $\tn=nk$, for $n \in
\{0,\ldots,\Nk+1\}$.

A finite volume scheme reads

\be
\bari
\dsp {\mcv \frac{\ucvnp-\ucvn}{k}+\sumv \medgecvcvv \Fcvcvvn} 
=\mcv \gcv, \\
\cv \in \T, \, n \in \{0,\ldots,\Nk\},
\ean
\label{rtgschema}
\ee

where

\begin{enumerate}

\item $\mcv$ (resp. $\medgecvcvv$) denotes the $\dim$ (resp. $\dim-1$) 
Lebesgue
measure of $\cv$ (resp. $\edgecvcvv$), 

\item 
the quantity $\gcv$, which depends on $\ucvn$ (or $\ucvnp$ or 
$\ucvn$ and
$\ucvnp$), for $\cv \in \T$, is some ``consistent" approximation of $g$ on element
$\cv$, between times $\tn$ and $\tnp$ (we do not discuss this 
approximation here). 

\item  the quantity $\Fcvcvvn$, which
depends on the set of discrete unknowns $u_M^n$ (or $u_M^{n+1}$ or $u_M^n$ and $u_M^{n+1}$) 
for $M \in \T$, is an
approximation of $F\cdot
\ncvcvv$ on $\edgecvcvv$ between times
$\tn$ and $\tnp$.
\end{enumerate}

 In order to obtain a ``good" scheme, this approximation  of $F\cdot
\ncvcvv$ has to be consistent,
conservative  (that is $\Fcvcvvn$=-$\Fcvvcvn$) 
and must ensure some stability properties  on the approximate 
solution given by the scheme
(indeed, one also needs some consistency with respect to entropies, when 
entropies exist\ldots). Except in the scalar case, it is not so easy to see what
kind of stability properties is needed\ldots. Indeed, in the scalar case, that is $m=1$,
taking $g=0$ and $\O=\R^\dim$ (for simplicity), it is essentially sufficient to have an
$L^{\infty}$  estimate
(that is a bound on $\ucvn$ independent of $\cv$, $n$, and of the time and space
discretizations) and a ``touch" of ``$BV$  estimate" (see, for instance,
chapters  \ref{hyper1d} and  \ref{hypmd} and \citet{Chainais1} for more precise assumptions). In the case $m
>1$, it is  not  generally possible  to give stability properties 
from which a mathematical proof
of convergence could be deduced.
However, it is advisable to require some stability properties such as
the positivity of some quantities depending on the unknowns; in the case of
flows, the required stability may be the positivity of the density, energy, pressure\ldots;
 the positivity of these quantities may be essential for the computation of
$F(u)$ or for its hyperbolicity.

~

The computation of $\Fcvcvvn$ is often  performed, at each ``interface", by solving
the following 1D (for the space variable) system (where, for simplicity,  the possible dependency of $F$
with respect to $x$ and $t$ is omitted): 

\be
\frac {\pa u}{\pa t} (z,t) + \frac {\pa f_{\cv,\cvv}(u)} {\pa z} (z,t)=0, \label{fnkl}
\ee

where $f_{\cv,\cvv}(u) (z,t)=   F\cdot \ncvcvv(u(z,t))$, for all $z \in \R$ and $t \in
(0,T)$, which gives consistency, conservativity  (and, hopefully,
stability) of  the final scheme (that
is \eqref{rtgschema}). To be more precise, in the case of lower order schemes, 
$\Fcvcvvn$ may be 
taken as: $\Fcvcvvn = F. \ncvcvv(w)$ where $w$ is the solution for $z = 0$ of  \eqref{fnkl} with
initial conditions $u(x,0) = u_\cv^n$ if $x <0$ and $u(x,0) = u_\cvv^n$ if $x >0$.
Note that the
variable $z$ lies in $\R$, so that the multidimensional problem has therefore been
transformed (as in chapter \ref{hypmd}) into a succession of one-dimensional problems. Hence,
in the following, we shall mainly keep to the case $\dim=1$.

~

Let us describe two classical schemes, namely the Godunov scheme and the Roe scheme, 
in the
case $\dim=1$, $\O=\R$, $F(x,t,u)=F(u)$ and
$g=0$ (but $m \ge 1$), in which case   System \eqref{rtgeq} becomes 

\be
\frac {\pa u}{\pa t} (x,t)
+  \frac {\pa F(u)}{\pa x} (x,t)= 0, 
\ \
 x \in \R, \, t \in (0,T).
\label{rtgeqs}
\ee

in order to complete this system, an initial condition must be specified, the discretization
of  which is standard.

Let $\mesh$ be an admissible mesh in the sense of Definition
\refp{meshhyp1d}, that is $\T=(\cv_i)_{i \in \Z}$, with 
$\cv_i$=$(x_{i-1/2}$,$x_{i+1/2})$, with $x_{i-1/2}<x_{i+1/2}$, $i \in 
\Z$. One sets $h_i= x_{i+1/2}-x_{i-1/2}$, $i \in \Z$. The discrete unknowns are $u_i^n$, $i \in \Z$, $n \in
\{0,\ldots,\Nk+1\}$ and the scheme \eqref{rtgschema} then reads 

\be
h_i\frac{\uinp-\uin}{k}+\fip-\fim=0, \ \
i \in \Z, \, n \in \{0,\ldots,\Nk\},
\label{rtgschemas}
\ee

where $\fipt$ is a consistent  approximation of $F(u(x_{i+1/2},\tn)$.
This scheme is clearly conservative  (in the
sense defined above).   Let us consider  explicit schemes, so that $\fipt$ is a function of
$u_j^n$, $j \in \Z$.  The principle of the Godunov scheme \citet{Go} is to take $\fipt=F(w)$ 
where $w$ is the solution, for $x=0$ (and any 
$t>0$), of the
following (Riemann) problem 

\be
\frac {\pa u}{\pa t} (x,t)
+  \frac {\pa F(u)}{\pa x}(x,t) = 0, 
\ \
 x \in \R, \, t \in \Rp,
\label{rtgrp}
\ee

\be
\bac
u(x,0)= u_i^n, \text{~if } x < 0, \\
u(x,0)= u_{i+1}^n, \text{~if } x > 0.
\ean
\label{rtgbcrp}
\ee

Then, $w$ depends on $u_i^n$, $u_{i+1}^n$ and $F$.

The time step is limited by the so called ``CFL condition", which reads 
$k \le L h_i$, for
all
$i \in \Z$, where $L$ is given by $F$ and the initial condition.
The quantity $\uinp$, given by the Godunov scheme, see \citet{Go}, is, for
all $i \in \Z$, the mean value on $\cv_i$ of the exact solution at time $k$ of \eqref{rtgeqs} 
with
the initial condition (at time $t=0$)
$u_0$ defined, a.e. on $\R$, by $u_0(x)=\uin$ if $x_{i-1/2} <x< x_{i+1/2}$.

\medskip

The Godunov scheme is an efficient scheme (consistent, 
conservative, 
stable), sometimes too
diffusive (especially if $k$ is far from $Lh_i$ defined above), but easy improvements are
possible, such as the MUSCL technique, see below and Section \ref{hyperhos}. Its
principal drawback is its difficult implementation for many problems, indeed the
computation of
$F(w)$ can be impossible or too expensive. For instance, this computation may need a 
non
trivial parametrization of the non linear waves. Note also that $F$ is generally not given
directly as a function of $u$ (the components of
$u$ are called ``conservative  unknowns") but as a function of some ``physical" 
unknowns
(for instance, pressure, velocity, energy\ldots), and the passage from $u$ to these 
physical
unknowns (or the converse) is often not so easy\ldots  it may be the consequence of
expensive and implicit calculations, using, for instance, Newton's 
algorithm\index{Newton's algorithm}.

Due to this difficulty of implementation, some ``Godunov type"
schemes were developed (see \citet{hlvl}). The idea is to take, for $\uinp$, the mean 
value
on $\cv_i$ of an \textit{approximate} solution at time $k$ of \eqref{rtgeqs} with the initial 
condition
(at time $t=0$), $u_0$, defined by $u_0(x)=\uin$, if $x_{i-1/2} <x< x_{i+1/2}$. 
In order for the scheme to be written under the conservative  form \eqref{rtgschemas},  with
a consistent  approximation of the fluxes, this approximate solution must  satisfy some
consistency  relation (another relation is needed for the 
consistency  with entropies). One of
the best known of this family of schemes is the Roe scheme (see 
\citet{roe} and \citet{roe2}), where this
approximate solution is computed by the solution of the following linearized Riemann problems:

\be
\frac {\pa u(x,t)}{\pa t}
+  A(\uin, u_{i+1}^n) \frac {\pa u(x,t)}{\pa
x} = 0, \, 
 x \in \R, \, t \in \Rp,
\label{rtglrp}
\ee

\be
\bac
u(x,0)= u_i^n, \text{~if } x<0, \\
u(x,0)= u_{i+1}^n, \text{~if } x>0,
\ean
\label{rtgbclrp}
\ee

where $A(\cdot,\cdot)$ is an $m \times m$ matrix, continuously depending on its two
arguments, with only real eigenvalues, diagonalizable and satisfying the so called ``Roe
condition":

\be
A(u,v)(u-v)=F(u)-F(v), \forall u,v \in \R^m.
\label{rtgrc}
\ee

Thanks to \eqref{rtgrc}, the Roe scheme \index{Finite volume scheme!-- for 
hyperbolic systems|(}\index{Roe scheme}\index{Scheme!Roe --}
can be written as \eqref{rtgschemas} with 

\be
\bac
\fip=F(\uin)+A^-(\uin,
u_{i+1}^n)(\uin-u_{i+1}^n) \\
(=F(u_{i+1}^n)+A^+(\uin,
u_{i+1}^n)(\uin-u_{i+1}^n)),
\ean
\label{rtgfr}
\ee

where $A^{\pm}$ are the classical nonnegative and nonpositive
parts of the matrix $A$: let $A$ be a matrix with only real
eigenvalues, $(\l_{p})_{p=1,\ldots,m}$, and diagonalizable, let
$(\varphi_p)_{p=1,\ldots,m}$ be a basis of $\R^m$ associated to these 
eigenvalues. Then, the
matrix $A^+$ is the matrix which has the same eigenvectors as $A$ and has
$(\max\{\l_{p},0\})_{p=1,\ldots,m}$ as corresponding
eigenvalues. The matrix $A^-$ is $(-A)^+$.

Roe's scheme was proved to be an efficient scheme, often less expensive than Godunov's scheme, with, 
more or
less  the same limitation on the time step, the same diffusion effect and some lack of
entropy consistency, which can be corrected. It has some properties of consistency  and
stability. Its main drawback is the difficulty of the computation of a matrix $A(u,v)$
satisfying
\eqref{rtgrc}. For instance, when it is possible to
compute and diagonalize the derivative of $F$, $DF(u)$, one can take 
$A(u,v)=DF(u^{\star})$,
but the difficulty is to find $u^{\star}$ such that \eqref{rtgrc} holds (note that this
condition is crucial in order to ensure conservativity  of Roe's scheme). In some difficult
cases, the Roe matrix is computed approximately  by using a  ``limited expansion" with
respect to   some small parameter.

\subsection{Rough schemes for complex hyperbolic systems}\label{rough}

The aim of this section is to present some discretization techniques for ``complex"
hyperbolic systems. 
In many applications, the expressions of $g$ and $F$ which appear in \eqref{rtgeq} are rather
``complex", and it is difficult or impossible to use classical schemes such as the 1D
Godunov or Roe schemes  or their standard extensions, for multidimensional problems, 
using
1D solvers on the interfaces of the mesh. This is the case of gas dynamics (Euler 
equations\index{Euler equations})
with real gas, for which the state law (pressure as a function of density and internal
energy) is tabulated or given by some complex analytical expressions. This is also  the
case when modelling   multiphase flows in pipe-lines: the function $F$ is difficult to
handle and highly depends on $x$ and $u$, because, for instance, of changes of the geometry
and  slope
of the pipe, of changes of the friction law or, more generally, of the varying nature of
the flow. Most of the attempts given below were developed for this last 
situation. 
Other interesting cases of ``complexity" are the treatment of boundary conditions
(mathematical literature is rather scarce on this subject, see Section 
\ref{hypbc} for a first insight), 
and the way to handle the case where the eigenvalues (of the derivative of
$F \cdot \n$ with respect to its third argument) are of very different magnitude, see
Section \ref{turbo}. Another case of complexity is the treatment of nonconservative
terms in the equations. One refers, for instance, to
\citet{h2} and references therein, for this important case.

\medskip

Possible modifications of Godunov and Roe schemes (including
``classical" improvements to avoid excessive artificial diffusion) are described now to
handle ``complex" systems. 
Because of the complexity of the models, the justification of the
schemes presented here is rather numerical than mathematical. Many variations have also
been developed, which are  not presented here. Note that other approaches 
are also possible, see e.g. \citet{Ghidaglia}.
For simplicity, one considers the case $\dim=1$, $\O=\R$, $F(x,t,u)=F(u)$
and
$g=0$ (but $m \ge 1$) described in Section \ref{godroe}, with the same notations.
The Godunov and Roe schemes can both be written under the form
\eqref{rtgschemas} with $\fipt$ computed as a function of $\uin$ and $\uipn$; both 
schemes are
consistent (in the sense of Section \ref{godroe}, i.e. consistency  of the
``fluxes") since $\fipt=F(u)$ if $\uin=\uipn=u$.

\medskip

Going further along this line of thought
yields (among other possibilities, see below) the ``VFRoe" 
scheme\index{Scheme!VFRoe --} which is
\eqref{rtgschemas}, that is:
\be
h_i\frac{\uinp-\uin}{k}+\fip-\fim=0, \ \
i \in \Z, \, n \in \{0,\ldots,\Nk\},
\label{rtgvfroe1}
\ee
 with
$\fipt=F(w)$, where $w$ is the solution of the linearized Riemann 
problem \eqref{rtglrp}, \eqref{rtgbclrp}, with
$A(\uin,\uipn)=DF(w^{\star})$, that is:

\be
\frac {\pa u(x,t)}{\pa t}
+ DF(w^{\star}) \frac {\pa u(x,t)}{\pa
x} = 0, \, 
 x \in \R, \, t \in \Rp,
\label{rtgvfroe2}
\ee

\be
\bac
u(x,0)= u_i^n, \text{~if } x<0, \\
u(x,0)= u_{i+1}^n, \text{~if } x>0,
\ean
\label{rtgvfroe3}
\ee

where $w^{\star}$ is some value between $\uin$ and 
$\uipn$
(for instance, $w^{\star}= (1/2) (\uin+\uipn)$). In this scheme, the 
Roe condition \eqref{rtgrc} 
is
not required (note that it is naturally conservative, thanks to its finite volume
origin). 
Hence, the VFRoe scheme \index{Scheme!VFRoe --} appears to be a
simplified version of the Godunov  and
Roe schemes. The study of the scalar case ($m=1$) shows that, in order to have some
stability, at least as much as in Roe's scheme, the choice of
$w^{\star}$ is essential. In practice, the choice $w^{\star}=(1/2) (\uin+\uipn)$ is
often adequate, at least for regular meshes.

\begin{remark} \normalfont
In Roe's scheme, the Roe condition \eqref{rtgrc} ensures conservativity\index{Conservativity}. 
The VFRoe scheme \index{Scheme!VFRoe --} is ``naturally" conservative, and therefore no such 
condition is needed. Also note that the VFRoe scheme \index{Scheme!VFRoe --} yields precise 
approximations of the shock velocities, without Roe's condition.
\end{remark}

Numerical tests show the good behaviour of the VFRoe \index{Scheme!VFRoe --} scheme. Its two main flaws are a 
lack
of entropy consistency (as in Roe's scheme) and a large diffusion effect (as in the 
Godunov
and Roe schemes). 
The first drawback can be corrected, as for Roe's scheme, with a
nonparametric entropy correction inspired from \citet{hh} (see \citet{mfg}).
The two
drawbacks can be corrected with a classical MUSCL technique, which consists in 
replacing, in
\refep{rtgbclrp}, $\uin$ and $\uipn$ by $u_{i+1/2,-}^n$ and $u_{i+1/2,+}^n$, which
depend on $\{u_j^n$, $j=i-1,i,i+1,i+2\}$ (see, for instance, Section
\refp{hyperhos} and  \citet{godlewski-springer} or  \citet{leveque}). For stability reasons, the
computation of the gradient of the unknown (cell by cell) and of the ``limiters" is
performed on some ``physical" quantities (such as density, pressure, velocity for 
Euler
equations\index{Euler equations}) instead of
$u$. The extension of the MUSCL technique to the case $\dim>1$ is more or less
straightforward.

\medskip

This MUSCL technique improves the space accuracy (in the truncation error) and the 
numerical
results are significantly better. However, stability  is sometimes lost. Indeed, considering
the linear scalar equation, one remarks that the scheme is antidiffusive when the limiters
are not active, this might lead to a loss of stability. The time step must then be
reduced (it is reduced by a factor
$10$ in severe situations\ldots). In order to allow larger time steps, the time accuracy
should be improved by using, for instance, an order 2 Runge-Kutta scheme (in the severe
situations suggested above, the time step is then multiplied by a factor $4$). 
Surprisingly,
this improvement of time accuracy is used to gain stability rather than precision\ldots

\medskip

Several numerical experiments (see \citet{mfg}) were performed which prove the 
efficiency of
the VFRoe scheme\index{Scheme!VFRoe --}, such as the classical Sod tests (\citet{sod}). The shock velocities are
exact, there are no oscillations\ldots. For these tests, the treatment of the
boundary conditions is straightforward. Throughout these experiments, the use of a
MUSCL technique yields a significant improvement, while the use of a higher order time
scheme is not necessary. In one of the Sod tests, the entropy correction is needed.

A comparison between the VFRoe scheme \index{Scheme!VFRoe --} and the Godunov scheme was performed by 
J. M.
H\'{e}rard (personal communication) for the Euler equations\index{Euler 
equations} on a Van Der Wals gas, for
which a matrix satisfying \eqref{rtgrc} seems difficult to find. The numerical results are
better with the VFRoe schem\index{Scheme!VFRoe --}, which is also much cheaper computationally.
An improvment of the VFRoe scheme is possible, using, instead of 
\eqref{rtgvfroe2}-\eqref{rtgvfroe3}, linearized Riemann 
problems associated to a nonconservative form of the initial system, 
namely System \eqref{rtgeqs} or more generally System \eqref{rtgeq}, for 
the computation of $w$ (which gives the flux $\fipt$ in \eqref{rtgvfroe1} 
by the formula $\fipt=F(w)$), see for instance \citet{bgh} for a simple 
example.


In some more complex cases, the flux $F$ may also highly, and not continuously, 
depend on the space variable
$x$. In the space discretization, it is ``natural" to set the 
discontinuities of $F$ with respect to $x$ on the boundaries of the mesh. 
The function $F$ may change drastically from $\cv_i$ to $\cv_{i+1}$. In this 
case, the implementation of  the VFRoe scheme \index{Scheme!VFRoe --} yields two
additional difficulties: 

\benum
\item[(i)] The matrix $A(\uin,\uipn)$ in the linearized Riemann problem 
\eqref{rtglrp}, \eqref{rtgbclrp} now depends on $x$:

$A(\uin,\uipn)$ $=$ $D_uF(x,w^{\star})$, where
$w^{\star}$ is some value between
$\uin$ and $\uipn$ and $D_uF$ denotes the derivative of $F$
with respect to its ``$u$" argument. 

\item[(ii)] once the solution, $w$, of the linearized problem \eqref{rtglrp}
\eqref{rtgbclrp}, for $x=0$ and any
$t>0$, is calculated, the choice $\fipt=F(x,w)$ again depends on
$x$.
\eenum

The choice of $\fipt$ (point (ii)) may be solved by remarking that, in Roe's scheme,
$\fipt$ may be written (thanks to \eqref{rtgrc}) as

\be
\fip=\frac 1 2 (F(\uin)+F(\uipn)) +\frac {1}{2} A_{i+\frac 1 2}^n (\uin-\uipn),
\label{rtgfrc}
\ee

where $A_{i+1/2}^n=\vert A(\uin,\uipn) \vert$, and $\vert A \vert=A^+ + A^-$.

Under this form, the second term of the right hand side of \eqref{rtgfrc} appears to be 
a stabilization
term, which does not affect the consistency\index{Consistency!-- of the 
fluxes}. Indeed, in the scalar case ($m=1$), one has
$A_{i+1/2}^n= \vert F(\uin)-F(\uipn) \vert/\vert \uin-\uipn
\vert$, which easily yields the $L^{\infty}$ stability   of the scheme (but not the
consistency with respect to the entropies). Moreover, the scheme is stable and consistent
with respect to the entropies, under a Courant-Friedrichs-Levy (CFL) condition, if $\fipt$ is nondecreasing with respect to
$\uin$ and nonincreasing with respect to $\uipn$, which holds if $A_{i+1/2}^n \ge
\sup\{\vert F'(s) \vert,$ $s \in [\uin,\uipn]$ or $[\uipn, \uin]\}$.
This remark suggests
a slightly different version of the VFRoe\index{Scheme!VFRoe --}scheme (closer to Roe's
scheme), which is the scheme \eqref{rtgvfroe1}-\eqref{rtgvfroe3}, taking  

\[
\fipt = \frac 1 2 (F(\uin)+F(\uipn)) +\frac {1}{2} \vert DF(w^{\star}) \vert (\uin-\uipn),
\]

in \eqref{rtgvfroe1}, instead of $\fipt=F(w)$.  Note that it is also possible  to take other
convex combinations of $F(\uin)$ and $F(\uipn)$ in the latter expression of $\fipt$, without 
modifying the consistency  of the scheme.

\medskip

When $F$ depends on $x$, the discontinuities of $F$ being on the boundaries of the
control volumes, the generalization of
\eqref{rtgfrc} is obvious, except for  the choice of $A_{i+1/2}^n$. The quantity
$F(\uin)$ is replaced by $F(x_i,\uin)$, where $x_i$ is the center of $\cv_i$.
Let us now turn to the choice of a convenient matrix $A_{i+1/2}^n$ for this 
modified
VFRoe scheme\index{Scheme!VFRoe --}, when $F$ highly depends on $x$. A first possible choice is
\[A_{i+1/2}^n=(1/2) (\vert D_u F(x_i,\uin)\vert + \vert 
D_u F(x_{i+1},\uipn)\vert).\]
The following slightly different choice for $A_{i+1/2}^n$ seems, however, to give better
numerical results (see \citet{fh}). Let us define
\[
A_i=D_uF(x_i,\uin), \forall \, i \in \Z
\]
(for the determination of $A_{i+1/2}^n$ 
the fixed index $n$ is omitted). Let $(\l^{(i)}_{p})_{p=1,\ldots,m}$ be the 
eigenvalues of
$A_i$ (with $\l^{(i)}_{p-1} \le
\l^{(i)}_{p}$, for all $p$) and $(\varphi^{(i)}_{p})_{p=1,\ldots,m}$ a basis of $\R^m$
associated to these eigenvalues. Then, the matrix $A_{i+1/2}^{(-)}$ [resp. 
$A_{i+1/2}^{(+)}$] is the matrix which has the same eigenvectors as $A_i$ [resp. $A_{i+1}$] and 
has
$(\max\{\vert\l^{(i)}_{p}\vert,\vert\l^{(i+1)}_{p}\vert\})_{p=1,\ldots,m}$ as 
corresponding
eigenvalues. The choice of $A_{i+1/2}^n$ is  

\be
A_{i+\frac 1 2}^n=\frac \lambda 2 (A_{i+\frac 1 2}^{(-)} +
A_{i+\frac 1 2}^{(+)}),
\label{rtggroomf}
\ee

where $\lambda$ is a parameter, the ``normal" value of which is $1$. Numerically, larger 
values
of $\lambda$, say $\lambda=2$ or $\lambda=3$, are sometimes needed, in severe situations, to obtain 
enough
stability. Too large values of $\lambda$ yield too much artificial diffusion.

~

The new scheme is then \eqref{rtgvfroe1}-\eqref{rtgvfroe3}, taking  
\be
\fipt = \frac 1 2 \Bigl(F(x_i,\uin)+F(x_i,\uipn)\Bigr) +\frac {1}{2} A_\ipud^n (\uin-\uipn).
\label{rtgfrcb}
\ee
where $ A_{i+1/2}^n $ is defined by \eqref{rtggroomf}.
It has, more or 
less,
the same properties as the Roe and VFRoe \index{Scheme!VFRoe --} schemes but allows the simulation of more 
complex
systems. It needs a MUSCL technique to reduce diffusion effects and order 
2 Runge-Kutta for
stability. It was implemented for the simulation of multiphase flows in pipe lines (see
\citet{fh}). The other difficulties encountered in this case are the treatment of
the boundary conditions and the different magnitude of the eigenvalues, which are
discussed
in the next sections.

\subsection{Partial implicitation of explicit scheme}\label{turbo}

In the modelling of flows, where ``propagation" 
phenomena and  ``convection" phenomena coexist, the Jacobian matrix of $F$ often has
eigenvalues of
different magnitude, the  ``large" eigenvalues (large meaning ``far from
$0$", positive or negative) corresponding to the  propagation  
phenomena and   ``small" eigenvalues corresponding to the ``convection" 
phenomena\index{Convection term} . 
Large and
small eigenvalues may differ by a factor $10$ or $100$.

~

With the explicit schemes described in the previous sections, the time step is limited by
the CFL condition corresponding to the large eigenvalues. Roughly speaking, with the
notations of Section \ref{godroe}, this condition is (for all $i \in \Z$) $k \le \vert \l 
\vert^{-1}h_i$, where $\l$
is the largest eigenvalue. In some cases, this limitation can be unsatisfactory for two
reasons. Firstly, the time step is too small and implies a prohibitive computational cost.
Secondly, the discontinuities in the solutions, associated to the small eigenvalues, are
not sharp because the time 
step is far from the CFL condition of the small eigenvalues (however, 
this can be somewhat corrected with a MUSCL method). This is in fact a major problem
when the discontinuities associated to the small eigenvalues need to be computed
precisely. It is the case of interest here.

~

A first method to avoid the time step limitation is to take a ``fully implicit" version of
the schemes developed in the previous sections, that is $\fipt$ function of $u_j^{n+1}$, 
$j
\in \Z$, instead of $u_j^{n}$, $j
\in \Z$ (the terminology ``fully implicit" is by opposition to ``linearly
implicit", see below and \citet{fer}). However, in order to be
competitive with  explicit schemes, the fully implicit scheme is used
with large time steps. In practice, this prohibits the use of a MUSCL technique in the
computation of the solution at time $\tnp$ by, for instance, a Newton algorithm\index{Newton's algorithm}. This
implicit scheme is therefore very diffusive and will smear discontinuities.

A second method consists in splitting the system into two systems, the first one
is associated with the ``small" eigenvalues, and the second one with the ``large"
eigenvalues (in the case of the Euler equations\index{Euler equations}, this splitting may correspond to a
``convection" system and a ``propagation" system). At each time step, the first system is
solved with an explicit scheme and the second one with an implicit scheme. Both use the
same time step, which is limited by the CFL condition of the small eigenvalues. Using
a MUSCL technique and an order 2 Runge-Kutta method for the first system yields sharp
discontinuities associated to the small eigenvalues. This method is often satisfactory, but
is difficult to handle in the case of severe boundary conditions, since the convenient
boundary conditions for each system may be difficult to determine.

Another method, developed by E. Turkel (see \citet{tur}), in connexion with Roe's
scheme, uses a change of variables in order to reduce the ratio between large and small
eigenvalues.

~

Let us now describe a partially linearly implicit method (``turbo" scheme) which
was successfully tested for multiphase flows in pipe lines (see \citet{fh}) and other cases
(see \citet{fer}). For the sake of simplicity, the method is described for the last scheme
of Section \ref{rough}, i.e. the scheme defined by \eqref{rtgvfroe1}-
\eqref{rtgvfroe3}, where $\fip$ is defined by \eqref{rtgfrcb} and \eqref{rtggroomf}
 (recall that $F$ may
depend on $x$).

\medskip

Assume that $I\subset \{1,\ldots,m\}$ is the set of index of large eigenvalues (and
does not depend on $i$). The aim here is to ``implicit" the unknowns coresponding to
the large  eigenvalues  only: let
$\Ab_{i}$, $\Ab_{i+1/2}^{(-)}$ and $\Ab_{i+1/2}^{(+)}$ be
the matrix having the same eigenvectors as $A_i$, $A_{i+1/2}^{(-)}$ and 
$A_{i+1/2}^{(+)}$, with the same large eigenvalues (i.e. corresponding to $p \in I$) and $0$ as small
eigenvalues. Let 

\[
\Ab_{i+1/2}^n= (\lambda/2) (\Ab_{i+1/2}^{(-)} +
\Ab_{i+1/2}^{(+)}).
\]

Then, the partially linearly implicit scheme is obtained by
replacing $\fipt$ in
\eqref{rtgschemas} by $\fipbt$ defined by

\[
\bac
\fipb=\fip + \frac 1 2 (\Ab_i(\uinp-\uin) + \Ab_{i+1}(u_{i+1}^{n+1}-\uipn))
\\
+ \frac 1 2 \Ab_{i+\frac 1 2}^n(\uinp-\uin+\uipn-u_{i+1}^{n+1}).
\ean
\]

In order to obtain sharp discontinuities corresponding to the small eigenvalues, a
MUSCL technique is used for the computation of $\fipt$. Then, again for stability 
reasons,
it is preferable to add an order 2 Runge-Kutta method for the time
discretization. Although it is not so easy to implement, the order 2
Runge-Kutta method is needed to enable the use of ``large" time steps. 
The time step is, in severe situations, very close to that given by the 
usual CFL condition corresponding to the small eigenvalues, and can be 
considerably larger than that given by the large eigenvalues (see 
\citet{fh} for several tests).

\subsection{Boundary conditions}\label{hypbc}

In many simulations of real situations, the treatment of the boundary conditions is not
easy (in particular in the case of sign change of eigenvalues). We give here 
a classical possible mean (see e.g. \citet{kum} and \citet{dlf}) of handling boundary conditions (a more
detailed  description may be found in \citet{jmm} for the case of multiphase flows in 
pipe lines).

~

Let us consider now the system \eqref{rtgeqs} where ``$x \in \R$" is replaced by 
``$x \in \O$" with $\O=(0,1)$. In order for the system to be well-posed,  an initial
condition (for $t=0$) and some convenient boundary  conditions for $x=0$ and $x=1$ are
needed;  these boundary conditions will appear later in the  discretization (we do not
detail here the mathematical analysis of the  problem of the adequacy of the boundary
conditions, see e.g. \citet{serre} and references therein). 
Let us now explain the numerical 
treatment of the boundary condition at $x=0$.

~

With the notations of Section \ref{godroe}, the space mesh is given by 
$\{\cv_i$, $i \in \{0,\ldots, N_\mesh\}\}$, 
with $\sum_{i=1}^{N_\mesh}h_i=1$. Using the finite volume 
scheme \eqref{rtgschemas} with $i \in \{1,\ldots,N_\mesh\}$ instead of $i 
\in \Z$ needs, for the computation of $u_1^{n+1}$, with $\{u_i^n$, $i \in 
\{1,\ldots, N_\mesh\}\}$ given, a value for $F_{1/2}^n$ (which 
corresponds to the flux at point $x=0$ and time $t=\tn$).

~

For the sake of simplicity, consider only the case of the Roe and VFRoe 
\index{Scheme!VFRoe --}
schemes. Then, the ``interior fluxes", that is 
$\fipt$ for $i \in \{1,\ldots, N_\mesh-1\}$, are determined by using 
matrices $A(\uin, \uipn)$ ($i \in \{1,\ldots, N_\mesh-1\}$). In the case 
of the Roe scheme, $\fipt$ is given by \eqref{rtgfr} or \eqref{rtgfrc} and 
$A(\cdot,\cdot)$ satisfies the Roe condition \eqref{rtgrc}. In the case of 
the VFRoe \index{Scheme!VFRoe --} scheme, $\fipt$ is given through the  resolution of the linearized 
Riemann problem \eqref{rtglrp}, 
\eqref{rtgbclrp} with e.g. $A(\uin,\uipn)=DF((1/2) (\uin+\uipn))$. In order to compute
$F_{1/2}^n$, a possibility is to take the same method as  
for  the interior fluxes; this requires
the determination of some $u_0^n$.
In some cases (e.g. when all the eigenvalues 
of $D_uF(u)$ are nonnegative), the given boundary conditions at $x=0$ are sufficient to 
determine the value $u_0^n$, or directly $F_{1/2}^n$, but this is not true in the 
general case\ldots. 
In the general case, there are not enough given boundary conditions to determine 
$u_0^n$ and   missing equations need to be introduced. The idea is to use an iterative
process. Since  $A(u_0^n,u_1^n)$ is diagonalizable and has only real eigenvalues, let 
$\l_1,\ldots,\l_m$ be the eigenvalues of $A(u_0^n,u_1^n)$ and 
$\varphi_1,\ldots, \varphi_m$ a basis of $\R^m$ associated to these 
eigenvalues. Then the vectors $u_0^n$ and $u_1^n$ may be decomposed on this basis, this yields

\[
u_0^n= \sum_{i=1}^m \a_{0,i}\  \varphi_i, \;  u_1^n= \sum_{i=1}^m 
\a_{1,i} \  \varphi_i.
\]

Assume that the number of  negative eigenvalues of $A(u_0^n,u_1^n)$  does not
depend on  $u_0^n$ (this is 
a  simplifying assumption); 
let   $p$ be the number of negative eigenvalues and $m-p$  the number
of positive eigenvalues  of $A(u_0^n,u_1^n)$.

Then, the number of (scalar) given boundary conditions is (hopefully \ldots) 
$m-p$. Therefore, one takes, for $u_0^n$, the solution of the (nonlinear) 
system of $m$ (scalar) unknowns,  and $m$ (scalar) equations. 
The $m$ unknowns are  the components of $u_0^n$ and 
the $m$ equations are obtained with the $m-p$
boundary  conditions and the $p$ following equations:

\be
\a_{0,i}=\a_{1,i}, \text{ if } \l_i < 0.
\label{rtgbc}
\ee

Note that the quantities $\a_{0,i}$ depend on  $A(u_0^n,u_1^n)$; the resulting system is 
therefore 
nonlinear and may be solved with, for instance, a Newton algorithm\index{Newton's algorithm}.

~

Other possibilities around this method are possible. For instance, 
another possibility, perhaps more natural, consists in writing the $m-p$
boundary  conditions on $u_{1/2}^n$ instead of $u_0^n$ and to take 
\eqref{rtgbc} with the components of $u_{1/2}^n$ instead of those of 
$u_0^n$, where $u_{1/2}^n$ is the solution at $x=0$ of \eqref{rtglrp}, 
\eqref{rtgbclrp} with $i=0$. With the VFRoe scheme\index{Scheme!VFRoe --}, the flux 
at the boundary $x=0$ is then $F_{1/2}^n=F(u_{1/2}^n)$. In the case of a 
linear system with linear boundary conditions and with the VFRoe scheme, 
this method gives the same flux $F_{1/2}^n$ as the preceding method, the value $u_{1/2}^n$ is 
completely determined although $u^n_0$ is not completely determined.

~

In the case of the scheme described in the second part of Section 
\ref{rough}, the following ``simpler" possibility was implemented. For this scheme, $\fipt$ is 
given, for $i \in \{1,\ldots, N_\mesh -1\}$, by \eqref{rtgfrc} with 
\eqref{rtggroomf}. Then, the idea is to take the same equation for the 
computation of $F_{1/2}^n$ but to compute $u_0^n$ as above (that is with 
$m-p$ boundary conditions and \eqref{rtgbc}) with the choice $A(u_0^n,u_1^n)=D_uF(x_1,u_1^n)$.

~

This method of computation of the boundary fluxes gives good results but 
is not adapted to all cases (for instance, if $p$ changes during the 
Newton iterations \index{Newton's algorithm} or if the number of boundary conditions is not equal to 
$m-p$\ldots). Some particular methods, depending on the 
problems under consideration, have to be developped.

~

We now  give an attempt for the justification of this treatment of the
boundary conditions, at least for a linear system with linear boundary 
conditions.

Consider the system

\be\bal
u_t(x,t) + u_x (x,t) = 0, \, x \in (0,1), \, t \in \Rp,
\\
v_t(x,t)-v_x(x,t)=0,\, x \in (0,1), \, t \in \Rp,
\ean\label{sl}\ee

with the boundary conditions

\be\bal
u(0,t)+\a v(0,t)=0,  \, t \in \Rp,
\\
v(1,t)+\b u(1,t)=0,  \, t \in \Rp,
\ean\label{slcl}\ee

and the initial conditions

\be\bal
u(x,0)=u_0(x),  \, x \in (0,1),
\\
v(x,0)=v_0(x),  \, x \in (0,1),

\ean\label{slci}\ee

where $\a \in \R^\star$, $\b \in \R^\star$, $u_0 \in L^\infty(\O)$ and $v_0 \in 
L^\infty(\O)$ are given.
It is well known that the problem \eqref{sl}-\eqref{slci} admits a unique 
weak solution (entropy conditions are not necessary to obtain uniqueness 
of the solution of this linear system).

~ 

A stable numerical scheme for the discretization of the problem 
\eqref{sl}-\eqref{slci} will add some numerical diffusion terms. It seems 
quite natural to assume that this diffusion does not lead a coupling 
between the two equations of \eqref{sl}. Then, roughly speaking, the 
numerical scheme will consist in an approximation of the following 
parabolic system:

\be\bal
u_t(x,t) + u_x (x,t) -\eps u_{xx}(x,t)= 0, \, x \in (0,1), \, t \in \Rp,
\\
v_t(x,t)-v_x(x,t)- \eta v_{xx}(x,t)=0,\, x \in (0,1), \, t \in \Rp,
\ean\label{slp}\ee

for some $\eps>0$ and $\eta>0$ depending on the mesh (and time step) and 
$\eps \tends 0$, $\eta \tends 0$ as the space and time steps tend to 0.

In order to be well posed, this parabolic system has to be completed 
with the initial conditions 
\eqref{slci} and (for all $t >0$) four boundary conditions, i.e. two conditions 
at $x=0$ and two  conditions at $x=1$. This is also the case for the 
numerical scheme which may be viewed as a discretization of \eqref{slp}.
There are two boundary conditions given by \eqref{slcl}. 
Hence  two  other boundary conditions
must be found,
one at $x=0$ and the other at $x=1$.

~ 
 
If these  two  additional conditions are, for instance, $v(0,t)=u(1,t)=0$, 
then the (unique) solution to  \eqref{slcl}-\eqref{slp} with these two additional 
conditions does not converge, as $\eps \tends 0$ and $\eta \tends 0$, to 
the weak solution of \eqref{sl}-\eqref{slci}. This negative result is also 
true for a large choice of other additional boundary conditions.
However, if the additional boundary conditions are (wisely) chosen to be
$v_x(0,t)=u_x(1,t)=0$, 
the solution to  \eqref{slcl}-\eqref{slp} with these two additional 
conditions converges  to the weak solution of \eqref{sl}-\eqref{slci}.

~

The numerical treatment of the boundary conditions described above may 
be viewed as a discretization of \eqref{slcl} and $v_x(0,t)=u_x(1,t)=0$;
this remark gives a formal
justification to such a choice.

\subsection{Staggered grids}

\index{Staggered grid}
For some systems of equations it may be ``natural" (in the sense that the discretization seems simpler) to associate different grids to different unknowns of the problem. 
To each unknown is associated an equation and this equation is
integrated over the elements (which are the control volumes) of the corresponding
mesh, and then discretized by using one discrete unknown per control volume (and time
step, for evolution problems). 
This is the case, for instance, of the well known discretization of the incompressible
Navier-Stokes equations\index{Navier-Stokes equations}  with staggered grids, 
see \citet{Pat} and Section 
\ref{sstokes}.

Let us now give an example in order to show that staggered grids should be avoided in
the case of nonlinear hyperbolic systems since they  may
yield some kind of ``instability". As an illustration, 
let us consider the following ``academic" problem:

\be
\bal
u_t (x,t) + (vu)_x (x,t) =  0, \, x \in \R, \, t \in \R_+,
\\
v_t (x,t) + (v^2)_x (x,t) =  0, \, x \in \R, \, t \in \R_+,
\\
u(x,0)  =  u_0(x), \, x \in \R,
\\
v(x,0)  =  u_0(x), \, x \in \R,
\ean
\label{sg}
\ee

where $u_0$ is a bounded function from $\R$ to $[0,1]$.
Taking $u=v$ equal to the weak entropy solution of the B\"{u}rgers 
equation (namely $u_t + (u^2)_x=0$), with initial
condition $u_0$, leads to a solution of the problem \eqref{sg}. One would
expect a numerical scheme to give an approximation of this solution. 
Note that the solution of the B\"{u}rgers equation, with initial
condition $u_0$,  also takes its values in $[0,1]$, and hence, a ``good"
numerical scheme can be expected to give approximate 
solutions taking values in $[0,1]$. Let us show that this property is not satisfied when
using staggered grids.

~

Let $k$ be the time step and $h$ be the (uniform) space step. Let $x_i=ih$ and
$x_{i+1/2}=(i+1/2)h$, for $i \in \Z$. Define, for $i \in \Z$,
$\cv_i=(x_{i-1/2},x_{i+1/2})$ and $\cv_{i+1/2}=(x_{i},x_{i+1})$.

The mesh associated to $u$ is $\{\cv_i, i \in \Z \}$ and the mesh associated to $v$ is
$\{\cv_{i+1/2}, i \in \Z \}$. Using the principle of staggered grids\index{Staggered grid}, the discrete
unknowns are $u_i^n$, $i \in \Z$, $n \in \N^{\star}$, and 
$v_{i+1/2}^n$, $i \in
\Z$,  $n \in \N^{\star}$. The discretization of the initial conditions is, for
instance,

\be
\bal
\dsp{u_i^0 = \frac 1 h \int_{\cv_i} u_0(x) dx}, \, i \in \Z,
\\
\dsp{v_{i+\frac 1 2}^0 = \frac 1 h \int_{\cv_{i+\frac 1 2}} u_0(x) dx}, \, 
i \in \Z.
\ean
\label{sgic}
\ee

The second equation of (\ref{sg}) does not depend on $u$. It seems reasonable to
discretize this equation with the Godunov scheme, which is here the upstream
scheme, since $u_0$ is nonnegative.  The discretization of the first equation of
(\ref{sg}) with the principle of staggered\index{Staggered grid} grids is easy. Since $v_{i+1/2}^n$ is
always nonnegative, we also take an upstream\index{Upstream} value for $u$ at the extremities of the
cell $\cv_i$. Then, with the explicit Euler scheme  in time, the scheme becomes 

\be
\bal
\dsp{\frac 1 k (u_i^{n+1}-u_i^n) + 
\frac 1 h(v_{i+\frac 1 2}^nu_i^n - v_{i-\frac 1 2}^nu_{i-1}^n)=0, 
\, i \in \Z, \, n \in \N,}
\\
\\
\dsp{\frac 1 k (v_{i+\frac 1 2}^{n+1}-v_{i+\frac 1 2}^n) + 
\frac 1 h((v_{i+\frac 1 2}^n)^2 -(v_{i-\frac 1 2}^n)^2)=0,
\, i \in \Z, \, n \in \N.}
\ean
\label{sgd}
\ee

It is easy to show that, whatever $k$ and $h$, there exists $u_0$ (function
from $\R$ to $[0,1]$) such that sup$\{u_i^1, i \in \Z \}$ is strictly larger than $1$.
In fact, it is possible to have, for instance, sup$\{u_i^1, i \in \Z 
\}=1+ k/(2h)$.
In this sense the scheme (\ref{sgd}) appears to be unstable. Note that the same 
phenomenon
exists with the implicit Euler scheme  instead of the explicit Euler scheme
\index{Finite volume scheme!-- for 
hyperbolic systems|)}. Hence
staggered grids\index{Staggered grid} do not seem to 
be the best choice for nonlinear hyperbolic systems.

\section{Incompressible Navier-Stokes  Equations}\label{nseq}

\hpage{Systems}
\index{Navier-Stokes equations}

The discretization of the stationary Navier-Stokes  equations by the finite volume method is
presented in this section. We first recall the classical discretization on cartesian
staggered grids\index{Staggered grid}. We then study, in the linear case of the Stokes equations, a finite volume
method on a staggered triangular grid, for which we show, in a particular case, the
convergence of the method.

\subsection{The continuous equation}\label{cstokes}

Let us consider here the stationary Navier-Stokes equations:

\be\bal        
\dsp{- \nu \Delta u^{(i)}(x) + \sum_{j=1}^d u^{(j)}(x) \dxj{ u^{(i)} }(x) +  \dxi{p} (x)= 
f^{(i)}(x), \, x \in \O, \, \forall i=1,\ldots,d,}
\\
\sumund \dxi{ u^{(i)} }(x) = 0, \, x \in \O.
\ean
\label{nstkc}
\ee  

with Dirichlet boundary condition 

\be        
u^{(i)}(x) = 0, \,  x \in \dr\O,  \, \forall   i=1,\ldots,d,
\label{stkb}\ee 

under the following assumption:
\begin{hyp} \label{stkH}
\[
\ban
\text{(i)} & \O \ \text{is an open bounded connected polygonal subset of} \ \R^d, 
d = 2, 3,\\ 
\text{(iii)} &\nu > 0,\\
\text{(iii)} & f^{(i)} \in L^2(\O), \, \forall i=1,\ldots,d.\\
\ean 
\]

\end{hyp} 

In the above equations, $u^{(i)}$ represents the $i$th component  of the velocity of a
fluid,
$\nu$ the kinematic viscosity and $p$ the pressure. The unknowns of the 
problem are $u^{(i)}$, $i \in \{1,\ldots,d\}$ and $p$. The number of 
unknown functions from $\O$ to $\R$ which are to be computed  is therefore 
$d+1$. Note that \eqref{nstkc} yields 
$d+1$ (scalar) equations.

~

We shall  also  consider the
Stokes equations, which are obtained by neglecting the nonlinear 
convection term\index{Convection term}.

\be\bal       
\dsp{- \nu \Delta u^{(i)}(x) + \dxi{p}  (x)= f^{(i)}(x), \,  x \in \O, \, 
\forall   i=1,\ldots,d,}
\\
\sumund \dxi{u^{(i)}} = 0, \,   x \in \O.
\ean
\label{stkc} \ee  

There exist several convenient mathematical formulations of \eqref{nstkc}-\eqref{stkb}
and  \eqref{stkc}-\eqref{stkb}, see e.g. \citet{temam}.  Let us give one of them for the
Stokes problem. Let

\[
V=\{ u=(u^{(1)},\ldots,u^{(d)})^t \in (H^1_0(\O))^d, \; \sumund\dxi { u^{(i)} }= 0
\}.
\]
Under assumption \ref{stkH},  there exists a unique function $u$ such that  

\be \ban
u \in V, \\
\nu \sumund \int_\O \grad u^{(i)}(x) \cdot \grad v^{(i)}(x) dx =   \sumund \int_\O
f^{(i)}(x) v^{(i)}(x) dx, \;     \forall v=(v^{(1)},\ldots,v^{(d)})^t\in V.
\ean
\label{stksolvar}
\ee

Equation \eqref{stksolvar} yields the existence of  $p\in L^2$ (unique if  $\int_\O p(x)
dx = 0$) such that

\be
-\nu \Delta u^{(i)} + \dxi{p} = f^{(i)} \text{ in } \D'(\O), \, \forall i 
\in \{1,\ldots,d\}.
\label{stkp}\ee

In the following, we shall study finite volume schemes for  the discretization 
 of Problem  \eqref{nstkc}-\eqref{stkb} and   \eqref{stkc}-\eqref{stkb}.  Note that the
Stokes equations may also be successfully discretized by the finite 
element method\index{Finite element method},
see e.g. \citet{giraultraviart} and references therein.

\subsection{Structured staggered grids}\label{sstokes}
\index{Staggered grid}The  discretization of the incompressible
Navier-Stokes equations with staggered grids is classical (see \citet{Pat}): the idea
is to associate different control volume grids to the different unknowns. In the two-dimensional
case, the meshes  consist  in rectangles. Consider, for instance, the mesh, say
$\mesh$, for the pressure $p$. Then, considering that the discrete 
unknowns are located at the centers of the elements of their associated 
mesh, the discrete unknowns for $p$
are, of course, located at the centers of the element of
$\mesh$. The meshes are staggered such that the discrete unknowns for the $x$-velocity
are located at the centers of the edges of $\mesh$ parallel to the $y$-axis, and the
discrete unknowns for the $y$-velocity are located at the centers of the edges of $\mesh$
parallel to the $x$-axis. 
The two equations of ``momentum" are associated to the $x$ and
$y$-velocity (and integrated over the control volumes of the considered 
mesh) and the ``divergence free" equation is associated to the pressure 
(and integrated over the control volume of $\mesh$). Then the
discretization of all the terms of the equations is straightforward, except for the
convection terms\index{Convection term} (in the momentum equations) which, eventually,  have to be discretized
according to the Reynolds number (upstream or centered discretization\ldots).
The convergence analysis of this so-called ``MAC" (Marker and Cell) is performed in 
\citet{nicol1} in the linear case
and \citet{nicol2} in the case of the Navier-Stokes equations. 

\subsection{A finite volume scheme on unstructured staggered grids\index{Staggered grid}}\label{ustokes}

Let us now turn to the case of unstructured grids; the scheme we shall study uses
the same control volumes for all the components of the velocity. The pressure  unknowns
are located at the vertices, and a Galerkin expansion \index{Galerkin expansion} is used for the approximation
of the pressure. Note that other finite volume schemes have been proposed for the
discretization of the Stokes and incompressible Navier-Stokes equations
on unstructured grids\index{Stokes equations} 
(\citet{bottahempel}), but, to our knowledge, no proof of convergence has been given yet.
 
\medskip

We again use the notion of admissible mesh, introduced in Definition
\refp{meshdirichlet}, in the particular case of triangles, if $d=2$, or tetrahedra, if
$d=3$.
We limit the description below to the case $d=2$ and to the Stokes equations. 
Let $\O$ be an open 
bounded polygonal connected subset $\O$ of
$\R^2$. Let $\mesh$ be a mesh of $\O$ consisting of triangles, satisfying the properties
required  for  the finite element method (see e.g. \citet{ciarlet}), 
with  acute angles only. Defining, for all $\cv \in \mesh$, the  point   $\centercv$  as the
intersection of the orthogonal bisectors of the sides of the triangle $\cv$  yields that
$\mesh$ is an admissible mesh in the sense of Definition \refp{meshdirichlet}.
Let $\S_\mesh$ be the 
set of vertices of $\mesh$. For $S \in \S_\mesh$, let 
$\phi_S$ be the shape function associated to $S$ in the piecewise linear  finite element method for the mesh $\mesh$. For all
$\cv\in\mesh$, let $\S_\cv\subset\S_\mesh$ be the set of the vertices of $\cv$.

A possible finite volume scheme \index{Finite volume scheme!-- for 
the Stokes system} using a Galerkin expansion \index{Galerkin expansion} for the pressure is defined
by the following equations, with the notations of Definition \refp{meshdirichlet}:

\be        
\bari
\dsp{\nu \sum_{\edge \in \edgescv} F^{(i)}_{\cv,\edge} 
 + \sum_{S\in\S_\cv} p_S \int_\cv
\dxi{\phi_S} (x) dx = } \mcv f_\cv^{(i)},
\\
\forall \cv\in\mesh, \, \forall i=1,\ldots,d, 
\ean
\label{stkschu}
\ee

\be
\bal
\dsp{F^{(i)}_{\cv,\edge}=\transedge (u_\cv^{(i)} - u_\cvv^{(i)}), 
\text{ if } \edge \in \edgesint, \, \edge=\edgecvcvv, \, i=1,\ldots,d,}
\\
\dsp{F^{(i)}_{\cv,\edge}=\transedge u_\cv^{(i)}, 
\text{ if } \edge \in \edgesext \cap \edgescv, \,  i=1,\ldots,d,}
\ean
\label{stkflu}
\ee

\be
\summesh \sumund u_\cv^{(i)} \int_\cv \dxi{\phi_S} (x)dx = 0, \, \forall S \in \S_\mesh,
\label{stkschd}\ee

\be
\int_\O \sum_{S\in\S_\mesh} p_S \phi_S(x)dx = 0, 
\label{stkschpp}\ee

\be         
f_\cv^{(i)} = \frac 1 {\mcv} \int_\cv f(x) dx, \,  \forall \cv\in\mesh.
\label{stkschf}\ee

The  discrete unknowns of  \eqref{stkschu}-\eqref{stkschf} 
are $ u^{(i)}_\cv$, $\cv \in \mesh$, $i=1,\ldots,d $ and
$ p_S$, $S \in \S_\mesh $. 

The approximate solution is defined by

\be
p_\mesh = \sum_{S\in\S_\mesh} p_S \phi_S,
\label{spapp}
\ee

\be
u_\mesh^{(i)}(x) = u_\cv^{(i)}, \, \ae x\in\cv, \, \forall \cv\in\mesh, 
\, \forall i=1,\ldots,d.
\label{suapp}
\ee

The proof of the
convergence of the scheme is not straightforward in the general case. 
We shall prove  in the following proposition the convergence of the discrete velocities 
given by the finite volume scheme
\eqref{stkschu}-\eqref{stkschf} in the simple case of a mesh consisting of 
equilateral triangles.

\begin{prop}\label{pstk}

Under Assumption \ref{stkH}, let $\mesh$ be a triangular finite element mesh of $\O$, with 
acute angles only, and let, for all $\cv \in \mesh$,  $\centercv$ be the intersection of the
orthogonal bisectors of the sides of the triangle $\cv$ (hence $\mesh$ is an
admissible mesh in the sense of Definition \refp{meshdirichlet}).
Then, there exists a unique 
solution to \eqref{stkschu}-\eqref{stkschf}, denoted by $\{u^{(i)}_\cv$, $\cv \in \mesh$, $i=1,\ldots,d\}$ and
$\{p_S$, $S \in \S_\mesh\}$.
Furthermore, if the elements of $\mesh$ are equilateral triangles, then $u_\mesh \tends u$ in
$(L^2(\O))^d$, as $\size(\mesh) \tends 0$, where  $u$ is the (unique) solution to
\eqref{stksolvar} and  $u_\mesh=(u_\mesh^{(1)},\ldots,u_\mesh^{(d)})^d$ is defined 
by \eqref{suapp}.

\end{prop}

\proof of Proposition \ref{pstk}.

\smallskip

\textit{Step 1 (estimate on $u_\mesh$)}

Let $\mesh$ be an admissible mesh, in the sense of Proposition \ref{pstk}, 
and $\{u^{(i)}_\cv$, $\cv \in \mesh$, $i=1,\ldots,d\}$, $\{p_S$, $S \in 
\S_\mesh\}$ be a solution of \eqref{stkschu}-\eqref{stkschd} with 
\eqref{stkschf}.

Multiplying the equations
\eqref{stkschu} by
$u_\cv^{(i)}$, summing over $i=1,\ldots,d$ and
$\cv\in\mesh$ and using \eqref{stkschd} yields

\be
\nu \sumund \sum_{\edge \in \edges} \transedge (D_\edge u^{(i)})^2
= \sumund \summesh \mcv u_\cv^{(i)}  f_\cv^{(i)},
\label{stkest1}\ee

with $D_\edge u^{(i)}=\vert u_\cvv^{(i)} - u_\cv^{(i)} \vert$ if $\edge \in 
\edgesint$, $\edge=\edgecvcvv$, $i \in \{1,\ldots,d\}$ and $D_\edge 
u^{(i)}=\vert u_\cv^{(i)} \vert$ if $\edge \in \edgesext \cap \edgescv$, 
$i \in \{1,\ldots,d\}$.

\smallskip

In step 2, the existence and the uniqueness of the solution of 
\eqref{stkschu}-\eqref{stkschf} will be essentially deduced from 
\eqref{stkest1}.

\smallskip

Using the discrete 
Poincar\'e inequality \eqref{ellinpoin} in \eqref{stkest1} gives an
$L^2$ estimate and an estimate on the ``discrete $H^1_0$ norm" on the 
component of the approximate 
velocities, as in Lemma \refp{ellexistu}, that is: 

\[
\norm{u^{(i)}_\mesh}{1,\mesh} \le C, \; \norm{u^{(i)}_\mesh}{L^2(\O)} \le 
C, \forall i \in \{1,\ldots,d\},
\]

where $C$ only depends on $\O$, $vu$ and $f^{(i)}$, $i=1,\ldots,d$.

As in Theorem \refp{ellcvgce} (thanks to Lemma \refp{h1dtot} and Theorem 
\refp{kolmh10}), this estimate gives the relative 
compactness in $(L^2(\O))^d$ of the set of approximate solutions 
$u_\mesh$, for $\mesh$
in the set of  admissible meshes in the sense of Proposition \ref{pstk}. It  also gives 
that if $u_{\mesh_n} \tends u$ in $(L^2(\O))^d$, as $n \tends \infty$, where 
$u_{\mesh_n}$ is the solution associated to the mesh $\mesh_n$, and 
$\size(\mesh_n) \tends 0$ as $n \tends \infty$, then $u \in 
(H^1_0(\O))^d$. This will be 
used in Step 3 in order to prove the convergence of $u_\mesh$ to the solution of 
\eqref{stksolvar}.

\smallskip

\textit{Step 2 (existence and uniqueness of $u_\mesh$ and $p_\mesh$)}

Let $\mesh$ be an admissible mesh, in the sense of Proposition \ref{pstk}.
Replace, in the right hand side of \eqref{stkschd}, ``$0$" by ``$g_S$" with some 
$\{g_S$, $S \in \S_\mesh\} \subset \R$. Eliminating 
$F^{(i)}_{\cv,\edge}$, the system 
\eqref{stkschu}-\eqref{stkschd} becomes a linear 
system with as many equations as unknowns. The sets of unknowns are $\{u^{(i)}_\cv$, $\cv \in \mesh$, $i=1,\ldots,d\}$ and
$\{p_S$, $S \in \S_\mesh\}$. Ordering the equations and the unknowns 
yields a matrix, say $A$, defining this system.

~

Let us determine the kernel of $A$; let $f_\cv^{(i)}=0$ and $g_S=0$ for all $\cv \in 
\mesh$, all $S \in \S_\mesh$ and all $i \in 
\{1,\ldots,d\}$. Then, \eqref{stkest1} leads to $u_\cv^{(i)}=0$ for all $\cv \in \mesh$ and all $i \in 
\{1,\ldots,d\}$. Turning back to \eqref{stkschu} yields that 
$p_\mesh$ (defined by \eqref{spapp}) is constant on $\cv$ for all $\cv \in \mesh$. Therefore, since 
$\O$ is connected, $p_\mesh$ is constant on $\O$. Hence, the dimension of 
the kernel of $A$ is $1$ and so is  the codimension of the range of $A$.
In order to determine the range of $A$, note that

\[
\sum_{S \in \S_\mesh} \varphi_S(x)=1, \, \forall x \in \O.
\]

Then, a necessary condition in 
order that the linear system \eqref{stkschu}-\eqref{stkschd} has a solution is 

\be
\sum_{S \in \S_\mesh} g_S=0
\label{csex}
\ee

and, since the codimension of the range of $A$ is $1$, this condition is 
also sufficient.
Therefore, under the condition \eqref{csex}, the linear system 
\eqref{stkschu}-\eqref{stkschd} has a solution, this 
solution is unique up to an additive constant for $p_\mesh$.
In the particular case $g_S=0$ for all $S \in \S_\mesh$, this yields that 
\eqref{stkschu}-\eqref{stkschf} has a unique solution.

\smallskip

\textit{Step 3 (convergence of $u_\mesh$ to $u$)}

In this step the convergence of  $u_\mesh$ towards $u$ in  $(L^2(\O))^d$  as $\sizemesh
\tends 0$ is shown for 
 meshes consisting of
equilateral triangles.
 Let $(\mesh_n)_{\nnn}$ be a sequence of meshes (such as defined in
Proposition \ref{pstk}) consisting of equilateral triangles 
and let  $(u_{\mesh_n})_{\nnn}$ be the associated solutions. Assume that 
$\size(\mesh_n) \tends 0$ and $u_{\mesh_n} \tends u$ in $(L^2(\O))^d$ as 
$n \tends \infty$. 
Thanks to the compactness result of Step 1, proving 
that $u$ is the solution of \eqref{stksolvar} is sufficient to conclude 
this step and to conclude Proposition \ref{pstk}.

By Step 1, $u \in (H^1_0(\O))^d$. It remains to show that $u \in V$ 
(which is the first part of \eqref{stksolvar}) and that $u$ satisfies the 
second part of \eqref{stksolvar}.

~

For the sake of simplicity of the notations, let us omit, from now on, the index $n$ in
$\mesh_n$ and let $h = \size(\mesh)$.   Note that $\centercv$ (which is the intersection of
the orthogonal bisectors of the sides of the triangle $\cv$) is the center of gravity of
$\cv$, for all $\cv\in\mesh$. 
Let $\varphi=(\varphi^{(1)},\ldots,\varphi^{(d)})^t \in V$ and assume that 
the functions $\varphi^{(i)}$ are regular functions with compact support 
in $\O$, 
say $\varphi^{(i)} \in C^\infty_c(\O)$ for all $i \in \{1,\ldots,d\}$. 
There exists $C >0$ only depending on $\varphi$ such that

\be
\vert \varphi^{(i)}(\centercv) - \frac 1 \mcv \int_\cv \varphi^{(i)}(x) dx\vert \le 
C h^2,
\label{stkapp}
\ee

for all $\cv\in\mesh$ and $i=1,\ldots,d$. Let us proceed as in the proof of  
convergence of  the finite volume scheme  for the  Dirichlet  problem (Theorem
\refp{ellcvgce}).

Assume that $h$ is small enough so that  $\varphi(x)=0$ for all $x$ such that $x \in \cv$,
$\cv  \in \mesh$ and 
$\edgescv \cap \edgesext \ne \emptyset$.

Note that $ (\dr \phi_S)/ 
(\dr  x_i)$ is constant in each
$\cv\in\mesh$ and that 

\[
\sumund \int_\O \frac {\dr \phi_S} 
{\dr  x_i} (x) \varphi^{(i)}(x) dx = - \int_\O \phi_S(x)
\sumund \frac  {\dr \varphi^{(i)}} 
{\dr  x_i} (x) dx = 0.
\]

Then,

\[
\bac
\sumund \summesh \sum_{S\in\S_\cv} p_S \int_\cv \frac {\dr \phi_S} 
{\dr  x_i} (x)  dx \frac 1 \mcv \int_\cv
\varphi^{(i)}(x) dx = 0.
\ean
\]

Therefore, multiplying the equations
\eqref{stkschu} by
$(1 / \mcv) \int_\cv \varphi^{(i)}(x) dx $,  for each $i=1,\ldots,d$,   
summing the results over $\cv \in \mesh$ and $i \in \{i\ldots,d\}$ yields

\be
\bac
\nu \sumund \sum_{\edgecvcvv\in\edgesint} \transcvcvv (u_\cvv^{(i)} - u_\cv^{(i)}) 
(\frac 1 \mcvv \int_\cvv \varphi^{(i)}(x) dx -  \frac 1 \mcv \int_\cv
\varphi^{(i)}(x) dx) = 
\\
 \sumund \summesh  f_\cv^{(i)} \int_\cv \varphi^{(i)}(x) dx. 
\ean
\label{stkeq1}
\ee

Passing to the limit in \eqref{stkeq1} as $n \tends \infty$ and using 
\eqref{stkapp} gives, in the same way as for
the Dirichlet  problem (see Theorem \refp{ellcvgce}), that $u$ satisfies 
the equation given in \eqref{stksolvar}, at least for 
$v \in V \cap (C^\infty_c(\O))^d$. Then, 
since  $V \cap (C^\infty_c(\O))^d$ is dense (for the 
$(H^1_0(\O))^d$-norm) in $V$ (see, for 
instance, \citet{lions} for a proof of this result), $u$ satisfies the 
equation given in \eqref{stksolvar}.

\smallskip

Since $u \in (H^1_0(\O))^d$, it remains to show that $u$ is divergence free.
Let $\varphi\in C_c^\infty(\O)$. Multiplying \eqref{stkschd} by 
$\varphi(S)$, summing over $S\in \S_\mesh$ and noting that the function $ 
\sum_{S\in\S_\mesh} \varphi(S) \phi_S $ converges to
$\varphi$ in $H^1(\O)$, one obtains that $u$ is divergence free and then 
belongs to $V$. This completes the proof that $u$ is the (unique) 
solution of \eqref{stksolvar} and concludes the proof of Proposition 
\ref{pstk}. \eof

\section{Flows in porous media}\label{flowpm}

\hpage{Systems}

\subsection{Two phase flow}\index{Two phase flow}

This section is devoted to the discretization of a system which may
be viewed as an elliptic equation coupled to a hyperbolic equation. This 
system appears in the modelling of a two phase flow in a porous medium.
Let $\O$ be an open bounded polygonal subset of $\R^\dim$, $\dim = 2$ or 
$3$, and  let $a$ and $b$ be functions of class $C^1$ from
$\R$ to $\R_+$. Assume that $a$ is nondecreasing and $b$ is nonincreasing. Let $g$
and $\overline{u}$ be bounded functions from $\partial \O \times \R_+$ to $\R$, and
$u_0$ be a bounded function from $\O$ to $\R$. Consider the following problem:

\be
\begin{array}{rcl}
u_t (x,t) - \div (a(u) \grad p)(x,t)  & = &  0, \ \
 (x,t) \in \O \times \R_+, \\
(1-u)_t (x,t) - \div (b(u) \grad p)(x,t)  & = &  0, \ \
 (x,t) \in \O \times \R_+, \\
\grad p(x,t) \cdot \bfn (x)  & = &  g(x,t), \ \  (x,t) \in \partial \O \times \R_+, \\
u(x,t)  & = &  \overline{u}(x,t), \ \  (x,t) \in \partial \O \times \R_+  \ ; \ g(x,t) \ge 0, \\
u(x,0)  & = &  u_0(x), \ \  x \in \O,
\end{array}
\label{egeqeeh2D}
\ee
where $\bfn$ is the normal to $\partial \O$, outward to $\O$.
The unknowns of this system are the functions $p$ and $u$ (from $\O \times \R_+$ to $\R$).
Adding the two first equations of (\ref{egeqeeh2D}), this system may be viewed as an elliptic equation with respect to the  unknown $p$, for a given $u$ (note that there is no time derivative in this equation), with a Neumann condition,  coupled to  a hyperbolic equation with respect to the  unknown  $u$ (for a given $p$). 
Note that, for the elliptic problem with the  Neumann condition, the compatibility condition on $g$ reads
\[
\int_{\partial \O}M(u(x,t))g(x,t)\dfrontiere(x)=0, \ \  t \in \R_+,
\]
where $M=a+b$.
It is not known whether the system (\ref{egeqeeh2D}) has a solution, except in the simple case where the function $M$ is a positive constant (which is, however, already  an interesting case for real  applications). 

~

In order to discretize \index{Finite volume scheme!-- for 
multiphase flow problems|(} (\ref{egeqeeh2D}), let $\mesh$ be an admissible mesh of $\O$ in the
sense of Definition  \refp{meshneuman} and $k>0$ be the time step. 
The discrete unknowns are $p_\cv^n$ and $u_\cv^n$ for $\cv\in\mesh$ and $n \in 
\N^{\star}$. 
The discretization of the initial condition is

\[
\ucvz = \frac 1 \mcv \int_\cv u_0(x) dx,
\ \cv\in\mesh.
\]

In order to take into account the boundary condition on $u$, define, with 
$t_n=nk$,

\[
\overline{u}_\cv^n = \frac 1 {k \, \meas(\dr\cv\cap\dr \O)} {\int_{\dr\cv\cap\dr
\O}\int_{t_n}^{t_{n+1}} \overline{u}(x,t) \dfrontiere (x) dt} ,
\ \cv\in\mesh, \ \  n \in \N. 
\]

The scheme will use an ``upstream\index{Upstream} choice" of $a(u)$ and
$b(u)$ on each ``interface" of the mesh, that is, for all $\cv\in\mesh$, 
$\cvv\in\Nv$,

\[
\begin{array}{l}
\dsp {(a(u))^n_{\cv, \cvv}} = a(\ucvn) \ \  \mathrm{if} \ \  p_\cv^{n+1} \ge
p_\cvv^{n+1}  \\
\dsp {(a(u))^n_{\cv, \cvv}} = a(\ucvvn) \ \  \mathrm{if} \ \  p_\cv^{n+1}
< p_\cvv^{n+1}, 
\\
\dsp {(b(u))^n_{\cv, \cvv}} = b(\ucvn) \ \  \mathrm{if} \ \  p_\cv^{n+1} \ge p_\cvv^{n+1} 
\\
\dsp {(b(u))^n_{\cv, \cvv}} = b(\ucvvn) \ \  \mathrm{if} \ \  p_\cv^{n+1}
< p_\cvv^{n+1}, 
\end{array} 
\]

The discrete equations are, for all $\cv\in\mesh$, $n \in \N$,

\[
\begin{array}{l}
\mcv \dsp{\frac {\ucvnp - \ucvn} k - \sumv \transcvcvv
(p_\cvv^{n+1} - p_\cv^{n+1}) (a(u))^n_{\cv, \cvv}}
\\
\dsp{-\frac {a(\overline{u}_\cv^n)}
k  \int_{\dr\cv\cap\dr\O}\int_{t_n}^{t_{n+1}}
g^+(x,t) \dfrontiere (x) dt + \frac {a(\ucvn)}
k  \int_{\dr\cv\cap\dr\O}\int_{t_n}^{t_{n+1}}
g^-(x,t) \dfrontiere (x) dt=0, }
\\
- \mcv \dsp{\frac {\ucvnp - \ucvn} k - \sumv \transcvcvv
(p_\cvv^{n+1} - p_\cv^{n+1}) (b(u))^n_{\cv, \cvv}}
\\
\dsp{-\frac {b(\overline{u}_\cv^n)}
k  \int_{\dr\cv\cap\dr\O}\int_{t_n}^{t_{n+1}}
g^+(x,t) \dfrontiere (x) dt + \frac {b(\ucvn)}
k  \int_{\dr\cv\cap\dr\O}\int_{t_n}^{t_{n+1}}
g^-(x,t) \dfrontiere (x) dt=0.}
\end{array}
\]

Recall that $g^+(x,t)=\max \{g(x,t),0\}$, $g^-=(-g)^+$ and 
$\transcvcvv=\medgecvcvv / \dcvcvv$ (see Definition \refp{meshdirichlet}).
This  finite volume scheme  \index{Finite volume scheme!-- for 
multiphase flow problems|)} gives very good numerical results under a 
usual stability  condition on the time step with respect to the space mesh. 
It can be generalized to more complicated systems (in
particular, for the simulation of multiphase flows in porous medium such as the ``black
oil" case of reservoir engineering, see \citet{E}). It is possible to 
prove the convergence of this scheme in the case where the function $M$ 
is constant and the function $g$ does not depend on $t$. 
In this case, the scheme may be written as a finite volume scheme for a 
stationary diffusion equation
with respect to the unknown $p$ (which does not depend on $t$) and an upstream\index{Upstream} finite
volume scheme for a hyperbolic equation with respect to the unknown $u$. 
The proof of this convergence is given below 
(Theorem \ref{egcvgce}) under the 
assumptions that  $a(u)=u$ and  $b(u)=1-u$ (see also \citet{Vignal1}). Note that the elliptic
equation with respect to the pressure may also be discretized with a
finite element method\index{Finite element method},
and coupled to the finite volume scheme for the hyperbolic equation. 
This coupling of finite elements and finite volumes was
introduced in \citet{F2}, where it is called ``CVFE" (Control Volume 
Finite Element\index{Control 
volume finite element method}), in \citet{ES} and in \citet{EG}, where the
convergence of the finite element-finite volume scheme is shown
under the same assumptions.

\subsection{Compositional multiphase flow}\label{thermo}

Let us now turn to the study of a system of partial differential equations which arises 
in the simulation of a multiphase flow in a porous medium (the
so called ``Black  Oil" case in petroleum engineering, see e.g. \citet{E}). This system consists in
a parabolic equation coupled with hyperbolic equations 
and algebraic equations and inequalities
(these  algebraic equations and inequalities are given by an assumption of thermodynamical
equilibrium). It may be written, for  $x \in \O$ and  $t\in \R_+$,  as:

\be
\frac {\partial } {\partial t} (\rho_1(p) u)(x,t) - 
\div (f_1(u,v,c)\grad p)(x,t)  = 
0,  
\label{triph1}
\ee
\be
\frac {\partial } {\partial t} (\rho_2(p,c) (1-u-v)(1-c))(x,t) - \div (f_2(u,v,c)\grad
p)(x,t)  =   0,
\label{triph2}
\ee
\be
\frac {\partial } {\partial t} (\rho_2(p,c) (1-u-v)c + \rho_3(p)v)(x,t) -
\div (f_3(u,v,c)\grad p)(x,t)  =  0,
\label{triph3}
\ee
\be
(v(x,t)=0 \text{ and } c(x,t) \le f(p(x,t)) \text{ or } (c(x,t) = f(p(x,t))
\text{ and } v(x,t)  \ge  0),
\label{triph4}
\ee

where $\O$ is a given open bounded polygonal subset of $\R^\dim$ 
($\dim=2$ or $3$), $f_1$, $f_2$, $f_3$ are given functions from
$\R^3$ to $\R_+$, $f$, $\rho_1$, $\rho_3$ are given functions from $\R$ to
$\R_+$ and  $\rho_2$ is a given function from $\R^2$ to $\R_+$. The 
problem is completed by
initial and boundary conditions which are omitted  here. The unknowns of
this problem are the functions $u,v,c,p$ from $\O \times \R_+$ to $\R$.

\medskip

In order to discretize this problem, let $k$ be the time step (as usual, $k$ may in
fact be variable) and $\mesh$ be a cartesian mesh of $\O$. Following the
ideas (and notations) of the previous chapters, the discrete unknowns 
are $u_\cv^n$,
$v_\cv^n$, $c_\cv^n$ and $p_\cv^n$, for $\cv \in \mesh$ and $n \in 
\N^{\star}$ and it is 
quite
easy to discretize (\ref{triph1})-(\ref{triph3}) with a classical finite
volume method.  Note that the time
discretization of the unknown $p$  must generally be implicit while the
time discretization of the unknowns  $u,v,c$ may be explicit or implicit. The explicit choice
requires a usual restriction on the time step (linearly with respect to the space step). The
only new problem is the discretization of (\ref{triph4}), which is now described.

\smallskip

Let $n \in \N$. The discrete unknowns at time $t_{n+1}$, namely $u_\cv^{n+1}$,
$v_\cv^{n+1}$, $c_\cv^{n+1}$ and $p_\cv^{n+1}$, $\cv \in \mesh$, have to be computed 
from the
discrete unknowns at time $t_{n}$, namely $u_\cv^{n}$,
$v_\cv^{n}$, $c_\cv^{n}$ and $p_\cv^{n}$, $\cv \in \mesh$. Even if the time discretization
of (\ref{triph1})-(\ref{triph3}) is explicit with respect to the unknowns $u,v$ and $c$, the
system of discrete equations (with unknowns  $u_\cv^{n+1}$, $v_\cv^{n+1}$, 
$c_\cv^{n+1}$ and
$p_\cv^{n+1}$, $\cv \in \mesh$) is nonlinear, whatever the discretization of
(\ref{triph4}). It can be solved by, say, a Newton process\index{Newton's algorithm}. Let $l \in \N$ be the index of the
 ``Newton iteration", and $u_\cv^{n+1,l}$, $v_\cv^{n+1,l}$, $c_\cv^{n+1,l}$ 
and
$p_\cv^{n+1,l}$ ($\cv \in \mesh$) be the computed unknowns at iteration $l$. As usual,
these unknowns are, for $l=0$, taken equal to  $u_\cv^{n}$, $v_\cv^{n}$, $c_\cv^{n}$ 
and
$p_\cv^{n}$. In order to discretize (\ref{triph4}), a
``phase index" is introduced; it is denoted by $i_\cv^n$, for all $\cv \in \mesh$ and $n \in
\N$ and it is defined by:

\[
\begin{array}{l}
\text{ if } i_\cv^n=0 \text{ then } v_\cv^n=0  \ \ (\text{ and } c_\cv^n \le f(p_\cv^n)),
\\
\text{ if } i_\cv^n=1 \text{ then } c_\cv^n = f(p_\cv^n) \ \ (\text{ and } v_\cv^n \ge 0).
\end{array}
\]

In the Newton process for the computation of the unknowns at time $t_{n+1}$, a 
``phase index", denoted by $i_\cv^{n+1,l}$ is also introduced, with $i_\cv^{n+1,0}=i_\cv^{n}$. 
This phase index is used in the computation of $u_\cv^{n+1,l+1}$,
$v_\cv^{n+1,l+1}$, $c_\cv^{n+1,l+1}$, $p_\cv^{n+1,l+1}$ and $i_\cv^{n+1,l+1}$ ($\cv 
\in
\mesh$), starting from $u_\cv^{n+1,l}$, $v_\cv^{n+1,l}$, $c_\cv^{n+1,l}$, $p_\cv^{n+1,l}$ and
$i_\cv^{n+1,l}$. Setting $v_\cv^{n+1,l+1}=0$ if $i_\cv^{n+1,l}=0$, and
$c_\cv^{n+1,l+1}=f(p_\cv^{n+1,l+1})$ if $i_\cv^{n+1,l}=1$, the computation of 
(intermediate)
values of $u_\cv^{n+1,l+1}$, $v_\cv^{n+1,l+1}$, $c_\cv^{n+1,l+1}$, 
$p_\cv^{n+1,l+1}$ is possible
with a ``Newton iteration" on (\ref{triph1}), (\ref{triph2}), (\ref{triph3}) (note that the
number of unknowns is equal to the number of equations). Then, for each $\cv \in 
\mesh$,
three cases are possible:

\benum
\item if $c_\cv^{n+1,l+1} \le f(p_\cv^{n+1,l+1})$ and $v_\cv^{n+1,l+1} \ge 0$, then 
set $i_\cv^{n+1,l+1}=i_\cv^{n+1,l}$,
\item if $c_\cv^{n+1,l+1} \,  > \,  f(p_\cv^{n+1,l+1})$ (and necessarily 
$i_\cv^{n+1,l} \  =0$), then  
set $c_\cv^{n+1,l+1} \,  = \, f(p_\cv^{n+1,l+1})$ and $i_\cv^{n+1,l+1}=1$,
\item if $v_\cv^{n+1,l+1} < 0$ (and necessarily $i_\cv^{n+1,l}=1$), then set
$v_\cv^{n+1,l+1} = 0$ and $i_\cv^{n+1,l+1}=0$.
\eenum

This yields the final values of $u_\cv^{n+1,l+1}$, $v_\cv^{n+1,l+1}$, 
$c_\cv^{n+1,l+1}$,
$p_\cv^{n+1,l+1}$ and $i_\cv^{n+1,l+1}$ ($\cv \in \mesh$).

\medskip

When the ``convergence" of the Newton process \index{Newton's algorithm} is achieved, say at 
iteration $l^\star$, the
values of the unknowns at time $t_{n+1}$ are found. They are taken equal to those
indexed by $(n+1,l^\star)$ (for $u,v,c,p,i$).  
It can be proved, under convenient hypotheses on the function $f$ (which are realistic in
the applications), that there is no ``oscillation" of the ``phase index" during the
Newton iterations \index{Newton's algorithm} performed from time $t_n$ to time $t_{n+1}$ (see 
\citet{EG1}). 
This method, using the phase index, was also successfully adapted for the treatment of the obstacle problem  and the Signorini problem, see \citet{HM,HM-imajna}.

\subsection{A simplified case}\label{casimple}

The aim of this section and of the following sections is the study of the convergence of 
two coupled finite volume schemes, for the system of equations
 $u_{t} - \div (u \grad p) = 0$ and $\Delta p = 0$, defined on an open set
$\O$. A finite volume mesh $\mesh$ is used for the discretization in space, 
together  with an explicit Euler time discretization. Similar results are 
in \citet{Vignal1}
and \citet{vigver} where the case of different space meshes for the two
equations is also studied.

\separ

We assume that the following assumption is satisfied.

\begin{hyp}\label{eghyp1}
Let $\O$ be an open polygonal bounded connected subset of $\R^d$, $d=2$ or $3$, and $\dr\O $ its boundary.
We denote by  $\bfn$ the normal vector to $\dr\O $ outward to $\O $.

Let $g\in L^{2}(\dr\O)$ be a function such that 
\[\dsp{\int_{\dr\O}} \ g(x)\dfrontiere (x)  = 0,\]
and let
$\dr\O^+=$$ \{x \in \dr\O$, $g(x)\ge 0\}$,
$\O^+ =$$ \O \cup \dr\O ^+$ and 
$\dr\O ^-= $ $\{x \in \dr\O$, $g(x)\le 0\}$.
Let $u_0\in L^{\infty }(\O )$ and
$\bar{u}\in L^{\infty }(\dr\O ^+\times\Rpast)$ 
represent respectively  the initial condition and the boundary condition
for the unknown $u$.
\end{hyp}

\separ

The set
\[
\mathcal{D}(\O ^+\times\Rp) = 
\{\varphi \in C^{\infty }_{c}(\R^\dim\times\R,\R), \ 
\varphi  = 0 \text{ on } \dr\O ^-\times\Rp\} 
\]
will be the set of test functions for Equation \eqref{egeq1} in
the weak formulation of the problem, which is given below.

\begin{defi} \normalfont
A pair $(u,p) \in  L^{\infty }(\O\times\Rpast)\times H^{1}(\O )$
($u$ is the saturation, $p$ is the pressure) is a weak solution of

\be
\left\{\begin{array}{ll}
\Delta p(x) = 0, & \forall x\in \O,  \\
\grad p(x)\cdot\bfn(x) = g(x), & \forall x \in \dr\O, \\
u_{t}(x,t) -\div(u\grad p) (x,t) = 0, & \forall x\in \O ,  
\forall t\in \Rp,   \\
u(x,0) = u_0(x), & \forall x \in \O, \\
u(x ,t) = \bar{u}(x ,t), & \forall x \in \dr\O ^+, 
\forall t \in \Rp.
\end{array}\right. \label{upfort}
\ee
if it verifies
\be
p \in  H^{1}(\O ), \label{egeq02}
\ee

\be
u \in  L^{\infty }(\O\times\Rpast), \label{egeq01}
\ee

\be
\int_\O \grad p(x)\cdot\grad X(x) dx -
\int_{\dr\O } X(x) g (x) \dfrontiere(x)  = 0, \forall X \in H^{1}(\O).
\label{egeq2}
\ee

and
\be
\begin{array} {l}
\dsp{ \int_{\Rp}\int_\O u(x,t) 
(\varphi_t(x,t) - 
\grad p(x)\cdot\grad \varphi (x,t) dx dt }
+ \dsp{ \int_\O u_{0}(x) \varphi (x,0) dx +} \\  
\dsp{\int_{\Rp}\int_{{\dr\O }^+} 
\bar{u}(x,t) \varphi (x,t) g(x)\dfrontiere (x) dt = 0, \forall \varphi 
\in \mathcal{D}(\O ^+\times\Rp)}.
\end{array}
\label{egeq1}
\ee
\end{defi}

Under Assumption \ref{eghyp1}, a classical result gives the  existence of $p\in H^{1}(\O )$ 
and the uniqueness of  $\grad p$ where $p$ is the 
solution of \eqref{egeq02},\eqref{egeq2}, which is a variational formulation of
the classical Neumann problem. 
Additional hypotheses on the function $g$
are necessary to get the uniqueness of $u\in L^{\infty }(\R^\dim\times\Rpast)$ solution of 
\eqref{egeq1}. 
The existence of $u$ results from the convergence of the scheme, but not its uniqueness,
which could be obtained thanks to regularity properties of  $\grad p$. We shall
assume such regularity, which   ensures  the uniqueness of
the function $u$ and allows an error estimate between  the finite volume scheme 
approximation of the pressure and the exact pressure. In fact, for the sake of 
simplicity, we assume (in Assumption \ref{eghyp4} below) that $p \in 
C^2({\overline \O})$. This is a rather ``strong" assumption which can be 
weakened. However, a convergence result (such as in Theorem \ref{egcvgce}) with the only assumption $p \in 
H^1(\O)$ seems not easy to obtain. Note also that similar results of 
convergence (for the ``pressure scheme" and for the ``saturation scheme") 
are possible with an open bounded connected subset of $\R^d$ with a $C^2$ 
boundary (instead of an open bounded connected polygonal subset of 
$\R^d$) using Definition \refp{parnlHT} of admissible 
meshes.

\begin{hyp}\label{eghyp4}
The pressure $p$, weak solution in $H^1(\O)$ to \eqref{egeq2}, belongs to
$C^2({\overline \O})$. \end{hyp}

\begin{remark} \normalfont
The solution $(u,p)$ of \eqref{egeq02}-\eqref{egeq1}  is also a
weak solution of

\[
(1 - u)_t(x,t) -\div ((1 - u)\grad p) (x,t)= 0.
\]

\end{remark}

\begin{remark} \normalfont
The finite volume scheme will ensure the conservation of
each of the quantities $u$ and $1-u$. It can be extended to
more complex phenomena such as compressibility, thermodynamic
equilibrium\ldots(see Section \ref{thermo})
\end{remark}

\begin{remark} \normalfont
The proof which is given here can easily be extended to the case of the existence of a source term which reads
\[
\begin{array} {ll}
- \Delta p(x)=v(x), & x \in \O,
\\
\grad p(x)\cdot\bfn(x) = g(x), &  x \in \dr\O,
\\
u_{t}(x,t) -\div (u \grad p) (x,t) + u(x,t) v^-(x) = s(x,t) v^+(x), &  x \in \O , \,  t \in \Rp,
\\
u(x,0) = u_0(x), &  x \in \O,
\\
u(x,t) = \bar{u}(x,t), &  x \in \dr\O ^+, 
 \, t \in \Rp,
\end{array}
\]
where $v \in  L^{2}(\O )$ with $\dsp{\int_{\dr\O }}g(x)\dfrontiere(x) + \dsp{\int_{\O }}v(x)dx = 0$  and $s \in  L^{\infty }(\O \times\Rpast)$. 
All modifications  which are connected to such terms will be stated in remarks.
\end{remark}

\subsection {The scheme for the simplified case}
Let $\O$ be an open polygonal bounded connected subset of $\R^\dim$.
Let $\Mesh$ be an admissible mesh, in the sense of Definition \refp{meshneuman}, and let $h=\size(\mesh)$. 
Assume furthermore that there exists $\a > 0$ such that $\dedge \ge \a h$ for all $\edge \in \edgesint$.

\subsubsection{The pressure finite volume scheme} \index{Finite volume 
scheme!-- for 
multiphase flow problems|(}

We first define the approximate pressure, 
using the finite volume scheme defined in section \refp{nll} (that is 
\eqref{nllcondbor}-\eqref{nllschemaO}).

\smallskip

(i) The values $G_\cv $, for $\cv \in\mesh$, are defined by

\be\bal     
\dsp{G_\cv  =  \int_{\dr\cv \cap\dr\O} g(x) \dfrontiere(x) \text{ if } 
\mcvbord \not = 0},
\\
G_\cv= 0,  \text{ if }  \mcvbord = 0.
\ean
\label{egschemapb}\ee

\separ

(ii) The scheme is defined by

\be        
- \sumv\transcvcvv  \Bigl(\pcvv - \pcv\Bigr) =  G_\cv , \
\forall \cv \in\mesh,  
\label{egschemap}\ee 

and 

\be        
\summesh  \mcv \pcv = 0.
\label{egschemap0}\ee 

We recall that, from lemma \refp{nllexistu}, there exists a  unique function
$p_\mesh \in X(\mesh)$ defined by $p_\mesh(x)=\pcv$ for a.e. $x \in 
\cv$, for all $\cv \in \mesh$, where  $(\pcv)_{\cv \in \mesh}$ satisfy equations
\eqref{egschemapb}-\eqref{egschemap0}. Then, using
Theorem \refp{nllesterr}, there exist $C_1$ and $C_2$, 
only depending on $p$ and $\O$, such that

\be
\Vert p_\mesh - p \Vert_{L^2(\O)} \le C_1 h
\label{egllesterr}\ee

and

\be
\sumE \medgecvcvv  \dcvcvv \big(\frac {\pcvv - \pcv} {\dcvcvv} -
\frac 1 {\medgecvcvv} \int_{\edgecvcvv} \grad p (x) \cdot \ncvcvv 
\dfrontiere (x) \big)^2  
\le (C_2 h)^2.
\label{egllesterr1}\ee

Last but not least, using   lemma \refp{nllest}, there exists $C_3$, 
only depending on $g$ and $\O$, such that

\be
\sumE\transcvcvv  ( \pcvv - \pcv )^2 \le (C_3)^2.
\label{egllesterr2}\ee

\subsubsection{The saturation finite volume scheme}
Let us now turn to the finite volume discretization of the hyperbolic 
equation
\eqref{egeq1}. 
In order to write the scheme, let us  introduce the following notations: 
let
\[G^{(+)}_\cv  =\int_{\dr\cv \cap\dr\O } g^+(x)\dfrontiere(x) \text{ and } G^{(-)}_\cv  =\int_{\dr\cv \cap\dr\O } g^-(x)\dfrontiere(x),
\]
so that  $G^{(+)}_\cv  - G^{(-)}_\cv  = G_\cv$. 
Let
\[G^{(+)} =\dsp{\int_{\dr\O }}g^+(x)\dfrontiere(x)  = \summesh G^{(+)}_\cv 
 \]
(note that $G^{(+)}$ does not depend on $\mesh$). 
The scheme \eqref{egschemap} may also be written

\be
\sumv\transcvcvv  \Bigl(\pcvv - \pcv\Bigr) + G^{(+)}_\cv  - G^{(-)}_\cv  = 0, \
\forall \cv \in\mesh.
\label{egeq3}
\ee

\begin{remark} \normalfont
In the case of the problem with source terms, the right hand side of the equation
\eqref{egschemap} is replaced by $G_\cv + V^{(+)}_\cv - V^{(-)}_\cv$ with
\[
V^{(\pm)}_\cv  = \int_{\cv} v^{\pm}(x)dx.
\]
Then, in the equation \eqref{egeq3} the quantities $G^{(\pm)}_\cv$ are replaced by  $G^{(\pm)}_\cv +V^{(\pm)}_\cv$.
\end{remark}

Let $\xi \in (0,1)$. 
Given an  admissible mesh $\mesh$,   the time step is defined by a real value $k>0$ such that

\be
k \le \inf_{\cv \in \mesh } \ \{ \frac{ \mcv \ (1-\xi)} {\dsp{\sum_{\cvv  \in \Nv}}\transcvcvv  (\pcvv - \pcv)^+ + G^{(+)}_\cv }.
\label{egeqS}
\ee

\begin{remark} \normalfont
Since the right hand side of \eqref{egeqS} has a strictly positive lower
bound, it is always possible to find values $k>0$ which satisfy \eqref{egeqS}.
Roughly speaking, the condition \eqref{egeqS} is a linear condition between the time step and
the size of the mesh. Let us explain  this point in more detail:  in most practical
cases, function $g$ is regular enough so that $\vert \pcvv-\pcv \vert / \dcvcvv$ is bounded
by some $C$ only depending on $g$ and  $\O$. Assume furthermore that the mesh $\mesh$ is
admissible in the sense of Definition \refp{meshneuman} and that, for some $\a >0$,
$\dcvedge \ge \a h$, for all $\cv \in \mesh$ and $\edge \in \edges$. 
Then the condition $k \le Dh$, with
$D=((1-\xi) \a) / (d(C+\norm{g}{L^\infty(\dr \O)}))$, implies the condition
\eqref{egeqS}.
Note also that for all  $g \in L^{2}(\dr \O )$ we already have a bound for 
$\snorm{p_\mesh}{1,\Mesh}$ (but this does not yield a bound on $\vert 
\pcvv-\pcv \vert / \dcvcvv$). Finally, note that condition \eqref{egeqS} 
is easy to implement in practise, since the values $\transcvcvv$ and 
$\pcv$ are available by the pressure scheme.
\end{remark}

\begin{remark} \normalfont
In the problem with source terms, the condition \eqref{egeqS} will
be modified as follows:
\[
k \le \inf_{\cv \in \mesh } \ \ \frac{\mcv \ (1-\xi)}{\dsp{\sum_{\cvv  \in \Nv}}\transcvcvv  (\pcvv - \pcv)^+ + G^{(+)}_\cv  + V^{(+)}_\cv }.
\]
\end{remark}

The initial condition is discretized by:  
 
\be
\ucvz = \dsp{\frac 1 {\mcv} \int_\cv } u_0(x)dx, \ \forall \cv  \in\mesh.
\label{egschema0}
\ee

We extend the definition of $\bar{u}$ by  $0$ on $\dr\O^-\times\Rp $, and we
define $\bar{u}^n_\cv $, for $\cv \in \mesh$ and $\nnn$, by

\be
\bal
\bar{u}^n_\cv  = \dsp{\frac 1 {k \, \meas(\dr\cv \cap\dr\O)}} \int^{(n+1)k}_{nk}
\int_{\dr\cv  \cap\dr\O} \bar{u}(x,t)\dfrontiere (x) dt, \text{ if } \meas(\dr\cv
\cap\dr\O)\neq 0,
\\
\bar{u}^n_\cv  = 0, \text{ if } \meas(\dr\cv \cap\dr\O)=0.
\ean
\label{egschemab}\ee

Hence the following function may be defined on $\dr\O \times \Rp$:

\[
\bar{u}_{\mesh,k}(x,t) = \bar{u}^n_\cv, \, \forall x\in\dr\cv \cap\dr\O, \forall \cv \in
\mesh, \forall t \in [nk,(n+1)k), \, \nnn.
\]

\separ

The finite volume discretization of the hyperbolic equation \eqref{egeq1}  
is then written as the following relation between  $u^{n+1}_\cv $  
and  all $u^n_\cvv $, $\cvv  \in \mesh $.

\be
\mcv ( u^{n+1}_\cv  - u^n_\cv  ) - 
k \ \ \Bigl[ \sumv \transcvcvv  u^n_{\cv, \cvv } (\pcvv - \pcv) + 
\bar{u}^n_\cv  G^{(+)}_\cv  - u^n_\cv  G^{(-)}_\cv  \Bigr] = 0, \ \forall \cv\in\mesh, \ \forall 
n\in\N,
\label{egschema}\ee

in which   the upstream\index{Upstream} value $u^n_{\cv, \cvv }$ is defined  by

\be\bal
u^n_{\cv, \cvv } = u^n_\cv, \text{ if } \pcv \ge \pcvv,
\\
u^n_{\cv, \cvv } = u^n_\cvv, \text{ if } \pcvv > \pcv.
\ean\label{egschemau}\ee

\separ

The approximate solution, denoted by  $\uapp$, is defined 
a.e. from $\O \times\Rp\rightarrow$ to $\R$ by

\be
\uapp (x,t) = \ucvn, \, \forall x \in \cv, \, \forall \cv \in \mesh, \, \forall t\in
[nk,(n+1)k), \, \forall n \in \N.
\label{eguapp}\ee

\begin{remark} \normalfont
In the case of source terms, the following term is defined:
\[
s^n_\cv  =  \frac 1 {\mcv k} \int^{(n+1)k}_{nk}\int_{\cv} s(x,t) dx dt
\]
and the term $k(s^n_\cv  V^{(+)}_\cv  - u^n_\cv  V^{(-)}_\cv)$ is added to  the right hand side of \eqref{egschema} \index{Finite volume scheme!-- for multiphase flow problems|)}.
\end{remark}

\subsection{Estimates on the approximate solution}

\subsubsection{Estimate in $L^\infty(\O\times\Rpast)$}

\index{Estimate on the approximate solution!-- for a hyperbolic equation}
\begin{lemma}\label{egestlinf}
Under the assumptions \ref{eghyp1} and \ref{eghyp4}, let $\mesh$ be an admissible mesh in the
sense of Definition \refp{meshneuman} and $k>0$ satisfying \eqref{egeqS}. Then, the
function  $\uapp$ defined by \eqref{egschemapb}-\eqref{egschemap0} and
\eqref{egschema0}-\eqref{eguapp} satisfies 

\be
\Vert \uapp\Vert_{L^\infty(\O\times\Rpast)} \le \max \{\Vert u_0 \Vert_{L^\infty(\O)},
\Vert \bar{u} \Vert_{L^\infty(\dr\O^+\times\Rpast)}\}.
\label{egeqestlinf}\ee

\end{lemma}

\proof of Lemma \ref{egestlinf}

\smallskip

Relation \eqref{egschema} can be written as

\[
\begin{array}{ll}
u^{n+1}_\cv  = & u^n_\cv  \dsp{\big[ 1 -  \frac k {\mcv} \big( \sum_{\cvv \in \Nv}}\transcvcvv  (\pcv - \pcvv)^- + 
G^{(-)}_\cv \big) \big] + \\ & \dsp{\frac k  {\mcv} \big(\sum_{\cvv \in \Nv}}\transcvcvv  u^n_\cvv  
(\pcvv - \pcv)^+ + G^{(+)}_\cv  \bar{u}^n_\cv  \big).
\end{array}
\]

Using 
\[\sum_{\cvv \in \Nv} \transcvcvv (\pcvv - \pcv)^+ + G^{(+)}_\cv  = \sum_{\cvv \in \Nv}\transcvcvv  (\pcv - \pcvv)^- + G^{(-)}_\cv ,
\]
and Inequality \eqref{egeqS},  the term  $u^{n+1}_\cv $  may be expressed as a linear combination of terms $u^n_\cvv $, $\cvv \in \mesh$, and $\bar{u}^n_\cv $, with positive coefficients. 
Thanks to relation \eqref{egeq3}, the sum of these coefficients is equal to $1$. The estimate  \eqref{egeqestlinf} 
follows by an easy induction. \eof

\begin{remark} \normalfont
In the case of source terms, Lemma \ref{egestlinf} remains true with the following estimate  instead of \eqref{egeqestlinf}:
\[
\Vert \uapp \Vert_{L^\infty(\O\times\Rpast)} \le \max\{\Vert u_0 \Vert_{L^\infty(\O)},
\Vert \bar{u} \Vert_{L^\infty(\dr\O^+\times\Rpast)}, \Vert s \Vert_{L^{\infty }(\O \times
\Rp^\ast)}\}.
\]

\end{remark}

\subsubsection{Weak $BV$ estimate}
\index{BV!Weak -- estimate}
\begin{lemma}\label{egestbvw}
Under the assumptions \ref{eghyp1} and \ref{eghyp4}, let $\mesh$ be an admissible mesh in the
sense of Definition \refp{meshneuman}. Let $h=\sizemesh$ and $\a > 0$ be such that $\dedge \ge \a
h$ for all $\edge \in \edgesint$. Let $k>0$ satisfying \eqref{egeqS}. Let $\{\ucvn$, $\cv
\in \mesh$, $\nnn\}$ be the solution to \eqref{egschema0}-\eqref{egschemau} with $\{\pcv$, $\cv
\in \mesh\}$ given by \eqref{egschemapb}-\eqref{egschemap0}. Let $T > k$ be a given real
value, and let  $\NTk$ be  the integer value such that $\NTk k <  T \le (\NTk+1) k$.  Then
there exists $H$, which only depends on $T$, $\O$, $u_0$, $\bar{u}$,  $g$, $\a$ and $\xi$,
such that the following inequality holds:

\be
k\sumTsk  \sumE  
\transcvcvv \vert \pcv - \pcvv \vert \vert u^n_\cv  - u^n_\cvv \vert  
+ k\sumTsk  \summesh 
G^{(+)}_\cv  \vert u^n_\cv  - \bar{u}^n_\cv \vert  \le  \frac {H} {\sqrt{h}}.
\label{egeqbvw}
\ee

\end{lemma}

\proof of Lemma \ref{egestbvw}

\smallskip

For $n\in {\N}$ and $\cv \in \mesh $,  multiplying \eqref{egschema} by 
$u^n_\cv $ yields
\be
\mcv (u^{n+1}_\cv  u^n_\cv  - u^n_\cv  u^n_\cv ) - 
k ( \sumv \transcvcvv  u^n_{\cv, \cvv }u^n_\cv  (\pcvv - \pcv) 
+ \bar{u}^n_\cv u^n_\cv G^{(+)}_\cv  - (u^n_\cv)^{2} G^{(-)}_\cv  ) = 0.
\label{egeq9}\ee

Writing
$u^{n+1}_\cv  u^n_\cv  - u^n_\cv  u^n_\cv  = 
- \half  (u^{n+1}_\cv  - u^n_\cv )^{2} 
- \half  (u^n_\cv )^{2} + \half  (u^{n+1}_\cv )^{2}$  and  summing 
\eqref{egeq9} on $\cv \in \mesh$ and $n\in\NN$ gives

\be
\begin{array}{ll}
\dsp{- \half  \sumTsk  \summesh \mcv 
(u^{n+1}_\cv  - u^n_\cv )^{2} + 
\half  \summesh \mcv
( (u^{\NTk+1}_\cv )^{2}-(u^{0}_\cv )^{2})} &
\\
- k\sumTsk  \summesh 
( \sumv \transcvcvv  u^n_{\cv, \cvv } u^n_\cv  
(\pcvv  - \pcv) + \bar{u}^n_\cv  u^n_\cv  G^{(+)}_\cv  
- (u^{n}_\cv)^2 G^{(-)}_\cv  ) & = 0.
\end{array}
\label{egeq10}\ee

Using \eqref{egschemau} gives, for all $\cv \in \mesh$,

\[
-\sumv \transcvcvv  u^n_{\cv, \cvv } u^n_\cv  
(\pcvv  - \pcv)=\sumv \transcvcvv  (u^n_\cv)^2  
(\pcv  - \pcvv)^+ -\sumv \transcvcvv  u^n_{\cvv } u^n_\cv  
(\pcvv  - \pcv)^+.
\]

Then,

\[
- \summesh \sumv \transcvcvv  u^n_{\cv, \cvv } u^n_\cv  
(\pcvv  - \pcv)=\summesh \sumv \transcvcvv ((u^n_\cv)^2 - u^n_{\cvv } u^n_\cv)
(\pcv  - \pcvv)^+.
\]

Therefore, since $(u^n_\cv )^{2}- u^n_\cv  u^n_\cvv  = 
\half (u^n_\cv  - u^n_\cvv )^{2} + 
\half ( (u^n_\cv )^{2} - (u^n_\cvv )^{2})$,

\[
\begin{array}{rl}
- \summesh \sumv \transcvcvv  u^n_{\cv, \cvv } u^n_\cv  
(\pcvv  - \pcv)= &  \half \summesh \sumv \transcvcvv (u^n_\cv-u^n_{\cvv })^2
(\pcv  - \pcvv)^+
\\
 & +\half \summesh \sumv \transcvcvv (u^n_\cv)^2
(\pcv  - \pcvv)^+
\\
 & -\half \summesh \sumv \transcvcvv (u^n_\cvv)^2
(\pcv  - \pcvv)^+
\\
= &  \half \summesh \sumv \transcvcvv (u^n_\cv-u^n_{\cvv })^2
(\pcv  - \pcvv)^+
\\
 & + \half \summesh \sumv \transcvcvv (u^n_\cv)^2
(\pcv  - \pcvv)
\end{array}
\]

and, using \eqref{egeq3},

\[
\begin{array}{rl}
- \summesh \sumv \transcvcvv  u^n_{\cv, \cvv } u^n_\cv  
(\pcvv  - \pcv)= &  \half \summesh \sumv \transcvcvv (u^n_\cv-u^n_{\cvv })^2
(\pcv  - \pcvv)^+
\\
 & +\half \summesh G_\cv^{(+)} (u^n_\cv)^2
-\half \summesh G_\cv^{(-)} (u^n_\cv)^2.
\end{array}
\]

Hence

\be
\begin{array}{ll}
- k \sumTsk  \summesh  
(\sumv  \transcvcvv  u^n_{\cv, \cvv }u^n_\cv  (\pcvv - \pcv) + 
\bar{u}^n_\cv u^n_\cv G^{(+)}_\cv  - (u^{n}_\cv)^2 G^{(-)}_\cv  ) & = \\
\half  k\sumTsk  
(\sumE \transcvcvv  \vert \pcv - \pcvv \vert (u^n_\cv  - u^n_\cvv )^{2} + 
\summesh  G^{(+)}_\cv (u^n_\cv  - \bar{u}^n_\cv )^{2})  & - \\
\half  k\sumTsk  \summesh  
( G^{(+)}_\cv  (\bar{u}^{n}_\cv)^2 - G^{(-)}_\cv  (u^{n}_\cv)^2  ).
\end{array}
\label{egeq11}
\ee
Using \eqref{egschema}, we get
\[
\sumTsk  \summesh \mcv (u^{n+1}_\cv  - u^n_\cv )^{2} = \sumTsk  \summesh \frac{k^{2}}{\mcv}\big( \sumv \transcvcvv  u^n_{\cv, \cvv } (\pcvv - \pcv) +  \bar{u}^n_\cv G^{(+)}_\cv  - u^n_\cv G^{(-)}_\cv  \big)^{2}.
\]
Then, for all $\cv \in \mesh $, using again \eqref{egeq3} and the definition \eqref{egschemau},
\[
\bal
\sumTsk  \summesh \mcv (u^{n+1}_\cv  - u^n_\cv )^{2}= \\ 
\sumTsk  \summesh  \frac{k^{2}}{\mcv} \big(\dsp{\sumv} \transcvcvv (\ucvvn-\ucvn ) (\pcvv - \pcv)^+ + G^{(+)}_\cv  (\bar{u}^n_\cv -u^n_\cv ) \big)^{2}.
\ean 
\]
The Cauchy-Schwarz inequality yields
\[
\begin{array}{ll}
\sumTsk \summesh \mcv (u^{n+1}_\cv  - u^n_\cv )^{2} \le \sumTsk \summesh  \frac{k^{2}}{\mcv} \Bigl( \dsp{\sumv}\transcvcvv  (\pcvv - \pcv)^+ + G^{(+)}_\cv  \Bigr) \\
\Bigl( \dsp{\sumv}\transcvcvv   (\pcvv - \pcv)^+ (u^n_\cvv -u^n_\cv )^{2} + G^{(+)}_\cv  (\bar{u}^n_\cv -u^n_\cv ) ^{2}\Bigr). &
\end{array}
\]
Using the stability  condition \eqref{egeqS} and reordering the summations gives
\be
\begin{array}{ll}
\sumTsk \summesh \mcv 
(u^{n+1}_\cv  - u^n_\cv )^{2} \le 
\sumTsk k \, (1 - \xi) \\ 
\Bigl( \sumE\transcvcvv   \vert \pcvv - \pcv \vert (u^n_\cvv -u^n_\cv )^{2} 
+ \summesh G^{(+)}_\cv  (\bar{u}^n_\cv -u^n_\cv ) ^{2}\Bigr). &
\end{array}
\label{egeq12}\ee

Using \eqref{egeq10}, \eqref{egeq11} and \eqref{egeq12}, we obtain

\be
\begin{array}{ll}
\summesh \mcv
( (u^{\NTk+1}_\cv )^{2}-(\ucvz)^{2} ) \\
+ \xi k\sumTsk  
\Bigl( \sumE  \transcvcvv  \vert \pcv - \pcvv \vert
(u^n_\cv  - u^n_\cvv )^{2} + 
 \summesh   G^{(+)}_\cv (u^n_\cv  - \bar{u}^n_\cv )^{2} \Bigr) \\
- k\sumTsk  \summesh  
( G^{(+)}_\cv  (\bar{u}^{n}_\cv)^2 - G^{(-)}_\cv  (u^{n}_\cv)^2  ) \le  0.
\end{array}
\label{egeq13}\ee

Then, setting $C_4=\meas(\O)\norm{u_0}{L^\infty(\O)}^2 +
2 T G^{(+)}\norm{\bar{u}}{L^{\infty}(\dr\O ^+\times\Rpast)}^{2}$ which only 
depends on $\O$,
$u_0$, $T$, $g$ and $\bar{u}$,

\[
\summesh \mcv (u^{\NTk+1}_\cv )^{2} 
+ k\sumTsk  \summesh G^{(-)}_\cv  
(u^{n}_\cv)^2  \le C_4
\]

(this inequality will not be used in the sequel) and

\be
k \sumTsk  \sumE 
\transcvcvv \vert \pcv - \pcvv \vert (u^n_\cv  - u^n_\cvv )^{2} + 
k\sumTsk  \summesh 
G^{(+)}_\cv (u^n_\cv  - \bar{u}^n_\cv )^{2} \le  \frac {C_4} {\xi}.
\label{egeq17}\ee

The Cauchy-Schwarz inequality yields

\be
\begin{array}{ll}
k\sumTsk  \sumE  
\transcvcvv \vert \pcv - \pcvv \vert \vert u^n_\cv  - u^n_\cvv \vert  + 
k\sumTsk  \summesh G^{(+)}_\cv  
\vert u^n_\cv  - \bar{u}^n_\cv \vert  & \le \\
( k \sumTsk   \sumE 
\transcvcvv \vert \pcv - \pcvv \vert (u^n_\cv  - u^n_\cvv )^{2} + 
k\sumTsk   \summesh   
G^{(+)}_\cv (u^n_\cv  - \bar{u}^n_\cv )^{2}  )^{\half } \\
\Bigl( k \sumTsk  ( \sumE 
\transcvcvv \vert \pcv - \pcvv \vert + \summesh G^{(+)}_\cv) \Bigr)^{\half }
\end{array}
\label{egeq18}\ee

The expression $W$, defined by $W = \sumE\transcvcvv \vert \pcv - \pcvv \vert$, verifies

\be
W \le ( \sumE \transcvcvv  )^{\half } 
( \sumE \transcvcvv (\pcv - \pcvv)^{2}) ^{\half } \le C_3  ( \sumE \transcvcvv  
)^{\half }
\label{egeq19}\ee

using \eqref{egllesterr2}. Recall that $C_3$ only depends on $g$ and $\O$.

Since

\be
\sumE\transcvcvv  \le (\sumE \medgecvcvv \dcvcvv) \frac {1} {\a^2 h^2} \le \frac {\dim
\meas(\O)} {\a^2 h^2}
\label{sumtau}
\ee

and

\[
\summesh G^{(+)}_\cv  = \dsp{\int_{\dr\O }} g^+(x) \dfrontiere(x),
\]

we finally conclude that \eqref{egeqbvw} holds.

\begin{remark} \normalfont
In the case of source terms, one adds the term $k\sumTsk  \summesh 
V^{(+)}_\cv  \vert u^n_\cv  - s^n_\cv \vert $ in the left hand side of \eqref{egeqbvw} (and
$H$ also  depends on $v$ and $s$).
\end{remark}

\subsection{Convergence result}\label{casimplec}

We already know, by the results of  section \refp{nll}, 
that the pressure scheme converges. Let us now prove 
the convergence of the saturation scheme \eqref{egschema}. 
Thanks to the estimate \eqref{egeqestlinf} in $L^{\infty 
}(\O\times\Rpast)$ (Lemma \ref{egestlinf}), for any sequence of meshes and time steps, such that 
the size of the mesh tends to $0$, we can extract a subsequence such that the 
approximate saturation converges to a function $u$ in $L^{\infty 
}(\O\times\Rpast)$ for the weak-$\star$ topology. We have to show that $u$
is the (unique) solution of \eqref{egeq01}, \eqref{egeq1} (the uniqueness of the solution
 is given by Assumption \ref{eghyp4}).

\begin{theo}\label{egcvgce}
Under assumptions \ref{eghyp1} and \ref{eghyp4}, let $\xi \in (0,1)$ and $\a >0$ be
given. For an admissible mesh $\mesh$, in the sense of Definition
\refp{meshneuman}, such that $\dedge \ge \a \ \size(\mesh)$ for all $\edge \in
\edgesint$ and for a time step $k>0$ satisfying \eqref{egeqS},
let $\uapp$ be defined by
\eqref{egschemapb}-\eqref{egschemap0} and \eqref{egschema0}-\eqref{eguapp}. Then
$\uapp$ converges to the solution $u$ of \eqref{egeq01}, \eqref{egeq1} 
in $L^{\infty }(\O\times\Rpast)$ for the
weak-$\star$ topology, as $\sizemesh \tends 0$.

\end{theo}

\proof of Theorem \ref{egcvgce}

\smallskip

In the case $g(x)=0$ for a.e. (for the ($d-1$)-dimensional Lebesgue measure) $x \in \dr \O$, the proof of
Theorem \ref{egcvgce} is easy. Indeed, $\grad p(x)= 0$ for a.e. $x \in \O$ and, for any
mesh and time step, $\pcv-\pcvv=0$ for all $\cv$, $\cvv \in \mesh$. Then, $u^n_\cv=u^0_\cv$
for all $\cv \in \mesh$ and all $\nnn$. Therefore, it is easy to prove 
that the sequence $\uapp$  converges,   as $\size(\mesh) \tends
0$ (for any $k$\ldots),  to $u$, defined by $u(x,t)=u_0(x)$ for 
a.e. $(x,t) \in \O \times \Rp$; note that  $u$ is the unique solution to \eqref{egeq01}, \eqref{egeq1}.

\smallskip

Let us now assume that $g$ is not the null function in $L^2(\dr \O)$.

\smallskip

Let $(\mesh_m,k_m)_{m\in\N}$ be a sequence of space meshes and time steps. For all $m
\in \N$, assume that $\mesh_m$ is an admissible mesh in the sense of Definition
\ref{meshneuman}, that $\dedge \ge \a \size(\mesh_m)$ for all $\edge \in
\edgesint$ and that $k_m>0$ satisfies \eqref{egeqS} (with $k=k_m$ and
$\mesh=\mesh_m$).
Assume also that $\size(\mesh_m)$
$\tends 0$ as $m \tends \infty$.

Let $u_m$ be the function $\uapp$ defined by
\eqref{egschemapb}-\eqref{egschemap0} and \eqref{egschema0}-\eqref{eguapp}, for 
$\mesh = \mesh_m$ and $k = k_m$.  By  Lemma \ref{egestlinf}, the sequence $(u_m)_{m \in
\N}$ is bounded in $L^{\infty
}(\O\times\Rpast)$.
In order to prove that the sequence $(u_m)_{m \in \N}$ converges in $L^{\infty
}(\O\times\Rpast)$ for the weak-$\star$ topology to the solution of 
\eqref{egeq01}, \eqref{egeq1}, using a classical contradiction argument, it 
is sufficient to prove that if $u_m \tends u$ in $L^{\infty 
}(\O\times\Rpast)$ for the weak-$\star$ topology then the
function $u$ is a solution of \eqref{egeq01}, \eqref{egeq1}.

\smallskip

Let us proceed in two steps. In the first step, it
is proved that $k_m \tends 0$ as $m \tends \infty$. Then, in the second step, it is proved
that the function $u$ is a solution of \eqref{egeq01}, \eqref{egeq1}.

From now on, the index ``$m$" is omitted.

\smallskip

\textit{Step 1 (proof of $k \tends 0$ as $m \tends \infty$)}

The proof that $k \tends 0$ (as $m \tends \infty$) uses \eqref{egeqS} and 
the fact that 
$\sizemesh \tends 0$.Indeed, define
\[
A_\mesh=\sumEll \medgecvcvv \vert \pcv-\pcvv \vert,
\]
and, for  $\edge \in \edgesint$, define $\chi_\edge$ from $\O \times \O$ to $\{0,1\}$ by
\[
\chi_\edge(x,y)=1, \text{ if } \edge \cap [x,y] \ne \emptyset,
\]
\[
\chi_\edge(x,y)=0, \text{ if } \edge \cap [x,y] = \emptyset.
\]
Let $\eta \in
\R^\dim \setminus \{0\}$ and $\comp \subset \O$ be a compact set such that $d(\comp,\O^c) \ge
\eta$.
Recall that $p_\mesh$ is defined by $p_\mesh(x)=\pcv$ for a.e. $x \in \cv$ and all $\cv \in
\mesh$.
For a.e. $x \in \comp$ one has
\[
\vert p_\mesh(x+\eta)-p_\mesh(x) \vert \le \sumedgesintkl \chi_\edge(x,x+\eta) \vert
\pcv-\pcvv \vert,
\]
integrating this inequality over $\comp$ yields, using $\int_\comp \chi_\edge(x,x+\eta)dx
\le \vert \eta \vert \medge$,
\be
\norm{p_\mesh(\cdot+\eta)-p_\mesh}{L^1(\comp)} \le \vert \eta \vert A_\mesh.
\label{atau}
\ee
Assume $A_\mesh \tends 0$ as $m \tends \infty$. Then, since $p_\mesh \tends p$ in
$L^1(\O)$, one deduces from \eqref{atau} that $\grad p =0$ a.e. on $\O$ which is impossible
(since $g$ is not the null function in $L^2(\dr \O)$). By the same way, it is also impossible
that $A_\mesh \tends 0$ for a subsequence. Then there exists $a >0$ 
(only depending on the sequence $(p_\mesh)_{m\in\N}$, recall that 
$p_\mesh=p_{\mesh_m}$ since we omit the index $m$) such that
$A_\mesh \ge a$ for all $m \in \N$.

Therefore, since $A_\mesh=\summesh \sumv \medgecvcvv (\pcvv-\pcv)^+ \ge a$, there exists $\cv
\in \mesh$ such that

\[
\sumv \medgecvcvv (\pcvv-\pcv)^+ \ge a \frac {\mcv} {\meas(\O)},
\]

Then, since $\transcvcvv=\medgecvcvv / \dcvcvv$ and $\dcvcvv \le 2h$,

\[
\sumv \transcvcvv (\pcvv-\pcv)^+ \ge a \frac {\mcv} {2h\meas(\O)},
\]

which yields, using  \eqref{egeqS},

\[
k \le (1-\xi) \meas(\O) \frac 2 a h.
\]

Hence $k \tends 0$ as $m \tends \infty$ (since $h \tends 0$ as $m \tends \infty$). This
concludes Step 1. 

\smallskip

\textit{Step 2 (proof of $u$ solution to \eqref{egeq1})}

Let $\varphi \in \mathcal{D}(\O^+\times\Rp)$. Let
$T>0$ such that, for all  $t>T-1$  and all $x\in \O$, $\varphi (x,t)=0$. 
Let $m\in\N$ such that $h<1$ and $k<1$ (thanks to Step 1, this is true 
for $m$ large enough).
Recall that we denote  $\mesh = \mesh_m$, $h= \size(\mesh_m)$ and $k = k_m$. Let $\NTk \in \N$ be such
that $\NTk k < T \le (\NTk + 1)k$. Multiplying equation \eqref{egschema} by 
$\varphi (\centercv,nk)$ and   summing the result on $\cv\in\mesh$ and $n\in\N$
yields 
\[
E_{1,m} + E_{2,m} = 0,
\]
with
\[
E_{1,m} =\sumTsk \sum_{\cv \in \mesh} \mcv ( u^{n+1}_\cv  - u^n_\cv  ) 
\varphi (\centercv,nk) 
\]
and
\[
E_{2,m} = - \sumTsk  k \sum_{\cv \in \mesh} \Bigl(\sumv \transcvcvv  u^n_{\cv, \cvv } (\pcvv
- \pcv) +  G^{(+)}_\cv  \bar{u}^n_\cv  - G^{(-)}_\cv  u^n_\cv \Bigr) 
\varphi (\centercv,nk).
\]

It is  shown below that

\be
\lim_{m\tends\infty} E_{1,m} = T_{1}, 
\label{egcvgce1}
\ee

where

\[
  T_{1} = -\int_{\Rp}\int_{\O }u(x,t) 
\varphi_t(x,t) dx dt -
\int_{\O }u_{0}(x) \varphi (x,0) dx,
\]

and that

\be
\lim_{m\tends\infty} E_{2,m} = T_{2},
\label{egcvgce2}\ee

where

\[
T_{2} =\int_{\Rp}\int_{\O }
u(x,t) \grad p(x)\cdot\grad \varphi (x,t) dx dt -
\int_{\Rp}\int_{\dr\O }\bar{u}(x,t) 
\varphi (x,t) g (x) \dfrontiere (x)  dt.
\]

Then, passing to the limit in $E_{1,m} + E_{2,m} = 0$ proves that $u$ is the (unique)
solution of \eqref{egeq01}, \eqref{egeq1} and concludes the proof of Theorem \ref{egcvgce}.

\smallskip

Let us first prove \eqref{egcvgce1}. Writing $E_{1,m}$ in the following way:
\[
E_{1,m} = \sum^{\NTk}_{n=1}\summesh \mcv \frac{\varphi (\centercv,(n-1)k) - \varphi (\centercv,nk)}{k}\ucvn  -
\summesh \mcv \ucvz \varphi (\centercv,0),
\]
the assertion \eqref{egcvgce1} is easily proved, in the same way as, for instance, in the proof of Theorem \refp{parnlcvgce}.

\smallskip

Let us prove now \eqref{egcvgce2}. 
To this purpose, we need auxiliary expressions, which make use of the convergence of the approximate pressure to the continuous one.
Define  $E_{3,m}$ and $E_{4,m}$ by

\[
\begin{array}{ll}
E_{3,m} & = \sumTsk k\sumE 
(u^n_\cv -u^n_\cvv )  \frac {\pcvv-\pcv} {\dcvcvv}
\int_{\edgecvcvv} \varphi (x,nk) \dfrontiere (x) \\
& + \sumTsk k \summesh  (\ucvn -\bar{u}^n_\cv ) \int_{\dr\cv\cap\dr\O} g(x) \varphi 
(x,nk) \dfrontiere(x)
\end{array}
\]

and

\[
E_{4,m} = \sumIN \int_{nk}^{(n+1)k} \big(\int_{\O }
u_{\mesh,k} (x,t) \grad p(x) \cdot \grad \varphi (x,nk) dx - 
\int_{\dr\O }\bar{u}_{\mesh,k} (x,t) \varphi (x,nk) g(x) \dfrontiere(x)\big)dt.
\]

We have $E_{4,m}\tends  T_{2}$  as $m\tends\infty$ thanks to the convergence of $\uapp $ to
$u$ in $L^\infty(\O \times \R)$ for the weak-$\star$ topology and to the convergence of
$\bar{u}_{\mesh,k}$ to $\bar{u}$ in $L^\infty(\dr\O ^+ \times \Rp)$ for the weak-$\star$
topology (the latter convergence holds also in $L^p(\dr\O ^+ \times (0,S))$ for all $1 \le p <
\infty$ and all $0<S<\infty$). 
Let us prove that $\vert E_{3,m} - E_{4,m}\vert \tends  0$  as $m\tends\infty$ (which gives
$E_{3,m}\tends  T_{2}$  as $m\tends\infty$).

using the equation satisfied by $p$ leads to

\[
\begin{array}{ll}
E_{4,m} & = \sumTsk k\sumE  (u^n_\cv -u^n_\cvv ) 
\int_{\edgecvcvv} \varphi (x,nk) \grad p(x) \cdot \ncvcvv \dfrontiere(x) \\
& + \sumTsk k \summesh  (\ucvn -\bar{u}^n_\cv ) \int_{\dr\cv\cap\dr\O} g(x) \varphi 
(x,nk) \dfrontiere(x).
\end{array}
\]

Therefore,

\[
\begin{array}{ll}
E_{3,m}-E_{4,m} & = \sumTsk k\sumE  (u^n_\cv -u^n_\cvv ) 
\int_{\edgecvcvv} (\frac {\pcvv-\pcv} {\dcvcvv}-\grad p(x) \cdot \ncvcvv) \varphi (x,nk)
\dfrontiere(x)
\\ 
 & = \sumTsk k \summesh  u^n_\cv \Bigl(\sumv \int_{\edgecvcvv} (\frac {\pcvv-\pcv}
{\dcvcvv}-\grad p(x) \cdot \ncvcvv) \varphi (x,nk) \dfrontiere(x)\Bigr).
\end{array}
\]

Using the equation satisfied by the pressure in \eqref{upfort} and the 
pressure scheme
\eqref{egschemap} yields
\[
E_{3,m}-E_{4,m} = \sumTsk k \summesh  u^n_\cv \Bigl(\sumv \int_{\edgecvcvv} (\frac {\pcvv-\pcv}
{\dcvcvv}-\grad p(x) \cdot \ncvcvv) (\varphi (x,nk)-\varphi(\centercv,nk))
\dfrontiere(x)\Bigr).
\]

Thanks to the regularity of $\varphi$ and $p$, there exists $C_5>0$, only depending on $p$,
and $C_6$, only depending on  $\varphi$, such that, for all $\edgecvcvv \in
\edgesint$,

\[
\vert \frac {\pcvv-\pcv} {\dcvcvv}-\grad p(x) \cdot \ncvcvv) \vert \le \vert
\frac {\pcvv-\pcv} {\dcvcvv}- \frac 1 {\medgecvcvv} \int_\edge \grad p(x) \cdot
\ncvcvv\dfrontiere(x) \vert + C_5h, \, \forall x \in \edgecvcvv
\]

and, for all $\cv \in \mesh$,

\[
\vert \varphi (x,nk)-\varphi(\centercv,nk) \vert \le C_6h, \, \forall x \in \cvbar, \, \forall \nnn.
\]

Thus,

\[
\begin{array}{ll}
\vert E_{3,m}-E_{4,m} \vert & \le \dsp{ \sumTsk k \summesh  \vert u^n_\cv \vert \Bigl(\sumv
\vert \transcvcvv (\pcvv-\pcv) - \int_\edgecvcvv \grad p(x) \cdot
\ncvcvv \dfrontiere(x) \vert \Bigr) C_6h}
\\
 & + \sumTsk k \summesh  \vert u^n_\cv \vert (\sumv
\medgecvcvv C_6C_5h^2),
\end{array}
\]

which leads to $\vert E_{3,m} - E_{4,m}\vert \tends  0$  as $m\tends\infty$,
using \eqref{egllesterr1}, \eqref{sumtau} and  the Cauchy-Schwarz inequality.

In order to prove that $E_{2,m}\tends  T_{2}$  as $m\tends\infty$ (which concludes the proof
of Theorem \ref{egcvgce}), let us show that $\vert E_{2,m} -  E_{3,m}\vert \tends  0$  as
$m\tends\infty$. 

We get, using \eqref{egeq3} and \eqref{egschemau}

\[
\begin{array}{ll}
E_{2,m} = &- \sumTsk k \sumE 
\transcvcvv (u^n_\cvv - u^n_\cv ) (\pcvv - \pcv)
\varphi (\centercv,nk) \\ & -
\sumTsk k\summesh 
(\bar{u}^n_\cv -u^n_\cv ) G^{(+)}_\cv \varphi (\centercv,nk).
\end{array}
\]

This yields

\be
\begin{array}{ll}
E_{3,m} - E_{2,m} = & \sumTsk k \sumE 
\transcvcvv (u^n_\cv -u^n_\cvv ) (\pcvv - \pcv) \phi ^n_{\cv, \cvv } + \\ &
\sumTsk k\summesh 
(u^n_\cv  - \bar{u}^n_\cv ) G^{(+)}_\cv  \phi ^n_\cv ,
\end{array}
\label{egeq28}\ee

where

\[
\phi ^n_{\cv, \cvv } =\frac 1 {\medgecvcvv} \int_{\edgecvcvv} \varphi 
(x,nk) \dfrontiere(x)  
- \varphi (\centercv,nk), \ \forall \cv\in\mesh, \ \forall\cvv\in\Nv
\]

and

\[
G^{(+)}_\cv \phi ^n_\cv  = \int_{\dr\cv\cap\dr\O} \varphi (x,nk) g(x)\dfrontiere(x) -
G^{(+)}_\cv \varphi (\centercv,nk).
\]

We recall that, for all $x\in\dr\O $, $\varphi (x,nk)g^+(x) = \varphi 
(x,nk)g(x)$,
by definition of $\mathcal{D}(\O^+ \times \Rp)$.
Therefore, there exists $C_7$, which only depends on 
$\varphi$,
such that $\vert \phi ^n_{\cv, \cvv }\vert \le C_7 h$ and $G^{(+)}_\cv
\vert \phi^n_{\cv} \vert \le G^{(+)}_\cv C_7 h$, for all $\cv \in \mesh$, $\cvv\in\Nv$ and
all $\nnn$. Therefore, using Lemma \ref{egestbvw}, we get $\vert E_{3,m} - E_{2,m} \vert \le
C_7 h \frac H {\sqrt{h}}$ which yields $\vert E_{2,m} -  E_{3,m}\vert \tends  0$ and
then $E_{2,m}\tends T_2$ as $m\tends\infty$. This concludes the proof of Theorem
\ref{egcvgce}.  \eof

\begin{remark} \normalfont
In the case of source terms, the convergence theorem \ref{egcvgce} still holds. There are
some minor modifications in the proof. The definitions of $E_{2,m}$, $E_{3,m}$ and
$E_{4,m}$ change. In the definition of $E_{2,m}$, the quantity $G^{(+)}_\cv  \bar{u}^n_\cv  -
G^{(-)}_\cv  u^n_\cv$ is replaced by $G^{(+)}_\cv  \bar{u}^n_\cv  - G^{(-)}_\cv  u^n_\cv +
V^{(+)}_\cv  s^n_\cv  - V^{(-)}_\cv  u^n_\cv$. In the definition of $E_{3,m}$ one adds

\[
\sumTsk k \summesh  (\ucvn -s^n_\cv ) \int_{\cv} v^+(x) \varphi 
(x,nk) dx.
\]

The quantity $E_{3,m}-E_{4,m}$ does not change and in order to prove $E_{3,m}-E_{2,m}
\tends 0$ it is sufficient to remark that there exists $C_8$, only depending on $\varphi$,
such that 
\[
\vert \int_{\cv} \varphi (x,nk) v^+(x)dx - V^{(+)}_\cv \varphi (\centercv,nk)
\vert \le V^{(+)}_\cv C_8 h.
\]

\end{remark}

\section{Boundary conditions}\label{hyperbcsys}\index{Boundary conditions for
hyperbolic equations}

In the industrial context, efficient numerical simulators are often developped after
a long  ``trial and error" procedure. The efficiency of the simulators may be
evaluated, for instance,
by the fact that the solution
satisfies some natural constraints and that it is in agreement with experimental data.
In some cases, estimates on the approximate solutions allow to obtain the convergence
of some sequences of approximate solutions as the discretization size tends to $0$.
However, it is not
easy to give the answer to the following question: ``Which problem is the
limit of the approximate solutions the unique solution to ?".

\vspace{.1cm}

This paper will focus on the problem of boundary conditions needed in the
discretization of non linear hyperbolic equations or systems of equations; this
 problem
is not yet clearly understood in many cases.
Two different cases will be presented: a two phase flow in a pipeline and
a two phase flow in a porous medium.

\subsection{A two phase flow in a pipeline}\label{pipe}

\paragraph{Description of the system}

A ``simple" model for a two phase flow in a pipeline (see \cite{fh}, for instance) leads to a $3 \times 3$ system
of conservations laws. The unknown $w$ is a function from $(0,1) \times \mathbf{R}_+$
in $\mathbf{R}^3$, solution of the following system:

\begin{equation}
w_t + (F(w))_x =0, \; x \in (0,1), \; t \in  \mathbf{R}_+,
\label{pipesys}
\end{equation}
where $(\cdot)_t$ and $(\cdot)_x$ denote the derivatives with respect to $t$ and $x$ variables.
The first two equations of (\ref{pipesys}) give the mass conservation of the 
$2$ phases (gas and liquid) and the third one is the momentum equation for the mixture. The expression of the
given function $F~:~\mathbf{R}^3 \rightarrow \mathbf{R}^3$ is quite complicated. It takes
into account thermodynamical laws and a hydrodynamical law. System (\ref{pipesys}) is
hyperbolic: for any $w\in \mathbf{R}^3$, the Jacobian matrix $DF(w)$ 
is diagonalizable in $\mathbf{R}$. The three eigenvalues
can be ordered: $\lambda_1(w) < \lambda_2(w) < \lambda_3(w)$. In real situations, the first eigenvalue, $\lambda_1(w)$ is
negative and the third, $\lambda_3(w)$, is positive (they correspond to some
``pressure waves" which are related to a ``sound velocity"). The
second eigenvalue, $\lambda_2(w)$, corresponds to some mean velocity between the two phases and
can change sign. One can also note that the field related to this second eigenvalue is
quite complicated because it is not, in general, a genuinely non linear field or
a linearly degenerate field. In petroleum engineering, the wave associated to this
second eigenvalue is a ``void fraction wave"; engineers require
a good representation of this wave in the numerical simulations.

\begin{remark}
In real situations, the function $F$ in System (\ref{pipesys})
also depends on $x$, in order to take into account, for
instance, the variation in the 
slope of the pipeline. Moreover,   some source terms have to be added to the system, 
in order to take into account, for
instance, some friction terms.
 \label{fdisc}
\end{remark}

In order to complete System (\ref{pipesys}), an initial condition is prescribed: 

\begin{equation}
w(x,0)=w_0(x), \; x \in (0,1),
\label{pipeic}
\end{equation}
and it is also necessary to give some boundary conditions. This appears to be
not so easy. Indeed, classically, a general principle is that the number of boundary conditions needs to
be equal to the number of positive eigenvalues of the Jacobian matrix at $x=0$ and to the number of
negative eigenvalues of the Jacobian matrix at $x=1$ (and these boundary conditions
have to satisfy some compatibility conditions). However, this principle is not so easy to
understand when an eigenvalue changes sign during the simulation (or in
the case of a null eigenvalue). A very interesting case is the so called ``severe slugging"  case
in  a pipeline. For this case, there are always two positive eigenvalues at $x=0$ and
two natural boundary conditions are prescribed at $x=0$, namely the fluxes of gas and
liquid; these
boundary conditions can be taken constant in time. At $x=1$, there is
one natural boundary condition, namely the pressure (which is the same for the two
phases, in this model), to be prescribed. It can also be constant in time. The
true physical solution, which is measured by experiments (and the aim is to modelize these
experiments), is
periodical in time and it appears that, at $x=1$, the first eigenvalue is always positive and the third one is always
negative but the second eigenvalue changes sign during the simulation. In the sequel, one presents
different ways to take into account the boundary conditions and one gives a convergence
result in a simplified case.

\paragraph{Discretization of the problem}\label{discsys}

In order to discretize Problem (\ref{pipesys}), (\ref{pipeic}) and some boundary conditions,
which will be introduced later, let $h=\frac 1 N$ (with $N \in \mathbf{N}^\star$) be the mesh size
and $k>0$ be the time step (assumed to be constant, for the sake of simplicity).
The discrete unknown are
the values $w_i^n \in \mathbf{R}^3$ for $i \in \{1,\ldots,N\}$ and $n \in \mathbf{N}$. 
The discretization of the initial condition leads to

\begin{equation}
w_i^0= \frac 1 h \int_{(i-1)h}^{ih} w_0(x)dx, \; i \in \{1,\ldots,N\}.
\label{pipeicd}
\end{equation}

For the computation of $w_i^n$ for $n >0$, one uses an explicit, $3$-points scheme:

\begin{equation}
\frac h k (w_i^{n+1}-w_i^n) + F_{i+\frac 1 2}^n-F_{i-\frac 1 2}^n= 0, \; i \in \{1,\ldots,N\}, \; n \in \mathbf{N}.
\label{pipesysd}
\end{equation}

For $i \in 1,\ldots,N-1$, one takes $F_{i+\frac 1 2}^n=g(w_i^n, w_{i+1}^n)$, where $g$ is the
numerical flux. It has to satisfy, in particular, the classical consistency condition, namely $g(a,a)=F(a)$,
and needs to be chosen in order to obtain some stability properties  for the numerical scheme
under a so called CFL condition on the time step
(see Sect. \ref{hyperbc1d} for the study of a scalar model).
In the case of two phase flow in a pipeline, the classical numerical fluxes such as
 the Godunov flux (see \cite{Go})
or the Roe flux (see \cite{roe2}) may not be implemented, because of computational difficulties.
 A convenient
choice is obtained with a simplified Roe flux, 
namely $g(a,b)= \frac {g(a)+g(b)} 2 + \frac 1 2
\vert  A(a,b) \vert (a-b)$, where $A(a,b)$ is some appoximation of the Jacobian matrix, 
depending on $a$
and $b$, but not satisfying the so called Roe condition, see \cite{fh}.

\begin{remark} In fact, for the simulation of a two phase flow in a pipeline, 
the magnitude  of the so-called fast 
eigenvalues, $\lambda_1$ and $\lambda_3$,  is much greater than that of $\lambda_2$;
 the choice in \cite{fh} is to use an implicit scheme
with respect to the fast eigenvalues, whereas the eingevalue $\lambda_2$, 
which corresponds to  the
void fraction wave,  is handled with an explicit  second order 
discretization, since the void fraction wave needs to be simulated precisely
 (see \cite{fh} for details).
\end{remark}
Let us now define the fluxes  $F_{\frac 1 2}^n$ and $F_{N + \frac 1 2}^n$ 
at the boundary.

\paragraph{Boundary conditions for the discretized problem}

In order to compute $F_{\frac 1 2}^n$ (and similarily $F_{N + \frac 1 2}^n$) a good way
is to know, or to determine, some artificial value $w_0^n \in \mathbf{R}^3$ (and $w_{N+1}^n \in \mathbf{R}^3$) and to take $F_{\frac 1 2}^n=g_0(w_0^n,w_1^n)$ (and $F_{N + \frac 1 2}^n=g_1(w_N^n,w_{N+1}^n)$).
The numerical fluxes $g_0$ and $g_1$ can be chosen equal to $g$, but this is not at all
necessary (see the convergence result of sections \ref{hyperbc1d} and 
\ref{hyperbcmd}); in fact,
 there are numerous situations where one should
take $g_0$ and $g_1$ different from $g$. 
Indeed, the scheme is often very sensitive to the
computation of the  boundary fluxes
 and it is often worthwhile  to use a more precise, but
also 
more expensive  numerical flux (such as the Godunov flux, for instance) for the computation
of the boundary fluxes than for the computation of the interior fluxes.
The difficulty is now to determine these artificial values,  $w_0^n$ and $w_{N+1}^n$.

\begin{remark} In some cases, the choice of $w_0^n$ and $w_{N+1}^n$ is quite easy. 
A well known example
is given by the wall-boundary condition for the Euler equations (with a perfect gas state
law or a more general state law). For the sake of simplicity, let us mention the
one-dimensional case;  the generalization to a multi-dimensional
 case is quite easy. The Euler equations may be written the form
(\ref{pipesys}), corresponding to conservation of mass, momentum and energy, 
with $w = (\rho,\rho u, E)^t$,  where $\rho$ is the density of the fluid,  $u$ its
 velocity, and $E$ its energy. The wall-boundary condition at $x=0$ 
 is $u=0$, and
 the only component to compute for the boundary condition is the second
component of $F_{\frac 1 2}^n$ which is equal here to the pressure at $x=0$ (since $u = 0$ at the wall), 
say $p_{\frac 1 2}^n$.
The value $w_1^n$ may be computed from the values $\rho_1^n$, $u_1^n$ and $p_1^n$.
 A natural choice for $w_0^n$ is to take $\rho_0^n=\rho_1^n$, $u_0^n=-u_1^n$ 
 and $p_0^n=p_1^n$. The flux $F_{\frac 1 2}^n$ (that is the value $p_{\frac 1 2}^n$) is then
obtained with $F_{\frac 1 2}^n=g_0(w_0^n,w_1^n)$ and a convenient choice of the numerical
flux $g_0$. We suggest to choose $g_0$ as  the Godunov flux (or as a linearized Godunov
flux, see \cite{bgh2} for instance).  Numerical tests which were performed in 
\cite{bgh2} show that this choice is very satisfactory, even in the difficult case of
a strong depressurization at the boundary.  These tests also show that the pressure obtained 
with   the Roe flux is not
so satisfactory and neither is  the choice $p_{\frac 1 2}^n=p_{1}^n$ which may seem natural
(in particular, in 2D simulations, using a dual mesh obtained with a 
finite element primal mesh).
\end{remark}

In most cases, however, the choice of $w_0^n$ and $w_{N+1}^n$ is not so easy. 
A possible method, which is described in \cite{egh}, is
now layed out, for a fixed $n$ and $g_0$ given:

\begin{enumerate}

\item Compute $DF(w_1^n)$, its eigenvalues $\{\lambda_1, \lambda_2,\lambda_3\}$  and
a basis of $\mathbf{R}^3$, $\{\varphi_1, \varphi_2,\varphi_3\}$, such that $DF(w_1^n) \varphi_i=
\lambda_i \varphi_i$, $i=1,2,3$.

\item Write $w_1^n$ on the basis $\{\varphi_1, \varphi_2,\varphi_3\}$, namely
$w_1^n=\alpha_1 \varphi_1 +\alpha_2 \varphi_2 + \alpha_3 \varphi_3$,

\item Let $p$ be the number of positive eigenvalues, compute 
$w_0^n=\beta_1 \varphi_1 +\beta_2 \varphi_2 + \beta_3 \varphi_3$ and
$F_{\frac 1 2}^n=g_0(w_0^n,w_1^n)$, where the three unknowns $\beta_1,$  $\beta_2,$
and  $\beta_3$ are determined by the 
 $p$ equations stating the boundary conditions 
 (note that these equations involve the components of $F_{\frac 1 2}^n$) 
 and by the $3-p$ equalities  $\beta_i=\alpha_i$ for $\lambda_i <0$.

\end{enumerate}

This method leads, at each time step, to a non linear system of
$3$ equations with $3$ unknowns (except if $\lambda_i=0$ for some $i$), 
namely $\beta_1$, $\beta_2$ and $\beta_3$; note that some compatibility conditions 
are needed in order that this non linear system has a
solution. Several variants of this method are possible.
For instance, a boundary condition may be imposed on $w_0^n$ rather than
 $F_{\frac 1 2}^n$. A similar method is, of course, possible at point $x=1$ 
 (changing the role of positive and negative eigenvalues).

\vspace{.1cm}

This method is not always satisfactory. In the case of severe slugging for
the simulation of two phase flow in a pipeline, the method seems to perform well at $x=0$, 
where the eigenvalues
$\lambda_1$ and $\lambda_2$ are always positive and the two boundary conditions
(gas and liquid fluxes) are convenient. However, at $x=1$, the second eigenvalue sometimes  becomes
negative and one needs a second boundary condition (the first one is
a condition on the pressure). A natural condition seems to be $Q_l=0$, 
where $Q_l$ is the second component
of the flux $F$, that is the liquid flux, but this condition does not lead to good results.
Other possible choices of this additional boundary condition at $x=1$ were tested and did not
give good results. A possible interpretation of this problem is the fact that
the sign of $\lambda_2$ is computed with $w_N^n$. Roughly speaking, it is ``too late"  when
$\lambda_2(w_N^n)$ becomes negative (see Sect. \ref{hyperbc1d} for the study of a simple
scalar case). Indeed,
good results (in agreement with experiments) are obtained with the unilateral condition
$Q_l \ge 0$ (whatever the sign of $\lambda_2(w_N^n)$). It consists in using the preceeding method 
(for the boundary condition at
$x=1$) and in replacing, in the numerical scheme (\ref{pipesysd}), the second component
of $F_{N + \frac 1 2}^n$ by its
positive part. Then, if $\lambda_2(w_N^n) <0$, two boundary conditions are given at $x=1$ (pressure
and $Q_l=0$) and if $\lambda_2(w_N^n) \ge 0$, one boundary condition is given at $x=1$ (pressure)
but, in (\ref{pipesysd}), the second component of $F_{N + \frac 1 2}^n$ is replaced by its
positive part.

\vspace{.1cm}

We studied in Section \refp{hyperbc1d} the sense of this  boundary condition
in  the simplified scalar case.

\subsection{Two phase flow in a porous medium}

A second example is given by the modelization of a two phase flow, oil and water (for instance),
in a porous medium. Phases are immiscible. Compressibility and capillarity effects are neglected. 
The model is obtained using the conservation of mass for each phase and Darcy's law.
This study is limited to the one dimensional case. In this case the pressure can be eliminated and
the problem is reduced to a single equation, namely (\ref{pipeeq}) with~:

\begin{equation}
f(u)=\frac {f_1(u) (\alpha + \beta f_2(u)) 
} {f_1(u) + f_2(u)}.
\label{fmp}
\end{equation}

The unknown is the saturation of one phase, say water, and is denoted by $u$.
 The quantity $\alpha$ is the total flux, which is constant in space, thanks to
 the incompressibility of the phases. One assumes also that it is constant in time and positive. 
 The quantity $\beta$ is the difference between the densities of the phases.
 The functions $f_1$ and $f_2$ are the mobilities of the phases. The function $f_1$ is nondecreasing,
 regular and satisfies $f_1(0)=0$. The function $f_2$ is nonincreasing,
 regular and satisfies $f_2(1)=0$. The function $f_1 +f_2$ is bounded from below by a
 positive number. 

 \begin{remark}
For the equivalent two or three dimensional model, the pressure cannot be eliminated
and the resulting model is a coupled system of two partial differential equations and two unknowns (pressure and saturation). The problem to which the
limit of the approximate solutions is solution is then much more complicated to
determine. See \cite{EG} for a partial study of this question.
\end{remark}

Here again, an initial condition is prescribed, namely (\ref{pipeeqic}), with
$u_0 \in L^\infty((0,1))$, $ 0 \le u_0 \le 1$ a.e.. The boundary condition will be given later.

\vspace{.1cm}

The numerical scheme is as in Sect. \ref{conv}; it is given by (\ref{pipeeqicd}) and 
(\ref{pipeeqd}) with (\ref{numf}). The choice of the numerical flux, $g$, satisfying (C1)-(C3),
is usually given, for this model, using an  ``upwinding phase by phase", that is
(see \cite{bj}, for instance)~:

\begin{equation}
\begin{array}{l}
\displaystyle  g(a,b) = \frac {f_1(a) (\alpha + \beta f_2(a)) } {f_1(a) + f_2(a)}
\text{ if } -\alpha + \beta f_1(a) \le 0 \\
\displaystyle g(a,b) = \frac {f_1(a) (\alpha + \beta f_2(b)) } {f_1(a) + f_2(b)}
\text{ if } -\alpha + \beta f_1(a) > 0.
\end{array}
\label{ap}
\end{equation}

Let us then define $f_{\frac 1 2}^n$ and $f_{N + \frac 1 2}^n$.
On considers here the case of an injection of pure water at $x=0$.
Then~:

\begin{equation}
f_{\frac 1 2}^n=\alpha, \; n \ge 0.
\label{fnpe}
\end{equation}

At $x=1$, The boundary condition is quite complicated. A simple example is
(see \cite{egv} for a more complete study):

\begin{equation}
f_{N+ \frac 1 2}^n=\frac {f_1(u_N^n) \alpha } {f_1(u_N^n)
 + f_2(u_N^n)}.
\label{fnps}
\end{equation}

Then, the approximate solution is given with (\ref{pipeeqicd})-(\ref{numf}), $g$ given
by (\ref{ap}), and (\ref{fnpe})-(\ref{fnps}).

\vspace{.3cm}

In order to prove that the approximate solutions converge, as $h$ and $k$ go to zero,
and to determine the problem which the limit of the approximate solutions is the unique solution to, one
proceeds as in Sect. \ref{simple}. One has to find $g_0$ and $g_1$
satisfying (C1)-(C3) and $\overline{u}, \overline{\overline{u}} \in L^\infty(\mathbf{R}_+)$ such that $f_{\frac 1 2}^n$ and $ f_{N+\frac 1 2}^n$, respectively defined by (\ref{fnpe}) and (\ref{fnps}),
satisfy (\ref{eqbc}). This is again performed in \cite{egv}.
The most interesting case is obtained for $\beta f_1(1) > \alpha$ and when the function $f$ is
increasing on $(0,u_M)$ and decreasing on $(u_M,1)$, as in Sect. \ref{simple}.
In fact, the main point is the existence of a unique $u_m \in (0,1)$such that $f(u_m)=f(1)=\alpha$ and
that $f$ is increasing on $[0,u_m]$ and greater or equal to $\alpha$ on $[u_m,1]$.
Then, it is quite easy to prove that (\ref{fnpe}) gives
\[
f_{\frac 1 2}^n=\alpha=g_G(u_m,u_1^n),
\]
where $g_G$ is the Godunov flux given in Sect. \ref{simple}.

\vspace{.1cm}

For the boundary condition at  $x=1$, it is possible to construct (see \cite{egv}) a function $g_1~:~[0,1]^2 \rightarrow \mathbf{R}$ satisfying (C1)-(C3) such that (\ref{fnps}) gives~:

\[
f_{N+ \frac 1 2}^n=g_1(u_N^n,1).
\]

It is now possible to use Theorem \ref{thm}.

Let $L$ be a common Lipschitz constant for $g$ (given by (\ref{ap})),
$g_G$ and $g_1$ (on $[0,1]^2$) and let $\zeta >0$.
If $k \le (1-\zeta) \frac h L$, the approximate solution $u_{h,k}$, that is the solution defined
by (\ref{pipeeqicd})-(\ref{numf}) (with $g$ given by \ref{ap}),  and by the boundary fluxes (\ref{fnpe})-(\ref{fnps}),
takes its values in $[0,1]$
and converges towards the unique solution of (\ref{hweak}) in
$L^p_{loc}([0,1] \times \mathbf{R}_+)$ for any $1\le p < \infty$, as $h \rightarrow 0$:

\begin{equation}
\begin{array}{l}
u \in L^\infty((0,1) \times (0,\infty)),
\\
\displaystyle \int_0^\infty \int_0^1 [ (u-\kappa)^{\pm}\varphi_t + 
\text{sign}_{\pm}(u-\kappa)(f(u)-f(\kappa)) \varphi_x ] dx dt
\\
\hfill \displaystyle  + M \int_0^\infty (
 u_m - \kappa)^{\pm} \varphi(0,t) dt + M \int_0^\infty (1 - 
 \kappa)^{\pm} \varphi(1,t) dt 
 \\
\hfill \displaystyle + \int_0^1 (u_0-\kappa)^{\pm}\varphi(x,0)dx \ge 0,
\\
\hfill  \forall \kappa \in [0,1], \; \forall \varphi \in C^1_c([0,1] \times 
 [0,\infty), \mathbf{R}_+),
\end{array}
\label{hweak}
\end{equation}
where $M$ is a bound for  $\vert 
f^\prime \vert$ on $[0,1]$ ($f$ is given by (\ref{fmp}). 
As in Sect. \ref{simple}. It is possible to give the sense of the boundary condition
if $u$ is regular enough.
Indeed, let $u$ be a regular solution of (\ref{hweak}). Then, $u$ satisfies the boundary conditions
in the sense given by \cite{bln}, that is~:

\[
\text{sign}(u(0,t)-u_m) (f(u(0,t)) - f(\kappa)) \le 0, \; 
\forall \kappa \in [u_m, u(0,t)], \text{ for a.e. } t \in \mathbf{R}_+,
\]
\[
\text{sign}(u(1,t)-1) (f(u(1,t)) - f(\kappa)) \ge 0, \; 
\forall \kappa \in [1, u(1,t)], \text{ for a.e. } t \in \mathbf{R}_+,
\]
with $[a,b]= \{ta +(1-t)b,$ $t \in [0,1]\}$ and $\text{sign}(s)=1$ for $s>0$,
 $\text{sign}(s)=-1$ for $s<0$, $\text{sign}(0)=0$.

This gives $u(0,t)=u_m$ or $u(0,t)= 
1$ and $u(1,t) \le u_m$ or $u(1,t) = 1$.
In particular, at $x=0$, one has $f(u(0,t))=\alpha$ (only water is injected) and, at $x=1$,
$f(u(1,t)) < \alpha$ if $u(1,t) < u_m$ (which states that there is some oil production).


\markboth{}{Bibliography}

\markboth{}{Index}
\cleardoublepage\phantomsection
\addcontentsline{toc}{chapter}{Index}
\printindex

\markboth{}{Index}


\begin{thebibliography}{00}

\addcontentsline{toc}{chapter}{Bibliography}

\markboth{}{Bibliography}

\bibitem[\textsc{Angermann}   (1996)]{angermann}
 \textsc{Angermann, L. }  (1996),
 Finite volume schemes as non-conforming Petrov-Galerkin
 approximations of primal-dual mixed formulations,  Report 181,
 Instit\"ut f\"ur Angewandte Mathematik, Universt\"at Erlangen-Nurnberg.


\
\bibitem[\textsc{Amiez} and \textsc{Gremaud} (1991)]{AG}
 \textsc{Amiez, G. } and  \textsc{P.A. Gremaud }  (1991),
On a numerical approach to Stefan-like
problems, \textit{Numer. Math.} \textbf{59}, 71-89.


\bibitem[\textsc{Angot} (1989)]{Angot}  \textsc{Angot, P.} (1989), Contribution
\`a l'etude des
transferts thermiques dans les syst\`emes complexes, application aux
composants
\'electroniques,
Thesis, Universit\'e de Bordeaux 1.


\bibitem[\textsc{Agouzal, Baranger, Maitre} and \textsc{Oudin} (1995)] {ABMO} \textsc{Agouzal, A., J. Baranger, J.-F.Maitre } and  \textsc{F. Oudin }(1995),
Connection between finite volume and mixed finite element methods
for a diffusion problem with non constant coefficients,
with application to Convection Diffusion,
\textit{East-West Journal on Numerical Mathematics.}, \textbf{3}, 4, 237-254.



\bibitem[\textsc{Arbogast, Wheeler} and \textsc{Yotov} (1997)] {arbogast} \textsc{Arbogast, T., M.F. Wheeler} and 
\textsc{I. Yotov}(1997), Mixed finite elements for elliptic problems with tensor coefficients as 
cell-centered finite differences, 
\textit{SIAM J. Numer. Anal.} \textbf{34}, 2, 828-852.

\bibitem[\textsc{Attey} (1974)] {ATT} \textsc{Atthey, D.R.} (1974), A Finite
Difference Scheme for Melting Problems,
\textit{J. Inst. Math. Appl.} \textbf{13}, 353-366.


\bibitem[\textsc{Bank} and \textsc{Rose} (1986)]{BR} 
\textsc{R.E.  Bank} and  \textsc{D.J. Rose} (1986),
Error estimates for the box method, \textit{SIAM J. Numer.  Anal.},  777-790

\bibitem[\textsc{Baranger, Maitre} and \textsc{Oudin} (1996)]{BMO}
\textsc{Baranger, J., J.-F. Maitre } and  \textsc{F.  Oudin} (1996),
Connection between finite volume and mixed finite element methods, \textit{Mod\'el. Math.  Anal. Num\'er.}, 30,  3,  4, 444-465.


\bibitem{bln} Bardos, C., LeRoux, A.Y., N\'{e}d\'{e}lec, J.C. (1979): First order
quasilinear equations with boundary conditions. Comm. Partial Differential
Equations, \textbf{9}, 1017--1034


\bibitem[\textsc{Barth} (1994)]{Barth} \textsc{Barth, T.J.} (1994), Aspects of unstructured
grids and
finite volume solvers for the Euler and Navier-Stokes equations,
\textit{Von Karmann Institute Lecture}.


\bibitem[\textsc{Belmouhoub} (1996)]{belmouhoub} \textsc{Belmouhoub R.} (1996) 
Mod\'elisation
tridimensionnelle de la
gen\`ese des bassins s\'edimentaires, Thesis,  Ecole Nationale
Sup\'erieure des
Mines de Paris, 1996.

\bibitem[\textsc{Berger, Brezis} and \textsc{Rogers} (1979)]{BBR}
\textsc{Berger, A.E., H. Brezis} and  \textsc{J.C.W. Rogers}
(1979),  A Numerical Method
for Solving the Problem $u_t -  \Delta f(u) = 0$, \textit{RAIRO Numerical
Analysis},
\textbf{13}, 4, 297-312.

\bibitem[\textsc{Bertsch, Kersner} and \textsc{Peletier} (1995)]{BKP} 
\textsc{Bertsch, M.,  R. Kersner} and  \textsc{L.A. Peletier}
(1995),
Positivity  versus localization in degenerate diffusion equations,
Nonlinear Analysis
TMA, \textbf{9}, 9, 987-1008.


\bibitem[\textsc{Botta} and \textsc{Hempel} (1996)]
{bottahempel}\textsc{Botta, N.} and  \textsc{D. Hempel} (1996),
A finite volume projection method for
the numerical solution of the incompressible Navier-Stokes equations on
triangular
grids,  in:  F.
Benkhaldoun and R. Vilsmeier eds,
\textit{Finite volumes for complex applications, Problems and
Perspectives}
(Hermes, Paris), 355-363.


\bibitem{bj} Brenier,  Y., Jaffr\'e, J. (1991): Upstream differencing 
for multiphase flow in reservoir simulation.
SIAM J. Num. Ana. \textbf{28}, 685--696




\bibitem[\textsc{Brezis} (1983)]{bre}
\textsc{Brezis,  H.} (1983),
\textit{Analyse Fonctionnelle:
Th\'eorie et Applications}
(Masson, Paris).

 
\bibitem[\textsc{Brun, H\'{e}rard, Leal De Sousa} and
\textsc{Uhlmann} (1996)]
{h2}\textsc{Brun, G., J. M. H\'{e}rard, L. Leal De Sousa} and
\textsc{M. Uhlmann} (1996),
Numerical modelling of turbulent compressible flows using second order
models,
 in:  F.
Benkhaldoun and R. Vilsmeier eds,
\textit{Finite volumes for complex applications, Problems and
Perspectives}
(Hermes, Paris),  338-346.

\bibitem[\textsc{Buffard, Gallou\"et} and
\textsc{H\'erard} (1998)]
{bgh}\textsc{Buffard, T., T. Gallou\"et} and
\textsc{J. M. H\'{e}rard} (1998),
Un sch\'ema simple pour les \'equations de Saint-Venant,
\textit{C. R. Acad. Sci., S\'erie I},\textbf{326}, 386-390.



\bibitem{bgh2} Buffard, T., Gallou\"{e}t, T., H\'erard, J.M. (2000): A sequel to a rough Godunov scheme : Application to real gases. 
Computers and Fluids,  \textbf{29}, 813--847

\bibitem[\textsc{Cai} (1991)]{Cai} \textsc{Cai, Z.} (1991), 
On the finite volume element method, 
\textit{Numer. Math.}, \textbf{58}, 713-735.  



\bibitem[\textsc{Cai, Mandel} and \textsc{Mc Cormick} (1991)]{CMMC} \textsc{Cai, Z., J. Mandel} and 
\textsc{S. Mc Cormick}
(1991), The finite volume
element method for diffusion equations on general triangulations, \textit{SIAM
J. Numer. Anal.}, \textbf{28}, 2, 392-402. 



\bibitem[\textsc{Chainais-Hillairet} (1996)]{Chainais1}  
\textsc{Chainais-Hillairet, C.} (1996), First and
second order schemes for a hyperbolic equation: convergence and error estimate,
in:  F. Benkhaldoun and R. Vilsmeier eds,
\textit{Finite volumes for complex applications, Problems and
Perspectives}
(Hermes, Paris),  137-144.

\bibitem[\textsc{Chainais-Hillairet} (1999)]{Chainais2}  
\textsc{Chainais-Hillairet, C.} (1999), 
Finite volume schemes for a nonlinear hyperbolic equation. Convergence towards the entropy solution and error
estimate,   \textit{Mod\'el. Math.  Anal. Num\'er.}, 33, 1, 129-156.

\bibitem[\textsc{Champier} and \textsc{Gallou\"{e}t} (1992)]{cg}  
\textsc{Champier, S.} and  \textsc{T. Gallou\"{e}t} (1992),
Convergence d'un sch\'{e}ma
d\'{e}centr\'{e} amont pour une \'{e}quation hyperbolique lin\'{e}aire
sur
un maillage triangulaire, \textit{Mod\'el. Math.  Anal. Num\'er.}
\textbf{26},   7, 835-853.

\bibitem[\textsc{Champier, Gallou\"{e}t}  and  \textsc{Herbin} (1993)]{cgh}  \textsc{Champier, S.,  T. Gallou\"{e}t}  and  \textsc{R.
Herbin} (1993), Convergence of
an Upstream finite volume Scheme on a Triangular Mesh for a Nonlinear
Hyperbolic Equation, \textit{Numer. Math. } 66, 139-157.

\bibitem[\textsc{Chavent} and \textsc{Jaffr\'e} (1990)]{CJ}  
\textsc{Chavent, G.} and  \textsc{J.  Jaffr\'e} (1990),
Mathematical
Models and Finite Element for Reservoir Simulation, \textit{Studies in
Mathematics and its Applications} (North Holland, Amsterdam).

\bibitem[\textsc{Chevrier} and \textsc{Galley} (1993)]{Che}  
\textsc{Chevrier, P.} and  \textsc{Galley, H.} (1993),
A Van Leer finite volume scheme for the Euler
equations on unstructured meshes,
\textit{Mod\'el. Math.  Anal. Num\'er.}, \textbf{27}, 2, 183-201.


\bibitem[\textsc{Chou} and \textsc{Li} (2000)]{chouli}  
\textsc{Chou, S. } and \textsc{Q. Li} Error estimates in $L^2$, $H^1$, and $L^\infty$ in covolume methods for 
elliptic and parabolic problems: a unified approach, 
\textit{Math. Comp.}, \textbf{69}, 2000, 103-120.

\bibitem[\textsc{Ciarlet} (1978)]{ciarlet}  
\textsc{Ciarlet, P.G.} (1978), \textit{The Finite Element
Method for Elliptic
Problems} (North-Holland, Amsterdam).

\bibitem[\textsc{Ciarlet} (1991)]{finitelement}  \textsc{Ciarlet, P.G.} (1991), \textit{Basic error
estimates
for elliptic problems}  in: \textit{Handbook of Numerical Analysis }\textbf{II}
(North-Holland,Amsterdam) 17-352.


\bibitem[\textsc{Ciavaldini} (1975)]{CIA}  \textsc{Ciavaldini, J.F.}  (1975),
 Analyse num\'erique d'un
probl\`eme de
Stefan  \`a deux phases par une m\'ethode d'\'el\'ements finis,
\textit{SIAM J. Numer. Anal.},
\textbf{12}, 464-488.

\bibitem[\textsc{Cockburn, Coquel} and \textsc{LeFloch} (1994)]{Coquel}
 \textsc{Cockburn, B., F. Coquel} and  \textsc{P. LeFloch} (1994), An error estimate for  finite volume
methods for multidimensional conservation laws.
\textit{Math. Comput.} \textbf{63},  207, 77-103.

\bibitem[\textsc{Cockburn, Coquel} and \textsc{LeFloch} (1995)]{Coquelconv} 
\textsc{Cockburn, B., F. Coquel} and  \textsc{P. LeFloch} (1995), Convergence of the finite
volume method  for multidimensional conservation laws,
\textit{SIAM J. Numer. Anal. } \textbf{32},  687-705.
s
\bibitem[\textsc{Cockburn} and \textsc{Gremaud} (1996)]{C} 
\textsc{Cockburn, B.} and  \textsc{P. A. Gremaud} (1996),
A priori error estimates for numerical methods for scalar
conservation laws. I. The general approach,
\textit{Math. Comput. } \textbf{65},522-554.

\bibitem[\textsc{Cockburn} and \textsc{Gremaud} (1996)]{cgr}
\textsc{Cockburn, B.} and \textsc{P. A. Gremaud} (1996),
 Error estimates for finite element methods
for scalar conservation laws,   \textit{SIAM J. Numer. Anal.} \textbf{33}, 2,
522-554.

\bibitem[\textsc{Conway} and \textsc{Smoller} (1966)]{CS}
  \textsc{Conway E.} and  \textsc{J. Smoller } (1966), Global
solutions of the Cauchy problem
for quasi-linear first-order equations in several space variables,\textit{Comm.  Pure Appl. Math.},
\textbf{19},
95-105.

\bibitem[\textsc{Coquel} and \textsc{LeFloch} (1996)]{clfmuscl}  \textsc{Coquel, F.} and \textsc{P. LeFloch} (1996), An entropy satisfying MUSCL scheme for systems 
 of conservation laws,
\textit{Numer. Math. } \textbf{74},  1-33.

\bibitem[\textsc{Cordes} and \textsc{Putti} (1998)]{Putti} 
\textsc{Cordes C.} and  \textsc{M. Putti} (1998),  Finite element
approximation of the diffusion operator on tetrahedra,
SIAM J. Sci. Comput, \textbf{19}, 4, 1154-1168.

\bibitem[\textsc{Coudi\`ere, Gallou\"et} and \textsc{Herbin} (1998)]{sobolev}
 \textsc{Coudi\`ere, Y., T. Gallou\"et} and \textsc{R. Herbin} (1998), Discrete Sobolev Inequlities and $L^p$
 error estimates for approximate finite volume solutions of
 convection diffusion equations, submitted.

\bibitem[\textsc{Coudi\`ere, Vila} and \textsc{Villedieu} (1996)]{rouencoudiere}  
 \textsc{Coudi\`ere, Y., J.P. Vila, } and  \textsc{P. Villedieu} (1996), Convergence of
a finite volume scheme for a diffusion problem,
in:  F. Benkhaldoun and R. Vilsmeier eds,
\textit{Finite volumes for complex applications, Problems and
Perspectives}
(Hermes, Paris),  161-168.


\bibitem[\textsc{Coudi\`ere, Vila} and \textsc{Villedieu} (1999)]{coudiere} \textsc{Y. Coudi\`ere, J.P. Vila} and  \textsc{P. Villedieu} (1999), Convergence rate of
a finite volume scheme for a two-dimensional convection 
diffusion problem, Math. Model. Numer. Anal. 33 (1999), no. 3, 493--516.


\bibitem[\textsc{Courbet} and  \textsc{Croisille} (1998)]{croisille} 
\textsc{Courbet, B.} and \textsc{J. P.
Croisille} (1998),
Finite-volume
box-schemes on triangular meshes, M2AN, vol 32, 5, 631-649, 1998.


\bibitem[\textsc{Crandall} and \textsc{Majda} (1980)]{CM} 
\textsc{Crandall, M.G.} and  \textsc{A.  Majda} (1980),
Monotone Difference Approximations
for  Scalar Conservation Laws,
 \textit{Math. Comput.} \textbf{34}, 149, 1-21.

\bibitem[\textsc{Dahlquist} and  \textsc{Bj\"orck} (1974)]{DB}  \textsc{Dahlquist, G.} and  \textsc{A. Bj\"orck} (1974), \textit{Numerical Methods}, Prentice Hall
Series in Automatic Computation.

\bibitem[\textsc{Deimling} (1980)]{deimling}  \textsc{Deimling, K.} (1980), 
\textit{Nonlinear Functional
Analysis}, (Springer, New York).

\bibitem[\textsc{DiPerna} (1985)]{DP}   
\textsc{DiPerna, R.} (1985), Measure-valued solutions to
conservation laws,
\textit{Arch. Rat. Mech. Anal.}, 88, 223-270.

\bibitem[\textsc{Dubois} and \textsc{LeFloch} (1988)]{dlf}   
\textsc{Dubois, F.} and \textsc{P. LeFloch} (1988), Boundary
conditions for nonlinear hyperbolic systems of conservation laws,
\textit{Journal of Differential Equations}, 71, 93-122.

\bibitem[\textsc{Eymard}(1992)]{E}   
\textsc{Eymard, R. }(1992),  Application \`{a} la simulation
de r\'{e}servoir des
m\'{e}thodes volumes-\'{e}l\'{e}ments finis; probl\`{e}mes de mise en
oeuvre, Cours
CEA-EDF-INRIA.

\bibitem[\textsc{Eymard} and \textsc{Gallou\"{e}t} (1993)]{EG}   
\textsc{Eymard, R.} and \textsc{T. Gallou\"{e}t}  (1993),
Convergence d'un sch\'ema de type
el\'ements finis - volumes finis pour un syst\`eme coupl\'e elliptique -
hyperbolique, \textit{Mod\'el. Math.  Anal. Num\'er.}
\textbf{27}, 7, 843-861.


\bibitem[\textsc{Eymard} and \textsc{Gallou\"{e}t} (1991)]{EG1}  
\textsc{Eymard, R.} and \textsc{T. Gallou\"{e}t} (1991),
Traitement des changements de
phase
dans la mod\'{e}lisation de gisements p\'{e}troliers,  \textit{Journ\'{e}es
num\'{e}riques de
Besan\c{c}on}.


\bibitem{eg} Eymard, R., Gallou\"{e}t, T. (2003):  H-convergence and numerical schemes for elliptic
equations.
SIAM Journal on Numerical Analysis,  \textbf{41}, Number 2,  539--562


\bibitem[\textsc{Eymard, Gallou\"{e}t, Ghilani} and
\textsc{Herbin} (1998)]{eggh}  \textsc{Eymard, R., T. Gallou\"{e}t,  M. Ghilani}
 and \textsc{R. Herbin} (1997),
Error estimates for the approximate solutions of a
nonlinear hyperbolic equation given by finite volume schemes,
\textit{IMA Journal of Numerical Analysis}, 18, 563-594.

 

\bibitem[\textsc{Eymard, Gallou\"{e}t,} and\textsc{Herbin} (2000)]{egh}   
\textsc{Eymard, R., T. Gallou\"{e}t, R. Herbin}  
(2000):
Finite volume methods. Handbook of numerical analysis, Vol. VII, 713--1020.
North-Holland, Amsterdam
  
\bibitem[\textsc{Eymard, Gallou\"{e}t} and  \textsc{Herbin}(1995)]{EGH2}
\textsc{Eymard, R., T. Gallou\"{e}t} and  \textsc{R. Herbin} (1995), 
Existence and uniqueness of
the entropy solution to a nonlinear hyperbolic equation, \textit{Chin. Ann.
of Math.},
\textbf{16B} 1, 1-14.


\bibitem[\textsc{Eymard, Gallou\"{e}t} and  \textsc{Herbin} (1999)]{cvnl}
\textsc{Eymard, R., T. Gallou\"{e}t} and  \textsc{R. Herbin} (1999),
Convergence  of finite volume approximations to the solutions of
semilinear convection  diffusion reaction equations, Numer. Math., 82, 91-116.


\bibitem{egv} Eymard, R., Gallou\"{e}t, T., Vovelle, J.:  Boundary conditions in the numerical approximation of some physical problems via finite volume schemes.
Accepted for publication in  ``Journal of CAM"


\bibitem[\textsc{Eymard, Ghilani} (1994)]{EG2}  \textsc{Eymard R., M.
Ghilani} (1994),
 Convergence d'un sch\'ema implicite de
type \'el\'ements finis-volumes finis pour un syst\`eme form\'e d'une
\'equation
elliptique et d'une \'equation hyperbolique, \textit{C.R. Acad. Sci. Paris
Serie I,}, \textbf{319}, 1095-1100.

\bibitem[\textsc{Faille} (1992a)]{Fa1} \textsc{Faille, I.} (1992),
Mod\'elisation bidimensionnelle de
la gen\`ese et la migration
des hydrocarbures dans un bassin s\'edimentaire,  Thesis,  Universit\'e
de Grenoble.

\bibitem[\textsc{Faille}  (1992b)]{Fa2}\textsc{Faille, I.}  (1992),  A
control volume method to solve an
elliptic equation on a
two-dimensional irregular meshing, \textit{Comp. Meth. Appl. Mech. Engrg.}, 100,
275-290.


\bibitem[\textsc{Faille} and \textsc{Heintz\'{e}} (1998)]
{fh} \textsc{Faille, I.} and
\textsc{Heintz\'{e} E.} (1998),
A rough finite volume scheme for modeling two-phase flow in a pipeline
Computers and Fluids 28, 2,  213--241.

\bibitem[\textsc{Feistauer, Felcman} and
\textsc{Lukacova-Medvidova} (1995)]{FFL1} \textsc{Feistauer M., J. Felcman }
and
\textsc{M. Lukacova-Medvidova} (1995),
Combined finite element-finite volume solution of compressible flow,
\textit{J. Comp.  Appl.  Math.} \textbf{63}, 179-199.


\bibitem[\textsc{Feistauer, Felcman} and \textsc{Lukacova-Medvidova} (1997)]{FFL2} \textsc{Feistauer M., J. Felcman } and \textsc{M. Lukacova-Medvidova} (1997)],
On the convergence of a combined finite volume-finite element method
for nonlinear convection-
diffusion problems, \textit{Numer. Meth. for P.D.E.'s}, \textbf{13}, 163-190.


\bibitem[\textsc{Fernandez} (1989)]{fer} \textsc{Fernandez G.} (1989),
Simulation num\'{e}rique
d'\'{e}coulements
r\'{e}actifs \`{a} petit nombre de Mach, Th\`{e}se de Doctorat,
Universit\'{e} de Nice.






\bibitem[\textsc{Fezoui, Lanteri, Larrouturou} and
\textsc{Olivier} (1989)]{La} \textsc{Fezoui L., S. Lanteri,  B. Larrouturou}
and
\textsc{C. Olivier} (1989),
R\'esolution num\'erique des \'equations de Navier-Stokes pour un Fluide
Compressible en
Maillage Triangulaire, INRIA report 1033.

\bibitem[\textsc{Fiard} (1994)]{thesejmf}\textsc{Fiard, J.M.} (1994),
Mod\'elisation math\'ematique et
simulation num\'erique des piles au gaz naturel \`a oxyde solide,
th\`ese de Doctorat, Universit\'e de Chamb\'ery.


\bibitem[\textsc{Fiard}, \textsc{Herbin} (1994)]{FH1} \textsc{Fiard, J.M.},
\textsc{R. Herbin} (1994),  Comparison between finite volume finite
element
methods for the numerical simulation of  an elliptic problem arising in
electrochemical engineering, Comput. Meth. Appl. Mech. Engin., 115,
315-338.


\bibitem[\textsc{Forsyth} (1989)]{F1} \textsc{Forsyth, P.A.} (1989),  A
control volume finite element method for local
mesh refinement, SPE 18415, 85-96.

\bibitem[\textsc{Forsyth} (1991)]{F2}\textsc{Forsyth, P.A.} (1991),  A control
volume finite element approach to NAPL
groundwater
contamination, SIAM J. Sci. Stat. Comput., 12, 5, 1029-1057.

\bibitem[\textsc{Forsyth} and \textsc{Sammon} (1988)]{FS} \textsc{Forsyth,
P.A.} and \textsc{P.H. Sammon}  (1988), Quadratic Convergence for
Cell-Centered Grids, Appl. Num. Math. 4, 377-394.

\bibitem[\textsc{Gallou\"{e}t} (1996)]{rouentg} \textsc{Gallou\"{e}t, T.}
(1996), Rough schemes for systems for complex
hyperbolic systems, in:  F. Benkhaldoun and R. Vilsmeier eds,
\textit{Finite volumes for complex applications, Problems and
Perspectives}
(Hermes, Paris), 10-21.

\bibitem[\textsc{Gallou\"{e}t} and \textsc{Herbin} (1994)]{GH}   \textsc{Gallou\"{e}t, T.} and  \textsc{R. Herbin} (1994), A Uniqueness Result for
Measure
Valued
Solutions of a Nonlinear Hyperbolic Equations. Int.
Diff. Int. Equations., 6,6, 1383-1394.



\bibitem[\textsc{Gallou\"{e}t}, \textsc{Herbin}  and \textsc{Vignal}(1999)]{ghv}
\textsc{Gallou\"{e}t, T.},  \textsc{R. Herbin}  and  \textsc{M.H. 
Vignal}, (1999), Error estimate for the approximate finite volume solutions of convection 
diffusion 
equations with general boundary conditions,   SIAM J. 
Numer. Anal., 37, 6, 1935-1972, 2000.  




\bibitem[\textsc{Girault} and \textsc{Raviart} (1986)]{giraultraviart} \textsc{Girault,
V.} and  \textsc{P. A.  Raviart}(1986), \textit{Finite
Element Approximation of the Navier-Stokes Equations} (Springer-Verlag).


\bibitem[\textsc{Ghidaglia, Kumbaro}
and \textsc{Le Coq} (1996)]{Ghidaglia}   \textsc{Ghidaglia J.M., A.
Kumbaro}
and \textsc{G. Le Coq} (1996), Une m\'ethode
``volumes finis" \`a flux caract\'eristiques pour
la r\'esolution num\'erique de lois
de conservation, \textit{C.R. Acad. Sci. Paris S\'er. I} \textbf{332},
981-988.



\bibitem[\textsc{Godlewski} and \textsc{Raviart} (1991)]{godlewski-ellipses}
\textsc{Godlewski E.} and  \textsc{P. A. Raviart} (1991), \textit{Hyperbolic
systems of conservation laws}, Ellipses.

\bibitem[\textsc{Godlewski} and \textsc{Raviart} (1996)]{godlewski-springer} \textsc{Godlewski E.} and  \textsc{P.
A.
Raviart} (1996), \textit{Numerical approximation of hyperbolic systems
of conservation laws}, \textit{Applied
Mathematical Sciences } \textbf{118} (Springer, New York).

\bibitem[\textsc{Godunov} (1976)]{Go} \textsc{Godunov S.} (1976), \textit{R\'{e}solution num\'{e}rique
des probl\`{e}mes multidimensionnels de la dynamique des gaz}
(Editions de Moscou).


\bibitem[\textsc{Guedda,  Hilhorst} and  \textsc{Peletier} (1997)]{GHP} \textsc{Guedda, M., D. Hilhorst} and  \textsc{M.A. Peletier} (1997),
Disappearing
interfaces in nonlinear diffusion, \textit{Adv. Math. Sciences and Applications}, 
7,  695-710.


\bibitem[\textsc{Guo} and  \textsc{Stynes} (1997)]{guo} \textsc{Guo, W.} and  \textsc{M. Stynes} (1997),
An analysis of a cell-vertex finite volume method for a parabolic convection-diffusion
problem, \textit{Math. Comput.} \textbf{66}, 217,
105-124.

\bibitem[\textsc{Harten} (1983)]{HA}  \textsc{Harten A.} (1983), On a class
of high resolution
total-variation-stable finite-difference schemes,  \textit{J. Comput.
Phys.}
\textbf{49}, 357-393.


\bibitem[\textsc{Harten, Lax} and \textsc{Van Leer } (1983)]{hlvl}\textsc{Harten A.,  P. D. Lax} and \textsc{B. Van Leer } (1983),
On upstream  differencing and
Godunov-type schemes for hyperbolic conservations laws,
\textit{SIAM review} \textbf{25}, 35-61.


\bibitem[\textsc{Harten, Hyman} and
\textsc{Lax} (1976)]{hh}\textsc{Harten A., J. M. Hyman} and  \textsc{P. D. Lax}
(1976), On finite
difference approximations and entropy conditions, \textit{Comm. Pure Appl.
Math.}
\textbf{29}, 297-322.


\bibitem[\textsc{Heinrich} (1986)]{Heinrich}  \textsc{Heinrich B.} (1986),
Finite diference methods on irregular networks, I.S.N.M. 82(Birkhauser).

\bibitem[\textsc{Herbin} (1995)]{VF4}  \textsc{Herbin R.} (1995),
An error estimate for a finite volume scheme for a
diffusion-convection problem on a triangular mesh, \textit{Num. Meth.
P.D.E.}
\textbf{11},
165-173.

\bibitem[\textsc{Herbin} (1996)]{rouenrh}  \textsc{Herbin R.} (1996), Finite
volume methods for
diffusion convection
equations  on general meshes, \textit{in Finite volumes for complex
applications, Problems and Perspectives}, F. Benkhaldoun and R.
Vilsmeier eds,
Hermes, 153-160.


\bibitem[\textsc{Herbin} and  \textsc{Labergerie} (1997)]{HL}  \textsc{Herbin
R.} and  \textsc{O. Labergerie} (1997), Finite volume
schemes for
elliptic and elliptic-hyperbolic problems on triangular meshes,
\textit{Comp. Meth. Appl. Mech. Engin.} \textbf{147},85-103.

\bibitem[\textsc{Herbin} and  \textsc{Marchand} (1997)]{HM}  \textsc{Herbin
R.} and  \textsc{E. Marchand} (1997),
\textit{Numerical approximation of a nonlinear problem
with a Signorini condition }, Third IMACS International
Symposium On Iterative Methods In Scientific Computation, Wyoming, USA,
July
97.

\bibitem[\textsc{Idelsohn} and  \textsc{Onate} (1994)]{ide-94-fin}  \textsc{Idelsohn S.} and  \textsc{Onate E.} (1997), Finite volumes and finite element: Two good friends,
\textit{Internat. J. Numer. Methods Engrg.} \textbf{37},3323-3341.


\bibitem[\textsc{Herbin} and  \textsc{Marchand} (2001)]{HM-imajna}  \textsc{Herbin R.} and  \textsc{E. Marchand} (2001),
Finite volume approximation of a class of variational inequalities, 
 \textit{IMA Journal of Numerical Analysis}, \textbf{37}, 553-585.

\bibitem[\textsc{Hughes, Engel, Mazzei } and  \textsc{Larson} (2000)]{hug-00-con}  \textsc{Herbin R.} and  \textsc{E. Marchand} (2000),
The continuous Galerkin method is locally conservative, 
 \textit{J. Comput. Phys.}, \textbf{37}, 467-488.

 \bibitem{monograph} Kagan, A.M., Linnik, Y.V., Rao, C.R. (1973):
Characterization Problems in Mathematical Statistics. Wiley, New York


\bibitem[\textsc{Kamenomostskaja, S.L.} (1995)]{KAM}\textsc{Kamenomostskaja, S.L.} (1995), On the Stefan
problem,
\textit{Mat. Sb.} 53, 489-514 (1961 in Russian).

\bibitem[\textsc{Keller} (1971)]{keller}  \textsc{Keller H. B.} (1971), A new
difference sheme for
parabolic problems, \textit{Numerical solutions of partial differential
equations, II},  B.
Hubbard ed., (Academic Press, New-York) 327-350.


\bibitem[\textsc{Kr\"oner} (1997) ]{kronerbook}  \textsc{Kr\"oner D.} (1997)
\textit{Numerical schemes for  conservation laws
in two dimensions},
Wiley-Teubner Series Advances in Numerical Mathematics. (John Wiley and
Sons,
Ltd., Chichester; B. G. Teubner, Stuttgart).

\bibitem[\textsc{Kr\"oner} and \textsc{Rokyta} (1994)]{kroner1}  \textsc{Kr\"oner D.} and \textsc{M. Rokyta} (1994),
Convergence of
upwind finite volume schemes
on unstructured grids for scalar conservation laws
in two dimensions,   \textit{SIAM J. Numer. Anal.} \textbf{31}, 2, 324-343.

\bibitem[\textsc{Kr\"oner, Noelle} and \textsc{Rokyta}  (1995)]{kroner2}
\textsc{Kr\"oner D., S. Noelle} and \textsc{M. Rokyta}  (1995),
Convergence of higher order
upwind finite volume schemes
on unstructured grids for scalar conservation laws
in several space dimensions,\textit{Numer. Math.} \textbf{71}, 527-560,.





\bibitem[\textsc{Krushkov} (1970)]{krushkov} \textsc{Krushkov S.N.}
(1970),
First Order quasilinear equations
with several space variables, \textit{Math. USSR. Sb.} \textbf{10}, 217-243.

\bibitem[\textsc{Kumbaro} (1992)]{kum} \textsc{Kumbaro A.} (1992),
Mod\'elisation, analyse math\'ematique
et num\'erique des mod\`eles
bi-fluides d'\'ecoulement diphasique, Thesis,   Universit\'e
Paris XI Orsay.


\bibitem[\textsc{Kuznetsov} (1976)]{Kuznetsov}  \textsc{Kuznetsov} (1976),
Accuracy of some approximate methods for computing the
weak solutions of a first-order quasi-linear equation,
\textit{USSR Comput. Math. and Math. Phys.} \textbf{16},  105-119.


\bibitem[\textsc{Lady\v{z}enskaja, Solonnikov} and \textsc{Ural'ceva} (1968)]{LSU}  \textsc{Lady\v{z}enskaja, O.A., V.A. Solonnikov} and
\textsc{Ural'ceva, N.N.} (1968),
Linear and Quasilinear Equations of Parabolic Type,
\textit{Transl. of Math. Monographs} \textbf{23}.


\bibitem[\textsc{Lazarov} and \textsc{Mishev} (1996)]{rouenlazarov} \textsc{Lazarov R.D.} and \textsc{I.D. Mishev} (1996),
Finite volume methods for
reaction diffusion problems
in:  F. Benkhaldoun and R. Vilsmeier eds,
\textit{Finite volumes for complex applications, Problems and
Perspectives}
(Hermes, Paris), 233-240.

\bibitem[\textsc{Lazarov, Mishev} and  \textsc{Vassilevski} (1996)]{lazarovsiam}  \textsc{Lazarov R.D.,  I.D. Mishev} and
\textsc{P.S
Vassilevski} (1996), Finite
volume methods for convection-diffusion problems,  SIAM J. Numer. Anal.,
33,
1996, 31-55.




\bibitem[\textsc{LeVeque} (1990)]{leveque} \textsc{LeVeque, R. J. } (1990),
\textit{Numerical methods for conservation laws}
(Birkhauser verlag).



\bibitem[\textsc{Li  , Chen}, and  \textsc{Wu}(2000)]{vfchinois} 
\textsc{Li, R. } \textsc{Z. Chen } and \textsc{W. Wu}(2000),
\textit{Generalized difference methods for partial differential equations}
Marcel Dekker, New York.

\bibitem[\textsc{Lions} (1996)]{lions} \textsc{Lions, P. L. } (1996), \textit{Mathematical Topics in Fluid
Mechanics; v.1: Incompressible Models}, Oxford Lecture Series in
Mathematics \& Its Applications, No.3, Oxford Univ. Press.


\bibitem[\textsc{Mackenzie}, and  \textsc{Morton} (1992)]{mackenzie} \textsc{Mackenzie, 
J.A.}, and  \textsc{K.W. Morton} (1992),
Finite volume solutions of convection-diffusion test problems,
\textit{Math. Comput.} \textbf{60}, 201,
189-220.


\bibitem[\textsc{Manteuffel} and  \textsc{White} (1986)]{WM} \textsc{Manteuffel,
T.}, and  \textsc{A.B.White} (1986),
The numerical solution of second  order
boundary value problem on non uniform meshes,\textit{Math. Comput.} \textbf{47},
511-536.


\bibitem[\textsc{Masella} (1997)]{jmm} \textsc{Masella, J. M.} (1997),
Quelques m\'ethodes num\'eriques
pour les \'ecoulements  diphasiques bi-fluide en conduites
petroli\`eres,
Thesis, Universit\'e Paris VI.

\bibitem[\textsc{Masella, Faille}, and  \textsc{Gallou\"{e}t}
(1996)] {mfg} \textsc{Masella, J. M.,  I. Faille}, and  \textsc{T.
Gallou\"{e}t}  (1996), On a rough Godunov scheme, accepted for
publication in Intl. J. Computational Fluid Dynamics.

\bibitem[\textsc{Meirmanov} (1992)]{MEI}  \textsc{Meirmanov, A.M.} (1992),
\textit{The Stefan Problem}
(Walter de Gruyter Ed, New
York).

\bibitem[\textsc{Meyer} (1973)]{MEY}  \textsc{Meyer, G.H.} (1973),
Multidimensional Stefan Problems,
\textit{SIAM J. Num. Anal.} \textbf{10}, 522-538.


\bibitem[\textsc{Mishev} (1998) ]{mishev}
\textsc{Mishev, I.D.} (1998),   Finite volume methods on Vorono\"{\i} meshes  
\textit{Num. Meth. P.D.E.}, 14, 2, 193-212.

\bibitem[\textsc{Morton} (1996)]{morton}
\textsc{Morton, K.W.} (1996), \textit{Numerical Solutions of
Convection-Diffusion
problems} (Chapman and Hall, London).


\bibitem[\textsc{Morton}  and  \textsc{S\"uli} (1991)]{mortonsuli}  \textsc{Morton, K.W.}  and  \textsc{E. S\"uli} (1991),  Finite
volume methods and their analysis,  {IMA J. Numer. Anal.} \textbf{11},
241-260.


\bibitem[\textsc{Morton, Stynes}  and  \textsc{S\"uli} (1997)]{mortonss}  \textsc{Morton, K.W., Stynes M.}  and  \textsc{E. S\"uli} (1997),  Analysis of a cell-vertex
finite volume method for a convection-diffusion problems,  \textit{Math. Comput.} \textbf{66},
220, 1389-1406.

\bibitem[\textsc{Ne\v{c}as} (1967)]{necas} \textsc{Ne\v{c}as, J. } (1967),
\textit{Les m\'ethodes directes en th\'eorie des
\'e quations elliptiques} (Masson, Paris).

\bibitem[\textsc{Nicolaides} (1992)]{nicol1} \textsc{Nicolaides, R.A.} (1992) Analysis and convergence
of the MAC scheme I. The linear problem, \textit{SIAM J. Numer. Anal.}\textbf{29}, 6,  1579-1591.


\bibitem[\textsc{Nicolaides} (1993)]{nicolbook} 
\textsc{Nicolaides, R.A.} (1993) The Covolume Approach to Computing Incompressible
Flows, in: \textit{Incompressible computational fluid dynamics} 
M.D Gunzburger, R. A. Nicolaides eds, 295--333.
 

\bibitem[\textsc{Nicolaides} and \textsc{Wu}  (1996)]{nicol2} \textsc{Nicolaides, R.A.} and \textsc{X. Wu}  (1996)
 Analysis and convergence
of the MAC scheme II. Navier-Stokes equations ,\textit{Math. Comput.}\textbf{65}, 213,  29-44.

\bibitem[\textsc{No\"elle} (1996)]{Noelle} \textsc{No\"elle, S.} (1996) A
note on entropy inequalities
and error  estimates for higher-order accurate finite volume schemes on
irregular grids, \textit{Math. Comput.} \textbf{65}, 1155-1163.


\bibitem[\textsc{Oden} (1991)]{oden}  \textsc{Oden, J.T.} (1991), 
\textit{Finite elements:
An Introduction}  in: \textit{Handbook of Numerical Analysis }\textbf{II}
(North-Holland,Amsterdam) 3-15.

\bibitem[\textsc{Oleinik} (1960)]{Ol2} \textsc{Oleinik,  O.A.} (1960),
A method of solution of the general Stefan Problem, \textit{Sov.
Math. Dokl.} \textbf{1}, 1350-1354.


\bibitem[\textsc{Oleinik} (1963)]{Ol}  \textsc{Oleinik, O. A.} (1963), On
discontinuous solutions of nonlinear
differential equations,\textit{Am. Math. Soc. Transl., Ser. 2} \textbf{26},
95-172.

\bibitem[\textsc{Osher} (1984)]{O}  \textsc{Osher S.} (1984),  Riemann
Solvers, the entropy Condition, and
difference approximations, \textit{SIAM J. Numer. Anal.} \textbf{21}, 217-235.


\bibitem{o} Otto, F. (1996): Initial-boundary value problem for a scalar conservation
law. C. R. Acad. Sci. Paris S\'er. I Math. \textbf{8},  729--734

\bibitem[\textsc{Patankar} (1980)]{Pat} \textsc{Patankar, S.V.} (1980), \textit{Numerical Heat Transfer and
Fluid Flow}, Series
in Computational Methods in Mechanics and  Thermal Sciences, Minkowycz
and
Sparrow Eds. (Mc Graw Hill).

\bibitem{tran} Patault, S., Q.-H. Tran, Q.H. (1996): Mod\`ele et sch\'ema num\'erique du code TACITE-NPW, tech. report,
rapport IFP 42415





\bibitem[\textsc{Pfertzel} (1987)]{pfertzel}  \textsc{Pfertzel, A.} (1987),
Sur quelques sch\'emas num\'eriques pour la r\'esolution des
\'ecoulements
diphasiques en milieux poreux, Thesis, Universit\'e de Paris 6.


\bibitem[\textsc{Roberts} and \textsc{Thomas}  (1991)]{mixedfinitelement}
\textsc{Roberts J.E. } and {J.M. Thomas}  (1991),
Mixed and hybrids methods, in: \textit{Handbook of Numerical Analysis }\textbf{II}
(North-Holland,Amsterdam) 523-640.



\bibitem[\textsc{Roe} (1980)]{roe} \textsc{Roe P. L.} (1980), The use of
Riemann problem
in finite difference
schemes, \textit{Lectures notes of physics} \textbf{141}, 354-359.



\bibitem[\textsc{Roe} (1981)]{roe2} \textsc{Roe P. L.} (1981), Approximate
Riemann solvers, parameter
vectors, and difference schemes, \textit{J.  Comp. Phys.} \textbf{43},
357-372.

\bibitem[\textsc{Rudin} (1987)]{rudin} \textsc{Rudin W.} (1987), \textit{Real and
Complex analysis}
(McGraw Hill).

\bibitem[\textsc{Samarskii} (1965)]{samarski1} \textsc{Samarski A.A.} (1965), 
On monotone difference schemes for elliptic 
and parabolic equations in the case of a non-selfadjoint elliptic operator,
\textit{Zh. Vychisl. Mat. i. Mat. Fiz.} \textbf{5},
548-551 (Russian).



\bibitem[\textsc{Samarskii }(1971)]{samarski2}\textsc{Samarskii A.A.} (1971), 
\textit{Introduction to the  Theory of Difference Schemes}, Nauka, Moscow (Russian).


\bibitem[\textsc{Samarskii, Lazarov} and \textsc{Makarov}(1987)]{samarskilaz} 
\textsc{Samarskii A.A., R.D. Lazarov} and \textsc{V.L. Makarov} (1987), 
\textit{Difference Schemes for differential equations having generalized solutions},
Vysshaya Shkola Publishers,  Moscow (Russian).


\bibitem[\textsc{Sanders}  (1983)]{S}  \textsc{Sanders R.}  (1983), On the
Convergence of
Monotone Finite Difference
Schemes with Variable Spatial Differencing, \textit{Math. Comput.} \textbf{40}, 161,
91-106.

\bibitem[\textsc{Selmin}  (1993)]{Sel}  \textsc{Selmin V.}  (1993),
The node-centered finite volume approach: Bridge
between finite differences and finite elements,
\textit{Comp. Meth. in Appl. Mech. Engin.} \textbf{102}, 107-
138.
 

\bibitem[\textsc{Serre} (1996)]{serre}  \textsc{Serre D.} (1996), \textit{Syst\`emes de lois de
conservation} (Diderot).

\bibitem[\textsc{Shashkov }(1987)]{shashkov} 
\textsc{Shashkov M.}  (1996), 
\textit{Conservative Finite-difference Methods on general grids},
CRC Press, New York.


\bibitem[\textsc{Sod} (1978)]{sod} \textsc{Sod G. A.} (1978), A survey of
several finite difference
methods for systems of nonlinear hyperbolic conservation laws,
\textit{J. Comput. Phys.},
\textbf{27}, 1-31.


\bibitem[\textsc{Sonier} and \textsc{Eymard} (1993)]{ES} \textsc{Sonier F.}
and  \textsc{R. Eymard} (1993),
Mathematical and numerical properties of control-volume finite-element
scheme for reservoir simulation,  \textit{paper SPE 25267,  12th
symposium on reservoir simulation, New Orleans}.

\bibitem[\textsc{S\"uli} (1989)]{suli} \textsc{S\"uli E.} (1992), 
The accuracy of cell vertex finite volume methods on 
quadrilateral meshes, \textit{Math. Comp.} \textbf{59}, 200,
359-382.

\bibitem[\textsc{Szepessy} (1989)]{Sz} \textsc{Szepessy A.} (1989), An
existence result for scalar
conservation laws
using measure valued solutions, \textit{Comm. P.D.E.} \textbf{14}, 10,
1329-1350.

\bibitem[\textsc{Temam} (1977)]{temam} \textsc{Temam R. } (1977),  \textit{Navier-Stokes Equations}
(North Holland,
Amsterdam).


\bibitem[\textsc{Tichonov} and \textsc{Samarskii} (1962)]{samarskitich} 
\textsc{Tichonov A.N.} and \textsc{A.A. Samarskii} (1962), 
Homogeneous  difference schemes on nonuniform nets \textit{Zh. Vychisl. Mat. i. Mat. Fiz.} \textbf{2},
812-832 (Russian).

\bibitem[\textsc{Thomas} (1977)]{T}  \textsc{Thomas J.M.} (1977),  Sur
l'analyse num\'erique
des m\'ethodes d'\'el\'ements finis hybrides
et mixtes,   Thesis, Universit\'e Pierre et Marie Curie.

\bibitem[\textsc{Thom\'ee} (1991)]{finitedifference}  \textsc{Thom\'ee V.}
(1991),  Finite differences for linear
parabolic equations, in: \textit{Handbook of Numerical Analysis }\textbf{I}
(North-Holland,Amsterdam) 5-196.

\bibitem[\textsc{Turkel} (1987)]{tur} \textsc{Turkel E.} (1987),
Preconditioned methods for solving
the imcompressible
and low speed compressible equations, \textit{J. Comput. Phys.} \textbf{72},
277-298.

\bibitem[\textsc{Van Leer} (1974)]{vanleer1} \textsc{Van Leer B.} (1974),
Towards the ultimate conservative
difference scheme, II. Monotonicity and conservation combined in a
second-order scheme
\textit{J. Comput. Phys.} \textbf{14}, 361-370.

\bibitem[\textsc{Van Leer} (1977)]{vanleer2} \textsc{Van Leer B.} (1977),
Towards the ultimate conservative
difference scheme, IV. a new approach to numerical convection.
\textit{J. Comput. Phys.} \textbf{23}, 276-299.

\bibitem[\textsc{Van Leer} (1979)]{vanleer3} \textsc{Van Leer B.} (1979),
Towards the ultimate conservative
difference scheme, V,
\textit{J. Comput. Phys.} \textbf{32}, 101-136.

\bibitem[\textsc{Vanselow} (1996)]{Vanselow} \textsc{Vanselow R.} (1996),
Relations between FEM and FVM, in:  F. Benkhaldoun and R. Vilsmeier eds,
\textit{Finite volumes for complex applications, Problems and
Perspectives}
(Hermes, Paris), 217-223.


\bibitem[\textsc{Vassileski, Petrova}  and  \textsc{Lazarov} (1992)]{vass}  \textsc{Vassileski, P.S., S.I. Petrova}  and  \textsc{R.D. Lazarov} (1992),  Finite
difference schemes on triangular cell-centered grids with local refinement,  
\textit{SIAM J. Sci. Stat. Comput.} \textbf{13}, 6,
1287-1313.


\bibitem[\textsc{Vignal} (1996a)]{Vignal1} \textsc{Vignal M.H.} (1996),
Convergence of a finite volume
scheme for a system of an elliptic equation and a hyperbolic equation
\textit{Mod\'el. Math.  Anal. Num\'er.} \textbf{30}, 7, 841-872.

\bibitem[\textsc{Vignal} (1996b)]{Vignal2} \textsc{Vignal M.H.} (1996),
Convergence of Finite Volumes Schemes for an elliptic hyperbolic
system with boundary conditions,in:  F. Benkhaldoun and R. Vilsmeier
eds,
\textit{Finite volumes for complex applications, Problems and
Perspectives}
(Hermes, Paris), 145-152.

\bibitem[\textsc{Vignal} and \textsc{Verdi\`ere} (1998)]{vigver}
\textsc{Vignal M.H.} and \textsc{S. Verdi\`ere} (1998), Numerical
and theoretical study of a dual mesh method using finite volume schemes
for two phase flow problems in porous media, \textit{Numer. Math.}, 80, 4, 601-639.

\bibitem[\textsc{Vila} (1986)]{vilathese}  \textsc{Vila J.P.} (1986), Sur la
th\'eorie et
l'approximation num\'erique de probl\`emes hyperboliques non
lin\'eaires.
Application aux \'equations de saint Venant et \`a la mod\'elisation des
avalanches de neige dense, Thesis, Universit\'e Paris VI.


\bibitem[\textsc{Vila} (1994)]{Vi}  \textsc{Vila J.P.} (1994), Convergence and
error estimate in finite
volume schemes for general multidimensional conservation laws, I.
explicit monotone schemes, \textit{Mod\'el. Math.  Anal. Num\'er.}
 \textbf{28}, 3, 267-285.




\bibitem[\textsc{Vol'pert}  (1967)]{Vo} \textsc{Vol'pert A. I.}  (1967),
The spaces $BV$ and quasilinear equations,
\textit{Math. USSR Sb. 2}, 225-267.

\bibitem{v} Vovelle, J. (2002): Convergence of finite volume monotone schemes for scalar
conservation laws on bounded domains.
Num. Math., \textbf{3}, 563--596


\bibitem[\textsc{Weiser} and \textsc{Wheeler} (1988)]{weiser} \textsc{Weiser A.} and \textsc{M.F. Wheeler} (1988),
On convergence of block-centered finite-differences for elliptic problems,
\textit{SIAM J. Numer. Anal.}, \textbf{25}, 351-375.


\end{thebibliography}
\end{document}